\documentclass[leqno,11pt]{amsart}
\usepackage{amssymb,amsmath,latexsym,amsfonts,amsbsy,amsthm,mathtools}
\usepackage{mathrsfs, comment}
\usepackage[T1]{fontenc}
\usepackage[utf8]{inputenc}
\usepackage{booktabs,array,longtable,microtype}
\usepackage{graphicx,tikz,pdflscape}
\usetikzlibrary{arrows.meta,calc,fit}
\usepackage{xcolor}
\usepackage[hidelinks]{hyperref}
\newcolumntype{P}[1]{>{\raggedright\arraybackslash}p{#1}}
\newcommand{\dbla}[1]{\langle #1\rangle}
\numberwithin{equation}{section}
\newtheorem{thm}{Theorem}[section]
\newtheorem{lem}[thm]{Lemma}
\newtheorem{rem}[thm]{Remark}
\newtheorem{prop}[thm]{Proposition}
\newtheorem{cor}[thm]{Corollary}

\allowdisplaybreaks[4]
\AtBeginDocument{%
  \setlength{\abovedisplayskip}{4.5pt plus 1.5pt minus 1pt}%
  \setlength{\belowdisplayskip}{4.5pt plus 1.5pt minus 1pt}%
  \setlength{\abovedisplayshortskip}{0pt plus 1pt}%
  \setlength{\belowdisplayshortskip}{3.5pt plus 1.5pt minus 1pt}%
  \setlength{\jot}{2pt}%
}

\hypersetup{pdftitle={Stability of Plasma Echo Growth and the Sharp Sobolev Stability Threshold for the Vlasov--Poisson--Landau Equation},
pdfsubject={Echo cascades before enhanced dissipation}}
\title[Plasma echo growth and the sharp Sobolev stability threshold]{{Stability of Plasma Echo Growth and the Sharp Sobolev Stability Threshold for the Vlasov--Poisson--Landau Equation}}

\author{ Weiren Zhao}
\address[W. Zhao]{Department of Mathematics, New York University Abu Dhabi, Saadiyat Island, P.O. Box 129188, Abu Dhabi, United Arab Emirates.}
\email{zjzjzwr@126.com, wz19@nyu.edu}

\author{Ruizhao Zi}
\address[R. Zi]{School of Mathematics and Statistics, and Key Laboratory of Nonlinear Analysis \& Applications (Ministry of Education), Central China Normal University, Wuhan,  430079,  P. R. China.}
\email{rzz@ccnu.edu.cn}

\date{}
\begin{document}
\begin{abstract}
We establish stability of growing plasma echo cascades for the 
Vlasov--Poisson--Landau equation with Coulomb collisions near a global
Maxwellian. The occurrence of plasma echoes, or even a growing cascade
in a reduced model, does not by itself establish this stability. We prove
that the amplification persists under changes in the linearized operator,
including collisions and the exact evolving background coefficients,
and under nonlinear perturbations. This yields optimality of the exponent $1/3$ obtained in \cite{ref4} in the
Sobolev stability threshold. Specifically, we construct solutions with
initial Sobolev perturbations of size $\nu^{1/3}[\log(1/\nu)]^p$ whose
transport profiles exhibit norm inflation before the enhanced-dissipation
time scale $\nu^{-1/3}$. Since collisions preclude the backward Sobolev
construction used in collisionless echo arguments, we prepare the initial
data at physical time zero and control the entire forward evolution.
\end{abstract}
\maketitle
\tableofcontents

\section{Introduction}\label{sec:introduction}

\subsection{Main results}\label{sec:mainResults}

We study the Vlasov--Poisson--Landau system on
$\mathbb T^3_x\times\mathbb R^3_v$, with spatial period $2\pi$
and collision frequency $\nu>0$. In the physical time variable
$t\ge0$, {the distribution $\widetilde F$ and the self-consistent electric}
field $\widetilde E$ satisfy
\begin{equation}\label{eq:vplPhysicalTime}
\begin{cases}
\displaystyle
\partial_t\widetilde F
+v\cdot\nabla_x\widetilde F
+\widetilde E\cdot\nabla_v\widetilde F
=\nu Q(\widetilde F,\widetilde F),\\[1mm]
\displaystyle
\widetilde E(t,x)
=-\nabla_x(-\Delta_x)^{-1}
\bigl(\rho[\widetilde F](t,x)-1\bigr),\\[1mm]
\widetilde F(0,x,v)=F_{\rm in}(x,v).
\end{cases}
\end{equation}
Here $\rho[H](t,x)=\int_{\mathbb R^3}H(t,x,v)\,dv$ denotes the
spatial density associated with a distribution $H$.
The initial data have mean density one, and the inverse Laplacian
acts on functions of zero spatial mean.
The Coulomb Landau collision operator is
\begin{equation}\label{eq:landauoperator}
\begin{aligned}
Q({F_1},{F_2})(v)
&=\partial_{v_i}\int_{\mathbb R^3}\Phi_{ij}(v-v_*)
\bigl(
{F_1}(v_*)\partial_{v_j}{F_2}(v)
-{F_2}(v)\partial_{v_{*j}}{F_1}(v_*)
\bigr)\,dv_*,
\\
\Phi_{ij}(z)
&=\frac1{|z|}
\left(\delta_{ij}-\frac{z_i z_j}{|z|^2}\right).
\end{aligned}
\end{equation}

For a distribution $F$ with spatially averaged density one, let
\[
E[F]=-\nabla_x(-\Delta_x)^{-1}\left(\int_{\mathbb R^3}F\,dv-1\right),
\]
where the inverse Laplacian acts on functions of zero spatial mean.
\par
The total energy per unit spatial volume includes the electrostatic energy:
\begin{equation}\label{eq:totalEnergyDefinition}
\mathcal E({F})=\frac1{(2\pi)^3}\int_{\mathbb T^3}\!\int_{\mathbb R^3}
\frac{|v|^2}{2}{F}\,dv\,dx
+\frac1{2(2\pi)^3}\int_{\mathbb T^3}|E[{F}]|^2\,dx.
\end{equation}
Smooth solutions with sufficient velocity decay conserve mass, momentum
and this total energy. The Landau operator preserves the velocity moments
$1,v,|v|^2$, while the kinetic work of the electric field cancels the
time derivative of the electrostatic energy.
\par

For a temperature parameter $\theta>0$, define the centered
Maxwellian of unit mass by
\[
\mu_\theta(v)=(2\pi\theta)^{-3/2}e^{-|v|^2/(2\theta)}.
\]
{The Maxwellian moments satisfy}
\[
\int_{\mathbb R^3}\mu_\theta\,dv=1,\qquad
\int_{\mathbb R^3}v\mu_\theta\,dv=0,\qquad
\int_{\mathbb R^3}\frac{|v|^2}{2}\mu_\theta\,dv
=\frac32\theta.
\]
We take {$\mu =\mu_1$} as the fixed reference equilibrium. {Then {$E[\mu ]=0$} and {$\mathcal E(\mu )=3/2$}.}
In the initial-data construction below, $\theta$ will be chosen
near one so that the initial kinetic plus electric energy equals
the energy of {$\mu $}. This choice may depend on $\nu$ through the
parameters of the initial data. Once those data are prescribed,
$\theta$ is a constant independent of time, position, and velocity.
It parametrizes the initial Maxwellian and imposes no constancy condition
on the temperature of the evolving solution.

We use the standard Maxwellian-normalized perturbation:
\begin{equation}\label{eq:normalizedPhysicalPerturbation}
\widetilde f=\frac{\widetilde F-\mu }{\sqrt \mu },\qquad
\widetilde F=\mu +\sqrt \mu \,\widetilde f,\qquad
\widetilde f(0)=\frac{F_{\rm in}-\mu }{\sqrt \mu }.
\end{equation}
We use the conventional collision matrix
\begin{equation}\label{eq:maxwellianCollisionMatrix}
\sigma_{ij}=\Phi_{ij}*\mu ,\qquad
\sigma=(\sigma_{ij})_{1\le i,j\le3}.
\end{equation}
Here and below repeated velocity indices are summed, and convolution is
in velocity. The normalized Maxwellian linearization is
\begin{equation}\label{eq:normalizedLandauExplicit}
\begin{aligned}
\mathcal L f
={}&-\partial_{v_i}(\sigma_{ij}\partial_{v_j}f)
 +\left(\frac14\sigma_{ij}v_iv_j
 -\frac12\partial_{v_i}(\sigma_{ij}v_j)\right)f-\mathcal K f,\\
\mathcal K f
={}&-\mu^{-1/2}(v)\partial_{v_i}\left\{\mu (v)
 \int_{\mathbb R^3}\Phi_{ij}(v-v_*)\sqrt \mu (v_*)
 \left(\partial_{v_{*j}}f(v_*)+\frac{v_{*j}}2f(v_*)\right)\,dv_*\right\}.
\end{aligned}
\end{equation}
The normalized bilinear collision operator is
\begin{equation}\label{eq:normalizedGammaExplicit}
\begin{aligned}
\Gamma(f,g)
={}&\partial_{v_i}\left[(\Phi_{ij}*(\sqrt \mu  f))\partial_{v_j}g\right]
 -\frac12\big(\Phi_{ij}*(v_i\sqrt \mu  f)\big)\partial_{v_j}g\\
&-\partial_{v_i}\left[(\Phi_{ij}*(\sqrt \mu \,\partial_{v_j}f))g\right]
 +\frac12\big(\Phi_{ij}*(v_i\sqrt \mu \,\partial_{v_j}f)\big)g.
\end{aligned}
\end{equation}

{We construct solutions depending on $x$} only through $x_1$,
while retaining all three velocity variables. In this class,
$\widetilde E=(\widetilde E_1,0,0)$.
We identify such distributions with functions of $(x_1,v)$ and use
this active spatial variable in the norms below. The theorem is stated
in the physical time $t$; {the translated time used in the proof}
is introduced in Section~\ref{sec:mainIdeas}.
{In this class,} equation~\eqref{eq:vplPhysicalTime} is equivalently
\begin{equation}\label{eq:normalizedVPLPhysicalTime}
\begin{cases}
\displaystyle
\partial_t\widetilde f+v_1\partial_{x_1}\widetilde f
+\widetilde E_1\left(\partial_{v_1}-\frac{v_1}{2}\right)\widetilde f
-\widetilde E_1v_1\sqrt \mu +\nu\mathcal L\widetilde f
=\nu\Gamma(\widetilde f,\widetilde f),\\[1mm]
\displaystyle
\widetilde E_1=-\partial_{x_1}\widetilde\phi,\qquad
-\partial_{x_1}^2\widetilde\phi
=\int_{\mathbb R^3}\sqrt \mu \,\widetilde f\,dv.
\end{cases}
\end{equation}
The normalized free-transport profile is
\begin{equation}\label{eq:normalizedPhysicalProfile}
\widetilde g(t,z,v)=\widetilde f(t,z+tv_1,v),\qquad
\sqrt \mu \,\widetilde g(t,z,v)
=\widetilde F(t,z+tv_1,v)-\mu (v).
\end{equation}
Normalization and the transport change of variables commute because {$\mu $}
depends only on $v$. We use the superscript $\nu$ for the solution family
constructed below.
For free transport,
\[
\partial_t\widetilde f+v_1\partial_{x_1}\widetilde f=0
\quad\Longrightarrow\quad
\widetilde g(t,z,v)=\widetilde f(0,z,v).
\]
Measuring growth in this profile therefore removes the velocity
derivatives generated by free transport alone.
\par

The weighted Sobolev spaces $\mathscr H^r_\ell$ and
$H^{{\rm s}_0}_{q,\ell}$ are defined in Section~\ref{sec:notation}.
The initial norm is applied to $\widetilde f^\nu(0)$. At later times
we measure both $\widetilde g^\nu$ and the unnormalized profile
{$\sqrt \mu \,\widetilde g^\nu$}.

{Nonlinear stability and damping hold for sufficiently small
Sobolev perturbations of size $O(\nu^{1/3})$, with enhanced dissipation
on the characteristic time scale $\nu^{-1/3}$ \cite{ref4}.} We test the sharpness of this exponent by
constructing classical Vlasov--Poisson--Landau solutions whose initial
perturbations exceed this size by a logarithmic factor, tend to zero
as $\nu\to0$, and produce Sobolev norm inflation before that time scale.
\par

\begin{thm}\label{thm:logNonlinearInflation}
There is a finite integer ${\rm s}_{\log}$ such that, for every fixed integer ${\rm s}_0\ge {\rm s}_{\log}$, $p>{\rm s}_0+7/2$, $\ell\ge4$ and $0<q<1/4$, the following holds. For all sufficiently small $\nu>0$, there is a nonnegative classical solution $\widetilde F^\nu$ of the original initial-value problem \eqref{eq:vplPhysicalTime} on $[0,t_\nu]$, depending only on $x_1$ and all three velocity variables, with the mass, momentum and total energy of {$\mu $}, such that
\begin{equation}\label{eq:logNonlinearInitialTime}
{
\|\widetilde f^\nu(0)\|_{H^{{\rm s}_0}_{q,\ell}} \mathrel{{\sim}}\nu^{1/3}[\log(1/\nu)]^p,\qquad t_\nu\mathrel{{\sim}}\nu^{-1/3}[\log(1/\nu)]^{{\rm s}_0+2-p}.
}
\end{equation}
For each fixed sufficiently small $\nu>0$, the normalized profile
$\widetilde g^\nu$ from \eqref{eq:normalizedPhysicalProfile} belongs to
$C([0,t_\nu];\mathscr H^{{\rm s}_0}_\ell)$, as does
{$\sqrt \mu \,\widetilde g^\nu$}.
There is a constant $c>0$, independent of $\nu$, such that
\begin{equation}\label{eq:logNonlinearInflation}
\left\|\widetilde F^\nu(t_\nu,z+t_\nu v_1,v)-\mu (v)
\right\|_{\mathscr H^{{\rm s}_0-4}_\ell}
\ge c\nu^{-1/3}.
\end{equation}
Consequently, after decreasing $c$ if necessary,
\[
\|\widetilde g^\nu(t_\nu)\|_{\mathscr H^{{\rm s}_0-4}_\ell}
\ge c\nu^{-1/3}.
\]
There are nonzero modes $k_\nu\mathrel{{\sim}}\log(1/\nu)$ for which
\begin{equation}\label{eq:logNonlinearDensity}
\langle k_\nu,k_\nu t_\nu e_1\rangle^{{\rm s}_0-3}
|\rho_{k_\nu}[\widetilde F^\nu-\mu ](t_\nu)|\ge c\nu^{-2/3},
\end{equation}
where $\rho_k[H](t)$ denotes the $k$th-Fourier
transform in $x$ of the density of $H(t,x_1,v)$, namely,
\[
\rho_k[H](t)=\frac1{2\pi}\int_0^{2\pi}e^{-ikx_1}
\rho[H](t,x_1)\,dx_1
=\frac1{2\pi}\int_0^{2\pi}\!\int_{\mathbb R^3}
e^{-ikx_1}H(t,x_1,v)\,dv\,dx_1,\qquad k\in\mathbb Z.
\]
The constant $c>0$ is independent of $t, \nu$. 
\end{thm}

\begin{rem}
\label{rem:echoRobustness}\label{rem:inflationScope}
Growth in an echo model does not by itself establish instability of
the collisional nonlinear equation. Two difficulties are essential.
First, collisions act over a long interval on high velocity
frequencies; observing before enhanced dissipation, with
$\nu t_\nu^3\ll1$, does not alone make the collision error small
on the prepared oscillations. Second, errors must be small relative
to the amplified density at the final resonance. A global energy
comparison can lose the phase localization needed for this lower
bound. We therefore compare the collision solution operators on
{functions with localized velocity Fourier transforms and estimate
the terminal density separately.}
The linear corrections and nonlinear remainder are then controlled
relative to the same terminal density coefficient; see
Sections~\ref{sec:roadmapLinearErrors}--\ref{sec:roadmapNonlinear}.

The logarithmic power is not claimed to be optimal. The theorem
concerns the displayed high-order profile and density observables;
no inflation of a fixed low-order electric-field norm is asserted.
\end{rem}

\begin{rem}[Linear Fokker--Planck collisions]
\label{rem:linearFokkerPlanckInflation}
One obvious difficulty is the distribution dependence of the Landau collision operator, which leads to complicated linear estimates, see Propositions \ref{prop:centeredpropagation}, \ref{prop:exactCollisionPropagation}. If the collision operator is simple as the linear Fokker-Planck, which has the explicit semi-group representation, then the proof will be much simpler. 
Indeed, replacing the Landau term in \eqref{eq:vplPhysicalTime} by
$\nu(\Delta_v\widetilde F+\nabla_v\cdot(v\widetilde F))$
gives the Vlasov--Poisson--Fokker--Planck equation with the fixed linear collision operator $\nu(\Delta_v+\nabla_v\cdot(v\,\cdot))$. A parallel construction gives nonnegative
classical solutions with the mass and momentum of $\mu$ such that,
for sufficiently large fixed ${\rm s}_0$, $\ell\ge4$ and $0<q<1/4$,
\[
\|\widetilde f^\nu(0)\|_{H^{{\rm s}_0}_{q,\ell}}
\sim\nu^{1/3}[\log(1/\nu)]^{{\rm s}_0+4},
\qquad t_\nu\sim\nu^{-1/3}[\log(1/\nu)]^{-2},
\]
and the profile lower bound \eqref{eq:logNonlinearInflation} holds.
Thus norm inflation again occurs before the enhanced-dissipation
time scale.

Here the approximate background solves the transport--Fokker--Planck
equation explicitly, with collisions included. The explicit Fourier
solution formula also puts collision damping directly into the
Volterra density kernels. Initial phase mixing makes the background
electric field exponentially small in the square of a mixing time
of order $\log(1/\nu)$; this smallness absorbs the fixed powers of
$\nu^{-1}$ in the error estimates. Background-field transport can
therefore be included by a single perturbation argument, without
finite background-field corrections. The approximate background
still has an electric residual in the full VPFP equation, but this
residual is absorbed into the nonlinear remainder, so no separate
exact nonlinear background is needed. Moreover, linear collisions
produce neither variable collision-coefficient errors nor quadratic
collision terms. Since this linear Fokker--Planck operator does not conserve total
energy, the construction does not impose energy matching with $\mu$
and requires no adjustment of the initial temperature.
\end{rem}

\begin{comment}

\begin{rem}[Nonlinear collisionless norm inflation below Gevrey--3]
\label{rem:collisionlessGevreyInflation}
By modifying the proof method of this paper, one can obtain the
following result for the collisionless Vlasov--Poisson equation on
$\mathbb T\times\mathbb R$.

There are positive global classical solutions $F_N$ near any fixed
Maxwellian $\mu$, for either interaction sign, with the same total
mass as $\mu$, such that, for every fixed
$r\ge0$ and $\lambda,\beta>0$, the same sequence satisfies
\[
\sum_{j=0}^{4}\sum_{k\in\mathbb Z}\int_{\mathbb R}
\langle k,\eta\rangle^{2r}
\exp\!\left\{\frac{2\lambda\langle k,\eta\rangle^{1/3}}
{\log^\beta(e+\langle k,\eta\rangle)}\right\}
\left|\partial_\eta^j\widehat{(F_N(0)-\mu)}_k(\eta)\right|^2\,d\eta
\longrightarrow0,
\]
{whereas the corresponding self-consistent electric field satisfies}
\[
\|E_N(T_N)\|_{H_x^{30}}\longrightarrow\infty,
\qquad T_N\mathrel{{\sim}} N^2\log\log(e^e+N)\longrightarrow\infty.
\]
Here $N\to\infty$ indexes the collisionless family,
$\eta\in\mathbb R$ is the velocity Fourier variable and
$\langle k,\eta\rangle=(1+k^2+\eta^2)^{1/2}$.
In particular, the initial perturbations tend to zero in every fixed
Gevrey class of order strictly larger than $3$, with four velocity
moments. The logarithmic loss does not include the unweakened
Gevrey--3 endpoint.
\end{rem}
\par
\end{comment}

\subsection{Historical results}\label{sec:collisionlessEchoConstructions}

{Landau damping is the decay of the self-consistent electric field
through phase mixing, even though the distribution develops fine
velocity structure. Landau \cite{Landau1946} identified this mechanism
in the linearized equation near a Maxwellian. Mouhot and Villani \cite{MouhotVillani2011} proved nonlinear damping
for small analytic perturbations of suitable stable homogeneous equilibria.}

{The periodic nonlinear damping theory extends to Gevrey classes
of order strictly below $3$ \cite{BedrossianMasmoudiMouhot2016}.
Paraproduct decompositions and controlled loss of regularity in time
replace the Newton iteration of the original argument. A proof near
Penrose-stable equilibria on the torus for analytic and Gevrey data in
this range uses pointwise resolvent estimates and a nonlinear bootstrap
\cite{GrenierNguyenRodnianski2021}.} {The theory reaches the critical weighted Gevrey--3 regularity,
including nonlinear wave and scattering operators, in
\cite{IonescuPausaderWangWidmayer2024}. These positive results control the accumulation of plasma
echoes through regularity assumptions on the perturbation.}
For a broader account of phase mixing, damping and echoes, see
\cite{ref3}.

On $\mathbb R_x^n\times\mathbb R_v^n$, spatial dispersion changes
the stability problem. For screened interactions, nonlinear damping
near Penrose-stable homogeneous equilibria holds for small localized
perturbations with finite regularity, with subsequent refinements
of the regularity assumptions and decay estimates for $n\ge3$
\cite{BedrossianMasmoudiMouhot2018,HanKwanNguyenRousset2021,HuangNguyenXu2025}.
For the unscreened Vlasov--Poisson equation in three dimensions,
nonlinear damping near the Poisson equilibrium is established in
\cite{IonescuPausaderWangWidmayerPoisson2024}; an alternative proof
is given in \cite{NguyenWeiZhang2025}.
The corresponding two-dimensional problem has recently been treated
for small smooth localized perturbations of the Poisson equilibrium
with zero total charge \cite{IonescuNguyen2026}.
These whole-space results concern different interaction laws or
specific equilibria and do not directly give a periodic Sobolev
stability theorem near a Maxwellian.

The fluid analogue of Landau damping is inviscid damping: mixing by
shear flows produces decay of the nonzero horizontal modes of the
velocity. Nonlinear damping for the two-dimensional Euler equation
is established near Couette flow on $\mathbb T\times\mathbb R$
\cite{BedrossianMasmoudi2015} and in finite channels near Couette
and suitable monotone shear flows
\cite{IonescuJia2020,IonescuJia2023,MasmoudiZhao2024,HeIyerZhao2026}.
Related results for the inhomogeneous Euler equation include
perturbations of constant density near Couette flow and nonconstant
density near monotone shear flows
\cite{ChenWeiZhangZhang2025,ZhaoInhomogeneous2025}.
At lower Sobolev regularity, cat's-eye structures obstruct a general
damping theorem, both near Couette flow for homogeneous Euler
\cite{LinZengEuler2011} and near monotone shear flows for
inhomogeneous Euler \cite{ZhaoZhao2026}.

{Weak collisions provide another mechanism for controlling echoes:
phase mixing creates fine velocity scales at which collisions act
more effectively. The resulting competition between resonant growth
and enhanced dissipation makes the perturbation size relative to
$\nu$ a second parameter in the stability problem.}
{Echo suppression and Sobolev Landau damping hold for a
Vlasov--Fokker--Planck equation \cite{ref7}.} The nonlinear
{Fokker--Planck problem, the collisionless limit and the tradeoff}
between regularity and perturbation size are studied in \cite{ref5}.
{For the Coulomb Vlasov--Poisson--Landau equation, the Sobolev
stability theorem \cite{ref4} is complemented by a Gevrey result
\cite{ref6}: uniform damping and a long-time
collisionless limit hold with a smallness threshold independent of
$\nu$.} The present construction uses this
{last result to construct the exact VPL background;} the additional oscillatory
perturbation is controlled in a fixed Sobolev norm.

For the two-dimensional Navier--Stokes equation on
$\mathbb T\times\mathbb R$, nonlinear inviscid damping and enhanced
dissipation are established near Couette flow for small Gevrey data
uniformly in viscosity \cite{BedrossianMasmoudiVicol2016}, for
Sobolev perturbations of size $O(\nu^{1/3})$
\cite{MasmoudiZhaoSobolev2022,WeiZhang2023}, and with quantitative
tradeoffs between regularity and perturbation size
\cite{LiMasmoudiZhaoViscous2025}.
In finite channels, stability and damping estimates are available
under no-slip \cite{ChenLiWeiZhang2020}, Navier-slip
\cite{WeiZhangChannel2026,BedrossianHeIyerWang2025,BedrossianHeIyerWang2024},
and Robin boundary conditions \cite{HeNiuZhao2026}.

Several different mechanisms obstruct stability or damping.
The Penrose instability criterion \cite{Penrose1960} detects
exponentially growing modes of the collisionless linearization at
certain homogeneous equilibria. Nonlinear instability is known for a class of double-humped
equilibria satisfying this criterion, in a cube with specular
reflection \cite{GuoStrauss1995}. This transfers spectral instability to the
nonlinear dynamics. In the present problem, the global Maxwellian
is linearly stable; growth instead occurs along a nearby
inhomogeneous background.

A different obstruction is the persistence of nontrivial coherent
structures. For the one-dimensional periodic Vlasov--Poisson equation,
traveling BGK waves exist arbitrarily close to homogeneous equilibria
in $W^{r,p}$ for $p>1$ and $r<1+1/p$ \cite{LinZeng2011}. Such waves prevent a general nonlinear damping theorem
in these neighborhoods, but their persistence does not by itself
give norm inflation.

{Plasma echoes were identified experimentally in
\cite{GouldONeilMalmberg1967}. The mechanism used here is their
cumulative amplification.}
Arbitrarily many nonlinear density oscillations can develop from
small Sobolev perturbations of suitable fixed backgrounds for
gravitational interactions \cite{ref2}. For electrostatic
interactions, the construction uses a sequence of small backgrounds
and still smaller perturbations. The distinction matters when
comparing these results with instability near a fixed Maxwellian.
A second iterate already contains a long cascade, and time-dependent
norms with an imbalance weight control the nonlinear error.
\par

{Secondary linear instability around traveling wave-like
backgrounds is studied in \cite{ref9}.} Compact velocity Fourier
support makes their field vanish after an initial interval, so that
the background is then an exact free-transport solution and successive
resonances can be separated. The instability theorem in that work
concerns the one-dimensional linearization with zero homogeneous
background, for either interaction sign, and identifies Gevrey--3
as a linear stability threshold. Sobolev growth there can coexist
with decay of the force in $L^2$. 
%A nonlinear collisionless norm-inflation result below Gevrey--3, obtainable by modifying the proof method of this paper, is stated in Remark~\ref{rem:collisionlessGevreyInflation}.
\par

At lower regularity, weakly collisional VPFP has a $\nu^{1/2}$
stability threshold in critical weighted spaces; larger perturbations
can violate enhanced-dissipation estimates before time $\nu^{-1/3}$
\cite{ZhaoZiVPFPLow}.
The regularity and the measured growth in that result differ from
the Sobolev profile inflation studied here. Thus the instability
mechanism, the perturbation space and the observable must all be
specified when comparing stability thresholds.
\par

\subsection{Organization of the paper}\label{sec:organization}
A road map of the proof, including the density decompositions and the
corresponding distribution equations, is given in Appendix~\ref{app:roadMap}.
Section~\ref{sec:mainIdeas} presents the main ideas and the reduction to the
key estimates. Section~\ref{sec:referenceEstimates} constructs and estimates
the prepared background. Section~\ref{sec:approximateDensity} studies the
approximate density equation and its growing solutions.
Section~\ref{sec:completeLinear} establishes the linear estimates and transfers
the growth to the complete linearized equation. Finally,
Section~\ref{sec:nonlinear} treats the nonlinear problem and proves the main theorem.
\par

\subsection{Notation}\label{sec:notation}

We fix the basic conventions and function spaces used throughout the
paper. The norms below apply to the initial data, the background
comparison, and the linear and nonlinear evolution estimates.

\noindent$\bullet$\quad\emph{Variables and derivatives.}
The spatial variable is $x_1\in\mathbb T=\mathbb R/(2\pi\mathbb Z)$ and
$v\in\mathbb R^3$; a subscript $x$ on a norm denotes the active spatial
variable. Physical time is $t\ge0$. In the construction of
Section~\ref{sec:mainIdeas}, we use $s=t+t_0$ and $I=[t_0,{T_1}]$, where
$t_0>0$ is selected with the initial data and ${T_1}=t_0+t_\nu$.
{We denote the identity operator by $\mathrm{Id}$;
$I$ denotes a time interval throughout.}
\begingroup
A finite upper time endpoint in a local estimate is denoted by $T$;
$T_1$ remains the fixed observation endpoint of the construction.
\par
The free-transport coordinate is $z=x_1-tv_1=x_1-(s-t_0)v_1$,
so the transport profile of a distribution $H$ written in shifted time is
$H(s,z+(s-t_0)v_1,v)$. The roman symbols ${\rm s}_0$ and ${\rm s}$
are fixed Sobolev orders, whereas the italic $s$ denotes time.
Initial data carry the subscript ${\rm in}$; a Fourier index is appended
in the same subscript, as in $H_{{\rm in},k}$. The background marker
${\rm b}$ is a superscript on distributions, as in $F^{\rm b}_{\rm in}$
and $F^{\rm b,\nu}$. Amplitude parameters on initial families are
superscripts, as in $F_{\rm in}$ and {$W_{\rm in}$};
the chosen family uses $\delta=\delta_\nu$.
For a scalar or vector $a$, we write $\langle a\rangle=(1+|a|^2)^{1/2}$;
in particular,
\[
\langle v\rangle=(1+|v|^2)^{1/2},\qquad
\langle k,\eta\rangle=(1+|k|^2+|\eta|^2)^{1/2}.
\]
We use $\mathbb N_0=\{0,1,2,\ldots\}$ and $e_1=(1,0,0)$.
For $\alpha\in\mathbb N_0^3$,
$|\alpha|=\alpha_1+\alpha_2+\alpha_3$ and
$\partial_v^\alpha=\partial_{v_1}^{\alpha_1}
\partial_{v_2}^{\alpha_2}\partial_{v_3}^{\alpha_3}$.
Multi-index inequalities are componentwise. In velocity calculations,
$\partial_i=\partial_{v_i}$, and repeated velocity indices are summed
from one to three. Spatial derivatives are written explicitly.
\par\endgroup
\noindent$\bullet$\quad\emph{Perturbations and profiles.}
We reserve $\widetilde f^\nu$ for the normalized nonlinear perturbation
in physical time and $\widetilde g^\nu$ for its free-transport profile,
as defined in \eqref{eq:normalizedPhysicalPerturbation}--\eqref{eq:normalizedPhysicalProfile}.

Their shifted-time versions $f^\nu,g^\nu$ use $s=t+t_0$.
\begingroup
Uppercase letters $F,H,W$ denote distributions or perturbations before
Maxwellian normalization, including all indexed variants. Lowercase
$f,g$ denote normalized perturbations and their transport profiles.
For example, $f=(F-\mu)/\sqrt\mu$, whereas an unnormalized linear
perturbation $H$ has normalized representative $H/\sqrt\mu$.
Time-dependent densities are denoted by $\rho$; general test functions are
identified locally.
\par
\par
\par\endgroup
\noindent$\bullet$\quad\emph{Measures and inner products.}
The notation $L^2_{x,v}$ means $L^2(\mathbb T\times\mathbb R^3)$,
with the ordinary measures $dx_1\,dv$; $L^2_v$ and $L^2_x$ denote the
corresponding spaces on $\mathbb R^3$ and $\mathbb T$.
Inner products are linear in the first argument and conjugate-linear
in the second. We omit their variable subscripts, or write
$\|\cdot\|_2$ for the norm, when the arguments determine the domain.
The same conventions apply with $z$ in place of $x_1$.
For a velocity weight $\ell$, we write
\begin{equation}\label{eq:weightedInnerProduct}
\langle f,g\rangle_\ell
=\int_{\mathbb T\times\mathbb R^3}\langle v\rangle^{2\ell}
f(x_1,v)\overline{g(x_1,v)}\,dx_1\,dv.
\end{equation}
For functions of $v$ alone, the same notation means integration only
over $\mathbb R^3$; a variable subscript is added only if needed to
specify the domain. We write ${\langle f,g\rangle}$ at weight zero.
For distributional sources, the weighted pairing is extended by the
collision duality defined below. The subscript $\ell$ specifies only
the velocity weight, not a sum of spatial derivatives.
\par

\noindent$\bullet$\quad\emph{Fourier transform and density.}
We denote the frequencies of the Fourier transform in $x$ by $k\in\mathbb Z$ and use the velocity frequency
$\eta\in\mathbb R^3$, with
\begin{equation}\label{eq:fourier}
H_k(v)=\frac1{2\pi}\int_0^{2\pi}H(x_1,v)e^{-ikx_1}\,dx_1,\qquad
\widehat H_k(\eta)=\int_{\mathbb R^3}H_k(v)e^{-i\eta\cdot v}\,dv.
\end{equation}
Spatial means include the factor $(2\pi)^{-1}$. With the ordinary
space--velocity measure fixed above, Parseval's identity is
\[
\|H\|_{L^2_{x,v}}^2
=\frac1{(2\pi)^2}\sum_{k\in\mathbb Z}
\int_{\mathbb R^3}|\widehat H_k(\eta)|^2\,d\eta.
\]
\begingroup
For a distribution $H$ with finite velocity mass and first moment, set
\[
\rho[H]=\int_{\mathbb R^3}H\,dv,\qquad
j[H]=\int_{\mathbb R^3}v_1H\,dv.
\]
Thus $\rho_k[H]=\widehat H_k(0)$; for a normalized perturbation $f$,
the physical density is {$\rho[\sqrt \mu  f]$}. For any such $H$, define
the mean-zero potential and electric field by
\[
\phi_k[H]=\frac{\rho_k[H]}{k^2},\qquad
E_k[H]=-\frac{i}{k}\rho_k[H]\quad(k\ne0),\qquad
\phi_0[H]=E_0[H]=0.
\]
No zero-total-mass assumption is required. For a distribution {${F}$}
with mean density one, including $F^{\rm b}_{\rm app}$ and
$F^{\rm b,\nu}$, this gives {$-\partial_{x_1}^2\phi[{F}]=\rho[{F}]-1$}
and {$E[{F}]=E[{F}-\mu ]$}, consistently with \eqref{eq:vplPhysicalTime}.
For a time-dependent density ${\rho}$, the notation $E[{\rho}]$ uses the same
Fourier multiplier.
For $q\ge0$ and a compact interval $I'\subset I$, the density space
$\ell^1_qL^\infty(I')$ consists of time-dependent density sequences
${\rho}=({\rho}_k)_{k\ne0}$ with finite norm
\begin{equation}\label{eq:fourierDensityNorm}
\|{\rho}\|_{\ell^1_qL^\infty(I')}
=\sum_{k\ne0}\langle k\rangle^q
\|{\rho}_k\|_{L^\infty(I')}.
\end{equation}
The time supremum is taken in each mode before summation; for continuous
functions $\rho_k$ it is the ordinary supremum. An omitted interval means $I$.
The norm $L^\infty(I';H_x^q)$ instead takes the spatial Sobolev norm
before the time supremum.
\par\endgroup
\noindent$\bullet$\quad\emph{Weighted Sobolev spaces.}
For a scalar periodic function $u$ and $q\ge0$, define
\[
\|u\|_{H_x^q}^2
=2\pi\sum_{k\in\mathbb Z}\langle k\rangle^{2q}|u_k|^2.
\]
In particular, Parseval's identity reads
$\|u\|_{L^2_x}^2=2\pi\sum_k|u_k|^2$, and $H_x^0=L^2_x$
with the same norm. For integer orders, this Fourier norm is
equivalent to the derivative-sum norm, with constants depending only
on the order. For functions of both space and velocity we use the
derivative-sum convention specified below.
For $\ell\ge0$, let
\[
L^2_\ell=L^2(\mathbb R^3,\langle v\rangle^{2\ell}dv),\qquad
\|f\|_{L^2_\ell}=\|\langle v\rangle^\ell f\|_{L^2_v}.
\]
For an integer $r\ge0$, the mixed spatial Sobolev space
$H^r_{x_1}L^2_\ell$ has norm
\[
\|f\|_{H^r_{x_1}L^2_\ell}^2
=\sum_{j=0}^r\|\langle v\rangle^\ell
\partial_{x_1}^jf\|_{L^2_{x,v}}^2.
\]
At $\ell=0$ we write $H^r_{x_1}L^2_v$.
For an integer $r\ge0$ and $\ell\ge0$, let $\mathscr H^r_\ell$ be
the completion of smooth functions, periodic in $z$ and compactly
supported in $v$, under
\begin{equation}\label{eq:Hnorm}
\|g\|_{\mathscr H^r_\ell}^2
=\sum_{j+|\alpha|\le r}\|\langle v\rangle^\ell
\partial_z^j\partial_v^\alpha g\|_{L^2_{z,v}}^2.
\end{equation}
For the initial data, $H^{{\rm s}_0}_{q,\ell}$ is defined by the
corresponding completion under
\begin{equation}\label{eq:initialGaussianNorm}
\|f\|_{H^{{\rm s}_0}_{q,\ell}}^2
=\sum_{j+|\alpha|\le {\rm s}_0}
\|e^{q|v|^2}\langle v\rangle^\ell
\partial_{x_1}^j\partial_v^\alpha f\|_2^2,
\qquad 0<q<1/4.
\end{equation}
We apply $\mathscr H^r_\ell$ both to the normalized profile
{$g^\nu$ and to the physical perturbation {$\sqrt \mu \,g^\nu$}, always}
displaying the argument. The same derivative norm in $(x_1,v)$ is
used for physical-coordinate functions. The initial norm is applied
to {$\widetilde f^\nu(0)=(F_{\rm in}-\mu )/\sqrt \mu $}.
The reference Maxwellian {$\mu =\mu_1$} is fixed throughout; auxiliary
Gaussian weights introduced later do not change this normalization.

\noindent$\bullet$\quad\emph{Collision norms.}
With the Maxwellian collision matrix $\sigma_{ij}$ from
\eqref{eq:maxwellianCollisionMatrix}, define, for $\ell\ge0$,
\begin{equation}\label{eq:backgroundCollisionNorm}
\|f\|_{\Sigma_\ell}^2
=\int_{\mathbb R^3}\langle v\rangle^{2\ell}\sigma_{ij}(v)
\left(\partial_i f\,\overline{\partial_j f}
+\frac{v_iv_j}{4}|f|^2\right)\,dv.
\end{equation}
The space $\Sigma_\ell$ is the completion of
$C_c^\infty(\mathbb R^3)$ under this norm. For a smooth velocity source
{$S$, define the dual norm through the weighted pairing:}
\begin{equation}\label{eq:quadraticDualNorm}
\|S\|_{\Sigma_\ell^*}
=\sup_{\varphi\in C_c^\infty(\mathbb R^3),\,\|\varphi\|_{\Sigma_\ell}\le1}
|\langle S,\varphi\rangle_\ell|.
\end{equation}
The completion defines the same weighted distributional dual for
nonsmooth sources. For functions depending also on $x_1$, we set
\[
\|f\|_{L^2_x\Sigma_\ell}^2
=\int_{\mathbb T}\|f(x_1,\cdot)\|_{\Sigma_\ell}^2\,dx_1,
\qquad
\|S\|_{L^2_x\Sigma_\ell^*}^2
=\int_{\mathbb T}\|S(x_1,\cdot)\|_{\Sigma_\ell^*}^2\,dx_1.
\]

\noindent$\bullet$\quad\emph{Energy and dissipation.}
For a fixed integer $r\ge0$ and a normalized distribution $f$, define
\begin{equation}\label{eq:fullLinearEnergyNorms}
\begin{aligned}
\mathcal E_r[f]&=\sum_{j=0}^r
\|\langle v\rangle^{2r+6-2j}\partial_{x_1}^jf\|_2^2,\qquad \mathcal D_r[f]=\nu\sum_{j=0}^r
\|\partial_{x_1}^jf\|_{L_x^2\Sigma_{2r+6-2j}}^2,\\
\|f\|_{\mathcal N_r^{I'}}&=\sup_{s\in I'}\mathcal E_r[f(s)]^{1/2}
+\left(\int_{I'}\mathcal D_r[f(s)]\,ds\right)^{1/2},
\qquad \mathcal N_r=\mathcal N_r^I.
\end{aligned}
\end{equation}
Thus $\mathcal D_r$ includes the collision frequency, and
the norm $\|\cdot\|_{\mathcal N_r^{I\prime}}$ controls energy and integrated dissipation.
All spatial modes, including zero, are retained.
For statements about continuity and weak convergence, $\mathscr E_r$
denotes the Hilbert completion of smooth periodic functions compactly
supported in velocity under $\mathcal E_r^{1/2}$.
Strong continuity means
${\mathcal E_r[f(s)-f(s')]^{1/2}\to0}$ as ${s\to s'}$.
The corresponding energy inner product is
\[
\langle f,g\rangle_{\mathscr E_r}
=\sum_{j=0}^r
\langle\partial_{x_1}^jf,\partial_{x_1}^jg\rangle_{2r+6-2j}.
\]
Numerical energy bounds are written using $\mathcal E_r^{1/2}$.
\par
For $r\ge2$, the comparison space $\mathscr V_r(I')$ consists of
physical perturbations with finite norm
\begin{equation}\label{eq:errorComparisonNorm}
\|H\|_{\mathscr V_r(I')}
=\|\rho[H]\|_{\ell^1_{r-1}L^\infty(I')}
+(1+{T_1})^{-1}\|H/\sqrt \mu\|_{\mathcal N_r^{I'}}.
\end{equation}

The density term sums nonzero modes, whereas the energy term retains
the zero mode. The factor $(1+T_1)^{-1}$ allows the energy and integrated
dissipation to carry one additional time-integration factor relative
to the density; it records a loss in the estimate, not a necessary
growth rate of the solution. This factor is restored when the energy
enters a nonlinear source estimate. On subintervals $I'$, the
denominator still uses the original final time. We write
$\mathscr V_r=\mathscr V_r(I)$.
\par

\noindent$\bullet$\quad\emph{Source integrals.}
A physical source $S$ enters the normalized equation as $S/\sqrt\mu$.
For the time-space $L^2$ norm with velocity norm ${\Sigma_\ell^*}$, set
\[
\|S/\sqrt\mu\|_{L^2_{s,x}(I')\Sigma_\ell^*}^{\,2}
=\int_{I'}\int_{\mathbb T}
\|S(s,x_1,\cdot)/\sqrt\mu\|_{\Sigma_\ell^*}^{\,2}\,dx_1\,ds.
\]
The two source conditions in the energy estimates are
\begin{equation}\label{eq:quadraticSourceSpaces}
\int_{I'}\mathcal E_r[S_0/\sqrt\mu]^{1/2}\,ds<\infty,
\qquad
\nu^{-1/2}\sum_{j=0}^r
\|\partial_{x_1}^j(S_1/\sqrt\mu)\|_{L^2_{s,x}(I')\Sigma_{2r+6-2j}^*}<\infty.
\end{equation}
For fixed $r$, the finite sum in the second condition is equivalent
to the square root of the sum of squares, with comparison constants
$1$ and ${\sqrt{r+1}}$. The first condition controls pairing with the
energy supremum; the second controls pairing with the integrated
collision dissipation, without requiring a velocity derivative of
the source.
For broader Gaussian normalizations, replace $S/\sqrt\mu$ by
${\mu_{{\varpi}}^{-1/2}S}$ in the same expressions.
\par

\noindent$\bullet$\quad\emph{Constants and asymptotic notation.}
The positive constants $C,c$ may change from line to line. Unless
stated otherwise, they may depend on fixed regularity orders, weights
and construction parameters, but are independent of $\nu$ and $N$.
{For nonnegative quantities, we write $X\lesssim Y$ for
$X\le CY$ and $X\sim Y$ for $cY\le X\le CY$, with constants as above.
Here $\sim$ denotes comparability; no convergence of the ratio to one
is required.} In the collisional construction, asymptotic notation
refers to $\nu\downarrow0$, with $N=K\log(1/\nu)\to\infty$ and $K$
fixed. All auxiliary regularity and expansion orders are fixed before
this limit.

\begingroup
\noindent$\bullet$ \emph{Gevrey norms.}
For fixed integers ${\rm s},L\ge0$, $\lambda>0$ and $0<\chi<1$, we use
\begin{equation}\label{eq:bggevreyinitialnorm}
\|f\|_{\mathscr G^{\lambda,\chi}_{{\rm s},L}}^2
=\sum_{j+|\alpha|\le {\rm s}}
\|\langle v\rangle^L\partial_{x_1}^j\partial_v^\alpha
 e^{\lambda\langle\nabla_{x_1},\nabla_v\rangle^\chi}f\|_2^2.
\end{equation}
The exponential denotes the Fourier multiplier with symbol
$e^{\lambda\langle k,\eta\rangle^\chi}$. The same definition in $(z,v)$
is used for transport profiles.
\par\endgroup
\addtocontents{toc}{\protect\setcounter{tocdepth}{2}}
\section{Main ideas and reduction to the key estimates}\label{sec:mainIdeas}\label{sec:proofOverview}\label{sec:keyEstimates}

We explain the construction in the order in which the comparisons
are needed, stating each key estimate with the mechanism it controls.
Section~\ref{sec:secondaryStrategy} identifies the secondary instability,
and Section~\ref{sec:roadmapBackground} prepares the initial data and
the exact background. Section~\ref{sec:echoMechanismOverview} explains
plasma echoes through a toy model and derives the initial frequency
selection. Section~\ref{sec:roadmapEcho} extracts the growing density
from the linearized VPL equation. Sections~\ref{sec:roadmapLinearErrors}
and~\ref{sec:roadmapFieldCompletion} estimate the remaining terms and
construct the complete linear solution. Section~\ref{sec:roadmapNonlinear}
constructs the nonlinear solution and proves physical-time norm
inflation. The proofs of the key estimates are given in
Sections~\ref{sec:referenceEstimates}--\ref{sec:nonlinear}.
\par

\subsection{Secondary instability and the construction strategy}\label{sec:secondaryStrategy}

The instability mechanism in this paper is a secondary instability
triggered by plasma echoes. Here secondary instability refers to the
growth of a perturbation along a
background that has itself been produced by perturbing an equilibrium.
To distinguish these two perturbations, consider an abstract equation
with a linear operator $\mathscr L$ and a bilinear operator $\mathscr B$:
\[
\partial_t F+\mathscr L F+\mathscr B(F,F)=0.
\]
Let {$F_{\rm eq}$} be an equilibrium, so that
{$\mathscr L F_{\rm eq}+\mathscr B(F_{\rm eq},F_{\rm eq})=0$}, and suppose
that the linearized evolution at this equilibrium is stable:
\[
\partial_t H+\mathscr L H
+\mathscr B(F_{\rm eq},H)+\mathscr B(H,F_{\rm eq})=0.
\]
We first choose a nonzero perturbation of {$F_{\rm eq}$} and construct the
exact solution {$F^{\rm b}(t)$} generated by those data. The difference
{$F^{\rm b}(t)-F_{\rm eq}$} is the \emph{primary perturbation}; its role
is to change the background along which a further perturbation evolves.
For a second solution {$F$}, write
\[
F=F^{\rm b}+H=F_{\rm eq}+(F^{\rm b}-F_{\rm eq})+H.
\]
The difference {$H=F-F^{\rm b}$} is the \emph{secondary perturbation}.
Subtracting the two exact equations gives
\[
\partial_t H+\mathscr L H
+\mathscr B(F^{\rm b},H)+\mathscr B(H,F^{\rm b})
=-\mathscr B(H,H).
\]
Dropping the quadratic term on the right gives the linearized evolution
along {$F^{\rm b}$}. The coefficients now contain the nonzero primary
perturbation, so stability of the evolution at {$F_{\rm eq}$} does not
by itself give a comparable bound along {$F^{\rm b}$} over long times.
The modified evolution may amplify perturbations through unstable
modes or transient growth, as in related constructions
\cite{ref2,ref1,LiMasmoudiZhao2024,ZhaoZiVPFPLow}.

In our VPL problem, the equilibrium is {$\mu $}. In physical time, the exact
background is $\widetilde F^{\rm b,\nu}$, and the two perturbations are
identified by
\[
\widetilde F^\nu-\mu 
=\underbrace{\widetilde F^{\rm b,\nu}-\mu }_{\text{primary perturbation}}
+\underbrace{\widetilde F^\nu-\widetilde F^{\rm b,\nu}}_{\text{secondary perturbation}}.
\]
The leading spatial modes $\pm1$ of the background couple adjacent
Fourier modes of the secondary perturbation and generate a sequence
of plasma echoes. We compute these interactions using an explicit
approximation and control the remaining background coefficients by
comparison. The exact background remains controlled while the
secondary perturbation grows.

The exact solution will be the sum of the prepared background, a small
multiple of a growing solution of the linearized equation, and a
nonlinear remainder. The background and full initial data {are chosen as perturbations of $\mu$,}
with the same mass, momentum and total energy as {$\mu $}.
We now specify these data before introducing the linear comparison.

\subsection{Prepared initial data and the exact background}\label{sec:roadmapBackground}

\begingroup
We first explain the choice of initially phase-mixed data at physical time zero. We then
fix the scales and select the temperatures in the initial data, and state the comparison
of the explicit approximation with the exact VPL background.
\par\endgroup
A main difficulty is to construct initial data at physical time zero
whose forward VPL evolution develops the prescribed echo amplification.
One of the main contributions of this paper is to prepare these data
and control the entire forward evolution up to the echo time.
Earlier collisionless instability constructions can prescribe a
configuration at a large positive time, adapted to the subsequent
resonant interactions, and control the solution on both sides of that
time. Deng and Masmoudi \cite{ref1} use this strategy in their
low-Gevrey instability construction for Euler near Couette flow.
Bedrossian \cite[Section~2]{ref2} prescribes a suitable configuration
at a positive time and recovers the initial datum through a nonlinear
final-value problem before analyzing the later echo cascade.
This separates the preparation of the configuration from the
subsequent growth analysis. Related final-data constructions occur
in dispersive equations: multi-soliton solutions are obtained from
approximate configurations prescribed at times tending to infinity,
with uniform control on earlier intervals
\cite{CoteMartelMerle2011,MartelMerle2006}.

The backward Sobolev estimates used in such a construction are not
available for the dissipative Landau equation considered here.
Our approach instead prepares the required phases in the initial
data. The next step explains how a time shift of an explicit
free-transport solution makes this preparation compatible with a
small comparison error throughout the forward collisional evolution.
\par

We first work in physical time $t\ge0$. The exact VPL background
$\widetilde F^{\rm b,\nu}(t)$ is not available in explicit form,
so the transfer between spatial modes and the accumulated echo
{amplification in the equation linearized about this background are difficult to compute}
directly. We therefore seek an explicit approximation whose error
is small compared with the background amplitude $a$ throughout the
echo evolution.

A natural first choice would be to start from {$\mu (v)(1+a\cos x_1)$}.
{For this preliminary calculation, denote the corresponding unshifted free-transport}
evolution by
\[
\widetilde F^{\rm b}_{\rm app,0}(t,x_1,v)
=\mu (v)\bigl(1+a\cos(x_1-tv_1)\bigr),\qquad t\ge0.
\]
Here the subscript $0$ denotes the absence of a time shift.
Direct differentiation and velocity integration give
\[
(\partial_t+v_1\partial_{x_1})\widetilde F^{\rm b}_{\rm app,0}=0,\qquad
E[\widetilde F^{\rm b}_{\rm app,0}](t,x_1)
=ae^{-t^2/2}\sin x_1.
\]
Subtracting this free-transport equation from VPL leaves the known source
\[
\underbrace{\nu Q(\widetilde F^{\rm b}_{\rm app,0},
\widetilde F^{\rm b}_{\rm app,0})}_{\text{collision term}}
-\underbrace{E[\widetilde F^{\rm b}_{\rm app,0}]
\partial_{v_1}\widetilde F^{\rm b}_{\rm app,0}}_{\text{electric term}}
\]
in the equation for the background difference. Noting that 
\[
-E[\widetilde F^{\rm b}_{\rm app,0}]\partial_{v_1}\widetilde F^{\rm b}_{\rm app,0}=-E[\widetilde F^{\rm b}_{\rm app,0}]\partial_{v_1}\big(\widetilde F^{\rm b}_{\rm app,0}-\mu\big)-E[\widetilde F^{\rm b}_{\rm app,0}]\partial_{v_1}\mu,
\]
the part of the electric source linear in $a$ is
\[
-E[\widetilde F^{\rm b}_{\rm app,0}](t)\partial_{v_1}\mu 
=ae^{-t^2/2}v_1\mu \sin x_1.
\]
The normalized time integral satisfies
\[
\int_0^\infty
\left\|\frac{E[\widetilde F^{\rm b}_{\rm app,0}](t)\partial_{v_1}\mu }{\sqrt \mu }\right\|_2\,dt
\sim a\int_0^\infty e^{-t^2/2}\,dt\sim a.
\]
Thus estimating the electric contribution by its time-integrated
norm gives only an $O(a)$ bound, with no additional small factor
relative to the background amplitude. 

To obtain the {additional} small factor, we choose $t_0\ge1$ and translate
the explicit free-transport solution by $t_0$ before prescribing
the physical initial data. The resulting approximation is
\[
\widetilde F^{\rm b}_{\rm app}(t,x_1,v)
=\widetilde F^{\rm b}_{\rm app,0}(t+t_0,x_1,v)
=\mu (v)\bigl(1+a\cos(x_1-(t+t_0)v_1)\bigr),\qquad t\ge0.
\]
{At physical time $t=0$, the phase is already} $x_1-t_0v_1$.
{The resulting field and source integral satisfy}
\begin{align}\label{small-E}
E[\widetilde F^{\rm b}_{\rm app}](t,x_1)
=ae^{-(t+t_0)^2/2}\sin x_1,\qquad
\int_0^\infty
\left\|\frac{E[\widetilde F^{\rm b}_{\rm app}](t)\partial_{v_1}\mu }{\sqrt \mu }\right\|_2\,dt
\lesssim\frac{a}{t_0}e^{-t_0^2/2}.
\end{align}
Here the time integral is a Gaussian tail starting at $t_0$.
The time shift therefore prepares the initial phase and supplies
the missing small factor.

We now introduce the shifted time variable and notation
for the exact background:
\[
s=t+t_0,\qquad
F^{\rm b,\nu}(s,x_1,v)
=\widetilde F^{\rm b,\nu}(s-t_0,x_1,v),\qquad s\ge t_0.
\]
Thus $s=t_0$ denotes physical time $t=0$, and the approximation becomes
\begin{equation}\label{eq:bgapp}
F^{\rm b}_{\rm app}(s,x_1,v)
=\widetilde F^{\rm b}_{\rm app}(s-t_0,x_1,v)
=\mu (v)\bigl(1+a\cos(x_1-sv_1)\bigr).
\end{equation}
\begingroup
The choice of initially phase-mixed data comes from translating the explicit solution
before choosing the initial data, not from merely relabeling time
in a fixed initial-value problem.
\par\endgroup
This preparation reduces the field without reducing the oscillatory
coefficients responsible for echoes. Indeed,
\[
(F^{\rm b}_{\rm app}-\mu )_{\pm1}(s,v)=\frac a2\mu (v)e^{\mp isv_1},
\qquad
\rho_{\pm1}[F^{\rm b}_{\rm app}](s)=\frac a2e^{-s^2/2}.
\]
The distribution coefficients have modulus {$a\mu /2$}; the Gaussian
decay arises only after integration in velocity. In particular,
the normalized distribution perturbation has $L^2$ size
comparable to $a$, and
\begin{equation}\label{eq:backgroundCouplingMechanism}
\partial_{v_1}(F^{\rm b}_{\rm app}-\mu )
=a\,\partial_{v_1}\mu \cos(x_1-sv_1)+as\mu \sin(x_1-sv_1).
\end{equation}
Hence the term $E[H]\partial_{v_1}F^{\rm b}_{\rm app}$ in the
{linearized equation contains adjacent-mode interactions whose coefficients}
can be computed explicitly.

For any planar solution, the same change of variables reads
\[
s=t+t_0,\qquad
F(s,x_1,v)=\widetilde F(s-t_0,x,v),\qquad
E[F](s,x_1)=\widetilde E_1(s-t_0,x).
\]
Equation~\eqref{eq:vplPhysicalTime} becomes
\begin{equation}\label{eq:vpl}
\begin{cases}
\displaystyle
\partial_sF+v_1\partial_{x_1}F
+E[F]\partial_{v_1}F=\nu Q(F,F),\\[1mm]
\displaystyle
E[F]=-\partial_{x_1}(-\partial_{x_1}^2)^{-1}
\left(\int_{\mathbb R^3}F\,dv-1\right),\\[1mm]
F(t_0,x_1,v)=F_{\rm in}(x_1,v).
\end{cases}
\end{equation}
{We use \eqref{eq:vpl} throughout the construction: $s=t_0$
corresponds to physical time $t=0$, and $s=T_1$ to $t=t_\nu=T_1-t_0$.
The parameter $t_0$ may depend on $\nu$ and enters the prescribed
initial phases.}

The normalized perturbation and its profile in shifted time are
\begin{equation}\label{eq:shiftedProfileCorrespondence}
\begin{aligned}
f^\nu(s,x_1,v)
&=\widetilde f^\nu(s-t_0,x_1,v)
=\frac{F^\nu(s,x_1,v)-\mu (v)}{\sqrt \mu (v)},\\
g^\nu(s,z,v)
&=f^\nu(s,z+(s-t_0)v_1,v)
=\widetilde g^\nu(s-t_0,z,v),\\
\sqrt \mu (v)\,g^\nu(s,z,v)
&=F^\nu(s,z+(s-t_0)v_1,v)-\mu (v).
\end{aligned}
\end{equation}
\begin{rem}
For comparison, consider the freely transported perturbation
\[
a\psi(v)\cos(x_1-(t+t_0)v_1),\qquad t\ge0,
\]
where $\psi\in\mathcal S(\mathbb R^3)$ is real-valued.
Here $\psi$ is an auxiliary profile used only in this comparison;
the background $F^{\rm b}_{\rm app}$ in our construction remains
the Gaussian approximation defined in \eqref{eq:bgapp}.
Velocity integration gives
\[
\rho_k\!\left[a\psi(v)\cos(x_1-(t+t_0)v_1)\right]
=\frac a2\mathbf1_{|k|=1}
\widehat\psi\bigl(k(t+t_0),0,0\bigr).
\]
Thus, if
\[
\operatorname{supp}\widehat\psi
\subset\{|\xi_1|\le\tfrac12\},
\qquad t_0\ge1,
\]
the density of this perturbation vanishes for every $t\ge0$.

This calculation does not provide an alternative nonnegative
background of the form used here. Indeed,
\[
\mu(v)+a\psi(v)\cos(x_1-(t+t_0)v_1)\ge0
\quad\text{for all }x_1,v
\quad\Longrightarrow\quad
a|\psi(v)|\le\mu(v).
\]
The resulting Gaussian decay makes $\widehat\psi$ analytic in
$\xi_1$, so the stated support condition would force $\psi=0$.

For the Gaussian background actually used in \eqref{eq:bgapp},
the nonzero density modes are small rather than identically zero:
\[
\rho_{\pm1}[\widetilde F^{\rm b}_{\rm app}-\mu](t)
=\frac a2e^{-(t+t_0)^2/2}.
\]
More importantly for the background comparison, the leading
electric source satisfies \eqref{small-E}.
Our argument therefore uses the small accumulated electric source
provided by the initially mixed phase, rather than exact vanishing
of the field. The collision defect and the difference between the
exact and approximate backgrounds are estimated separately below.
\end{rem}

\begingroup
\subsubsection{Parameters for the construction}
The prepared initial phases are now fixed. We specify their amplitudes
and frequencies before defining the complete initial distribution.
The common parameter relations are as follows.
\par\endgroup For the main theorem,
the fixed mode parameters and loss constants are chosen first, then $K$, and
then a sufficiently large Sobolev order ${\rm s}_0$; the complete
order of choices is given at the end of Section~\ref{sec:roadmapNonlinear}.
The number $a$ is the background amplitude, $\zeta$ is the frequency
parameter in the initial phase, and $m_0,m_1$ are the first and last
spatial modes of the descending chain. A large parameter $N$
determines these modes.
Fix an initial Sobolev order ${\rm s}_0\ge11$ and a velocity weight
$\ell\ge4$, and choose the mode ratios by
\begin{equation}\label{eq:parameters}
\Lambda=\frac{\sqrt{2\pi}}2,\quad
0<{\bf b}_1<{\bf b}_0=\Lambda^{1/3},\quad
m_i=\lfloor{\bf b}_iN\rfloor,\quad 0<a\le1,\quad \zeta=\frac{N^3}{a}.
\end{equation}
For fixed $A_0\ge1$ and $\gamma_{\rm a}>0$, the construction uses
\begin{equation}\label{eq:mainparameters}
e^{-\gamma_{\rm a}N}\le a\le N^{-{\rm s}_0-2},\qquad
t_0=A_0N,\qquad T_1=\frac{\zeta}{m_1}.
\end{equation}
Here $\gamma_{\rm a}$ bounds how rapidly the amplitude may decrease
with $N$; $A_0$ is independent of $a,N$, with no further largeness
requirement. For the collisional construction choose fixed $K>0$,
$p>{\rm s}_0+7/2$ and small $c_0>0$, and set
\begin{equation}\label{eq:comparisonLogParameters}
N=K\log(1/\nu),\qquad
a=c_0\nu^{1/3}[\log(1/\nu)]^p(1+t_0)^{-{\rm s}_0}.
\end{equation}
Choose $\gamma_{\rm a}>1/(3K)$. For sufficiently small $\nu$, these
choices satisfy \eqref{eq:mainparameters}, $0<a\le1/2$, $m_1\ge1$,
$m_0\ge2$ and $m_0t_0\le\zeta/2$. 
These relations specify the proposed construction. The reasons for the
frequency choices are given in Section~\ref{sec:echoMechanismOverview}.
 Lemma~\ref{lem:logarithmicScales}
computes the scale restrictions and the condition on $p$.
\par
Unless a statement
explicitly concerns only the collisionless equation or allows a weaker
condition, the collisional estimates use these common choices.
We write $I=[t_0,T_1]$ and $t_\nu=T_1-t_0$.

The initial Sobolev size and the high-frequency amplitude are
\begin{equation}\label{eq:nonlinearTwoScales}
\epsilon_\nu=\nu^{1/3}[\log(1/\nu)]^p,\qquad
\delta_\nu=c_1\epsilon_\nu\zeta^{-{\rm s}_0},
\end{equation}
with fixed small $c_1>0$. Thus
$a(1+t_0)^{{\rm s}_0}=c_0\epsilon_\nu$ and
$\delta_\nu\zeta^{{\rm s}_0}=c_1\epsilon_\nu$:
the background and the secondary perturbation have the same prescribed
Sobolev scale, despite their different amplitudes and frequencies.
Lemma~\ref{lem:initialGaussianBounds} proves the corresponding initial
norm bounds.
\par

The prepared oscillations cost $(1+t_0)^{{\rm s}_0}$ in the initial
Sobolev norm; this is the reason for the factor in
\eqref{eq:comparisonLogParameters}. At fixed derivative and moment
orders, their Gevrey cost satisfies
$a(1+t_0)^C e^{Ct_0^\chi}=\nu^{1/3-o(1)}\to0$.
Also, $t_0=A_0N$ makes the Gaussian tail smaller than every fixed
power of $\nu$.

\begingroup
\subsubsection{The full initial datum and its properties}
We now prescribe the complete distribution at physical time zero.
The background and high-frequency oscillations have the amplitudes and
frequencies fixed above. We adjust only the temperature of the first
Maxwellian to match the total energy of $\mu$. The following lemma
gives the datum, its decomposition, and the estimates used later;
the proof is in Appendix~\ref{app:initialNormalization}.

\begin{lem}\label{lem:data}\label{lem:initialGaussianBounds}\label{lem:bginitialsize}
Let $0<a\le1/2$, $t_0\ge1$, $m_0\in\mathbb N$ with $m_0\ge2$, and $\zeta\in\mathbb R$.
For the amplitude $\delta_\nu$ in \eqref{eq:nonlinearTwoScales},
provided $|\delta_\nu|$ is sufficiently small, the complete initial datum is
\begin{equation}\label{eq:positivefamily}
F_{\rm in}
=\mu_{\theta_{\delta_\nu}}(1+a\cos(x_1-t_0v_1))
+\frac{\delta_\nu}{\sqrt\pi}\mu_{1/2}
\cos(m_0x_1+(\zeta-m_0t_0)v_1).
\end{equation}
For all sufficiently small $|\delta|$, the temperature $\theta_\delta$
is the unique root near one of
\begin{equation}\label{eq:thetadelta}
3(\theta_\delta-1)+\frac{a^2}{2}e^{-\theta_\delta t_0^2}
+\frac{\delta^2}{2\pi m_0^2}e^{-(\zeta-m_0t_0)^2/2}=0.
\end{equation}
The data have the following properties.

\emph{(i) Decomposition, positivity and conserved quantities.}
Set $\theta=\theta_0$ and define
\begin{equation}\label{eq:maininitial}
\begin{aligned}
F^{\rm b}_{\rm in}(x_1,v)&=\mu_\theta(v)(1+a\cos(x_1-t_0v_1)),\\
H_{\rm in}(x_1,v)&=\frac{\mu_{1/2}(v)}{\sqrt\pi}
\cos(m_0x_1+(\zeta-m_0t_0)v_1).
\end{aligned}
\end{equation}
For the auxiliary amplitude family
$F^\delta_{\rm in}=\mu_{\theta_\delta}(1+a\cos(x_1-t_0v_1))+\delta H_{\rm in}$,
the two distinct field modes give
\begin{equation}\label{eq:initialFullEnergy}
\mathcal E(F^\delta_{\rm in})
=\frac32\theta_\delta+\frac{a^2}{4}e^{-\theta_\delta t_0^2}
+\frac{\delta^2}{4\pi m_0^2}e^{-(\zeta-m_0t_0)^2/2}
=\mathcal E(\mu).
\end{equation}
Thus \eqref{eq:thetadelta} imposes equality of total energies while
leaving the amplitudes and phases of both oscillations unchanged.
Writing
\begin{equation}\label{eq:initialTemperatureCorrection}
\begin{aligned}
W^\delta_{\rm in}&=(\mu_{\theta_\delta}-\mu_\theta)(1+a\cos(x_1-t_0v_1)),
\qquad W_{\rm in}=W^{\delta_\nu}_{\rm in},\\
F_{\rm in}&=F^{\rm b}_{\rm in}+\delta_\nu H_{\rm in}+W_{\rm in},
\end{aligned}
\end{equation}
we have $F_{\rm in}=F^{\delta_\nu}_{\rm in}$ and
$\partial_\delta F^\delta_{\rm in}|_{\delta=0}=H_{\rm in}$.
The background temperature is the unique root in $[1-a^2/6,1]$ of
\begin{equation}\label{eq:bgtemp}
3(\theta-1)+\frac{a^2}{2}e^{-\theta t_0^2}=0,
\end{equation}
and
\begin{equation}\label{eq:thetabound}
0\le1-\theta\le Ca^2e^{-t_0^2},\qquad
|\theta_\delta-\theta|\le C\delta^2m_0^{-2}e^{-(\zeta-m_0t_0)^2/2}.
\end{equation}
Both $F^{\rm b}_{\rm in}$ and $F^\delta_{\rm in}$ have mass one,
momentum zero and total energy $\mathcal E(\mu)$; moreover,
$F^{\rm b}_{\rm in}>0$ and $F^\delta_{\rm in}\ge\mu_{\theta_\delta}/4$.

\emph{(ii) The initial remainder.}
For fixed $n\ge0$, $\ell\ge0$ and $0<q<1/4$,
\begin{equation}\label{eq:initialTemperatureMixedBound}
\|W^\delta_{\rm in}/\sqrt\mu\|_{H^n_{q,\ell}}
\le C_{n,q,\ell}\delta^2m_0^{-2}e^{-(\zeta-m_0t_0)^2/2}
[1+a(1+t_0)^n],
\end{equation}
whereas spatial derivatives alone satisfy
\begin{equation}\label{eq:initialTemperatureSpatialBound}
\sum_{j=0}^n\|\langle v\rangle^\ell\partial_{x_1}^j(W^\delta_{\rm in}/\sqrt\mu)\|_2
\le C_{n,\ell}\delta^2m_0^{-2}e^{-(\zeta-m_0t_0)^2/2}.
\end{equation}

\emph{(iii) Sobolev size.}
Under \eqref{eq:parameters}--\eqref{eq:mainparameters}, for fixed
$0<q<1/4$, $\ell\ge4$ and all sufficiently large $N$,
\begin{equation}\label{eq:maininitialnorm}
\|(F^{\rm b}_{\rm in}-\mu)/\sqrt\mu\|_{H^{{\rm s}_0}_{q,\ell}}
\sim a(1+t_0)^{{\rm s}_0},\qquad
\|H_{\rm in}/\sqrt\mu\|_{H^{{\rm s}_0}_{q,\ell}}\le C\zeta^{{\rm s}_0}.
\end{equation}
The first norm tends to zero as $N\to\infty$.

\emph{(iv) Gevrey size of the background datum.}
Assume also \eqref{eq:comparisonLogParameters} and fix $1/2<\chi<2/3$.
For every fixed ${\rm s},L,\lambda$, the
normalized datum $f^{\rm b}_{\rm in}=(F^{\rm b}_{\rm in}-\mu)/\sqrt\mu$
satisfies
\begin{equation}\label{eq:bginitialsize}
\|f^{\rm b}_{\rm in}\|_{\mathscr G^{\lambda,\chi}_{{\rm s},L}}
\le aL_N,\qquad L_N=C(1+N)^{r_L}e^{CN^\chi},
\end{equation}
with constants independent of $\nu$. In particular,
\begin{equation}\label{eq:epsNsmall}
\log L_N=o(N),\qquad aL_N=\nu^{1/3-o(1)}.
\end{equation}
After a fixed reduction of radius, the same bound holds in the
time-zero vector-field Gevrey norms of \cite{ref6}, with any
prescribed finite polynomial velocity moments.
\end{lem}

\par\endgroup

\begingroup
\subsubsection{The exact background and its explicit approximation}
Chaturvedi, Luk and Nguyen \cite{ref4} prove nonlinear stability
and Landau damping for sufficiently small Sobolev perturbations of
\emph{the global Maxwellian {$\mu $}}, of size $O(\nu^{1/3})$.
Our background data are also perturbations of {$\mu $}, but we choose
\[
\left\|\frac{\widetilde F^{\rm b,\nu}(0)-\mu }{\sqrt \mu }\right\|_{\mathscr H^{{\rm s}_0}_\ell}
\sim\epsilon_\nu=\nu^{1/3}[\log(1/\nu)]^p.
\]
The Sobolev theorem does not supply the uniform control needed for these data.
\par

We choose the background data to be small in a fixed-radius Gevrey
class $1/\chi$, with $1/2<\chi<2/3$. The uniform-in-$\nu$ theorem
of Bedrossian, Zhao and Zi \cite{ref6} then supplies a global exact
background, uniform Gevrey control and Landau damping. These estimates
control the coefficients of the linearized equation. The accurate
comparison of this background with an explicit free-transport solution,
needed to compute the echo amplification, is a separate estimate
proved here. In contrast, the much higher velocity frequency of the
secondary perturbation makes that perturbation small in the prescribed
Sobolev norm but not in the same fixed-radius Gevrey norm. There is
therefore no conflict between control and damping of the background
and Sobolev growth of the secondary perturbation.

\par
\par\endgroup
The profile \eqref{eq:bgapp} has only spatial modes $0,\pm1$;
the exact VPL background need not retain this support.

{The low-order comparison supplies smallness, while the
higher-order Gevrey bounds supply regularity. }

Fix $\chi\in(1/2,2/3)$. The coefficient estimates use a factor $L_N$
with $\log L_N=o(N)$ and the small size $\varepsilon_N=aL_N$;
{the precise Gevrey definition and dependence on fixed derivative orders}
are given in Proposition~\ref{lem:fixedCoefficientOrders}.

{We distinguish the background error at the distribution,
normalized-perturbation and transport-profile levels:}
\begin{equation}\label{eq:backgroundRemainder}
{
\begin{aligned}
{F^{\rm b}_{\rm err}}(s,x_1,v)
&=F^{\rm b,\nu}(s,x_1,v)-F^{\rm b}_{\rm app}(s,x_1,v),\\
{f^{\rm b}_{\rm err}}(s,x_1,v)
&=\frac{{F^{\rm b}_{\rm err}}(s,x_1,v)}{\sqrt \mu (v)}
=f^{\rm b,\nu}(s,x_1,v)-f^{\rm b}_{\rm app}(s,x_1,v),\\
{g^{\rm b}_{\rm err}}(t,z,v)
&={f^{\rm b}_{\rm err}}(t_0+t,z+tv_1,v).
\end{aligned}}
\end{equation}
{Here {$f^{\rm b,\nu}=(F^{\rm b,\nu}-\mu )/\sqrt \mu $} and
{$f^{\rm b}_{\rm app}=(F^{\rm b}_{\rm app}-\mu )/\sqrt \mu $}.
The error quantities depend on $\nu$, which is suppressed in their
superscripts. Density, current and electric field are evaluated on
{${{F^{\rm b}_{\rm err}}=\sqrt \mu  {f^{\rm b}_{\rm err}}}$}.}
{The normalized error {${f^{\rm b}_{\rm err}}$} includes the initial temperature discrepancy and satisfies} \eqref{eq:direct-r-eq}, whose
known sources are the electric and collision terms omitted by the
explicit approximation. The Gaussian terms below absorb the
initial discrepancy and the electric defect; the terms proportional
to $\nu$ measure the collision error.
\par

\begin{prop}\label{thm:bgcomparison}\label{prop:directBackground}
{Assume \eqref{eq:comparisonLogParameters}, with fixed ${\rm s}_0\ge11$,
$p>{\rm s}_0+7/2$, and $A_0\ge1$. For all sufficiently small $\nu>0$,
the VPL initial-value problem with datum $F^{\rm b}_{\rm in}$ in
\eqref{eq:maininitial} admits a global smooth nonnegative solution
$F^{\rm b,\nu}$ with the mass, momentum and total energy of {$\mu $}.
For every fixed velocity moment $\ell\ge0$ and sufficiently small $\nu$ depending on $\ell$, the remainder
{${{f^{\rm b}_{\rm err}}=\frac{F^{\rm b,\nu}-F^{\rm b}_{\rm app}}{\sqrt \mu }}$} satisfies}
\begin{align}
\|{f^{\rm b}_{\rm err}}(s)\|_2+\|E[{F^{\rm b}_{\rm err}}](s)\|_2
&\le Ca\{\nu(s-t_0)\langle s\rangle^2+e^{-t_0^2/4}\},
\label{eq:polynomialBackgroundLow}\\
\sup_{s\in I}\|\langle v\rangle^\ell {f^{\rm b}_{\rm err}}(s)\|_2
&\le C_\ell a\{\nu {T_1}^3+e^{-t_0^2/8}\}.
\label{eq:polynomialBackgroundMoments}
\end{align}
The first estimate holds for every $s\in I$. Neither estimate contains
$L_N$. The weighted $L^2$ norm in the second estimate is invariant
under spatial shears in the active variable.

\end{prop}

The proof is in Section~\ref{sec:backgroundresponse}. The first estimate
{keeps the dependence on the elapsed time;} the second supplies the
velocity moments needed to replace the electric coefficient
$E[H]\partial_{v_1}F^{\rm b}_{\rm app}$ by
$E[H]\partial_{v_1}F^{\rm b,\nu}$. 

\subsection{Plasma echoes in a toy model and the choice of initial frequencies}\label{sec:echoMechanismOverview}

\begingroup
Before treating the complete linearized equation, we isolate one
interaction that creates a plasma echo. The toy model below explains
both the delayed density growth and its amplification, and thereby
motivates the spatial and velocity frequencies in {$H_{\rm in}$}.
For this purpose only, we omit collisions, background-field transport,
the diagonal Maxwellian term in the density equation and the upward interactions.
We keep free transport and the interaction of the perturbation field
with the background {Fourier mode}
\endgroup
\[
(F^{\rm b}_{\rm app})_{-1}(\tau,v)=\frac a2\mu (v)e^{i\tau v_1}.
\]
The field of mode $m\ge2$ has coefficient $-i\rho_m/m$.
Multiplying by this background {Fourier mode} produces spatial mode $m-1$.
The density created in that mode then generates a field which can
produce mode $m-2$ through the same background {Fourier mode}.
Thus the background mode stays fixed while the perturbation moves
through the sequence $m_0,m_0-1,\ldots,m_1$.

A distribution-level version of this descending toy model is
\begin{equation}\label{eq:toyDistributionSystem}
\begin{aligned}
(\partial_s+im_0v_1)H_{m_0}^{\rm toy}&=0,
& H_{m_0}^{\rm toy}(t_0)&=H_{{\rm in},m_0},\\
(\partial_s+i(m-1)v_1)H_{m-1}^{\rm toy}
&=\frac{i\rho_m^{\rm toy}(s)}{m}
\partial_{v_1}\left(\frac a2\mu e^{isv_1}\right),
& H_{m-1}^{\rm toy}(t_0)&=0,
\end{aligned}
\end{equation}
where $m_1<m\le m_0$ and
$\rho_k^{\rm toy}(s)=\int_{\mathbb R^3}H_k^{\rm toy}(s,v)\,dv$.
 It is an auxiliary approximation, not the full linearized
VPL equation; related echo-chain models are studied in \cite{ref9}.
Duhamel's formula for the second equation, followed by velocity
integration, gives the kernel computed next.
\par

After transport from $\tau$ to $s$, integration by parts gives
\begin{equation}\label{eq:toyEchoKernel}
\begin{aligned}
\frac{i}{m}\int_{\mathbb R^3}e^{-i(m-1)(s-\tau)v_1}
\partial_{v_1}\left(\frac a2\mu e^{i\tau v_1}\right)\,dv
&=-\frac{a(m-1)}{2m}(s-\tau)
\int_{\mathbb R^3}\mu  e^{-i((m-1)s-m\tau)v_1}\,dv\\
&=-\frac{a(m-1)}{2m}(s-\tau)
e^{-((m-1)s-m\tau)^2/2}.
\end{aligned}
\end{equation}
{Taking the densities in \eqref{eq:toyDistributionSystem} therefore gives the triangular system}
\begin{equation}\label{eq:toyEchoRecursion}
\begin{aligned}
\rho_{m_0}^{\rm toy}(s)
&=\int_{\mathbb R^3}e^{-im_0(s-t_0)v_1}H_{{\rm in},m_0}(v)\,dv,\\
\rho_{m-1}^{\rm toy}(s)
&=-\frac{a(m-1)}{2m}\int_{t_0}^s(s-\tau)
e^{-((m-1)s-m\tau)^2/2}\rho_m^{\rm toy}(\tau)\,d\tau,
\qquad m_1<m\le m_0.
\end{aligned}
\end{equation}
 The function $\rho_{m_0}^{\rm toy}$ is determined by the initial perturbation; we now explain which frequencies make
this recurrence useful.

\begingroup
If the density in mode $m$ is concentrated near $\tau=\zeta/m$,
the Gaussian in \eqref{eq:toyEchoRecursion} selects the later time
$s=\zeta/(m-1)$ in mode $m-1$:
\[
\tau=\frac{\zeta}{m},\quad (m-1)s=m\tau
\quad\Longrightarrow\quad s=\frac{\zeta}{m-1},\qquad
\frac{\zeta}{m-1}-\frac{\zeta}{m}=\frac{\zeta}{m(m-1)}>0.
\]
This newly generated density at a later time is the plasma echo. For fixed $s$, the Gaussian kernel restricts the integration variable $\tau$ to an interval of width
$O(1/m)$. Freezing the prefactor at the resonant times and integrating
the Gaussian over the real line gives the single-step scale
\begin{equation}\label{eq:toyEchoGain}
\frac{a(m-1)}{2m}
\left(\frac{\zeta}{m-1}-\frac{\zeta}{m}\right)
\int_{\mathbb R}e^{-(m\tau-\zeta)^2/2}\,d\tau
=\frac{\sqrt{2\pi}\,a\zeta}{2m^3}.
\end{equation}
This calculation concerns the mass of the localized transfer kernel;
successive Gaussian convolutions also increase the width of the density profile. A terminal
lower bound therefore requires control of that spreading, not only
multiplication of the gains.
\par\endgroup
Formula~\eqref{eq:toyEchoGain} explains the spatial frequency range.
The gain is at least one for $m^3\le\Lambda a\zeta$, where
$\Lambda=\sqrt{2\pi}/2$. Setting $a\zeta=N^3$ and starting at
$m_0=\lfloor\Lambda^{1/3}N\rfloor$ places the first mode at the
upper end of this amplifying range. Descending to
$m_1=\lfloor{\bf b}_1N\rfloor$, with
$0<{\bf b}_1<\Lambda^{1/3}$, gives order $N$ amplifying interactions:
\[
\log\prod_{m=m_1+1}^{m_0}\frac{\Lambda N^3}{m^3}
=\sum_{m=m_1+1}^{m_0}\log\frac{\Lambda N^3}{m^3}\sim N.
\]
\begingroup
These are the choices in \eqref{eq:parameters}. They fix the first
density maximum near $\zeta/m_0$ and the observation time $T_1=\zeta/m_1$.
Thus order $N$ amplifying interactions are available without following
the chain to a fixed low mode. 
\par
\par\endgroup

\begingroup
The initial velocity frequency must place the top density maximum
at that first target time. For a Gaussian with modulation $e^{i\xi v_1}$,
where $\xi\in\mathbb R$, free transport gives
\endgroup
\[
\int_{\mathbb R^3}e^{-im_0(s-t_0)v_1}\mu_{1/2}(v)e^{i\xi v_1}\,dv
=e^{-(m_0(s-t_0)-\xi)^2/4}.
\]
This Gaussian has its maximum at $s=t_0+\xi/m_0$. To put the maximum at
$\zeta/m_0$, we choose
\[
t_0+\frac{\xi}{m_0}=\frac{\zeta}{m_0}
\quad\Longleftrightarrow\quad \xi=\zeta-m_0t_0.
\]
This determines the phase of the positive Fourier mode of
\eqref{eq:maininitial}:
\begin{equation}\label{eq:highdatum}
H_{{\rm in},m_0}(v)=\frac{\mu_{1/2}(v)}{\sqrt{4\pi}}
e^{i(\zeta-m_0t_0)v_1}.
\end{equation}
Set {$H_{{\rm in},-m_0}=\overline{H_{{\rm in},m_0}}$} and set all other
spatial modes to zero. In velocity Fourier variables this choice reads
\begin{equation}\label{eq:initialEchoFrequencyBand}
\widehat H_{{\rm in},m_0}(\eta)
=\frac1{\sqrt{4\pi}}
e^{-|\eta-(\zeta-m_0t_0)e_1|^2/4}.
\end{equation}
\begingroup
Thus the spatial frequencies are exactly $\pm m_0$, and the velocity
frequencies are concentrated in Gaussian Fourier profiles of width $O(1)$
around $\pm(\zeta-m_0t_0)e_1$. This is Gaussian localization, not
compact Fourier support. The condition $m_0t_0\le\zeta/2$ puts
the first free-transport density maximum after the initial time and gives a velocity frequency
comparable to $\zeta$. An unmodulated Gaussian would instead have
the largest free-transport density at $s=t_0$.
\par\endgroup

\begingroup
Using only the two conjugate spatial modes starts one descending chain
in each sign block, with $H_{{\rm in},k}=0$ for $0<|k|<m_0$. The fixed velocity-frequency width gives a first density maximum with width $O(1/m_0)$ in time, matching the width of one adjacent-mode interaction.
The Gaussian factor in $\widehat H_{{\rm in},\pm m_0}$ makes the subsequent convolutions explicit.
\par\endgroup
{Evaluating the velocity Fourier transform in
\eqref{eq:initialEchoFrequencyBand} at $\eta=m_0(s-t_0)e_1$
directly gives the free-transport density $\mathfrak g^0$:}
\begin{equation}\label{eq:freeGaussianDensity}
{
\begin{aligned}
\mathfrak g^0_{m_0}(s)
&=\int_{\mathbb R^3}e^{-im_0(s-t_0)v_1}H_{{\rm in},m_0}(v)\,dv
=\widehat H_{{\rm in},m_0}\big(m_0(s-t_0)e_1\big)\\
&=\frac{1}{\sqrt{4\pi}}
\exp\!\left(-\frac{[m_0(s-t_0)-(\zeta-m_0t_0)]^2}{4}\right)
=\frac{e^{-(m_0s-\zeta)^2/4}}{\sqrt{4\pi}}.
\end{aligned}}
\end{equation}
\begingroup
In particular, $\mathfrak g^0_{m_0}(s)$ is exponentially small at $s=t_0$ but
has size of order one near $s=\zeta/m_0$ before multiplication by
$\delta_\nu$. The frequency choice therefore prepares a density maximum at a positive time, rather than a large initial electric field.
The cooler Maxwellian {$\mu_{1/2}$} also permits pointwise domination
by the background in Lemma~\ref{lem:data}. The high-frequency
Sobolev cost of {$H_{\rm in}/\sqrt \mu $} is of order $\zeta^{{\rm s}_0}$, which explains the
factor $\zeta^{-{\rm s}_0}$ in the amplitude $\delta_\nu$ in
\eqref{eq:nonlinearTwoScales}.
\par\endgroup

\begingroup
The actual linearized VPL equation has
the exact, nonexplicit background, both Landau collision terms,
the Maxwellian term in the density equation, upward interactions and
background-field transport. We next extract a growing density from
that equation and identify the terms whose effects must be controlled.
\par
\par\endgroup
\subsection{Extracting the growing density from the linearized VPL equation}\label{sec:roadmapOperators}\label{sec:roadmapEcho}

We now return from the toy model to the complete equation linearized
about the exact background:
\begin{equation}\label{eq:fullRangePhysicalEquation}
\begin{split}
\partial_sH^\nu+v_1\partial_{x_1}H^\nu
+E[F^{\rm b,\nu}]\partial_{v_1}H^\nu
+E[H^\nu]\partial_{v_1}F^{\rm b,\nu}
={}&\nu\big(Q(F^{\rm b,\nu},H^\nu)
+Q(H^\nu,F^{\rm b,\nu})\big),\\
H^\nu(t_0)={}&H_{\rm in}.
\end{split}
\end{equation}
\begingroup
Unlike \eqref{eq:toyEchoRecursion}, this equation contains collisions
and both electric terms, and the background is not explicit.  In this subsection we derive the density equation and
define the descending density with the Maxwellian Volterra term included.
Section~\ref{sec:roadmapLinearErrors} controls the collision,
background-coefficient and upward-interaction errors;
Section~\ref{sec:roadmapFieldCompletion} restores background-field
transport.
\par\endgroup
We begin by omitting only background-field transport and defining
the operator
\par

\begin{equation}\label{eq:fieldExpansionPhysicalOperator}
\mathcal L_{\rm C}H=\partial_sH+v_1\partial_{x_1}H
+E[H]\partial_{v_1}F^{\rm b,\nu}
-\nu\big(Q(F^{\rm b,\nu},H)+Q(H,F^{\rm b,\nu})\big).
\end{equation}

For the initial perturbation {$H_{\rm in}$} in \eqref{eq:maininitial}, we first seek
\begin{equation}\label{eq:exactCollisionModel}
\mathcal L_{\rm C}H^{\rm C,\nu}=0,\qquad H^{\rm C,\nu}(t_0)=H_{\rm in}.
\end{equation}

For the density reduction, define the normalized background collision
operator by
\begin{equation}\label{eq:backgroundNormalizedLinearization}
\begin{aligned}
\mathcal C_{\rm b}(s)f
&=\mu^{-1/2}\big[Q(F^{\rm b,\nu}(s),\sqrt \mu  f)
+Q(\sqrt \mu  f,F^{\rm b,\nu}(s))\big]\\
&=-\mathcal L f
+\Gamma\left(f^{\rm b,\nu}(s),f\right)
+\Gamma\left(f,f^{\rm b,\nu}(s)\right).
\end{aligned}
\end{equation}
{Here $f^{\rm b,\nu}=(F^{\rm b,\nu}-\mu)/\sqrt\mu$ was defined in \eqref{eq:backgroundRemainder}.}
The second equality follows from
\eqref{eq:normalizedLandauExplicit}--\eqref{eq:normalizedGammaExplicit}.
For {$f=H^{\rm C,\nu}/\sqrt \mu $}, equation \eqref{eq:exactCollisionModel}
therefore reads
\begin{equation}\label{eq:normalizedExactCollisionModel}
\partial_sf+v_1\partial_{x_1}f-\nu\mathcal C_{\rm b}(s)f
=-\mu^{-1/2}E[\rho^{\rm C,\nu}]\partial_{v_1}F^{\rm b,\nu},
\qquad \rho^{\rm C,\nu}=\int_{\mathbb R^3}\sqrt \mu  f\,dv.
\end{equation}

To derive a Volterra equation for this density, use the linear
initial-value problem without the self-consistent electric term,
as in \cite[Section~7.1]{ref4} and \cite[Section~6.2]{ref6}:
\begin{equation}\label{eq:exactCollisionEvolution}
\partial_sf+v_1\partial_{x_1}f=\nu\mathcal C_{\rm b}(s)f,
\qquad f(\tau)=f_\tau.
\end{equation}

Denote its solution operator by $\mathbb U_{\rm b}^\nu(s,\tau)$, so
that $f(s)=\mathbb U_{\rm b}^\nu(s,\tau)f_\tau$. This equation uses only
the known background. Because the coefficients depend on time, both
time arguments are needed; spatially dependent coefficients also couple
Fourier modes. Section~\ref{sec:centeredcollisions} constructs this
solution operator and proves the estimates used below.

For the normalized equation \eqref{eq:normalizedExactCollisionModel},
Duhamel's formula gives
\begin{equation}\label{eq:exactCollisionDuhamel}
\frac{H^{\rm C,\nu}(s)}{\sqrt \mu }
=\mathbb U_{\rm b}^\nu(s,t_0)(H_{\rm in}/\sqrt \mu )
-\int_{t_0}^s\mathbb U_{\rm b}^\nu(s,\tau)
\big[\mu^{-1/2}E[\rho^{\rm C,\nu}](\tau)\partial_{v_1}F^{\rm b,\nu}(\tau)\big]\,d\tau,
\end{equation}
where $E[\rho^{\rm C,\nu}]$ is
given by Poisson's equation. The initial datum contributes the known
density
\begin{equation}\label{eq:exactDensityInitialTrace}
\mathfrak g^{\rm C,\nu}(s)
=\int_{\mathbb R^3}\sqrt \mu \,
\mathbb U_{\rm b}^\nu(s,t_0)(H_{\rm in}/\sqrt \mu )\,dv.
\end{equation}
The function $\mathfrak g^{\rm C,\nu}$ is known once the datum and the
solution operator are fixed; it is distinct from the unknown
self-consistent density $\rho^{\rm C,\nu}$.
Multiplying \eqref{eq:exactCollisionDuhamel} by {$\sqrt \mu $} and
integrating in $v$ yields the exact density equation
\begin{equation}\label{eq:exactDensityVolterra}
\rho^{\rm C,\nu}(s)=\mathfrak g^{\rm C,\nu}(s)-\int_{t_0}^s\!\int_{\mathbb R^3}\sqrt \mu \,
\mathbb U_{\rm b}^\nu(s,\tau)
\big[\mu^{-1/2}E[\rho^{\rm C,\nu}](\tau)\partial_{v_1}F^{\rm b,\nu}(\tau)\big]\,dv\,d\tau.
\end{equation}

The background and the solution operator in this equation are not
explicit, which makes the accumulated echo amplification difficult
to compute.

\emph{We approximate this equation by replacing $F^{\rm b,\nu}$ with
$F^{\rm b}_{\rm app}$ and the collision solution operator with free
transport.}
\par At the distribution level, the approximate equation is
\begin{equation}\label{eq:overviewApproximateLinearEquation}
\partial_sH^{\rm app}+v_1\partial_{x_1}H^{\rm app}
+E[H^{\rm app}]\partial_{v_1}F^{\rm b}_{\rm app}=0,\qquad H^{\rm app}(t_0)=H_{\rm in}.
\end{equation}
Write $\rho^{\rm app}=\int H^{\rm app}\,dv$ for the density of this
approximate solution. For $k\ne0$, Poisson's equation gives
$E_k[H^{\rm app}]=-i\rho_k^{\rm app}/k$, and
Duhamel's formula for this equation reads
\[
H_k^{\rm app}(s,v)=e^{-ik(s-t_0)v_1}H_{{\rm in},k}(v)
+\sum_{m\ne0}\frac{i}{m}\int_{t_0}^s
e^{-ik(s-\tau)v_1}\rho_m^{\rm app}(\tau)
\partial_{v_1}(F^{\rm b}_{\rm app})_{k-m}(\tau,v)\,d\tau.
\]

The explicit background has only the spatial modes
\[
(F^{\rm b}_{\rm app})_0=\mu ,\qquad
(F^{\rm b}_{\rm app})_{\pm1}(\tau,v)=\frac a2\mu (v)e^{\mp i\tau v_1}.
\]
Integration by parts in $v_1$ and Gaussian integration give
\[
\begin{aligned}
\frac{i}{k}\int_{\mathbb R^3}e^{-ik(s-\tau)v_1}
\partial_{v_1}\mu \,dv
&=-(s-\tau)e^{-k^2(s-\tau)^2/2},\\
\frac{i}{m}\int_{\mathbb R^3}e^{-ik(s-\tau)v_1}
\partial_{v_1}\left(\frac a2\mu  e^{-i(k-m)\tau v_1}\right)\,dv
&=-\frac{ak}{2m}(s-\tau)e^{-(ks-m\tau)^2/2},
\qquad |k-m|=1.
\end{aligned}
\]
The density $\mathfrak g^0$ and the computed kernels are
\begin{equation}\label{eq:overviewApproximateDensityKernels}
\begin{aligned}
\mathfrak g_k^0(s)&=\widehat H_{{\rm in},k}\bigl(k(s-t_0)e_1\bigr),\qquad K_k(t)=-t e^{-k^2t^2/2}\mathbf1_{t\ge0}\quad(k\ge1),\\
P_{k,m}(s,\tau)&=-\frac{ak}{2m}(s-\tau)e^{-(ks-m\tau)^2/2}
\quad(k,m\ne0,\ |k-m|=1).
\end{aligned}
\end{equation}

\begingroup
Velocity integration therefore gives the approximate density equation,
for $k\ne0$. The braces define the density integral operators through their
action on the actual density $\rho^{\rm app}$:
\endgroup
\begin{equation}\label{eq:overviewApproximateDensityEquation}
\begin{aligned}
{\rho_k^{\rm app}}(s)={}&\mathfrak g_k^0(s)
+\underbrace{\int_{t_0}^sK_{|k|}(s-\tau){\rho_k^{\rm app}}(\tau)\,d\tau}
_{(\mathbb K{\rho^{\rm app}})_k(s)}
+\underbrace{\sum_{\substack{m\ne0,\ |k-m|=1\\|m|=|k|+1}}
\int_{t_0}^sP_{k,m}(s,\tau){\rho_m^{\rm app}}(\tau)\,d\tau}
_{(\mathbb P_\downarrow{\rho^{\rm app}})_k(s)}\\
&+\underbrace{\sum_{\substack{m\ne0,\ |k-m|=1\\|m|=|k|-1}}
\int_{t_0}^sP_{k,m}(s,\tau){\rho_m^{\rm app}}(\tau)\,d\tau}
_{(\mathbb P_\uparrow{\rho^{\rm app}})_k(s)}.
\end{aligned}
\end{equation}
The first integral comes from the Maxwellian mode of the background.
The other two come from the background modes $\pm1$ and map $\rho_m$ to a contribution to $\rho_k$ with $|k|<|m|$ or $|k|>|m|$, respectively. These are the three contributions of the
negative electric Duhamel term. For an arbitrary density $\rho$,
the corresponding physical-space operator identity is
\begin{equation}\label{eq:approximateElectricIntegralIdentity}
\begin{aligned}
&-\int_{t_0}^s\!\int_{\mathbb R^3}
E[{\rho}]\bigl(\tau,x_1-(s-\tau)v_1\bigr)
(\partial_{v_1}F^{\rm b}_{\rm app})
\bigl(\tau,x_1-(s-\tau)v_1,v\bigr)\,dv\,d\tau\\
&\qquad=\big[(\mathbb K+\mathbb P_\downarrow+\mathbb P_\uparrow){\rho}\big](s,x_1).
\end{aligned}
\end{equation}

\emph{We make a second approximation by omitting the upward term in
\eqref{eq:overviewApproximateDensityEquation}, keeping the Maxwellian
and descending terms, and retaining the same function $\mathfrak g^0$.}
\par
This defines the \emph{descending density} $\rho^\downarrow$.
Use the fixed modes $m_1<m_0$ from \eqref{eq:parameters}. Since the datum has only modes
$\pm m_0$, $\mathfrak g_k^0=0$ for $k>0$, $k\ne m_0$.
The defining equations on $m_1\le k\le m_0$ are
\begin{equation}\label{eq:reduced}
\begin{aligned}
\rho_{m_0}^\downarrow(s)
&=\mathfrak g_{m_0}^0(s)
+\int_{t_0}^sK_{m_0}(s-\tau)\rho_{m_0}^\downarrow(\tau)\,d\tau,\\
\rho_k^\downarrow(s)
&=\int_{t_0}^sK_k(s-\tau)\rho_k^\downarrow(\tau)\,d\tau
-\frac{ak}{2(k+1)}\int_{t_0}^s(s-\tau)
 e^{-(ks-(k+1)\tau)^2/2}\rho_{k+1}^\downarrow(\tau)\,d\tau,
\end{aligned}
\end{equation}
where $m_1\le k<m_0$ in the second line and $s\ge t_0$.
The top equation determines $\rho_{m_0}^\downarrow$; the second then
successively determines $\rho_{m_0-1}^\downarrow,\ldots,
\rho_{m_1}^\downarrow$. Each step is a scalar Volterra equation with
a known inhomogeneous term determined by $\rho_{k+1}^\downarrow$. We continue the same recursion to $k=m_1-1,\ldots,1$,
set $\rho_k^\downarrow=0$ for $k>m_0$, and set
$\rho_{-k}^\downarrow=\overline{\rho_k^\downarrow}$ and
$\rho_0^\downarrow=0$. Modes below $m_1$ do not enter the equations
on $m_1\le k\le m_0$, so this continuation preserves that chain.
The defining equations now give, on all nonzero modes,
\begin{equation}\label{eq:overviewDescendingDensityEquation}
\rho^\downarrow=\mathfrak g^0+(\mathbb K+\mathbb P_\downarrow)\rho^\downarrow.
\end{equation}

\begingroup
The free-transport solution operator is
$[\mathbb U^0(s,\tau)f]_k=e^{-ik(s-\tau)v_1}f_k$.
Each diagonal Volterra equation is inverted through the scalar equation
\begin{equation}\label{eq:scalar}
\rho_k(s)=J_k(s)+\int_{t_0}^sK_k(s-\tau)\rho_k(\tau)\,d\tau.
\end{equation}
\par\endgroup

\begingroup
Here $J_k$ is a prescribed inhomogeneous term. Lemma~\ref{prop:inverse}
proves uniform and weighted inverse bounds, including the signed
kernel at $k=1$.
\par
\par\endgroup
For the datum above, $\mathfrak g^0_{m_0}$ is the Gaussian
$\mathfrak g^0_{m_0}(s)={\frac{e^{-(m_0s-\zeta)^2/4}}{\sqrt{4\pi}}}$ in
\eqref{eq:freeGaussianDensity}. In particular,
\[
\rho_{m_0}^\downarrow(t_0)=\frac{e^{-(m_0t_0-\zeta)^2/4}}{\sqrt{4\pi}},
\qquad \rho_k^\downarrow(t_0)=0\quad(1\le k<m_0).
\]
The diagonal kernel is signed; we do not assume positivity of the
descending density.

\begingroup
The resonance times and the single-step scale were derived in
Section~\ref{sec:echoMechanismOverview}. The density $\rho^\downarrow$
defined here differs from $\rho^{\rm toy}$: at each step it also solves
the signed Maxwellian Volterra term in \eqref{eq:reduced}.
The following estimates show that this additional Volterra term does not
destroy the Gaussian chain growth. The single-step gains and their
products are normalized by
\par
\par\endgroup
\begin{equation}\label{eq:alphagain}
\alpha_m=\frac{a\sqrt{2\pi}\zeta}{2m^3}
=\frac{\Lambda N^3}{m^3},\qquad
B_m=\prod_{j=m+1}^{m_0}\alpha_j.
\end{equation}

An empty product equals one; this defines $B_m$ for every integer $m\ge1$, with $B_m=1$ for $m\ge m_0$.

The density recursion uses the gain
\begin{equation}\label{eq:echoGainParameters}
B_*=B_{m_1}=\prod_{j=m_1+1}^{m_0}\frac{\Lambda N^3}{j^3},\qquad c_{\rm e}=\int_{{\bf b}_1}^{{\bf b}_0}\log\frac{\Lambda}{y^3}\,dy>0,
\qquad \log B_*=c_{\rm e}N+O(1).
\end{equation}

To control cancellation, we compare $(-1)^{m_0-k}\rho_k^\downarrow$
with a positive Gaussian convolution. The diagonal correction is small:
\[
\int_0^\infty|K_k(t)|\,dt=k^{-2}=O(N^{-2}),\qquad m_1\le k\le m_0.
\]
In the rescaled variable $w=ms-\zeta$ for mode $m$,
Lemma~\ref{lem:step} shows that, after division by the single-step gain,
one transfer differs from convolution with $(2\pi)^{-1/2}e^{-w^2/2}$ by
$O(a/N+N^{-2})$. After
$m_0-m_1$ transfers, the comparison profile is
$\dfrac{e^{-w^2/[2(m_0-m_1+2)]}}{\sqrt{2\pi(m_0-m_1+2)}}$, whose value at zero is
$[2\pi(m_0-m_1+2)]^{-1/2}$. Thus the product $B_*$ gives exponential
amplification, while convolution produces the factor $N^{-1/2}$.
\par
{The following proposition makes this comparison quantitative
and justifies the continuation to the modes needed in the later inverse.}

\begin{prop}\label{thm:chain}
{Assume \eqref{eq:parameters}, with sufficiently large $N$, $m_1\ge1$,
$t_0\ge1$ and $m_0t_0\le\zeta/2$. Let $\rho^\downarrow$ solve the
reduced density system \eqref{eq:reduced}, and set $T_1=\zeta/m_1$.
Then}
\begin{equation}\label{eq:chaincomparison}
\left|\frac{(-1)^{m_0-m_1}\rho_{m_1}^\downarrow({T_1})}{B_{m_1}}
-\frac1{\sqrt{2\pi(m_0-m_1+2)}}\right|
\le C\left(a+\frac1N+e^{-c\zeta^2}\right).
\end{equation}
Consequently, if $a\sqrt N\to0$, then for all sufficiently large $N$,
\begin{equation}\label{eq:chainlower}
|\rho_{m_1}^\downarrow({T_1})|\ge\frac{c B_{m_1}}{\sqrt N},
\qquad
{\log B_{m_1}=c_{\rm e}N+O(1).}
\end{equation}
{For the real datum \eqref{eq:maininitial}, the continuation
specified after \eqref{eq:reduced} satisfies
\eqref{eq:overviewDescendingDensityEquation} and
\[
\rho^\downarrow\in C(I;H_x^r)\quad\text{for every fixed }r.
\]
This continuation preserves the estimates on $m_1\le k\le m_0$.}
\end{prop}

The mode $m_0$ fixed in \eqref{eq:parameters} is the largest integer
with $\alpha_m\ge1$. Consequently,
sources above $m_0$ do not accumulate a further exponential gain before
entering the descending chain. Write
\[
\widetilde B=\sup_{m\ge1}B_m,\qquad C_{\rm B}=3{\bf b}_0+1.
\]
Lemma~\ref{lem:extendedGain} proves $\widetilde B\le e^{C_{\rm B}N}$.
This upper bound and $B_*=e^{c_{\rm e}N+O(1)}$ have different roles:
the former controls arbitrary modes in an error estimate; the latter
is the growth factor in the lower bound for $\rho^\downarrow_{m_1}(T_1)$.
Proposition~\ref{thm:chain} is proved in Section~\ref{sec:chain}.

\begingroup
The time shift also makes the upward-interaction term perturbatively small.
For an upward interaction on positive modes, $k=m+1$ and
$t_0\le\tau\le s$, so
\[
ks-m\tau=s+m(s-\tau)\ge s\ge t_0.
\]
Unlike the descending interaction, this Gaussian cannot reach its
zero of the kernel phase in the integration region. Since
$k(s-\tau)\le ks-m\tau$, the change of variable
$y=ks-m\tau$ gives
\begin{equation}\label{eq:upwardTimeShiftSmallness}
\int_{t_0}^s|P_{m+1,m}(s,\tau)|\,d\tau
\le\frac{a}{2m^2}\int_s^\infty y e^{-y^2/2}\,dy
=\frac{a}{2m^2}e^{-s^2/2}
\le\frac{a}{2m^2}e^{-t_0^2/2}.
\end{equation}
The negative modes satisfy the same estimate after conjugation.
Thus the choice of initially phase-mixed data has two uses: it suppresses the electric
source in the background comparison, and it separates all upward
interactions from their resonances. With $t_0=A_0N$, the latter
effect supplies a Gaussian-small factor in $N$. At $t_0=0$, the displayed
bound supplies no additional Gaussian factor uniformly over the
initial time interval.
\par
\par\endgroup
\par

\begingroup
We have therefore identified both the growing density and the small
upward term omitted from its equation. The next subsection compares
this density with the exact collisional density.
\par
\par\endgroup
\subsection{Collision and coefficient errors}\label{sec:roadmapScales}\label{sec:roadmapLinearErrors}

The descending density now has a lower bound. We compare it with
$\rho^{\rm C,\nu}$, the density solving
\eqref{eq:exactDensityVolterra}. That equation includes both Landau
collision terms and the exact background in $E[\rho]\partial_{v_1}F^{\rm b,\nu}$;
background-field transport is left to
Section~\ref{sec:roadmapFieldCompletion}. The required comparison is
\par

\[
\frac{\sqrt N}{B_*}
\left|\rho_{m_1}^{\rm C,\nu}(T_1)-\rho_{m_1}^\downarrow(T_1)\right|\longrightarrow0.
\]
\begingroup
The defining equation \eqref{eq:overviewDescendingDensityEquation}
already identifies the error due to omitting the upward-interaction term:
\[
({\mathrm{Id}}-\mathbb K-\mathbb P_\downarrow-\mathbb P_\uparrow)\rho^\downarrow-\mathfrak g^0
=-\mathbb P_\uparrow\rho^\downarrow.
\]
It remains to account for the changes in collision dynamics and
background coefficients, as well as the evolution of the same initial
datum under these different equations.
\par\endgroup
For an arbitrary prescribed density ${\rho}$, subtract the electric
Duhamel term of the approximate equation from that of the exact
equation. Writing the free transport explicitly gives
\begin{equation}\label{eq:overviewDensityResponseDifference}
\begin{aligned}
(\Delta\mathbb T{\rho})(s,x_1)
={}&-\int_{t_0}^s\!\int_{\mathbb R^3}\sqrt \mu (v)\,
\Big[\mathbb U_{\rm b}^\nu(s,\tau)
\big(\mu^{-1/2}E[{\rho}](\tau)\partial_{v_1}F^{\rm b,\nu}(\tau)\big)\Big]
(x_1,v)\,dv\,d\tau\\
&+\int_{t_0}^s\!\int_{\mathbb R^3}
E[{\rho}]\bigl(\tau,x_1-(s-\tau)v_1\bigr)\,
(\partial_{v_1}F^{\rm b}_{\rm app})
\bigl(\tau,x_1-(s-\tau)v_1,v\bigr)\,dv\,d\tau.
\end{aligned}
\end{equation}
\begingroup
In the second integral, $\partial_{v_1}F^{\rm b}_{\rm app}$ is
computed at fixed spatial position before evaluation at
$x_1-(s-\tau)v_1$. Free transport leaves {$\mu (v)$} unchanged, so the
normalizing factors {$\sqrt \mu $} and {$\mu^{-1/2}$} cancel in this term.
By \eqref{eq:approximateElectricIntegralIdentity}, the second,
positive integral equals
$-[(\mathbb K+\mathbb P_\downarrow+\mathbb P_\uparrow){\rho}](s,x_1)$.
Thus it subtracts the approximate electric term in the density equation. This formula
displays both changes:
the exact collision dynamics replace free transport, and
$F^{\rm b,\nu}$ replaces $F^{\rm b}_{\rm app}$ in the electric term $E[\rho]\partial_{v_1}F^{\rm b}_{\rm app}$.
\par\endgroup
By the definition of $\Delta\mathbb T$, the exact density equation
\eqref{eq:exactDensityVolterra} takes the form
\begin{equation}\label{eq:overviewExactCoefficientDensityEquation}
\rho^{\rm C,\nu}=\mathfrak g^{\rm C,\nu}
+(\mathbb K+\mathbb P_\downarrow+\mathbb P_\uparrow+\Delta\mathbb T)\rho^{\rm C,\nu}.
\end{equation}
Subtracting the defining equation
\eqref{eq:overviewDescendingDensityEquation} from this equation gives
\[
\rho^{\rm C,\nu}-\rho^\downarrow
=\mathfrak g^{\rm C,\nu}-\mathfrak g^0
+(\mathbb K+\mathbb P_\downarrow)(\rho^{\rm C,\nu}-\rho^\downarrow)
+(\mathbb P_\uparrow+\Delta\mathbb T)\rho^{\rm C,\nu}.
\]
In the last term, use
\[
(\mathbb P_\uparrow+\Delta\mathbb T)\rho^{\rm C,\nu}
=(\mathbb P_\uparrow+\Delta\mathbb T)(\rho^{\rm C,\nu}-\rho^\downarrow)
+(\mathbb P_\uparrow+\Delta\mathbb T)\rho^\downarrow.
\]
Moving every term containing the density difference to the left yields
\begin{equation}\label{eq:overviewDensityDifferenceEquation}
({\mathrm{Id}}-\mathbb K-\mathbb P_\downarrow-\mathbb P_\uparrow-\Delta\mathbb T)
(\rho^{\rm C,\nu}-\rho^\downarrow)
=\mathfrak g^{\rm C,\nu}-\mathfrak g^0
+\Delta\mathbb T\rho^\downarrow+\mathbb P_\uparrow\rho^\downarrow.
\end{equation}
The left-hand operator is precisely the one in the exact density
equation. For general inhomogeneous problems $\mathfrak g$ also contains the density obtained from the prescribed source. The exact density equation then reads
\begin{equation}\label{eq:exactCoefficientDensityEquation}
{\rho}=\mathfrak g+(\mathbb K+\mathbb P_\downarrow+\mathbb P_\uparrow){\rho}+\Delta\mathbb T{\rho}.
\end{equation}

\begingroup
{For each fixed spatial order $r_{\rm x}\ge0$, we work on
the Banach space $\mathscr X_N^{r_{\rm x}}(I)$ of continuous,
mean-zero time-dependent densities. Its norm is defined precisely in
\eqref{eq:echoDensityNorm} in Section~\ref{sec:gainDensityInverse},
using the phase weight in \eqref{eq:echoPhaseWeights} and
the early--late time split in \eqref{eq:electricTimeSplit}.
It measures the small early-time density without a gain denominator
and divides the late-time density by the chain gains $B_{|k|}$.}
\endgroup
To solve \eqref{eq:exactCoefficientDensityEquation}, we first
invert the full approximate density equation and estimate the change:
\par
\[
\mathbb V_{\rm e}^0=({\mathrm{Id}}-\mathbb K-\mathbb P_\downarrow-\mathbb P_\uparrow)^{-1},
\qquad {
\|\mathbb V_{\rm e}^0\Delta\mathbb T\|_
{\mathscr X_N^{r_{\rm x}}(I)\to\mathscr X_N^{r_{\rm x}}(I)}=o(1).}
\]
\begingroup
Consequently the exact density solution operator satisfies
\[
\mathbb V^{\rm C,\nu}
=(\mathrm{Id}-\mathbb K-\mathbb P_\downarrow-\mathbb P_\uparrow-\Delta\mathbb T)^{-1} =(\mathrm{Id}-\mathbb V_{\rm e}^0\Delta\mathbb T)^{-1}\mathbb V_{\rm e}^0.
\]
\par
Propositions~\ref{prop:approximateDensityInverse} and~\ref{prop:exactDensityInverse}
prove these assertions on the same weighted density space. 
Applying this inverse to \eqref{eq:overviewDensityDifferenceEquation}
gives the three-error decomposition
\begin{equation}\label{eq:directChainComparisonIdentity}
\rho^{\rm C,\nu}-\rho^\downarrow
=\mathbb V^{\rm C,\nu}\Bigl(
\underbrace{\mathfrak g^{\rm C,\nu}-\mathfrak g^0}
_{\substack{\text{initial-data}\\\text{evolution error}}}
+\underbrace{\Delta\mathbb T\rho^\downarrow}
_{\substack{\text{collision and}\\\text{background-coefficient error}}}
+\underbrace{\mathbb P_\uparrow\rho^\downarrow}
_{\substack{\text{omitted upward}\\\text{interaction}}}
\Bigr).
\end{equation}
The same density solution operator acts on each of the three functions in parentheses.
\par\endgroup
Once $\rho^{\rm C,\nu}=\mathbb V^{\rm C,\nu}\mathfrak g^{\rm C,\nu}$
has been constructed, substitution into \eqref{eq:exactCollisionDuhamel}
gives $H^{\rm C,\nu}$. Taking the density of that formula and using
\eqref{eq:exactDensityVolterra} verifies
$\rho[H^{\rm C,\nu}]=\rho^{\rm C,\nu}$.

\begingroup
{
Constants may depend on the fixed order $r_{\rm x}$.}
 The bound $CN$ for $\mathbb V^{\rm C,\nu}$
controls the full density function and permits the reconstruction of
$H^{\rm C,\nu}$. A sharper bound controls only the terminal
coefficient
\[
B_*^{-1}(\mathbb V^{\rm C,\nu}\mathfrak g)_{m_1}(T_1)
\]
by $C$ times the norm of $\mathfrak g$, in the
same space. Applying this terminal bound to the three functions in
\eqref{eq:directChainComparisonIdentity} avoids the extra factor
$N$ from estimating the whole density function. 
Proposition~\ref{prop:exactDensityInverse} states both estimates.
\par\endgroup
The collision and coefficient changes can now be quantified.
The collision terms in \eqref{eq:exactCollisionModel} are
\[
\nu\bigl(Q(F^{\rm b,\nu},H)+Q(H,F^{\rm b,\nu})\bigr).
\]
Their effect depends on velocity frequency as well as elapsed time.
For physical time $t=s-t_0$, a free-transport phase satisfies
\[
\partial_{v_1}^2e^{ikx_1+i(\xi-kt)v_1}
=-(\xi-kt)^2e^{ikx_1+i(\xi-kt)v_1}.
\]
At the prepared frequency scale $|k|\lesssim N$, $|\xi|\lesssim\zeta$,
and with $NT_1\sim\zeta$, the corresponding accumulated cost obeys
\[
\nu\int_0^{t_\nu}(1+|\xi-kt|^2)\,dt
\lesssim\nu T_1(1+\zeta^2)\sim N^2\nu T_1^3.
\]
Thus observation before the usual enhanced-dissipation time,
$\nu T_1^3\ll1$, does not alone make the collision error small on
functions with velocity Fourier frequency of size $\zeta$. Section~\ref{sec:centeredcollisions}
turns this frequency calculation into localized comparisons for the
Landau equation, accounting for the velocity-dependent coefficients,
the nonlocal collision term and the Gaussian Fourier tails of the
data.

\begingroup
There are also two distinct changes of coefficients. In the collision
terms, the Maxwellian coefficients are replaced by $F^{\rm b,\nu}$;
in the electric coefficient, $\partial_{v_1}F^{\rm b}_{\rm app}$ is replaced by
$\partial_{v_1}F^{\rm b,\nu}$. For the latter difference,
Proposition~\ref{thm:bgcomparison} supplies the relative size
$\nu T_1^3$ in weighted low norms (the Gaussian term is smaller under
\eqref{eq:comparisonLogParameters}). The higher background
bounds and the density-kernel estimates in
Section~\ref{sec:densityComparisonProofs} control the required
velocity derivative and the density trace. Smallness of the
background difference alone is not the terminal density estimate.
\par
\par\endgroup

\begingroup
To separate the collision effects, let $\mathbb U_\mu^\nu(s,\tau)$
be the solution operator with Maxwellian collision coefficients:
\[
\partial_sf+v_1\partial_{x_1}f+\nu\mathcal Lf=0,
\qquad f(\tau)=f_\tau.
\]
\begingroup
\begin{rem}[Shifted commuting derivatives]
For the function $e^{i(\zeta-ks)v_1}V(v)$, the derivative
$\partial_{v_1}+iks-i\zeta$ differentiates $V$ without a factor $\zeta$.
We use $Z_{k,c}(s)=\partial_{v_1}+i(ks-c)$ and establish estimates uniform
in $c$ in Proposition~\ref{prop:centeredpropagation}.
Chaturvedi, Luk and Nguyen use frequency-shifted vector fields in their
collision solution-operator and nonlinear density estimates
\cite[Sections 6.1--6.2]{ref4}. Here the shifts also allow comparison with
the exact inhomogeneous background. For a source with factor
$e^{-il\tau v_1}$ in mode $k$, the appropriate shift is
$c=(k-l)\tau$: it depends on the insertion time $\tau$, but remains fixed
as $s$ evolves. The physical-time version replaces $ks$ by $k(s-t_0)$
and chooses $c$ from the same source phase.
\end{rem}
\par\endgroup
\par\endgroup
The decomposition of $\Delta\mathbb T$ in
\eqref{eq:directResponseDecomposition} separates
\par

\[
\begin{array}{ccl}
\mathbb U^0\longrightarrow\mathbb U_\mu^\nu
&:&\text{Maxwellian collision evolution},\\
\mathbb U_\mu^\nu\longrightarrow\mathbb U_{\rm b}^\nu
&:&\text{exact background collision coefficients},\\
F^{\rm b}_{\rm app}\longrightarrow F^{\rm b,\nu}
&:&\text{coefficient of }E[\rho].
\end{array}
\]

\begingroup
The three source estimates are
\eqref{eq:electricResponseDifferenceBound},
\eqref{eq:exactCollisionInitialComparison}, and
\eqref{eq:descendingReferenceDefect}. With
$\mathcal O_N{\rho}=B_*^{-1}{\rho}_{m_1}({T_1})$, their combination is
\[
\begin{aligned}
\sqrt N\big|\mathcal O_N\mathbb V^{\rm C,\nu}
 (\mathfrak g^{\rm C,\nu}-\mathfrak g^0)\big|
&\lesssim\sqrt N\big(\nu {T_1}(1+\zeta^2)+aL_N\big),\\
\sqrt N\big|\mathcal O_N\mathbb V^{\rm C,\nu}
 \Delta\mathbb T\rho^\downarrow\big|
&\lesssim N^{7/2}\nu {T_1}^3+N^{7/2}aL_N
 +N^{9/2}\nu T_1^3,\\
\sqrt N\big|\mathcal O_N\mathbb V^{\rm C,\nu}
 \mathbb P_\uparrow\rho^\downarrow\big|
&\lesssim e^{-c'N^2}.
\end{aligned}
\]
\par
\par\endgroup
The last bound is the weighted form of the phase-separation estimate
\eqref{eq:upwardTimeShiftSmallness}; the time-shift mechanism was
explained in Section~\ref{sec:roadmapOperators}. In the second line,
the three terms correspond, respectively, to Maxwellian collisions,
the change of collision coefficients, and the background error in the
electric coefficient. 
Combining these bounds gives the following proposition.
\par

\begingroup
\begin{prop}\label{thm:exactCollisionEcho}
{Assume \eqref{eq:comparisonLogParameters}, with $p>{\rm s}_0+7/2$,
and \eqref{eq:parameters}. Fix an integer $r\ge0$ and a velocity weight
$\ell\ge4$. Let $H^{\rm C,\nu}$ be a mild solution of
\eqref{eq:exactCollisionModel} with
{$H^{\rm C,\nu}/\sqrt \mu \in C(I;H^r_{x_1}L^2_\ell)$} and finite
integrated collision dissipation. For sufficiently small $\nu$,
the terminal density satisfies}
\begin{equation}\label{eq:unifiedExactTerminal}
\frac{\sqrt N}{B_*}|\rho_{m_1}[H^{\rm C,\nu}]({T_1})-\rho_{m_1}^{\downarrow}({T_1})|
\le C\big(N^{9/2}\nu T_1^3+N^{7/2}aL_N\big).
\end{equation}

In particular,
\begin{equation}\label{eq:logExactCollisionTerminal}
{
\frac{\sqrt N}{B_*}|\rho_{m_1}[H^{\rm C,\nu}](T_1)-\rho_{m_1}^{\downarrow}(T_1)| \le C[\log(1/\nu)]^{3{\rm s}_0+21/2-3p} +\nu^{1/3-o(1)}\longrightarrow0,
}
\end{equation}
and hence
\begin{equation}\label{eq:logExactCollisionLower}
|\rho_{m_1}[H^{\rm C,\nu}]({T_1})|\ge cB_*/\sqrt N,
\qquad {T_1}=o(\nu^{-1/3}).
\end{equation}
{For sufficiently small $\nu$, a unique solution in the stated class
exists; no restriction is imposed on the support of the Fourier transform in $x$.
The constants and the smallness threshold may depend on the prescribed
fixed orders. The solutions constructed at different orders agree on
their common interval.}
\end{prop}
\par
\par\endgroup

\begingroup
We now check that the right side of \eqref{eq:unifiedExactTerminal}
tends to zero for the scales chosen in
Section~\ref{sec:roadmapBackground}. The coefficient term
$N^{7/2}aL_N=\nu^{1/3-o(1)}$ tends to zero.
The largest remaining contribution is $N^{9/2}\nu T_1^3$;
the next lemma computes this quantity together with the collision
cost on the initial oscillations. All constants, including $K$,
are fixed independently of $\nu$.
\par
\par\endgroup
\begin{lem}\label{lem:logarithmicScales}\label{rem:inflationTime}
Use \eqref{eq:comparisonLogParameters} with fixed ${\rm s}_0\ge11$, $K>0$ and $\gamma_{\rm a}>1/(3K)$; for this scale calculation allow any fixed $p\in\mathbb R$. For sufficiently small $\nu$, the restrictions in \eqref{eq:mainparameters} hold and
\begin{equation}\label{eq:bgtimewindow}
{
{T_1}\mathrel{{\sim}}\nu^{-1/3}[\log(1/\nu)]^{{\rm s}_0+2-p},\qquad \nu {T_1}^3\mathrel{{\sim}}[\log(1/\nu)]^{3{\rm s}_0+6-3p},\qquad \log {T_1}=O(N).
}
\end{equation}
Moreover $t_\nu={T_1}-t_0\mathrel{{\sim}} {T_1}$, and
\begin{equation}\label{eq:prescribedLandau}
\varepsilon_{\rm coll}:=\nu {T_1}(1+\zeta^2)
\mathrel{{\sim}} N^2\nu {T_1}^3\mathrel{{\sim}}[\log(1/\nu)]^{3{\rm s}_0+8-3p}.
\end{equation}
The two terminal costs satisfy
\begin{equation}\label{eq:logObjectiveCollision}
{
\sqrt N\varepsilon_{\rm coll}\mathrel{{\sim}}[\log(1/\nu)]^{3{\rm s}_0+17/2-3p},\qquad N^{9/2}\nu {T_1}^3\mathrel{{\sim}}[\log(1/\nu)]^{3{\rm s}_0+21/2-3p}.
}
\end{equation}
Consequently,
\begin{align}
{T_1}=o(\nu^{-1/3})&\quad\Longleftrightarrow\quad p>{\rm s}_0+2,
\label{eq:preparationTimeRange}\\
\varepsilon_{\rm coll}\to0&\quad\Longleftrightarrow\quad p>{\rm s}_0+8/3,
\label{eq:logCollisionPerturbativeConstraint}\\
\sqrt N\varepsilon_{\rm coll}+N^{9/2}\nu {T_1}^3\to0
&\quad\Longleftrightarrow\quad p>{\rm s}_0+7/2.
\label{eq:retainedComparisonRange}
\end{align}

In particular, writing $p={\rm s}_0+7/2+\eta$ with $\eta>0$,
the largest terminal contribution satisfies
\begin{equation}\label{eq:introTerminalCollisionEstimates}
N^{9/2}\nu T_1^3\sim[\log(1/\nu)]^{-3\eta}.
\end{equation}
\par

For $|k|\le C N$ and $|\xi|\le C\zeta$,
\[
\nu\int_0^{T_1}(1+|\xi-kt|^2)\,dt\le C\varepsilon_{\rm coll}.
\]
\end{lem}
\par

The proof of Lemma~\ref{lem:logarithmicScales} is in
Section~\ref{sec:background}. The condition $p>{\rm s}_0+7/2$
makes the terminal collision errors small; merely observing before
enhanced dissipation requires only $p>{\rm s}_0+2$.
These are sufficient conditions for the present comparison, not
necessary conditions for arbitrary VPL perturbations. 

We have now constructed $H^{\rm C,\nu}$ with the descending-chain
lower bound at $T_1$. The remaining linear term is
$E[F^{\rm b,\nu}]\partial_{v_1}H$; controlling its contribution requires
the separate argument of the next subsection.
\par

\subsection{Background-field transport and the complete linear solution}\label{sec:roadmapFieldCompletion}

We restore the term omitted in \eqref{eq:exactCollisionModel},
while preserving the lower bound just obtained. Define
\begin{equation}\label{eq:fullPhysicalOperatorDefinition}
\mathcal T_{\rm b}H=E[F^{\rm b,\nu}]\partial_{v_1}H,\qquad
\mathcal A_\nu H=\mathcal L_{\rm C}H+\mathcal T_{\rm b}H.
\end{equation}

With this notation, the complete equation
\eqref{eq:fullRangePhysicalEquation} is {$\mathcal A_\nu H^\nu=0$}.
\par

The distribution $H^{\rm C,\nu}$ now has the required density growth,
and we set $H_0=H^{\rm C,\nu}$. The next comparison is with the solution
{$H^\nu$} of \eqref{eq:fullRangePhysicalEquation}.
Since $\mathcal L_{\rm C}H_0=0$, direct comparison gives
\[
\mathcal A_\nu(H^\nu-H_0)=-E[F^{\rm b,\nu}]\partial_{v_1}H_0,
\qquad (H^\nu-H_0)(t_0)=0.
\]
The required error is relative to a specified density coefficient:
\[
|\rho_{m_1}[H^\nu-H_0](T_1)|=o(B_*/\sqrt N).
\]
We explain first why the available general source estimate is too
coarse for this purpose, then which properties of this particular
source make finite corrections effective, and finally how those
corrections are constructed.

The complete inhomogeneous estimate of
Proposition~\ref{thm:fullForcedCoarse} gives, for zero initial data,
\[
\|\rho[W]\|_{\ell^1_{r-1}L^\infty(I)}
\le C_rN^4e^{C_*N}\nu^{-1/2}\sum_{j=0}^{r}\left\|\partial_{x_1}^j\!\left(\mu^{-1/2}\bigl(S\bigr)\right)\right\|_{L^2_{s,x}(I)\Sigma_{2r+6-2j}^*},
\qquad \mathcal A_\nu W=S.
\]
It applies to finite-regularity sources without phase localization.
The exponential is a loss in this estimate, not the lower-bound
amplification $B_*$ or an asserted growth rate of the error.
Removing the mode weights used in the density estimate costs at most
$e^{C\sqrt{a(1+T_1)}}$. Our endpoint gives
\begin{equation}\label{eq:observationGainBalance}
a(1+T_1)\sim\frac{a\zeta}{m_1}=\frac{N^3}{m_1}\sim N^2,
\qquad e^{C\sqrt{a(1+T_1)}}\le e^{C'N}.
\end{equation}
The product calculation is proved in
Proposition~\ref{prop:coarseDensityInverse};
Lemma~\ref{lem:coarseLossVersusEcho} compares this loss with the chain gain.

For an equation such as \eqref{eq:fieldFiniteUnknownError}, ordinary velocity derivatives are difficult to propagate, but the coarse density estimate does not require their pointwise-in-time control. Accordingly, $\mathcal E_r$ contains spatial derivatives and velocity weights only. The velocity derivative in the source is transferred to the test function in the dual collision norm. Since $\mathcal D_r$ includes $\nu$, pairing with integrated dissipation costs $\nu^{-1/2}$ and an additional velocity moment. The hierarchy
\[
\mathcal E_r[f]=\sum_{j=0}^r
\|\langle v\rangle^{2r+6-2j}\partial_{x_1}^jf\|_2^2
\]
provides two additional velocity powers at the same spatial derivative when $r$ is increased by one. We use higher energy levels to control these moment and derivative losses, and estimate the finite iteration at successively lower spatial orders. This uses the joint hierarchy of spatial derivatives and velocity moments; it does not replace a velocity weight by a spatial derivative.

Fix $1/2<\sigma_{\rm int}<1$. Integration by parts in velocity,
as in Lemma~\ref{lem:electricFormDual}, gives
\[
\begin{aligned}
&\nu^{-1/2}\sum_{j=0}^{r}\left\|\partial_{x_1}^j\!\left(\mu^{-1/2}\bigl(E[F^{\rm b,\nu}]\partial_{v_1}H_0\bigr)\right)\right\|_{L^2_{s,x}(I)\Sigma_{2r+6-2j}^*}\\
&\quad\le C_r\nu^{-1/2}\|E[F^{\rm b,\nu}]\|_{L^2_sC_x^r}
\sup_I\mathcal E_{r+1}[H_0/\sqrt\mu]^{1/2}.
\end{aligned}
\]
\begingroup
The field has $L^2_sC_x^r$ norm at most
$a\nu^{\sigma_{\rm int}-o(1)}$, whereas the distribution bound is
$(1+T_1)e^{C_{\rm B}N+o(N)}$.
Since $a(1+T_1)\sim N^2$, the source bound is
$\nu^{\sigma_{\rm int}-1/2-KC_{\rm B}-o(1)}$.
With $N=K\log(1/\nu)$ and $B_*=e^{c_{\rm e}N+O(1)}$, direct comparison therefore gives only
\begin{equation}\label{eq:overviewDirectFieldError}
\frac{\sqrt N}{B_*}|\rho_{m_1}[H^\nu-H_0](T_1)|
\lesssim\nu^{\sigma_{\rm int}-1/2-K(C_*+C_{\rm B}-c_{\rm e})-o(1)}.
\end{equation}
Taking $K$ small makes this bound tend to zero, but also reduces the
amplification. Already at the initial Sobolev order, the quantity required to diverge is
\[
\delta_\nu\zeta^{{\rm s}_0}\frac{B_*}{\sqrt N}
\sim\nu^{1/3-Kc_{\rm e}}[\log(1/\nu)]^{p-1/2},
\qquad \delta_\nu\zeta^{{\rm s}_0}\sim\epsilon_\nu.
\]
Its divergence requires $Kc_{\rm e}\ge1/3$, including the logarithmic
possibility at equality. Lemma~\ref{lem:coarseLossVersusEcho} gives
$C_*>\pi\sqrt3\,c_{\rm e}$ and hence
\[
Kc_{\rm e}\ge\tfrac13
\quad\Longrightarrow\quad
\sigma_{\rm int}-\tfrac12-K(C_*+C_{\rm B}-c_{\rm e})<0.
\]
Thus this general estimate does not both justify the direct comparison
and establish divergence of the quantity required to diverge.
\par
\par\endgroup
The useful structure of the source is a relation between its spatial
modes and velocity oscillations. To see it explicitly, use physical
time $t=s-t_0$ and consider one free-transport component
$e^{imx_1+i(c-mt)v_1}V(v)$, with $V$ smooth and rapidly decreasing.
Write $E_l^{\rm b,\nu}(s)$ for the $l$th Fourier coefficient of
$E[F^{\rm b,\nu}](s)$. A field mode $l$ acting at time $u$ changes the
spatial mode to $k=m+l$. After free transport to time $t$, the phase is
\[
e^{-ik(t-u)v_1}e^{i(c-mu)v_1}=e^{i(c+lu-kt)v_1}.
\]
The corresponding density contribution has the exact expression
\begin{equation}\label{eq:overviewSingleFieldDensity}
\begin{aligned}
&-\int_0^t E_l^{\rm b,\nu}(t_0+u)
\int_{\mathbb R^3}e^{-ik(t-u)v_1}
\partial_{v_1}\big(e^{i(c-mu)v_1}V(v)\big)\,dv\,du\\
&\qquad=-ik\int_0^t(t-u)E_l^{\rm b,\nu}(t_0+u)
\widehat V\big((kt-c-lu)e_1\big)\,du.
\end{aligned}
\end{equation}
This follows by integration by parts in $v_1$, with
$\widehat V(\eta)=\int e^{-i\eta\cdot v}V(v)\,dv$ and
$e_1=(1,0,0)$. Thus the velocity integral is small away from
$kt=c+lu$ when $\widehat V$ decays away from zero. The factor
$k(t-u)$ also shows that the velocity derivative has zero density
at the instant of application; subsequent transport can produce a
nonzero density. This identity illustrates the frequency structure,
not a formula for the collisional solution.

\begingroup
The field shifts the velocity frequency by $lu$. Crucially,
Proposition~\ref{prop:backgroundFieldBounds} proves the matching bound
\begin{equation}\label{eq:overviewFieldEnvelope}
\sum_{l\ne0}e^{\varrho_{\rm F}(|l|^\chi+|lu|^\chi)}
|E_l^{\rm b,\nu}(t_0+u)|
\lesssim a\nu^{\sigma_{\rm int}-o(1)}e^{-c_{\rm dec}u^\chi},
\qquad 0\le u\le t_\nu.
\end{equation}
The constants $\varrho_{\rm F},c_{\rm dec}$ are fixed and positive.
This follows from interpolation between the time-dependent background
comparison and the uniform Gevrey bound, followed by the density trace.
The Gaussian-small explicit field is included in this estimate.
The comparison must be used at each elapsed time $u$, before time
integration; replacing $u$ by $t_\nu$ would lose the power of $\nu$.
\par\endgroup

\begingroup
For $0<\chi<1$ and sufficiently small $\lambda'>0$,
\[
{
|kt-c|^\chi\le|kt-c-lu|^\chi+|lu|^\chi,\qquad e^{-\varrho_{\rm F}|lu|^\chi}e^{-\lambda|kt-c-lu|^\chi} \le e^{-\lambda'|kt-c|^\chi}.
}
\]
The factor $e^{-\varrho_{\rm F}|lu|^\chi}$ therefore compensates for the shift $lu$ in the velocity frequency.
This is the localization used here: the quantitative Fourier estimate in the preceding display, not preservation of a compact
support. A small $L^1_tC^r_x$ norm of the field alone would not give
this conclusion. With collisions, the centered derivative estimates
of Section~\ref{sec:centeredcollisions} and the finite field trace
estimate of Proposition~\ref{prop:finiteFieldTrace} replace the
explicit free-transport calculation. They give the corresponding
spatial-mode and phase decay at smaller positive Gevrey radii.
\par\endgroup

\begingroup
For the quantitative construction, record the integrated field size
\begin{equation}\label{eq:fieldExpansionSizes}
\varepsilon_{\rm F}
=aL_N\nu^{\sigma_{\rm int}},\qquad
\zeta\varepsilon_{\rm F}
=N^3L_N\nu^{\sigma_{\rm int}}
\le\nu^{\sigma_{\rm int}-o(1)}.
\end{equation}
\par\endgroup
The fixed polynomial time moments of the Fourier sum in
\eqref{eq:overviewFieldEnvelope} are bounded by $\varepsilon_{\rm F}$;
see Proposition~\ref{prop:backgroundFieldBounds}.

\begingroup
\emph{These estimates allow us to use the density inverse
$\mathbb V^{\rm C,\nu}$ with the chain weights, where the inverse costs
only $CN$.}
\par The prescribed phases are followed separately for the
initial datum and for each coefficient $\mathcal F_m^{\rm b}(\tau)$; no single velocity frequency is
assigned to the whole distribution. On the prepared datum, a velocity
derivative costs at most $\zeta$ after the weighted time integration.
Since $a\zeta=N^3$, one further field operation contributes
\[
\zeta a\nu^{\sigma_{\rm int}-o(1)}
=N^3\nu^{\sigma_{\rm int}-o(1)}
=\nu^{\sigma_{\rm int}-o(1)}.
\]
The point is to use the phase and mode information before replacing
the source by an ordinary energy norm. Finite corrections do not
create this information; they allow it to be used repeatedly.
\par\endgroup
We implement this procedure through the zero-data increments
\begin{equation}\label{eq:fieldCorrectionRecursion}
\mathcal L_{\rm C}H_n=-E[F^{\rm b,\nu}]\partial_{v_1}H_{n-1},
\qquad H_n(t_0)=0,\qquad 1\le n\le J.
\end{equation}
\par
Each field insertion changes the oscillatory factor from
$e^{i(c-kt)v_1}$ to a sum of factors with different values of $c+lu$, as in \eqref{eq:overviewSingleFieldDensity}. Consequently, a single fixed parameter $c$ does not describe all components of $H_n$ without an additional cost. Applying a coarse estimate at that same parameter at every step would discard the frequency information needed to obtain a new power of $\nu$.
We therefore expand the transport terms in ordered Duhamel integrals,
keeping the frequency shifts at their actual integration times.
The first two corrections explain how this leads to the density
operators used in Section~\ref{sec:fieldExpansion}.

Set $f_n=H_n/\sqrt\mu$, $f_{\rm in}=H_{\rm in}/\sqrt\mu$ and
$\rho^{(n)}=\rho[H_n]$, so that $\rho^{(0)}=\rho^{\rm C,\nu}$.
Recall the collision evolution $\mathbb U_{\rm b}^\nu(s,\tau)$ from
\eqref{eq:exactCollisionEvolution}. The normalized transport term and
the self-consistent electric source are
\[
\begin{aligned}
\widetilde{\mathcal T}_{\rm b}(s)f
&=-\mu^{-1/2}\mathcal T_{\rm b}(s)(\sqrt\mu f)
=-E[F^{\rm b,\nu}](s)\left(\partial_{v_1}-\frac{v_1}{2}\right)f,\\
-\mu^{-1/2}E[\rho]\partial_{v_1}F^{\rm b,\nu}
&=\sum_{m\ne0}\mathcal F_m^{\rm b}\rho_m,\qquad
[\mathcal F_m^{\rm b}]_j
=\frac{i}{m}\mu^{-1/2}\partial_{v_1}F_{j-m}^{\rm b,\nu}.
\end{aligned}
\]
Thus \eqref{eq:fieldCorrectionRecursion}, with zero initial datum,
gives
\[
f_n(s)=\int_{t_0}^s\mathbb U_{\rm b}^\nu(s,u)\widetilde{\mathcal T}_{\rm b}(u)f_{n-1}(u)\,du
+\sum_{m\ne0}\int_{t_0}^s\mathbb U_{\rm b}^\nu(s,\tau)
\mathcal F_m^{\rm b}(\tau)\rho_m^{(n)}(\tau)\,d\tau.
\]
The term with $\rho^{(n)}$ must be retained: the self-consistent
electric term belongs to $\mathcal L_{\rm C}$, whereas
$\mathbb U_{\rm b}^\nu$ contains only transport and collisions.

For $q\ge0$, use the ordered evolution operators
\[
\begin{aligned}
\mathbb U_{\rm b}^{(0)}(s,\tau)&=\mathbb U_{\rm b}^\nu(s,\tau),\\
\mathbb U_{\rm b}^{(q)}(s,\tau)f
&=\int_\tau^s\mathbb U_{\rm b}^\nu(s,u)\widetilde{\mathcal T}_{\rm b}(u)
\mathbb U_{\rm b}^{(q-1)}(u,\tau)f\,du,\qquad q\ge1.
\end{aligned}
\]
In particular, the first two are
\[
\begin{aligned}
\mathbb U_{\rm b}^{(1)}(s,\tau)f
&=\int_\tau^s\mathbb U_{\rm b}^\nu(s,u_1)\widetilde{\mathcal T}_{\rm b}(u_1)
\mathbb U_{\rm b}^\nu(u_1,\tau)f\,du_1,\\
\mathbb U_{\rm b}^{(2)}(s,\tau)f
&=\int_{\tau\le u_2\le u_1\le s}
\mathbb U_{\rm b}^\nu(s,u_1)\widetilde{\mathcal T}_{\rm b}(u_1)
\mathbb U_{\rm b}^\nu(u_1,u_2)\widetilde{\mathcal T}_{\rm b}(u_2)
\mathbb U_{\rm b}^\nu(u_2,\tau)f\,du_2\,du_1.
\end{aligned}
\]
The chronological order of the factors is essential; the coefficients
and the shifted derivatives depend on their respective times.

\emph{The first correction.}
The base solution satisfies
\[
f_0(s)=\mathbb U_{\rm b}^\nu(s,t_0)f_{\rm in}
+\sum_{m\ne0}\int_{t_0}^s\mathbb U_{\rm b}^\nu(s,\tau)
\mathcal F_m^{\rm b}(\tau)\rho_m^{(0)}(\tau)\,d\tau.
\]
Substituting this expression for $f_0(u)$ in the equation for $f_1$
and interchanging the ordered integrals gives
\[
\begin{aligned}
\frac{H_1(s)}{\sqrt\mu}=f_1(s)
={}&\mathbb U_{\rm b}^{(1)}(s,t_0)f_{\rm in} +\sum_{m\ne0}\int_{t_0}^s\mathbb U_{\rm b}^{(1)}(s,\tau)
\mathcal F_m^{\rm b}(\tau)\rho_m^{(0)}(\tau)\,d\tau\\
&+\sum_{m\ne0}\int_{t_0}^s\mathbb U_{\rm b}^\nu(s,\tau)
\mathcal F_m^{\rm b}(\tau)\rho_m^{(1)}(\tau)\,d\tau.
\end{aligned}
\]
For example, the second term results from the identity
\[
\int_{t_0}^s\!\int_{t_0}^{u}
\mathbb U_{\rm b}^\nu(s,u)\widetilde{\mathcal T}_{\rm b}(u)
\mathbb U_{\rm b}^\nu(u,\tau)
\mathcal F_m^{\rm b}(\tau)\rho_m^{(0)}(\tau)\,d\tau\,du =\int_{t_0}^s\mathbb U_{\rm b}^{(1)}(s,\tau)
\mathcal F_m^{\rm b}(\tau)\rho_m^{(0)}(\tau)\,d\tau.
\]

\emph{The second correction.}
Substituting the displayed formula for $f_1(u)$ into the equation
for $f_2$ adds one transport insertion to each of its terms:
\[
\begin{aligned}
\frac{H_2(s)}{\sqrt\mu}=f_2(s)
={}&\mathbb U_{\rm b}^{(2)}(s,t_0)f_{\rm in} +\sum_{m\ne0}\int_{t_0}^s\mathbb U_{\rm b}^{(2)}(s,\tau)
\mathcal F_m^{\rm b}(\tau)\rho_m^{(0)}(\tau)\,d\tau\\
&+\sum_{m\ne0}\int_{t_0}^s\mathbb U_{\rm b}^{(1)}(s,\tau)
\mathcal F_m^{\rm b}(\tau)\rho_m^{(1)}(\tau)\,d\tau\\
&+\sum_{m\ne0}\int_{t_0}^s\mathbb U_{\rm b}^\nu(s,\tau)
\mathcal F_m^{\rm b}(\tau)\rho_m^{(2)}(\tau)\,d\tau.
\end{aligned}
\]

\emph{General correction and density recursion.}
Induction gives the exact finite identity
\[
\frac{H_n(s)}{\sqrt\mu}
=\mathbb U_{\rm b}^{(n)}(s,t_0)f_{\rm in}
+\sum_{q=0}^n\sum_{m\ne0}\int_{t_0}^s
\mathbb U_{\rm b}^{(q)}(s,\tau)\mathcal F_m^{\rm b}(\tau)
\rho_m^{(n-q)}(\tau)\,d\tau.
\]
Here $q$ counts the transport insertions after the electric source
at time $\tau$; the remaining $n-q$ insertions belong to the density
$\rho^{(n-q)}$. Multiplying by $\sqrt\mu$ and integrating in velocity
leads to the functions and density operators
\[
\begin{aligned}
\mathfrak g_{q,k}(s)
&=\int_{\mathbb R^3}\sqrt\mu\,
[\mathbb U_{\rm b}^{(q)}(s,t_0)f_{\rm in}]_k\,dv,\\
(\mathbb T_q\rho)_k(s)
&=\sum_{m\ne0}\int_{t_0}^s\!\int_{\mathbb R^3}\sqrt\mu\,
[\mathbb U_{\rm b}^{(q)}(s,\tau)\mathcal F_m^{\rm b}(\tau)]_k
\rho_m(\tau)\,dv\,d\tau,\qquad q\ge1.
\end{aligned}
\]
These are the operators in \eqref{eq:fieldDensityCoefficients}.
The $q=0$ integral is already
$(\mathbb K+\mathbb P_\downarrow+\mathbb P_\uparrow+\Delta\mathbb T)\rho^{(n)}$.
Therefore the same inverse $\mathbb V^{\rm C,\nu}$ used for $H_0$ gives
\[
\begin{aligned}
\rho^{(1)}&=\mathbb V^{\rm C,\nu}
 (\mathfrak g_1+\mathbb T_1\rho^{(0)}),\qquad \rho^{(2)}=\mathbb V^{\rm C,\nu}
 (\mathfrak g_2+\mathbb T_1\rho^{(1)}+\mathbb T_2\rho^{(0)}),\\
\rho^{(n)}&=\mathbb V^{\rm C,\nu}
 \left(\mathfrak g_n+\sum_{q=1}^n\mathbb T_q\rho^{(n-q)}\right),
 \qquad 1\le n\le J.
\end{aligned}
\]
Thus one first determines $\rho^{(n)}$ from the preceding densities
and then reconstructs $H_n$ by the distribution formula above.
The ordered integral estimates and the justification of these identities
by smooth approximation and energy convergence are given in
Section~\ref{sec:fieldExpansion}.
Lemma~\ref{lem:oneFieldOperation} controls each change of the shifted
vector field; Proposition~\ref{prop:finiteFieldTrace} records the losses
after $q$ insertions. Its proof retains the times $u_1,\ldots,u_q$
before applying the Fourier bound \eqref{eq:fieldFourierEnvelope}.
Section~\ref{sec:fieldExpansion} estimates the ordered integrals using
\eqref{eq:overviewFieldEnvelope} and the shifted-vector-field collision bounds,
and proves
\[
\frac{\sqrt N}{B_*}|\rho_{m_1}[H_n](T_1)|
\le\nu^{n\sigma_{\rm int}-o(1)},\qquad 1\le n\le J.
\]

The derivative losses are controlled using a finite hierarchy of spaces.
In normalized variables,
\[
\frac{E[F^{\rm b,\nu}]\partial_{v_1}H_{n-1}}{\sqrt \mu }
=E[F^{\rm b,\nu}]\left(\partial_{v_1}-\frac{v_1}{2}\right)
\frac{H_{n-1}}{\sqrt \mu }.
\]
Decreasing the Gevrey radius controls the velocity derivative, and one
additional velocity moment controls multiplication by $v_1$.
For a desired spatial density order $r_{\rm x}$, we use the orders
\[
r_{\rm x}+J,\quad r_{\rm x}+J-1,\quad\ldots,\quad r_{\rm x}.
\]
The Gaussian datum and the background estimates supply these finitely
many stronger norms. 

The corrected approximation is $H^{[J]}=\sum_{n=0}^JH_n$;
$H_n$ alone is an increment. Summing the equations gives
\[
\mathcal A_\nu H^{[J]}
=E[F^{\rm b,\nu}]\partial_{v_1}
\left(\sum_{n=0}^JH_n-\sum_{n=1}^JH_{n-1}\right)
=E[F^{\rm b,\nu}]\partial_{v_1}H_J.
\]
The residual contains $J+1$ field operations. We estimate its last
velocity derivative by integration by parts, using
\begin{equation}\label{eq:overviewFieldResidualSource}
\nu^{-1/2}\sum_{j=0}^{r}\left\|\partial_{x_1}^j\!\left(\mu^{-1/2}\bigl(E[F^{\rm b,\nu}]\partial_{v_1}H_J\bigr)\right)\right\|_{L^2_{s,x}(I)\Sigma_{2r+6-2j}^*}
\le e^{C_{\rm B}N+o(N)}\frac{1+T_1}{\sqrt\nu\,\zeta}
(\zeta\varepsilon_{\rm F})^{J+1}.
\end{equation}
Proposition~\ref{prop:fieldQuantitativeInputs} proves this bound from
the field's $L^2$ time estimate and the spatial energy of $H_J$.
Only this residual is passed to the complete inhomogeneous estimate:
\[
\mathcal A_\nu W_J=-E[F^{\rm b,\nu}]\partial_{v_1}H_J,
\qquad W_J(t_0)=0,\qquad H^\nu=H^{[J]}+W_J.
\]
Since $(1+T_1)/\zeta=\nu^{o(1)}$, we obtain
\[
\|\rho[W_J]\|_{\ell^1_{r-1}L^\infty(I)}
\le\nu^{(J+1)\sigma_{\rm int}-1/2-K(C_*+C_{\rm B})-o(1)}.
\]
After fixing $K,C_*,C_{\rm B}$ and $\sigma_{\rm int}$, choose $J$ independently of $\nu$ so that
\begin{equation}\label{eq:logFieldOrderChoice}
\iota_{\rm F}:=(J+1)\sigma_{\rm int}-\frac12-K(C_*+C_{\rm B})>1.
\end{equation}
The constants $C_*,C_{\rm B}$ are independent of $J$ and of the fixed
regularity orders. Thus finitely many corrections compensate for the
fixed power $e^{C_*N}=\nu^{-KC_*}$.
Continuing the chain to a fixed low mode would instead give
$T_1\sim\zeta$ and the cost $e^{CN^{3/2}}$ in the same argument;
no fixed positive power of $\nu$ compensates for that cost.
All radii, moments and spatial orders needed through $J$ are fixed
before $\nu\downarrow0$. The following proposition records the construction.
\par

\begin{prop}\label{thm:finiteFieldExpansion}\label{lem:fieldErrorIdentity}\label{prop:quantitativeFieldResidual}
Assume the parameters of Proposition~\ref{thm:exactCollisionEcho}.
Fix $J\ge1$ and $1/2<\sigma_{\rm int}<1$, independently of $\nu$.
Set $H_0=H^{\rm C,\nu}$ and let $H_n$, $1\le n\le J$,
solve the recurrence \eqref{eq:fieldCorrectionRecursion}.
\par
{Each initial-value problem has a unique finite-energy solution.
Define the corrected approximations by}
\[
H^{[n]}=\sum_{j=0}^nH_j,\qquad 0\le n\le J.
\]
They satisfy
\begin{equation}\label{eq:fieldPartialSumIdentities}
\begin{aligned}
\mathcal L_{\rm C}H^{[n]}&=-\mathcal T_{\rm b}H^{[n-1]}
&& (1\le n\le J),\\
\mathcal A_\nu H^{[n]}&=\mathcal T_{\rm b}H_n,\qquad
H^{[n]}(t_0)=H_{\rm in}
&& (0\le n\le J).
\end{aligned}
\end{equation}
At the final order,
\begin{equation}\label{eq:fieldApproximationToBase}
\frac{\sqrt N}{B_*}|\rho_{m_1}[H^{[J]}-H^{\rm C,\nu}]({T_1})|
\le\nu^{{\sigma_{\rm int}}-o(1)},\qquad
|\rho_{m_1}[H^{[J]}]({T_1})|\ge cB_*/\sqrt N.
\end{equation}
If {$W_J=H^\nu-H^{[J]}$}, where {$H^\nu$} solves
\eqref{eq:fullRangePhysicalEquation} with the same datum, then
\begin{equation}\label{eq:fieldFiniteUnknownError}
\mathcal A_\nu W_J=-\mathcal T_{\rm b}H_J,
\qquad W_J(t_0)=0.
\end{equation}
All assertions concern a fixed finite $J$. If
\eqref{eq:logFieldOrderChoice} holds, the remainder satisfies, for every
fixed $r\ge2$ and sufficiently small $\nu$ depending on that order,
\begin{equation}\label{eq:commonFieldCompletionError}
\|\rho[H^\nu-H^{[J]}]\|_{\ell^1_{r-1}L^\infty(I)}
\le\nu^{\iota_{\rm F}-o(1)}.
\end{equation}
The complete solution is taken in the order-$r$ energy class.
\end{prop}

Section~\ref{sec:fieldExpansion} proves the finite-order estimates and
applies Proposition~\ref{thm:fullForcedCoarse} only to the final residual.

 {Finite-order approximate-solution constructions in other settings include nonlinear}
instability constructions \cite{Grenier2000,DesjardinsGrenier2003,GrenierNguyen2024}
and truncated Hilbert expansions for kinetic equations
\cite{Caflisch1980,GuoJang2010,LeiLiuXiaoZhao2022}. 

Here the correction is entirely at the \emph{linear} level: with the
background fixed, every $H_n$ is linear in $H_{\rm in}$, and $n$ counts
insertions of $E[F^{\rm b,\nu}]\partial_{v_1}$. It does not count powers
of the perturbation amplitude in a nonlinear instability expansion,
or powers of a Knudsen number in a hydrodynamic Hilbert expansion.
The issue specific to this construction is to obtain the factor
$\zeta\varepsilon_{\rm F}\le\nu^{\sigma_{\rm int}-o(1)}$ per insertion
on a long time interval, despite the large velocity frequency and
the derivative losses, while proving that each correction has terminal
density $o(B_*/\sqrt N)$. The phase-dependent field decay and the
finite regularity hierarchy provide these estimates. 
\par

The complete comparison can now be displayed in one identity:
\begin{equation}\label{eq:linearConstructionThreeErrors}
\begin{aligned}
\rho[H^\nu]=\rho^\downarrow
&+\underbrace{\mathbb V^{\rm C,\nu}
  (\mathfrak g^{\rm C,\nu}-\mathfrak g^0)}_{\text{initial-data evolution}} +\underbrace{\mathbb V^{\rm C,\nu}\mathbb P_\uparrow
  \rho^\downarrow}_{\text{upward interaction}}
+\underbrace{\mathbb V^{\rm C,\nu}\Delta\mathbb T
  \rho^\downarrow}_{\text{collision and electric coefficients}}\\
&+\underbrace{\sum_{n=1}^J\rho[H_n]}_{\text{finite background-field corrections}}
+\underbrace{\rho[W_J]}_{\text{complete linear remainder}}.
\end{aligned}
\end{equation}
Every term after $\rho^\downarrow$ is $o(B_*/\sqrt N)$ at $(m_1,{T_1})$.

The notation {$W_J$} is reserved for this linear completion; {$W$} is the
nonlinear remainder in \eqref{eq:nonlinearFirstOrderExpansion}.

We state the complete linear result in the energy and comparison
spaces defined in Section~\ref{sec:notation}. The energy norms include
all spatial modes; the density norms below concern the nonzero modes
that generate the electric field.

\begingroup
\begin{thm}[The complete linear echo at logarithmic amplitudes]\label{thm:logFullLinearEcho}
Fix ${\rm s}_0\ge11$, $K>0$, $p>{\rm s}_0+7/2$, and $1/2<\chi<2/3$. Use \eqref{eq:comparisonLogParameters} and \eqref{eq:parameters}, and let $F^{\rm b,\nu}$ be the exact background of Proposition~\ref{prop:bguniform}. For all sufficiently small $\nu>0$, the complete physical equation \eqref{eq:fullRangePhysicalEquation}, with initial datum {$H_{\rm in}$}, has a unique finite-energy solution {$H^\nu$}. More precisely, for every fixed spatial order $r\ge2$ there is $\nu_r>0$ such that, for $0<\nu<\nu_r$, {{$H^\nu/\sqrt \mu $} is continuous} in $\mathscr E_r$ and has finite integrated collision dissipation at that order. The solutions obtained at different orders agree by energy uniqueness. {The terminal density satisfies}
\begin{equation}\label{eq:logFullLinearComparison}
{
\frac{\sqrt N}{B_*} \left|\rho_{m_1}[H^\nu](T_1)-\rho_{m_1}^{\downarrow}(T_1)\right| \le{}C[\log(1/\nu)]^{3{\rm s}_0+21/2-3p} +\nu^{1/3-o(1)}\longrightarrow0,
}
\end{equation}
In particular,
\begin{equation}\label{eq:logFullLinearLower}
|\rho_{m_1}[H^\nu]({T_1})|\ge cB_*/\sqrt N,
\qquad {T_1}-t_0=o(\nu^{-1/3}).
\end{equation}
Both exact collision coefficients, the electric term in the density equation to {$H^\nu$}, and background field transport are included.

{If also $Kc_{\rm e}>4/3$, the amplitudes
$\epsilon_\nu,\delta_\nu$ in \eqref{eq:nonlinearTwoScales} satisfy}
\begin{equation}\label{eq:logLinearObservable}
\left\|\delta_\nu H_{\rm in}/\sqrt \mu \right\|_{H^{{\rm s}_0}_{q,\ell}}
\le C\epsilon_\nu,\qquad
\langle m_1,m_1({T_1}-t_0)\rangle^{{\rm s}_0-3}
|\delta_\nu\rho_{m_1}[H^\nu]({T_1})|
\ge\nu^{4/3-Kc_{\rm e}+o(1)}\longrightarrow\infty.
\end{equation}
Here $0<q<1/4$ and $\ell\ge4$ are fixed weights. The nonlinear conclusion for the normalized initial datum is established in Theorem~\ref{thm:logNonlinearInflation}.
\end{thm}
\par\endgroup
The proof combines the preceding propositions in
Section~\ref{sec:fieldExpansion}. The nonlinear comparison also requires
upper bounds for the density and spatial energy of the same solution.

\begin{prop}\label{prop:linearNonlinearInput}
Use \eqref{eq:comparisonLogParameters} and the fixed correction order
from \eqref{eq:logFieldOrderChoice}. Let $H^\nu$ be the complete linear
solution of Theorem~\ref{thm:logFullLinearEcho}. For every fixed integer
${\rm s}\ge3$ and sufficiently small $\nu$, depending on that order,
\begin{equation}\label{eq:linearHierarchyInput}
\|H^\nu\|_{\mathscr V_{\rm s}(I)}\le A_{\rm s}(N)\widetilde B,
\qquad A_{\rm s}(N)\ge1,\qquad \log A_{\rm s}(N)=o(N).
\end{equation}
\end{prop}

{Section~\ref{sec:fieldExpansion} proves this estimate.} The factor
$A_{\rm s}$ may depend on the fixed order, whereas the exponential
constant $C_{\rm B}$ does not. The Sobolev order can therefore be chosen
after the losses in the nonlinear comparison have been fixed.
\par

\subsection{Nonlinear construction and norm inflation}\label{sec:roadmapNonlinear}

The complete linear solution {$H^\nu$} has now been constructed and
satisfies the echo lower bound. We next construct an exact nonlinear
solution with the amplitude $\delta_\nu$ fixed in \eqref{eq:nonlinearTwoScales}, of the form
\begin{equation}\label{eq:nonlinearFirstOrderExpansion}
F^\nu=F^{\rm b,\nu}+\delta_\nu H^\nu+W,
\qquad \mathcal A_\nu H^\nu=0.
\end{equation}
 The comparison must give
\[
|\rho_{m_1}[W](T_1)|=o(\delta_\nu B_*/\sqrt N).
\]

The initial correction was fixed in Lemma~\ref{lem:data}. Thus
\begingroup
\begin{equation}\label{eq:roadmapInitialDecomposition}
F_{\rm in}
=F^{\rm b}_{\rm in}+\delta_\nu H_{\rm in}+W_{\rm in},
\qquad W(t_0)=W_{\rm in},
\end{equation}
\endgroup
\begingroup
where $W_{\rm in}$ is the change from $\theta$ to ${\theta_{\delta_\nu}}$
\eqref{eq:initialTemperatureCorrection}. The bound
\eqref{eq:initialTemperatureSpatialBound} makes this initial error
negligible even after the complete inhomogeneous estimate.
\par
\par\endgroup
To derive the equation for {$W$}, subtract the exact background VPL
equation from the equation for $F^\nu$. Poisson's law gives
{$E[F^\nu]-E[F^{\rm b,\nu}]=E[\delta_\nu H^\nu+W]$}, and the
bilinearity of $Q$ gives
\[
\mathcal A_\nu(\delta_\nu H^\nu+W)
={}-E[\delta_\nu H^\nu+W]\partial_{v_1}(\delta_\nu H^\nu+W) +\nu Q(\delta_\nu H^\nu+W,\delta_\nu H^\nu+W).
\]
Using {$\mathcal A_\nu H^\nu=0$} and expanding the right side now yields
\[
\begin{aligned}
\mathcal A_\nu W
={}&\underbrace{\delta_\nu^2\big(\nu Q(H^\nu,H^\nu)
-E[H^\nu]\partial_{v_1}H^\nu\big)}_{\text{known quadratic source}} -E[\delta_\nu H^\nu+W]\partial_{v_1}W
-E[W]\partial_{v_1}(\delta_\nu H^\nu)\\
&+\nu Q(\delta_\nu H^\nu+W,W)+\nu Q(W,\delta_\nu H^\nu).
\end{aligned}
\]

For the subsequent estimates, write {$H_{\rm a}=\delta_\nu H^\nu$} and
use the bilinear expression
\begin{equation}\label{eq:nonlinearBilinearOperator}
\mathbb Q_\nu({F_1},{F_2})=-E[{F_1}]\partial_{v_1}{F_2}+\nu Q({F_1},{F_2}).
\end{equation}

The expanded equation above is equivalently
\begingroup
\begin{equation}\label{eq:nonlinearExpansionRoadmap}
\begin{aligned}
\mathcal R_{\rm a}
&=\mathcal A_\nu H_{\rm a}-\mathbb Q_\nu(H_{\rm a},H_{\rm a})
=-\delta_\nu^2\mathbb Q_\nu(H^\nu,H^\nu),\\
\mathcal A_\nu W
&=-\mathcal R_{\rm a}+\mathbb Q_\nu(H_{\rm a},W)
+\mathbb Q_\nu(W,H_{\rm a})+\mathbb Q_\nu(W,W),\\
W(t_0)&=F_{\rm in}-F^{\rm b,\nu}(t_0)-\delta_\nu H_{\rm in}.
\end{aligned}
\end{equation}
\par\endgroup
We integrate the electric and collision terms by parts in velocity,
then use Cauchy--Schwarz and Young's inequalities to absorb the derivative
of the energy test function into the collision dissipation.
{The resulting bounds use the {source integrals} in \eqref{eq:quadraticSourceSpaces},}
defined in Section~\ref{sec:notation}. Lemma~\ref{lem:electricFormDual}
and the collision estimate give
\[
\begin{aligned}
\nu^{-1/2}\sum_{j=0}^{\rm s}\left\|\partial_{x_1}^j\!\left(\mu^{-1/2}\bigl(\delta_\nu^2E[H^\nu]\partial_{v_1}H^\nu\bigr)\right)\right\|_{L^2_{s,x}(I)\Sigma_{2{\rm s}+6-2j}^*}
&\lesssim\delta_\nu^2\nu^{-1/2}(1+T_1)^{3/2}e^{2C_{\rm B}N+o(N)},\\
\nu^{-1/2}\sum_{j=0}^{\rm s}\left\|\partial_{x_1}^j\!\left(\mu^{-1/2}\bigl(\nu\delta_\nu^2Q(H^\nu,H^\nu)\bigr)\right)\right\|_{L^2_{s,x}(I)\Sigma_{2{\rm s}+6-2j}^*}
&\lesssim\delta_\nu^2(1+T_1)^2e^{2C_{\rm B}N+o(N)}.
\end{aligned}
\]
The first estimate uses the spatial energy of $H^\nu$ at order
${\rm s}+1$, already supplied by Proposition~\ref{prop:linearNonlinearInput}.
Proposition~\ref{prop:quadraticApproximation} proves the resulting residual bound.
\par

The scale of $\delta_\nu$ explains why a first-order approximation
can suffice. For the prescribed initial Sobolev size
$\epsilon_\nu=\nu^{1/3}[\log(1/\nu)]^p$, the chosen frequencies give
\[
\delta_\nu\sim\epsilon_\nu\zeta^{-{\rm s}_0}
=\nu^{({\rm s}_0+1)/3+o(1)},\qquad
\zeta=\nu^{-1/3+o(1)},\qquad 1+T_1=\nu^{-1/3+o(1)}.
\]
Applying the same complete inhomogeneous estimate as in
Section~\ref{sec:roadmapFieldCompletion}, now to the quadratic source,
and dividing by the linear density lower-bound scale $\delta_\nu B_*/\sqrt N$ gives
\par

\[
\frac{{\delta_\nu^2\nu^{-1/2}(1+T_1)^{3/2}e^{(C_*+2C_{\rm B})N+o(N)}}}
{\delta_\nu B_*/\sqrt N}
=\nu^{({\rm s}_0-2)/3-K(C_*+2C_{\rm B}-c_{\rm e})-o(1)}.
\]
The source is quadratic in $\delta_\nu$, while the density lower-bound scale $\delta_\nu B_*/\sqrt N$ is linear in $\delta_\nu$. The remaining factor of $\delta_\nu$
{absorbs the dual-source, time and amplification losses when}
\begin{equation}\label{eq:nonlinearRoadmapMargin}
\Theta=\frac{{\rm s}_0-2}{3}-K(C_*+2C_{\rm B}-c_{\rm e})>0.
\end{equation}
The constants on the right are fixed before the Sobolev order.
This is the reason the nonlinear approximation here stops at first
order, whereas the echo construction of Deng and Masmoudi
\cite{ref1} uses a deep nonlinear Taylor expansion. The finite
field corrections in the preceding part all belong to the linear
solution {$H^\nu$} and do not increase this nonlinear expansion order.

It remains to solve the equation for {$W$}, since the quadratic-source
bound alone does not control the terms containing {$W$}. We keep the
derivative of the unknown remainder in the transport operator and
keep the collision terms in the linear equation at each step.
Starting with {$W^{(0)}=0$}, solve successively
\begingroup
\[
\begin{aligned}
\mathcal A_\nu W^{(n+1)}
&+E[\delta_\nu H^\nu+W^{(n)}]\partial_{v_1}W^{(n+1)}
+E[W^{(n+1)}]\partial_{v_1}(\delta_\nu H^\nu)\\
&-\nu\big(Q(\delta_\nu H^\nu+W^{(n)},W^{(n+1)})
+Q(W^{(n+1)},\delta_\nu H^\nu)\big)\\
&=\delta_\nu^2\big(\nu Q(H^\nu,H^\nu)
-E[H^\nu]\partial_{v_1}H^\nu\big),
\qquad W^{(n+1)}(t_0)=W_{\rm in}.
\end{aligned}
\]
\endgroup
For a fixed {$W^{(n)}$}, this is a linear equation for {$W^{(n+1)}$}.

Proposition~\ref{thm:fullForcedCoarse} solves this equation.

The construction is bounded at order ${\rm s}$ and contracts at
order ${\rm s}-1$. Subtracting consecutive equations leaves, on the
right, the difference source
\[
-E[W^{(n)}-W^{(n-1)}]\partial_{v_1}W^{(n)}
+\nu Q(W^{(n)}-W^{(n-1)},W^{(n)}),\qquad n\ge1.
\]
The order-${\rm s}$ energy of the bounded iterates controls
the difference source at order ${\rm s}-1$.
\par
Section~\ref{sec:nonlinearConstruction} proves, uniformly in $n$,
\[
\begin{aligned}
\|W^{(n)}\|_{\mathscr V_{\rm s}}
&\le\delta_\nu e^{C_{\rm B}N+o(N)},\\
\|W^{(n+1)}-W^{(n)}\|_{\mathscr V_{{\rm s}-1}}
&\le{\nu^{\Theta+K(C_{\rm B}-c_{\rm e})-o(1)}}
\|W^{(n)}-W^{(n-1)}\|_{\mathscr V_{{\rm s}-1}}.
\end{aligned}
\]
The second inequality holds for $n\ge1$. The common subexponential
factors may depend on the fixed order, but not on $n$.
Since $C_{\rm B}>c_{\rm e}$, the preceding condition
for the relative source error also makes the contraction factor tend
to zero. Passing to the limit gives the exact remainder; the sharper
source bound then gives
\[
\frac{|\rho_{m_1}[W](T_1)|}{\delta_\nu B_*/\sqrt N}
\le{\nu^{\Theta-o(1)}}\longrightarrow0.
\]
The Gaussian-small initial correction is included in these estimates.
Higher fixed regularity orders and the exact VPL equation yield the
classical solution, nonnegativity and conservation; these steps are
completed in Sections~\ref{sec:classicalConclusion}
and~\ref{sec:nonlinearNormInflation}.

For each fixed integer ${\rm s}\ge4$, choose $A_{\rm s}(N)\ge1$ with
$\log A_{\rm s}(N)=o(N)$, dominating the factor in
Proposition~\ref{thm:fullForcedCoarse}, the order-$({\rm s}+1)$ factor in
Proposition~\ref{prop:linearNonlinearInput}, and the estimates below, and put
\begin{equation}\label{eq:additionalSpatialFactors}
\eta_{\nu,{\rm s}}=2\delta_\nu A_{\rm s}(N)^2\widetilde B e^{C_*N}
\big((1+T_1)^2+\nu^{-1/2}(1+T_1)^{3/2}\big).
\end{equation}
\par

All enlargements of $A_{\rm s}$ are made simultaneously in the ball radius and in \eqref{eq:additionalSpatialFactors}, and the factors are chosen nondecreasing in the fixed integer ${\rm s}$. They include only fixed constants and subexponential factors supplied by the linear estimates.

\begingroup
\begin{prop}\label{thm:exactNonlinearSolution}\label{prop:logNonlinearCalibration}
Use \eqref{eq:comparisonLogParameters} and \eqref{eq:nonlinearTwoScales}.
Fix ${\rm s}\ge4$, ${\rm s}_0\ge11$, $p>{\rm s}_0+7/2$ and $K>0$, with
$\Theta>0$ as in \eqref{eq:nonlinearRoadmapMargin}.
For the initial datum $F_{\rm in}$ from
\eqref{eq:positivefamily} and sufficiently small $\nu$, there is a real
distributional VPL solution
\[
F^\nu=F^{\rm b,\nu}+H_{\rm a}+W
\]
on all of $I$. {The remainder {$W$} satisfies}
\begin{equation}\label{eq:exactNonlinearClass}
W/\sqrt \mu \in C(I;\mathscr E_{\rm s}),\qquad
\int_I\mathcal D_{\rm s}[W/\sqrt \mu ]\,ds<\infty,\qquad
\rho[W]\in\ell^1_{{\rm s}-1}L^\infty(I),
\end{equation}
and
\begin{equation}\label{eq:nonlinearComparisonClosed}
\|W\|_{\mathscr V_{\rm s}}
\le C_{\rm s}\delta_\nu A_{\rm s}(N)\widetilde{B}\eta_{\nu,{\rm s}}.
\end{equation}
It is unique among solutions in this class with
{$\|W\|_{\mathscr V_{\rm s}}\le 2\delta_\nu A_{\rm s}(N)\widetilde{B}$}.
For the same datum, the constructions at any two fixed orders agree once
$\nu$ is sufficiently small for both orders. With $\varepsilon_N$
from Proposition~\ref{lem:fixedCoefficientOrders}, moreover,
\[
\sup_I\mathcal E_{\rm s}[f^\nu]^{1/2}
\le C_{\rm s}\varepsilon_N+2C_{\rm s}(1+{T_1})\delta_\nu A_{\rm s}(N)\widetilde{B}=o(1),\qquad \int_I\mathcal D_{\rm s}[f^\nu]\,ds<\infty.
\]
For every fixed $0<\omega<\Theta$, the terminal comparison is
\begin{equation}\label{eq:nonlinearRelativeComparison}
\frac{|\rho_{m_1}[W]({T_1})|}{\delta_\nu B_*/\sqrt N}
\le\nu^{\Theta-o(1)}\le C_\omega\nu^\omega\longrightarrow0.
\end{equation}
Consequently,
\begin{equation}\label{eq:exactNonlinearEcho}
\begin{aligned}
|\rho_{m_1}[F^\nu-F^{\rm b,\nu}]({T_1})|&\ge c\delta_\nu B_*/\sqrt N,\\
\frac{\langle m_1,m_1({T_1}-t_0)\rangle^{{\rm s}_0-3}
|\rho_{m_1}[F^\nu-F^{\rm b,\nu}]({T_1})|}{\epsilon_\nu}
&\ge c\zeta^{-3}B_*/\sqrt N.
\end{aligned}
\end{equation}
The second lower bound tends to infinity if $Kc_{\rm e}>1$.
Constants and the upper bound on $\nu$ may depend on the fixed orders,
$p$ and $K$.
\end{prop}
\par\endgroup

\begingroup
The proof is in Section~\ref{sec:nonlinearConstruction}; the initial
change from $\theta$ to $\theta_\delta$ is included through
\eqref{eq:initialTemperatureMixedBound}--\eqref{eq:initialTemperatureSpatialBound}.
\par\endgroup
We can now compare the full solution with the descending density.
The background contribution is small by damping:
\[
|\rho_{m_1}[F^{\rm b,\nu}-\mu ](T_1)|
\lesssim m_1a e^{o(N)}e^{-c(T_1-t_0)^\chi}
=o(\delta_\nu B_*/\sqrt N).
\]
Together with the linear and nonlinear comparisons, this gives
\begin{equation}\label{eq:roadmapTerminalTarget}
\begin{aligned}
\rho_{m_1}[F^\nu-\mu ](T_1)
={}&\delta_\nu\rho_{m_1}^\downarrow(T_1)
+\underbrace{\delta_\nu\big(\rho_{m_1}[H^\nu](T_1)
-\rho_{m_1}^\downarrow(T_1)\big)}_{\text{linear comparison error}}\\
&+\underbrace{\rho_{m_1}[W](T_1)}_{\text{nonlinear remainder}}
+\underbrace{\rho_{m_1}[F^{\rm b,\nu}-\mu ](T_1)}_{\text{background density}}\\
={}&\delta_\nu\rho_{m_1}^\downarrow(T_1)
+o(\delta_\nu B_*/\sqrt N).
\end{aligned}
\end{equation}
Since {$|\rho_{m_1}^\downarrow(T_1)|\ge cB_*/\sqrt N$}, it follows that
\[
|\rho_{m_1}[F^\nu-\mu ](T_1)|\ge(c-o(1))\delta_\nu B_*/\sqrt N.
\]

Finally, return to physical time $t=s-t_0$. The profile used in the
main theorem satisfies
\[
\sqrt \mu (v)\,\widetilde g^\nu(t,z,v)
=\widetilde F^\nu(t,z+tv_1,v)-\mu (v),\qquad
\rho_k[\widetilde F^\nu-\mu ](t)
=\widehat{\sqrt \mu \,\widetilde g^\nu}_k(t,kt e_1).
\]
Lemma~\ref{lem:nonlinearSobolevTrace} converts this density trace
into a weighted profile Sobolev lower bound.
\par
At $t=t_\nu=T_1-t_0$ and $k=m_1$, we have
$m_1t_\nu=\zeta-m_1t_0\ge\zeta/2$. Consequently,
\[
\|\sqrt \mu \,\widetilde g^\nu(t_\nu)\|_{\mathscr H^{{\rm s}_0-4}_\ell}
\gtrsim\zeta^{{\rm s}_0-4}\delta_\nu B_*/\sqrt N \sim\epsilon_\nu\zeta^{-4}B_*/\sqrt N
=\nu^{5/3-Kc_{\rm e}+o(1)}.
\]
Choosing $Kc_{\rm e}>2$ gives the lower bound $c\nu^{-1/3}$ in
Theorem~\ref{thm:logNonlinearInflation}. 
The proof of Theorem~\ref{thm:logNonlinearInflation} in
Section~\ref{sec:nonlinearNormInflation} verifies the normalized-profile
and density conclusions with the same choice of $K$.
\par

\begingroup
The physical time $t_\nu$ satisfies
\[
t_\nu\sim\nu^{-1/3}[\log(1/\nu)]^{{\rm s}_0+2-p}
=o(\nu^{-1/3}).
\]
\par\endgroup
All choices are made before the limit $\nu\to0$, in the following order. The distinction between an exponential coefficient and a smallness threshold is essential.
\begin{center}
\small
\begin{tabular}{@{}P{.12\linewidth}P{.39\linewidth}P{.41\linewidth}@{}}
\toprule
Order & Choice & Dependence and purpose\\
\midrule
1 & Mode parameters, $\chi\in(1/2,2/3)$, coarse gain constant & Fix $c_{\rm e},C_{\rm B},C_*$ independently of the Sobolev and expansion orders.\\
2 & $Kc_{\rm e}>2$ & Gives the quantitative terminal lower bounds in Theorem~\ref{thm:logNonlinearInflation}.\\
3 & ${\rm s}_0\ge {\rm s}_{\log}$ with $\Theta>0$ & Absorbs the fixed exponential losses in the nonlinear comparison.\\
4 & $p>{\rm s}_0+7/2$ and the finite field order $J$ & Gives small terminal collision errors and a sufficiently small field residual.\\
5 & Remaining fixed derivative, moment and radius choices & Proposition~\ref{lem:fixedCoefficientOrders} supplies the background bounds at these orders.\\
6 & Sufficiently small $\nu$ & This threshold may depend on every preceding fixed choice.\\
\bottomrule
\end{tabular}
\end{center}
The coefficient $C_*$ and the order-dependent smallness thresholds are
specified in \eqref{eq:coarseConstantDependence}. For qualitative profile
inflation, $Kc_{\rm e}>5/3$ suffices, with the same initial-size and
observation-time exponents.
\par

\par

The density decompositions and the corresponding distribution equations are
summarized in Appendix~\ref{app:roadMap},
Figures~\ref{fig:roadMapDistributions}--\ref{fig:roadMapDensityComparison}.
\par

\begingroup
\addtocontents{toc}{\protect\setcounter{tocdepth}{1}}
\section{The prepared background}\label{sec:referenceEstimates}\label{part:linear}
The data and their norms are specified in Lemma~\ref{lem:data}; their
proof and the scale calculation of Lemma~\ref{lem:logarithmicScales}
are given in Appendix~\ref{app:initialNormalization}. We first construct
the global exact background. We then estimate the difference from the
explicit free-transport approximation, the electric field, and the
Fourier coefficients of the background error at a reduced Gevrey radius.
These are the coefficient bounds used in the linear comparison of
Section~\ref{sec:completeLinear}.

\subsection{The exact background and uniform estimates}\label{sec:background}\label{sec:collisionalbackground}\label{sec:collisionalLinear}
We apply the global background theorem to the data in Lemma~\ref{lem:data}.
The Gevrey norm is defined in \eqref{eq:bggevreyinitialnorm}.
The endpoint $T_1$ is used for the comparison estimates below, not
as a restriction on the existence time of the background.
\par\endgroup
Denote the exact VPL background in physical time by
\[
{
\widetilde F^{\rm b,\nu}(t,x_1,v)=F^{\rm b,\nu}(t_0+t,x_1,v),\qquad \widetilde f^{\rm b,\nu}=\frac{\widetilde F^{\rm b,\nu}-\mu }{\sqrt \mu },\qquad E^\nu=E[\widetilde F^{\rm b,\nu}].
}
\]
The equation for $\widetilde F^{\rm b,\nu}$ is
\begin{equation}\label{eq:exactVPLbackground}
\begin{aligned}
\partial_t\widetilde F^{\rm b,\nu}+v_1\partial_{x_1}\widetilde F^{\rm b,\nu}
+E^\nu\partial_{v_1}\widetilde F^{\rm b,\nu}
&=\nu Q(\widetilde F^{\rm b,\nu},\widetilde F^{\rm b,\nu}),\\
\widetilde F^{\rm b,\nu}(0)&=F^{\rm b}_{\rm in}=\mu +\sqrt \mu  f^{\rm b}_{\rm in}.
\end{aligned}
\end{equation}
The Coulomb operator $Q$ has the normalization in \eqref{eq:landauoperator}. Define $Y_t=\nabla_v+t e_1\partial_{x_1}$ and
$g^{\rm b,\nu}(t,z,v)=\widetilde f^{\rm b,\nu}(t,z+tv_1,v)$.

The required normalization and fixed spatial lattice are checked in
Lemma~\ref{lem:fixedPeriodNormalization}. We use \cite[Theorems 1.1 and 2.1]{ref6}
for this collisional background. No fixed-radius Gevrey smallness is imposed
on the high-frequency perturbation.

\begin{prop}\label{prop:bguniform}\label{lem:fixedCoefficientOrders}
Use \eqref{eq:comparisonLogParameters} with ${\rm s}_0\ge11$, $K>0$ and $p>{\rm s}_0+7/2$, and any fixed $A_0\ge1$.
Fix derivative and moment orders ${\rm s},L$. Choose $L_N=C_{{\rm s},L}N^{P_{{\rm s},L}}e^{C_{{\rm s},L}N^\chi}$ with constants independent of $\nu$. For all sufficiently small $\nu$, \eqref{eq:exactVPLbackground} has a global smooth nonnegative solution with the same conservation quantities as {$\mu $}. It depends only on $x_1$. On $0\le t\le {T_1}-t_0$, the exact collisional background satisfies
\begin{align}
\sum_{j+|\alpha|\le {\rm s}}
\|\langle v\rangle^L\partial_{x_1}^jY_t^\alpha \widetilde f^{\rm b,\nu}(t)\|_2
&\le C\varepsilon_N,\qquad \varepsilon_N=aL_N,
\label{eq:bguniformsobolev}\\
\|E^\nu(t)\|_{H^{\rm s}_x}
&\le C_{\rm s}\varepsilon_N e^{-c_{\rm s}t^\chi}.
\label{eq:bguniformfield}
\end{align}
For some fixed $\lambda_*>0$ one also has
\begin{equation}\label{eq:bguniformgevrey}
\|e^{\lambda_*\langle\nabla_z,\nabla_v\rangle^\chi}
g^{\rm b,\nu}(t)\|_2\le C\varepsilon_N.
\end{equation}
All constants are independent of $\nu$. The choice of initial moment order may depend on the prescribed $L,{\rm s}$.

The smallness threshold for $\nu$ may depend on the fixed orders ${\rm s},L$.
Changing these orders changes only $L_N$ and that threshold, not the
construction scales in \eqref{eq:comparisonLogParameters}.
\end{prop}

\begin{proof}
Fix the desired orders ${\rm s},L$ and choose initial indices ${\rm s}_{\rm in},L_0$ in \cite[Theorem 2.1]{ref6} larger than these orders, with velocity moment index at least $4(L+{\rm s}+8)$. Lemma~\ref{lem:bginitialsize} gives
\[
\|f^{\rm b}_{\rm in}\|_{\mathscr G^{\lambda,\chi}_{{\rm s}_{\rm in},L_0}}
\le aL_N=\nu^{1/3-o(1)}\to0,
\qquad {T_1}=o(\nu^{-1/3})<\nu^{-1/2}.
\]
Lemma~\ref{lem:fixedPeriodNormalization} transfers the equation and norms to the fixed-period convention of that theorem.

The reference norm controls the homogeneous part and the transported
nonzero modes. In the present symmetry class, the vector fields
$Y_t=\nabla_v+t e_1\partial_{x_1}$ and the already defined profile
$g^{\rm b,\nu}$ satisfy
\[
{
\partial_z^j\partial_v^\alpha g^{\rm b,\nu}(t,z,v) =(\partial_{x_1}^jY_t^\alpha\widetilde f^{\rm b,\nu}) (t,z+tv_1,v),\qquad dx_1\,dv=dz\,dv.
}
\]

In the fixed normalization of Lemma~\ref{lem:fixedPeriodNormalization},
write {$f'=(F'-\mu')/\sqrt{\mu'}$} and $t'=t/c_t$. {Here $\mu'(v')=e^{-|v'|^2}$ is the Maxwellian in the rescaled variables.}
Equation~\eqref{eq:fixedLatticeReferenceOutputs}, at the chosen initial
orders, gives
\[
\begin{aligned}
&\|f'_0(t')\|_{\mathcal G_{c_x;\chi,L_0}^{\lambda_\infty,{\rm s}_{\rm in}}}
+\langle(\nu')^{1/3}t'\rangle^{L_0/6}
\|f'_{\ne}(t')\|_{\mathcal G_{c_x;\chi,L_0/4}^{\lambda_\infty,{\rm s}_{\rm in}}}
\le C\varepsilon_N.
\end{aligned}
\]
The multiplier on the nonzero-mode norm is at least one, and the
associated decay factor is at most one. Dropping the multiplier gives
the uniform bound on the sum of the two norms used below. The fixed
normalization changes their derivative and moment factors only by fixed
constants.
For each fixed $j+|\alpha|\le {\rm s}$, the corresponding
finite-derivative term in that norm has a fixed positive coefficient.
Since $L_0/4\ge L$,
\[
\|\langle v\rangle^L\partial_{x_1}^jY_t^\alpha \widetilde f^{\rm b,\nu}(t)\|_2 =\|\langle v\rangle^L\partial_z^j\partial_v^\alpha g^{\rm b,\nu}(t)\|_2 \le C_{{\rm s},L}\varepsilon_N.
\]
\par
Summing the finitely many terms proves \eqref{eq:bguniformsobolev}.

For the Fourier radius, use the factorial-series comparison
\eqref{eq:gevreyRadiusComparison} proved in Lemma~\ref{lem:bginitialsize}.
The four pure derivatives in $(z,v_1,v_2,v_3)$ occur in the reference
norm. Since the largest component of $(k,\eta)$ is at least
$|(k,\eta)|/2$, for sufficiently small fixed $\lambda_*>0$ and
$\lambda_{\rm der}>0$,
\[
{
\|e^{\lambda_*\langle\nabla_z,\nabla_v\rangle^\chi}g^{\rm b,\nu}(t)\|_2^2 \le C\sum_{n\ge0}\frac{\lambda_{\rm der}^{2n}}{(n!)^{2/\chi}} \left(\|\partial_z^ng^{\rm b,\nu}(t)\|_2^2 +\sum_{i=1}^3\|\partial_{v_i}^ng^{\rm b,\nu}(t)\|_2^2\right) \le C\varepsilon_N^2.
}
\]
The radii are reduced only by fixed factors, independently of $\nu$.
This proves \eqref{eq:bguniformgevrey}.

The reference damping estimate and $E_0=0$ imply, for $t\ge1$,
\[
\|e^{c|t\partial_{x_1}|^\chi}E^\nu(t)\|_2\le C\varepsilon_N,
\qquad
\sup_{k\ne0}\langle k\rangle^{\rm s} e^{-c|kt|^\chi}
\le C_{\rm s}e^{-c_{\rm s}t^\chi}.
\]
Thus \eqref{eq:bguniformfield} holds for $t\ge1$. For $0\le t\le1$,
\[
\begin{aligned}
|\rho[\widetilde F^{\rm b,\nu}-\mu ]_k(t)|
&=\left|\int\sqrt \mu  \widetilde f_k^{\rm b,\nu}(t,v)\,dv\right|
\le\|\sqrt \mu \|_2\|\widetilde f_k^{\rm b,\nu}(t)\|_2,\\
\|E^\nu(t)\|_{H^{\rm s}_x}
&\le C\|\rho[\sqrt \mu  \widetilde f^{\rm b,\nu}](t)\|_{H^{{\rm s}-1}_x} \le C\|\widetilde f^{\rm b,\nu}(t)\|_{H^{{\rm s}-1}_xL^2_v}
\le C_{\rm s}\varepsilon_N.
\end{aligned}
\]
Global existence and conservation follow from \cite[Theorem 1.1]{ref6}; classical positivity and uniqueness preserve nonnegativity and independence of $x_2,x_3$.

For any further fixed orders ${\rm s},L$, Lemma~\ref{lem:bginitialsize}
supplies a factor $L_{N,{\rm s},L}=C_{{\rm s},L}N^{P_{{\rm s},L}}e^{C_{{\rm s},L}N^\chi}$;
{then the corresponding amplitude bound is $aL_{N,{\rm s},L}$.} Since $\log L_{N,{\rm s},L}=o(N)$,
these fixed orders do not change the powers of $\nu$. For $1/2<{\sigma_{\rm int}}<1$, the required small quantities are
\begin{equation}\label{eq:logAdditionalCoefficientSmallness}
\begin{aligned}
\nu {T_1}^3&\mathrel{{\sim}} [\log(1/\nu)]^{3{\rm s}_0+6-3p}=o(1),&
aL_{N,{\rm s},L}&=\nu^{1/3-o(1)},\\
\nu aL_{N,{\rm s},L}{T_1}^3&=\nu^{1/3-o(1)},&
a^2L_{N,{\rm s},L}\nu^{2{\sigma_{\rm int}}-1}&=\nu^{2{\sigma_{\rm int}}-1/3-o(1)}.
\end{aligned}
\end{equation}
Here $p>{\rm s}_0+7/2$ and ${\sigma_{\rm int}}>1/2$, so every exponent has a strict smallness margin. This proves the assertion, including the coefficient smallness needed at any fixed energy order. We henceforth write $L_N$ and $\varepsilon_N$ also for these enlarged factors; each enlargement involves only finitely many fixed orders.
\end{proof}

\subsection{Direct comparison and the background field}\label{sec:backgroundresponse}\label{app:directBackground}

\begingroup
The explicit defects provide
small source integrals, but the unknown difference also has a
self-consistent electric term. {We close the remainder energy estimate using}
field-energy cancellation and collision coercivity; the exact
background damping controls the remaining electric terms over the long
comparison interval.
We use physical time $t=s-t_0$, with $x=x_1$ and $t_\nu=T_1-t_0$.
\par
\par\endgroup
We use the normalized remainder defined in
\eqref{eq:backgroundRemainder}.
\par
In physical-coordinate identities below, background distributions and
the normalized remainder without a displayed time argument are evaluated at $s=t_0+t$;
$\partial_t$ differentiates these compositions.
The electric fields of the exact background, explicit background and
physical difference are $E[F^{\rm b,\nu}]$,
$E[F^{\rm b}_{\rm app}]$ and {$E[{F^{\rm b}_{\rm err}}]$}, respectively.
We first estimate the exact current, derive the difference equation,
{and prove the low-norm energy estimate.} Interpolation followed by a
density trace gives the smaller background field; that small field is
{not assumed in the proof of that estimate.}

Recall the normalized background profile $g^{\rm b,\nu}$ introduced after
\eqref{eq:exactVPLbackground}. The exact and explicit profiles satisfy
\begin{equation}\label{eq:physicalCollisionProfiles}
{
\begin{aligned}
g^{\rm b,\nu}(t,z,v)
&=f^{\rm b,\nu}(t_0+t,z+tv_1,v),\\
g^{\rm b}_{\rm app}(t,z,v)
&=f^{\rm b}_{\rm app}(t_0+t,z+tv_1,v)
=a\sqrt \mu \cos(z-t_0v_1),\\
{g^{\rm b}_{\rm err}}(t,z,v)
&=g^{\rm b,\nu}(t,z,v)-g^{\rm b}_{\rm app}(t,z,v)
={f^{\rm b}_{\rm err}}(t_0+t,z+tv_1,v).
\end{aligned}}
\end{equation}
We estimate this error profile {as a whole; the required smallness is not obtained}
by estimating the two profiles separately.
Write $\rho^\nu=\int\widetilde F^{\rm b,\nu}\,dv$ and
$j^\nu=\int v_1\widetilde F^{\rm b,\nu}\,dv$.
We use the Sobolev and Gevrey bounds
\eqref{eq:bguniformsobolev} and \eqref{eq:bguniformgevrey}, with
$\varepsilon_N=aL_N$ as in Proposition~\ref{prop:bguniform}.
\par

Lemma~\ref{lem:logarithmicScales} also gives $t_\nu\sim T_1$,
$\nu T_1^3\to0$ and $\log L_N=o(N)$. In addition,
\begin{equation}\label{eq:backgroundTimeLowerBound}
T_1\sim N^2/a\ge cN^{{\rm s}_0+4},
\end{equation}
since $a\le N^{-{\rm s}_0-2}$.
\par
\subsubsection{The current, explicit defect and difference equation}

We first obtain the time-integrable potential bound. Then the exact
subtraction separates the prescribed collision and electric defects from
the terms acting on the unknown difference.

\begin{lem}\label{lem:direct-current}
Let $\phi^\nu$ be the mean-zero potential, $E^\nu=-\partial_x\phi^\nu$.
For every fixed $j\ge0$, there are $C_j,c_j>0$ such that
\begin{equation}\label{eq:direct-coarse-field}
\|E^\nu(t)\|_{C^j_x}+\|\phi^\nu(t)\|_{C^j_x}
+\|\partial_t\phi^\nu(t)\|_{C^j_x}
\le C_j\varepsilon_N e^{-c_jt^\chi}.
\end{equation}
Also
\begin{equation}\label{eq:direct-Er-coarse}
\|E[{F^{\rm b}_{\rm err}}](t_0+t)\|_2\le C\varepsilon_N e^{-ct^\chi}.
\end{equation}
Writing $\phi^{\rm b,\nu}(s)=\phi^\nu(s-t_0)$, one has
\begin{equation}\label{eq:backgroundPotentialIntegral}
\sup_I\|\phi^{\rm b,\nu}\|_{C^{j+2}_x}
+\int_I\|\partial_s\phi^{\rm b,\nu}\|_{C^j_x}\,ds
\le C_j\varepsilon_N.
\end{equation}
\end{lem}
\begin{proof}
For {$q(v)=\sqrt \mu (v)$} or {$q(v)=v_1\sqrt \mu (v)$},
$|\widehat q(\xi)|\le C(1+|\xi|)e^{-|\xi|^2}$.
For $k\ne0$ and $0<\lambda\le\lambda_*$, Fourier convolution gives
\begin{align*}
I_{q,k}(t)
&=\int q(v)e^{-iktv_1}g^{\rm b,\nu}_k(t,v)\,d v\\
&=(2\pi)^{-3}\int\widehat q(kt e_1-\eta)
       \widehat g^{\rm b,\nu}_k(t,\eta)\,d\eta,\\
\lambda\dbla{k,kt e_1}^\chi
&\le\lambda_*\dbla{k,\eta}^\chi+
       \lambda|kt e_1-\eta|^\chi,\\
|I_{q,k}(t)|
&\le C e^{-\lambda\dbla{k,kt e_1}^\chi}
\|e^{\lambda_*\dbla{k,\eta}^\chi}\widehat g^{\rm b,\nu}_k\|_2
\|e^{\lambda|\xi|^\chi}\widehat q(\xi)\|_2\\
&\le C\varepsilon_Ne^{-c(|k|^\chi+|kt|^\chi)}.
\end{align*}
These are the density and current coefficients. Local mass conservation
and Poisson give
\[
\partial_t\rho^\nu+\partial_xj^\nu=0,\quad
\phi^\nu_k=\rho_k^\nu/k^2,\quad
\partial_t\phi^\nu_k=-ij^\nu_k/k .
\]
For $k\ne0$, $|kt|^\chi\ge t^\chi$. Summing the exponentially decaying
spatial modes proves \eqref{eq:direct-coarse-field}.
Direct integration of \eqref{eq:bgapp} gives
\begin{equation}\label{eq:direct-EB}
\rho[F^{\rm b}_{\rm app}]-1=a e^{-s^2/2}\cos x,\qquad
E[F^{\rm b}_{\rm app}]=a e^{-s^2/2}\sin x.
\end{equation}
Since $s=t_0+t\ge t$, {the $L^2$ norm of this source is bounded by}
$Ca e^{-ct^\chi}$. This and \eqref{eq:direct-coarse-field} imply
\eqref{eq:direct-Er-coarse}. Applying the potential estimate at order $j+2$
and integrating the time-derivative estimate at order $j$ gives
\[
\sup_I\|\phi^{\rm b,\nu}\|_{C^{j+2}_x}
+\int_I\|\partial_s\phi^{\rm b,\nu}\|_{C^j_x}\,ds
\le C_j\varepsilon_N\left(1+\int_0^\infty e^{-c_jt^\chi}\,dt\right)
\le C_j\varepsilon_N,
\]
which proves \eqref{eq:backgroundPotentialIntegral}.
\end{proof}

To estimate the terms in \eqref{eq:direct-r-eq} and the later
Gaussian-weighted collision forms, set
\par
\begin{equation}\label{eq:gaussianConjugatedDerivatives}
D_{{\varpi},i}^{\pm}
=\partial_{v_i}\pm\frac{v_i}{2{\varpi}},\qquad
D_{{\varpi}}^{\pm}
=\nabla_v\pm\frac{v}{2{\varpi}},\qquad
1\le{\varpi}\le\frac32.
\end{equation}

\begingroup
These are the derivatives obtained by conjugating with a Gaussian:
\[
\mu_{\varpi}^{-1/2}\partial_{v_i}(\mu_{\varpi}^{1/2}f)
=D^-_{\varpi,i}f,\qquad
\mu_{\varpi}^{1/2}\partial_{v_i}(\mu_{\varpi}^{-1/2}f)
=D^+_{\varpi,i}f.
\]
Thus $D^-_{1,i}=\partial_{v_i}-v_i/2$ differentiates a distribution
normalized by $\sqrt\mu$, and integration by parts places $D^+_{1,i}$
on the other factor. The subscript $1$ denotes the variance of
$\mu=\mu_1$.
\par\endgroup

\begin{lem}\label{lem:direct-equation}
Use the backgrounds and remainder defined above, with $s=t_0+t$.
Then the normalized error {${f^{\rm b}_{\rm err}}$} satisfies
\begin{equation}\label{eq:direct-r-eq}
{
\begin{aligned}
&\left(\partial_t+v_1\partial_x+E^\nu\left(\partial_1-\frac{v_1}{2}\right)\right){f^{\rm b}_{\rm err}}
+\mu^{-1/2}E[{F^{\rm b}_{\rm err}}]\partial_1F^{\rm b}_{\rm app}\\
&\quad=\nu \mu^{-1/2}\bigl(Q(\widetilde F^{\rm b,\nu},{F^{\rm b}_{\rm err}})
+Q({F^{\rm b}_{\rm err}},F^{\rm b}_{\rm app})\bigr)\\
&\qquad\quad+\mu^{-1/2}\bigl(\nu Q(F^{\rm b}_{\rm app},F^{\rm b}_{\rm app})
-E[F^{\rm b}_{\rm app}]\partial_1F^{\rm b}_{\rm app}\bigr).
\end{aligned}}
\end{equation}
{The two known sources are} the collision and electric defects of the
explicit approximation. For every fixed $\ell\ge0$, they obey the following
bounds, with $D_{1,i}^-$ from \eqref{eq:gaussianConjugatedDerivatives}:
\begin{equation}\label{eq:direct-source}
\begin{gathered}
\sup_x\left\|\dbla{v}^\ell D_{1,1}^-\frac{F^{\rm b}_{\rm app}(s,x,\cdot)-\mu }{\sqrt \mu }\right\|_{L^2_v}
\le C_\ell a\dbla{s},\\
\|\dbla{v}^\ell \mu^{-1/2}Q(F^{\rm b}_{\rm app},F^{\rm b}_{\rm app})(s)\|_2
\le C_\ell a\dbla{s}^2,\\
{
\left\|\dbla{v}^\ell \mu^{-1/2} \big(\nu Q(F^{\rm b}_{\rm app},F^{\rm b}_{\rm app}) -E[F^{\rm b}_{\rm app}]\partial_1F^{\rm b}_{\rm app}\big)\right\|_2 \le C_\ell a\{\nu\dbla{s}^2+(1+a\dbla{s})e^{-s^2/2}\}.
}
\end{gathered}
\end{equation}

In particular, for $t_0\ge1$,
\[
\int_{t_0}^{\infty}
\left\|\mu^{-1/2}E[F^{\rm b}_{\rm app}](s)
\partial_1F^{\rm b}_{\rm app}(s)\right\|_2\,ds
\le Ca(t_0^{-1}+a)e^{-t_0^2/2}.
\]
Both $\widetilde F^{\rm b,\nu}$ and $F^{\rm b}_{\rm app}$ satisfy local mass conservation; hence
\begin{equation}\label{eq:direct-field-energy-id}
\begin{aligned}
\partial_tE[{F^{\rm b}_{\rm err}}]
&=-j[{F^{\rm b}_{\rm err}}]+\frac1{2\pi}\int_0^{2\pi}j[{F^{\rm b}_{\rm err}}](t_0+t,y_1)\,dy_1,\\
\frac12\frac{d}{dt}\|E[{F^{\rm b}_{\rm err}}]\|_2^2
&=-\int E[{F^{\rm b}_{\rm err}}]j[{F^{\rm b}_{\rm err}}]\,dx.
\end{aligned}
\end{equation}
\end{lem}
\begin{proof}
As $(\partial_t+v_1\partial_x)F^{\rm b}_{\rm app}=0$,
\[
(\partial_t+v_1\partial_x+E^\nu\partial_1)\bigl({F^{\rm b}_{\rm err}}\bigr)
=\nu Q(\widetilde F^{\rm b,\nu},\widetilde F^{\rm b,\nu})-E^\nu\partial_1F^{\rm b}_{\rm app}.
\]
The algebraic identity
\[
Q(\widetilde F^{\rm b,\nu},\widetilde F^{\rm b,\nu})
=Q(\widetilde F^{\rm b,\nu},{F^{\rm b}_{\rm err}})+Q({F^{\rm b}_{\rm err}},F^{\rm b}_{\rm app})+Q(F^{\rm b}_{\rm app},F^{\rm b}_{\rm app})
\]
and {$E^\nu=E[{F^{\rm b}_{\rm err}}]+E[F^{\rm b}_{\rm app}]$} give \eqref{eq:direct-r-eq} after division by {$\sqrt \mu $}.

Direct differentiation gives
\begin{align*}
D_{1,1}^-\big(\mu^{-1/2}(F^{\rm b}_{\rm app}-\mu )\big)&=a\sqrt \mu (-v_1\cos(x-sv_1)+s\sin(x-sv_1)),\\
\partial_{ij}(a\mu \cos(x-sv_1))
&=a\mu (v_iv_j-\delta_{ij}-s^2\delta_{i1}\delta_{j1})\cos(x-sv_1)\\
&\quad-as\mu (v_i\delta_{j1}+v_j\delta_{i1})\sin(x-sv_1).
\end{align*}
The first estimate follows immediately. Here {$F^{\rm b}_{\rm app}-\mu =a\mu \cos(x-sv_1)$}.
The cancellation identity of Lemma~\ref{lem:LandauCancellation} is
\[
Q({F_1},{F_2})=(\Phi_{ij}*{F_1})\partial_{ij}{F_2}+8\pi {F_1} {F_2}.
\]
Since {$Q(\mu ,\mu )=0$},
\begin{align*}
\mu^{-1/2}Q(F^{\rm b}_{\rm app},F^{\rm b}_{\rm app})={}&(\Phi_{ij}*\mu )\mu^{-1/2}\partial_{ij}(F^{\rm b}_{\rm app}-\mu )\\
&+(\Phi_{ij}*(F^{\rm b}_{\rm app}-\mu ))(v_iv_j-\delta_{ij})\sqrt \mu \\
&+(\Phi_{ij}*(F^{\rm b}_{\rm app}-\mu ))\mu^{-1/2}\partial_{ij}(F^{\rm b}_{\rm app}-\mu )\\
&+16\pi\sqrt \mu  (F^{\rm b}_{\rm app}-\mu )+8\pi \mu^{-1/2}(F^{\rm b}_{\rm app}-\mu )^2.
\end{align*}
The convolution kernel {is locally integrable, and convolution with}
a Gaussian is bounded; thus
{$|\Phi*\mu |\le C$}, {$|\Phi*(F^{\rm b}_{\rm app}-\mu )|\le Ca$}.
The displayed derivative formula and Gaussian integration give
{$\|\dbla{v}^\ell \mu^{-1/2}Q(F^{\rm b}_{\rm app},F^{\rm b}_{\rm app})\|_2\le C_\ell(a+a^2)\dbla{s}^2$}.
Furthermore
{$\mu^{-1/2}\partial_1F^{\rm b}_{\rm app}=-v_1\sqrt \mu +D_{1,1}^-\big(\mu^{-1/2}(F^{\rm b}_{\rm app}-\mu )\big)$}.
Together with \eqref{eq:direct-EB}, this proves \eqref{eq:direct-source}.

For the integrated electric bound, the same derivative formula gives
\[
\begin{aligned}
\|\mu^{-1/2}E[F^{\rm b}_{\rm app}]\partial_1\mu \|_2
&\le Ca e^{-s^2/2},\\
\|\mu^{-1/2}E[F^{\rm b}_{\rm app}]
\partial_1(F^{\rm b}_{\rm app}-\mu )\|_2
&\le Ca^2\langle s\rangle e^{-s^2/2}.
\end{aligned}
\]
Since $s\ge t_0\ge1$, integration uses only
\[
\int_{t_0}^{\infty}e^{-s^2/2}\,ds
\le t_0^{-1}e^{-t_0^2/2},\qquad
\int_{t_0}^{\infty}\langle s\rangle e^{-s^2/2}\,ds
\le\sqrt2\,e^{-t_0^2/2}.
\]
Integrating free transport for $F^{\rm b}_{\rm app}$ and the VPL equation for $\widetilde F^{\rm b,\nu}$
in velocity gives the two continuity equations.
Their difference, together with {$\partial_xE[{F^{\rm b}_{\rm err}}]=\rho[{F^{\rm b}_{\rm err}}]$} and
{$\int E[{F^{\rm b}_{\rm err}}]\,d x=0$}, proves \eqref{eq:direct-field-energy-id}.
\end{proof}
\par

The difference equation now requires a coercive collision form.
{We first record the Maxwellian part of the energy estimate,} which will also be used in the
localized propagation argument of Section~\ref{sec:centeredcollisions}.

Recall the normalized Maxwellian operator $\mathcal L$ from
\eqref{eq:normalizedLandauExplicit}.
We use the collision matrix {$\sigma_{ij}=\Phi_{ij}*\mu $} and the norm
$\Sigma_\ell$ defined in \eqref{eq:backgroundCollisionNorm}.

\begingroup
\begin{lem}\label{lem:weightedMaxwellianForm}
For each fixed $\ell\ge0$ there are $c_\ell,C_\ell>0$ such that every
smooth rapidly decreasing function $y(v)$ satisfies
\begin{equation}\label{eq:allWeightMaxwellianCoercivity}
\operatorname{Re}\langle\mathcal Ly,y\rangle_\ell
\ge c_\ell\|y\|_{\Sigma_\ell}^2-C_\ell\|y\|_{L_\ell^2}^2.
\end{equation}
\end{lem}
\begin{proof}
The weighted lower bound is \cite[Lemma 4.4]{ref4}, following the
classical weighted Landau estimates; see also
\cite[Lemma 5]{Guo2002}. Set the exponential-weight parameter in
\cite[Lemma 4.4]{ref4} to $\vartheta=0$ and choose the small parameter
there so that the coefficient of the collision norm is positive.
The fixed velocity scaling $v=\sqrt2v'$ in
Lemma~\ref{lem:fixedPeriodNormalization} gives
\[
\operatorname{Re}\langle\mathcal Ly,y\rangle_\ell
\ge c_\ell\|y\|_{\Sigma_\ell}^2
-C_\ell\|\mathbf1_{|v|\le C_\ell}y\|_2^2
\ge c_\ell\|y\|_{\Sigma_\ell}^2-C_\ell\|y\|_{L_\ell^2}^2.
\]
The polynomial weights before and after this fixed scaling are
equivalent. For complex $y$, apply the real estimate separately to
$\operatorname{Re}y$ and $\operatorname{Im}y$ and add.
\end{proof}
For the differentiated calculation in Section~\ref{sec:centeredcollisions},
we record the explicit decomposition following from
\eqref{eq:normalizedLandauExplicit}:
\begin{equation}\label{eq:centeredlocaldecomposition}
\begin{aligned}
\mathcal L y&=-\partial_i(\sigma_{iq}\partial_qy)+\mathcal V y-\mathcal K y,\quad
\mathcal V=\frac14 \sigma_{iq}v_iv_q-4\pi \mu ,\\
\mathcal K y&=\sqrt \mu (v)(v_iv_q-\delta_{iq})
\int\Phi_{iq}(v-v_*)\sqrt \mu (v_*)y(v_*)\,dv_*+8\pi \mu  y.
\end{aligned}
\end{equation}
Indeed, $\partial_i\sigma_{iq}=-\sigma_{iq}v_i$ and
$\partial_i(\sigma_{iq}v_q)=8\pi\mu$. These identities give the local
part, while the cancellation formula for $Q(\sqrt\mu y,\mu)$ gives
$\mathcal K$.
\par\endgroup

\begin{lem}
\label{lem:direct-form}
Use {$\sigma=\Phi*\mu $} and the norm $\Sigma_\ell$ from
{\eqref{eq:backgroundCollisionNorm}.}
Define
\begin{equation}\label{eq:direct-A}
\mathcal C_{\rm b,app}(t) h=\mu^{-1/2}
\{Q(\widetilde F^{\rm b,\nu},\sqrt \mu  h)+Q(\sqrt \mu  h,F^{\rm b}_{\rm app})\}.
\end{equation}
For sufficiently small $\nu$ (depending on the fixed $\ell$),
\begin{equation}\label{eq:direct-form}
\operatorname{Re}\langle\mathcal C_{\rm b,app}(t) h,h\rangle_{\ell}
\le-c_\ell\|h\|_{L^2_x\Sigma_\ell}^2
+C_\ell(1+a\dbla{s}^2)\|h\|_{L_x^2L_\ell^2}^2.
\end{equation}
In particular, no physical velocity derivative of {$\widetilde F^{\rm b,\nu}-\mu $} is required
in the first collision ordering.
\end{lem}
\begin{proof}
Decompose {$\widetilde F^{\rm b,\nu}=\mu +{F^{\rm c}}$} and use {$F^{\rm b}_{\rm app}-\mu =a\mu \cos(x-sv_1)$}.
The weighted Maxwellian estimate \eqref{eq:allWeightMaxwellianCoercivity},
with the sign convention of \eqref{eq:normalizedLandauExplicit}, gives
\begin{equation}\label{eq:direct-maxwell-form}
\operatorname{Re}\langle\mu^{-1/2}
[Q(\mu ,\sqrt \mu  h)+Q(\sqrt \mu  h,\mu )],h\rangle_{\ell}
\le-c_\ell^0\|h\|_{L^2_x\Sigma_\ell}^2
+C_\ell\|h\|_{L_x^2L_\ell^2}^2.
\end{equation}

Sobolev embedding in $(z,v)$ in {\eqref{eq:bguniformsobolev}} and the shear give
{$|{F^{\rm c}}(t,x,v)|\le C\varepsilon_Ne^{-|v|^2/4}$}.
For {$\sigma[{F^{\rm c}}]=\Phi*{F^{\rm c}}$} and {$b_{i}[{F^{\rm c}}]=(\partial_j\Phi_{ij})*{F^{\rm c}}$},
Lemma~\ref{lem:gaussianCollisionCoefficients} therefore yields
\begin{equation}\label{eq:direct-matrix}
|\overline X \sigma[{F^{\rm c}}] Z|
\le C\varepsilon_N|\sigma^{1/2}X||\sigma^{1/2}Z|,\qquad
|\sigma^{-1/2}b[{F^{\rm c}}]|\le C\varepsilon_N\dbla{v}^{-1/2}.
\end{equation}
With $D_{1,i}^{\pm}$ from \eqref{eq:gaussianConjugatedDerivatives}, integration by parts and
\eqref{eq:physicalFormFactors} give
\begin{align*}
\langle\mu^{-1/2}Q({F^{\rm c}},\sqrt \mu  h),h\rangle_{\ell}
&=-\int(\sigma[{F^{\rm c}}]D_1^-h-b[{F^{\rm c}}]h)\cdot
\overline{D_1^+(\langle v\rangle^{2\ell}h)}\,dx\,dv,\\
\big|\langle\mu^{-1/2}Q({F^{\rm c}},\sqrt \mu  h),h\rangle_{\ell}\big|
&\le C\varepsilon_N
\big(\|\sigma^{1/2}D_1^-h\|_{L_x^2L_\ell^2}
+\|\langle v\rangle^{\ell-1/2}h\|_2\big)
\|\langle v\rangle^{-\ell}\sigma^{1/2}D_1^+(\langle v\rangle^{2\ell}h)\|_2\\
&\le C_\ell\varepsilon_N\|h\|_{L^2_x\Sigma_\ell}^2.
\end{align*}

For the second ordering, cancellation and the explicit derivative
formula in Lemma~\ref{lem:direct-equation} give
\[
\begin{aligned}
\mu^{-1/2}Q(\sqrt \mu  h,F^{\rm b}_{\rm app}-\mu )
&=(\Phi_{ij}*(\sqrt \mu  h))\mu^{-1/2}\partial_{ij}(F^{\rm b}_{\rm app}-\mu )+8\pi (F^{\rm b}_{\rm app}-\mu )h,\\
|\mu^{-1/2}\partial_{ij}(F^{\rm b}_{\rm app}-\mu )|
&\le Ca\dbla{s}^2\dbla{v}^2\sqrt \mu .
\end{aligned}
\]
{After conjugation by $\langle v\rangle^\ell$,} the integral kernel satisfies
\[
\begin{aligned}
|K_\ell(v,v_*)|
&\le Ca\dbla{s}^2\dbla{v}^{\ell+2}\sqrt \mu (v)
|v-v_*|^{-1}\sqrt \mu (v_*)\dbla{v_*}^{-\ell}\\
&\le C_\ell a\dbla{s}^2\dbla{v}^{-4}
|v-v_*|^{-1}e^{-|v_*|^2/4}.
\end{aligned}
\]
Thus \eqref{eq:weightedCoulombSchur}, with $d=4$, and the multiplication
term {$8\pi (F^{\rm b}_{\rm app}-\mu )h$} yield
\[
\|\mu^{-1/2}Q(\sqrt \mu  h,F^{\rm b}_{\rm app}-\mu )\|_{L_x^2L_\ell^2}
\le C_\ell a\dbla{s}^2\|h\|_{L_x^2L_\ell^2}.
\]
Adding the three pairings gives
\[
\operatorname{Re}\langle\mathcal C_{\rm b,app}(t) h,h\rangle_{\ell}
\le-(c_\ell^0-C_\ell\varepsilon_N)\|h\|_{L^2_x\Sigma_\ell}^2
+C_\ell(1+a\dbla{s}^2)\|h\|_{L_x^2L_\ell^2}^2.
\]
Taking $C_\ell\varepsilon_N\le c_\ell^0/2$ proves \eqref{eq:direct-form}.
\end{proof}

\subsubsection{The low norm and polynomial velocity moments}

\begingroup
We now prove the two low-norm bounds needed downstream. Field energy cancels
the Maxwellian term in the density equation. A split at a polylogarithmic time controls
the remaining electric terms without an exponential in $a{T_1}^2$.
The weighted estimate uses the explicit velocity dependence of $F^{\rm b}_{\rm app}$.
\par\endgroup
\begin{proof}[Proof of Proposition~\ref{thm:bgcomparison}]
Proposition~\ref{prop:bguniform} supplies the global smooth nonnegative
background. Lemma~\ref{lem:data} {ensures Maxwellian mass, momentum and}
total energy, which are conserved by the VPL evolution. We prove
\eqref{eq:polynomialBackgroundLow}--\eqref{eq:polynomialBackgroundMoments}
in elapsed time $t=s-t_0$, on $0\le t\le t_\nu={T_1}-t_0$.

The initial discrepancy is determined explicitly:
\begin{equation}\label{eq:directInitialMismatch}
\begin{aligned}
{f^{\rm b}_{\rm err}}(t_0)&=\frac{\mu_\theta-\mu }{\sqrt \mu }
\big(1+a\cos(x_1-t_0v_1)\big),\\
\|\langle v\rangle^\ell {f^{\rm b}_{\rm err}}(t_0)\|_2+\|E[{F^{\rm b}_{\rm err}}](t_0)\|_2
&\le C_\ell|\theta-1|\le C_\ell a^2e^{-t_0^2}.
\end{aligned}
\end{equation}
Lemma~\ref{lem:data} and the mean-value formula give the second line.

Set {$y=e^{\phi^\nu/2}{f^{\rm b}_{\rm err}}$}. Lemma~\ref{lem:potentialTransport},
with physical time $t$ and Gaussian width one, gives
\[
e^{\phi^\nu/2}\left(\partial_t+v_1\partial_x+E^\nu\left(\partial_1-\frac{v_1}{2}\right)\right){f^{\rm b}_{\rm err}}
=(\partial_t+v_1\partial_x+E^\nu\partial_1)y
-\tfrac12(\partial_t\phi^\nu)y.
\]
Since $e^{\phi^\nu/2}$ is independent of velocity, it commutes with both terms of
$\mathcal C_{\rm b,app}(t)$. Equation \eqref{eq:direct-r-eq} becomes
\begin{equation}\label{eq:direct-y}
\begin{aligned}
&(\partial_t+v_1\partial_x+E^\nu\partial_1)y =\nu\mathcal C_{\rm b,app}(t)y+\tfrac12(\partial_t\phi^\nu)y
+e^{\phi^\nu/2}E[{F^{\rm b}_{\rm err}}]v_1\sqrt \mu\\
&\qquad-e^{\phi^\nu/2}E[{F^{\rm b}_{\rm err}}]D_{1,1}^-\big(\mu^{-1/2}(F^{\rm b}_{\rm app}-\mu )\big)\\
&\qquad+e^{\phi^\nu/2}\mu^{-1/2}\big(\nu Q(F^{\rm b}_{\rm app},F^{\rm b}_{\rm app})
-E[F^{\rm b}_{\rm app}]\partial_1F^{\rm b}_{\rm app}\big).
\end{aligned}
\end{equation}
By \eqref{eq:direct-coarse-field},
$\|\phi^\nu\|_\infty\le C\varepsilon_N=o(1)$, so multiplication by $e^{\phi^\nu/2}$ and
${e^{-\phi^\nu/2}}$ is uniformly bounded in every norm used in this proof.

Define
\[
Y(t)^2=\|y(t)\|_2^2+\|E[{F^{\rm b}_{\rm err}}](t_0+t)\|_2^2.
\]
The transport vector field $(v_1,E^\nu)$ is divergence free in $(x,v_1)$.
Using \eqref{eq:direct-field-energy-id}, the Maxwellian electric term contributes
exactly
\begin{equation}\label{eq:direct-cancellation}
\begin{aligned}
&\int e^{\phi^\nu/2}E[{F^{\rm b}_{\rm err}}]v_1\sqrt \mu  y\,dx\,dv
-\int E[{F^{\rm b}_{\rm err}}]v_1{F^{\rm b}_{\rm err}}\,dx\,dv\\
&\qquad=\int E[{F^{\rm b}_{\rm err}}]v_1(e^{\phi^\nu}-1){F^{\rm b}_{\rm err}}\,dx\,dv.
\end{aligned}
\end{equation}
{The absolute value of this pairing is at most}
$C\|\phi^\nu\|_\infty Y^2$.
The explicit first bound in \eqref{eq:direct-source} gives
\[
\begin{aligned}
\|e^{\phi^\nu/2}E[{F^{\rm b}_{\rm err}}]D_{1,1}^-\big(\mu^{-1/2}(F^{\rm b}_{\rm app}-\mu )\big)\|_2
&\le Ca\dbla{s}\|E[{F^{\rm b}_{\rm err}}]\|_2,\\
\|E[{F^{\rm b}_{\rm err}}](t_0+t)\|_2
&\le C\min\{Y(t),\varepsilon_Ne^{-ct^\chi}\}.
\end{aligned}
\]
Pair \eqref{eq:direct-y} with $y$ and add the field-energy identity
\eqref{eq:direct-field-energy-id}. The cancellation
\eqref{eq:direct-cancellation}, the collision form
\eqref{eq:direct-form} and the explicit source bounds give
\[
\begin{aligned}
\frac12\frac{d}{dt}Y^2+c\nu\|y\|_{L_x^2\Sigma_0}^2
\le{}&C\big[\|\partial_t\phi^\nu\|_\infty+\|\phi^\nu\|_\infty
+\nu(1+a\langle s\rangle^2)\big]Y^2\\
&+Ca\langle s\rangle\|E[F^{\rm b}_{\rm err}]\|_2Y +Ca\big[\nu\langle s\rangle^2
+(1+a\langle s\rangle)e^{-s^2/2}\big]Y.
\end{aligned}
\]
Using \eqref{eq:direct-coarse-field} and the preceding minimum bound,
then dividing by $Y$ when $Y>0$, yields
\par
\begin{equation}\label{eq:direct-scalar}
\begin{split}
Y'\le{}&C\big[\varepsilon_Ne^{-ct^\chi}+\nu(1+a\dbla{s}^2)\big]Y
+Ca\dbla{s}\min\{Y,\varepsilon_Ne^{-ct^\chi}\}\\
&+Ca\nu\dbla{s}^2+Ca(1+a\dbla{s})e^{-s^2/2}.
\end{split}
\end{equation}
The coefficient of $Y$ is time-integrable, with
\begin{equation}\label{eq:direct-b}
\int_0^{t_\nu}\big[\varepsilon_Ne^{-ct^\chi}
+\nu(1+a\dbla{t_0+t}^2)\big]\,dt
\le C(\varepsilon_N+\nu {T_1}+a\nu {T_1}^3)=o(1).
\end{equation}
{The inequality at zeros of $Y$ follows by first using}
{$\sqrt{Y^2+\eta^2}$ and then letting $\eta\downarrow0$.}

Choose a fixed $C_1$ sufficiently large and put
\begin{equation}\label{eq:direct-R}
R_\nu=(C_1\log(1/\nu))^{1/\chi}.
\end{equation}
For all small $\nu$,
\begin{equation}\label{eq:direct-R-properties}
\begin{aligned}
1<R_\nu<t_\nu,\qquad (1+t_0+R_\nu)&\le CN^{1/\chi},&
a(1+t_0+R_\nu)^2&=o(1),\\
\int_{R_\nu}^\infty(1+t_0+t)^j e^{-ct^\chi}\,d t
&\le\nu^2,\qquad j=0,1 .
\end{aligned}
\end{equation}
To verify the tail estimate, substitute $u=t^\chi$. A fixed polynomial
factor is absorbed into $e^{cu/2}$, leaving
$C(1+t_0)^j e^{-cR_\nu^\chi/2}\le\nu^2$
when $C_1$ is sufficiently large. The other assertions follow from
{\eqref{eq:backgroundTimeLowerBound} and \eqref{eq:comparisonLogParameters}}.

Use $Y$ in the minimum for $t\le R_\nu$ and the damping bound thereafter.
By \eqref{eq:directInitialMismatch} and the bounded potential multiplier,
\[
Y(0)\le C\big(\|{f^{\rm b}_{\rm err}}(t_0)\|_2+\|E[{F^{\rm b}_{\rm err}}](t_0)\|_2\big)
\le Ca^2e^{-t_0^2}.
\]
The explicit source satisfies
\[
\int_0^{t_\nu}(1+a\dbla{t_0+t})e^{-(t_0+t)^2/2}\,d t
\le Ce^{-t_0^2/4}.
\]
The full Gronwall inequality from \eqref{eq:direct-scalar} is therefore
\begin{align}
Y(t)
&\le \exp\left(C(\varepsilon_N+\nu {T_1}+a\nu {T_1}^3)+
Ca\int_0^{R_\nu}\dbla{t_0+\tau}\,d\tau\right)\nonumber\\
&\quad\times\bigg[Y(0)+Ca\nu\int_0^t\dbla{t_0+\tau}^2\,d\tau
+Cae^{-t_0^2/4}\nonumber\\
&\hspace{42mm}
+Ca\varepsilon_N\mathbf1_{\{t>R_\nu\}}
\int_{R_\nu}^t\dbla{t_0+\tau}e^{-c\tau^\chi}\,d\tau\bigg]\nonumber\\
&\le C\{ae^{-t_0^2/4}+a\nu t\dbla{t_0+t}^2+
a\varepsilon_N\nu^2\mathbf1_{\{t>R_\nu\}}\}.
\label{eq:direct-low-gronwall}
\end{align}
For $t>R_\nu>1$, $\varepsilon_N\nu^2\le\nu t\dbla{t_0+t}^2$.
This proves \eqref{eq:polynomialBackgroundLow}.

{For the moment estimate set ${Z_\ell=\|y\|_{L_x^2L_\ell^2}}$.}
The only contribution of field transport to the weighted energy is
\[
\operatorname{Re}\langle E^\nu\partial_1y,y\rangle_{\ell}
=-\ell\int E^\nu v_1\langle v\rangle^{2\ell-2}|y|^2\,dx\,dv,
\qquad
\big|\operatorname{Re}\langle E^\nu\partial_1y,y\rangle_\ell\big|
\le C_\ell\|E^\nu\|_\infty Z_\ell^2.
\]
{Take the weighted inner product of \eqref{eq:direct-y} with $y$. Equations}
\eqref{eq:direct-source}, \eqref{eq:direct-form}, and \eqref{eq:direct-coarse-field} yield
\begin{equation}\label{eq:direct-moment-scalar}
\begin{split}
Z_\ell'\le{}&C_\ell\big[\varepsilon_Ne^{-ct^\chi}
+\nu(1+a\dbla{s}^2)\big]Z_\ell
+C_\ell(1+a\dbla{s})\|E[{F^{\rm b}_{\rm err}}](t_0+t)\|_2\\
&+C_\ell a\nu\dbla{s}^2
+C_\ell a(1+a\dbla{s})e^{-s^2/2}.
\end{split}
\end{equation}
The coefficient of $Z_\ell$ has integral $o(1)$ by \eqref{eq:direct-b}.
Only the $L^2_x$ norm of {$E[{F^{\rm b}_{\rm err}}]$} is needed: the weighted $L^2_v$ norm of
{$D_{1,1}^-\big(\mu^{-1/2}(F^{\rm b}_{\rm app}-\mu )\big)$} is uniformly bounded in $x$ by \eqref{eq:direct-source}.

For $0\le t\le R_\nu$, $1+a\dbla{s}\le C$ and \eqref{eq:polynomialBackgroundLow} give
\begin{align}
\int_0^{R_\nu}(1+a\dbla{s})\|E[{F^{\rm b}_{\rm err}}](t_0+t)\|_2\,d t
&\le Ca\int_0^{R_\nu}\{\nu t\dbla{t_0+t}^2+e^{-t_0^2/4}\}\,d t\nonumber\\
&\le Ca\nu (1+t_0+R_\nu)^4+CaR_\nu e^{-t_0^2/4}\nonumber\\
&\le Ca\nu N^{4/\chi}+Ca e^{-t_0^2/8}.
\label{eq:direct-moment-early}
\end{align}
For the remaining interval \eqref{eq:direct-Er-coarse} and
\eqref{eq:direct-R-properties} give
\begin{equation}\label{eq:direct-moment-late}
\int_{R_\nu}^{t_\nu}(1+a\dbla{s})\|E[{F^{\rm b}_{\rm err}}](t_0+t)\|_2\,d t
\le C\varepsilon_N\nu^2=CaL_N\nu^2\le Ca\nu .
\end{equation}
Finally
\[
N^{4/\chi}\le C{T_1}^3,\qquad
\int_0^{t_\nu}\dbla{t_0+t}^2\,d t\le C{T_1}^3,\qquad
Z_\ell(0)\le C_\ell a^2e^{-t_0^2}.
\]
The first inequality follows from ${T_1}\ge cN^{{\rm s}_0+4}$,
${\rm s}_0\ge11$, and $4/\chi<8$.
The integrating factor in \eqref{eq:direct-moment-scalar} is bounded
by \eqref{eq:direct-b}. The two field integrals
\eqref{eq:direct-moment-early}--\eqref{eq:direct-moment-late} thus give
\[
\sup_{0\le t\le t_\nu}Z_\ell(t)
\le C_\ell\big[Z_\ell(0)+a\nu T_1^3+a\nu N^{4/\chi}
+ae^{-t_0^2/8}+\varepsilon_N\nu^2\big] \le C_\ell a\big(\nu T_1^3+e^{-t_0^2/8}\big),
\]
where $N^{4/\chi}\le CT_1^3$ and
$\varepsilon_N\nu^2=aL_N\nu^2\le Ca\nu$. The bounded inverse
potential multiplier proves \eqref{eq:polynomialBackgroundMoments}.
\par
All integrations by parts can first be made with velocity cutoffs.
The fixed higher weighted regularity in {\eqref{eq:bguniformsobolev}} and the Gaussian
profile justify passage to the limit on each finite interval.
For every real $\alpha$, the shear has unit Jacobian and leaves the
velocity weight unchanged:
\[
\left\|\langle v\rangle^\ell
{f^{\rm b}_{\rm err}}(t_0+t,z+\alpha v_1,v)\right\|_{L^2_{z,v}}
=\left\|\langle v\rangle^\ell
{f^{\rm b}_{\rm err}}(t_0+t,x_1,v)\right\|_{L^2_{x_1,v}}.
\]
\end{proof}

\subsubsection{Interpolation, the density trace and the exact field}

The pointwise-in-time estimate \eqref{eq:polynomialBackgroundLow} must be used before taking
the density trace. This gives a positive power of $\nu$ in the weighted Fourier bound for $E^{\rm b,\nu}$ in \eqref{eq:fieldFourierEnvelope}, even though the observation time is long. We prove this bound explicitly, since it also supplies every integrated field bound
used by the later linear and nonlinear arguments.

\begingroup
\begin{lem}\label{lem:direct-gevrey}
Under \eqref{eq:comparisonLogParameters}, for every fixed $0<{\sigma_{\rm int}}<1$, there is a fixed $\lambda_{\sigma_{\rm int}}>0$ such that
\begin{equation}\label{eq:bgdifferencegevrey}
\|e^{\lambda_{\sigma_{\rm int}}\dbla{\nabla_z,\nabla_v}^\chi}{g^{\rm b}_{\rm err}}(t)\|_2
\le CaL_N\nu^{\sigma_{\rm int}}\dbla{t}^{3{\sigma_{\rm int}}},
\qquad 0\le t\le t_\nu .
\end{equation}
\end{lem}
\begin{proof}
Before absorbing powers of $t_0$ and weakening the Gaussian exponent,
the more precise bound is
\begin{equation}\label{eq:direct-precise-gevrey}
\|e^{\lambda_{\sigma_{\rm int}}\dbla{\nabla_z,\nabla_v}^\chi}{g^{\rm b}_{\rm err}}(t)\|_2
\le CaL_N^{1-{\sigma_{\rm int}}}
\{\nu t\dbla{t_0+t}^2+e^{-t_0^2/4}\}^{{\sigma_{\rm int}}}.
\end{equation}

Use the explicit normalized profile $g^{\rm b}_{\rm app}$ from
\eqref{eq:physicalCollisionProfiles}.
{Gaussian Fourier transformation gives}
$\|e^{\lambda_*\dbla{\nabla}^\chi}g^{\rm b}_{\rm app}\|_2\le aL_N$
at a fixed radius, after enlarging $L_N$. Together with Proposition~\ref{prop:bguniform}, this yields
\begin{equation}\label{eq:direct-high-D}
\|e^{\lambda_*\dbla{\nabla}^\chi}{g^{\rm b}_{\rm err}}(t)\|_2\le CaL_N.
\end{equation}
The shear preserves $L^2$, so \eqref{eq:polynomialBackgroundLow} gives
{$\|{g^{\rm b}_{\rm err}}(t)\|_2\le Ca(\nu t\dbla{t_0+t}^2+e^{-t_0^2/4})$}.
Set $\lambda_{\sigma_{\rm int}}=(1-{\sigma_{\rm int}})\lambda_*$. Fourier H\"older interpolation gives
\[
\|e^{\lambda_{\sigma_{\rm int}}\dbla{\nabla}^\chi}{g^{\rm b}_{\rm err}}\|_2
\le\|{g^{\rm b}_{\rm err}}\|_2^{\sigma_{\rm int}}
\|e^{\lambda_*\dbla{\nabla}^\chi}{g^{\rm b}_{\rm err}}\|_2^{1-\sigma_{\rm int}}.
\]
This proves \eqref{eq:direct-precise-gevrey}.
Use
\[
(x+y)^{\sigma_{\rm int}}\le x^{\sigma_{\rm int}}+y^{\sigma_{\rm int}},\quad
t\dbla{t_0+t}^2\le C\dbla{t_0}^2\dbla{t}^3,\quad
\log L_N=o(N^2),\quad t_0=A_0N
\]
to obtain the right-hand side of \eqref{eq:bgdifferencegevrey}
plus $CaL_Ne^{-c_{\sigma_{\rm int}}N^2}$. Finally,
\[
\frac{e^{-c_{\sigma_{\rm int}}N^2}}
{\nu^{\sigma_{\rm int}}\langle t\rangle^{3\sigma_{\rm int}}}
\le e^{-c_{\sigma_{\rm int}}N^2+\sigma_{\rm int}N/K}\longrightarrow0.
\]
This proves \eqref{eq:bgdifferencegevrey}. The Gaussian term in
\eqref{eq:direct-precise-gevrey} is kept because the collision term
vanishes at $t=0$.
\end{proof}
\par\endgroup

\begingroup
\begin{prop}
\label{prop:backgroundFieldBounds}
Under the hypotheses of Proposition~\ref{thm:bgcomparison}, fix
$0<\sigma_{\rm int}<1$. There are fixed $\varrho_{\rm F},c_{\rm dec},c>0$ and a factor $L_N$
{of the form specified in Proposition~\ref{prop:bguniform}} such that the following holds.
\begingroup
Recall the integrated field size
\[
\varepsilon_{\rm F}=aL_N\nu^{\sigma_{\rm int}}
\]
from \eqref{eq:fieldExpansionSizes}, with the same
$\sigma_{\rm int}$, $L_N$ and $c$.
\par\endgroup The estimates here hold for the
full range $0<\sigma_{\rm int}<1$; the finite field correction uses the
{choice $\sigma_{\rm int}>1/2$ made in Section~\ref{sec:roadmapFieldCompletion}.}
Write $E^{\rm b,\nu}(s)=E^\nu(s-t_0)=E[F^{\rm b,\nu}](s)$. Then
\begin{equation}\label{eq:fieldFourierEnvelope}
\mathfrak f(s)
:=\sum_{l\ne0}e^{\varrho_{\rm F}(|l|^\chi+|l(s-t_0)|^\chi)}
|E_l^{\rm b,\nu}(s)|
\le C\varepsilon_{\rm F}e^{-c_{\rm dec}(s-t_0)^\chi},\qquad t_0\le s\le {T_1}.
\end{equation}
For every fixed $q\ge0$, integer $j\ge0$, and $t_0\le\tau\le {T_1}$,
\begin{align}
\int_\tau^{T_1}\langle s\rangle^q\mathfrak f(s)\,ds
&\le C_q\varepsilon_{\rm F}\langle t_0\rangle^q
e^{-c_{{\rm dec},q}(\tau-t_0)^\chi},\label{eq:fieldFourierMoments}\\
\int_{t_0}^{T_1}\langle s\rangle^q
\|E^{\rm b,\nu}(s)\|_{C^j_x}\,ds
&\le aL_N\nu^{\sigma_{\rm int}},
\label{eq:exactFieldIntegral}\\
\int_{t_0}^{T_1}\|E^{\rm b,\nu}(s)\|_{C^j_x}^2\,ds
&\le a^2L_N\nu^{2\sigma_{\rm int}}.
\label{eq:directFieldSquare}
\end{align}
All constants are independent of $\nu,N$, but may depend on the fixed
indices and on $\sigma_{\rm int}$. The factor $L_N$ and the Gaussian
constant may be enlarged and reduced, respectively, for the fixed
indices in the last two bounds.
\end{prop}
\begin{proof}
Under \eqref{eq:comparisonLogParameters}, for every fixed $c>0$ and
$\sigma_{\rm int}\in(0,1)$,
\[
\frac{e^{-cN^2}}{\nu^{2\sigma_{\rm int}}}
=\exp\left(-cN^2+\frac{2\sigma_{\rm int}}K N\right)\longrightarrow0.
\]
Thus Gaussian errors can be absorbed into both displayed powers of $\nu$;
the estimate below still shows separately where these errors arise.

Use physical time $t=s-t_0$, so that $E^\nu(t)=E^{\rm b,\nu}(t_0+t)$.
The distribution-level error {$F^{\rm b}_{\rm err}$} has the Gaussian-weighted trace
\[
\begin{aligned}
\rho_k[{F^{\rm b}_{\rm err}}](t_0+t)
&=\int\sqrt \mu (v)e^{-iktv_1}{g}_{{\rm err},k}^{\rm b}(t,v)\,dv\\
&=(2\pi)^{-3}\int\widehat{\sqrt \mu }(kt e_1-\eta)
\widehat {g}_{{\rm err},k}^{\rm b}(t,\eta)\,d\eta.
\end{aligned}
\]
For a fixed sufficiently small $c>0$,
\[
c\dbla{k,kt e_1}^\chi
\le\lambda_{\sigma_{\rm int}}\dbla{k,\eta}^\chi+
c|kt e_1-\eta|^\chi.
\]
Cauchy--Schwarz, as in Lemma \ref{lem:direct-current}, gives
\begin{align}
|(E[{F^{\rm b}_{\rm err}}])_k(t_0+t)|
&\le \frac{C}{|k|}e^{-c\dbla{k,kt e_1}^\chi}
\|e^{\lambda_{\sigma_{\rm int}}\dbla{k,\eta}^\chi}\widehat {g}_{{\rm err},k}^{\rm b}(t)\|_2\nonumber\\
&\le\frac{Ca}{|k|}
\{L_N\nu^{\sigma_{\rm int}}\dbla{t}^{3{\sigma_{\rm int}}}+e^{-c_2N^2}\}
e^{-c_3(|k|^\chi+|kt|^\chi)}.
\label{eq:direct-Er-mode}
\end{align}
The second line follows from \eqref{eq:direct-precise-gevrey}, using
$t\langle t_0+t\rangle^2\le C\langle t_0\rangle^2\langle t\rangle^3$
and enlarging $L_N$ to absorb the fixed power of $t_0$.
For sufficiently large $N$, $t_0=A_0N$ and $\log L_N=o(N^2)$ give
\[
L_N^{1-\sigma_{\rm int}}e^{-\sigma_{\rm int}t_0^2/4}
\le e^{-c_2N^2},\qquad 0<c_2<\frac{\sigma_{\rm int}A_0^2}{4}.
\]
No pointwise bound on {${\widehat {g}_{{\rm err},k}^{\rm b}}$} is used.

The explicit field has only the modes $k=\pm1$, with
$|(E[F^{\rm b}_{\rm app}])_{\pm1}|=(a/2)e^{-(t_0+t)^2/2}$.
For any fixed sufficiently small $c_4>0$,
\begin{equation}\label{eq:direct-EB-envelope}
e^{-(t_0+t)^2/2}
\le C e^{-t_0^2/4}e^{-c_4t^\chi}.
\end{equation}
Choose $\varrho_{\rm F}<c_3/4$, and also small enough for \eqref{eq:direct-EB-envelope}.
Since $|k|\ge1$,
\[
\dbla{t}^{3{\sigma_{\rm int}}}e^{-c_3|kt|^\chi/2}
\le C_{\sigma_{\rm int}} e^{-c_{\rm dec}t^\chi},\qquad
\sum_{k\ne0}|k|^{-1}e^{-c_3|k|^\chi/2}<\infty.
\]
Adding {$E^\nu=E[{F^{\rm b}_{\rm err}}]+E[F^{\rm b}_{\rm app}]$} in \eqref{eq:direct-Er-mode} proves
\eqref{eq:fieldFourierEnvelope}.

For fixed $q\ge0$ and every $x\ge0$, substitution $u=t^\chi$ and
absorption of a fixed polynomial into the exponential imply
\[
\int_x^\infty\dbla{t}^q e^{-c_{\rm dec}t^\chi}\,d t
\le C_q e^{-c_{{\rm dec},q}x^\chi}.
\]
Taking $x=\tau-t_0$, changing variables by $s=t_0+t$, and using
$\dbla{t_0+t}^q\le C_q\dbla{t_0}^q\dbla{t}^q$ gives
\eqref{eq:fieldFourierMoments}.
Moreover,
\[
\|E^\nu(t)\|_{C^j_x}
\le\sum_{k\ne0}|k|^j|E^\nu_k(t)|
\le C_{j,\varrho_{\rm F}}\mathfrak f(t_0+t).
\]
Integrating proves \eqref{eq:exactFieldIntegral}, after absorbing the fixed
power of $t_0$ into $L_N$ and reducing the Gaussian exponent.
Squaring and integrating \eqref{eq:fieldFourierEnvelope} gives
\[
\int_0^{t_\nu}\|E^\nu\|_{C^j_x}^2\,d t
\le C_j a^2(L_N^2\nu^{2{\sigma_{\rm int}}}+e^{-2cN^2})
\int_0^\infty e^{-2c_{\rm dec}t^\chi}\,d t .
\]
The square of a factor {of the form specified in Proposition~\ref{prop:bguniform}} has the same form.
This proves \eqref{eq:directFieldSquare}. Finally,
\[
\zeta\varepsilon_{\rm F}
=N^3L_N\nu^{\sigma_{\rm int}}
\le\nu^{\sigma_{\rm int}-o(1)},
\]
since $a\zeta=N^3$, $\log L_N=o(N)$ and $N\mathrel{{\sim}}\log(1/\nu)$.
This also verifies the smallness stated in \eqref{eq:fieldExpansionSizes}.
\end{proof}
\par\endgroup

\begingroup
{The fixed-radius bound}
\eqref{eq:bgdifferencegevrey} is not a sufficiently small relative
error on the whole interval; Proposition~\ref{prop:decreasingRadiusBackground}
decreases the coefficient radius and retains the full factor
$a\nu T_1^3$.
\par\endgroup

\begingroup
\subsection{The background error at a decreasing Gevrey radius}\label{sec:backgroundDecreasingRadius}
The collision comparison will require the background error in $(z,v)$
coordinates with its small amplitude, rather than the larger fixed-radius
bound $aL_N$. We interpolate the low moment estimate with the fixed-radius
Gevrey estimate. For a function $V=V(v)$ independent of the spatial variable,
define
\begin{equation}\label{eq:velocityInputGevrey}
\|V\|_{\mathscr A_\ell^{\lambda,\chi}}^2
=\sum_{n\ge0}\frac{\lambda^{2n}(n+1)^{24}}{(n!)^{2/\chi}}
\|\langle v\rangle^\ell\partial_{v_1}^nV\|_2^2.
\end{equation}
We apply this norm to individual coefficients of the Fourier transform
in $x$ of the transport profile $g^{\rm b}_{\rm err}$ from
\eqref{eq:physicalCollisionProfiles}. The velocity radius decreases like
$N^{-1}$ and the Gevrey radius for the Fourier transform in $x$ like $N^{-\chi}$.
\begin{prop}\label{prop:decreasingRadiusBackground}
{Assume \eqref{eq:comparisonLogParameters}, and fix an integer
$\ell\ge4$. For a sufficiently small fixed velocity radius
$\lambda_{\rm v}>0$, there are constants $c,C>0$, independent of
$N,\nu$, with the following property. Define}
\begin{equation}\label{eq:decreasingRadiusParameters}
\lambda_{{\rm v},N}=\lambda_{\rm v}/N,\qquad
\varrho_N=cN^{-\chi}.
\end{equation}
{Then the error profile {${g^{\rm b}_{\rm err}}$} in
\eqref{eq:physicalCollisionProfiles} satisfies}
\begin{equation}\label{eq:decreasingRadiusBackground}
\begin{split}
&\|{g}_{{\rm err},l}^{\rm b}(t)\|_{\mathscr A_\ell^{\lambda_{{\rm v},N},\chi}}
+\|v_1{g}_{{\rm err},l}^{\rm b}(t)\|_{\mathscr A_\ell^{\lambda_{{\rm v},N},\chi}} +N^{-1}\|\partial_{v_1}{g}_{{\rm err},l}^{\rm b}(t)\|_{\mathscr A_\ell^{\lambda_{{\rm v},N},\chi}}\\
&\hspace{20mm}\le Ca\nu T_1^3e^{-\varrho_N\langle l\rangle^\chi},
\qquad 0\le t\le {T_1}-t_0,\quad l\in\mathbb Z.
\end{split}
\end{equation}
In particular, at $\tau=t_0+t$, the coefficient involving $g^{\rm b}_{\rm err}$ obeys
\begin{equation}\label{eq:decreasingRadiusElectricInput}
\left\|(\partial_{v_1}-v_1/2-ilt){g}_{{\rm err},l}^{\rm b}(t)\right\|_{\mathscr A_\ell^{\lambda_{{\rm v},N},\chi}}
\le Ca\nu T_1^3(N+|l|\tau)
e^{-\varrho_N\langle l\rangle^\chi}.
\end{equation}
\end{prop}

\begin{proof}
The error profile {${g^{\rm b}_{\rm err}}$} is {${f^{\rm b}_{\rm err}}$} from
Proposition~\ref{thm:bgcomparison}, expressed in free-transport coordinates. {The small moment bounds} \eqref{eq:polynomialBackgroundMoments} therefore give, for every fixed integer $J_{\rm mom}$,
\begin{equation}\label{eq:decreasingRadiusLowMoments}
\sum_{|\alpha|\le J_{\rm mom}}\|v^\alpha {g^{\rm b}_{\rm err}}(t)\|_2
\le C_{J_{\rm mom}}a\nu T_1^3.
\end{equation}
The Gaussian error in \eqref{eq:polynomialBackgroundMoments} is smaller:
\[
\frac{e^{-cN^2}}{\nu T_1^3}
\le \frac{e^{-cN^2}}{\nu}
=\exp(-cN^2+N/K)\longrightarrow0.
\]
Thus the low-moment error is bounded by ${Ca\nu T_1^3}$. We next obtain a high norm with a fixed radius. Proposition~\ref{prop:bguniform} and the explicit profile give {$\|e^{\lambda_*\langle\nabla_z,\nabla_v\rangle^\chi}{g^{\rm b}_{\rm err}}\|_2\le aL_N$} and every required fixed velocity moment bounded by $aL_N$. Apply Lemma~\ref{lem:weightedGevreyInterpolation} {with smallness parameter equal to one. The resulting per-mode bounds,} summed after a further fixed reduction of both radii, imply that for some $\lambda_*>0$,
\begin{equation}\label{eq:decreasingRadiusHighMoments}
\sum_{|\alpha|\le J_{\rm mom}}
\|e^{\lambda_*\langle\nabla_z,\nabla_v\rangle^\chi}(v^\alpha {g^{\rm b}_{\rm err}}(t))\|_2
\le aL_N.
\end{equation}
Here $J_{\rm mom}$ is fixed, and the constants and radius may depend on $J_{\rm mom}$.

For each $|\alpha|\le J_{\rm mom}$, H\"older's inequality in the Fourier variables, using \eqref{eq:decreasingRadiusLowMoments} and \eqref{eq:decreasingRadiusHighMoments}, gives
\begin{equation}\label{eq:decreasingRadiusInterpolation}
\begin{split}
\|e^{N^{-\chi}\lambda_*\langle\nabla_z,\nabla_v\rangle^\chi}
(v^\alpha {g^{\rm b}_{\rm err}}(t))\|_2
&\le(C_{J_{\rm mom}}a\nu T_1^3)^{1-N^{-\chi}}(aL_N)^{N^{-\chi}}\\
&\le Ca(\nu T_1^3)^{1-N^{-\chi}}\le Ca\nu T_1^3.
\end{split}
\end{equation}
Indeed $L_N^{N^{-\chi}}\le C$, because $N^{-\chi}\log L_N\le C+CN^{-\chi}\log(1+N)$. Since $\nu T_1^3\mathrel{{\sim}} N^{-3(p-{\rm s}_0-2)}$,
\[
(\nu T_1^3)^{-N^{-\chi}}=\exp\big(O(N^{-\chi}\log N)\big)\le C.
\]

Choose $c_3>0$ so that, for all $l,\eta$,
\[
c_3\langle l\rangle^\chi+c_3\langle\eta\rangle^\chi
\le\lambda_*\langle l,\eta\rangle^\chi.
\]
Applying this Fourier multiplier inequality to
\eqref{eq:decreasingRadiusInterpolation} gives, for $|\alpha|\le J_{\rm mom}$,
\par

\[
\|e^{c_3N^{-\chi}\langle\nabla_v\rangle^\chi}(v^\alpha {g}_{{\rm err},l}^{\rm b}(t))\|_2
\le Ca\nu T_1^3e^{-c_3N^{-\chi}\langle l\rangle^\chi}.
\]
Since $\ell$ is an integer, choose $J_{\rm mom}\ge\ell+2$.
For $m\ge0$ and $|\alpha|\le J_{\rm mom}$, Fourier multiplication gives
\begin{align*}
\sup_\eta|\eta_1|^m e^{-c_3N^{-\chi}\langle\eta\rangle^\chi/2}
&\le (CN)^m(m!)^{1/\chi},\\
\|\partial_1^m(v^\alpha {g}_{{\rm err},l}^{\rm b})\|_2
&\le Ca\nu T_1^3e^{-c_3N^{-\chi}\langle l\rangle^\chi}(CN)^m(m!)^{1/\chi}.
\end{align*}
The exact commutation formula is
\[
v^\alpha\partial_1^ng_{{\rm err},l}^{\rm b}
=\sum_{h=0}^{\min(n,\alpha_1)}(-1)^h\binom nh
\frac{\alpha_1!}{(\alpha_1-h)!}
\partial_1^{n-h}(v^{\alpha-he_1}{g}_{{\rm err},l}^{\rm b}).
\]
Since $N\ge1$, $h\le J_{\rm mom}$, and
$\langle v\rangle^\ell\le C_\ell\sum_{|\alpha|\le \ell}|v^\alpha|$,
this yields
\[
\begin{aligned}
\|\langle v\rangle^\ell\partial_1^ng_{{\rm err},l}^{\rm b}\|_2
&\le Ca\nu T_1^3e^{-c_3N^{-\chi}\langle l\rangle^\chi}(CN)^n(n!)^{1/\chi}(n+1)^{J_{\rm mom}},\\
N^{-1}\|\langle v\rangle^\ell\partial_1^{n+1}{g}_{{\rm err},l}^{\rm b}\|_2
&\le Ca\nu T_1^3e^{-c_3N^{-\chi}\langle l\rangle^\chi}(CN)^n(n!)^{1/\chi}(n+1)^{J_{\rm mom}+1/\chi},\\
\|\langle v\rangle^\ell\partial_1^n(v_1{g}_{{\rm err},l}^{\rm b})\|_2
&\le Ca\nu T_1^3e^{-c_3N^{-\chi}\langle l\rangle^\chi}(CN)^n(n!)^{1/\chi}(n+1)^{J_{\rm mom}+1}.
\end{aligned}
\]
For the last line use
{${\partial_1^n(v_1{g}_{{\rm err},l}^{\rm b})=v_1\partial_1^ng_{{\rm err},l}^{\rm b}+n\partial_1^{n-1}{g}_{{\rm err},l}^{\rm b}}$};
the first term uses moments through $\ell+1\le J_{\rm mom}$.
Take $C\lambda_{\rm v}<1/2$. Each of the three norms in
\eqref{eq:decreasingRadiusBackground}, with the factor $N^{-1}$
attached to the velocity-derivative term, has square bounded by
\[
\begin{aligned}
&Ca^2(\nu T_1^3)^2e^{-2c_3N^{-\chi}\langle l\rangle^\chi}\sum_{n\ge0}
\frac{(\lambda_{\rm v}/N)^{2n}(CN)^{2n}(n!)^{2/\chi}
(n+1)^{2J_{\rm mom}+2+2/\chi}}{(n!)^{2/\chi}(n+1)^{-24}}\\
&\quad=Ca^2(\nu T_1^3)^2e^{-2c_3N^{-\chi}\langle l\rangle^\chi}\sum_{n\ge0}(C\lambda_{\rm v})^{2n}(n+1)^{2J_{\rm mom}+26+2/\chi}\\
&\le C_{J_{\rm mom}}a^2(\nu T_1^3)^2e^{-2c_3N^{-\chi}\langle l\rangle^\chi}.
\end{aligned}
\]
The constant is independent of $N,l,t$. Choose $c\le c_3$ and
take square roots to prove \eqref{eq:decreasingRadiusBackground}.
Finally, the triangle inequality gives
\[
\begin{aligned}
&\|(\partial_1-v_1/2-ilt){g}_{{\rm err},l}^{\rm b}\|_{\mathscr A_\ell^{\lambda_{\rm v}/N,\chi}} \le Ca\nu T_1^3e^{-c_3N^{-\chi}\langle l\rangle^\chi}(N+1+|l|t)\\
&\quad\le Ca\nu T_1^3(N+|l|\tau)
 e^{-\varrho_N\langle l\rangle^\chi},
\end{aligned}
\]
which is \eqref{eq:decreasingRadiusElectricInput}.
\end{proof}
\par\endgroup
\section{\texorpdfstring{{The approximate density equation}}{The approximate density equation}}\label{sec:approximateDensity}
{We analyze the density equation obtained from free transport
and the explicit background $F^{\rm b}_{\rm app}$.
Section~\ref{sec:chain} constructs the signed Maxwellian inverse and
proves the descending-chain lower bound.
Section~\ref{sec:gainDensityInverse} then inverts the full approximate
density equation, including the upward-interaction term, and establishes the
estimates in $\mathscr X_N^{r_{\rm x}}(I)$ and for $B_*^{-1}\rho_{m_1}(T_1)$ needed for the exact comparison. Collision propagation and the complete linearized
VPL equation are treated in Section~\ref{sec:completeLinear}.}

\begingroup
\subsection{\texorpdfstring{{The signed Maxwellian kernel and the descending chain}}{The signed Maxwellian kernel and the descending chain}}\label{sec:collisionless}\label{sec:estimates}\label{sec:chain}\label{sec:linear}
We first construct the scalar inverse used in
Section~\ref{sec:roadmapOperators}, with the weighted bounds needed for
the echo chain. {We then identify the density} and follow the
causal transfers in \eqref{eq:reduced}. Normalizing each transfer
isolates a Gaussian convolution; {the accumulated comparison error}
yields Proposition~\ref{thm:chain}.
\par\endgroup
In the variable $w=ks$, the kernel is
\begin{equation}\label{eq:D}
D_k(w)=k^{-1}K_k(w/k)=-k^{-2}w e^{-w^2/2}\mathbf1_{w\ge0}.
\end{equation}
Convolution below is on $\mathbb R$, with kernels supported in $[0,\infty)$.

\begingroup
\begin{lem}\label{prop:inverse}
There exist $r_*>0$ and $C$ such that, for every integer $k\ge1$, a real continuous kernel $\widetilde R_k$ supported in $[0,\infty)$ satisfies
\begin{equation}\label{eq:inverse}
\widetilde R_k=D_k+D_k*\widetilde R_k,
\quad |\widetilde R_k(w)|\le Ck^{-2}e^{-r_*w}\quad(w\ge0),
\quad \mathbb V_k U=U+\widetilde R_k*U=({\mathrm{Id}}-D_k*)^{-1}U.
\end{equation}
The operator identity holds on $L^\infty(\mathbb R)$ and on the spaces with finite norm
$\operatorname*{ess\,sup}_w e^{\lambda\langle w\rangle}|U(w)|$, for $0\le\lambda\le r_*/2$.
For these $\lambda$, put
\[
q_\lambda=\int_0^\infty w e^{-w^2/2+\lambda w}\,dw.
\]
If $k^2>q_\lambda$, then
\begin{equation}\label{eq:moment}
\int_0^\infty e^{\lambda w}|\widetilde R_k(w)|\,dw
\le\frac{q_\lambda}{k^2-q_\lambda}.
\end{equation}
In particular, for $k\mathrel{{\sim}} N$ and $\lambda\le C_0/N$, this norm is at most $C/N^2$. The kernel is not nonnegative.
For $J_k\in C([t_0,{T_1}])$, the unique solution
$\rho_k\in C([t_0,{T_1}])$ of \eqref{eq:scalar} is given by
\[
\rho_k(s)=J_k(s)+\int_{t_0}^s k\widetilde R_k(k(s-\tau))J_k(\tau)\,d\tau.
\]

\end{lem}

\begin{proof}[Proof of Lemma~\ref{prop:inverse}]
The Volterra resolvent framework is given in
\cite[Appendix D, Lemma D.2]{ref3}. We verify here the uniform dependence
on the spatial frequency, the exponential decay, and the weighted bounds
needed below.
Define the entire function
\[
\mathfrak k(z)=\int_0^\infty w e^{-w^2/2}e^{-zw}\,dw.
\]
For $\operatorname{Re}z\ge0$, $|\mathfrak k(z)|\le1$. Equality requires $\operatorname{Re}z=0$ and a constant phase on the positive half-line, so equality occurs only at $z=0$. There $\mathfrak k(0)=1$. Consequently
\[
1+\lambda \mathfrak k(z)\ne0\qquad
(0\le\lambda\le1,\ \operatorname{Re}z\ge0).
\]
For $z\ne0$, integration by parts twice gives
\[
\mathfrak k(z)=\frac1{z^2}\left[1+\int_0^\infty
(w^3-3w)e^{-w^2/2}e^{-zw}\,dw\right].
\]
Thus $|\mathfrak k(z)|\le C|z|^{-2}$ for $\operatorname{Re}z\ge-1$ and
$|z|\ge1$. Combining this decay with the preceding nonvanishing
statement and compactness gives
\[
\inf_{\substack{0\le\lambda\le1\\\operatorname{Re}z\ge0}}
|1+\lambda\mathfrak k(z)|\ge c>0.
\]
\par
Also
\[
|\mathfrak k(z)-\mathfrak k(z+r)|\le
\int_0^\infty w e^{-w^2/2+rw}(1-e^{-rw})\,dw\le Cr
\]
when $\operatorname{Re}z\ge-r$ and $0<r\le1$.
Choose $0<r_*\le1$ so that $Cr_*\le c/2$. Then
\[
\inf_{\substack{k\ge1\\\operatorname{Re}z\ge-r_*}}
|1+k^{-2}\mathfrak k(z)|\ge c/2.
\]
\par

Now define
\begin{equation}\label{eq:laplaceR}
\widetilde R_k(w)=\frac1{2\pi i}
\int_{-r_*-i\infty}^{-r_*+i\infty}
e^{zw}\frac{-k^{-2}\mathfrak k(z)}{1+k^{-2}\mathfrak k(z)}\,dz.
\end{equation}
On this line the fraction in \eqref{eq:laplaceR} is bounded by
$Ck^{-2}(1+|\operatorname{Im}z|^2)^{-1}$, so
\[
|\widetilde R_k(w)|\le Ck^{-2}e^{-r_*w}
\int_{\mathbb R}\frac{dy}{1+y^2}\le Ck^{-2}e^{-r_*w}.
\]
Dominated convergence gives continuity. For $w<0$, move the contour
to $\operatorname{Re}z=R\ge1$. The horizontal integrals tend to zero
by the $O(|z|^{-2})$ bound, and
\[
|\widetilde R_k(w)|\le Ck^{-2}e^{Rw}
\int_{\mathbb R}\frac{dy}{R^2+y^2}
\le Ck^{-2}\frac{e^{Rw}}R\longrightarrow0\qquad(R\to\infty).
\]
Hence $\widetilde R_k(w)=0$ for $w<0$. Since
$\mathfrak k(\overline z)=\overline{\mathfrak k(z)}$, conjugating
\eqref{eq:laplaceR} and replacing $y$ by $-y$ gives
$\overline{\widetilde R_k(w)}=\widetilde R_k(w)$.
Multiplication of the Laplace transform of $\widetilde R_k$ by
$1+k^{-2}\mathfrak k$ yields the transform of $D_k$, proving the
first identity in \eqref{eq:inverse}. The exponential bounds justify
the convolution products.
\par

The exponential bound makes the left side of \eqref{eq:moment} finite. Take absolute values in $\widetilde R_k=D_k+D_k*\widetilde R_k$ and integrate with weight $e^{\lambda w}$. This gives
\[
\|e^{\lambda w}\widetilde R_k\|_1
\le k^{-2}q_\lambda
\left(1+\|e^{\lambda w}\widetilde R_k\|_1\right),
\]
which proves \eqref{eq:moment}. The inverse identity follows by multiplication of the convolution operators. Near zero, the convolution term is $O(w^3)$, whereas $D_k(w)=-w/k^2+O(w^3)$; thus $\widetilde R_k$ is negative for small positive $w$.

Substitution and the resolvent identity verify the formula. If
$e=K_{|k|}*e$ on $I$, then
\[
|e(s)|\le\|K_{|k|}\|_{L^\infty(0,{T_1}-t_0)}
\int_{t_0}^s|e(\tau)|\,d\tau.
\]
Gronwall's inequality gives $e=0$ on the stated finite interval.

\end{proof}
\par\endgroup

The explicit density calculation in
Section~\ref{sec:roadmapOperators} extends to an equation with a
prescribed source. We record the source contribution and the scalar
inverse before estimating the descending transfers.
\par

\begingroup
Recall the kernel mapping $\rho_m$ to the contribution in mode $k$,
\[
P_{k,m}(s,\tau)=-\frac{ak}{2m}(s-\tau)e^{-(ks-m\tau)^2/2},
\qquad k,m\ne0,\quad |k-m|=1,
\]
from \eqref{eq:overviewApproximateDensityKernels}.
\par\endgroup

\begingroup
\begin{lem}\label{lem:fullDensityIdentity}
Suppose that a smooth, zero-mass physical perturbation $H$, with $H(t_0)=H_{\rm in}$, satisfies
\begin{equation}\label{eq:applinear}
\partial_sH+v_1\partial_{x_1}H+
E[H]\partial_{v_1}F^{\rm b}_{\rm app}=S.
\end{equation}
{Then the actual density $\rho=\rho[H]$ satisfies, for $k\ne0$,}
\begin{equation}\label{eq:densityexact}
\rho_k(s)=\mathfrak g_k(s)+J_k(s)+\int_{t_0}^sK_{|k|}(s-\tau)\rho_k(\tau)\,d\tau
+\sum_{\substack{m\ne0\\|k-m|=1}}\int_{t_0}^sP_{k,m}(s,\tau)\rho_m(\tau)\,d\tau.
\end{equation}
\par
where
\begin{equation}\label{eq:forcing}
{
\mathfrak g_k(s)=\widehat H_{{\rm in},k}(k(s-t_0)e_1),\qquad J_k(s)=\int_{t_0}^s\widehat S_k(\tau,k(s-\tau)e_1)\,d\tau.
}
\end{equation}
For the scalar inversion and individual mode estimates, write, for $k,m\ne0$ with $|k-m|=1$,
\[
(\mathbb K_k\rho_k)(s)=\int_{t_0}^sK_{|k|}(s-\tau)\rho_k(\tau)\,d\tau,
\qquad
(\mathbb P_{k,m}\rho_m)(s)=\int_{t_0}^sP_{k,m}(s,\tau)\rho_m(\tau)\,d\tau.
\]
The scalar solution map is
$\mathbb W_k=({\mathrm{Id}}-\mathbb K_k)^{-1}$ on $C(I)$, with
$\sup_{k\ne0}\|\mathbb W_k\|_{C(I)\to C(I)}\le C$.
Let $\mathbb W=({\mathrm{Id}}-\mathbb K)^{-1}$ act mode by mode.
Using the density integral operators already defined in
\eqref{eq:overviewApproximateDensityEquation}, no contribution arises
from the zero field mode.
Thus \eqref{eq:densityexact} is equivalently
\begin{equation}\label{eq:allmodeidentity}
\rho=\mathfrak g+J+\mathbb K\rho+\mathbb P_\downarrow\rho+\mathbb P_\uparrow\rho,\qquad \rho=\mathbb W(\mathfrak g+J)+\mathbb W\mathbb P_\downarrow\rho
+\mathbb W\mathbb P_\uparrow\rho.
\end{equation}
Section~\ref{sec:gainDensityInverse} strengthens these inverse bounds
in the density weights used for comparison with the exact-background
solution operator.
\end{lem}
\par
\par\endgroup

\begingroup
\begin{proof}[Proof of Lemma~\ref{lem:fullDensityIdentity}]
The Fourier coefficients of the explicit background are
{$(F^{\rm b}_{\rm app})_0=\mu $} and
{$(F^{\rm b}_{\rm app})_{\pm1}=(a/2)\mu e^{\mp i\tau v_1}$}.
{The mode equation and the corresponding free Duhamel formula are}
\begin{align*}
(\partial_s+ikv_1)H_k
 &=S_k+\sum_{m\ne0}\frac{i\rho_m}{m}\partial_{v_1}(F^{\rm b}_{\rm app})_{k-m},\\
H_k(s)&=e^{-ik(s-t_0)v_1}H_k(t_0)
\\
&\quad+\int_{t_0}^se^{-ik(s-\tau)v_1}
 \left(S_k(\tau)+\sum_{m\ne0}\frac{i\rho_m(\tau)}m
 \partial_{v_1}(F^{\rm b}_{\rm app})_{k-m}(\tau)\right)d\tau.
\end{align*}
The velocity integrals are those computed in
Section~\ref{sec:roadmapEcho}, giving the kernels in
\eqref{eq:overviewApproximateDensityKernels}.
\par
Integration of the Duhamel formula in $v$ gives
\eqref{eq:densityexact}--\eqref{eq:forcing}, including their signs.
Lemma~\ref{prop:inverse} gives the scalar time inverse
\[
(\mathbb W_k {\rho})(s)={\rho}(s)+
\int_{t_0}^s |k|\widetilde R_{|k|}(|k|(s-\tau)){\rho}(\tau)\,d\tau.
\]
Extending $\rho$ by zero outside $I$, the change of variable $w=|k|s$ gives
\par
\[
(\mathbb W_k {\rho})(s)=(\mathbb V_{|k|}[{\rho}(\cdot/|k|)])(|k|s).
\]
Consequently,
\[
\|\mathbb W_k\|_{C(I)\to C(I)}
\le1+\int_0^\infty|k|\,|\widetilde R_{|k|}(|k|t)|\,dt
=1+\|\widetilde R_{|k|}\|_1\le C.
\]
Applying $\mathbb W_k=({\mathrm{Id}}-K_{|k|}*)^{-1}$ gives
\[
\rho_k=\mathbb W_k(\mathfrak g_k+J_k)
+\sum_{\substack{m\ne0\\|k-m|=1}}\mathbb W_k\mathbb P_{k,m}\rho_m.
\]
Separating $|m|=|k|+1$ and $|m|=|k|-1$ proves
\eqref{eq:allmodeidentity}. The index $m=0$ is absent because $E_0=0$.
\end{proof}
\par\endgroup

\begingroup
Fix $\vartheta>0$, and write
\begin{equation}\label{eq:weights}
\lambda_m=\frac{\vartheta}{N}\left(1+\frac{m-m_1}{N}\right),\qquad
\|U\|_{\mathscr C_\lambda}
=\operatorname*{ess\,sup}_{w\in\mathbb R}e^{\lambda\langle w\rangle}|U(w)|.
\end{equation}
Take $N$ large enough that $\lambda_{m_0}<r_*/2$, with $r_*$ from Lemma~\ref{prop:inverse}. For bounded measurable functions the following integral is defined pointwise:
\begin{equation}\label{eq:Bplus}
[\mathbb B_m^\zeta U](w)=\frac{a}{2m^3}\int_{\mathbb R}
[\zeta+u+m(w-u)]_+e^{-(w-u)^2/2}U(u)\,du,
\qquad [y]_+=\max(y,0).
\end{equation}
The positive part enforces the time ordering; it is not an assumption that the density is positive. The integral need not define a bounded function for an arbitrary bounded function $U$. The weighted domains just defined will justify each use of $\mathbb V_{m-1}$.
\par\endgroup
\begin{lem}\label{lem:identity}
Under \eqref{eq:parameters} and the time and weight restrictions above, the finite system \eqref{eq:reduced} has a unique continuous solution on $[t_0,{T_1}]$. {The source in the top mode is} the free-transport density of \eqref{eq:highdatum}.

Set $w_m=mt_0-\zeta$ and
\begin{equation}\label{eq:initialtrace}
U_{m_0}=\mathbb V_{m_0}\big({\left(\frac{e^{-(\cdot)^2/4}}{\sqrt{4\pi}}\right)}\mathbf1_{[w_{m_0},\infty)}\big),
\qquad U_{m-1}=\mathbb V_{m-1}\mathbb B_m^\zeta U_m.
\end{equation}
Then $U_m(w)=0$ for $w<w_m$, and
\begin{equation}\label{eq:reducedidentity}
\rho_m^\downarrow(s)=(-1)^{m_0-m}U_m(ms-\zeta),\qquad s\in[t_0,{T_1}].
\end{equation}
All modes below $m_0$ are continuous. The value of $U_{m_0}$ at $w_{m_0}$ is taken from the right, consistently with \eqref{eq:reduced}.
\end{lem}

\begingroup
\begin{proof}
The top equation has a unique continuous solution by Lemma~\ref{prop:inverse}. Once a higher mode is known, the forcing of the next lower equation is continuous, so the same lemma constructs that mode uniquely. This finite induction gives existence and uniqueness. {The function $\mathfrak g^0_{m_0}$ is given by the direct Gaussian calculation \eqref{eq:freeGaussianDensity}.}

At the top mode, substitute $w=m_0s-\zeta$ and $u=m_0\tau-\zeta$ into the first equation of \eqref{eq:reduced}. It becomes
\[
U_{m_0}(w)={\frac{e^{-w^2/4}}{\sqrt{4\pi}}}+\int_{w_{m_0}}^wD_{m_0}(w-u)U_{m_0}(u)\,du
\quad(w\ge w_{m_0}).
\]
Lemma~\ref{prop:inverse} proves {the top identity and the support assertion. The exponential moment bound for the resolvent} also shows that $U_{m_0}\in\mathscr C_{\lambda_{m_0}}$, since the truncated Gaussian has finite norm there.

\begingroup
Suppose that \eqref{eq:reducedidentity} and
$U_m\in\mathscr C_{\lambda_m}$ hold. Put
$s=(w+\zeta)/(m-1)$. In the adjacent-mode term, set
$u=m\tau-\zeta$, so $d\tau=du/m$. Then
\begin{align*}
&(-1)^{m_0-m+1}\int_{t_0}^{s}P_{m-1,m}(s,\tau)
\rho_m^\downarrow(\tau)\,d\tau\\
&\quad=\frac{a(m-1)}{2m}\int_{t_0}^{s}
(s-\tau)e^{-((m-1)s-m\tau)^2/2}U_m(m\tau-\zeta)\,d\tau\\
&\quad=\frac{a}{2m^3}
\int_{w_m}^{\frac m{m-1}(w+\zeta)-\zeta}
[\zeta+u+m(w-u)]e^{-(w-u)^2/2}U_m(u)\,du.
\end{align*}
The integrand coefficient and the upper endpoint satisfy
\[
\zeta+u+m(w-u)=m(m-1)(s-\tau),\qquad
u\le\frac m{m-1}(w+\zeta)-\zeta
\iff\zeta+u+m(w-u)\ge0.
\]
Since $U_m(u)=0$ for $u<w_m$, the same integral equals
\[
\frac{a}{2m^3}\int_{\mathbb R}
[\zeta+u+m(w-u)]_+e^{-(w-u)^2/2}U_m(u)\,du
=\mathbb B_m^\zeta U_m(w).
\]
For the diagonal integral use $u=(m-1)\tau-\zeta$ and
$d\tau=du/(m-1)$. With
$U_{m-1}(w)=(-1)^{m_0-m+1}\rho_{m-1}^\downarrow(s)$ this gives
\[
U_{m-1}(w)=\int_{w_{m-1}}^w D_{m-1}(w-u)U_{m-1}(u)\,du
+\mathbb B_m^\zeta U_m(w).
\]
\par\endgroup
Before applying the inverse, put $z=w-u$ and use $\lambda_m-\lambda_{m-1}>0$. Since
\[
e^{\lambda_{m-1}\langle w\rangle-\lambda_m\langle u\rangle}
\le e^{\lambda_{m-1}|z|}e^{-(\lambda_m-\lambda_{m-1})\langle u\rangle},
\]
the inequality $[\zeta+u+mz]_+\le\zeta+|u|+m|z|$ gives
\[
\|\mathbb B_m^\zeta U_m\|_{\mathscr C_{\lambda_{m-1}}}
\le \frac{Ca}{m^3}\big(\zeta+(\lambda_m-\lambda_{m-1})^{-1}+m\big)
\|U_m\|_{\mathscr C_{\lambda_m}}<\infty.
\]
Here $\sup_u|u|e^{-(\lambda_m-\lambda_{m-1})\langle u\rangle}\le C/(\lambda_m-\lambda_{m-1})$ and the remaining Gaussian $z$ integrals are finite. Lemma~\ref{prop:inverse} therefore applies in this weighted space and defines $U_{m-1}\in\mathscr C_{\lambda_{m-1}}$. This argument proves well-definedness for the finite recursion; the gain-normalized bounds uniform in $N$ are proved next.

For $w<w_{m-1}$,
\[
\frac{m}{m-1}(w+\zeta)-\zeta<mt_0-\zeta=w_m,
\qquad \mathbb B_m^\zeta U_m(w)=0.
\]
The causal support of $\widetilde R_{m-1}$ therefore gives
$U_{m-1}(w)=0$ for $w<w_{m-1}$. Applying the inverse identity recovers
the transformed Volterra equation, and uniqueness proves
\eqref{eq:reducedidentity}. The weighted bound above supplies a locally
uniform integrable majorant for the Gaussian integral; dominated
convergence gives continuity of $\mathbb B_m^\zeta U_m$ and of every subsequent mode.
\par
\end{proof}
\par\endgroup
The indicator in \eqref{eq:initialtrace} extends a time trace by zero before the initial time. It does not cut off the Fourier transform of the datum \eqref{eq:highdatum}. In particular, the initial distribution {still has Gaussian velocity decay.}

\begingroup
For every variance $\alpha>0$, completing the square gives
\begin{equation}\label{eq:palpha}
\int_{\mathbb R}\frac{e^{-(w-u)^2/2}}{\sqrt{2\pi}}
\frac{e^{-u^2/(2\alpha)}}{\sqrt{2\pi\alpha}}\,du
=\frac{e^{-w^2/(2(1+\alpha))}}{\sqrt{2\pi(1+\alpha)}}.
\end{equation}
Taking $\alpha=1$ gives the top Gaussian on the right; each subsequent convolution adds one to the variance.
\par\endgroup

To estimate the recursion we compare each transfer with convolution by ${(2\pi)^{-1/2}e^{-w^2/2}}$, {after dividing by the single-step gain $\alpha_m$.} A decreasing exponential weight absorbs the unbounded integration variable in the actual coefficient. A weighted supremum norm is sufficient and also accommodates the initial jump of the extended trace.

\begin{lem}\label{lem:step}
Uniformly for $m_1<m\le m_0$ and sufficiently large $N$,
\begin{align}
\|{\left(\frac{e^{-(\cdot)^2/2}}{\sqrt{2\pi}}\right)}*U\|_{\mathscr C_{\lambda_{m-1}}}
&\le(1+C/N)\|U\|_{\mathscr C_{\lambda_m}},\label{eq:Pbound}\\
\|\alpha_m^{-1}\mathbb B_m^\zeta U-{\left(\frac{e^{-(\cdot)^2/2}}{\sqrt{2\pi}}\right)}*U\|_{\mathscr C_{\lambda_{m-1}}}
&\le\frac{Ca}{N}\|U\|_{\mathscr C_{\lambda_m}},\label{eq:Berror}\\
\|(\mathbb V_{m-1}-{\mathrm{Id}})U\|_{\mathscr C_{\lambda_{m-1}}}
&\le\frac{C}{N^2}\|U\|_{\mathscr C_{\lambda_{m-1}}}.\label{eq:Verror}
\end{align}
The constants depend on ${\bf b}_0,{\bf b}_1,\vartheta$, not on $a,N,t_0$.
\end{lem}

\begin{proof}
Write $z=w-u$. Since $\lambda_m-\lambda_{m-1}=\vartheta/N^2$,
\[
e^{\lambda_{m-1}\langle w\rangle-\lambda_m\langle u\rangle}
\le e^{C|z|/N}e^{-\vartheta\langle u\rangle/N^2}.
\]
Dropping the last factor and using
$\int {\frac{e^{-z^2/2}}{\sqrt{2\pi}}}e^{C|z|/N}\,dz\le1+C/N$ proves \eqref{eq:Pbound}.

To estimate the difference in \eqref{eq:Berror}, observe that
\[
\alpha_m^{-1}\mathbb B_m^\zeta U(w)
=\int {\frac{e^{-z^2/2}}{\sqrt{2\pi}}}\left[1+\frac{u+mz}{\zeta}\right]_+U(u)\,du.
\]
For every real $x$, $|[1+x]_+-1|\le|x|$. Multiplication by the weight ratio bounds the difference by the norm of $U$ times
\[
\int {\frac{e^{-z^2/2}}{\sqrt{2\pi}}}e^{C|z|/N}
\left(\frac{CN^2}{\zeta}+\frac{m|z|}{\zeta}\right)\,dz
\le\frac{Ca}{N},
\]
where we used $\sup_u|u|e^{-\vartheta\langle u\rangle/N^2}\le CN^2$ and $\zeta=N^3/a$. This argument includes the region where the positive part is zero; no Fourier-support assumption has been made.

Finally $\langle w\rangle\le\langle w-x\rangle+x$ for $x\ge0$. Apply \eqref{eq:moment} with $\lambda=\lambda_{m-1}$ to the convolution defining $\mathbb V_{m-1}-{\mathrm{Id}}$. Since $m-1\mathrel{{\sim}} N$, this proves \eqref{eq:Verror}.
\end{proof}

\begingroup
\begin{proof}[Proof of Proposition~\ref{thm:chain}]
First compare the actual top mode with ${\left(\frac{e^{-(\cdot)^2/4}}{\sqrt{4\pi}}\right)}$. Since $w_{m_0}\le-\zeta/2$, the Gaussian tail gives
\[
\|{\left(\frac{e^{-(\cdot)^2/4}}{\sqrt{4\pi}}\right)}\mathbf1_{(-\infty,w_{m_0})}\|_{\mathscr C_{\lambda_{m_0}}}
\le Ce^{-c\zeta^2}.
\]
Using \eqref{eq:moment} at $m_0$ in \eqref{eq:initialtrace}, we obtain
\begin{equation}\label{eq:toperror}
\|U_{m_0}-{\left(\frac{e^{-(\cdot)^2/4}}{\sqrt{4\pi}}\right)}\|_{\mathscr C_{\lambda_{m_0}}}
\le C N^{-2}+Ce^{-c\zeta^2},\qquad
\|U_{m_0}\|_{\mathscr C_{\lambda_{m_0}}}\le C.
\end{equation}

Set $Z_m=U_m/B_m$. Lemma~\ref{lem:step} gives, between the indicated successive weighted spaces,
\[
\begin{aligned}
\|\alpha_m^{-1}\mathbb V_{m-1}\mathbb B_m^\zeta\|&\le1+C/N,\\
\|\alpha_m^{-1}\mathbb V_{m-1}\mathbb B_m^\zeta U-{\left(\frac{e^{-(\cdot)^2/2}}{\sqrt{2\pi}}\right)}*U\|_{\mathscr C_{\lambda_{m-1}}}
&\le C(a/N+N^{-2})\|U\|_{\mathscr C_{\lambda_m}}.
\end{aligned}
\]

Since $m_0-m_1\le C_0N$,
\[
\begin{aligned}
\|Z_{m-1}\|_{\mathscr C_{\lambda_{m-1}}}&\le(1+C/N)\|Z_m\|_{\mathscr C_{\lambda_m}},\\
\|Z_m\|_{\mathscr C_{\lambda_m}}&\le(1+C/N)^{m_0-m}\|Z_{m_0}\|_{\mathscr C_{\lambda_{m_0}}}
\le e^{CC_0}\|Z_{m_0}\|_{\mathscr C_{\lambda_{m_0}}}\le C_2,\\
\|U_m\|_{\mathscr C_{\lambda_m}}&=B_m\|Z_m\|_{\mathscr C_{\lambda_m}}
\le C_2B_m,\qquad m_1\le m\le m_0.
\end{aligned}
\]

For $m_1\le m\le m_0$, the comparison profile is explicitly
\[
q_m(w)=\frac{1}{\sqrt{2\pi(m_0-m+2)}}
\exp\!\left(-\frac{w^2}{2(m_0-m+2)}\right).
\]
In particular, $q_{m_0}(w)={\frac{e^{-w^2/4}}{\sqrt{4\pi}}}$. Completing the square in the
convolution integral gives, for $m_1<m\le m_0$,
\[
\begin{aligned}
({\left(\frac{e^{-(\cdot)^2/2}}{\sqrt{2\pi}}\right)}*q_m)(w)
&=\frac{1}{2\pi\sqrt{m_0-m+2}}\int_{\mathbb R}
\exp\!\left(-\frac{(w-u)^2}{2}-\frac{u^2}{2(m_0-m+2)}\right)\,du\\
&=\frac{1}{\sqrt{2\pi(m_0-m+3)}}
\exp\!\left(-\frac{w^2}{2(m_0-m+3)}\right)=q_{m-1}(w).
\end{aligned}
\]
Each ideal transfer therefore increases the Gaussian variance by one.
\par

{Define} $e_m=\|Z_m-q_m\|_{\mathscr C_{\lambda_m}}$. By \eqref{eq:Pbound},
\begin{align*}
Z_{m-1}-q_{m-1}
&={\left(\frac{e^{-(\cdot)^2/2}}{\sqrt{2\pi}}\right)}*(Z_m-q_m)
+\alpha_m^{-1}\mathbb V_{m-1}\mathbb B_m^\zeta Z_m-{\left(\frac{e^{-(\cdot)^2/2}}{\sqrt{2\pi}}\right)}*Z_m,\\
e_{m-1}&\le(1+C/N)e_m+C(a/N+N^{-2})\|Z_m\|_{\mathscr C_{\lambda_m}}\\
&\le(1+C/N)e_m+C_3(a/N+N^{-2}),\\
e_m&\le(1+C/N)^{m_0-m}e_{m_0}
+C_3(a/N+N^{-2})\sum_{j=0}^{m_0-m-1}(1+C/N)^j\\
&\le e^{CC_0}e_{m_0}+C_3C_0Ne^{CC_0}(a/N+N^{-2}).
\end{align*}
Together with \eqref{eq:toperror}, this gives
\begin{equation}\label{eq:fullcomparison}
\|Z_m-q_m\|_{\mathscr C_{\lambda_m}}
\le C\left(a+N^{-1}+e^{-c\zeta^2}\right).
\end{equation}
There is no factor $e^{CN}$: the product of at most $CN$ factors $1+C/N$ is bounded independently of $N$.

Lemma~\ref{lem:identity} and \eqref{eq:fullcomparison} give
\eqref{eq:chaincomparison}. Since $m_0-m_1\sim N$, if $a\sqrt N\to0$,
then for all sufficiently large $N$,
\[
\frac{|\rho_{m_1}^\downarrow(T_1)|}{B_{m_1}}
\ge\frac1{\sqrt{2\pi(m_0-m_1+2)}}
-C\left(a+N^{-1}+e^{-c\zeta^2}\right)\ge\frac c{\sqrt N}.
\]
\par
Finally,
\[
\log B_{m_1}=\sum_{j=m_1+1}^{m_0}\log\frac{\Lambda}{(j/N)^3}
=N\int_{{\bf b}_1}^{{\bf b}_0}\log\frac{\Lambda}{y^3}\,dy+O(1),
\]
by the Riemann-sum estimate on a fixed interval away from zero. The integrand is positive on $[{\bf b}_1,{\bf b}_0)$.

To extend the chain, use the integral operator $\mathbb P_{k,m}$ from Lemma~\ref{lem:fullDensityIdentity}.
For $k=m_1-1,\ldots,1$ define successively
\[
\rho_k^\downarrow=\mathbb W_k\mathbb P_{k,k+1}\rho_{k+1}^\downarrow,
\qquad
\|\rho_k^\downarrow\|_{C(I)}
\le C\sup_{s\in I}\int_{t_0}^s|P_{k,k+1}(s,\tau)|\,d\tau
\|\rho_{k+1}^\downarrow\|_{C(I)}<\infty.
\]
Keep the already constructed modes $m_1\le k\le m_0$, put
$\rho_k^\downarrow=0$ for $k>m_0$, and reflect by complex conjugation.
Since the datum has only modes $\pm m_0$, the definitions give, for
every $k>0$,
\[
\rho_k^\downarrow=\mathbb W_k\mathfrak g_k
+\mathbb W_k\mathbb P_{k,k+1}\rho_{k+1}^\downarrow.
\]
Each retained mode is continuous. Hence, for every fixed $r$ and
$s_*\in I$,
\[
\|\rho^\downarrow(s)-\rho^\downarrow(s_*)\|_{H_x^r}^2
=4\pi\sum_{k=1}^{m_0}\langle k\rangle^{2r}
|\rho_k^\downarrow(s)-\rho_k^\downarrow(s_*)|^2
\longrightarrow0\qquad(s\to s_*).
\]
This proves the stated time continuity. The added equations
depend only on higher modes, so they preserve every previously proved
estimate on $[m_1,m_0]$.
\end{proof}
\par\endgroup
\subsection{\texorpdfstring{{Solving the approximate density equation}}{Solving the approximate density equation}}\label{sec:gainDensityInverse}

\begingroup
We solve the approximate density equation directly in the two-interval
space used for the collision comparison. The early equation has a
uniformly bounded inverse. Its contribution to the late equation is
then the known function $\mathfrak g^+$ in \eqref{eq:lateDensityInput}, and the descending gains give the late inverse
bound. Adding the upward-interaction term gives $\mathbb V_{\rm e}^0$.
Finally we estimate the single terminal coefficient without the factor
$N$ arising from summation over the spatial indices $k$.
\par
\par\endgroup
For the fixed mode parameters \eqref{eq:parameters}, interactions starting at modes $m>m_0$ add no further exponential amplification before
reaching $m_0$. The estimates below use this property to control
sources whose Fourier transforms in $x$ have unrestricted support.

Fix $\chi\in(1/2,2/3)$. {We retain the gain $B_k$ and the}
empty-product convention from \eqref{eq:alphagain}, and the maximum
$\widetilde B$ introduced in Section~\ref{sec:roadmapEcho}.
\begin{lem}\label{lem:extendedGain}
Assume \eqref{eq:parameters}. For all sufficiently large $N$, the sequence $B_k$ is nonincreasing, $B_{m_1}=B_*$, and
\begin{equation}\label{eq:maximumEchoGain}
1\le B_k\le\widetilde{B}\le e^{C_{\rm B}N}.
\end{equation}
\begin{equation}\label{eq:maximumEchoGainAsymptotic}
{
\log \widetilde{B} =3(m_0-1)\log({\bf b}_0N)-3\log(m_0!) =3{\bf b}_0N-\frac92\log N+O(1).
}
\end{equation}
One may take $C_{\rm B}=3{\bf b}_0+1$, independently of all fixed derivative and expansion orders; also $c_{\rm e}<3{\bf b}_0$.
\end{lem}

\begin{proof}
For $k<m_0$, $B_k/B_{k+1}=\alpha_{k+1}\ge1$ by \eqref{eq:parameters}. For $k\ge m_0$ the extended gain is one. Hence the maximum is attained at $k=1$, and
\[
\widetilde{B}=\prod_{j=2}^{m_0}\alpha_j=\frac{(\Lambda N^3)^{m_0-1}}{(m_0!)^3}.
\]
Stirling's formula and $0\le {\bf b}_0N-m_0<1$ give
\[
\begin{aligned}
\log\widetilde B
&=3(m_0-1)\log({\bf b}_0N)
-3(m_0+\tfrac12)\log m_0+3m_0+O(1)\\
&=3m_0\log\frac{{\bf b}_0N}{m_0}+3m_0
-3\log({\bf b}_0N)-\tfrac32\log m_0+O(1),\\
m_0\log\frac{{\bf b}_0N}{m_0}
&={\bf b}_0N-m_0+O(N^{-1}).
\end{aligned}
\]
Substitution proves \eqref{eq:maximumEchoGainAsymptotic}, including the
bounded effect of the floor, and hence \eqref{eq:maximumEchoGain} with
$C_{\rm B}=3{\bf b}_0+1$. Finally
\[
0<c_{\rm e}=\int_{{\bf b}_1}^{{\bf b}_0}3\log({\bf b}_0/y)\,dy
<\int_0^{{\bf b}_0}3\log({\bf b}_0/y)\,dy=3{\bf b}_0.
\]
The chosen $C_{\rm B}$ depends only on the fixed mode parameters in \eqref{eq:parameters}.
\end{proof}

\begingroup
We estimate the density contributions in \eqref{eq:allmodeidentity}, with
the function $\mathfrak g^0$ defined in
\eqref{eq:overviewApproximateDensityKernels}.
\par
\par\endgroup
Set
\begin{equation}\label{eq:electricTimeSplit}
T_{\rm a}=t_0+N^4,\qquad I_-=[t_0,T_{\rm a}],\qquad I_+=[T_{\rm a},T_1].
\end{equation}
Choose $\vartheta>0$ sufficiently small and define the phase weights
\begin{equation}\label{eq:echoPhaseWeights}
\lambda_N(s)=\vartheta N^{-2\chi}\left(2-\frac{s-t_0}{T_1-t_0}\right),
\qquad w_{k,N}(s)=e^{\lambda_N(s)\langle ks-\zeta\rangle^\chi}.
\end{equation}
For $I=[t_0,T_1]$ and $r_{\rm x}\ge0$, use the single density norm
\begin{equation}\label{eq:echoDensityNorm}
\|\rho\|_{\mathscr X_N^{r_{\rm x}}(I)}
=\sum_{k\ne0}\langle k/N\rangle^{r_{\rm x}}
\left(\sup_{I_-}w_{|k|,N}|\rho_k|
+B_{|k|}^{-1}\sup_{I_+}w_{|k|,N}|\rho_k|\right).
\end{equation}

\begingroup
The space consists of continuous mean-zero functions on $I$ with
finite norm \eqref{eq:echoDensityNorm}. The positive and negative
spatial modes are preserved separately by the approximate equation;
the two blocks have the same kernels after replacing $(k,m)$ by
$(|k|,|m|)$. We estimate the positive block and sum the two bounds.
No reality assumption is needed, and the zero mode generates no field.

\begin{rem}[The weights and the division of the time interval]\label{rem:densitySpaceMotivation}
The weight in \eqref{eq:weights} controls the density away from
$ms=\zeta$. For the comparison with the exact collision equation we
use the same resonant variable and a decreasing Gevrey radius:
\[
\frac{w_{k,N}(s)}{w_{m,N}(\tau)}
\le\exp\left(\lambda_N(s)|ks-m\tau|^\chi
-\frac{\vartheta N^{-2\chi}(s-\tau)}{T_1-t_0}
\langle m\tau-\zeta\rangle^\chi\right).
\]
The decrease absorbs polynomial factors in $m\tau-\zeta$.
The scale $N^{-2\chi}$ is compatible with the coefficient estimate
of Proposition~\ref{prop:decreasingRadiusBackground}: with
$n=k-m-l$, the identity
\[
ks-m\tau=k(s-t_0)-m(\tau-t_0)+(l+n)t_0
\]
shows that the weight costs factors with exponents
$\lambda_N(s)|k(s-t_0)-m(\tau-t_0)|^\chi$ and
$\lambda_N(s)t_0^\chi(|l|^\chi+|n|^\chi)$.
Since $t_0=A_0N$, $\lambda_N=O(N^{-2\chi})$ gives
$\lambda_Nt_0^\chi=O(N^{-\chi})$, which fits the coefficient Fourier decay.

The early endpoint $T_{\rm a}=t_0+N^4$ precedes the prepared
resonances. For the kernel in \eqref{eq:overviewApproximateDensityKernels},
$z=ks-m\tau$ and $|k-m|=1$ give
$|k|(s-\tau)\le\tau+|z|$. Gaussian integration therefore gives
\[
\sup_{s\in I_-}\int_{t_0}^s |P_{k,m}(s,\tau)|
\frac{w_{|k|,N}(s)}{w_{|m|,N}(\tau)}\,d\tau
\le\frac{Ca(1+T_{\rm a})}{m^2},\qquad a(1+T_{\rm a})=o(1).
\]
No gain denominator is needed on this interval. At later times the
Gaussian does not suppress the descending resonances $ks=m\tau$;
the mode-dependent growth is normalized by
\[
\alpha_m\frac{B_m}{B_{m-1}}=1\qquad(2\le m\le m_0).
\]
For the nonresonant configurations of the exact coefficient kernels,
$\tau-t_0\ge N^4$ supplies phase separation. The decay at frequency
width $N$ then bounds the large gain ratios by
$e^{CN-cN^{3\chi}}=o(1)$, since $3\chi>1$.
Thus the time split permits direct smallness before the resonances,
gain normalization on the descending chain, and exponential decay
for the other late-time configurations.
\end{rem}
\par\endgroup

\begingroup
\begin{prop}\label{prop:approximateDensityInverse}
Assume \eqref{eq:mainparameters} and \eqref{eq:parameters}, fix
$r_{\rm x}\ge0$, $1/2<\chi<2/3$, and choose the constant in
\eqref{eq:echoPhaseWeights} sufficiently small.
The descending equation
$\rho=\mathfrak g+(\mathbb K+\mathbb P_\downarrow)\rho$
has a unique solution for every $\mathfrak g\in\mathscr X_N^{r_{\rm x}}(I)$, and
\begin{equation}\label{eq:timeGdowninverse}
\mathbb V_\downarrow=(\mathrm{Id}-\mathbb K-\mathbb P_\downarrow)^{-1},\qquad
\|\mathbb V_\downarrow\|_{\mathscr X_N^{r_{\rm x}}(I)\to\mathscr X_N^{r_{\rm x}}(I)}
\le C_{r_{\rm x}}N.
\end{equation}
The equation including the upward term has the inverse
\begin{equation}\label{eq:referenceElectricInverse}
\mathbb V_{\rm e}^0
=(\mathrm{Id}-\mathbb K-\mathbb P_\downarrow-\mathbb P_\uparrow)^{-1}
=(\mathrm{Id}-\mathbb V_\downarrow\mathbb P_\uparrow)^{-1}\mathbb V_\downarrow,\qquad \|\mathbb V_{\rm e}^0\|_{\mathscr X_N^{r_{\rm x}}(I)\to\mathscr X_N^{r_{\rm x}}(I)}
\le C_{r_{\rm x}}N.
\end{equation}
Both restricted equations on $I_-$ satisfy the sharper estimate
\begin{equation}\label{eq:earlyApproximateInverse}
\sum_{k\ne0}\langle k/N\rangle^{r_{\rm x}}\sup_{I_-}w_{|k|,N}|\rho_k|
\le C_{r_{\rm x}}\sum_{k\ne0}\langle k/N\rangle^{r_{\rm x}}\sup_{I_-}w_{|k|,N}|\mathfrak g_k|.
\end{equation}
For the descending density extended to mode one,
\begin{equation}\label{eq:descendingTwoIntervalBounds}
\|\mathfrak g^0\|_{\mathscr X_N^{r_{\rm x}}}\le C_{r_{\rm x}},
\qquad \|\rho^\downarrow\|_{\mathscr X_N^{r_{\rm x}}}\le C_{r_{\rm x}}N,
\end{equation}
and
\begin{equation}\label{eq:descendingReferenceDefect}
(\mathrm{Id}-\mathbb K-\mathbb P_\downarrow-\mathbb P_\uparrow)\rho^\downarrow-\mathfrak g^0
=-\mathbb P_\uparrow\rho^\downarrow,\qquad \|\mathbb P_\uparrow\rho^\downarrow\|_{\mathscr X_N^{r_{\rm x}}}\le e^{-cN^2}.
\end{equation}
With
\[
\|\rho\|_{\ell^1_{r_{\rm x}}L^\infty(I)}
=\sum_{k\ne0}\langle k\rangle^{r_{\rm x}}\|\rho_k\|_{L^\infty(I)},
\]
the norm satisfies
\begin{equation}\label{eq:echoDensityEmbeddings}
\|\rho\|_{\ell^1_{r_{\rm x}}L^\infty(I)}
\le C_{r_{\rm x}}N^{r_{\rm x}}\widetilde B\|\rho\|_{\mathscr X_N^{r_{\rm x}}(I)}.
\end{equation}
For large $N$, $t_0<T_{\rm a}<T_1$, $m_0T_{\rm a}/\zeta\to0$ and $a(1+T_{\rm a})\to0$.
\end{prop}
\par\endgroup

\begingroup
\begin{proof}
\emph{1. Phase weights and the scalar inverse.}
The kernels are those in \eqref{eq:overviewApproximateDensityKernels}.
On positive modes the equation reads
\begin{equation}\label{eq:referenceDensityResponse}
\rho_k=\mathfrak g_k+\mathbb K_k\rho_k+\mathbb P_{k,k+1}\rho_{k+1}
+\mathbf1_{k\ge2}\mathbb P_{k,k-1}\rho_{k-1}.
\end{equation}
Subadditivity and the decreasing radius give, for $s\ge\tau$,
\begin{equation}\label{eq:timeGweightratio}
\frac{w_{k,N}(s)}{w_{m,N}(\tau)}
\le\exp\left(\lambda_N(s)|ks-m\tau|^\chi
-\frac{\vartheta N^{-2\chi}(s-\tau)}{T_1-t_0}\langle m\tau-\zeta\rangle^\chi\right).
\end{equation}
Since $\big||k|s-|m|\tau\big|\le|ks-m\tau|$, this also implies
\begin{equation}\label{eq:oscDweightratio}
\frac{w_{|k|,N}(s)}{w_{|m|,N}(\tau)}\le e^{\lambda_N(s)|ks-m\tau|^\chi},\qquad k,m\ne0.
\end{equation}
On either interval $I_\pm$, the signed Maxwellian resolvent of
Lemma~\ref{prop:inverse} defines $\mathbb W_k=(\mathrm{Id}-\mathbb K_k)^{-1}$ and gives
\begin{equation}\label{eq:timeGinverse}
\|w_{k,N}(\mathbb W_k-\mathrm{Id})\mathfrak g_k\|_{L^\infty(I_\pm)}
\le\frac C{k^2}\|w_{k,N}\mathfrak g_k\|_{L^\infty(I_\pm)}.
\end{equation}
Indeed, the weighted convolution norm is bounded by
$\int_0^\infty|\widetilde R_k(x)|e^{2\vartheta N^{-2\chi}x^\chi}\,dx\le C/k^2$.
This includes $k=1$ without an $L^1$ contraction assumption.

\emph{2. The early equation.}
The parameter relations give
$T_{\rm a}/T_1+m_0T_{\rm a}/\zeta=O(aN^2)$ and
$a(1+T_{\rm a})=O(aN^4)=o(1)$.
For $|k-m|=1$, put $z=ks-m\tau$. Since
$|k|(s-\tau)\le\tau+|z|$, the Gaussian kernel gives
\begin{equation}\label{eq:earlyAdjacentKernel}
\sup_{s\in I_-}\int_{t_0}^s|P_{k,m}(s,\tau)|
\frac{w_{|k|,N}(s)}{w_{|m|,N}(\tau)}\,d\tau
\le\frac{Ca(1+T_{\rm a})}{|m|^2}.
\end{equation}
For $|k-m|=1$, $\langle k/N\rangle^{r_{\rm x}}\le
C_{r_{\rm x}}\langle m/N\rangle^{r_{\rm x}}$. Thus
\eqref{eq:timeGinverse} and \eqref{eq:earlyAdjacentKernel} give
\[
\begin{aligned}
&\sum_{k\ne0}\langle k/N\rangle^{r_{\rm x}}
\sup_{I_-}w_{|k|,N}
\left|\big[\mathbb W(\mathbb P_\downarrow+\mathbb P_\uparrow)\rho\big]_k\right|\\
&\quad\le C_{r_{\rm x}}a(1+T_{\rm a})
\sum_{m\ne0}\langle m/N\rangle^{r_{\rm x}}
\sup_{I_-}w_{|m|,N}|\rho_m|.
\end{aligned}
\]
For $N$ sufficiently large, $C_{r_{\rm x}}a(1+T_{\rm a})<1/2$.
Consequently the series
\[
\rho=\sum_{n=0}^\infty
\big[\mathbb W(\mathbb P_\downarrow+\mathbb P_\uparrow)\big]^n
\mathbb W\mathfrak g
\]
converges in the early norm and proves \eqref{eq:earlyApproximateInverse}.
Omitting $\mathbb P_\uparrow$ gives the descending early inverse with
the same bound.
\par

\begingroup
\emph{3. Descending transfers on the late interval.}
The time integrals in this step start at $T_{\rm a}$. We first estimate
$\mathbb P_{m-1,m}$ and then apply
$\mathbb W_{m-1}=(\mathrm{Id}-\mathbb K_{m-1})^{-1}$.
For $\mathbb A_m=\mathbb W_{m-1}\mathbb P_{m-1,m}$, the required bound is
\[
\frac{\langle(m-1)/N\rangle^{r_{\rm x}}B_m}
{\langle m/N\rangle^{r_{\rm x}}B_{m-1}}
\frac{\|w_{m-1,N}\mathbb A_m\rho_m\|_{L^\infty(I_+)}}
{\|w_{m,N}\rho_m\|_{L^\infty(I_+)}}
\le\begin{cases}(1+C/(m-1)^2)(1+C/N),&2\le m\le C_0N,\\1/2,&m>C_0N.
\end{cases}
\]
Here $C_0>{\bf b}_0$ will be fixed sufficiently large. When $\rho_m=0$, one has $\mathbb A_m\rho_m=0$; hence the displayed quotient is used only when the
denominator is positive.

For $m\ge2$ and $s\in I_+$, the defining time integral is
\[
[\mathbb P_{m-1,m}\rho_m](s)
=-\frac{a(m-1)}{2m}\int_{T_{\rm a}}^s
(s-\tau)e^{-((m-1)s-m\tau)^2/2}\rho_m(\tau)\,d\tau.
\]
For $2\le m\le m_0$, positivity of $w_{m,N}$ and
$B_{m-1}=\alpha_mB_m$, where
$\alpha_m=a\sqrt{2\pi}\zeta/(2m^3)=\Lambda N^3/m^3$, give
\begin{align*}
&\frac{w_{m-1,N}(s)}{B_{m-1}}
|[\mathbb P_{m-1,m}\rho_m](s)|\\
&\quad\le\frac{w_{m-1,N}(s)}{B_{m-1}}
\int_{T_{\rm a}}^s|P_{m-1,m}(s,\tau)|
\frac{w_{m,N}(\tau)|\rho_m(\tau)|}{w_{m,N}(\tau)}\,d\tau\\
&\quad\le\left(\frac1{B_m}\sup_{\tau\in I_+}
w_{m,N}(\tau)|\rho_m(\tau)|\right)
\int_{T_{\rm a}}^s\frac{|P_{m-1,m}(s,\tau)|}{\alpha_m}
\frac{w_{m-1,N}(s)}{w_{m,N}(\tau)}\,d\tau.
\end{align*}
We now estimate the last integral, uniformly in $s$, for
$2\le m\le C_0N$. This integral estimate itself does not require $m\le m_0$.
Set $k=m-1$ and make the change of variables
\[
u=m\tau-\zeta,\qquad z=ks-m\tau=ks-\zeta-u,\qquad
d\tau=\frac{du}{m},\qquad mT_{\rm a}-\zeta\le u\le ms-\zeta.
\]
The identities
\[
k(s-\tau)=\tau+z,\qquad mk(s-\tau)=\zeta+u+mz,
\qquad \frac{|P_{m-1,m}(s,\tau)|}{\alpha_m}\,d\tau
=\frac{e^{-z^2/2}}{\sqrt{2\pi}}
\left(1+\frac{u+mz}{\zeta}\right)du
\]
hold on this interval, where the bracket is nonnegative. Hence
\begin{equation}\label{eq:timeGprincipal}
\begin{aligned}
&\int_{T_{\rm a}}^s\frac{|P_{m-1,m}(s,\tau)|}{\alpha_m}
\frac{w_{m-1,N}(s)}{w_{m,N}(\tau)}\,d\tau\\
&\quad=\int_{mT_{\rm a}-\zeta}^{ms-\zeta}
\frac{e^{-z^2/2}}{\sqrt{2\pi}}
\left[1+\frac{u+mz}{\zeta}\right]_+
\frac{w_{m-1,N}(s)}{w_{m,N}((u+\zeta)/m)}\,du.
\end{aligned}
\end{equation}
In the integrals below, $z=ks-\zeta-u$ and $\tau=(u+\zeta)/m$.
The weight estimate \eqref{eq:timeGweightratio} reads
\[
\frac{w_{m-1,N}(s)}{w_{m,N}((u+\zeta)/m)}
\le\exp\left(\lambda_N(s)|z|^\chi
-\frac{\vartheta N^{-2\chi}(s-\tau)}{T_1-t_0}\langle u\rangle^\chi\right).
\]

On the part with $u<0$, the first integral is bounded by
\[
\begin{aligned}
&\int_{\substack{mT_{\rm a}-\zeta\le u\le ms-\zeta\\u<0}}
\frac{e^{-z^2/2}}{\sqrt{2\pi}}
\left[1+\frac{u+mz}{\zeta}\right]_+
\frac{w_{m-1,N}(s)}{w_{m,N}((u+\zeta)/m)}\,du\\
&\quad\le\int_{\substack{mT_{\rm a}-\zeta\le u\le ms-\zeta\\u<0}}
\frac{e^{-z^2/2+\lambda_N(s)|z|^\chi}}{\sqrt{2\pi}}
\left(1+\frac{m|z|}{\zeta}\right)du.
\end{aligned}
\]
If $u\ge0$ and $z\ge-\tau/2$, then
\[
\tau\ge\frac\zeta m,\qquad
s-\tau=\frac{\tau+z}{k}\ge\frac\zeta{2mk},\qquad
\frac{\vartheta N^{-2\chi}(s-\tau)}{T_1-t_0}
\ge\frac{\vartheta m_1N^{-2\chi}}{2mk}
\ge cN^{-1-2\chi}.
\]
The part of the bracket containing $u$ therefore satisfies
\begin{equation}\label{eq:timeGfrequencyabsorb}
\sup_{u\ge0}\frac u\zeta e^{-cN^{-1-2\chi}\langle u\rangle^\chi}
\le\frac{(c\chi e)^{-1/\chi}N^{2+1/\chi}}{\zeta}
\le CaN^{1/\chi-1}.
\end{equation}
Consequently the integral on this second region obeys
\[
\begin{aligned}
&\int_{\substack{mT_{\rm a}-\zeta\le u\le ms-\zeta\\u\ge0,\ z\ge-\tau/2}}
\frac{e^{-z^2/2}}{\sqrt{2\pi}}
\left[1+\frac{u+mz}{\zeta}\right]_+
\frac{w_{m-1,N}(s)}{w_{m,N}((u+\zeta)/m)}\,du\\
&\quad\le\int_{\substack{mT_{\rm a}-\zeta\le u\le ms-\zeta\\u\ge0,\ z\ge-\tau/2}}
\frac{e^{-z^2/2+\lambda_N(s)|z|^\chi}}{\sqrt{2\pi}}
\left(1+\frac{m|z|}{\zeta}+CaN^{1/\chi-1}\right)du.
\end{aligned}
\]
On the remaining region $u\ge0$, $z<-\tau/2$, we have
\[
|z|\ge\frac\tau2\ge\frac\zeta{2m}\ge cT_1,\qquad
0\le1+\frac{u+mz}{\zeta}
=\frac{m(\tau+z)}{\zeta}\le\frac{mT_1}{2\zeta}\le C.
\]
Since $dz=-du$ for fixed $s$, the third integral satisfies
\[
\begin{aligned}
&\int_{\substack{mT_{\rm a}-\zeta\le u\le ms-\zeta\\u\ge0,\ z<-\tau/2}}
\frac{e^{-z^2/2}}{\sqrt{2\pi}}
\left[1+\frac{u+mz}{\zeta}\right]_+
\frac{w_{m-1,N}(s)}{w_{m,N}((u+\zeta)/m)}\,du\\
&\quad\le C\int_{|z|\ge cT_1}e^{-z^2/2+\lambda_N(s)|z|^\chi}\,dz
\le Ce^{-c'T_1^2}.
\end{aligned}
\]

The first two regions are disjoint, and the map $u\mapsto z$ has
Jacobian of absolute value one. Their leading integrals therefore
combine into a single Gaussian integral. Using
$0<\lambda_N(s)\le2\vartheta N^{-2\chi}$, we obtain
\[
\begin{aligned}
\int_{\mathbb R}\frac{e^{-z^2/2+\lambda_N(s)|z|^\chi}}{\sqrt{2\pi}}\,dz
&=1+\int_{\mathbb R}\frac{e^{-z^2/2}}{\sqrt{2\pi}}
\big(e^{\lambda_N(s)|z|^\chi}-1\big)\,dz\\
&\le1+2\vartheta N^{-2\chi}\int_{\mathbb R}
\frac{|z|^\chi e^{-z^2/2+2\vartheta|z|^\chi}}{\sqrt{2\pi}}\,dz
\le1+CN^{-2\chi},\\
\int_{\mathbb R}\frac{|z|e^{-z^2/2+\lambda_N(s)|z|^\chi}}{\sqrt{2\pi}}\,dz
&\le\int_{\mathbb R}\frac{|z|e^{-z^2/2+2\vartheta|z|^\chi}}{\sqrt{2\pi}}\,dz\le C.
\end{aligned}
\]
Adding the three region estimates in \eqref{eq:timeGprincipal} gives
\[
\begin{aligned}
&\int_{T_{\rm a}}^s\frac{|P_{m-1,m}(s,\tau)|}{\alpha_m}
\frac{w_{m-1,N}(s)}{w_{m,N}(\tau)}\,d\tau\\
&\quad\le\int_{\mathbb R}\frac{e^{-z^2/2+\lambda_N(s)|z|^\chi}}{\sqrt{2\pi}}
\left(1+\frac{m|z|}{\zeta}+CaN^{1/\chi-1}\right)dz+Ce^{-c'T_1^2}\\
&\quad\le1+CN^{-2\chi}+C\frac m\zeta+CaN^{1/\chi-1}+Ce^{-c'T_1^2}
\le1+\frac CN.
\end{aligned}
\]
Here $2\chi>1$ and the parameter conditions give
\[
\frac m\zeta=\frac{am}{N^3}\le\frac C{N^2},\qquad
aN^{1/\chi}\le N^{-{\rm s}_0-2+1/\chi}\le1,\qquad
e^{-c'T_1^2}\le N^{-1}
\]
for sufficiently large $N$. Extracting $\sup_{\tau\in I_+}w_{m,N}(\tau)|\rho_m(\tau)|$, with or
without the gain factors, now proves
\begin{equation}\label{eq:timeGdownzero}
\|w_{m-1,N}\mathbb P_{m-1,m}\rho_m\|_{L^\infty(I_+)}
\le\alpha_m(1+C/N)\|w_{m,N}\rho_m\|_{L^\infty(I_+)},
\qquad 2\le m\le C_0N.
\end{equation}

For $m>C_0N$, return to the original time integral. With $k=m-1$,
$z=ks-m\tau$, $d\tau=-dz/m$ and $k(s-\tau)=\tau+z$, we have
\begin{equation}\label{eq:timeGhighmodes}
\begin{aligned}
&\frac{\|w_{m-1,N}\mathbb P_{m-1,m}\rho_m\|_{L^\infty(I_+)}}
{\|w_{m,N}\rho_m\|_{L^\infty(I_+)}} \le\sup_{s\in I_+}\frac{a}{2m^2}
\int_{-s}^{ks-mT_{\rm a}}(\tau+z)
e^{-z^2/2+\lambda_N(s)|z|^\chi}\,dz\\
&\quad\le\frac{a}{2m^2}\int_{\mathbb R}(T_1+|z|)
e^{-z^2/2+2\vartheta N^{-2\chi}|z|^\chi}\,dz
\le\frac{CN^2}{m^2},
\end{aligned}
\end{equation}
where $\tau=(ks-z)/m$ in the first integral and $aT_1\sim N^2$.

Finally, \eqref{eq:timeGinverse} gives
\[
\begin{aligned}
\|w_{m-1,N}\mathbb A_m\rho_m\|_{L^\infty(I_+)}
&=\|w_{m-1,N}\mathbb W_{m-1}\mathbb P_{m-1,m}\rho_m\|_{L^\infty(I_+)}\\
&\le\left(1+\frac{C}{(m-1)^2}\right)
\|w_{m-1,N}\mathbb P_{m-1,m}\rho_m\|_{L^\infty(I_+)}.
\end{aligned}
\]
The remaining weight ratios satisfy
\[
\frac{\langle(m-1)/N\rangle^{r_{\rm x}}}{\langle m/N\rangle^{r_{\rm x}}}\le1,
\qquad
\alpha_m\frac{B_m}{B_{m-1}}
=\begin{cases}1,&2\le m\le m_0,\\\alpha_m\le1,&m>m_0.\end{cases}
\]
Thus \eqref{eq:timeGdownzero} proves the required bound for
$2\le m\le C_0N$. For $m>C_0N>m_0$, $B_m/B_{m-1}=1$ and
\[
\left(1+\frac{C}{(m-1)^2}\right)\frac{CN^2}{m^2}\le\frac12
\]
once $C_0$ is sufficiently large. This proves the high-mode bound.
\par\endgroup

For $1\le k\le m$, the low-mode factors satisfy
\[
\prod_{j=k+1}^{\min\{m,\lfloor C_0N\rfloor\}}
\left(1+\frac{C}{(j-1)^2}\right)\left(1+\frac CN\right)
\le e^{CC_0}\prod_{j=2}^\infty\left(1+\frac{C}{(j-1)^2}\right)\le C.
\]
Each remaining factor is at most $1/2$. Combining this product with
the bound for $\mathbb W_m$ gives
\begin{equation}\label{eq:approximateDescendingPathBound}
\begin{aligned}
&\frac{\langle k/N\rangle^{r_{\rm x}}}{B_k}
\|w_{k,N}\mathbb A_{k+1}\cdots\mathbb A_m
\mathbb W_m\mathfrak g_m\|_{L^\infty(I_+)}\\
&\quad\le C_{r_{\rm x}}
2^{-(m-\max\{k,\lfloor C_0N\rfloor\})_+}
\frac{\langle m/N\rangle^{r_{\rm x}}}{B_m}
\|w_{m,N}\mathfrak g_m\|_{L^\infty(I_+)}.
\end{aligned}
\end{equation}
Empty products of operators are the identity. Summation over $k$ gives the column bound
\[
\sup_{m\ge1}\sum_{k=1}^m
2^{-(m-\max\{k,\lfloor C_0N\rfloor\})_+}
\le C_0N+\sum_{h=0}^\infty2^{-h}\le CN.
\]

For a continuous density $\mathfrak g$ with finite late norm, define
\[
\rho_k^{(M)}=\sum_{m=k}^{M}
\mathbb A_{k+1}\cdots\mathbb A_m\mathbb W_m\mathfrak g_m
\quad(1\le k\le M),\qquad \rho_k^{(M)}=0\quad(k>M).
\]
Cancellation of the common terms gives, for every $k\ge1$,
\[
\rho_k^{(M)}-\mathbb A_{k+1}\rho_{k+1}^{(M)}
=\mathbb W_k\mathfrak g_k\mathbf1_{k\le M},\qquad (\mathrm{Id}-\mathbb K_k)\rho_k^{(M)}
-\mathbb P_{k,k+1}\rho_{k+1}^{(M)}
=\mathfrak g_k\mathbf1_{k\le M}.
\]
The path and column bounds imply
\[
\sum_{k\ge1}\frac{\langle k/N\rangle^{r_{\rm x}}}{B_k}
\sup_{I_+}w_{k,N}|\rho_k^{(M)}|
\le C_{r_{\rm x}}N
\sum_{m=1}^{M}\frac{\langle m/N\rangle^{r_{\rm x}}}{B_m}
\sup_{I_+}w_{m,N}|\mathfrak g_m|.
\]
For $M>L$, the difference satisfies
\[
\begin{aligned}
&\sum_{k\ge1}\frac{\langle k/N\rangle^{r_{\rm x}}}{B_k}
\sup_{I_+}w_{k,N}|\rho_k^{(M)}-\rho_k^{(L)}|\\
&\quad\le C_{r_{\rm x}}N
\sum_{m=L+1}^{M}\frac{\langle m/N\rangle^{r_{\rm x}}}{B_m}
\sup_{I_+}w_{m,N}|\mathfrak g_m|\longrightarrow0
\qquad(L\to\infty).
\end{aligned}
\]
Thus $\rho^{(M)}$ converges in the late norm to a continuous density
$\rho$. To pass to the equation, note that
\[
\int_0^\infty |D_k(x)|e^{2\vartheta N^{-2\chi}x^\chi}\,dx
=\frac1{k^2}\int_0^\infty x e^{-x^2/2+2\vartheta N^{-2\chi}x^\chi}\,dx
\le\frac C{k^2}.
\]
This bounds $\mathbb K$ in the late norm. The one-step bounds control
$\mathbb W\mathbb P_\downarrow$, and
$\mathbb P_\downarrow=(\mathrm{Id}-\mathbb K)\mathbb W\mathbb P_\downarrow$
is consequently bounded there. Passing to the limit in the finite
equations gives $\rho=\mathfrak g+\mathbb K\rho+\mathbb P_\downarrow\rho$ and
\begin{equation}\label{eq:lateDescendingInverse}
\sum_{k\ge1}\frac{\langle k/N\rangle^{r_{\rm x}}}{B_k}\sup_{I_+}w_{k,N}|\rho_k|
\le C_{r_{\rm x}}N\sum_{m\ge1}\frac{\langle m/N\rangle^{r_{\rm x}}}{B_m}\sup_{I_+}w_{m,N}|\mathfrak g_m|.
\end{equation}
For the other inverse identity, let $\rho$ be any continuous density
supported in $1\le m\le M$. Then
\[
\begin{aligned}
&\sum_{m=k}^{M}\mathbb A_{k+1}\cdots\mathbb A_m\mathbb W_m
\big[(\mathrm{Id}-\mathbb K_m)\rho_m-\mathbb P_{m,m+1}\rho_{m+1}\big]\\
&\quad=\sum_{m=k}^{M}\mathbb A_{k+1}\cdots\mathbb A_m
\big[\rho_m-\mathbb A_{m+1}\rho_{m+1}\big] =\rho_k-\mathbb A_{k+1}\cdots\mathbb A_{M+1}\rho_{M+1}
=\rho_k.
\end{aligned}
\]
Finite-mode truncations converge in the late norm. The boundedness
of the solution map, $\mathbb K$ and $\mathbb P_\downarrow$ therefore
extends this identity to every density in that space. This proves
both inverse identities and uniqueness.
\par

\emph{4. Joining the two intervals.}
First solve the descending early equation. For $s\in I_+$ the known function in the descending late equation is
\begin{equation}\label{eq:lateDensityInput}
\mathfrak g_k^+(s)=\mathfrak g_k(s)
+\int_{t_0}^{T_{\rm a}}K_k(s-\tau)\rho_k(\tau)\,d\tau
+\int_{t_0}^{T_{\rm a}}P_{k,k+1}(s,\tau)\rho_{k+1}(\tau)\,d\tau.
\end{equation}
The early-to-late bounds are
\[
\begin{aligned}
\sup_{s\in I_+}\frac1{B_k}\int_{t_0}^{T_{\rm a}}
|K_k(s-\tau)|\frac{w_{k,N}(s)}{w_{k,N}(\tau)}\,d\tau&\le C/k^2,\\
\sup_{s\in I_+}\frac1{B_k}\int_{t_0}^{T_{\rm a}}
|P_{k,m}(s,\tau)|\frac{w_{k,N}(s)}{w_{m,N}(\tau)}\,d\tau
&\le\frac{Ca(1+T_{\rm a})}{m^2},\quad m=k+1.
\end{aligned}
\]
For the second inequality use $k(s-\tau)=\tau+z$,
$d\tau=-dz/m$, and $B_k\ge1$; the first follows from the scalar
Gaussian kernel. Therefore
\begin{equation}\label{eq:lateDensityInputBound}
\sum_{k\ge1}\frac{\langle k/N\rangle^{r_{\rm x}}}{B_k}
\sup_{I_+}w_{k,N}|\mathfrak g_k^+|
\le C_{r_{\rm x}}\|\mathfrak g\|_{\mathscr X_N^{r_{\rm x}}}.
\end{equation}
Apply \eqref{eq:lateDescendingInverse} to $\mathfrak g^+$ and add the early
bound. At $T_{\rm a}$ the late integrals vanish and
$\mathfrak g_k^+(T_{\rm a})=\rho_k(T_{\rm a})$ by the early equation,
so the solutions join continuously. This constructs the descending
inverse and proves \eqref{eq:timeGdowninverse}.

\emph{5. The upward-interaction term and the chain defect.}
For $k=m+1$, time ordering implies
$z=ks-m\tau=s+m(s-\tau)\ge s\ge t_0$ and $k(s-\tau)\le z$.
Thus on every allowed time block,
\begin{equation}\label{eq:referenceUpwardKernel}
\sup_s\int|P_{m+1,m}(s,\tau)|\frac{w_{m+1,N}(s)}{w_{m,N}(\tau)}\,d\tau
\le\frac{Ca}{m^2}\int_{t_0}^\infty z e^{-z^2/2+2\vartheta N^{-2\chi}z^\chi}\,dz
\le\frac{Ca}{m^2}e^{-ct_0^2}.
\end{equation}
The gain ratios on the early, crossing and late blocks are respectively
$1$, $B_{m+1}^{-1}\le1$ and $B_m/B_{m+1}\le CN^3$.
Hence
\begin{equation}\label{eq:timeGupsmall}
\|\mathbb P_\uparrow\|_{\mathscr X_N^{r_{\rm x}}\to\mathscr X_N^{r_{\rm x}}}
\le C_{r_{\rm x}}aN^3e^{-cN^2},\qquad \|\mathbb V_\downarrow\mathbb P_\uparrow\|
\le C_{r_{\rm x}}N\|\mathbb P_\uparrow\|
\le C_{r_{\rm x}}aN^4e^{-cN^2}<1/2.
\end{equation}
The early, crossing and late bounds also show that $\mathbb K$ and
$\mathbb P_\downarrow$ are bounded on the full two-interval space.
\par

\begingroup
The factorization follows from
\[
(\mathrm{Id}-\mathbb K-\mathbb P_\downarrow)(\mathrm{Id}-\mathbb V_\downarrow\mathbb P_\uparrow)
=\mathrm{Id}-\mathbb K-\mathbb P_\downarrow-\mathbb P_\uparrow.
\]
The preceding bound on $\mathbb V_\downarrow\mathbb P_\uparrow$ proves
\eqref{eq:referenceElectricInverse}.
\par\endgroup
The negative block obeys identical bounds, proving the two-sided claims.
The function $\mathfrak g^0$ on positive modes is
\begin{equation}\label{eq:collisionalzerokernels}
\mathfrak g_k^0(s)=\mathbf1_{k=m_0}\frac{e^{-(m_0s-\zeta)^2/4}}{\sqrt{4\pi}}.
\end{equation}
Since $B_{m_0}=1$, its Gaussian absorbs the phase weight and proves
$\|\mathfrak g^0\|_{\mathscr X_N^{r_{\rm x}}}\le C_{r_{\rm x}}$.
The descending recursion gives $\rho^\downarrow=\mathbb V_\downarrow\mathfrak g^0$.
Equations \eqref{eq:timeGdowninverse} and \eqref{eq:timeGupsmall} yield
\eqref{eq:descendingTwoIntervalBounds} and \eqref{eq:descendingReferenceDefect}.
Finally $w_{|k|,N}\ge1$, $B_{|k|}\le\widetilde B$ and
$\langle k\rangle^{r_{\rm x}}\le C N^{r_{\rm x}}\langle k/N\rangle^{r_{\rm x}}$
give \eqref{eq:echoDensityEmbeddings} directly.
\end{proof}
\par
\par\endgroup
The full inverse estimates a sum over the spatial indices $k$; the lower bound
uses only $(m_1,{T_1})$. The following functional observes that density
coefficient, {normalized by the echo gain $B_*$.} Since
{$|\rho_{m_1}^\downarrow({T_1})|\ge cB_*/\sqrt N$}, the required normalized
comparison error is $o(N^{-1/2})$. We estimate the observed inverse
before bounding the perturbation source, avoiding an unnecessary factor $N$.

\begingroup
\begin{lem}\label{lem:terminalObservation}
\begingroup
For ${r_{\rm x}}\ge0$, define the normalized terminal evaluation functional
\par\endgroup
\begin{equation}\label{eq:terminalObservationDefinition}
\mathcal O_N{\rho}=B_*^{-1}{\rho}_{m_1}({T_1}).
\end{equation}
On $\mathscr X_N^{r_{\rm x}}$, the collisionless inverse satisfies
\begin{equation}\label{eq:terminalObservationInverse}
\|\mathcal O_N\mathbb V_{\rm e}^0\|_{\mathscr X_N^{r_{\rm x}}\to\mathbb C}
\le C_{r_{\rm x}}.
\end{equation}
More generally, let $\Delta\mathbb T$ be a bounded perturbation of the density integral operator and
suppose that
$\mathbb V=({\mathrm{Id}}-\mathbb K-\mathbb P_\downarrow-\mathbb P_\uparrow-\Delta\mathbb T)^{-1}$
exists on the same space with $\|\mathbb V\|\le C_{r_{\rm x}}N$. Then
\begin{equation}\label{eq:terminalObservationPerturbation}
\|\mathcal O_N\mathbb V\|
\le C_{r_{\rm x}}\big(1+N\|\Delta\mathbb T\|\big).
\end{equation}
In particular this observation remains uniformly bounded whenever
$N\|\Delta\mathbb T\|$ stays bounded.
\end{lem}
\par\endgroup
\begin{proof}

Restrict first to the descending late equation of
Proposition~\ref{prop:approximateDensityInverse}. {For $m\ge m_1$, \eqref{eq:approximateDescendingPathBound} gives}
\begin{equation}\label{eq:terminalSinglePath}
\frac{\langle m_1/N\rangle^{r_{\rm x}}}{B_*}
\|w_{m_1,N}\mathbb A_{m_1+1}\cdots\mathbb A_m\mathbb W_m\mathfrak g_m^+\|_{L^\infty(I_+)}
\le C_{r_{\rm x}}\frac{\langle m/N\rangle^{r_{\rm x}}}{B_m}
\|w_{m,N}\mathfrak g_m^+\|_{L^\infty(I_+)}.
\end{equation}
There is only one descending path to mode $m_1$ at time $T_1$, so summing over $m$ costs no factor $N$. Since $w_{m_1,N}(T_1)\ge1$,
\eqref{eq:lateDensityInputBound} gives
\[
|\mathcal O_N\mathbb V_\downarrow\mathfrak g|
\le C_{r_{\rm x}}\sum_{m\ge1}\frac{\langle m/N\rangle^{r_{\rm x}}}{B_m}
\sup_{I_+}w_{m,N}|\mathfrak g_m^+|
\le C_{r_{\rm x}}\|\mathfrak g\|_{\mathscr X_N^{r_{\rm x}}}.
\]
The negative block does not contribute to this positive-mode observation.
The identity
\[
\mathbb V_{\rm e}^0=\mathbb V_\downarrow+
\mathbb V_\downarrow\mathbb P_\uparrow\mathbb V_{\rm e}^0
\]
and \eqref{eq:timeGupsmall} therefore imply
\[
\|\mathcal O_N\mathbb V_{\rm e}^0\|
\le C_{r_{\rm x}}+C_{r_{\rm x}}aN^4e^{-cN^2}\le C_{r_{\rm x}}.
\]
This proves \eqref{eq:terminalObservationInverse} directly in the
two-interval space.
\par

For the perturbed inverse, multiplication by ${\mathrm{Id}}-\mathbb K-\mathbb P_\downarrow-\mathbb P_\uparrow$ gives the exact identity
\[
\mathbb V=\mathbb V_{\rm e}^0
+\mathbb V_{\rm e}^0\Delta\mathbb T\mathbb V .
\]
Consequently,
\[
\|\mathcal O_N\mathbb V\|
\le\|\mathcal O_N\mathbb V_{\rm e}^0\|
\bigl(1+\|\Delta\mathbb T\|\,\|\mathbb V\|\bigr)
\le C_{r_{\rm x}}\bigl(1+N\|\Delta\mathbb T\|\bigr),
\]
which proves \eqref{eq:terminalObservationPerturbation}.
\par
\end{proof}

\section{\texorpdfstring{{The complete linear echo and inhomogeneous estimates}}{The complete linear echo and inhomogeneous estimates}}\label{sec:completeLinear}

\begingroup
We transfer the growing density from Section~\ref{sec:approximateDensity}
to the complete linearization about $F^{\rm b,\nu}$.
Section~\ref{sec:gainenergy} compares the two collision solution operators
and applies them to the prescribed initial datum.
Section~\ref{sec:exactElectricResponse} estimates the three kernel differences
and constructs the exact density inverse and $H^{\rm C,\nu}$.
Section~\ref{sec:threeDensityErrors} combines the initial, upward-interaction
and kernel errors, proving Proposition~\ref{thm:exactCollisionEcho}.
\par\endgroup
Section~\ref{sec:fullRangeClosure} establishes the complete forced estimate
for sources of finite regularity, with the order-independent exponential
constant fixed before the field correction order $J$ is chosen.
Section~\ref{sec:fieldExpansion} constructs the increments $H_n$, proves
the ordered-operation and derivative estimates, and completes the
linear solution. The exact identities are
\[
\begin{aligned}
H_0&=H^{\rm C,\nu},\qquad H^{[J]}=\sum_{n=0}^JH_n,
&H^\nu&=H^{[J]}+W_J,\\
\mathcal A_\nu H^{[J]}&=\mathcal T_{\rm b}H_J,
&\mathcal A_\nu W_J&=-\mathcal T_{\rm b}H_J,\qquad W_J(t_0)=0.
\end{aligned}
\]
The comparison is therefore organized by
\[
\rho[H^\nu]-\rho^\downarrow
=(\rho^{\rm C,\nu}-\rho^\downarrow)
+\sum_{n=1}^J\rho[H_n]+\rho[W_J].
\]
The first error is controlled by the precise density inverse and terminal
observation; the second by the finite correction estimates; the third
by the complete forced estimate applied to the reduced residual.
This proves Theorem~\ref{thm:logFullLinearEcho}.
{The same section concludes with the spatial energy and density
upper bounds of Proposition~\ref{prop:linearNonlinearInput}. Together with
the electric source estimate in Lemma~\ref{lem:electricFormDual}, these bounds control the quadratic
nonlinear residual in Section~\ref{sec:nonlinear} for the exact solution $H^\nu$.}

\begingroup
\subsection{\texorpdfstring{$c$-shifted vector fields and estimates for the solution operators}{c-shifted vector fields and estimates for the solution operators}}\label{sec:gainenergy}\label{sec:centeredcollisions}\label{sec:bidirectional}\label{sec:exactCollisionEvolution}
The initial data contain high velocity frequencies. An unshifted commuting
derivative would count those frequencies as a loss even for free transport:
\[
(\partial_{v_1}+iks)\big(e^{i(\zeta-ks)v_1}V(v)\big)
=e^{i(\zeta-ks)v_1}(\partial_{v_1}V+i\zeta V).
\]
Subtracting $i\zeta$ removes this loss. More generally, we use a fixed
frequency $c$ and prove collision coercivity uniformly in $c$ before
estimating the corresponding solution operators.

The time translation in the data requires a second distinction. At physical
time $t=s-t_0$, the Fourier coefficients of the data are
\begin{align*}
\left(\frac{F^{\rm b}_{\rm in}-\mu}{\sqrt\mu}\right)_j
&=\mathbf1_{j=0}\frac{\mu_\theta-\mu}{\sqrt\mu}
 +\mathbf1_{|j|=1}\frac a2e^{-ijt_0v_1}\frac{\mu_\theta}{\sqrt\mu},\\
\left(\frac{H_{\rm in}}{\sqrt\mu}\right)_{\pm m_0}
&=e^{\pm i(\zeta-m_0t_0)v_1}\frac{\mu_{1/2}}{\sqrt{4\pi}\sqrt\mu}.
\end{align*}
For the last line, the physical-time frequency is
$c=\pm(\zeta-m_0t_0)$; for a source inserted at time $\tau$, $c$ will be
chosen from the Fourier phase of that source and kept fixed during its
subsequent evolution. The norm $\mathscr A_\ell^{\lambda,\chi}$ in
\eqref{eq:velocityInputGevrey} measures the nonoscillatory function $V(v)$.
The background estimates in this norm were proved in $(z,v)$ coordinates
in Section~\ref{sec:backgroundDecreasingRadius}.

We first treat the Maxwellian collision operator, then the exact background
coefficient. These solution operators contain neither the self-consistent
electric term nor background-field transport. The density equation will be
solved in the next subsection.
\par\endgroup
\begingroup
For a spatial mode $k$ and a fixed real frequency $c$, set
\par\endgroup
\begin{equation}\label{eq:centeredGdefinitions}
\begin{aligned}
Z_{k,c}(s)&=\partial_{v_1}+i(ks-c),\qquad
g_\chi(n)=\frac{(n!)^{1/\chi}}{(n+1)^{12}},\qquad
 a_n=\frac{\lambda^n}{g_\chi(n)},\\
\|f\|_{\mathscr G^{\lambda,\chi}_{\ell;k,c}(s)}^2
&=\sum_{n\ge0}a_n^2\|\langle v\rangle^\ell Z_{k,c}(s)^n f\|_2^2,\qquad \mathfrak D^{\lambda,\chi}_{\ell;k,c}(s;f)
=\sum_{n\ge0}a_n^2\|Z_{k,c}(s)^nf\|_{\Sigma_\ell}^2.
\end{aligned}
\end{equation}
The weighted collision norm $\Sigma_\ell$ is defined in \eqref{eq:backgroundCollisionNorm}. Only the derivatives in the direction of transport occur in this Gevrey norm. The dissipation still contains all three velocity derivatives.
{Indeed,
\[
\mathcal F_v[Z_{k,c}(s)f](\eta)
=i(\eta_1+ks-c)\widehat f(\eta).
\]
Thus the derivative sum controls velocity Fourier decay away from
$\eta_1=c-ks$, with a Gevrey radius determined by $\lambda$.
It does not assert a compact Fourier support. The density trace
\eqref{eq:centeredGtrace} will convert this centered control into
decay away from the resonance $ks=c$.}

\begingroup
\begin{lem}\label{lem:centeredcoercivity}
Fix $\ell\ge0$ and $1/2<\chi<2/3$, and use the coefficients and norms in \eqref{eq:centeredGdefinitions}.

There exist $\lambda_*>0$ and $c_{\mathrm{coer},\ell},C_{\mathrm{coer},\ell}>0$ such that, for $0<\lambda\le\lambda_*$,
\begin{equation}\label{eq:centeredcoercivity}
{
\operatorname{Re}\sum_{n\ge0}a_n^2 \langle Z_{k,c}(s)^n\mathcal Lf, Z_{k,c}(s)^nf\rangle_{\ell} \ge{}c_{\mathrm{coer},\ell}\mathfrak D^{\lambda,\chi}_{\ell;k,c}(s;f) -C_{\mathrm{coer},\ell}\|f\|_{\mathscr G^{\lambda,\chi}_{\ell;k,c}(s)}^2.
}
\end{equation}
The constants are independent of $k,c,s$. The estimate holds also when both sums are truncated at any integer $n={J_{\rm der}}$.
The centered derivative also satisfies \begin{equation}\label{eq:centeredGcommutation}
[\partial_s+ikv_1,Z_{k,c}(s)]=0
\end{equation}
The Maxwellian operator is symmetric and nonnegative:
\[
\langle\mathcal Lf,w\rangle=\langle f,\mathcal Lw\rangle,
\qquad \langle\mathcal Lf,f\rangle\ge0
\]
for smooth rapidly decreasing functions.
\end{lem}
\par\endgroup

\begingroup
\begin{proof}
The symmetry and nonnegativity of $\mathcal L$ follow from
\cite[Lemma 4.5]{ref4}, after the fixed normalization of
Lemma~\ref{lem:fixedPeriodNormalization}. The centered weighted estimate
\eqref{eq:centeredcoercivity} is proved below.

The commutators $[\partial_s,Z_{k,c}]=ik$ and $[ikv_1,\partial_{v_1}]=-ik$ cancel, proving \eqref{eq:centeredGcommutation}.

We first separate the frequency shift from the differentiated coefficients. The shift $ks-c$ in $Z_{k,c}(s)$ is constant in $v$. Set
$\mathcal C_0=\mathcal L$ and define, for $j\ge0$,
\[
\mathcal C_{j+1}f
=\partial_{v_1}(\mathcal C_jf)-\mathcal C_j(\partial_{v_1}f).
\]
The resulting commutator identity is
\begin{equation}\label{eq:centeredcommutatorexpansion}
Z_{k,c}(s)^n\mathcal Lf
=\mathcal LZ_{k,c}(s)^nf+\sum_{j=1}^n\binom nj\mathcal C_jZ_{k,c}(s)^{n-j}f.
\end{equation}
Each $\mathcal C_j$ is independent of $ks-c$.
\par

{Lemma~\ref{lem:weightedMaxwellianForm}}
applied to $Z_{k,c}(s)^nf$ controls the undifferentiated term. We now estimate
the commutators, using the decomposition \eqref{eq:centeredlocaldecomposition}.

The coefficient estimate of \cite[Lemma~A.3]{ref6} gives, for $j\ge1$,
\begin{equation}\label{eq:centeredcoefficientbounds}
|\partial_{v_1}^j \sigma_{iq}(v)|+|\partial_{v_1}^j\mathcal V(v)|
\le C^{j+1}g_\chi(j)\langle v\rangle^{-1-j}.
\end{equation}
We also use the standard consequences of the Maxwellian eigenvalue bounds,
\[
\|\langle v\rangle^{\ell-3/2}\nabla_v y\|_2
+\|\langle v\rangle^{\ell-1/2}y\|_2
\le C\|y\|_{\Sigma_\ell};
\]
see \cite[Appendix~A, equation (A.6)]{ref6}.

For $j\ge2$, integration by parts, \eqref{eq:centeredcoefficientbounds}, and ${|\nabla\langle v\rangle^{2\ell}|\le C_\ell\langle v\rangle^{2\ell-1}}$ give
\begin{equation}\label{eq:centeredhigherlocal}
\left|\left\langle\left(-\partial_i
((\partial_{v_1}^j\sigma_{iq})\partial_q y)
+(\partial_{v_1}^j\mathcal V)y\right),z\right\rangle_{\ell}\right|
\le C^{j+1}g_\chi(j)\|y\|_{\Sigma_\ell}\|z\|_{\Sigma_\ell}.
\end{equation}
For $j=1$ this bound applies to the potential term and the term where a derivative falls on ${\langle v\rangle^{2\ell}}$. The remaining principal term {is estimated after taking the real part,} with $y=Z_{k,c}(s)^{n-1}f$ and $z=Z_{k,c}(s)y$:
\begin{align*}
\operatorname{Re}\langle(\partial_{v_1}\sigma_{iq})\partial_qy,
\partial_i(Z_{k,c}(s)y)\rangle_\ell
&=-\frac12\operatorname{Re}\int
\partial_{v_1}(\langle v\rangle^{2\ell}\partial_{v_1}\sigma_{iq})
\partial_qy\,\overline{\partial_i y}\,dv.
\end{align*}
The contribution of $i(ks-c)$ is zero, since the quadratic form with the real symmetric matrix $\partial_{v_1}\sigma$ is real. The coefficient on the right is bounded by ${C\langle v\rangle^{2\ell-3}}$, so the absolute value of this real part is at most $C\|Z_{k,c}(s)^{n-1}f\|_{\Sigma_\ell}^2$. This is the only term for which a general bilinear bound with two different functions would be insufficient.

For the nonlocal part, set $\mathcal K_0=\mathcal K$ and define
\[
\mathcal K_{j+1}f
=\partial_{v_1}(\mathcal K_jf)-\mathcal K_j(\partial_{v_1}f),\qquad j\ge0.
\]
At the kernel level,
\[
(\partial_{v_1}+\partial_{v_{*1}})\Phi_{iq}(v-v_*)=0.
\]
Thus all $j$ derivatives fall on Gaussian factors and polynomial
prefactors. After conjugation by $\langle v\rangle^\ell$, their
absolute value is bounded by
\[
C^{j+1}g_\chi(j)e^{-c(|v|^2+|v_*|^2)}|v-v_*|^{-1}.
\]
The Schur bounds are
\[
\sup_v\int_{\mathbb R^3}
\frac{e^{-c(|v|^2+|v_*|^2)}}{|v-v_*|}\,dv_*
+\sup_{v_*}\int_{\mathbb R^3}
\frac{e^{-c(|v|^2+|v_*|^2)}}{|v-v_*|}\,dv\le C.
\]
Together with the Gaussian multiplication term, this proves
\par
\begin{equation}\label{eq:centerednonlocalbound}
\|\langle v\rangle^\ell\mathcal K_j y\|_2
\le C^{j+1}g_\chi(j)\|\langle v\rangle^\ell y\|_2.
\end{equation}
This proves the required nonlocal estimate without differentiating the function $y$ on which $\mathcal K_j$ acts.

We now sum exactly the terms just estimated. 
The coefficients satisfy
\begin{equation}\label{eq:centeredfactorialratio}
\frac{a_n\binom nj C^jg_\chi(j)}{a_{n-j}}
=(C\lambda)^j\binom nj^{1-1/\chi}
\left(\frac{n+1}{(j+1)(n-j+1)}\right)^{12}
\le(C\lambda)^j,\qquad 1\le j\le n.
\end{equation}
\par
For a finite cutoff ${J_{\rm der}}$, set
\[
r_n=a_n\|Z_{k,c}(s)^nf\|_{\Sigma_\ell},\qquad
 e_n=a_n\|\langle v\rangle^\ell Z_{k,c}(s)^nf\|_2\quad(0\le n\le {J_{\rm der}}),
\]
and extend both sequences by zero outside this range.
The terms covered by \eqref{eq:centeredhigherlocal} and the order-one lower-order terms are bounded in absolute value by
\[
C\sum_{n=0}^{J_{\rm der}}\sum_{j=1}^n(C\lambda)^j r_{n-j}r_n
\le C\lambda\sum_{n=0}^{J_{\rm der}}r_n^2.
\]
The last inequality is the $\ell^1*\ell^2\to\ell^2$ convolution inequality, with $\sum_{j\ge1}(C\lambda)^j\le C'\lambda$.
Equation \eqref{eq:centerednonlocalbound} gives, in exactly the same way, $C\lambda\sum_{n=0}^{J_{\rm der}}e_n^2$.
The exceptional principal term with $j=1$ is bounded by
\[
C\sum_{n=1}^{J_{\rm der}} n a_n^2\|Z_{k,c}(s)^{n-1}f\|_{\Sigma_\ell}^2
\le C\lambda^2\sum_{n=0}^{J_{\rm der}}r_n^2,
\]
because $n a_n^2/a_{n-1}^2\le C\lambda^2$.
Combining these bounds with the undifferentiated coercivity yields
\[
\operatorname{Re}\sum_{n=0}^{J_{\rm der}} a_n^2
\langle Z_{k,c}(s)^n\mathcal Lf,Z_{k,c}(s)^nf\rangle_{\ell}
\ge(c_\ell-C\lambda)\sum_{n=0}^{J_{\rm der}}r_n^2
-(C_\ell+C\lambda)\sum_{n=0}^{J_{\rm der}}e_n^2.
\]
Choose $\lambda_*$ small. All constants are independent of ${J_{\rm der}}$ and $ks-c$.
This proves the finite-sum assertion. The infinite-sum assertion follows whenever the pairings are defined, or as the limit of these finite-sum form estimates; this is the form used in the evolution argument below.
\end{proof}
\par\endgroup
\begingroup
For an integer $j$ and a function $V=V(v)$ independent of $x_1$,
$\mathbf e_jV$ denotes the sequence with value $V$ at index $j$ and zero
elsewhere. Recall $\mathscr A_\ell^{\lambda,\chi}$ from
\eqref{eq:velocityInputGevrey}. In the notation just introduced,
\[
\|V\|_{\mathscr A_\ell^{\lambda,\chi}}
=\|V\|_{\mathscr G^{\lambda,\chi}_{\ell;0,0}(s)}.
\]
\par\endgroup

Use \eqref{eq:comparisonLogParameters} and fix $\ell\ge0$.
Recall the normalized transport profiles $g^{\rm b,\nu}$,
$g^{\rm b}_{\rm app}$ and {${g^{\rm b}_{\rm err}}$} from
\eqref{eq:physicalCollisionProfiles}, in the coordinates
$t=s-t_0$, $z=x_1-tv_1$.
To express the collision coefficients in physical coordinates, write
{${F^{\rm c}}_l(s,v)=\sqrt \mu  e^{-il(s-t_0)v_1}g_l(s-t_0,v)$}.
Let {$\mathcal C_{F^{\rm c}}=\mu^{-1/2}\big(Q({F^{\rm c}},\sqrt \mu \,\cdot)+Q(\sqrt \mu \,\cdot,{F^{\rm c}})\big)$}.
In particular,
\[
\mathcal C_{\rm b}=-\mathcal L+\mathcal C_{F^{\rm b,\nu}-\mu }.
\]
This operator acts on distributions normalized by {$\sqrt \mu $}.
The coefficient amplitude is $\varepsilon_N=aL_N=\nu^{1/3-o(1)}$.
{The amplitude factor $a$ allows} fixed coefficient radii in the collision comparison.

Use a hat to distinguish the physical-time centered derivative from \eqref{eq:centeredGdefinitions}:
\begin{equation}\label{eq:physicalCenteredNorms}
\begin{aligned}
\widehat Z_{k,c}(s)&=\partial_{v_1}+i\big(k(s-t_0)-c\big),\qquad \kappa_{k,j}=e^{\varrho_0|k-j|^\chi},\qquad
a_n=\lambda_{\rm v}^n/g_\chi(n),\\
\|f\|_{\mathfrak G_{c,j}^{\varrho,\lambda,\ell}(s)}^2
&=\sum_{k\in\mathbb Z}e^{2\varrho|k-j|^\chi}
\sum_{n\ge0}\frac{\lambda^{2n}}{g_\chi(n)^2}
\|\langle v\rangle^\ell\widehat Z_{k,c}(s)^nf_k\|_2^2,\\
\mathcal E^{\rm ph}_{c,j}(s;f)
&=\|f\|_{\mathfrak G_{c,j}^{\varrho_0,\lambda_{\rm v},\ell}(s)}^2,\qquad \mathcal D^{\rm ph}_{c,j}(s;f)
=\sum_{k\in\mathbb Z}\kappa_{k,j}^2\sum_{n\ge0}a_n^2
\|\widehat Z_{k,c}(s)^nf_k\|_{\Sigma_\ell}^2.
\end{aligned}
\end{equation}
At physical elapsed time $t=s-t_0$, the relation with
\eqref{eq:centeredGdefinitions} is
\begin{equation}\label{eq:physicalSingleModeRelations}
\begin{aligned}
\widehat Z_{k,c}(s)&=Z_{k,c}(t),\qquad \|f\|_{\mathfrak G_{c,j}^{\varrho,\lambda,\ell}(s)}^2
=\sum_{k\in\mathbb Z}e^{2\varrho|k-j|^\chi}
\|f_k\|_{\mathscr G^{\lambda,\chi}_{\ell;k,c}(t)}^2,\\
\mathcal D^{\rm ph}_{c,j}(s;f)
&=\sum_{k\in\mathbb Z}e^{2\varrho_0|k-j|^\chi}
\mathfrak D^{\lambda_{\rm v},\chi}_{\ell;k,c}(t;f_k).
\end{aligned}
\end{equation}
Thus the physical norms are spatially weighted sums of the existing
single-mode norms evaluated at $s-t_0$.
The radii $\varrho,\lambda>0$ are displayed when they vary.
In $\mathcal E^{\rm ph}$ and $\mathcal D^{\rm ph}$ we use the
sufficiently small fixed radii $\varrho_0,\lambda_{\rm v}>0$.

The coefficient bounds and the two collision-ordering estimates are proved in
Lemmas~\ref{lem:physicalCoefficientDerivatives} and \ref{lem:physicalCollisionForm}.
Lemma~\ref{lem:localizedConstruction} constructs the evolutions in the indicated
energy spaces, including strong time traces. We apply these results first to the
Maxwellian coefficient and then to their difference; the quantitative estimates
needed for the echo are proved below.

Let $\mathbb S_k^\nu(t)$ denote the solution operator of
\[
\partial_tf+ikv_1f+\nu\mathcal Lf=0,
\qquad \mathbb S_k^0(t)f=e^{-iktv_1}f.
\]
This is the Maxwellian transport--Landau evolution used in \cite[Section~6]{ref4}; it has no electric term. We next record {both the centered energy estimate and} the density estimate extracted from that energy.

\begingroup
\begin{prop}\label{prop:centeredpropagation}
{Fix $\ell\ge4$, $1/2<\chi<2/3$, $k\in\mathbb Z$, $c\in\mathbb R$,
$\nu\ge0$, and a sufficiently small $\lambda>0$.
Let $\tau<T<\infty$, and let $f\in C([\tau,T];L^2_v)$ be a mild
solution of}
\[
\partial_sf+ikv_1f+\nu\mathcal Lf={S},
\qquad f(\tau)=f_\tau{.}
\]
{Assume that the initial datum and source satisfy
\[
\|f_\tau\|_{\mathscr G^{\lambda,\chi}_{\ell;k,c}(\tau)}
+\int_\tau^T\|{S}(s)\|_{\mathscr G^{\lambda,\chi}_{\ell;k,c}(s)}\,ds<\infty.
\]
Then, for every $s\in[\tau,T]$,}
\begin{equation}\label{eq:centeredpropagation}
\|f(s)\|_{\mathscr G^{\lambda,\chi}_{\ell;k,c}(s)}
\le e^{C\nu(s-\tau)}
\left(\|f_\tau\|_{\mathscr G^{\lambda,\chi}_{\ell;k,c}(\tau)}
+\int_\tau^s\|{S}(\tau')\|_{\mathscr G^{\lambda,\chi}_{\ell;k,c}(\tau')}\,d\tau'\right).
\end{equation}
{For every initial datum and source satisfying these conditions,
such a mild solution exists and is unique.
Moreover, there is $C\ge1$, independent of $k,c,s$, such that}
\begin{equation}\label{eq:centeredGtrace}
\left|\int\sqrt \mu  f(s,v)\,dv\right|
\le Ce^{-C^{-1}|ks-c|^\chi}
\|f(s)\|_{\mathscr G^{\lambda,\chi}_{\ell;k,c}(s)}.
\end{equation}

In particular, let $V$ be any one of
\begin{equation}\label{eq:centeredGaussianlist}
\sqrt \mu ,\qquad v_1\sqrt \mu ,\qquad
\frac{\mu_{1/2}}{\sqrt{4\pi}\sqrt \mu }.
\end{equation}
For $t\ge0$ and $\eta\in\mathbb R$,
\begin{align}
\left|\int\sqrt \mu \,\mathbb S_k^\nu(t)(e^{i\eta v_1}V)\,dv\right|
&\le Ce^{C\nu t}e^{-C^{-1}|kt-\eta|^\chi},
\label{eq:centeredpackettrace}\\
\left|\int\sqrt \mu \,(\mathbb S_k^\nu(t)-\mathbb S_k^0(t))
(e^{i\eta v_1}V)\,dv\right|
&\le C\nu t(1+\eta^2+k^2t^2)
 e^{C\nu t}e^{-C^{-1}|kt-\eta|^\chi}.
\label{eq:centeredpacketdifference}
\end{align}
The same bounds hold for linear combinations of these three functions, multiplied by the sum of the absolute values of their coefficients.
\end{prop}
\par\endgroup
\begin{proof}
We use the Maxwellian semigroup defined in \cite[(6.2)--(6.3)]{ref4}, with the fixed normalization of Lemma~\ref{lem:fixedPeriodNormalization}. At $k=0$ it is the spectral semigroup of the nonnegative self-adjoint operator $\nu\mathcal L$. For fixed $k$, $\nu\ge0$ and $I'=[\tau,{T}]$, {the $L^2$ contraction of $\mathbb S_k^\nu$ gives} the unique mild solution
\begin{equation}\label{eq:centeredMildConstruction}
f(s)=\mathbb S_k^\nu(s-\tau)f_\tau
+\int_\tau^s\mathbb S_k^\nu(s-\tau'){S}(\tau')\,d\tau',
\qquad f\in C(I';L^2_v).
\end{equation}
In particular,
\[
\|f-\widetilde f\|_{C(I';L^2)}
\le\|f_\tau-\widetilde f_\tau\|_2
+\|{S}-\widetilde {S}\|_{L^1(I';L^2)}.
\]

If $\nu=0$, the formula is free transport and $[\partial_s+ikv_1,Z_{k,c}(s)]=0$; Minkowski's inequality directly proves \eqref{eq:centeredpropagation}. Hence suppose $\nu>0$ in the form argument.

Apply Lemma~\ref{lem:localizedConstruction} with the parameters
\[
{F^{\rm c}}=0,\quad \varepsilon_{\rm c}=0,\quad s_*=0,\quad
\kappa_j=e^{\varrho_0|j-k|^\chi},\quad
f_\tau^{\rm seq}=\mathbf e_k f_\tau,\quad
S^{\rm seq}(s)=\mathbf e_k {S}(s).
\]
The coefficient and perturbation-form hypotheses are then identically
zero. Since $\mathcal L$ is spatially diagonal, uniqueness implies
\[
f_j^{\rm seq}(s)=0\quad(j\ne k),\qquad
\mathcal Z_k(s)=\partial_{v_1}+i(ks-c)=Z_{k,c}(s).
\]
In particular the finite sums of that lemma become exactly
\[
E_{J_{\rm der}}(s;f^{\rm seq})
=\sum_{n=0}^{J_{\rm der}} a_n^2\|\langle v\rangle^\ell Z_{k,c}(s)^nf\|_2^2,
\qquad
D_{J_{\rm der}}(s;f^{\rm seq})=\sum_{n=0}^{J_{\rm der}}a_n^2\|Z_{k,c}(s)^nf\|_{\Sigma_\ell}^2.
\]
{The difference-quotient construction and regularization limit in that lemma}
therefore justify, for every finite ${J_{\rm der}}$,
\begin{equation}\label{eq:centeredFiniteEnergy}
\begin{aligned}
\tfrac12\partial_sE_{J_{\rm der}}(s;f^{\rm seq})+c\nu D_{J_{\rm der}}
&\le C\nu E_{J_{\rm der}}(s;f^{\rm seq})+E_{J_{\rm der}}(s;f^{\rm seq})^{1/2} E_{J_{\rm der}}(s;S^{\rm seq})^{1/2},\\
E_{J_{\rm der}}(s;f^{\rm seq})^{1/2}&\le e^{C\nu(s-\tau)}
\left(E_{J_{\rm der}}(\tau;f_\tau^{\rm seq})^{1/2}+\int_\tau^sE_{J_{\rm der}}(\tau';S^{\rm seq})^{1/2}\,d\tau'\right),
\end{aligned}
\end{equation}
This application includes $k=0$: the regularized commutator
\eqref{eq:constructionRegularizerCommutator} and the difference estimate
\eqref{eq:constructionDifferenceEnergy} require no division by $k$.
The constructed solution agrees with \eqref{eq:centeredMildConstruction}
by $L^2$ energy uniqueness. Finally,
\[
E_{J_{\rm der}}(s;f^{\rm seq})\uparrow\|f\|_{\mathscr G^{\lambda,\chi}_{\ell;k,c}}^2,
\qquad D_{J_{\rm der}}\uparrow\mathfrak D^{\lambda,\chi}_{\ell;k,c}(f),
\qquad {J_{\rm der}}\to\infty.
\]
Monotone convergence in \eqref{eq:centeredFiniteEnergy} proves
\eqref{eq:centeredpropagation}, with constants independent of $k,c,{J_{\rm der}}$.

To obtain the density, repeated integration by parts gives
\[
[i(ks-c)]^n\int\sqrt \mu  f\,dv
=\sum_{j=0}^n\binom nj
\int(\partial_{v_1}^j\sqrt \mu ) Z_{k,c}(s)^{n-j}f\,dv.
\]
The Gaussian derivative bounds and Cauchy--Schwarz imply
\[
|ks-c|^n\left|\int\sqrt \mu  f\,dv\right|
\le C(C_1/\lambda)^n(n!)^{1/\chi}\|f(s)\|_{\mathscr G^{\lambda,\chi}_{\ell;k,c}(s)}.
\]
If $(\lambda|ks-c|/(2C_1))^\chi\ge2$, take
$n=\lfloor(\lambda|ks-c|/(2C_1))^\chi\rfloor$. Since $n!\le n^n$,
\[
\left(\frac{C_1}{\lambda|ks-c|}\right)^n(n!)^{1/\chi}
\le\left(\frac{C_1n^{1/\chi}}{\lambda|ks-c|}\right)^n
\le2^{-n}\le C\exp\big(-C^{-1}(\lambda|ks-c|)^\chi\big).
\]
For bounded $\lambda|ks-c|$ use $n=0$. This proves
\eqref{eq:centeredGtrace} at the fixed radius $\lambda$.

For $c=\eta$ at initial time zero,
\[
Z_{k,\eta}(0)^n(e^{i\eta v_1}V)=e^{i\eta v_1}\partial_{v_1}^nV,
\qquad
\|e^{i\eta v_1}V\|_{\mathscr G^{\lambda,\chi}_{\ell;k,\eta}(0)}
=\|V\|_{\mathscr A_\ell^{\lambda,\chi}}.
\]
Equations \eqref{eq:centeredpropagation} and \eqref{eq:centeredGtrace} give
\[
\left|\int\sqrt\mu\,\mathbb S_k^\nu(t)(e^{i\eta v_1}V)\,dv\right|
\le Ce^{C\nu t-C^{-1}|kt-\eta|^\chi}
\|V\|_{\mathscr A_\ell^{\lambda,\chi}},
\]
which proves \eqref{eq:centeredpackettrace} for the three fixed Gaussians.
\par
For the difference estimate, write
\[
{f^{\rm free}}(\tau,v)=e^{i(\eta-k\tau)v_1}V(v),\qquad
 e(\tau)=\mathbb S_k^\nu(\tau)(e^{i\eta v_1}V)-{f^{\rm free}}(\tau).
\]
Then $e(0)=0$ and
\begin{equation}\label{eq:centeredpacketerror}
\partial_\tau e+ikv_1e+\nu\mathcal Le=-\nu\mathcal L{f^{\rm free}}.
\end{equation}
We now estimate this specified source; no derivatives of the unknown $e$ are put on the right-hand side.

For the listed Gaussians and a sufficiently small fixed $\lambda$,
\begin{equation}\label{eq:centeredGaussianforcing}
\|\mathcal L{f^{\rm free}}(\tau)\|_{\mathscr G^{\lambda,\chi}_{\ell;k,\eta}(\tau)}
\le C\big(1+|\eta-k\tau|^2\big).
\end{equation}
Use $Z_{k,\eta}(\tau)=\partial_{v_1}+i(k\tau-\eta)$ and the commutators
$\mathcal C_j$ defined in the proof of Lemma~\ref{lem:centeredcoercivity}.
The exact identities are
\begin{equation}\label{eq:gaussianSourceCommutator}
Z_{k,\eta}(\tau)^n\mathcal L{f^{\rm free}}
=\sum_{j=0}^n\binom nj\mathcal C_j
\big(e^{i(\eta-k\tau) v_1}\partial_{v_1}^{n-j}V\big),\qquad
Z_{k,\eta}(\tau)^m{f^{\rm free}}=e^{i(\eta-k\tau) v_1}\partial_{v_1}^mV.
\end{equation}
The local and nonlocal contributions to these commutators are, explicitly,
\[
\mathcal C_j^{\rm loc}f
=-\partial_{v_i}\big((\partial_{v_1}^j\sigma_{iq})\partial_{v_q}f\big)
 +(\partial_{v_1}^j\mathcal V)f,\qquad \mathcal C_j=\mathcal C_j^{\rm loc}-\mathcal K_j,\qquad j\ge0.
\]
Direct differentiation gives
\begin{equation}\label{eq:gaussianSourceLocalConjugation}
\begin{aligned}
e^{-i(\eta-k\tau) v_1}\mathcal C_j^{\rm loc}(e^{i(\eta-k\tau) v_1}\psi )
={}&-\partial_{v_i}((\partial_{v_1}^j \sigma_{iq})\partial_{v_q}\psi )
+(\partial_{v_1}^j\mathcal V)\psi +(\eta-k\tau)^2 (\partial_{v_1}^j \sigma_{11})\psi \\
&-i(\eta-k\tau)\big((\partial_{v_i}(\partial_{v_1}^j \sigma_{i1}))\psi 
+2(\partial_{v_1}^j \sigma_{i1})\partial_{v_i}\psi \big).
\end{aligned}
\end{equation}
For the three functions in \eqref{eq:centeredGaussianlist}, Gaussian
differentiation and integration against a Gaussian with a smaller fixed
exponent give, for every $m\ge0$,
\[
\sum_{|\alpha|\le2}
\|\langle v\rangle^{\ell+2}\partial_v^\alpha
\partial_{v_1}^mV\|_2
\le C^{m+1}\sqrt{(m+3)!}.
\]
The differentiated coefficient bounds in Lemma~\ref{lem:centeredcoercivity},
\eqref{eq:gaussianSourceLocalConjugation}, and the nonlocal estimate
\eqref{eq:centerednonlocalbound} consequently yield
\[
\begin{aligned}
&\|\langle v\rangle^\ell\mathcal C_j^{\rm loc}
(e^{i(\eta-k\tau)v_1}\partial_{v_1}^mV)\|_2 +\|\langle v\rangle^\ell\mathcal K_j
(e^{i(\eta-k\tau)v_1}\partial_{v_1}^mV)\|_2\\
&\qquad\le C^{j+m+2}g_\chi(j)
\big(1+(\eta-k\tau)^2\big)\sqrt{(m+3)!}.
\end{aligned}
\]
The frequency cost is quadratic, independently of $j,m$.
With $a_n=\lambda^n/g_\chi(n)$, equations
\eqref{eq:gaussianSourceCommutator} and \eqref{eq:centeredfactorialratio}
give
\begin{align*}
&a_n\|\langle v\rangle^\ell Z_{k,\eta}(\tau)^n\mathcal L{f^{\rm free}}\|_2\\*
&\quad\le C\big(1+(\eta-k\tau)^2\big)
\sum_{j=0}^n(C\lambda)^j a_{n-j}
C^{n-j+1}\sqrt{(n-j+3)!},\\
&\|\mathcal L{f^{\rm free}}\|_{\mathscr G^{\lambda,\chi}_{\ell;k,\eta}(\tau)}\\*
&\quad\le C\big(1+(\eta-k\tau)^2\big)
\left(\sum_{j\ge0}(C\lambda)^j\right)
\left(\sum_{m\ge0}
\frac{C^{2m+2}\lambda^{2m}(m+3)!}{g_\chi(m)^2}\right)^{1/2}\\*
&\quad\le C\big(1+(\eta-k\tau)^2\big).
\end{align*}
Here $C\lambda<1/2$, and consecutive summands of the last series
have ratio
\[
(C\lambda)^2\frac{m+4}{(m+1)^{2/\chi}}
\left(\frac{m+2}{m+1}\right)^{24}\longrightarrow0.
\]
This proves \eqref{eq:centeredGaussianforcing}.

Apply \eqref{eq:centeredpropagation} to \eqref{eq:centeredpacketerror}. It gives
\[
\|e(t)\|_{\mathscr G^{\lambda,\chi}_{\ell;k,\eta}(t)}
\le C\nu e^{C\nu t}\int_0^t(1+|\eta-k\tau|^2)\,d\tau
\le C\nu t(1+\eta^2+k^2t^2)e^{C\nu t}.
\]
Equation \eqref{eq:centeredGtrace} now supplies the frequency decay in
\eqref{eq:centeredpacketdifference}. Linearity proves the assertion for finite linear combinations.
\end{proof}

The form estimate now permits the coefficient difference to act on the evolved distribution itself. We compare exact evolutions, retaining the dissipation generated by each. This is the step that avoids placing two uncontrolled velocity derivatives on the source.

\begingroup
{Use the solution operator $\mathbb U_{\rm b}^\nu(s,\tau)$ of
{\eqref{eq:exactCollisionEvolution}, defined in Section~\ref{sec:roadmapOperators}.}
The Maxwellian solution operator is
{$[\mathbb U_\mu^\nu(s,\tau)f]_k=\mathbb S_k^\nu(s-\tau)f_k$}.
Neither operator includes an electric term in the density equation or background-field
transport.} Set
\begin{equation}\label{eq:physicalCollisionTimeCost}
\Theta_{\rm b}(s,\tau)=\nu(s-\tau)
+\nu\varepsilon_N\int_\tau^s(1+(\tau'-t_0)^2)\,d\tau'.
\end{equation}
\par\endgroup

\begingroup
\begin{prop}\label{prop:exactCollisionPropagation}
Use \eqref{eq:comparisonLogParameters} with $p>{\rm s}_0+7/2$, fix $\ell\ge4$, and take the radii in \eqref{eq:physicalCenteredNorms} sufficiently small for Lemma~\ref{lem:physicalCollisionForm}.

{Let $\tau\in I$, $j\in\mathbb Z$, $c\in\mathbb R$, and let
$f_\tau$ satisfy $\mathcal E^{\rm ph}_{c,j}(\tau;f_\tau)<\infty$.
\begingroup
Let $f$ solve one of the two initial-value problems
\[
\begin{cases}
\partial_sf+v_1\partial_{x_1}f+\nu\mathcal Lf=0,\\
f(\tau)=f_\tau,
\end{cases}
\qquad\text{or}\qquad
\begin{cases}
\partial_sf+v_1\partial_{x_1}f-\nu\mathcal C_{\rm b}(s)f=0,\\
f(\tau)=f_\tau,
\end{cases}
\]
where
$\mathcal C_{\rm b}(s)f=\mu^{-1/2}\big[Q(F^{\rm b,\nu},\sqrt\mu f)
+Q(\sqrt\mu f,F^{\rm b,\nu})\big]$.
\par\endgroup Write $f=\mathbb U_\sharp^\nu(\cdot,\tau)f_\tau$, where
$\mathbb U_\sharp^\nu$ denotes {$\mathbb U_\mu^\nu$} or
$\mathbb U_{\rm b}^\nu$, respectively. Then, for $\tau\le s\le T_1$,}
\begin{equation}\label{eq:physicalCollisionEnergy}
\sup_{\tau\le s'\le s}\mathcal E^{\rm ph}_{c,j}(s';f)
+\nu\int_\tau^s\mathcal D^{\rm ph}_{c,j}(s';f)\,ds'
\le C e^{C\Theta_{\rm b}(s,\tau)}
\mathcal E^{\rm ph}_{c,j}(\tau;f_\tau).
\end{equation}
{Let {$w=(\mathbb U_{\rm b}^\nu-\mathbb U_\mu^\nu)(\cdot,\tau)f_\tau$}.
Then}
\begin{equation}\label{eq:exactCollisionEvolutionDifference}
\sup_{\tau\le s'\le s}\mathcal E^{\rm ph}_{c,j}(s';w)^{1/2}
+\left(\nu\int_\tau^s\mathcal D^{\rm ph}_{c,j}(s';w)\,ds'\right)^{1/2}
\le C\varepsilon_N(1+\nu {T_1}^3)
\mathcal E^{\rm ph}_{c,j}(\tau;f_\tau)^{1/2}.
\end{equation}
{Both initial-value problems admit unique solutions in each fixed
polynomially weighted finite-energy class.
In particular, let $V\in\mathscr A_\ell^{\lambda_{\rm v},\chi}$,
$\eta\in\mathbb R$, $f_\tau=\mathbf e_j(e^{i\eta v_1}V)$, and
$c=j(\tau-t_0)+\eta$. Then, for every $k\in\mathbb Z$,}
\begin{align}
\left|\int\sqrt \mu [\mathbb U_\sharp^\nu(s,\tau)f_\tau]_k\,dv\right|
&\le Ce^{-\varrho_0|k-j|^\chi}e^{-C^{-1}|k(s-t_0)-c|^\chi}
\|V\|_{\mathscr A_\ell^{\lambda_{\rm v},\chi}},
\label{eq:exactCollisionTrace}\\
\left|\int\sqrt \mu [(\mathbb U_{\rm b}^\nu-\mathbb U_\mu^\nu)(s,\tau)f_\tau]_k\,dv\right|
&\le C\varepsilon_Ne^{-\varrho_0|k-j|^\chi}e^{-C^{-1}|k(s-t_0)-c|^\chi}
\|V\|_{\mathscr A_\ell^{\lambda_{\rm v},\chi}}.
\label{eq:exactCollisionTraceDifference}
\end{align}
{These estimates include $j=0$ and $k=0$.}
The energy estimate is uniform for $0<\lambda_{\rm v}\le\lambda_*$, where $\lambda_*>0$ is sufficiently small and fixed. In particular, at the radius $\lambda_{\rm v}/N$ it also gives
\begin{equation}\label{eq:exactCollisionVariableRadiusTrace}
\begin{split}
\left|\int\sqrt \mu [\mathbb U_\sharp^\nu(s,\tau)
\mathbf e_j(e^{i\eta v_1}V)]_k\,dv\right|
\le{}&C e^{-\varrho_0|k-j|^\chi}
e^{-C^{-1}(|k(s-t_0)-c|/N)^\chi}\\
&\quad\times\|V\|_{\mathscr A_\ell^{\lambda_{\rm v}/N,\chi}},
\qquad c=j(\tau-t_0)+\eta.
\end{split}
\end{equation}
The spatial radius $\varrho_0$ remains fixed. The constants in this estimate are independent of $N$.
Under \eqref{eq:comparisonLogParameters} with $p>{\rm s}_0+7/2$, $\Theta_{\rm b}({T_1},t_0)\le C\nu {T_1}+C\varepsilon_N\nu {T_1}^3=o(1)$.
\end{prop}
\par\endgroup
\begin{proof}
By Lemma~\ref{lem:logarithmicScales},
\[
\begin{aligned}
\Theta_{\rm b}(T_1,t_0)
&=\nu(T_1-t_0)+\nu\varepsilon_N
\left(T_1-t_0+\frac{(T_1-t_0)^3}{3}\right)\longrightarrow0,\\
\nu T_1^3&\sim[\log(1/\nu)]^{3{\rm s}_0+6-3p}=o(1),
\qquad \varepsilon_N=\nu^{1/3-o(1)}.
\end{aligned}
\]
We apply Lemma~\ref{lem:localizedConstruction} with
\[
s_*=t_0,\qquad c_{\rm col}(s)=1+(s-t_0)^2,\qquad
(F^{\rm c},\varepsilon_{\rm c})=
\begin{cases}
(0,0),&\mathbb U_\mu^\nu,\\
(F^{\rm b,\nu}-\mu,\varepsilon_N),&\mathbb U_{\rm b}^\nu.
\end{cases}
\]
For one fixed $J_v>\ell+8$, Lemma~\ref{lem:physicalCoefficientDerivatives}
gives the coefficient hypothesis. The spatial weight hypotheses follow from
\[
\begin{aligned}
\frac{\kappa_{k,j}}{\kappa_{m,j}}e^{-\varrho_{\rm b}|k-m|^\chi}
&\le e^{-(\varrho_{\rm b}-\varrho_0)|k-m|^\chi},\\
\max\left\{\frac{\langle k\rangle^r}{\langle m\rangle^r},
\frac{\langle k/N\rangle^r}{\langle m/N\rangle^r}\right\}
e^{-\varrho_{\rm b}|k-m|^\chi}
&\le C_r\langle k-m\rangle^r e^{-\varrho_{\rm b}|k-m|^\chi}
\le C_r e^{-\varrho_{\rm b}|k-m|^\chi/2}.
\end{aligned}
\]
For the finite sums in \eqref{eq:constructionFiniteNorms},
Lemma~\ref{lem:physicalCollisionForm} gives
\[
\begin{aligned}
\left|\sum_k\kappa_{k,j}^2\sum_{n=0}^{J_{\rm der}}a_n^2
\langle\widehat Z_{k,c}^n(\mathcal C_{F^{\rm c}}f)_k,
\widehat Z_{k,c}^nf_k\rangle_\ell\right|
\le C\varepsilon_N\big[D_{J_{\rm der}}(s;f)
+(1+(s-t_0)^2)E_{J_{\rm der}}(s;f)\big].
\end{aligned}
\]
Since $[\partial_s+ikv_1,\widehat Z_{k,c}]=0$, Lemma~\ref{lem:centeredcoercivity}
and the regularization limits of Lemma~\ref{lem:localizedConstruction} yield
\[
\frac12\partial_s\mathcal E^{\rm ph}_{c,j}
+(c_{\mathrm{coer},\ell}-C\varepsilon_N)\nu\mathcal D^{\rm ph}_{c,j}
\le C\nu[1+\varepsilon_N(1+(s-t_0)^2)]\mathcal E^{\rm ph}_{c,j}.
\]
Take $C\varepsilon_N\le c_{\mathrm{coer},\ell}/2$. As
$\partial_s\Theta_{\rm b}=\nu[1+\varepsilon_N(1+(s-t_0)^2)]$,
increasing the fixed constant $C$ gives
\[
\frac{d}{ds}\left(e^{-C\Theta_{\rm b}(s,\tau)}\mathcal E^{\rm ph}_{c,j}(s;f)\right)
+c_{\mathrm{coer},\ell}\nu e^{-C\Theta_{\rm b}(s,\tau)}\mathcal D^{\rm ph}_{c,j}(s;f)\le0.
\]
Integration proves \eqref{eq:physicalCollisionEnergy}.
Lemma~\ref{lem:localizedConstruction} also supplies existence, uniqueness,
strong time continuity and the Duhamel formula in the stated spaces.
\par
For the difference, put {${F^{\rm c}}=F^{\rm b,\nu}-\mu $}. First take finite-mode Schwartz initial data and use
the same projection $P_{J_x}$ and regularizer $\epsilon\mathcal H_v$
for both evolutions in \eqref{eq:constructionRegularization}. Write
{$f^{J_x,\epsilon}=f_\mu^{J_x,\epsilon}$} and
{$w^{J_x,\epsilon}=f^{{\rm b},J_x,\epsilon}-f_\mu^{J_x,\epsilon}$}. Subtraction gives
\begin{equation}\label{eq:collisionDifferenceRegularized}
\big(\partial_s+v_1\partial_{x_1}-\nu P_{J_x}\mathcal C_{\rm b}P_{J_x}
+\epsilon\mathcal H_v\big)w^{J_x,\epsilon}
=\nu P_{J_x}\mathcal C_{F^{\rm c}} P_{J_x} f^{J_x,\epsilon},\qquad w^{J_x,\epsilon}(\tau)=0.
\end{equation}
At each finite derivative order ${J_{\rm der}}$, the regularized solutions have
the regularity proved in Steps 1--3 of Lemma~\ref{lem:localizedConstruction}.
Thus the collision source can be tested in the form pairing.
Abbreviate the finite energies, dissipations and regularization forms by
$E_{f,{J_{\rm der}}},D_{f,{J_{\rm der}}}$ and $E_{w,{J_{\rm der}}},D_{w,{J_{\rm der}}},\mathfrak R_{w,{J_{\rm der}}}$.
Lemma~\ref{lem:physicalCollisionForm} and
\eqref{eq:constructionRegularizerCommutator} give
\begin{equation}\label{eq:collisionDifferenceFiniteEnergy}
\begin{split}
&\tfrac12 E_{w,{J_{\rm der}}}'+c\nu D_{w,{J_{\rm der}}}+c\epsilon \mathfrak R_{w,{J_{\rm der}}} \le C\big(\nu[1+\varepsilon_N(1+(s-t_0)^2)]+\epsilon\big)E_{w,{J_{\rm der}}}\\
&\qquad+C\nu\varepsilon_ND_{f,{J_{\rm der}}}^{1/2}D_{w,{J_{\rm der}}}^{1/2}
+C\nu\varepsilon_N(1+(s-t_0)^2)E_{f,{J_{\rm der}}}^{1/2}E_{w,{J_{\rm der}}}^{1/2}.
\end{split}
\end{equation}
The constants are independent of $J_x,\epsilon,{J_{\rm der}}$ and the parameter $c$.
For $\epsilon({T_1}-\tau)\le1$, the homogeneous bound
\eqref{eq:constructionUniformEstimate} and $\Theta_{\rm b}=o(1)$ imply
\[
\sup_{\tau\le s'\le s}E_{f,{J_{\rm der}}}(s')+\nu\int_\tau^sD_{f,{J_{\rm der}}}\,ds'
\le CE_{J_{\rm der}}(\tau;P_{J_x}f_\tau),\qquad \int_\tau^s (1+(s'-t_0)^2)\,ds'\le C{T_1}^3.
\]
Young's inequality in \eqref{eq:collisionDifferenceFiniteEnergy}, followed
{by the corresponding integrating factor,} yields
\begin{align*}
&\sup_{\tau\le s'\le s}E_{w,{J_{\rm der}}}(s')
+\nu\int_\tau^sD_{w,{J_{\rm der}}}\,ds'\\
&\quad\le C\varepsilon_N^2E_{J_{\rm der}}(\tau;P_{J_x}f_\tau)\\
&\qquad+C\nu\varepsilon_N{T_1}^3E_{J_{\rm der}}(\tau;P_{J_x}f_\tau)^{1/2}
\sup_{\tau\le s'\le s}E_{w,{J_{\rm der}}}(s')^{1/2}\\
&\quad\le\tfrac12\sup_{\tau\le s'\le s}E_{w,{J_{\rm der}}}(s')
+C\varepsilon_N^2[1+(\nu {T_1}^3)^2]E_{J_{\rm der}}(\tau;P_{J_x}f_\tau).
\end{align*}
Consequently,
\begin{equation}\label{eq:collisionDifferenceUniformBound}
\sup_{\tau\le s'\le s}E_{J_{\rm der}}(s';w^{J_x,\epsilon})
+\nu\int_\tau^sD_{J_{\rm der}}(s';w^{J_x,\epsilon})\,ds'
\le C\varepsilon_N^2(1+\nu {T_1}^3)^2E_{J_{\rm der}}(\tau;P_{J_x}f_\tau).
\end{equation}

First let $\epsilon\downarrow0$ at fixed $J_x$, and then let $J_x\to\infty$.
Step 5 of Lemma~\ref{lem:localizedConstruction} identifies the weak limits as
\[
f_\sharp^{J_x,\epsilon}\rightharpoonup\mathbb U_\sharp^\nu(\cdot,\tau)f_\tau,
\qquad
w^{J_x,\epsilon}\rightharpoonup
w=(\mathbb U_{\rm b}^\nu-\mathbb U_\mu^\nu)(\cdot,\tau)f_\tau.
\]
For fixed $J_{\rm der}$, distributional differentiation identifies
$\widehat Z^nw^{J_x,\epsilon}\rightharpoonup\widehat Z^nw$,
$0\le n\le J_{\rm der}$. Lower semicontinuity and
\eqref{eq:collisionDifferenceUniformBound} give
\[
\begin{aligned}
&\sup_{\tau\le s'\le s}E_{J_{\rm der}}(s';w)
+\nu\int_\tau^sD_{J_{\rm der}}(s';w)\,ds'\\
&\quad\le\liminf_{J_x\to\infty}\liminf_{\epsilon\downarrow0}
\left(\sup_{\tau\le s'\le s}E_{J_{\rm der}}(s';w^{J_x,\epsilon})
+\nu\int_\tau^sD_{J_{\rm der}}(s';w^{J_x,\epsilon})\,ds'\right)\\
&\quad\le C\varepsilon_N^2(1+\nu T_1^3)^2E_{J_{\rm der}}(\tau;f_\tau).
\end{aligned}
\]
For general data, choose finite-mode Schwartz $f_{\tau,h}$ with
$E_{J_{\rm der}}(\tau;f_{\tau,h}-f_\tau)\to0$.
The homogeneous estimate and \eqref{eq:localizedSourceCompletion} give
\[
\begin{aligned}
\sup_{\tau\le s'\le s}
E_{J_{\rm der}}(s';\mathbb U_\sharp^\nu(s',\tau)(f_{\tau,h}-f_{\tau,h'}))
&\le C E_{J_{\rm der}}(\tau;f_{\tau,h}-f_{\tau,h'}),\\
\mathbb U_\sharp^\nu(\cdot,\tau)f_{\tau,h}
&\longrightarrow\mathbb U_\sharp^\nu(\cdot,\tau)f_\tau
\quad\hbox{in }C([\tau,s];L^2_{x,v}).
\end{aligned}
\]
Thus the preceding bound holds for the prescribed datum. Finally,
\[
E_{J_{\rm der}}(s';w)\uparrow\mathcal E^{\rm ph}_{c,j}(s';w),\qquad
\int_\tau^sD_{J_{\rm der}}(s';w)\,ds'
\uparrow\int_\tau^s\mathcal D^{\rm ph}_{c,j}(s';w)\,ds'.
\]
Monotone convergence and square roots prove
\eqref{eq:exactCollisionEvolutionDifference}.
\par
For $f_\tau=\mathbf e_j(e^{i\eta v_1}V)$ and $c=j(\tau-t_0)+\eta$,
\begin{align*}
\widehat Z_{j,c}(\tau)^n(e^{i\eta v_1}V)
&=e^{i\eta v_1}\partial_1^nV,\\
\mathcal E^{\rm ph}_{c,j}(\tau;f_\tau)^{1/2}
&=\|V\|_{\mathscr A_\ell^{\lambda_{\rm v},\chi}}.
\end{align*}
Equation~\eqref{eq:centeredGtrace}, evaluated at time $s-t_0$ and
with parameter $c$ in mode $k$, and \eqref{eq:physicalSingleModeRelations} give
\begin{align*}
\left|\int\sqrt \mu  z_k(s)\,dv\right|
&\le Ce^{-C^{-1}|k(s-t_0)-c|^\chi-\varrho_0|k-j|^\chi}
\mathcal E^{\rm ph}_{c,j}(s;z)^{1/2}.
\end{align*}
For $z=\mathbb U_\sharp^\nu f_\tau$ use
$\mathcal E^{\rm ph}_{c,j}(s;z)^{1/2}\le C\|V\|_{\mathscr A_\ell^{\lambda_{\rm v},\chi}}$; for
{$z=(\mathbb U_{\rm b}^\nu-\mathbb U_\mu^\nu)f_\tau$} use
$\mathcal E^{\rm ph}_{c,j}(s;z)^{1/2}\le C\varepsilon_N(1+\nu {T_1}^3)\|V\|_{\mathscr A_\ell^{\lambda_{\rm v},\chi}}$.
Since $\nu {T_1}^3\to0$, these substitutions prove both stated traces.

Uniformity for $0<\lambda_{\rm v}\le\lambda_*$ follows from
\[
\sup_{0<\lambda_{\rm v}\le\lambda_*}
\sum_{h\ge1}(C\lambda_{\rm v})^h
\le\frac{C\lambda_*}{1-C\lambda_*},\qquad C\lambda_*<\frac12.
\]
The spatial radius is fixed. At velocity radius $\lambda_{\rm v}/N$,
the Gaussian pairing gives, for every integer $n\ge0$,
\[
\begin{aligned}
&|k(s-t_0)-c|^n
\left|\int\sqrt\mu[\mathbb U_\sharp^\nu(s,\tau)
\mathbf e_j(e^{i\eta v_1}V)]_k\,dv\right|\\
&\qquad\le Ce^{-\varrho_0|k-j|^\chi}(CN)^n(n!)^{1/\chi}
\|V\|_{\mathscr A_\ell^{\lambda_{\rm v}/N,\chi}}.
\end{aligned}
\]
If $(|k(s-t_0)-c|/(2CN))^\chi\ge2$, choose
$n=\lfloor(|k(s-t_0)-c|/(2CN))^\chi\rfloor$. Then
\[
\left(\frac{CN}{|k(s-t_0)-c|}\right)^n(n!)^{1/\chi}
\le\left(\frac{CNn^{1/\chi}}{|k(s-t_0)-c|}\right)^n
\le2^{-n}\le Ce^{-C^{-1}(|k(s-t_0)-c|/N)^\chi}.
\]
Otherwise use $n=0$. This proves \eqref{eq:exactCollisionVariableRadiusTrace}.
\par
\end{proof}

\begingroup
\begin{cor}[The density generated by the initial datum]\label{cor:initialDensityComparison}
Assume \eqref{eq:comparisonLogParameters} with $p>{\rm s}_0+7/2$,
and fix $r_{\rm x}\ge0$. For $H_{\rm in}$ in \eqref{eq:maininitial}, set
\[
\mathfrak g_k^{\rm C,\nu}(s)=\int\sqrt\mu\,
[\mathbb U_{\rm b}^\nu(s,t_0)(H_{\rm in}/\sqrt\mu)]_k\,dv,
\qquad
\mathfrak g_k^0(s)=\int\sqrt\mu\,
[\mathbb U^0(s,t_0)(H_{\rm in}/\sqrt\mu)]_k\,dv.
\]
Then, with $\varepsilon_{\rm coll}=\nu T_1(1+\zeta^2)$,
\begin{equation}\label{eq:exactCollisionInitialComparison}
\|\mathfrak g^{\rm C,\nu}-\mathfrak g^0\|_{\mathscr X_N^{r_{\rm x}}(I)}\le C_{r_{\rm x}}(\varepsilon_{\rm coll}+\varepsilon_N),
\qquad \|\mathfrak g^{\rm C,\nu}\|_{\mathscr X_N^{r_{\rm x}}(I)}\le C_{r_{\rm x}}.
\end{equation}
\end{cor}
\begin{proof}

For the initial trace split the difference as
\[
(\mathbb U_{\rm b}^\nu-\mathbb U^0)(H_{\rm in}/\sqrt \mu )
=(\mathbb U_{\rm b}^\nu-\mathbb U_\mu^\nu)(H_{\rm in}/\sqrt \mu )
+(\mathbb U_\mu^\nu-\mathbb U^0)(H_{\rm in}/\sqrt \mu ).
\]
The initial modes and frequency shifts are $j=\pm m_0$ and
$c=\eta=\pm(\zeta-m_0t_0)$. Their shifted-derivative Gaussian norms are
uniformly bounded. The weight inequality
\[
\big||k|s-\zeta\big|
\le|k(s-t_0)\mp(\zeta-m_0t_0)|+|k\mp m_0|t_0,
\qquad \lambda_N(s)t_0^\chi=o(1),
\]
allows the time weight to be absorbed by the phase and spatial decays
in \eqref{eq:exactCollisionTraceDifference}. Since $B_{|k|}^{-1}\le1$
and $m_0/N\le C$, Proposition~\ref{prop:centeredpropagation} and
Proposition~\ref{prop:exactCollisionPropagation} yield
\[
\begin{aligned}
\|\mathfrak g^{\rm C,\nu}-\mathfrak g^0\|_{\mathscr X_N^{r_{\rm x}}}
&\le C_{r_{\rm x}}\varepsilon_N
\sum_{\pm}\sum_{k\ne0}\langle k/N\rangle^{r_{\rm x}}
e^{-c|k\mp m_0|^\chi}+C_{r_{\rm x}}\nu T_1(1+\zeta^2)\\
&\le C_{r_{\rm x}}(\varepsilon_N+\varepsilon_{\rm coll}),\\
\|\mathfrak g^{\rm C,\nu}\|_{\mathscr X_N^{r_{\rm x}}}
&\le\|\mathfrak g^0\|_{\mathscr X_N^{r_{\rm x}}}
+C_{r_{\rm x}}(\varepsilon_N+\varepsilon_{\rm coll})\le C_{r_{\rm x}}.
\end{aligned}
\]
For the Maxwellian difference only $k=\pm m_0$ occurs.
This proves \eqref{eq:exactCollisionInitialComparison}.
\par
\end{proof}
\par\endgroup

\subsection{{\texorpdfstring{Treatment of $\Delta\mathbb T$ for Proposition~\ref{thm:exactCollisionEcho}}{Estimates for the density-operator difference and the exact density inverse}}}\label{sec:exactElectricResponse}\label{sec:densityComparisonProofs}\label{sec:exactBackgroundCollisions}

\begingroup
We estimate $\Delta\mathbb T$ to prove Proposition~\ref{thm:exactCollisionEcho}.
The exact background and both collision solution operators are already
available, and Corollary~\ref{cor:initialDensityComparison} controls the
initial-density difference. We split $\Delta\mathbb T$ into the three
kernels in \eqref{eq:directResponseDecomposition}, estimate them pointwise,
and then prove the weighted column bounds needed on
$\mathscr X_N^{r_{\rm x}}$. This gives the inverse of
\eqref{eq:exactCoefficientDensityEquation}, with both a norm bound and
an estimate of the coefficient at $(m_1,T_1)$. We then reconstruct
$H^{\rm C,\nu}$ from its density. Section~\ref{sec:threeDensityErrors}
combines the three errors and proves the terminal lower bound.
\par\endgroup
We begin with the coefficient $\mathcal F_m^{\rm b}(\tau)$ in
\eqref{eq:exactDensityVolterra}. For a density coefficient $\rho_m$, $m\ne0$,
Poisson's equation gives the normalized source coefficient
\begin{equation}\label{eq:exactCollisionElectricInput}
[\mathcal F_m^{\rm b}(\tau)]_j
=\frac{i}{m}\mu^{-1/2}\partial_{v_1}F_{j-m}^{\rm b,\nu}(\tau),
\qquad
-\mu^{-1/2}E[{\rho}]\partial_{v_1}F^{\rm b,\nu}
=\sum_{m\ne0}\mathcal F_m^{\rm b}{\rho}_m.
\end{equation}
For the comparison, split
{$F^{\rm b,\nu}=F^{\rm b}_{\rm app}+{F^{\rm b}_{\rm err}}$} and define
\begin{equation}\label{eq:directElectricInputs}
[\mathcal F_{{\rm app},m}(\tau)]_j
=\frac{i}{m}\mu^{-1/2}\partial_{v_1}(F^{\rm b}_{\rm app})_{j-m}(\tau),
\qquad
[{\mathcal F_{{\rm err},m}}(\tau)]_j
=\frac{i}{m}\left(\partial_{v_1}-\frac{v_1}{2}\right){f_{{\rm err},j-m}^{\rm b}}(\tau).
\end{equation}
Then $\mathcal F_m^{\rm b}=\mathcal F_{{\rm app},m}+{\mathcal F_{{\rm err},m}}$,
and \eqref{eq:overviewDensityResponseDifference} becomes
\[
(\Delta\mathbb T{\rho})_k(s)
=\sum_{m\ne0}\int_{t_0}^s{\rho}_m(\tau)\int\sqrt \mu \,
\big[\mathbb U_{\rm b}^\nu(s,\tau)\mathcal F_m^{\rm b}(\tau)
-\mathbb U^0(s,\tau)\mathcal F_{{\rm app},m}(\tau)\big]_k\,dv\,d\tau.
\]
\begingroup
{The approximate density integral equals}
$(\mathbb K+\mathbb P_\downarrow+\mathbb P_\uparrow){\rho}$ on both
signs of the spatial modes. This formula will be used to estimate the
difference of density integral operators, while \eqref{eq:exactDensityVolterra} remains
the equation to be solved. Background-field transport is still absent.
\par\endgroup
\begingroup
The bound for $g^{\rm b}_{\rm err}$ at the decreasing radius was
proved in Proposition~\ref{prop:decreasingRadiusBackground}.
It will be used for the third kernel below; the first two kernels use the
solution-operator comparisons of Section~\ref{sec:gainenergy}.
\par\endgroup

\begingroup
The exact difference of the two density kernels is
\begin{equation}\label{eq:directResponseDecomposition}
\begin{aligned}
(\Delta\mathbb T{\rho})_k(s)
&=\sum_{j=1}^3\sum_{m\ne0}\int_{t_0}^s
\Delta\mathcal K^{(j)}_{k,m}(s,\tau){\rho}_m(\tau)\,d\tau,\\
\Delta\mathcal K^{(1)}_{k,m}(s,\tau)
&=\int\sqrt \mu \,\big[(\mathbb U_\mu^\nu-\mathbb U^0)(s,\tau)
\mathcal F_{{\rm app},m}(\tau)\big]_k\,dv,\\
\Delta\mathcal K^{(2)}_{k,m}(s,\tau)
&=\int\sqrt \mu \,\big[(\mathbb U_{\rm b}^\nu-\mathbb U_\mu^\nu)(s,\tau)
\mathcal F_{{\rm app},m}(\tau)\big]_k\,dv,\\
\Delta\mathcal K^{(3)}_{k,m}(s,\tau)
&=\int\sqrt \mu \,\big[\mathbb U_{\rm b}^\nu(s,\tau)
{\mathcal F_{{\rm err},m}}(\tau)\big]_k\,dv.
\end{aligned}
\end{equation}
The first term changes free transport to Maxwellian collision evolution;
the second changes the collision coefficient from $\mu$ to $F^{\rm b,\nu}$;
the third replaces $\partial_{v_1}F^{\rm b}_{\rm app}$ by
$\partial_{v_1}F^{\rm b,\nu}$. We estimate the kernels first and then
integrate them in the density norm.
\begin{lem}\label{lem:exactElectricInputs}
Under \eqref{eq:comparisonLogParameters}, let $m\ne0$, $l\in\mathbb Z$
and $t=\tau-t_0\in[0,T_1-t_0]$.
Then $\mathcal F_m^{\rm b}(\tau)$ satisfies
\begin{equation}\label{eq:completeElectricProfile}
{
[\mathcal F_m^{\rm b}(\tau)]_{m+l} ={}-\mathbf1_{l=0}\frac{i}{m}v_1\sqrt \mu  +\frac{i}{m}e^{-iltv_1} (\partial_{v_1}-v_1/2-ilt)g_l^{\rm b,\nu}(t,v).
}
\end{equation}
For a sufficiently small fixed radius and $J_v$ large enough for the fixed weight,
\begin{equation}\label{eq:completeElectricInputBound}
\|(\partial_{v_1}-v_1/2-ilt)g_l^{\rm b,\nu}(t)\|_{\mathscr A_\ell^{\lambda_{\rm v},\chi}}
\le C\varepsilon_N(1+|l|\tau)e^{-\varrho_{\rm b}|l|^\chi}.
\end{equation}
For every function $V=V(v)$ independent of $x_1$ and integer $n\ge0$,
\[
\widehat Z_{m+l,\,m(\tau-t_0)}(\tau)^n
\big(e^{-il(\tau-t_0)v_1}V(v)\big)
=e^{-il(\tau-t_0)v_1}\partial_{v_1}^nV(v).
\]
The case $l=0$ also applies to the Gaussian term in \eqref{eq:completeElectricProfile}.

\end{lem}
\begin{proof}
At $t=\tau-t_0$, direct differentiation gives
\begin{align*}
F_l^{\rm b,\nu}&=\mathbf1_{l=0}\mu +\sqrt \mu  e^{-iltv_1}g_l^{\rm b,\nu},\\
\mu^{-1/2}\partial_{v_1}F_l^{\rm b,\nu}
&=-\mathbf1_{l=0}v_1\sqrt \mu 
+e^{-iltv_1}(\partial_{v_1}-v_1/2-ilt)g_l^{\rm b,\nu}.
\end{align*}
Multiplication by $i/m$ proves \eqref{eq:completeElectricProfile}.
At a fixed reduced radius, Lemmas~\ref{lem:weightedGevreyInterpolation}
and \ref{lem:physicalCoefficientDerivatives} yield
\begin{align*}
\|g_l^{\rm b,\nu}\|_{\mathscr A_\ell^{\lambda_{\rm v},\chi}}
+\|\partial_1g_l^{\rm b,\nu}\|_{\mathscr A_\ell^{\lambda_{\rm v},\chi}}
+\|v_1g_l^{\rm b,\nu}\|_{\mathscr A_\ell^{\lambda_{\rm v},\chi}}
&\le C\varepsilon_Ne^{-\varrho_{\rm b}|l|^\chi},\\
\|(\partial_1-v_1/2-ilt)g_l^{\rm b,\nu}\|_{\mathscr A_\ell^{\lambda_{\rm v},\chi}}
&\le C\varepsilon_N(1+|l|\tau)e^{-\varrho_{\rm b}|l|^\chi}.
\end{align*}
For $j=m+l$ and any function $V=V(v)$ independent of $x_1$,
\[
\widehat Z_{j,mt}(\tau)(e^{-iltv_1}V)
=e^{-iltv_1}\big(\partial_1+i(jt-mt-lt)\big)V
=e^{-iltv_1}\partial_1V.
\]
Iteration gives the stated identity for every $n\ge0$, including $l=0$.

\end{proof}
\begin{cor}[The three kernel errors]\label{cor:threeKernelErrors}
Under the hypotheses of Corollary~\ref{cor:initialDensityComparison},
let $k,m\ne0$ and $t_0\le\tau\le s\le T_1$. The kernels in
\eqref{eq:directResponseDecomposition} satisfy the following estimates.
The first vanishes unless $k=m+l$ for $l=0,\pm1$, and then
\begin{equation}\label{eq:maxwellianKernel}
|\Delta\mathcal K^{(1)}_{k,m}(s,\tau)|
\le \frac{C\nu {T_1}^3}{|m|}
\big(\mathbf1_{l=0}+a(1+{T_1})\mathbf1_{|l|=1}\big)e^{-c|ks-m\tau|^\chi}.
\end{equation}
\begin{equation}\label{eq:directKernel}
|\Delta\mathcal K^{(2)}_{k,m}(s,\tau)|
\le \frac{C\varepsilon_N}{|m|}
\sum_{l=0,\pm1}
\big(\mathbf1_{l=0}+a(1+\tau)\mathbf1_{|l|=1}\big)
e^{-\varrho_0|k-m-l|^\chi}e^{-c|k(s-t_0)-m(\tau-t_0)+lt_0|^\chi}.
\end{equation}
\begin{equation}\label{eq:backgroundRemainderKernel}
\begin{split}
|\Delta\mathcal K^{(3)}_{k,m}(s,\tau)|
\le{}&\frac{Ca\nu T_1^3}{|m|}
\sum_{l\in\mathbb Z}(N+|l|\tau)e^{-c(|l|/N)^\chi}
 e^{-\varrho_0|k-m-l|^\chi}\\
&\hspace{15mm}\times e^{-c_2(|k(s-t_0)-m(\tau-t_0)|/N)^\chi}.
\end{split}
\end{equation}
The constants are independent of $k,m,s,\tau,N,\nu$.
\end{cor}
\begin{proof}
For $l=0,\pm1$, the explicit background gives
\[
[\mathcal F_{{\rm app},m}(\tau)]_{m+l}
=\begin{cases}
-\dfrac{i}{m}v_1\sqrt\mu,&l=0,\\[1mm]
\dfrac{ia}{2m}e^{-il\tau v_1}(-v_1-il\tau)\sqrt\mu,&|l|=1.
\end{cases}
\]
The Maxwellian operator acts diagonally on the Fourier transform in $x$. For
$k=m+l$, the estimate \eqref{eq:centeredpacketdifference}, with
$\eta=-l\tau$, gives a factor $e^{-c|k(s-\tau)+l\tau|^\chi}$.
Moreover,
\[
\begin{aligned}
k(s-\tau)+l\tau&=ks-m\tau,\\
\nu(s-\tau)[1+l^2\tau^2+k^2(s-\tau)^2]
&\le C\nu T_1[1+T_1^2+|ks-m\tau|^2],\\
[1+|ks-m\tau|^2]e^{-c|ks-m\tau|^\chi}
&\le C e^{-c|ks-m\tau|^\chi/2}.
\end{aligned}
\]
Multiplying by the coefficients of the three functions displayed above proves
\eqref{eq:maxwellianKernel}.

For the second kernel, apply \eqref{eq:exactCollisionTraceDifference}
to $[\mathcal F_{{\rm app},m}(\tau)]_{m+l}$ for $l=0,\pm1$. At mode $j=m+l$,
\[
c=j(\tau-t_0)-l\tau=m(\tau-t_0)-lt_0,
\qquad k-j=k-m-l.
\]
The corresponding $\mathscr A_\ell^{\lambda_{\rm v},\chi}$ norms are bounded by $C/|m|$ for $l=0$ and
$Ca(1+\tau)/|m|$ for $|l|=1$. Substitution in
\eqref{eq:exactCollisionTraceDifference} proves \eqref{eq:directKernel}.

For the third kernel, \eqref{eq:directElectricInputs} yields
\[
[\mathcal F_{{\rm err},m}(\tau)]_{m+l}
=\frac{i}{m}e^{-il(\tau-t_0)v_1}
 (\partial_{v_1}-v_1/2-il(\tau-t_0))g_{{\rm err},l}^{\rm b}(\tau-t_0).
\]
Here
\[
c=(m+l)(\tau-t_0)-l(\tau-t_0)=m(\tau-t_0).
\]
Using \eqref{eq:decreasingRadiusElectricInput} in
\eqref{eq:exactCollisionVariableRadiusTrace} and summing over $l$ gives
\[
\begin{aligned}
|\Delta\mathcal K^{(3)}_{k,m}(s,\tau)|
&\le\frac C{|m|}\sum_l e^{-\varrho_0|k-m-l|^\chi}
 e^{-c_2(|k(s-t_0)-m(\tau-t_0)|/N)^\chi}\\
&\qquad\times\| (\partial_{v_1}-v_1/2-il(\tau-t_0))
 g_{{\rm err},l}^{\rm b}(\tau-t_0)\|_{\mathscr A_\ell^{\lambda_{\rm v}/N,\chi}}\\
&\le\frac{Ca\nu T_1^3}{|m|}\sum_l(N+|l|\tau)
 e^{-c(|l|/N)^\chi}e^{-\varrho_0|k-m-l|^\chi}
 e^{-c_2(|k(s-t_0)-m(\tau-t_0)|/N)^\chi}.
\end{aligned}
\]
This proves \eqref{eq:backgroundRemainderKernel}.
\end{proof}
\begin{lem}\label{lem:generalGainColumn}
{Fix $r_{\rm x}\ge0$, $1/3<\chi<1$, ${\bf b}_2>0$,
$0<q_{\rm g}<1$ and $C_{\rm g}>0$. Let $N\ge1$, $t_0\ge0$,
$T_{\rm a}=t_0+N^4<T_1$, $I_+=[T_{\rm a},T_1]$, and assume
$a>0$ and $a(1+T_1)\le C_0N^2$.
Suppose that $\mathfrak w_n>0$ satisfies}
\begin{equation}\label{eq:generalGainHypotheses}
\begin{aligned}
\mathfrak w_m/\mathfrak w_k&\le1 &&(1\le k\le m),\\
\mathfrak w_m/\mathfrak w_k&\le q_{\rm g}^{m-k} &&(1\le k\le m\le {\bf b}_2N),\\
\mathfrak w_m/\mathfrak w_k&\le e^{C_{\rm g}N} &&(k,m\ge1).
\end{aligned}
\end{equation}
Let $\varepsilon_{\rm ker}\ge0$, and let a kernel on $T_{\rm a}\le\tau\le s\le {T_1}$ obey,
for fixed $c,\varrho_0,c_2>0$,
\begin{equation}\label{eq:remainderKernel}
\begin{split}
|\mathcal K_{k,m}(s,\tau)|
&\le\frac{Ca\varepsilon_{\rm ker}}{|m|}\sum_{l\in\mathbb Z}(N+|l|\tau)
e^{-c(|l|/N)^\chi}e^{-\varrho_0|k-m-l|^\chi}e^{-c_2(|k(s-t_0)-m(\tau-t_0)|/N)^\chi},\\
k,m\ne0.
\end{split}
\end{equation}
{Then, for sufficiently large $N$,}
\begin{equation}\label{eq:generalGainColumn}
\sup_{m\ne0}\sum_{k\ne0}
\frac{\langle k/N\rangle^{r_{\rm x}}\mathfrak w_{|m|}}{\langle m/N\rangle^{r_{\rm x}}\mathfrak w_{|k|}}
\sup_{s\in I_+}\int_{T_{\rm a}}^s|\mathcal K_{k,m}(s,\tau)|\,d\tau
\le C_{{r_{\rm x}},\rm g}N^3\varepsilon_{\rm ker}.
\end{equation}
Here $C_{{r_{\rm x}},\rm g}$ and the lower threshold for $N$ may depend on the fixed
constants in \eqref{eq:generalGainHypotheses}--\eqref{eq:remainderKernel}.

\end{lem}
\begin{proof}

Put $n=k-m-l$. The polynomial ratio is bounded by
$C_{r_{\rm x}}\langle l/N\rangle^{r_{\rm x}}\langle n\rangle^{r_{\rm x}}$ and is absorbed by fixed
fractions of the two spatial exponentials in \eqref{eq:remainderKernel}.
For each fixed $k$, $\partial_\tau\big(k(s-t_0)-m(\tau-t_0)\big)=-m$, so on any subinterval
with left endpoint $\tau_*$ the change of variable $u=\big(k(s-t_0)-m(\tau-t_0)\big)/N$ gives
\begin{equation}\label{eq:phaseIntegral}
\sup_s\int_{\tau_*}^s e^{-c'(|k(s-t_0)-m(\tau-t_0)|/N)^\chi}\,d\tau
\le\frac{N}{|m|}\int_{\mathbb R}e^{-c'|u|^\chi}\,du
\le CN/|m|.
\end{equation}
Together with
\[
\sum_ne^{-r'|n|^\chi}\le C,\qquad
\sum_l(N+|l|T')e^{-c'(|l|/N)^\chi}\le CN^2(1+T'),
\]
this bounds a column whose gain ratio is at most one, at times at most $T'$, by
\begin{equation}\label{eq:baselineColumn}
\frac{C_{r_{\rm x}}a\varepsilon_{\rm ker} N^3(1+T')}{m^2}.
\end{equation}
The time supremum in \eqref{eq:phaseIntegral} is taken separately for
each spatial index $k$, as required in \eqref{eq:generalGainColumn}.

First suppose that $km>0$ and $1\le|k|\le|m|$.
For $|m|\ge {\bf b}_2N$, the first gain hypothesis and
\eqref{eq:baselineColumn} give
\[
\frac{C_{r_{\rm x}}a\varepsilon_{\rm ker} N^3(1+{T_1})}{m^2}
\le C_{{r_{\rm x}},{\bf b}_2}N^3\varepsilon_{\rm ker}.
\]
For $|m|<{\bf b}_2N$, put $h=|m|-|k|$. Then
\[
\frac{\mathfrak w_{|m|}}{\mathfrak w_{|k|}}\le q_{\rm g}^h,\qquad
|l|\le h+|n|,\qquad
\frac{\langle k/N\rangle^{r_{\rm x}}}{\langle m/N\rangle^{r_{\rm x}}}\le1.
\]
Discarding the remaining $l$ decay yields
\begin{equation}\label{eq:smallColumnSum}
\sum_{h\ge0}q_{\rm g}^h\sum_{n\in\mathbb Z}
\big(N+(h+|n|){T_1}\big)e^{-r'|n|^\chi}
\le C_{q_{\rm g}}(N+{T_1}).
\end{equation}
By \eqref{eq:phaseIntegral}, this part of the column is at most
\[
\frac{C_{{r_{\rm x}},\rm g}a\varepsilon_{\rm ker} N(N+{T_1})}{m^2}\le C_{{r_{\rm x}},\rm g}N^3\varepsilon_{\rm ker}.
\]
For the remaining configurations, $\tau-t_0\ge N^4$ and
\[
|k(s-t_0)-m(\tau-t_0)|=
\begin{cases}
|k|(s-\tau)+(|k|-|m|)(\tau-t_0),&km>0,\ |k|>|m|,\\
|k|(s-t_0)+|m|(\tau-t_0),&km<0,
\end{cases}
\quad\ge N^4.
\]
Reserving a fraction of the phase gives
\[
\frac{\mathfrak w_{|m|}}{\mathfrak w_{|k|}}e^{-c'(|k(s-t_0)-m(\tau-t_0)|/N)^\chi}
\le e^{C_{\rm g}N-c'N^{3\chi}}\le N^{-2}.
\]
Apply \eqref{eq:baselineColumn} to the remaining positive spatial
and phase factors:
\[
N^{-2}\frac{C_{r_{\rm x}}a\varepsilon_{\rm ker} N^3(1+{T_1})}{m^2}
\le C_{r_{\rm x}}N^3\varepsilon_{\rm ker}.
\]
These three estimates prove \eqref{eq:generalGainColumn}.

\end{proof}
\begin{rem}[The weights used in the applications]\label{rem:gainWeightsApplications}
The chain weights $\mathfrak w_k=B_k$ satisfy the hypotheses of
Lemma~\ref{lem:generalGainColumn}.
Fix $0<{\bf b}_2<{\bf b}_0/4$. The gain $B$ satisfies
\begin{equation}\label{eq:smallGain}
\begin{aligned}
\frac{B_m}{B_k}
&=\prod_{j=k+1}^m\frac{j^3}{\Lambda N^3}
\le ({\bf b}_2/{\bf b}_0)^{3(m-k)}, &&1\le k\le m\le {\bf b}_2N,\\
B_k&\ge e^{c_{{\bf b}_2}N}, &&1\le k\le2{\bf b}_2N,
\end{aligned}
\end{equation}
The first base is less than one. For the second bound,
\[
\log B_{\lfloor2{\bf b}_2N\rfloor}
\ge N\int_{2{\bf b}_2}^{{\bf b}_0}\log\frac{\Lambda}{y^3}\,dy-O(1)
\ge c_{{\bf b}_2}N,
\]
and $B_k$ is nonincreasing. The same monotonicity and the
bound $\widetilde{B}\le e^{CN}$ verify the other two hypotheses in
\eqref{eq:generalGainHypotheses}.

For the general finite-regularity density estimate in
Section~\ref{sec:fullLinearClosure}, the same lemma will be applied to
$\mathfrak w_k=G_k$. Those weights are defined from the kernel bound in the proof of Proposition~\ref{prop:coarseDensityInverse} below, where their three hypotheses are checked before use.
\end{rem}
\begin{lem}\label{prop:exactCollisionResponse}
Assume \eqref{eq:comparisonLogParameters} with $p>{\rm s}_0+7/2$.
On $\mathscr X_N^{r_{\rm x}}$, for every fixed ${r_{\rm x}}\ge0$,
\begin{equation}\label{eq:electricResponseDifferenceBound}
\|\Delta\mathbb T\|_{\mathscr X_N^{r_{\rm x}}(I)\to\mathscr X_N^{r_{\rm x}}(I)}\le C_{r_{\rm x}}\big(
N^2\nu {T_1}^3+N^2\varepsilon_N+N^3\nu T_1^3\big).
\end{equation}
The three terms on the right bound the operators given by the three
kernel summands of \eqref{eq:directResponseDecomposition}, in that order.

\end{lem}
\begin{proof}
We estimate the weighted column indexed by $m$, taking the supremum over $s$ separately for each $k$. On $I_-\times I_-$ the gain ratio is
$1$; on $I_+\times I_-$ it is $B_{|k|}^{-1}\le1$; on
$I_+\times I_+$ it is $B_{|m|}/B_{|k|}$.
\par

For the first kernel set $z=ks-m\tau$ and apply \eqref{eq:maxwellianKernel}. For $k=m+l\ne0$, $l=0,\pm1$, the required weight bounds are
\[
\begin{aligned}
\frac{w_{|k|,N}(s)}{w_{|m|,N}(\tau)}
&\le e^{\lambda_N(s)|z|^\chi},\qquad
\frac{\langle k/N\rangle^{r_{\rm x}}}{\langle m/N\rangle^{r_{\rm x}}}
\le C_{r_{\rm x}},\qquad \frac{B_{|m|}}{B_{|k|}}\le1\quad (|k|\le|m|),\\
\frac{B_{|m|}}{B_{|k|}}e^{-c|z|^\chi/4}
&\le CN^3e^{-ct_0^\chi/4}\le C
\quad (|k|=|m|+1),
\end{aligned}
\]
where the last line uses $|z|\ge s\ge t_0$.
Consequently, on each nonempty time block,
\[
\begin{aligned}
&\sum_{\substack{k\ne0\\|k-m|\le1}}
\frac{\langle k/N\rangle^{r_{\rm x}}}{\langle m/N\rangle^{r_{\rm x}}}
\max\left\{1,\frac{B_{|m|}}{B_{|k|}}\right\}
\sup_s\int_{t_0}^s
\frac{w_{|k|,N}(s)}{w_{|m|,N}(\tau)}
|\Delta\mathcal K^{(1)}_{k,m}(s,\tau)|\,d\tau\\
&\quad\le\frac{C_{r_{\rm x}}\nu T_1^3[1+a(1+T_1)]}{|m|}
\sup_{s,k}\int_{t_0}^s e^{-c'|ks-m\tau|^\chi}\,d\tau\\
&\quad\le\frac{C_{r_{\rm x}}\nu T_1^3[1+a(1+T_1)]}{m^2}
\le C_{r_{\rm x}}N^2\nu T_1^3.
\end{aligned}
\]
\par

For the second kernel use \eqref{eq:directKernel} and write $n=k-m-l$. Since $ks-m\tau=k(s-t_0)-m(\tau-t_0)+lt_0+nt_0$,
\[
\begin{aligned}
\frac{w_{|k|,N}(s)}{w_{|m|,N}(\tau)}
&\le\exp\left(\lambda_N(s)|k(s-t_0)-m(\tau-t_0)+lt_0|^\chi
+\lambda_N(s)t_0^\chi|n|^\chi\right),\\
\lambda_N(s)&\le c/4,\qquad
\lambda_N(s)t_0^\chi\le\varrho_0/4,\\
\frac{\langle k/N\rangle^{r_{\rm x}}}{\langle m/N\rangle^{r_{\rm x}}}
e^{-\varrho_0|n|^\chi/4}
&\le C_{r_{\rm x}}\langle n\rangle^{r_{\rm x}}
e^{-\varrho_0|n|^\chi/4}\le C_{r_{\rm x}}.
\end{aligned}
\]
The first two bounds leave positive fractions of both kernel decays.
\par
It remains to check the possibly large gain ratio. If
$s\ge\tau\ge T_{\rm a}$ and $|k|>|m|$, then $|k(s-t_0)-m(\tau-t_0)|\ge N^4$.
For equal signs this follows from
$|k(s-t_0)-m(\tau-t_0)|=|k|(s-\tau)+(|k|-|m|)(\tau-t_0)$; for opposite signs both terms have the same sign in $\big(k(s-t_0)-m(\tau-t_0)\big)$. Since $|l|\le1$ and $t_0=O(N)$,
$|k(s-t_0)-m(\tau-t_0)+lt_0|\ge N^4/2$ for large $N$. A reserved phase factor consequently satisfies
\[
\frac{B_{|m|}}{B_{|k|}}
e^{-c'|k(s-t_0)-m(\tau-t_0)+lt_0|^\chi}
\le e^{CN-c'N^{4\chi}/2^\chi}\le1 .
\]
After these weight bounds, each column is at most
\[
\begin{aligned}
&\frac{C_{r_{\rm x}}\varepsilon_N[1+a(1+T_1)]}{|m|}
\sum_{l=0,\pm1}\sum_n e^{-c'|n|^\chi}
\sup_s\int_{t_0}^s
e^{-c'|k(s-t_0)-m(\tau-t_0)+lt_0|^\chi}\,d\tau\\
&\quad\le\frac{C_{r_{\rm x}}\varepsilon_N[1+a(1+T_1)]}{m^2}
\le C_{r_{\rm x}}N^2\varepsilon_N.
\end{aligned}
\]
This gives the second contribution to \eqref{eq:electricResponseDifferenceBound}.
\par

For the third kernel apply \eqref{eq:backgroundRemainderKernel}.
First insert the frequency weight. Since
$ks-m\tau=\big(k(s-t_0)-m(\tau-t_0)\big)+(l+n)t_0$, the time weight ratio is at most
\[
\exp\!\left(\lambda_N(s)|k(s-t_0)-m(\tau-t_0)|^\chi+
\lambda_N(s)t_0^\chi(|l|^\chi+|n|^\chi)\right).
\]
Here $\lambda_N=O(\vartheta N^{-2\chi})$ and
$\lambda_Nt_0^\chi=O(\vartheta N^{-\chi})$. \begingroup
For sufficiently small fixed $\vartheta$ and then sufficiently large $N$,
\[
\lambda_N(s)N^\chi\le\frac{c_2}{2},\qquad
\lambda_N(s)t_0^\chi N^\chi\le\frac c2,\qquad
\lambda_N(s)t_0^\chi\le\frac{\varrho_0}{2}.
\]
Multiplication of \eqref{eq:backgroundRemainderKernel} by the time-weight ratio therefore gives
\par\endgroup
\begingroup
\begin{equation}\label{eq:thirdKernelTimeWeightBound}
\begin{aligned}
&\frac{w_{|k|,N}(s)}{w_{|m|,N}(\tau)}
|\Delta\mathcal K^{(3)}_{k,m}(s,\tau)|\\
&\quad\le\frac{Ca\nu T_1^3}{|m|}
\sum_{l\in\mathbb Z}(N+|l|\tau)
 e^{-\frac c2(|l|/N)^\chi}
 e^{-\frac{\varrho_0}{2}|k-m-l|^\chi}
 e^{-\frac{c_2}{2}(|k(s-t_0)-m(\tau-t_0)|/N)^\chi},\\
&\hspace{30mm} t_0\le\tau\le s\le T_1,\qquad k,m\ne0.
\end{aligned}
\end{equation}
\par\endgroup

Thus the weighted kernel satisfies \eqref{eq:remainderKernel} with
$\varepsilon_{\rm ker}=\nu T_1^3$ and smaller fixed
exponents. The weight hypotheses with $\mathfrak w=B$ were verified in
Remark~\ref{rem:gainWeightsApplications}.
On $I_-\times I_-$, no gain ratio occurs. Equation
\eqref{eq:baselineColumn}, with $T'=T_{\rm a}$, gives
\[
\frac{C_{r_{\rm x}}a\nu T_1^3N^3(1+T_{\rm a})}{m^2}
\le C_{r_{\rm x}}N^3\nu T_1^3,
\qquad a(1+T_{\rm a})=o(1).
\]
For $\tau\in I_-$ and $s\in I_+$, the gain ratio is
$B_{|k|}^{-1}\le1$. If $|m|\ge {\bf b}_2N$, use
\eqref{eq:baselineColumn} with $T'={T_1}$ and $m^{-2}\le C_{{\bf b}_2}N^{-2}$.
If $|m|<{\bf b}_2N$ and $|k|\le2{\bf b}_2N$, use \eqref{eq:smallGain}:
\[
\frac{C_{r_{\rm x}}a\nu T_1^3N^3(1+{T_1})}{m^2}
e^{-c_{{\bf b}_2}N}\le C_{r_{\rm x}}N^3\nu T_1^3.
\]
For $|m|<{\bf b}_2N$ and $|k|>2{\bf b}_2N$,
\[
|k(s-t_0)-m(\tau-t_0)|\ge |k|(s-t_0)-|m|(\tau-t_0)
\ge(|k|-|m|)N^4\ge {\bf b}_2N^5,
\]
so a reserved phase fraction is bounded by
$e^{-c_{{\bf b}_2}N^{4\chi}}\le N^{-2}$. Equation
\eqref{eq:baselineColumn} again gives $C_{r_{\rm x}}N^3\nu T_1^3$.
Finally, on $I_+\times I_+$ apply \eqref{eq:generalGainColumn}
with $\mathfrak w=B$. The remaining time block is empty by causality.
Adding the three blocks proves the third contribution to \eqref{eq:electricResponseDifferenceBound}.
The column sums are absolutely convergent, so finite continuous sources
extend to the full density space; dominated convergence preserves
time continuity. No finite-support condition on the density is used.

\end{proof}
\par\endgroup

\begingroup
\begin{prop}\label{prop:exactDensityInverse}
{Assume \eqref{eq:comparisonLogParameters} with $p>{\rm s}_0+7/2$ and \eqref{eq:parameters}.
Fix $r_{\rm x}\ge0$ and let $\mathfrak g\in\mathscr X_N^{r_{\rm x}}$,
with the density space defined in Section~\ref{sec:gainDensityInverse}.
Then \eqref{eq:exactCoefficientDensityEquation} has a unique solution
${\rho}\in\mathscr X_N^{r_{\rm x}}$. The solution map
$\mathbb V^{\rm C,\nu}:\mathfrak g\mapsto{\rho}$ satisfies}
\begin{equation}\label{eq:logExactCollisionInverse}
\|\mathbb V^{\rm C,\nu}\|_{\mathscr X_N^{r_{\rm x}}\to\mathscr X_N^{r_{\rm x}}}\le C_{r_{\rm x}}N,
\qquad \|\rho[H^{\rm C,\nu}]\|_{\mathscr X_N^{r_{\rm x}}}\le C_{r_{\rm x}}N.
\end{equation}
Evaluation of this solution map at the terminal mode has the stronger estimate
\begin{equation}\label{eq:exactTerminalInverse}
\|\mathcal O_N\mathbb V^{\rm C,\nu}\|_{\mathscr X_N^{r_{\rm x}}\to\mathbb C}\le C_{r_{\rm x}}.
\end{equation}
For $\mathfrak g=\mathfrak g^{\rm C,\nu}$, the density reconstructs a
unique solution $H^{\rm C,\nu}$ of \eqref{eq:exactCollisionModel} in every
fixed polynomially weighted finite-energy class, and
$\rho[H^{\rm C,\nu}]=\mathbb V^{\rm C,\nu}\mathfrak g^{\rm C,\nu}$.
\end{prop}
\begin{proof}
\emph{1. Density construction.}
{We solve the equation with prescribed $\mathfrak g$
\eqref{eq:exactCoefficientDensityEquation}.}
For $\Delta\mathbb T$ in \eqref{eq:overviewDensityResponseDifference},
Lemma~\ref{prop:exactCollisionResponse} and
Proposition~\ref{prop:approximateDensityInverse} give
\[
\|\mathbb V_{\rm e}^0\Delta\mathbb T\|
\le C_{r_{\rm x}}\big(N^3\nu {T_1}^3+N^3\varepsilon_N+
N^4\nu T_1^3\big)=o(1).
\]
Indeed $N^4\nu {T_1}^3\mathrel{{\sim}}[\log(1/\nu)]^{3{\rm s}_0+10-3p}$,
whose exponent is negative for $p>{\rm s}_0+7/2$, and
$N^3\varepsilon_N=\nu^{1/3-o(1)}$. The other terms are smaller.
The convergent operator series therefore defines
\begin{equation}\label{eq:exactCollisionDensityInverse}
\begin{aligned}
\mathbb V^{\rm C,\nu}
&=({\mathrm{Id}}-\mathbb V_{\rm e}^0\Delta\mathbb T)^{-1}\mathbb V_{\rm e}^0,
&\|\mathbb V^{\rm C,\nu}\|&\le C_{r_{\rm x}}N,\\
\rho^{\rm C,\nu}&=\mathbb V^{\rm C,\nu}\mathfrak g^{\rm C,\nu},
&\|\rho^{\rm C,\nu}\|_{\mathscr X_N^{r_{\rm x}}}&\le C_{r_{\rm x}}N.
\end{aligned}
\end{equation}
The last bound uses \eqref{eq:exactCollisionInitialComparison}.
Indeed,
\[
{
{\mathrm{Id}}-\mathbb K-\mathbb P_\downarrow-\mathbb P_\uparrow-\Delta\mathbb T =({\mathrm{Id}}-\mathbb K-\mathbb P_\downarrow-\mathbb P_\uparrow) ({\mathrm{Id}}-\mathbb V_{\rm e}^0\Delta\mathbb T),
}
\]
{Thus, for every $\mathfrak g\in\mathscr X_N^{r_{\rm x}}$,
${\rho}=\mathbb V^{\rm C,\nu}\mathfrak g$ is the unique solution of
\eqref{eq:exactCoefficientDensityEquation}. Taking
$\mathfrak g=\mathfrak g^{\rm C,\nu}$ gives $\rho^{\rm C,\nu}$.}
The resolvent identity and Lemma~\ref{lem:terminalObservation} give
\[
\begin{aligned}
\mathbb V^{\rm C,\nu}
&=\mathbb V_{\rm e}^0+\mathbb V_{\rm e}^0\Delta\mathbb T\mathbb V^{\rm C,\nu},\\
\|\mathcal O_N\mathbb V^{\rm C,\nu}\|
&\le\|\mathcal O_N\mathbb V_{\rm e}^0\|
\big(1+\|\Delta\mathbb T\|\,\|\mathbb V^{\rm C,\nu}\|\big)
\le C_{r_{\rm x}}(1+N\|\Delta\mathbb T\|)\le C_{r_{\rm x}}.
\end{aligned}
\]
This is \eqref{eq:exactTerminalInverse}.
\par

\emph{2. Distribution reconstruction.}
Fix $r,\ell$ and take ${r_{\rm x}}\ge r+1$. The electric source is
\[
\sum_{m\ne0}\mathcal F_m^{\rm b}(\tau)\rho_m^{\rm C,\nu}(\tau)
=-\mu^{-1/2}E[\rho^{\rm C,\nu}]\partial_{v_1}F^{\rm b,\nu}.
\]
The coefficient bounds in Lemma~\ref{lem:fullBackgroundCoefficients} and the spatial Leibniz rule give, for fixed $N,{T_1}$,
\[
\begin{split}
\sum_{m\ne0}\int_I
\|\mathcal F_m^{\rm b}(\tau)\rho_m^{\rm C,\nu}(\tau)\|_{H^r_xL^2_\ell}\,d\tau
&\le C_{r,\ell,N,{T_1}}\sum_{m\ne0}\langle m\rangle^r\sup_I|\rho_m^{\rm C,\nu}|\\
&\le C_{r,\ell,N,{T_1},{r_{\rm x}}}\|\rho^{\rm C,\nu}\|_{\mathscr X_N^{r_{\rm x}}}<\infty.
\end{split}
\]
The second inequality follows from \eqref{eq:echoDensityEmbeddings} and $r_{\rm x}\ge r$. {The constants in this reconstruction bound are used only} for finite-interval reconstruction.

{The solution operator in Proposition~\ref{prop:exactCollisionPropagation} gives the reconstruction}
\begin{equation}\label{eq:exactCollisionReconstruction}
\frac{H^{\rm C,\nu}(s)}{\sqrt \mu }
=\mathbb U_{\rm b}^\nu(s,t_0)(H_{\rm in}/\sqrt \mu )
-\int_{t_0}^s\mathbb U_{\rm b}^\nu(s,\tau)
\left[\mu^{-1/2}E[\rho^{\rm C,\nu}](\tau)
\partial_{v_1}F^{\rm b,\nu}(\tau)\right]d\tau.
\end{equation}
{Source convergence and the evolution energy estimate yield} continuity in $H^r_xL^2_\ell$ and finite integrated dissipation. The velocity density pairing is bounded on this energy space. Passing it through the convergent source sum and time integral gives, for $k\ne0$,
\[
{
\begin{aligned}
\rho_k[H^{\rm C,\nu}](s)
&=\mathfrak g_k^{\rm C,\nu}(s)-\int_{t_0}^s\!\int\sqrt \mu \,
\Big[\mathbb U_{\rm b}^\nu(s,\tau)
\big[\mu^{-1/2}E[\rho^{\rm C,\nu}](\tau)
\partial_{v_1}F^{\rm b,\nu}(\tau)\big]\Big]_k\,dv\,d\tau\\
&=\rho_k^{\rm C,\nu}(s),
\end{aligned}
}
\]
The last equality is \eqref{eq:exactDensityVolterra}.
For the zero density mode, local collision mass conservation and
velocity integration by parts give
\[
\partial_s\rho_0[H^{\rm C,\nu}]=0,\qquad
\rho_0[H^{\rm C,\nu}](t_0)=0,
\qquad \rho_0[H^{\rm C,\nu}](s)=0.
\]
Thus \eqref{eq:exactCollisionReconstruction} solves
\eqref{eq:exactCollisionModel} with the prescribed initial value,
including the homogeneous distribution mode.
\par
For uniqueness, Poisson's equation and the fixed-interval coefficient bounds give
\[
\|\mu^{-1/2}E[H]\partial_{v_1}F^{\rm b,\nu}\|_{H^r_xL^2_\ell}
\le C_{r,\ell,N,{T_1}}\|H/\sqrt \mu \|_{H^r_xL^2_\ell}.
\]
For two solutions with the same datum, put
$w=(H-\widetilde H)/\sqrt\mu$. The collision energy inequality gives
\[
\|w(s)\|_{H_x^rL_\ell^2}^2
\le C_{r,\ell,N,T_1}\int_{t_0}^s\|w(\tau)\|_{H_x^rL_\ell^2}^2\,d\tau,
\qquad w(t_0)=0.
\]
Gronwall yields $w=0$, so constructions at larger fixed orders agree.
\par
This proves the reconstruction and uniqueness assertions. The comparison
with the descending density, including the terminal lower bound, is
proved in the next subsection.
\end{proof}
\par\endgroup

\begingroup
The energy estimates will use the polynomial density norm
$\ell^1_qL^\infty$, without the weights $w_{k,N}$ and $B_k$ of
$\mathscr X_N^{r_{\rm x}}$. Before $T_{\rm a}$, the prepared high frequency has not reached its first critical time; the generic adjacent-mode kernel is also small because $a(1+T_{\rm a})=o(1)$. The following inverse gives a uniform early bound and permits an exponential loss on the late interval.

The density estimate may lose $e^{C_*N}$, provided $C_*$ is independent of the
fixed regularity order. Lemma~\ref{lem:collisionSourceEstimates} supplies the
collision density estimate for $L^2_v$ initial data and sources. We now estimate the electric inverse on the
two intervals already defined in \eqref{eq:electricTimeSplit}.

Use the split \eqref{eq:electricTimeSplit}. \begingroup
For $I'=[\tau_*,\tau^*]$, let $\mathbb V^{\rm C,\nu}_{I'}\mathfrak g$
denote the solution of the restricted density equation
\[
\rho_k(s)=\mathfrak g_k(s)+\sum_{m\ne0}\int_{\tau_*}^s
\left(\int_{\mathbb R^3}\sqrt\mu\,
[\mathbb U_{\rm b}^\nu(s,\tau)\mathcal F_m^{\rm b}(\tau)]_k\,dv\right)
\rho_m(\tau)\,d\tau,\qquad s\in I'.
\]
Equivalently, if
\[
(\mathbb T_{\nu,I'}\rho)_k(s)
=\sum_{m\ne0}\int_{\tau_*}^s
\left(\int\sqrt\mu\,[\mathbb U_{\rm b}^\nu(s,\tau)\mathcal F_m^{\rm b}(\tau)]_k\,dv\right)
\rho_m(\tau)\,d\tau,
\]
then $\mathbb V^{\rm C,\nu}_{I'}=(\mathrm{Id}-\mathbb T_{\nu,I'})^{-1}$.
The next proposition constructs this inverse on the polynomial density spaces.
These spaces do not impose the Fourier localization of
$\mathscr X_N^{r_{\rm x}}$, so the bounds are proved separately.
\par\endgroup
\begin{prop}\label{prop:coarseDensityInverse}
For each fixed $q\ge0$ and sufficiently small $\nu$,
\begin{equation}\label{eq:coarseDensityInverse}
\|\mathbb V^{\rm C,\nu}_{I_-}{\rho}\|_{\ell^1_qL^\infty(I_-)}
\le C_q\|{\rho}\|_{\ell^1_qL^\infty(I_-)},\qquad \|\mathbb V^{\rm C,\nu}_{I_+}{\rho}\|_{\ell^1_qL^\infty(I_+)}
\le C_qe^{C_*N}\|{\rho}\|_{\ell^1_qL^\infty(I_+)}.
\end{equation}
These estimates hold under \eqref{eq:comparisonLogParameters} for $p>{\rm s}_0+7/2$. The coefficient $C_*$ depends only on the fixed mode parameters, and is independent of fixed $q$, ${\rm s}_0$, $p$, $K$, and the expansion order. The constants $C_q$ and the smallness threshold may depend on the fixed orders and parameters.
\end{prop}

\begingroup
\begin{proof}
First invert the Maxwellian term in the density equation with Maxwellian collisions. {The scalar density kernel is}
\[
K_k^\nu(t)=\int\sqrt \mu \,\mathbb S_k^\nu(t)(-iv_1\sqrt \mu /k)\,dv,
\qquad k\ne0.
\]
By \eqref{eq:centeredpacketdifference} with $\eta=0$ and $\nu {T_1}=o(1)$,
\[
|K_k^\nu(t)-K_{|k|}(t)|
\le\frac{C\nu t}{|k|}(1+k^2t^2)e^{-c|kt|^\chi},
\qquad
\int_0^{T_1}|K_k^\nu-K_{|k|}|\,dt\le C\nu/|k|^3.
\]
By Lemma~\ref{prop:inverse}, uniformly in $k\ne0$,
\[
\begin{aligned}
\mathrm{Id}-K_k^\nu*
&=(\mathrm{Id}-K_{|k|}*)
\left[\mathrm{Id}-(\mathrm{Id}-K_{|k|}*)^{-1}
(K_k^\nu-K_{|k|})*\right],\\
\| (\mathrm{Id}-K_{|k|}*)^{-1}(K_k^\nu-K_{|k|})*\|_{L^\infty\to L^\infty}
&\le C\|K_k^\nu-K_{|k|}\|_{L^1}\le C\nu/|k|^3\le\tfrac12,\\
\|(\mathrm{Id}-K_k^\nu*)^{-1}\|_{L^\infty\to L^\infty}
&\le2\|(\mathrm{Id}-K_{|k|}*)^{-1}\|_{L^\infty\to L^\infty}\le C_0.
\end{aligned}
\]
The inversion is modewise, so the same $C_0$ applies to any positive
constant weights in the mode sum.
\par
For the late interval choose a fixed large $C_{\rm maj}>0$ and set
\begin{equation}\label{eq:coarseGainWeights}
G_k=\prod_{m=k+1}^{\infty}\left(1+\frac{C_{\rm maj}a(1+{T_1})}{m^2}\right),\quad k\ge1,
\qquad
\|{\rho}\|_{G,q}
=\sum_{k\ne0}\frac{\langle k/N\rangle^q}{G_{|k|}}\sup_{I_+}|{\rho}_k|.
\end{equation}
The notation $\|\cdot\|_{G,q}$ is local to this proof, with $N$ and
$I_+$ fixed; the same subscript on an operator {denotes the induced operator norm.}
The polynomial weight in this auxiliary norm uses the scale $N$;
the final density norm $\ell^1_qL^\infty(I_+)$ retains the weight
$\langle k\rangle^q$. The product converges and
\begin{equation}\label{eq:coarseGainRelations}
\frac{G_{m-1}}{G_m}=1+\frac{C_{\rm maj}a(1+{T_1})}{m^2},\qquad
1\le G_k\le e^{C\sqrt{a(1+{T_1})}}\le e^{CN}.
\end{equation}
For the second bound, monotonicity of the summand and integration by
parts give
\[
\begin{aligned}
\sum_{m\ge1}\log\left(1+\frac{C_{\rm maj}a(1+T_1)}{m^2}\right)
&\le\int_0^\infty\log\left(1+\frac{C_{\rm maj}a(1+T_1)}{x^2}\right)\,dx\\
&=2C_{\rm maj}a(1+T_1)\int_0^\infty
\frac{dx}{x^2+C_{\rm maj}a(1+T_1)}\\
&=\pi\sqrt{C_{\rm maj}a(1+T_1)}.
\end{aligned}
\]
The boundary terms vanish at zero and infinity. By
\eqref{eq:observationGainBalance}, this logarithmic cost is $O(N)$.
For a longer endpoint $T\sim\zeta$, the same upper bound and the
logarithmic summands with $2\le m\le\lfloor\sqrt{C_{\rm maj}a(1+T)}\rfloor$,
each at least $\log2$, give
\[
\log\prod_{m=2}^{\infty}\left(1+\frac{C_{\rm maj}a(1+T)}{m^2}\right)
\sim\sqrt{a(1+T)}\sim N^{3/2}.
\]
This is the longer-time weight-removal cost discussed in
Section~\ref{sec:roadmapFieldCompletion}.
\par
\par

For any matrix of time kernels $\mathcal K_{k,m}(s,\tau)$, the following column
estimate controls the associated integral operator:
\begin{equation}\label{eq:coarseColumnCriterion}
\begin{aligned}
&\left\|\left\{\sum_{m\ne0}\int_{T_{\rm a}}^s
\mathcal K_{k,m}(s,\tau){\rho}_m(\tau)\,d\tau\right\}_{k\ne0}\right\|_{G,q}\\
&\quad\le\|{\rho}\|_{G,q}\sup_{m\ne0}\sum_{k\ne0}
\frac{\langle k/N\rangle^qG_{|m|}}{\langle m/N\rangle^qG_{|k|}}
\sup_{s\in I_+}\int_{T_{\rm a}}^s|\mathcal K_{k,m}(s,\tau)|\,d\tau.
\end{aligned}
\end{equation}
It follows by taking the supremum over $s$ separately for each $k$ before summing, and then interchanging two nonnegative sums. 

For the adjacent kernel in \eqref{eq:densityexact},
$|P_{k,m}(s,\tau)|\le Ca(|ks-m\tau|+\tau)e^{-|ks-m\tau|^2/2}/|m|$.
Thus the substitution $z=ks-m\tau$ gives
\[
\sup_s\int_{T_{\rm a}}^s|P_{k,m}(s,\tau)|\,d\tau
\le\frac{Ca}{m^2}\int_{\mathbb R}(|z|+T_1)e^{-z^2/2}\,dz
\le\frac{Ca(1+T_1)}{m^2}.
\]
For $|k|=|m|-1$, the polynomial weight ratio is at most one and
\[
\begin{aligned}
\frac{G_{|m|}}{G_{|m|-1}}
\sup_s\int_{T_{\rm a}}^s|P_{k,m}(s,\tau)|\,d\tau
&\le\frac{Ca(1+T_1)/m^2}{1+C_{\rm maj}a(1+T_1)/m^2}
\le\frac{C}{C_{\rm maj}},\\
C_0\|\mathbb P_\downarrow\|_{G,q}
&\le\frac{C_0C}{C_{\rm maj}}\le\frac14
\qquad(C_{\rm maj}\ge4C_0C).
\end{aligned}
\]
This choice is independent of $q$.
For $|k|=|m|+1$, $|ks-m\tau|\ge s\ge t_0$, and
\[
\begin{aligned}
\frac{\langle k/N\rangle^qG_{|m|}}{\langle m/N\rangle^qG_{|k|}}
&\le C_qN^2,\\
\|\mathbb P_\uparrow\|_{G,q}
&\le C_qN^2e^{-ct_0^2}
\sup_{m\ne0}\frac{a(1+T_1)}{m^2}
\le C_qN^4e^{-ct_0^2}=o(1).
\end{aligned}
\]
After the scalar inverse, the sum of these two norms is at most
$1/4+o(1)<1/2$. Hence the Maxwellian collisional density operator together
with the two adjacent-mode integral operators for free transport has a uniformly bounded
inverse in $\|\cdot\|_{G,q}$.
\par
We now recover the exact density equation \eqref{eq:exactCoefficientDensityEquation} on $I_+$.
The three differences are the same as in
\eqref{eq:directResponseDecomposition}, {except that the Maxwellian}
diagonal collisional part is already included in {the approximate equation used here}.

For $[\mathcal F_{{\rm app},m}]_{m\pm1}$, \eqref{eq:maxwellianKernel} bounds the
change from free transport to {$\mathbb U_\mu^\nu$} by
\[
\frac{Ca(1+{T_1})\nu {T_1}^3}{|m|}e^{-c|ks-m\tau|^\chi}.
\]
When $|k|=|m|-1$, the gain ratio is at most one, {and the time integral}
is at most $C_qa(1+{T_1})\nu {T_1}^3/m^2\le C_qN^2\nu {T_1}^3$.
When $|k|=|m|+1$, the gain ratio is at most $CN^2$, but
$|ks-m\tau|\ge s\ge T_{\rm a}$. Reserve a fraction of the phase
decay to absorb $CN^2$; the remaining phase again integrates to
$C/|m|$. The same bound follows.

For the change from {$\mathbb U_\mu^\nu$} to $\mathbb U_{\rm b}^\nu$
acting on $[\mathcal F_{{\rm app},m}]_{m\pm1}$ and $[\mathcal F_{{\rm app},m}]_m$, use
\eqref{eq:directKernel}. Put
$n=k-m-l$, with $l=0,\pm1$.
The ratio $G_{|m|}/G_{|k|}$ is at most one for $|k|\le |m|$.
Otherwise $|k(s-t_0)-m(\tau-t_0)|\ge N^4$ and $|k(s-t_0)-m(\tau-t_0)+lt_0|\ge N^4/2$.
A reserved fraction of the phase absorbs $e^{CN}$.
The polynomial ratio is absorbed by a fraction of the spatial decay.
The remaining factors have time integral $C/|m|$ and a finite sum
in $n$. Thus the column criterion \eqref{eq:coarseColumnCriterion}
gives the bound
\[
\sup_{m\ne0}\frac{C_q\varepsilon_N[1+a(1+{T_1})]}{m^2}
\le C_qN^2\varepsilon_N.
\]

For the electric remainder, \eqref{eq:remainderKernel} holds with
$\varepsilon_{\rm ker}=\nu T_1^3$. To apply the column bound
\eqref{eq:generalGainColumn}, take $\mathfrak w=G$. Equation
\eqref{eq:coarseGainRelations} gives
\[
G_m/G_k\le1\quad(1\le k\le m),\qquad
G_m/G_k\le e^{CN}\quad(k,m\ge1).
\]
Fix ${\bf b}_2\in(0,{\bf b}_0/4)$. Since $a(1+{T_1})\ge cN^2$,
\[
\frac{G_m}{G_k}
=\prod_{j=k+1}^m(1+C_{\rm maj}a(1+{T_1})/j^2)^{-1}
\le(1+C_{\rm maj}c/{\bf b}_2^2)^{-(m-k)},\qquad
1\le k\le m\le {\bf b}_2N.
\]
Since $(1+C_{\rm maj}c/{\bf b}_2^2)^{-1}<1$, all three hypotheses in
\eqref{eq:generalGainHypotheses} hold.
Equations \eqref{eq:generalGainColumn} and \eqref{eq:coarseColumnCriterion}
give directly
\begin{equation}\label{eq:coarseReducedRemainder}
\left\|\left\{\sum_{m\ne0}\int_{T_{\rm a}}^s
\Delta\mathcal K^{(3)}_{k,m}(s,\tau){\rho}_m(\tau)\,d\tau
\right\}_{k\ne0}\right\|_{G,q}
\le C_qN^3\nu T_1^3\|{\rho}\|_{G,q}.
\end{equation}

{The perturbation after applying the inverse just constructed is bounded by}
\begin{equation}\label{eq:coarseUniformPerturbation}
C_q\big(N^2\nu {T_1}^3+N^2\varepsilon_N
+N^3\nu T_1^3\big)=o(1).
\end{equation}
Writing $p={\rm s}_0+7/2+\eta$ with $\eta>0$, the first and third terms
are bounded respectively by
$C[\log(1/\nu)]^{-5/2-3\eta}$ and
$C[\log(1/\nu)]^{-3/2-3\eta}+Ce^{-cN^2}$.
The middle term is $\nu^{1/3-o(1)}$.
Absorption gives a uniformly bounded exact inverse in the norm
$\|\cdot\|_{G,q}$.

To return to the original density norm, the scaled polynomial weights give
\begin{equation}\label{eq:coarseScaledNormConversion}
\|{\rho}\|_{G,q}
\le\|{\rho}\|_{\ell^1_qL^\infty(I_+)},
\qquad
\|{\rho}\|_{\ell^1_qL^\infty(I_+)}
\le C_qN^q e^{CN}\|{\rho}\|_{G,q}.
\end{equation}
Here $C$ depends only on the fixed gain parameter $C_{\rm maj}$ and the chain parameters. For every fixed $q$, $N^q\le e^N$ once $N$ is sufficiently large. Taking $C_*=C+1$ therefore proves the late bound in \eqref{eq:coarseDensityInverse}, with $C_*$ independent of the spatial order and of the amplitude calibration. The smallness threshold is allowed to depend on those fixed orders.

On $I_-$ use $\langle k\rangle^q$ without the gain weight.
The corresponding perturbations of the scalar inverse satisfy
\[
\begin{array}{ll}
\text{free adjacent-mode operators:}&C_q a(1+T_{\rm a}),\\
\text{adjacent collision difference:}&C_q a(1+T_{\rm a})\nu T_{\rm a}^3,\\
\text{collision-coefficient difference:}&C_q\varepsilon_N[1+a(1+T_{\rm a})],\\
\text{electric-coefficient remainder:}&C_q\varepsilon_N(1+T_{\rm a}).
\end{array}
\]
The last bound uses the fixed-radius coefficient estimate in
Lemma~\ref{lem:physicalCoefficientDerivatives} and
\eqref{eq:exactCollisionTrace}. Since $T_{\rm a}=t_0+N^4$,
\[
a(1+T_{\rm a})+a(1+T_{\rm a})\nu T_{\rm a}^3
+\varepsilon_N[1+a(1+T_{\rm a})]+
\varepsilon_N(1+T_{\rm a})=\nu^{1/3-o(1)}\longrightarrow0.
\]
Multiplication by the fixed $C_qC_0$ is therefore less than $1/2$
for sufficiently small $\nu$, proving the early inverse.
The column estimates extend both inverses from finite mode sums
to the completed density spaces.
\par
The choices can be recorded without reference to a derivative order:
\begin{equation}\label{eq:coarseConstantDependence}
\begin{aligned}
C_{\rm maj}&=C_{\rm maj}({\bf b}_0,{\bf b}_1,\chi),& C_*&=C+1,\\
N\ge N_q&\ \Longrightarrow\ N^q\le e^N,&
0<\nu<\nu_q({\rm s}_0,p,K)&\ \Longrightarrow\ \eqref{eq:coarseDensityInverse}.
\end{aligned}
\end{equation}

\end{proof}
\par\endgroup
\begingroup
For the full solution on $I_+$, restricting the equation leaves the
earlier integral in the prescribed density:
\[
\begin{aligned}
\rho_k(s)
&=\underbrace{\mathfrak g_k(s)+\sum_{m\ne0}\int_{t_0}^{T_{\rm a}}
\left(\int\sqrt\mu\,[\mathbb U_{\rm b}^\nu(s,\tau)\mathcal F_m^{\rm b}(\tau)]_k\,dv\right)
\rho_m(\tau)\,d\tau}_{\displaystyle\mathfrak g_{+,k}(s)}\\
&\quad+\sum_{m\ne0}\int_{T_{\rm a}}^s
\left(\int\sqrt\mu\,[\mathbb U_{\rm b}^\nu(s,\tau)\mathcal F_m^{\rm b}(\tau)]_k\,dv\right)
\rho_m(\tau)\,d\tau.
\end{aligned}
\]
Hence $\rho|_{I_+}=\mathbb V^{\rm C,\nu}_{I_+}\mathfrak g_+$.
In the complete inhomogeneous equation below we instead restart the
distribution equation at $T_{\rm a}$ with datum $H(T_{\rm a})$;
the early estimate controls this datum.
\par\endgroup

\subsection{Proof of Proposition~\ref{thm:exactCollisionEcho}}\label{sec:threeDensityErrors}
The estimates above now control the second, third and fourth terms of
\eqref{eq:linearConstructionThreeErrors}. They correspond to the initial
collision evolution, the upward interaction omitted from the descending
chain, and the difference of the density kernels. We state each contribution
separately before combining them. Throughout this subsection assume
\eqref{eq:comparisonLogParameters} with $p>{\rm s}_0+7/2$, fix
$r_{\rm x}\ge0$, and use $\mathcal O_N\rho=B_*^{-1}\rho_{m_1}(T_1)$.

\begin{lem}\label{lem:initialEvolutionTerminal}
The initial density difference satisfies
\[
\begin{aligned}
\|\mathbb V^{\rm C,\nu}(\mathfrak g^{\rm C,\nu}-\mathfrak g^0)\|_{\mathscr X_N^{r_{\rm x}}}
&\le C_{r_{\rm x}}N(\varepsilon_{\rm coll}+\varepsilon_N),\\
|\mathcal O_N\mathbb V^{\rm C,\nu}(\mathfrak g^{\rm C,\nu}-\mathfrak g^0)|
&\le C_{r_{\rm x}}(\varepsilon_{\rm coll}+\varepsilon_N).
\end{aligned}
\]
\end{lem}
\begin{proof}
Corollary~\ref{cor:initialDensityComparison} and
\eqref{eq:logExactCollisionInverse}--\eqref{eq:exactTerminalInverse} give
\[
\begin{aligned}
\|\mathbb V^{\rm C,\nu}(\mathfrak g^{\rm C,\nu}-\mathfrak g^0)\|_{\mathscr X_N^{r_{\rm x}}}
&\le\|\mathbb V^{\rm C,\nu}\|\,
\|\mathfrak g^{\rm C,\nu}-\mathfrak g^0\|_{\mathscr X_N^{r_{\rm x}}}
\le C_{r_{\rm x}}N(\varepsilon_{\rm coll}+\varepsilon_N),\\
|\mathcal O_N\mathbb V^{\rm C,\nu}(\mathfrak g^{\rm C,\nu}-\mathfrak g^0)|
&\le\|\mathcal O_N\mathbb V^{\rm C,\nu}\|\,
\|\mathfrak g^{\rm C,\nu}-\mathfrak g^0\|_{\mathscr X_N^{r_{\rm x}}}
\le C_{r_{\rm x}}(\varepsilon_{\rm coll}+\varepsilon_N).
\end{aligned}
\]
\end{proof}

\begin{lem}\label{lem:upwardTerminalError}
There are fixed $c,C_{r_{\rm x}}>0$ such that
\[
\|\mathbb V^{\rm C,\nu}\mathbb P_\uparrow\rho^\downarrow\|_{\mathscr X_N^{r_{\rm x}}}
+|\mathcal O_N\mathbb V^{\rm C,\nu}\mathbb P_\uparrow\rho^\downarrow|
\le C_{r_{\rm x}}e^{-cN^2}.
\]
\end{lem}
\begin{proof}
For some fixed $c>0$, \eqref{eq:descendingReferenceDefect} gives
$\|\mathbb P_\uparrow\rho^\downarrow\|_{\mathscr X_N^{r_{\rm x}}}
\le C_{r_{\rm x}}e^{-cN^2}$. Hence
\[
\|\mathbb V^{\rm C,\nu}\mathbb P_\uparrow\rho^\downarrow\|_{\mathscr X_N^{r_{\rm x}}}
+|\mathcal O_N\mathbb V^{\rm C,\nu}\mathbb P_\uparrow\rho^\downarrow|
\le C_{r_{\rm x}}(N+1)e^{-cN^2} \le C_{r_{\rm x}}e^{-cN^2/2}.
\]
\end{proof}

\begin{lem}\label{lem:kernelTerminalError}
The density-kernel difference satisfies
\[
\begin{aligned}
\|\mathbb V^{\rm C,\nu}\Delta\mathbb T\rho^\downarrow\|_{\mathscr X_N^{r_{\rm x}}}
&\le C_{r_{\rm x}}(N^4\nu T_1^3+N^4\varepsilon_N+N^5\nu T_1^3),\\
|\mathcal O_N\mathbb V^{\rm C,\nu}\Delta\mathbb T\rho^\downarrow|
&\le C_{r_{\rm x}}(N^3\nu T_1^3+N^3\varepsilon_N+N^4\nu T_1^3).
\end{aligned}
\]
\end{lem}
\begin{proof}
By \eqref{eq:descendingTwoIntervalBounds} and
Lemma~\ref{prop:exactCollisionResponse},
\[
\|\Delta\mathbb T\rho^\downarrow\|_{\mathscr X_N^{r_{\rm x}}}
\le\|\Delta\mathbb T\|\,\|\rho^\downarrow\|_{\mathscr X_N^{r_{\rm x}}} \le C_{r_{\rm x}}N(N^2\nu T_1^3+N^2\varepsilon_N+N^3\nu T_1^3).
\]
Applying respectively
$\|\mathbb V^{\rm C,\nu}\|\le C_{r_{\rm x}}N$ and
$\|\mathcal O_N\mathbb V^{\rm C,\nu}\|\le C_{r_{\rm x}}$
proves the two estimates.
\end{proof}

\begin{proof}[Proof of Proposition~\ref{thm:exactCollisionEcho}]
Proposition~\ref{prop:exactDensityInverse} constructs the unique distribution
$H^{\rm C,\nu}$ and identifies its density. The exact difference is
\[
\rho[H^{\rm C,\nu}]-\rho^\downarrow
=\mathbb V^{\rm C,\nu}(\mathfrak g^{\rm C,\nu}-\mathfrak g^0)
+\mathbb V^{\rm C,\nu}\mathbb P_\uparrow\rho^\downarrow
+\mathbb V^{\rm C,\nu}\Delta\mathbb T\rho^\downarrow.
\]
\begingroup
Apply $\mathcal O_N$ to the preceding identity and use
Lemmas~\ref{lem:initialEvolutionTerminal}--\ref{lem:kernelTerminalError}:
\begin{equation}\label{eq:exactTerminalSource}
\begin{aligned}
|\mathcal O_N(\rho^{\rm C,\nu}-\rho^\downarrow)|
&\le |\mathcal O_N\mathbb V^{\rm C,\nu}(\mathfrak g^{\rm C,\nu}-\mathfrak g^0)|
+|\mathcal O_N\mathbb V^{\rm C,\nu}\mathbb P_\uparrow\rho^\downarrow| +|\mathcal O_N\mathbb V^{\rm C,\nu}\Delta\mathbb T\rho^\downarrow|\\
&\le C_{r_{\rm x}}\big(\varepsilon_{\rm coll}+\varepsilon_N
+N^3\nu T_1^3+N^3\varepsilon_N+N^4\nu T_1^3+e^{-cN^2}\big).
\end{aligned}
\end{equation}
\par\endgroup
Since
$\zeta=m_1{T_1}$, $m_1\mathrel{{\sim}} N$, and $\varepsilon_N=aL_N$,
\[
\varepsilon_{\rm coll}=\nu {T_1}(1+\zeta^2)\le CN^2\nu T_1^3.
\]
Consequently,
\[
\begin{aligned}
\sqrt N\,|\mathcal O_N(\rho^{\rm C,\nu}-\rho^\downarrow)|
&\le C_{r_{\rm x}}\big(\sqrt N\varepsilon_{\rm coll}
+N^{7/2}\nu {T_1}^3+N^{7/2}\varepsilon_N
+N^{9/2}\nu T_1^3+e^{-c'N^2}\big)\\
&\le C_{r_{\rm x}}\big(N^{9/2}\nu T_1^3+N^{7/2}aL_N\big).
\end{aligned}
\]
Here the Gaussian term is absorbed into $N^{7/2}aL_N$, since
$a\ge e^{-\gamma_{\rm a}N}$ and $L_N\ge1$.
This proves \eqref{eq:unifiedExactTerminal}. Moreover,
\[
N^{9/2}\nu {T_1}^3\mathrel{{\sim}}[\log(1/\nu)]^{3{\rm s}_0+21/2-3p},
\qquad N^{7/2}aL_N=\nu^{1/3-o(1)}.
\]
The logarithmic exponent is negative for
$p>{\rm s}_0+7/2$, proving \eqref{eq:logExactCollisionTerminal}.
Proposition~\ref{thm:chain} therefore gives
\[
|\rho_{m_1}^{\rm C,\nu}({T_1})|
\ge |\rho_{m_1}^\downarrow({T_1})|
-|\rho_{m_1}^{\rm C,\nu}({T_1})-\rho_{m_1}^\downarrow({T_1})|
\ge(c-o(1))\frac{B_*}{\sqrt N}.
\]
Also ${T_1}=o(\nu^{-1/3})$ by \eqref{eq:comparisonLogParameters}.

\end{proof}
The last two terms of \eqref{eq:linearConstructionThreeErrors} are
$\sum_{n=1}^J\rho[H_n]$ and $\rho[W_J]$. They arise from adding the
background-field transport term $E[F^{\rm b,\nu}]\partial_{v_1}H$.
The next subsection establishes the complete inhomogeneous estimate;
Section~\ref{sec:fieldExpansion} then constructs the finite corrections
and estimates the remaining source.
\par\endgroup
\subsection{\texorpdfstring{Solution operators associated with $\mathcal A_\nu$}{Solution operators associated with A-nu}}\label{sec:fullLinearClosure}\label{sec:fullRangeClosure}

\begingroup
We prove the complete nonhomogeneous estimate used for the final linear
remainder and the nonlinear iteration. The equations are
\[
\mathcal A_\nu H=0,\qquad
\mathcal A_\nu H-\mathbb Q_\nu({\bf H}+{\bf W},H)
-\mathbb Q_\nu(H,{\bf H})=S_0+S_1.
\]
Here $\mathcal A_\nu$ is defined in
\eqref{eq:fullPhysicalOperatorDefinition}. Its coefficients depend on the
evolving background solution; after these coefficients are fixed,
the associated solution operators depend on both the initial and final
times and form an evolution family rather than a one-parameter semigroup.
The homogeneous equation uses
the prepared datum $H_{\rm in}$; the second problem allows arbitrary
finite-energy data and sources satisfying \eqref{eq:quadraticSourceSpaces}.
Taking ${\bf H}={\bf W}=0$ gives the equation for the remainder in
\eqref{eq:fieldFiniteUnknownError}. Section~\ref{sec:nonlinear} verifies
the coefficient hypotheses for the scaled linear solution and the
nonlinear remainder.

These sources need not have the phase localization required by
$\mathscr X_N^{r_{\rm x}}$. We therefore use the inverse in the polynomially weighted density norm on the early and late intervals, constructed in Proposition~\ref{prop:coarseDensityInverse}, and combine it with the energy estimate. Background-field transport
remains on the left; its contribution is controlled by late-time field
decay. Proposition~\ref{thm:fullForcedCoarse} gives a loss $e^{C_*N}$,
with $C_*$ independent of the fixed derivative order and chosen before
the correction order. The individual increments solve
\eqref{eq:fieldCorrectionRecursion}; the more precise estimates in
Section~\ref{sec:fieldExpansion} make their final residual small enough
for this complete estimate. The underlying solution constructions and
source-to-density estimates are proved in Appendix~\ref{app:evolution}.
\par\endgroup
We take $\sigma_{\rm int}\in(1/2,1)$ and retain a fixed velocity moment at least
$2{\rm s}_0+6$ in the direct background comparison.

Write
\[
{
{f^{\rm b,\nu}(s)=\widetilde f^{\rm b,\nu}(s-t_0)},\qquad E^{\rm b,\nu}=-\partial_{x_1}\phi^{\rm b,\nu},\qquad \int_0^{2\pi}\phi^{\rm b,\nu}(s,x_1)\,dx_1=0.
}
\]
Thus $-\partial_{x_1}^2\phi^{\rm b,\nu}=\int F^{\rm b,\nu}\,dv-1$.
Use $\mathcal E_r,\mathcal D_r,\mathcal N_r$ and $\mathscr E_r$ from
\eqref{eq:fullLinearEnergyNorms}, with
the explicit velocity weights and full-interval convention stated there.
For the density we use $\ell^1_qL^\infty(I')$ from
\eqref{eq:fourierDensityNorm}. {The polynomial weight in this norm is} $\langle k\rangle^q$;
it has neither the phase weight nor the gain denominator of
$\mathscr X_N^{r_{\rm x}}$. {The dissipation $\mathcal D_r$ includes $\nu$, whereas the collision form
$\mathcal D^{\rm ph}_{c,j}$ is defined without this factor.}

{The two source conditions are given in
\eqref{eq:quadraticSourceSpaces}. The energy-integrable part pairs with
the energy supremum; the collision-form part pairs with
${(\int_{I'}\mathcal D_r)^{1/2}}$.}

{For prescribed real distributions {${\bf H},{\bf W}$} of zero
 total mass, we use the prescribed-coefficient equation in
 \eqref{eq:fullForcedEquation}. The coefficient
 {$F^{\rm b,\nu}+{\bf H}+{\bf W}$} occurs in transport and the first
 collision ordering, whereas {$F^{\rm b,\nu}+{\bf H}$} occurs in the
 differentiated electric coefficient and the second collision ordering.
 In Section~\ref{sec:nonlinear} we take {${\bf H}=\delta_\nu H^\nu$}
 and {${\bf W}$} to be a prescribed remainder.}

\begingroup
The energy and density orders used below are
\[
\sup_{s\in I'}\mathcal E_r[H(s)/\sqrt\mu]^{1/2}
+\left(\int_{I'}\mathcal D_r[H/\sqrt\mu]\,ds\right)^{1/2},
\qquad \|\rho[H]\|_{\ell^1_{r-1}L^\infty(I')}.
\]
\par\endgroup

For collision sources we use $\Gamma$ from
\eqref{eq:normalizedGammaExplicit} and the coefficients
\begin{equation}\label{eq:normalizedBilinearNotation}
\sigma_{ij}[{F}]=\Phi_{ij}*{F},\qquad
b_i[{F}]=(\partial_j\Phi_{ij})*{F}.
\end{equation}
The argument of $\sigma[\cdot]$ and $b[\cdot]$ is always the physical
distribution. Thus a normalized coefficient $f$ contributes
{$\sigma[\sqrt \mu  f]$} and {$b[\sqrt \mu  f]$}; at Gaussian width
${\varpi}$ the physical distribution is {$\sqrt{\mu_{{\varpi}}}f$}.
The unadorned matrix {$\sigma=\sigma[\mu ]$} is the Maxwellian matrix used in the form norms.
{Lemma~\ref{lem:quadraticLandauForm} proves the source estimate for both
collision orderings by integration by parts.}

Define the normalized field operator
\begin{equation}\label{eq:fieldNormalizedOperator}
\widetilde{\mathcal T}_{\rm b}(s)f
=-E[F^{\rm b,\nu}](s,x_1)
\big(\partial_{v_1}-v_1/2\big)f.
\end{equation}
Thus {$\sqrt \mu \,\widetilde{\mathcal T}_{\rm b}(s)f=-\mathcal T_{\rm b}(\sqrt \mu  f)$}, with the sign corresponding to placing background field transport on the right.

The coefficient bounds, weighted multiplier estimates and potential identity used
next are proved in Appendix~\ref{app:energyCoefficients}. In particular, the
potential identity removes the unbounded $v_1E^{\rm b,\nu}$ factor created by
Gaussian normalization. {We give the application to the present equation and the full energy estimate here.}

For every integer $r\ge1$ and mean-zero density $\rho$ on $I^\prime\subset I$, the Fourier convention and Poisson equation give
\begin{equation}\label{eq:densityNormComparison}
\|\rho\|_{L^\infty(I';H_x^{r-1})}
\le\sqrt{2\pi}\|\rho\|_{\ell^1_{r-1}L^\infty(I')},\qquad \|E[\rho]\|_{L^\infty(I';H_x^r)}
\le\sqrt2\|\rho\|_{L^\infty(I';H_x^{r-1})}.
\end{equation}
At each time, the Fourier norm in Section~\ref{sec:notation} gives
\[
\|E[{\rho}](s)\|_{H_x^r}
\le\sqrt2\|{\rho}(s)\|_{H_x^{r-1}}
\le2\sqrt\pi\sum_{k\ne0}\langle k\rangle^{r-1}|{\rho}_k(s)|.
\]
Here $E_k=-i{\rho}_k/k$ for $k\ne0$ and the field zero mode vanishes.
Taking the time supremum proves \eqref{eq:densityNormComparison}.

\begin{prop}\label{prop:fullLinearEnergy}\label{lem:energyDensityEmbeddings}
Fix an integer $r\ge1$ and use the parameters of Proposition~\ref{lem:fixedCoefficientOrders}.

{Let {$H^\nu$} be a finite-energy solution of
\eqref{eq:fullRangePhysicalEquation}, and set}
\begin{equation}\label{eq:fullPotentialUnknown}
f=H^\nu/\sqrt \mu ,\qquad
y=e^{\phi^{\rm b,\nu}/2}f.
\end{equation}
Then, for sufficiently small $\nu$,
\begingroup
\begin{equation}\label{eq:fullLinearEnergyClosed}
\begin{aligned}
&\sup_{s\in I}\mathcal E_r[y(s)]^{1/2}
+\left(\int_I\mathcal D_r[y(s)]\,ds\right)^{1/2}\\
&\quad\le C_r\left(\mathcal E_r[y(t_0)]^{1/2}
+T_1[1+\varepsilon_N(1+T_1)]
\|\rho[H^\nu]\|_{L^\infty(I;H_x^{r-1})}\right).
\end{aligned}
\end{equation}
\par\endgroup

The same energy and dissipation norms of $f$ are bounded by a constant times those of $y$.
\end{prop}

\begingroup
\begin{proof}

Apply Lemma~\ref{lem:potentialTransport} with
${\varpi}=1$, $E=E^{\rm b,\nu}$ and $\phi=\phi^{\rm b,\nu}$.
The collision operator is local in $x$ and linear in the unknown, so
\[
e^{\phi^{\rm b,\nu}/2}\mathcal C_{\rm b}
(e^{-\phi^{\rm b,\nu}/2}y)=\mathcal C_{\rm b}y.
\]
Thus \eqref{eq:generalPotentialConjugation} gives the exact equation
\begin{equation}\label{eq:fullPotentialEquation}
\partial_sy+v_1\partial_{x_1}y+E^{\rm b,\nu}\partial_{v_1}y
=\nu\mathcal C_{\rm b}y
+\tfrac12\partial_s\phi^{\rm b,\nu}y
-e^{\phi^{\rm b,\nu}/2}E[H^\nu]
\mu^{-1/2}\partial_{v_1}F^{\rm b,\nu}.
\end{equation}
\begingroup
The contribution of $E^{\rm b,\nu}\partial_{v_1}y$ to the differentiated energy satisfies

\[
\left|\operatorname{Re}\sum_{j=0}^r
\left\langle \partial_{x_1}^j(E^{\rm b,\nu}\partial_{v_1}y),
\partial_{x_1}^jy\right\rangle_{2r+6-2j}\right|\le
\tfrac14c_r\mathcal D_r[y]
+C_r\big(\|E^{\rm b,\nu}\|_\infty+
\nu^{-1}\|E^{\rm b,\nu}\|_{C^r_x}^2\big)\mathcal E_r[y]
\]
\par\endgroup
by \eqref{eq:generalFieldEnergy}. Lemma~\ref{lem:weightedDualMultipliers}
and \eqref{eq:backgroundPotentialIntegral} give
\[
\mathcal E_r[e^{\pm\phi^{\rm b,\nu}/2}f]
+\mathcal D_r[e^{\pm\phi^{\rm b,\nu}/2}f]
\le C_r(\mathcal E_r[f]+\mathcal D_r[f]).
\]

{The contribution of $\tfrac12\partial_s\phi^{\rm b,\nu}y$ to the differentiated energy is bounded by}
$C_r\|\partial_s\phi^{\rm b,\nu}\|_{C^r_x}\mathcal E_r$.
{For the term $-e^{\phi^{\rm b,\nu}/2}E[H^\nu]\mu^{-1/2}\partial_{v_1}F^{\rm b,\nu}$ in \eqref{eq:fullPotentialEquation}, use}
\[
\mu^{-1/2}\partial_{v_1}F^{\rm b,\nu}
=-v_1\sqrt \mu +(\partial_{v_1}-v_1/2)f^{\rm b,\nu}.
\]
Choosing $L>2r+11$ in
\eqref{eq:fullBackgroundFourierEnvelope} gives
\[
\max_{0\le j_1\le j\le r}
\|\langle v\rangle^{2r+6-2j}\partial_{x_1}^{j-j_1}(\mu^{-1/2}\partial_{v_1}F^{\rm b,\nu})\|_{L^\infty_xL^2_v}
\le C_r(1+\varepsilon_N\langle s\rangle).
\]
The finite product rule and \eqref{eq:densityNormComparison} give
\begin{equation}\label{eq:fullElectricSourceProduct}
\begin{aligned}
\mathcal E_r[E[{\rho}]\mu^{-1/2}\partial_{v_1}F^{\rm b,\nu}]^{1/2}
&\le\sum_{j=0}^r\sum_{j_1=0}^j\binom j{j_1}
\|\partial_{x_1}^{j_1}E[{\rho}]\|_{L^2_x}
\|\langle v\rangle^{2r+6-2j}\partial_{x_1}^{j-j_1}(\mu^{-1/2}\partial_{v_1}F^{\rm b,\nu})\|_{L^\infty_xL^2_v}\\
&\le C_r(1+\varepsilon_N\langle s\rangle)\|E[{\rho}]\|_{H^r_x}\\
&\le C_r(1+\varepsilon_N\langle s\rangle)\|{\rho}\|_{H^{r-1}_x}.
\end{aligned}
\end{equation}
The multiplier bound \eqref{eq:spatialMultiplierNorms} applies to this
source with multiplier $e^{\phi^{\rm b,\nu}/2}$; hence
\eqref{eq:fullElectricSourceProduct} also bounds the conjugated source.
Together with Lemma~\ref{lem:fullCollisionForm}, these estimates yield
\begin{equation}\label{eq:fullEnergyDifferential}
\frac{d}{ds}\mathcal E_r+c_r\mathcal D_r
\le C_r c_{\rm en}^{\rm b}(s)\mathcal E_r
+C_r(1+\varepsilon_N\langle s\rangle)
\|\rho[H^\nu](s)\|_{H^{r-1}_x}\mathcal E_r^{1/2},
\end{equation}
where
\[
c_{\rm en}^{\rm b}(s)=\nu(1+\varepsilon_N\langle s\rangle^2)
+\|E^{\rm b,\nu}\|_{C^r_x}
+\|\partial_s\phi^{\rm b,\nu}\|_{C^r_x}
+\nu^{-1}\|E^{\rm b,\nu}\|_{C^r_x}^2.
\]
The coefficient estimates give
\begin{equation}\label{eq:fullEnergyIntegratedCoefficient}
\int_Ic_{\rm en}^{\rm b}(s)\,ds
\le C_r\big(\nu {T_1}+\nu\varepsilon_N{T_1}^3+\varepsilon_N
+a^2L_N\nu^{2{\sigma_{\rm int}}-1}\big)=o(1).
\end{equation}
Multiply \eqref{eq:fullEnergyDifferential} {by the corresponding integrating factor.} Since
$\int_I c_{\rm en}^{\rm b}=o(1)$, integration gives
\[
\begin{aligned}
\sup_{I}\mathcal E_r[y]+\int_I\mathcal D_r[y](s)\,ds
&\le C_r\mathcal E_r[y(t_0)]\\
&\quad+C_rT_1[1+\varepsilon_N(1+T_1)]
\|\rho[H^\nu]\|_{L^\infty(I;H_x^{r-1})}\sup_{I}\mathcal E_r[y]^{1/2}.
\end{aligned}
\]
Young's inequality gives
\[
\frac12\sup_I\mathcal E_r[y]+\int_I\mathcal D_r[y]\,ds
\le{}C_r\mathcal E_r[y(t_0)] +C_rT_1^2[1+\varepsilon_N(1+T_1)]^2
\|\rho[H^\nu]\|_{L^\infty(I;H_x^{r-1})}^2.
\]
Taking square roots and using the multiplier equivalence proves
\eqref{eq:fullLinearEnergyClosed}, also for $f$.
\par

\end{proof}
\par\endgroup
Dividing the energy estimate by $1+T_1$ explains the normalization in
\eqref{eq:errorComparisonNorm}. The remaining coefficient
$1+\varepsilon_N(1+T_1)$ is subexponential in $N$, since
$\varepsilon_N=aL_N$, $a(1+T_1)\sim N^2$ and $\log L_N=o(N)$.
\par

\begin{lem}\label{lem:coarseLossVersusEcho}
Assume \eqref{eq:parameters} and \eqref{eq:comparisonLogParameters}.
Let $C_*$ be the exponent obtained by the weighted contraction and
norm conversion in the proof of
Proposition~\ref{prop:coarseDensityInverse}. Then
\begin{equation}\label{eq:coarseLossVersusEcho}
C_*\ge\pi\sqrt{C_{\rm maj}/{\bf b}_1},\qquad
C_{\rm maj}\ge4\Lambda,\qquad C_*>\pi\sqrt3\,c_{\rm e}.
\end{equation}
In particular, for $1/2<\sigma_{\rm int}<1$ and $C_{\rm B}>c_{\rm e}$,
there is no $K>0$ satisfying both
\begin{equation}\label{eq:directFieldIncompatibleConditions}
Kc_{\rm e}\ge\tfrac13,\qquad
\sigma_{\rm int}-\tfrac12-K(C_*+C_{\rm B}-c_{\rm e})\ge0.
\end{equation}
These assertions concern the constant delivered by the stated proof;
they do not give a lower bound on the norm of the actual inverse.
\end{lem}
\begingroup
The proof is given in Appendix~\ref{app:coarseLossVerification}.
\par\endgroup

\par

\begingroup
{Lemma~\ref{lem:collisionSourceEstimates} constructs the solution of the collision equation to
{sources satisfying the second condition in \eqref{eq:quadraticSourceSpaces}} and bounds the resulting density by an energy
estimate applied separately to each Fourier mode of the initial datum and sources. We combine this bound and the
strong-source trace with the energy inequality.}
The size of the background field on an interval $I'$ is measured by
\begin{equation}\label{eq:coarseFieldSize}
\alpha_r(I')=
\int_{I'}\|E[F^{\rm b,\nu}]\|_{C^r_x}\,ds
+\left(\nu^{-1}\int_{I'}\|E[F^{\rm b,\nu}]\|_{C^r_x}^2\,ds\right)^{1/2}.
\end{equation}
{Recall from \eqref{eq:nonlinearBilinearOperator} that
\[
\mathbb Q_\nu(F_1,F_2)
=-E[F_1]\partial_{v_1}F_2+\nu Q(F_1,F_2).
\]
We use this notation below, with $\mathcal A_\nu$ denoting the complete
background linearized operator defined in
\eqref{eq:fullPhysicalOperatorDefinition}.}
\endgroup
\begingroup
We estimate the electric source by velocity integration by parts,
then apply Lemma~\ref{lem:collisionSourceEstimates} to the zero-data
collision equation. This gives the source contribution to the
density equation before the self-consistent density equation is solved.
\par\endgroup
\begingroup
The remainder equation will contain $-\mathcal T_{\rm b}H_J$.
The weighted dual norm is
\[
\|S\|_{\Sigma_\ell^*}
=\sup_{\varphi\in C_c^\infty,\ \|\varphi\|_{\Sigma_\ell}\le1}
|\langle S,\varphi\rangle_\ell|.
\]
The energy $\mathcal E_r$ contains spatial derivatives and velocity
weights, but no ordinary velocity derivative. We integrate the velocity
derivative in $E\partial_{v_1}H_2$ onto the test function and use the
collision form. The factor $\nu$ in $\mathcal D_r$ gives a cost
$\nu^{-1/2}$, while the additional velocity moment is supplied by
$\mathcal E_{r+1}$. The proof below records both costs explicitly.
\par\endgroup
The equations \eqref{eq:fieldPartialSumIdentities} and \eqref{eq:fieldFiniteUnknownError} lead to the prescribed-coefficient problem
\[
\mathcal A_\nu H-\mathbb Q_\nu({\bf H}+{\bf W},H)
-\mathbb Q_\nu(H,{\bf H})=S_0+S_1,\qquad H(t_0)=H_{\rm in}.
\]
For the linear remainder, ${\bf H}={\bf W}=0$, $S_0=0$, and
$S_1=-\mathcal T_{\rm b}H_J$ contains a velocity derivative. The next lemma estimates this source in the dual collision norm.

\begin{lem}\label{lem:electricFormDual}
Let $r\ge2$ be an integer and $I'\subset I$ an interval. Suppose that
$E\in L^2(I';C_x^r)$ is independent of $v$ and
${\sup_{I'}\mathcal E_{r+1}[H_2/\sqrt\mu]<\infty}$. Then
\begin{equation}\label{eq:electricL2SourceBound}
\nu^{-1/2}\sum_{j=0}^r
\|\partial_{x_1}^j(\mu^{-1/2}E\partial_{v_1}H_2)\|_{L^2_{s,x}(I')\Sigma_{2r+6-2j}^*}
\le C_r\nu^{-1/2}\|E\|_{L^2(I';C_x^r)}
\sup_{I'}\mathcal E_{r+1}[H_2/\sqrt\mu]^{1/2}.
\end{equation}
In particular, if $\rho[H_1]\in\ell^1_{r-1}L^\infty(I')$, then
\begin{equation}\label{eq:electricFormDualBound}
\begin{aligned}
&\nu^{-1/2}\sum_{j=0}^r
\|\partial_{x_1}^j(\mu^{-1/2}E[H_1]\partial_{v_1}H_2)\|_{L^2_{s,x}(I')\Sigma_{2r+6-2j}^*}\\
&\qquad\le C_r\sqrt{|I'|/\nu}\,
\|\rho[H_1]\|_{\ell^1_{r-1}L^\infty(I')}
\sup_{I'}\mathcal E_{r+1}[H_2/\sqrt\mu]^{1/2}.
\end{aligned}
\end{equation}
The physical sources ${E\partial_{v_1}H_2}$ and
${E[H_1]\partial_{v_1}H_2}$ are defined distributionally and have zero
velocity integral. For the corresponding normalized source, this means
\[
\int_{\mathbb R^3}E\partial_{v_1}H_2\,dv
=\int_{\mathbb R^3}\sqrt\mu\big(\mu^{-1/2}E\partial_{v_1}H_2\big)\,dv=0,
\]
and the same identity holds with $E=E[H_1]$.

Under the background and parameter assumptions of
Lemma~\ref{lem:collisionSourceEstimates}, let $I'=[\tau_*,\tau^*]$ and
let $g\in C(I';\mathscr E_r)$ be the energy solution of
\begin{equation}\label{eq:electricSourceCollisionEquation}
(\partial_s+v_1\partial_{x_1})g
=\nu\mathcal C_{\rm b}(s)g+\frac{E\partial_{v_1}H_2}{\sqrt\mu},
\qquad g(\tau_*)=0.
\end{equation}
Such a solution exists uniquely and satisfies
\begin{equation}\label{eq:electricSourceDensityResponse}
\begin{aligned}
&\sum_{k\ne0}\langle k\rangle^{r-1}\sup_{s\in I'}
\left|\int\sqrt\mu g_k(s)\,dv\right| \le C_r\nu^{-1/2}\sum_{j=0}^r
\|\partial_{x_1}^j(\mu^{-1/2}E\partial_{v_1}H_2)\|_{L^2_{s,x}(I')\Sigma_{2r+6-2j}^*}\\
&\qquad\le C_r\nu^{-1/2}\|E\|_{L^2(I';C_x^r)}
\sup_{I'}\mathcal E_{r+1}[H_2/\sqrt\mu]^{1/2}.
\end{aligned}
\end{equation}
For $E=E[H_1]$, the last expression is bounded by the
right-hand side of \eqref{eq:electricFormDualBound}.
\end{lem}

\begingroup
\begin{proof}
\emph{1. The weighted pairing and the source integral.}
Put $f=H_2/\sqrt\mu$. Then
${\mu^{-1/2}E\partial_{v_1}H_2=E(\partial_{v_1}-v_1/2)f}$.
First take $f(s,x_1,\cdot)\in\mathcal S(\mathbb R^3)$ and $\varphi\in C_c^\infty(\mathbb R^3)$, with $(s,x_1)$ fixed. Then
\[
\langle(\partial_{v_1}-v_1/2)f,\varphi\rangle_\ell
=-\left\langle f,
\left(\partial_{v_1}+\frac{v_1}{2}
+\frac{2\ell v_1}{\langle v\rangle^2}\right)\varphi\right\rangle_\ell.
\]
Here the functions depend only on $v$, with ${(s,x_1)}$ fixed.
The lower bound $\sigma(v)\ge c\langle v\rangle^{-3}\mathrm{Id}$ gives
\[
\begin{aligned}
|\langle(\partial_{v_1}-v_1/2)f,\varphi\rangle_\ell|
&\le C_\ell\|\langle v\rangle^{\ell+3/2}f\|_2
\left(\int\langle v\rangle^{2\ell-3}
\big(|\nabla_v\varphi|^2+|v|^2|\varphi|^2\big)\,dv\right)^{1/2}\\
&\le C_\ell\|\langle v\rangle^{\ell+3/2}f\|_2\|\varphi\|_{\Sigma_\ell}.
\end{aligned}
\]
Taking the supremum in \eqref{eq:quadraticDualNorm} yields
\begin{equation}\label{eq:electricDualSupremum}
\|(\partial_{v_1}-v_1/2)f\|_{\Sigma_\ell^*}
=\sup_{\substack{\varphi\in C_c^\infty(\mathbb R^3)\\
\|\varphi\|_{\Sigma_\ell}\le1}}
|\langle(\partial_{v_1}-v_1/2)f,\varphi\rangle_\ell| \le C_\ell\|\langle v\rangle^{\ell+3/2}f\|_2.
\end{equation}
Since $E$ is independent of $v$,
\[
\partial_{x_1}^j\big(E(\partial_{v_1}-v_1/2)f\big)
=\sum_{q=0}^j\binom jq(\partial_{x_1}^qE)
(\partial_{v_1}-v_1/2)\partial_{x_1}^{j-q}f.
\]
For $0\le j\le r$, \eqref{eq:electricDualSupremum} therefore gives
\[
\begin{aligned}
&\|\partial_{x_1}^j(\mu^{-1/2}E\partial_{v_1}H_2)\|_{L_x^2\Sigma_{2r+6-2j}^*} \le C_r\sum_{q=0}^j\|\partial_{x_1}^qE\|_{L_x^\infty}
\|\langle v\rangle^{2r+6-2j+3/2}\partial_{x_1}^{j-q}f\|_2\\
&\quad\le C_r\|E\|_{C_x^r}\sum_{q=0}^j
\|\langle v\rangle^{2(r+1)+6-2(j-q)}\partial_{x_1}^{j-q}f\|_2
\le C_r\|E\|_{C_x^r}\mathcal E_{r+1}[f]^{1/2},
\end{aligned}
\]
where
\[
(2(r+1)+6-2(j-q))-(2r+6-2j+3/2)=\tfrac12+2q>0.
\]
The same estimates and weighted duality give, at each time,
\begin{equation}\label{eq:electricWeightedEnergyPairing}
\left|\sum_{j=0}^r\left\langle
\partial_{x_1}^j\{E(\partial_{v_1}-v_1/2)f\},
\partial_{x_1}^jg\right\rangle_{2r+6-2j}\right| \le C_r\nu^{-1/2}\|E\|_{C_x^r}
\mathcal E_{r+1}[f]^{1/2}\mathcal D_r[g]^{1/2}
\end{equation}
for every test function $g$ with finite $\mathcal D_r[g]$.
For each $j$, time integration gives
\[
\begin{aligned}
&\|\partial_{x_1}^j(\mu^{-1/2}E\partial_{v_1}H_2)\|_{L^2_{s,x}(I')\Sigma_{2r+6-2j}^*}^{\,2} =\int_{I'}
\|\partial_{x_1}^j(\mu^{-1/2}E\partial_{v_1}H_2)(s)\|_{L_x^2\Sigma_{2r+6-2j}^*}^{\,2}\,ds\\
&\quad\le C_r\int_{I'}\|E(s)\|_{C_x^r}^2\mathcal E_{r+1}[f(s)]\,ds \le C_r\|E\|_{L^2(I';C_x^r)}^2
\sup_{I'}\mathcal E_{r+1}[H_2/\sqrt\mu].
\end{aligned}
\]
Taking square roots, summing over $j$ and multiplying by $\nu^{-1/2}$
proves \eqref{eq:electricL2SourceBound}.
The Poisson multiplier gives
\[
\begin{aligned}
\|E[H_1](s)\|_{C_x^r}
&\le C_r\sum_{k\ne0}\langle k\rangle^{r-1}|\rho_k[H_1](s)|,\\
\|E[H_1]\|_{L^2(I';C_x^r)}
&\le C_r|I'|^{1/2}\sum_{k\ne0}\langle k\rangle^{r-1}
\sup_{s\in I'}|\rho_k[H_1](s)|
=C_r|I'|^{1/2}\|\rho[H_1]\|_{\ell^1_{r-1}L^\infty(I')}.
\end{aligned}
\]
Substitution proves \eqref{eq:electricFormDualBound}.
The weighted pairing defines the source for functions $f$ and $E$ in the stated spaces.
Approximate $f$ in ${x_1,v}$ with uniformly bounded
${\mathcal E_{r+1}}$ and approximate $E$ in ${L^2(I';C_x^r)}$.
The displayed product bounds and dominated convergence give convergence
of each differentiated source in
${L^2_{s,x}(I')\Sigma_{2r+6-2j}^*}$.
Velocity cutoffs and $H_2=\sqrt\mu f$ also give
$\int E\partial_{v_1}H_2\,dv=0$ in distributions.

\emph{2. The density produced by the source.}
Apply Lemma~\ref{lem:collisionSourceEstimates} with zero initial datum,
$S_0=0$ and ${S_1=\mu^{-1/2}E\partial_{v_1}H_2}$.
The required source integrals are finite by
\eqref{eq:electricL2SourceBound}.
For smooth sources, Duhamel's formula and velocity integration give
\begin{equation}\label{eq:electricSourceDensityDuhamel}
\begin{aligned}
g(s)&=\int_{\tau_*}^s\mathbb U_{\rm b}^\nu(s,\tau)
\left[\frac{E(\tau)\partial_{v_1}H_2(\tau)}{\sqrt\mu}\right]\,d\tau,\\
\rho_k[\sqrt\mu g](s)
&=\int\sqrt\mu\left[\int_{\tau_*}^s\mathbb U_{\rm b}^\nu(s,\tau)
\left[\frac{E(\tau)\partial_{v_1}H_2(\tau)}{\sqrt\mu}\right]\,d\tau\right]_k dv.
\end{aligned}
\end{equation}
In the energy identity, weighted duality and Cauchy--Schwarz give
\par

\[
\begin{aligned}
&\left|\sum_{j=0}^r\int_{I'}
\left\langle\partial_{x_1}^j(\mu^{-1/2}E\partial_{v_1}H_2),
\partial_{x_1}^jg\right\rangle_{2r+6-2j}\,ds\right|\\
&\quad\le\sum_{j=0}^r\int_{I'}
\|\partial_{x_1}^j(\mu^{-1/2}E\partial_{v_1}H_2)\|_{L_x^2\Sigma_{2r+6-2j}^*}
\|\partial_{x_1}^jg\|_{L_x^2\Sigma_{2r+6-2j}}\,ds\\
&\quad\le\nu^{-1/2}\sum_{j=0}^r
\|\partial_{x_1}^j(\mu^{-1/2}E\partial_{v_1}H_2)\|_{L^2_{s,x}(I')\Sigma_{2r+6-2j}^*}
\left(\int_{I'}\mathcal D_r[g]\,ds\right)^{1/2}\\
&\quad\le\frac{c_r}{2}\int_{I'}\mathcal D_r[g]\,ds
+\frac{C_r}{\nu}\left(\sum_{j=0}^r
\|\partial_{x_1}^j(\mu^{-1/2}E\partial_{v_1}H_2)\|_{L^2_{s,x}(I')\Sigma_{2r+6-2j}^*}\right)^2.
\end{aligned}
\]
The first term is absorbed by collision dissipation. The energy and
density bounds in Lemma~\ref{lem:collisionSourceEstimates} therefore give
\par

\[
\begin{aligned}
&\|g\|_{\mathcal N_r^{I'}}
+\sum_{k\ne0}\langle k\rangle^{r-1}\sup_{s\in I'}
\left|\int\sqrt\mu g_k(s)\,dv\right|\\
&\qquad\le C_r\nu^{-1/2}\sum_{j=0}^r
\|\partial_{x_1}^j(\mu^{-1/2}E\partial_{v_1}H_2)\|_{L^2_{s,x}(I')\Sigma_{2r+6-2j}^*}.
\end{aligned}
\]
Steps 2--3 of that lemma's proof give the time supremum inside the
mode sum. For distributional sources, Step 4 supplies convergence
of smooth-source solutions in the norm $\|\cdot\|_{\mathcal N_r^{I'}}$ and in the
displayed density norm, and identifies the density limit with
$\rho[\sqrt\mu g]$. Thus \eqref{eq:electricSourceDensityDuhamel}
is the actual density produced by the source, and \eqref{eq:electricSourceDensityResponse}
follows. No pointwise action of $\mathbb U_{\rm b}^\nu(s,\tau)$ on
a distributional source is required.
\par
\end{proof}
\par
\par\endgroup
\begin{prop}\label{thm:fullForcedCoarse}\label{lem:splitFieldSize}\label{prop:prescribedRemainder}
Use the logarithmic parameters of this subsection. Fix ${\rm s}\ge3$ and let
$r\in\{{\rm s}-1,{\rm s}\}$. There is a factor $A_{\rm s}(N)\ge1$, with
$\log A_{\rm s}(N)=o(N)$, such that the following holds. Let $\mathfrak r\ge0$
and let {${{\bf H}},{{\bf W}}$} be prescribed real distributions of zero total mass satisfying
\begin{equation}\label{eq:prescribedRemainderClass}
\begin{aligned}
&{{\bf H}}/\sqrt \mu ,{{\bf W}}/\sqrt \mu \in C_{\rm w}(I;\mathscr E_{\rm s}),\qquad
\max\{\|{{\bf H}}\|_{\mathscr V_{\rm s}},\|{{\bf W}}\|_{\mathscr V_{\rm s}}\}\le\mathfrak r,\\
&{\sup_I\mathcal E_{{\rm s}+1}[{{\bf H}}/\sqrt\mu]^{1/2}
\le(1+T_1)\mathfrak r,}\\
&\partial_s\rho[{{\bf H}}]+\partial_{x_1}j[{{\bf H}}]=0,\qquad
\partial_s\rho[{{\bf W}}]+\partial_{x_1}j[{{\bf W}}]=0.
\end{aligned}
\end{equation}
The continuity equations are distributional. Assume that
\begin{equation}\label{eq:unifiedCoefficientSmallness}
A_{\rm s}e^{C_*N}\mathfrak r
\big((1+T_1)^2+\nu^{-1/2}(1+T_1)^{3/2}\big)=o(1).
\end{equation}
\par

{Let {$H_{\rm in}$}, $S_0$ and $S_1$ have zero total mass, and assume that}
\[
Y_r=\mathcal E_r[H_{\rm in}/\sqrt \mu]^{1/2}
+\int_{I}\mathcal E_r[\mu^{-1/2}S_0]^{1/2}\,ds+\nu^{-1/2}\sum_{j=0}^{r}\left\|\partial_{x_1}^j\!\left(\mu^{-1/2}S_1\right)\right\|_{L^2_{s,x}(I)\Sigma_{2r+6-2j}^*}<\infty{.}
\]
{Let {$H$} be a finite-energy solution of}
\begin{equation}\label{eq:fullForcedEquation}
{
\mathcal A_\nu H-\mathbb Q_\nu({\bf H}+{\bf W},H)
-\mathbb Q_\nu(H,{\bf H})=S_0+S_1,\qquad H(t_0)=H_{\rm in}.
}
\end{equation}
{Then}
\begin{equation}\label{eq:prescribedRemainderBound}
\|H\|_{\mathscr V_r(I)}\le A_{\rm s}(N)e^{C_*N}Y_r.
\end{equation}
{For every initial datum and pair of sources satisfying the stated
conditions, the initial-value problem has a unique finite-energy
solution, with {$H/\sqrt \mu \in C(I;\mathscr E_r)$}.
For the physical source $S_1$, zero total mass means, for almost every $s$,}
\[
\int_{\mathbb T}
\big\langle \mu^{-1/2}S_1(s,x_1,\cdot),\sqrt \mu \big\rangle\,dx_1=0.
\]
This unweighted pairing is defined by continuous extension: in the
weighted duality \eqref{eq:quadraticDualNorm}, it uses the test function
{$\langle v\rangle^{-4r-12}\sqrt \mu \in\Sigma_{2r+6}$}.
The estimate holds on all $[t_0,T']\subset I$ with the original final-time factors, uniformly
in the coefficient class. If the total source has zero velocity integral
at each spatial point, then
\begin{equation}\label{eq:prescribedRemainderContinuity}
\partial_s\rho[H]+\partial_{x_1}j[H]=0.
\end{equation}
For {${{\bf H}}={{\bf W}}=0$}, at every fixed $r\ge2$, the sharper bounds are
\begin{equation}\label{eq:fullForcedDensityBound}
\begin{aligned}
\|\rho[H]\|_{\ell^1_{r-1}L^\infty(I)}
\le C_rN^4e^{C_*N}\Bigg(&\mathcal E_r[H_{\rm in}/\sqrt\mu]^{1/2}
+\int_I\mathcal E_r[S_0/\sqrt\mu]^{1/2}\,ds\\
&+\nu^{-1/2}\sum_{j=0}^r
\|\partial_{x_1}^j(S_1/\sqrt\mu)\|_{L^2_{s,x}(I)\Sigma_{2r+6-2j}^*}\Bigg).
\end{aligned}
\end{equation}
\begin{equation}\label{eq:fullForcedEnergyBound}
\begin{aligned}
&\|H/\sqrt\mu\|_{\mathcal N_r} \le C_rN^4(1+T_1)[1+\varepsilon_N(1+T_1)]e^{C_*N}
\Bigg(\mathcal E_r[H_{\rm in}/\sqrt\mu]^{1/2}\\
&\qquad\quad+\int_I\mathcal E_r[S_0/\sqrt\mu]^{1/2}\,ds
+\nu^{-1/2}\sum_{j=0}^r
\|\partial_{x_1}^j(S_1/\sqrt\mu)\|_{L^2_{s,x}(I)\Sigma_{2r+6-2j}^*}\Bigg).
\end{aligned}
\end{equation}
The coefficient $C_*$ is the one in Proposition~\ref{prop:coarseDensityInverse};
it is independent of all fixed regularity and correction orders.
\end{prop}

\begingroup
\begin{proof}

Use $\varepsilon_N=aL_N$ as in Proposition~\ref{lem:fixedCoefficientOrders}.
Put {${F^{\rm c}_1}=F^{\rm b,\nu}+{{\bf H}}+{{\bf W}}$} and {${F^{\rm c}_2}=F^{\rm b,\nu}+{{\bf H}}$}.
Expanding \eqref{eq:nonlinearBilinearOperator} gives the equation
\begin{equation}\label{eq:prescribedRemainderExpanded}
(\partial_s+v_1\partial_{x_1}+E[{F^{\rm c}_1}]\partial_{v_1})H
=\nu\big(Q({F^{\rm c}_1},H)+Q(H,{F^{\rm c}_2})\big)
-E[H]\partial_{v_1}{F^{\rm c}_2}+S_0+S_1.
\end{equation}
In particular, {${{\bf W}}$} does not occur in the differentiated coefficient
{$\partial_{v_1}{F^{\rm c}_2}$}.

We verify the coefficient norms first. The density and continuity assumptions imply
\begin{equation}\label{eq:prescribedPotentialBounds}
\sup_I(\|E[{{\bf W}}]\|_{C^{\rm s}_x}+\|\phi[{{\bf W}}]\|_{C^{{\rm s}+1}_x})\le C_{\rm s}\mathfrak r,
\qquad
\int_{I'}\|\partial_s\phi[{{\bf W}}]\|_{C^{\rm s}_x}\,ds\le C_{\rm s}|I'|(1+{T_1})\mathfrak r.
\end{equation}
For $k\ne0$, the continuity equation and Poisson formula give
\[
\begin{aligned}
\partial_s\phi_k[{\bf W}]&=-ij_k[{\bf W}]/k,\qquad \|\partial_s\phi[{\bf W}]\|_{C_x^{\rm s}}
\le C_{\rm s}\sum_{k\ne0}\langle k\rangle^{{\rm s}-1}|j_k[{\bf W}]|
\le C_{\rm s}\|j[{\bf W}]\|_{H_x^{\rm s}}\\
&\le C_{\rm s}\mathcal E_{\rm s}[{\bf W}/\sqrt\mu]^{1/2}
\le C_{\rm s}(1+T_1)\mathfrak r.
\end{aligned}
\]
Time integration gives \eqref{eq:prescribedPotentialBounds}; the same
calculation with ${\bf W}$ replaced by ${\bf H}$ gives its potential bound.
\par

We first construct the energy solution and justify the density norm required below. Apply Lemma~\ref{lem:physicalDerivativeEnergy} with ${\varpi}=1$ and
\begin{equation}\label{eq:remainderCoefficientRealization}
E=E[{F^{\rm c}_1}],\qquad
\mathcal J_sU=-E[U]\partial_{v_1}{F^{\rm c}_2}.
\end{equation}
The coefficient assumptions follow from
\[
\begin{aligned}
\max_{i=1,2}\sup_I
\mathcal E_{\rm s}[(F_i^{\rm c}-\mu)/\sqrt\mu]^{1/2}
&\le C_{\rm s}[\varepsilon_N+(1+T_1)\mathfrak r]=o(1),\\
\sum_{i=1}^2\int_I\mathcal D_{\rm s}[(F_i^{\rm c}-\mu)/\sqrt\mu]\,ds
&\le C_{\rm s}\int_I\mathcal D_{\rm s}[f^{\rm b,\nu}]\,ds
+C_{\rm s}(1+T_1)^2\mathfrak r^2<\infty,
\end{aligned}
\]
and from \eqref{eq:prescribedPotentialBounds} and the background
potential estimates.
\par
 
For the operator in \eqref{eq:remainderCoefficientRealization}, take
$\mathcal J_s^0U=-E[U]\partial_{v_1}F^{\rm b,\nu}$ and
$\mathcal J_s^1U=-E[U]\partial_{v_1}{\bf H}$ in
Lemma~\ref{lem:physicalDerivativeEnergy}. The smooth background bounds
the normalized first term in $\mathscr E_r$. For the second, integration by parts
as in Lemma~\ref{lem:electricFormDual} gives pointwise
\[
\sum_{j=0}^{r}\left\|\partial_{x_1}^j\!\left(\frac{E[U]\partial_{v_1}{\bf H}}{\sqrt\mu}\right)\right\|_{L_x^2\Sigma_{2r+6-2j}^*}
\le C_{\rm s}\mathcal E_{r+1}[{\bf H}/\sqrt\mu]^{1/2}
\mathcal E_r[U/\sqrt\mu]^{1/2}.
\]
The coefficient on the right belongs to $L^2(I)$ by
\eqref{eq:prescribedRemainderClass}. Thus the lemma constructs a unique
strongly continuous energy solution for each fixed $\nu$.
\par

For this solution put $f=H/\sqrt\mu$. Gaussian integration and the
one-dimensional Poisson estimate give only the fixed-time bound
\[
\sup_I\big(\|\rho[H]\|_{H^r_x}+\|E[H]\|_{C^r_x}\big)
\le C_r\sup_I\mathcal E_r[f]^{1/2}<\infty.
\]
The pointwise estimate just proved therefore implies
\[
\begin{aligned}
&\nu^{-1/2}\sum_{j=0}^{r}\left\|\partial_{x_1}^j\!\left(\mu^{-1/2}\bigl(E[H]\partial_{v_1}{\bf H}\bigr)\right)\right\|_{L^2_{s,x}(I)\Sigma_{2r+6-2j}^*}\\
&\quad\le C_r\sqrt{|I|/\nu}\,
\sup_I\mathcal E_r[f]^{1/2}
\sup_I\mathcal E_{r+1}[{\bf H}/\sqrt\mu]^{1/2}<\infty.
\end{aligned}
\]
The smooth background also gives
${\int_{I}\mathcal E_r[\mu^{-1/2}\bigl(E[H]\partial_{v_1}F^{\rm b,\nu}\bigr)]^{1/2}\,ds<\infty}$.
For the transport and collision sources, the spatial product rule,
$2r+6-2j\ge6$ and Lemma~\ref{lem:quadraticLandauForm} give
\[
\begin{aligned}
&\int_I\|E[F^{\rm c}_1](\partial_{v_1}-v_1/2)f\|_{H^r_xL^2_v}\,ds \le C_r\sup_I\|E[F^{\rm c}_1]\|_{C^r_x}
\big(|I|+\sqrt{|I|/\nu}\big)\|f\|_{\mathcal N_r}<\infty,\\
&\nu^{-1/2}\sum_{j=0}^{r}\left\|\partial_{x_1}^j\!\left(\mu^{-1/2}\bigl(\nu Q({\bf H}+{\bf W},H)+\nu Q(H,{\bf H})\bigr)\right)\right\|_{L^2_{s,x}(I)\Sigma_{2r+6-2j}^*}\\
&\quad\le C_r\big(\|{\bf H}/\sqrt\mu\|_{\mathcal N_r}
+\|{\bf W}/\sqrt\mu\|_{\mathcal N_r}\big)\|f\|_{\mathcal N_r}<\infty.
\end{aligned}
\]
These bounds require no time supremum inside a sum over density modes.
Rewrite \eqref{eq:prescribedRemainderExpanded} with only the background
collision operator on the left. Duhamel's formula with
$\mathbb U_{\rm b}^\nu$, interpreted by source completion, and
{Lemma~\ref{lem:collisionSourceEstimates}} yield
\begin{equation}\label{eq:prescribedDensityMembership}
\begin{aligned}
&\|\rho[H]\|_{\ell^1_{r-1}L^\infty(I)}\\
&\quad\le C_r\Bigg(Y_r
+\int_I\mathcal E_r[\mu^{-1/2}E[H]\partial_{v_1}F^{\rm b,\nu}]^{1/2}\,ds\\
&\qquad+\nu^{-1/2}\sum_{j=0}^r
\|\partial_{x_1}^j(\mu^{-1/2}E[H]\partial_{v_1}{\bf H})\|_{L^2_{s,x}(I)\Sigma_{2r+6-2j}^*}\\
&\qquad+\int_I\|E[F^{\rm c}_1](\partial_{v_1}-v_1/2)f\|_{H_x^rL_v^2}\,ds\\
&\qquad+\sqrt\nu\sum_{j=0}^r
\|\partial_{x_1}^j\{\mu^{-1/2}[Q({\bf H}+{\bf W},H)+Q(H,{\bf H})]\}\|_{L^2_{s,x}(I)\Sigma_{2r+6-2j}^*}\Bigg)<\infty.
\end{aligned}
\end{equation}
The strong trace estimate requires only $H^r_xL^2_v$ for the transport
{source; the other sources satisfy the two integral conditions in \eqref{eq:quadraticSourceSpaces}.}
Thus the coupled estimates below apply to the constructed solution.
The qualitative existence constant is not used in their quantitative bounds.
\par

{We now derive the quantitative estimate.} Use the unknown
\[
y=e^{\psi/2}H/\sqrt \mu ,\qquad
\psi=\phi^{\rm b,\nu}+\phi[{{\bf H}}]+\phi[{{\bf W}}].
\]
Since $\partial_{x_1}\psi=-E[F_1^{\rm c}]$, direct substitution gives
\[
\begin{aligned}
(\partial_s+v_1\partial_{x_1}+E[F_1^{\rm c}]\partial_{v_1})y
={}&\nu\mu^{-1/2}\big[Q(F_1^{\rm c},\sqrt\mu y)
+Q(\sqrt\mu y,F_2^{\rm c})\big]+\tfrac12\partial_s\psi\,y\\
&-e^{\psi/2}E[H]\mu^{-1/2}\partial_{v_1}F_2^{\rm c}
+e^{\psi/2}(S_0+S_1)/\sqrt\mu.
\end{aligned}
\]
Here the two multiplication terms cancel:
\[
-\tfrac12v_1\partial_{x_1}\psi\,y
-\tfrac12v_1E[F_1^{\rm c}]y=0.
\]
Proposition~\ref{prop:fullLinearEnergy} bounds the background terms.
At $r={\rm s}-1,{\rm s}$, Lemma~\ref{lem:quadraticLandauForm}
bounds the additional collision pairings by
\par

\[
\begin{aligned}
&C_{\rm s}\big(\mathcal E_r[{{\bf H}}/\sqrt \mu ]^{1/2}
+\mathcal E_r[{{\bf W}}/\sqrt \mu ]^{1/2}\big)\mathcal D_r[y]\\
&\quad+C_{\rm s}\big(\mathcal D_r[{{\bf H}}/\sqrt \mu ]^{1/2}
+\mathcal D_r[{{\bf W}}/\sqrt \mu ]^{1/2}\big)
\mathcal E_r[y]^{1/2}\mathcal D_r[y]^{1/2}.
\end{aligned}
\]

The energy and form bounds in the coefficient ball imply
\begin{align*}
\sup_I\big(\mathcal E_r[{{\bf H}}/\sqrt \mu ]^{1/2}+\mathcal E_r[{{\bf W}}/\sqrt \mu ]^{1/2}\big)
&\le C_{\rm s}(1+{T_1})\mathfrak r=o(1),\\
&C_{\rm s}\big(\mathcal D_r[{{\bf H}}/\sqrt \mu ]^{1/2}+\mathcal D_r[{{\bf W}}/\sqrt \mu ]^{1/2}\big)
\mathcal E_r[y]^{1/2}\mathcal D_r[y]^{1/2}\\
&\quad\le\frac{c_r}{8}\mathcal D_r[y]
+C_{\rm s}\big(\mathcal D_r[{{\bf H}}/\sqrt \mu ]+\mathcal D_r[{{\bf W}}/\sqrt \mu ]\big)\mathcal E_r[y],\\
\int_I \big(\mathcal D_r[{{\bf H}}/\sqrt \mu ]+\mathcal D_r[{{\bf W}}/\sqrt \mu ]\big)\,ds
&\le C_{\rm s}\big(\|{{\bf H}}/\sqrt \mu\|_{\mathcal N_r}^2+\|{{\bf W}}/\sqrt \mu\|_{\mathcal N_r}^2\big)
\le C_{\rm s}(1+{T_1})^2\mathfrak r^2.
\end{align*}
Thus the first pairing is at most ${c_r\mathcal D_r/8}$ for small
$\nu$. The explicit sum in
\eqref{eq:fieldCommutatorPairing} gives
\begin{align*}
&\left|\sum_{j=1}^r\sum_{q=1}^j\binom jq
\langle(\partial_{x_1}^qE[{{\bf H}}+{{\bf W}}])
\partial_{v_1}\partial_{x_1}^{j-q}y,\partial_{x_1}^jy\rangle_{2r+6-2j}\right|\\
&\quad\le\frac{c_r}{8}\mathcal D_r[y]
+C_{\rm s}\nu^{-1}\|E[{{\bf H}}+{{\bf W}}]\|_{C^r_x}^2\mathcal E_r[y].
\end{align*}
Moreover,
\begin{align*}
\int_I\nu^{-1}\|E[{{\bf H}}+{{\bf W}}]\|_{C^r_x}^2\,ds
&\le C_{\rm s}\nu^{-1}(1+{T_1})\mathfrak r^2,\\
\int_I\big(\|E[{{\bf H}}+{{\bf W}}]\|_\infty+
\|\partial_s(\phi[{{\bf H}}]+\phi[{{\bf W}}])\|_{C^r_x}\big)\,ds
&\le C_{\rm s}(1+{T_1})^2\mathfrak r.
\end{align*}
The weight inequality used in the first line is precisely
$2r+6-2j+2q-3/2\ge2r+6-2j$ for $q\ge1$.
The total additional integrated energy coefficient is therefore
\[
C_{\rm s}\big((1+{T_1})^2\mathfrak r+\nu^{-1}(1+{T_1})\mathfrak r^2+(1+{T_1})^2\mathfrak r^2\big)=o(1).
\]
Each smallness follows from \eqref{eq:unifiedCoefficientSmallness}; for example $\nu^{-1}(1+{T_1})\mathfrak r^2=(\mathfrak r\sqrt{(1+{T_1})/\nu})^2$. The background energy coefficient is $o(1)$ by \eqref{eq:fullEnergyIntegratedCoefficient} at the present fixed order.

Lemma~\ref{lem:weightedDualMultipliers} bounds multiplication by
$e^{\psi/2}$ and $e^{-\psi/2}$ in both the energy and dual source norms. The source
pairings therefore satisfy
\[
\begin{aligned}
|\langle e^{\psi/2}S_0/\sqrt \mu ,y\rangle_{\mathscr E_r}|
&\le C_{\rm s}\mathcal E_r[S_0/\sqrt \mu]^{1/2}\mathcal E_r[y]^{1/2},\\
\left|\sum_{j=0}^r\langle\partial_{x_1}^j(e^{\psi/2}S_1/\sqrt\mu),\partial_{x_1}^jy\rangle_{2r+6-2j}\right|
&\le\frac{c_r}{8}\mathcal D_r[y]
+C_{\rm s}\nu^{-1}\sum_{j=0}^{r}\left\|\partial_{x_1}^j\!\left(S_1/\sqrt \mu\right)\right\|_{L_x^2\Sigma_{2r+6-2j}^*}^2.
\end{aligned}
\]
On $I'=[\tau_*,\tau^*]$, set
\[
Y_{r,I'}=\mathcal E_r[H(\tau_*)/\sqrt \mu]^{1/2}
+\int_{I'}\mathcal E_r[\mu^{-1/2}S_0]^{1/2}\,ds+\nu^{-1/2}\sum_{j=0}^{r}\left\|\partial_{x_1}^j\!\left(\mu^{-1/2}S_1\right)\right\|_{L^2_{s,x}(I')\Sigma_{2r+6-2j}^*}.
\]
The background electric term contributes
$C_{\rm s}|I'|[1+\varepsilon_N(1+\tau^*)]\|\rho[H]\|_{\ell^1_{r-1}L^\infty(I')}$
to the energy estimate. Lemma~\ref{lem:electricFormDual} bounds the
additional electric term by
\begin{equation}\label{eq:prescribedElectricDual}
\nu^{-1/2}\sum_{j=0}^{r}\left\|\partial_{x_1}^j\!\left(\mu^{-1/2}\bigl(E[H]\partial_{v_1}{\bf H}\bigr)\right)\right\|_{L^2_{s,x}(I')\Sigma_{2r+6-2j}^*}
\le C_{\rm s}\sqrt{|I'|/\nu}\,(1+T_1)\mathfrak r
\|\rho[H]\|_{\ell^1_{r-1}L^\infty(I')}.
\end{equation}
The energy inequality therefore yields
\begin{equation}\label{eq:prescribedCoupledEnergy}
\begin{aligned}
\|H/\sqrt\mu\|_{\mathcal N_r^{I'}}
&\le C_{\rm s}\big(Y_{r,I'}+\mathfrak p_{I'}\|\rho[H]\|_{\ell^1_{r-1}L^\infty(I')}\big),\\
\mathfrak p_{I'}&=|I'|[1+\varepsilon_N(1+\tau^*)]
+\sqrt{|I'|/\nu}\,(1+T_1)\mathfrak r.
\end{aligned}
\end{equation}
By \eqref{eq:unifiedCoefficientSmallness}, the added term in
$\mathfrak p_{I'}$ is $o(1)$. Thus $\mathfrak p_{I_-}\le C_{\rm s}N^4$
and $\mathfrak p_{I_+}\le C_{\rm s}(1+T_1)[1+\varepsilon_N(1+T_1)]$.
\par

{We next bound the density contributions of the sources using the background collision solution operator.} The extra electric transport satisfies
\begin{equation}\label{eq:prescribedTransportSource}
\int_{I'}\|E[{{\bf H}}+{{\bf W}}](\partial_{v_1}-v_1/2)(H/\sqrt \mu )\|_{H^r_xL^2_v}\,ds
\le C_{\rm s}\mathfrak r(|I'|+\sqrt{|I'|/\nu})\|H/\sqrt \mu\|_{\mathcal N_r^{I'}}.
\end{equation}

{In fact,}
\begin{align*}
&\partial_{x_1}^j\big(E[{{\bf H}}+{{\bf W}}](\partial_{v_1}-v_1/2)f\big)\\
&\quad=\sum_{q=0}^j\binom jq(\partial_{x_1}^qE[{{\bf H}}+{{\bf W}}])
(\partial_{v_1}-v_1/2)\partial_{x_1}^{j-q}f,\\
&\|E[{{\bf H}}+{{\bf W}}](\partial_{v_1}-v_1/2)f\|_{H^r_xL^2_v}\\
&\quad\le C_{\rm s}\mathfrak r\sum_{j=0}^r
\big(\|\partial_{v_1}\partial_{x_1}^jf\|_2+\|v_1\partial_{x_1}^jf\|_2\big)\\
&\quad\le C_{\rm s}\mathfrak r\big(\nu^{-1/2}\mathcal D_r[f]^{1/2}+\mathcal E_r[f]^{1/2}\big),\\
&\int_{I'}\big(\nu^{-1/2}\mathcal D_r[f]^{1/2}+\mathcal E_r[f]^{1/2}\big)\,ds\\
&\quad\le\sqrt{|I'|/\nu}\left(\int_{I'}\mathcal D_r[f]\,ds\right)^{1/2}
+|I'|\sup_{I'}\mathcal E_r[f]^{1/2}\\
&\quad\le(\sqrt{|I'|/\nu}+|I'|)\|H/\sqrt \mu\|_{\mathcal N_r^{I'}}.
\end{align*}
Here $2r+6-2j\ge6$ controls the unweighted derivative and multiplication
by $v_1$. Only this unweighted velocity source norm is used for the trace.

The collision source obeys
\begin{equation}\label{eq:prescribedCollisionSource}
\begin{aligned}
&\nu^{-1/2}\sum_{j=0}^{r}\left\|\partial_{x_1}^j\!\left(\mu^{-1/2}\bigl(\nu Q({{\bf H}}+{{\bf W}},H)+\nu Q(H,{{\bf H}})\bigr)\right)\right\|_{L^2_{s,x}(I')\Sigma_{2r+6-2j}^*}\\
&\quad\le C_{\rm s}\big(\|{{\bf H}}/\sqrt \mu\|_{\mathcal N_r}+\|{{\bf W}}/\sqrt \mu\|_{\mathcal N_r}\big)
\|H/\sqrt \mu\|_{\mathcal N_r^{I'}} \le C_{\rm s}(1+{T_1})\mathfrak r\|H/\sqrt \mu\|_{\mathcal N_r^{I'}},
\end{aligned}
\end{equation}
by Lemma~\ref{lem:quadraticLandauForm} and time Cauchy--Schwarz.
For the remaining electric source, the corresponding density is obtained by
solving the zero-data collision equation with normalized source
$E[H]\partial_{v_1}{\bf H}/\sqrt\mu$ and then integrating in velocity.
Applying \eqref{eq:electricSourceDensityResponse}, which is the
electric-source case of Lemma~\ref{lem:collisionSourceEstimates}, and
then \eqref{eq:prescribedElectricDual}, we obtain
\begin{equation}\label{eq:prescribedElectricDensityResponse}
\begin{aligned}
&\sum_{k\ne0}\langle k\rangle^{r-1}\sup_{s\in I'}
\left|\int\sqrt\mu
\left[\int_{\tau_*}^s\mathbb U_{\rm b}^\nu(s,\tau)
\left[\frac{E[H](\tau)\partial_{v_1}{\bf H}(\tau)}{\sqrt\mu}\right]
d\tau\right]_k dv\right|\\
&\qquad\le C_r\nu^{-1/2}\sum_{j=0}^{r}\left\|\partial_{x_1}^j\!\left(\mu^{-1/2}\bigl(E[H]\partial_{v_1}{\bf H}\bigr)\right)\right\|_{L^2_{s,x}(I')\Sigma_{2r+6-2j}^*}\\
&\qquad\le C_{\rm s}\sqrt{|I'|/\nu}\,(1+T_1)\mathfrak r
\|\rho[H]\|_{\ell^1_{r-1}L^\infty(I')}.
\end{aligned}
\end{equation}
{Thus the {source estimate} controls this density,}
before the self-consistent density equation is solved. The same
source-to-density estimate in Lemma~\ref{lem:collisionSourceEstimates}
applies to the collision source; the $L^1H_x^rL_v^2$ part of that lemma
applies to \eqref{eq:prescribedTransportSource}.
Combining these source and density estimates and applying
Proposition~\ref{prop:coarseDensityInverse} gives, with $A_-=1$ and
$A_+=e^{C_*N}$,
\par
\par

\begin{equation}\label{eq:approximationCoupledDensity}
\begin{aligned}
\|\rho[H]\|_{\ell^1_{r-1}L^\infty(I')}
\le C_{\rm s}A_\pm\Big[&Y_{r,I'}+\big\{\alpha_r(I')
+\mathfrak r(1+{T_1}+\sqrt{|I'|/\nu})\big\}\|H/\sqrt \mu\|_{\mathcal N_r^{I'}}\\
&+{\sqrt{|I'|/\nu}}(1+{T_1})\mathfrak r\|\rho[H]\|_{\ell^1_{r-1}L^\infty(I')}\Big].
\end{aligned}
\end{equation}

For fixed $r$ and $\sigma_{\rm int}\in(1/2,1)$,
Proposition~\ref{prop:backgroundFieldBounds}, after {reserving a fixed fraction of the background-field}
spatial radius, gives
\[
\|E[F^{\rm b,\nu}](t_0+t)\|_{C^r_x}
\le C_r\varepsilon_{\rm F}e^{-c_{\rm dec}t^\chi}.
\]
For $L\ge0$, the two required tail estimates are
\begin{align*}
\int_L^\infty e^{-c_{\rm dec} t^\chi}\,dt
&\le e^{-c_{\rm dec}L^\chi/2}\int_0^\infty e^{-c_{\rm dec} t^\chi/2}\,dt
\le C_{c_{\rm dec},\chi}e^{-c_{\rm dec}L^\chi/2},\\
\left(\int_L^\infty e^{-2c_{\rm dec} t^\chi}\,dt\right)^{1/2}
&\le e^{-c_{\rm dec}L^\chi/2}\left(\int_0^\infty e^{-c_{\rm dec} t^\chi}\,dt\right)^{1/2}
\le C_{c_{\rm dec},\chi}e^{-c_{\rm dec}L^\chi/2}.
\end{align*}
Taking $L=0$ on $I_-$ and $L=N^4$ on $I_+$ yields
\begin{equation}\label{eq:coarseFieldSplit}
\begin{aligned}
\alpha_r(I_-)&\le C_r\varepsilon_{\rm F}(1+\nu^{-1/2})
\le C_r\nu^{-1/2}\varepsilon_{\rm F},\\
\alpha_r(I_+)&\le C_r\varepsilon_{\rm F}(1+\nu^{-1/2})e^{-c_{\rm dec}N^{4\chi}/2}
\le C_r\nu^{-1/2}\varepsilon_{\rm F}e^{-cN^{4\chi}}.
\end{aligned}
\end{equation}
The Gaussian-small term satisfies the time-decay bound in \eqref{eq:fieldFourierEnvelope}, so the same tail integral applies to it.
Using \eqref{eq:fieldExpansionSizes}, $a=\nu^{1/3-o(1)}$,
$L_N=\nu^{-o(1)}$ and $N=K\log(1/\nu)$, we obtain
\begin{equation}\label{eq:logFieldSplitSmallness}
\nu^{-1/2}\varepsilon_{\rm F}
\le C\nu^{-1/2}aL_N\nu^{\sigma_{\rm int}}
\le\nu^{\sigma_{\rm int}-1/6-o(1)}\longrightarrow0.
\end{equation}

Substitution of \eqref{eq:prescribedCoupledEnergy} into
\eqref{eq:approximationCoupledDensity} gives
\begin{align*}
\|\rho[H]\|_{\ell^1_{r-1}L^\infty(I')}
\le{}&C_{\rm s}A_\pm\big[1+\alpha_r(I')
+\mathfrak r(1+{T_1}+\sqrt{|I'|/\nu})\big]Y_{r,I'}\\
&+C_{\rm s}A_\pm\Big[\alpha_r(I')\mathfrak p_{I'}
+{\mathfrak r\sqrt{|I'|/\nu}(1+T_1)}\\
&\hspace{20mm}+\mathfrak r(1+{T_1}+\sqrt{|I'|/\nu})\mathfrak p_{I'}\Big]
\|\rho[H]\|_{\ell^1_{r-1}L^\infty(I')}.
\end{align*}
The background part of the coefficient of the last norm is small by
\eqref{eq:coarseFieldSplit} and \eqref{eq:logFieldSplitSmallness}:
\begin{align*}
A_-\alpha_r(I_-)\mathfrak p_{I_-}&\le C_{\rm s}N^4\nu^{-1/2}\varepsilon_{\rm F}=o(1),\\
A_+\alpha_r(I_+)\mathfrak p_{I_+}
&\le C_{\rm s}e^{C_*N-c_{\rm s}N^{4\chi}}
\nu^{\sigma_{\rm int}-1/2-o(1)}=o(1).
\end{align*}
For the additional terms, $|I'|\le1+T_1$ and the preceding bound on
$\mathfrak p_{I'}$ give
\[
\begin{aligned}
&A_\pm\mathfrak r\big\{\sqrt{|I'|/\nu}(1+T_1)
+(1+T_1+\sqrt{|I'|/\nu})\mathfrak p_{I'}\big\}\\
&\quad\le C_{\rm s}e^{C_*N}[1+\varepsilon_N(1+T_1)]\mathfrak r
\big((1+T_1)^2+\nu^{-1/2}(1+T_1)^{3/2}\big)=o(1).
\end{aligned}
\]
\par

Here $1+\varepsilon_N(1+{T_1})=e^{o(N)}$ is absorbed in the fixed
factor $A_{\rm s}$ in \eqref{eq:unifiedCoefficientSmallness}.
Also $\alpha_r(I')+\mathfrak r(1+{T_1}+\sqrt{|I'|/\nu})=o(1)$ on either interval.
For small $\nu$, the full coefficient of
{$\|\rho[H]\|_{\ell^1_{r-1}L^\infty(I')}$} is at most $1/2$.
Absorbing it and returning to \eqref{eq:prescribedCoupledEnergy} yields
\[
{
\|\rho[H]\|_{\ell^1_{r-1}L^\infty(I')}\le C_{\rm s}A_\pm Y_{r,I'},\qquad \|H/\sqrt \mu\|_{\mathcal N_r^{I'}}\le C_{\rm s}(1+\mathfrak p_{I'}A_\pm)Y_{r,I'}.
}
\]
On the early and late intervals, respectively,
\[
\begin{aligned}
\|\rho[H]\|_{\ell^1_{r-1}L^\infty(I_-)}&\le C_{\rm s}Y_{r,I_-},&
\|H/\sqrt \mu\|_{\mathcal N_r^{I_-}}&\le C_{\rm s}N^4Y_{r,I_-},\\
\|\rho[H]\|_{\ell^1_{r-1}L^\infty(I_+)}&\le C_{\rm s}e^{C_*N}Y_{r,I_+},&
\|H/\sqrt \mu\|_{\mathcal N_r^{I_+}}&\le C_{\rm s}(1+{T_1})
[1+\varepsilon_N(1+{T_1})]e^{C_*N}Y_{r,I_+}.
\end{aligned}
\]
At the common endpoint,
\[
\begin{aligned}
Y_{r,I_+}
&\le\mathcal E_r[H(T_{\rm a})/\sqrt\mu]^{1/2}
+\int_{I_+}\mathcal E_r[S_0/\sqrt\mu]^{1/2}\,ds\\
&\quad+\nu^{-1/2}\sum_{j=0}^r
\|\partial_{x_1}^j(S_1/\sqrt\mu)\|_{L^2_{s,x}(I_+)\Sigma_{2r+6-2j}^*}
\le C_{\rm s}N^4Y_r.
\end{aligned}
\]
Adding the two interval estimates therefore gives
\[
\begin{aligned}
\|\rho[H]\|_{\ell^1_{r-1}L^\infty(I)}
&\le C_{\rm s}N^4e^{C_*N}Y_r,\\
\|H/\sqrt\mu\|_{\mathcal N_r}
&\le C_{\rm s}N^4(1+T_1)[1+\varepsilon_N(1+T_1)]e^{C_*N}Y_r,\\
\|H\|_{\mathscr V_r(I)}
&\le C_{\rm s}N^4[1+\varepsilon_N(1+T_1)]e^{C_*N}Y_r
\le A_{\rm s}(N)e^{C_*N}Y_r.
\end{aligned}
\]
Only the estimate on $I_+$ contains $e^{C_*N}$.
All bounds remain valid on shorter intervals with the original $T_1$.
For ${\bf H}={\bf W}=0$, set $\mathfrak r=0$ throughout; the first two
lines are \eqref{eq:fullForcedDensityBound}--\eqref{eq:fullForcedEnergyBound}
with constants depending only on $r$.

For the continuity equation, velocity cutoff approximation in the
weighted pairing is allowed because
\[
\langle v\rangle^{-4r-12}\sqrt\mu,
\quad v_1\langle v\rangle^{-4r-12}\sqrt\mu
\in\Sigma_{2r+6}.
\]
Local collision conservation and integration by parts give
\[
\int Q(F_1^{\rm c},H)\,dv=\int Q(H,F_2^{\rm c})\,dv=0,
\qquad
\int\partial_{v_1}H\,dv=\int\partial_{v_1}F_2^{\rm c}\,dv=0.
\]
Integrating \eqref{eq:prescribedRemainderExpanded} in velocity thus yields
\[
\partial_s\rho[H]+\partial_{x_1}j[H]
=\int(S_0+S_1)\,dv=0
\]
when the local source integral vanishes.
\par
\end{proof}
\par\endgroup
\begingroup
For the linear completion we use only ${\bf H}={\bf W}=0$. It is enough to show that $S_1=-\mathcal T_{\rm b}H_J$ is small in the dual norm on the right of \eqref{eq:fullForcedDensityBound}, including its factor $\nu^{-1/2}$, so that its size compensates the loss $e^{C_*N}$. This gives the term $\rho[W_J]$ in \eqref{eq:linearConstructionThreeErrors}. The full prescribed-coefficient version is also used for the nonlinear remainder. Proposition~\ref{thm:fullForcedCoarse} alone does not supply the smallness of $H_n$: that smallness is proved by Lemma~\ref{lem:oneFieldOperation} through Lemma~\ref{lem:fieldSpatialInputs}.

Recall the ordered evolution operators introduced after
\eqref{eq:fieldCorrectionRecursion}, with
$\widetilde{\mathcal T}_{\rm b}(s)f=-E[F^{\rm b,\nu}](s)(\partial_{v_1}-v_1/2)f$
as in \eqref{eq:fieldNormalizedOperator}:
\begin{equation}\label{eq:fieldEvolutionCoefficients}
\begin{aligned}
\mathbb U_{\rm b}^{(0)}(s,\tau)&=\mathbb U_{\rm b}^\nu(s,\tau),\\
\mathbb U_{\rm b}^{(q)}(s,\tau)V
&=\int_\tau^s\mathbb U_{\rm b}^\nu(s,\tau_1)
\widetilde{\mathcal T}_{\rm b}(\tau_1)\mathbb U_{\rm b}^{(q-1)}(\tau_1,\tau)V\,d\tau_1,
\qquad q\ge1.
\end{aligned}
\end{equation}
The correction equations \eqref{eq:fieldCorrectionRecursion} give, by Duhamel substitution,
\[
\frac{H_n(s)}{\sqrt\mu}=\mathbb U_{\rm b}^{(n)}(s,t_0)(H_{\rm in}/\sqrt\mu)
+\sum_{q=0}^n\sum_{m\ne0}\int_{t_0}^s
\mathbb U_{\rm b}^{(q)}(s,\tau)\mathcal F_m^{\rm b}(\tau)\rho_m[H_{n-q}](\tau)\,d\tau.
\]
The identity first holds for smooth finite approximations; the estimates below justify its convergence and reconstruction in \eqref{eq:fieldDistributionRecursion}.
\par\endgroup

\subsection{Finite field corrections and completion of the linear echo}\label{sec:fieldExpansion}\label{app:fieldIteration}

\begingroup
We now control the last two errors in
\eqref{eq:linearConstructionThreeErrors}.  We estimate finitely many field corrections with
the shifted-vector-field collision estimates and the chain-weighted density
inverse, and apply the complete estimate only to the remaining source.
\par\endgroup
The phase shifts and the reason for finite corrections were explained in
Section~\ref{sec:roadmapFieldCompletion}. Here we prove, in order,
the single-field bound, the ordered density traces, the density recursion,
and the final source estimate. The fixed order $J$ is chosen by
\eqref{eq:logFieldOrderChoice}; all smaller Gevrey radii, higher moments
and spatial orders are chosen before $\nu\to0$.
\par
\par

\begingroup
\subsubsection{One field operation and ordered propagation}
{We use the shifted-vector-field norm
$\mathfrak G_{c,j}^{\varrho,\lambda,\ell}$ from
\eqref{eq:physicalCenteredNorms}. The following estimate controls both
the derivative cost and the change of spatial mode; the factor
$e^{\varrho_{\rm F}|l(s-t_0)|^\chi}$ in the weighted Fourier estimate for the background field controls
the corresponding change of the shifted vector field.}
\begin{lem}\label{lem:oneFieldOperation}
{There exist fixed admissible radii
$0<\varrho_-<\varrho_+$ and $0<\lambda_-<\lambda_+$ in the range of
Proposition~\ref{prop:exactCollisionPropagation} with the following
property. Let $c\in\mathbb R$, $j\in\mathbb Z$, $s\in I$,
$t=s-t_0$, and let $f$ satisfy
$\|f\|_{\mathfrak G_{c,j}^{\varrho_+,\lambda_+,\ell+1}(s)}<\infty$.
Then}
\begin{equation}\label{eq:oneFieldOperation}
\|\widetilde{\mathcal T}_{\rm b}(s)f\|_{\mathfrak G_{c,j}^{\varrho_-,\lambda_-,\ell}(s)}
\le{}C\mathfrak f(s)\big(1+|c-jt|+t\big) \times
\|f\|_{\mathfrak G_{c,j}^{\varrho_+,\lambda_+,\ell+1}(s)}.
\end{equation}
The constant depends on the fixed differences of radii. In particular, this is a bound between two different regularity spaces.
\end{lem}
\par\endgroup
\begin{proof}
Write $k=k'+l$ and $t=s-t_0$.
By \eqref{eq:physicalSingleModeRelations}, the single-mode norm is
$\mathscr G^{\lambda,\chi}_{\ell;k',c}(t)$ and
$\widehat Z_{k',c}(s)=\partial_{v_1}+i(k't-c)$.
Use $a_n^\lambda=\lambda^n/g_\chi(n)$.
Since
\[
\binom nq\frac{g_\chi(q)g_\chi(n-q)}{g_\chi(n)}
=\binom nq^{1-1/\chi}
\left(\frac{n+1}{(q+1)(n-q+1)}\right)^{12}\le1,
\]
the identity $\widehat Z_{k,c}(s)=\widehat Z_{k',c}(s)+ilt$ and the
binomial formula give
\begin{equation}\label{eq:fieldDerivativeConvolution}
\begin{aligned}
&a_n^{\lambda_-}\|\langle v\rangle^\ell\widehat Z_{k,c}(s)^nw\|_2 \le\sum_{q=0}^n\frac{(\lambda_-|lt|)^q}{g_\chi(q)}
a_{n-q}^{\lambda_-}
\|\langle v\rangle^\ell\widehat Z_{k',c}(s)^{n-q}w\|_2,\\
&\left(\sum_{n\ge0}(a_n^{\lambda_-})^2
\|\langle v\rangle^\ell\widehat Z_{k,c}(s)^nw\|_2^2\right)^{1/2} \le\left(\sum_{q\ge0}\frac{(\lambda_-|lt|)^q}{g_\chi(q)}\right)
\|w\|_{\mathscr G^{\lambda_-,\chi}_{\ell;k',c}(t)}\\
&\quad\le Ce^{C\lambda_-^\chi|lt|^\chi}
\|w\|_{\mathscr G^{\lambda_-,\chi}_{\ell;k',c}(t)}.
\end{aligned}
\end{equation}
For $lt=0$ the series equals $1$. For $lt\ne0$, the summands satisfy
\[
\begin{aligned}
\frac{(\lambda_-|lt|)^{h+1}/g_\chi(h+1)}
{(\lambda_-|lt|)^h/g_\chi(h)}
&=\frac{\lambda_-|lt|}{(h+1)^{1/\chi}}
\left(\frac{h+2}{h+1}\right)^{12}\le\frac12
\quad\big(h\ge C[1+(\lambda_-|lt|)^\chi]\big),\\
\frac{(\lambda_-|lt|)^h}{g_\chi(h)}
&\le(h+1)^{12}
\left(\frac{e(\lambda_-|lt|)^\chi}{h}\right)^{h/\chi}
\le(h+1)^{12}e^{(\lambda_-|lt|)^\chi/\chi}\quad(h\ge1).
\end{aligned}
\]
The finite part and geometric tail consequently obey
\[
\sum_{h\ge0}\frac{(\lambda_-|lt|)^h}{g_\chi(h)}
\le C[1+(\lambda_-|lt|)^\chi]^{13}e^{(\lambda_-|lt|)^\chi/\chi}
\le Ce^{C\lambda_-^\chi|lt|^\chi}.
\]
\par
The derivative and moment costs follow from the explicit coefficients
\begin{align*}
\|\widehat Z_{k',c}(s)f_{k'}\|_{\mathscr G^{\lambda_-,\chi}_{\ell;k',c}(t)}
&\le\left(\sup_{n\ge0}\frac1{\lambda_+}
\left(\frac{\lambda_-}{\lambda_+}\right)^n
\frac{g_\chi(n+1)}{g_\chi(n)}\right)
\|f_{k'}\|_{\mathscr G^{\lambda_+,\chi}_{\ell;k',c}(t)},\\
\widehat Z_{k',c}(s)^n(v_1f_{k'})&=v_1\widehat Z_{k',c}(s)^nf_{k'}+n\widehat Z_{k',c}(s)^{n-1}f_{k'},\\
\|v_1f_{k'}\|_{\mathscr G^{\lambda_-,\chi}_{\ell;k',c}(t)}
&\le\|f_{k'}\|_{\mathscr G^{\lambda_-,\chi}_{\ell+1;k',c}(t)}\\
&\quad+\left(\sup_{n\ge1}\frac{n a_n^{\lambda_-}}{a_{n-1}^{\lambda_-}}\right)
\|f_{k'}\|_{\mathscr G^{\lambda_-,\chi}_{\ell;k',c}(t)}\\
&\le C\|f_{k'}\|_{\mathscr G^{\lambda_+,\chi}_{\ell+1;k',c}(t)}.
\end{align*}
Both suprema are finite. With $D_{1,1}^-=\partial_{v_1}-v_1/2$ from \eqref{eq:gaussianConjugatedDerivatives},
\[
\begin{aligned}
D_{1,1}^-f_{k'}&=\widehat Z_{k',c}(s)f_{k'}-i(k't-c)f_{k'}-\tfrac12v_1f_{k'},\\
\|D_{1,1}^-f_{k'}\|_{\mathscr G^{\lambda_-,\chi}_{\ell;k',c}(t)}
&\le C(1+|k't-c|)\|f_{k'}\|_{\mathscr G^{\lambda_+,\chi}_{\ell+1;k',c}(t)}.
\end{aligned}
\]
Since $\widehat Z_{k,c}=\widehat Z_{k',c}(s)+ilt$, applying
\eqref{eq:fieldDerivativeConvolution} to $w=D_{1,1}^-f_{k'}$ gives
\[
\|E_lD_{1,1}^-f_{k'}\|_{\mathscr G^{\lambda_-,\chi}_{\ell;k,c}(t)}
\le C|E_l|e^{C\lambda_-^\chi|lt|^\chi}
(1+|k't-c|)\|f_{k'}\|_{\mathscr G^{\lambda_+,\chi}_{\ell+1;k',c}(t)}.
\]
The spatial weights and frequency factor satisfy
\begin{align*}
e^{\varrho_-|k'+l-j|^\chi-\varrho_+|k'-j|^\chi}
&\le e^{\varrho_-|l|^\chi-(\varrho_+-\varrho_-)|k'-j|^\chi},\\
(1+|k't-c|)e^{-(\varrho_+-\varrho_-)|k'-j|^\chi}
&\le(1+|jt-c|+t|k'-j|)e^{-(\varrho_+-\varrho_-)|k'-j|^\chi}\\
&\le C_{\varrho_+,\varrho_-}(1+|jt-c|+t).
\end{align*}
Multiplying the single-mode estimate by the spatial weight, summing
over $l$, and applying the discrete Young inequality therefore gives
\[
\begin{aligned}
&\|\widetilde{\mathcal T}_{\rm b}(s)f\|_{\mathfrak G_{c,j}^{\varrho_-,\lambda_-,\ell}(s)} \le C(1+|jt-c|+t)
\left(\sum_l|E_l|e^{\varrho_-|l|^\chi+C\lambda_-^\chi|lt|^\chi}\right)\\
&\qquad\times
\left(\sum_{k'}e^{2\varrho_+|k'-j|^\chi}
\|f_{k'}\|_{\mathscr G^{\lambda_+,\chi}_{\ell+1;k',c}(t)}^2\right)^{1/2} \le C(1+|jt-c|+t)\mathfrak f(s)
\|f\|_{\mathfrak G_{c,j}^{\varrho_+,\lambda_+,\ell+1}(s)}.
\end{aligned}
\]
Here $\varrho_-\le\varrho_{\rm F}$ and $C\lambda_-^\chi\le\varrho_{\rm F}$,
so the displayed $l$-sum is at most $\mathfrak f(s)$; $E_0=0$.
This is \eqref{eq:oneFieldOperation}.
Truncation of the spatial and derivative sums, followed by completion,
extends the inequality {from smooth data to the stated domain.}
\end{proof}

\begingroup
{We now estimate the finite field corrections defined in
\eqref{eq:fieldEvolutionCoefficients}, first on smooth data and then by completion.}

The parenthesized superscript $(q)$ counts background-field insertions; in particular,
$\mathbb U_{\rm b}^{(0)}=\mathbb U_{\rm b}^\nu$ includes collisions and is different
from the free-transport operator $\mathbb U^0$.
Only finitely many $\mathbb U_{\rm b}^{(q)}$ will be used. No convergence of their
infinite sum is asserted.
\par\endgroup
\label{app:orderedFieldTraces}

\begingroup
{The next proposition estimates these field corrections and their
densities using the weighted Fourier estimate for the background field of
Proposition~\ref{prop:backgroundFieldBounds}.}
\par\endgroup
\begin{prop}\label{prop:finiteFieldTrace}
{Fix a velocity weight $\ell\ge4$ and an integer $q\ge1$.
There exist fixed admissible radii $\varrho_q,c_q,\lambda_q>0$ and
$0<\lambda'_q<\lambda_q$ for which the following estimates hold.
Let $t_0\le\tau\le s\le T_1$, $j\in\mathbb Z$, $\eta\in\mathbb R$,
and $V\in\mathscr A_{\ell+q}^{\lambda_q,\chi}$.
Set $c=j(\tau-t_0)+\eta$, and define}
\begin{equation}\label{eq:fieldOrderedIntegral}
\begin{split}
\Xi_q(j,\eta;s,\tau)
={}&\int_{\tau\le\tau_1\le\cdots\le\tau_q\le s}
\prod_{i=1}^q\Big[
\mathfrak f(\tau_i)
\big(1+|\eta-j(\tau_i-\tau)|+\tau_i-t_0\big)
\Big]\,d\tau_1\cdots d\tau_q.
\end{split}
\end{equation}
{Then}
\begin{equation}\label{eq:finiteFieldDistribution}
\begin{split}
&\left(\sum_k e^{2\varrho_q|k-j|^\chi}\sum_{n\ge0}
\frac{(\lambda'_q)^{2n}}{g_\chi(n)^2}
\|\langle v\rangle^\ell\widehat Z_{k,c}(s)^n
[\mathbb U_{\rm b}^{(q)}(s,\tau)\mathbf e_j(e^{i\eta v_1}V)]_k\|_2^2\right)^{1/2}\\
&\qquad\le C_q\Xi_q(j,\eta;s,\tau)
\|V\|_{\mathscr A_{\ell+q}^{\lambda_q,\chi}}.
\end{split}
\end{equation}
In particular, for every $k\in\mathbb Z$,
\begin{equation}\label{eq:finiteFieldTrace}
\begin{split}
\left|\int\sqrt \mu 
\big[\mathbb U_{\rm b}^{(q)}(s,\tau)\mathbf e_j(e^{i\eta v_1}V)\big]_k\,dv\right|
\le{}&C_q e^{-\varrho_q|k-j|^\chi}
e^{-c_q|k(s-t_0)-c|^\chi}\\
&\quad\times\Xi_q(j,\eta;s,\tau)
\|V\|_{\mathscr A_{\ell+q}^{\lambda_q,\chi}}.
\end{split}
\end{equation}
{The constants, radii and required coefficient moment orders may
depend on $\ell,q$ and the fixed background parameters, but are
independent of $\nu,N,s,\tau,j,\eta$.
Here $\mathscr A$ is the velocity derivative norm of
\eqref{eq:velocityInputGevrey}. The bound includes all intermediate
distribution modes, including mode zero.}
\end{prop}

\begingroup
\begin{proof}
Fix the number $q$ of field operations. Choose strictly
decreasing admissible radii $\varrho^{(0)}>\cdots>\varrho^{(q)}>0$ and
$\lambda^{(0)}>\cdots>\lambda^{(q)}>0$ as in
Lemma~\ref{lem:oneFieldOperation}.
The moment orders and radii are fixed before $\nu\to0$. Since
\[
\ell+q-i\ge\ell\ge4,\qquad 0\le i\le q,
\]
Proposition~\ref{prop:exactCollisionPropagation} applies at every level.
Equation~\eqref{eq:physicalCollisionEnergy}, with
$\Theta_{\rm b}({T_1},t_0)=o(1)$ at these finitely many orders, and
Lemma~\ref{lem:oneFieldOperation} give
\begin{equation}\label{eq:fieldLevelBounds}
\begin{aligned}
&\|\mathbb U_{\rm b}^\nu(s,\tau')f\|_{\mathfrak G_{c,j}^{\varrho^{(i)},\lambda^{(i)},\ell+q-i}(s)} \le C_q\|f\|_{\mathfrak G_{c,j}^{\varrho^{(i)},\lambda^{(i)},\ell+q-i}(\tau')},
\qquad 0\le i\le q,\\
&\|\widetilde{\mathcal T}_{\rm b}(\tau')f\|_{\mathfrak G_{c,j}^{\varrho^{(i)},\lambda^{(i)},\ell+q-i}(\tau')}\\
&\quad\le C_q\mathfrak f(\tau')(1+|c-j(\tau'-t_0)|+\tau'-t_0)
\|f\|_{\mathfrak G_{c,j}^{\varrho^{(i-1)},\lambda^{(i-1)},\ell+q-i+1}(\tau')},
\quad 1\le i\le q.
\end{aligned}
\end{equation}
For $V_{\rm in}=\mathbf e_j(e^{i\eta v_1}V)$ and
$c=j(\tau-t_0)+\eta$,
\[
\widehat Z_{j,c}(\tau)^n(e^{i\eta v_1}V)
=e^{i\eta v_1}\partial_{v_1}^nV,
\qquad
\|V_{\rm in}\|_{\mathfrak G_{c,j}^{\varrho^{(0)},\lambda^{(0)},\ell+q}(\tau)}
=\|V\|_{\mathscr A_{\ell+q}^{\lambda^{(0)},\chi}}.
\]
Minkowski's inequality, \eqref{eq:fieldEvolutionCoefficients}, and
\eqref{eq:fieldLevelBounds} yield successively
\begin{equation}\label{eq:fieldOrderedInduction}
\begin{aligned}
&\|\mathbb U_{\rm b}^{(i)}(s,\tau)V_{\rm in}\|_{\mathfrak G_{c,j}^{\varrho^{(i)},\lambda^{(i)},\ell+q-i}(s)}\\
&\quad\le C_q^2\int_\tau^s\mathfrak f(\tau')(1+|c-j(\tau'-t_0)|+\tau'-t_0)\\
&\hspace{38mm}\times\|\mathbb U_{\rm b}^{(i-1)}(\tau',\tau)V_{\rm in}\|_{\mathfrak G_{c,j}^{\varrho^{(i-1)},\lambda^{(i-1)},\ell+q-i+1}(\tau')}\,d\tau'\\
&\quad\le C_q^{2i+1}\|V_{\rm in}\|_{\mathfrak G_{c,j}^{\varrho^{(0)},\lambda^{(0)},\ell+q}(\tau)}\\
&\qquad\times\int_{\tau\le \tau_1\le\cdots\le \tau_i\le s}
\prod_{h=1}^i\big[\mathfrak f(\tau_h)(1+|c-j(\tau_h-t_0)|+\tau_h-t_0)\big]\,d\tau_1\cdots d\tau_i,
\end{aligned}
\end{equation}
where the initial step $i=0$ is the first bound in
\eqref{eq:fieldLevelBounds}. Since
\[
c-j(\tau_h-t_0)=\eta-j(\tau_h-\tau),
\]
the integral at $i=q$ equals $\Xi_q(j,\eta;s,\tau)$.
This proves \eqref{eq:finiteFieldDistribution}, with $\varrho_q=\varrho^{(q)}$, $\lambda'_q=\lambda^{(q)}$,
and $\lambda_q=\lambda^{(0)}$.

For $f=\mathbb U_{\rm b}^{(q)}(s,\tau)V_{\rm in}$, {extract the $k$th mode} and
apply the Gaussian trace estimate \eqref{eq:centeredGtrace} at the
shift parameter:
\[
\begin{aligned}
\left|\int\sqrt \mu  f_k(s)\,dv\right|
&\le C_qe^{-c_q|k(s-t_0)-c|^\chi}
\left(\sum_{n\ge0}\frac{(\lambda^{(q)})^{2n}}{g_\chi(n)^2}
\|\langle v\rangle^\ell\widehat Z_{k,c}^nf_k(s)\|_2^2\right)^{1/2}\\
&\le C_qe^{-\varrho^{(q)}|k-j|^\chi-c_q|k(s-t_0)-c|^\chi}
\|f\|_{\mathfrak G_{c,j}^{\varrho^{(q)},\lambda^{(q)},\ell}(s)}.
\end{aligned}
\]
Together with \eqref{eq:fieldOrderedInduction}, this proves
\eqref{eq:finiteFieldTrace}. For smooth $V_h\to V$ in
$\mathscr A_{\ell+q}^{\lambda_q,\chi}$, the same induction applied to
the difference gives
\[
\|\mathbb U_{\rm b}^{(q)}(s,\tau)
\mathbf e_j(e^{i\eta v_1}(V_h-V_{h'}))\|_{\mathfrak G_{c,j}^{\varrho_q,\lambda'_q,\ell}(s)} \le C_q\Xi_q(j,\eta;s,\tau)
\|V_h-V_{h'}\|_{\mathscr A_{\ell+q}^{\lambda_q,\chi}}\longrightarrow0.
\]
The integrable majorant in \eqref{eq:fieldOrderedInduction} passes each
ordered integral to this limit, proving the assertions on the stated domain.
\par
\end{proof}
\par\endgroup

\begingroup
\subsubsection{Density operators and their spatial derivative loss}
When acting on $\mathcal F_m^{\rm b}(\tau)$, each field operation contributes a factor involving $|m|$.
Poisson's equation supplies $|m|^{-1}$ and integration with respect to $\tau$
supplies a second $|m|^{-1}$. The following estimates record the
resulting loss $\max\{q-2,0\}$ of spatial derivatives.
\par\endgroup
Recall, for $q\ge1$, the functions and density operators introduced
after \eqref{eq:fieldCorrectionRecursion}:
\begin{equation}\label{eq:fieldDensityCoefficients}
\begin{aligned}
\mathfrak g_{q,k}(s)
&=\int\sqrt \mu [\mathbb U_{\rm b}^{(q)}(s,t_0)(H_{\rm in}/\sqrt \mu )]_k\,dv,\\
\mathcal K_{q;k,m}(s,\tau)
&=\int\sqrt \mu [\mathbb U_{\rm b}^{(q)}(s,\tau)\mathcal F_m^{\rm b}(\tau)]_k\,dv,\\
(\mathbb T_q{\rho})_k(s)
&=\sum_{m\ne0}\int_{t_0}^s
\mathcal K_{q;k,m}(s,\tau){\rho}_m(\tau)\,d\tau.
\end{aligned}
\end{equation}
\begingroup
The coefficient $\mathcal F_m^{\rm b}$ is defined in
\eqref{eq:exactCollisionElectricInput}. Set $\mathfrak g_0=\mathfrak g^{\rm C,\nu}$.
Only $q\ge1$ is used for the density operators $\mathbb T_q$:
the operator with no field insertion is already the sum
$\mathbb K+\mathbb P_\downarrow+\mathbb P_\uparrow+\Delta\mathbb T$
in \eqref{eq:exactCoefficientDensityEquation}.
The indices $n,q$ count correction orders, whereas $k,m$ denote spatial
Fourier modes; $\rho_{k}^{(n)}$ is the $k$th mode of the $n$th correction.
Use $[y]_+=\max\{y,0\}$ as in \eqref{eq:Bplus}, and the common density norm
\eqref{eq:echoDensityNorm}.
\par\endgroup

\begingroup
\begin{prop}\label{thm:finiteFieldResponse}
{Fix $q\ge1$ and $r_{\rm x}\ge0$, and take the density-radius constant
$\vartheta>0$ sufficiently small. Let
${\rho}\in\mathscr X_N^{r_{\rm x}+[q-2]_+}$. Then}
\begin{align}
\|\mathfrak g_q\|_{\mathscr X_N^{r_{\rm x}}}
&\le C_{{r_{\rm x}},q}(\zeta\varepsilon_{\rm F})^q,
\label{eq:fieldInitialCoefficient}\\
\|\mathbb T_q{\rho}\|_{\mathscr X_N^{r_{\rm x}}}
&\le C_{{r_{\rm x}},q}N^{2q+[q-2]_+}\varepsilon_{\rm F}^{q}
\|{\rho}\|_{\mathscr X_N^{{r_{\rm x}}+[q-2]_+}}.
\label{eq:fieldDensityCoefficient}
\end{align}
For $q=1,2$, the domain and range in \eqref{eq:fieldDensityCoefficient} have the same spatial polynomial weight. In particular,
\begin{equation}\label{eq:firstFieldInvertedOperator}
\|\mathbb V^{\rm C,\nu}\mathbb T_1\|_{\mathscr X_N^{r_{\rm x}}\to\mathscr X_N^{r_{\rm x}}}
\le C_{r_{\rm x}}N^3\varepsilon_{\rm F}.
\end{equation}
This is at most $\nu^{1/3+{\sigma_{\rm int}}-o(1)}$.

No restriction on the spatial support of ${\rho}$ is imposed.
\end{prop}
\par\endgroup

\begingroup
\begin{proof}
Use Proposition~\ref{prop:finiteFieldTrace} with a fixed velocity weight
$\ell\ge4$; both $H_{\rm in}/\sqrt\mu$ and $\mathcal F_m^{\rm b}$ have the
required higher moments.
For $\mathfrak g_q$, the spatial indices are $j=\pm m_0$ and the velocity frequencies are $\eta=\pm(\zeta-m_0t_0)$. Since $m_0\mathrel{{\sim}} N\ll\zeta$, every factor in \eqref{eq:fieldOrderedIntegral} is bounded by
\[
C\zeta\langle\tau_i-t_0\rangle\mathfrak f(\tau_i).
\]
The symmetric product integrand and \eqref{eq:fieldFourierEnvelope} give
\[
\Xi_q\le\frac{(C\zeta)^q}{q!}
\left(\int_{t_0}^s\langle u-t_0\rangle\mathfrak f(u)\,du\right)^q
\le C_q(\zeta\varepsilon_{\rm F})^q.
\]
For the two coefficients $(H_{\rm in}/\sqrt\mu)_{\pm m_0}$,
\[
\big||k|s-\zeta\big|
\le|k(s-t_0)\mp(\zeta-m_0t_0)|+|k\mp m_0|t_0,
\qquad \lambda_N(s)t_0^\chi\le C\vartheta N^{-\chi}.
\]
Choose $\vartheta$ small enough to leave half of the phase and spatial
decays in \eqref{eq:finiteFieldTrace}. Since $B_{|k|}^{-1}\le1$,
\[
\|\mathfrak g_q\|_{\mathscr X_N^{r_{\rm x}}}
\le C_q(\zeta\varepsilon_{\rm F})^q
\sum_\pm\sum_{k\ne0}\langle k/N\rangle^{r_{\rm x}}
e^{-\varrho_q|k\mp m_0|^\chi/2} \le C_{r_{\rm x},q}(\zeta\varepsilon_{\rm F})^q,
\]
because $m_0/N\le C$ and the norms of $\mu_{1/2}/\sqrt{4\pi\mu}$ are fixed.
This proves \eqref{eq:fieldInitialCoefficient}.
\par
For the electric kernel, use the coefficient $\mathcal F_m^{\rm b}(\tau)$ in \eqref{eq:completeElectricProfile}. {The component in mode $j=m+l$ has} phase $\eta=-l(\tau-t_0)$ and shift parameter
\[
c=m(\tau-t_0),\qquad
c-j(\tau_i-t_0)=m(\tau-\tau_i)-l(\tau_i-t_0).
\]
The Maxwellian term corresponds to $l=0$. For every field time $\tau_i\ge\tau$,
\[
1+|c-j(\tau_i-t_0)|+\tau_i-t_0
\le C(1+|m|)(1+|l|)\langle\tau_i\rangle.
\]
By \eqref{eq:completeElectricInputBound}, the non-Maxwellian part of $[\mathcal F_m^{\rm b}(\tau)]_{m+l}$ has norm at most
$C\varepsilon_N(1+|l|\tau)e^{-\varrho_{\rm b}|l|^\chi}/|m|$.
For $q\ge1$,
\begin{equation}\label{eq:orderedFieldTimeIntegral}
{
\int_{\tau\le t_1\le\cdots\le t_q\le s} \prod_{i=1}^q\langle t_i\rangle^2\mathfrak f(t_i)\,dt_1\cdots dt_q =\frac1{q!}\left(\int_\tau^s\langle t\rangle^2\mathfrak f(t)\,dt\right)^q \le C_qN^{2q}\varepsilon_{\rm F}^q e^{-c_{{\rm dec},q}(\tau-t_0)^\chi}.
}
\end{equation}
The equality follows by partitioning the product cube {into the $q!$}
{ordered regions; the diagonals have measure zero.} The last bound is
\eqref{eq:fieldFourierMoments} with the time weight $\langle s\rangle^2$ and $t_0\mathrel{{\sim}} N$.
For every ordered tuple,
\[
1+|l|\tau\le C\langle l\rangle\langle\tau\rangle,
\qquad \langle\tau\rangle\prod_{i=1}^q\langle t_i\rangle
\le\prod_{i=1}^q\langle t_i\rangle^2.
\]
Write $[\mathcal F_m^{\rm b}(\tau)]_{m+l}$ as
$e^{-il(\tau-t_0)v_1}V_{m,l}(\tau,v)$, {combining the two summands} if $l=0$.
{The corresponding estimate for $V_{m,l}(\tau)$ is}
\begin{equation}\label{eq:electricOrderedInputBound}
\begin{aligned}
&\Xi_q(m+l,-l(\tau-t_0);s,\tau)
\|V_{m,l}(\tau)\|_{\mathscr A_{\ell+q}^{\lambda_q,\chi}}\\
&\quad\le\frac{C_q(1+|m|)^q}{|m|}
\left(\mathbf1_{l=0}+\varepsilon_N\langle l\rangle^{q+1}
e^{-\varrho_{\rm b}|l|^\chi}\right) C_qN^{2q}\varepsilon_{\rm F}^q e^{-c_{{\rm dec},q}(\tau-t_0)^\chi}.
\end{aligned}
\end{equation}
All fixed powers of $\langle l\rangle$ are absorbed {by decreasing the}
fixed Gevrey radius. Using the physical phase,
Proposition~\ref{prop:finiteFieldTrace} yields
\begin{equation}\label{eq:fieldCompleteKernel}
\begin{split}
|\mathcal K_{q;k,m}(s,\tau)|
\le{}&\frac{C_qN^{2q}\varepsilon_{\rm F}^{q}(1+|m|)^q}{|m|}
e^{-c_q|k(s-t_0)-m(\tau-t_0)|^\chi}\\
&\quad\times\left(e^{-\varrho_q|k-m|^\chi}
+\varepsilon_N\sum_l e^{-\varrho_q'|l|^\chi}
e^{-\varrho_q|k-m-l|^\chi}\right).
\end{split}
\end{equation}
This argument uses only $\varepsilon_N=o(1)$ and the full background bound. It does not require a fixed-radius relative error to tend to zero. In particular, no factor ${T_1}$ has been used to integrate a field time.

We now insert the density weights. Set $n=k-m-l$, with $l=0$ for the Maxwellian term. Since $ks-m\tau=\big(k(s-t_0)-m(\tau-t_0)\big)+(n+l)t_0$, the decreasing time radius gives
\[
\frac{{w_{|k|,N}}(s)}{{w_{|m|,N}}(\tau)}
\le e^{\lambda_N(s)|k(s-t_0)-m(\tau-t_0)|^\chi+
\lambda_N(s)t_0^\chi(|n|^\chi+|l|^\chi)}.
\]
Choose $\vartheta$ so that
\[
\lambda_N(s)\le c_q/2,\qquad
\lambda_N(s)t_0^\chi\le\tfrac12\min\{\varrho_q,\varrho_q'\}.
\]
Then, with $n=k-m-l$,
\[
\begin{aligned}
&\frac{w_{|k|,N}(s)}{w_{|m|,N}(\tau)}
e^{-c_q|k(s-t_0)-m(\tau-t_0)|^\chi-
\varrho_q'|l|^\chi-\varrho_q|n|^\chi}\\
&\qquad\le e^{-\frac12c_q|k(s-t_0)-m(\tau-t_0)|^\chi-
\frac12\varrho_q'|l|^\chi-\frac12\varrho_q|n|^\chi}.
\end{aligned}
\]
\par
The only unfavorable gain ratio occurs when $s\ge\tau\ge T_{\rm a}$ and $|k|>|m|$. In that case
\[
|k(s-t_0)-m(\tau-t_0)|\ge |k|(s-t_0)-|m|(\tau-t_0)
\ge (|k|-|m|)(s-t_0)\ge N^4.
\]
Consequently,
\[
\frac{B_{|m|}}{B_{|k|}}
e^{-c_q|k(s-t_0)-m(\tau-t_0)|^\chi/4}
\le e^{CN-c_qN^{4\chi}/4}\le1.
\]
All other gain ratios are at most one; positive fixed fractions
of the phase and spatial decays remain for integration and summation.
\par

For fixed $s,k,m$, the remaining phase integrates to $C_q/|m|$. For the spatial polynomial ratio use
\[
\frac{\langle k/N\rangle^{r_{\rm x}}}{\langle m/N\rangle^{{r_{\rm x}}+[q-2]_+}}
\le C_{r_{\rm x}}\langle l\rangle^{r_{\rm x}}\langle n\rangle^{r_{\rm x}}
\langle m/N\rangle^{-[q-2]_+}.
\]
The first two factors are summable against the remaining spatial decays, uniformly in $N$. Thus each weighted column indexed by $m$ is bounded by
\[
C_{{r_{\rm x}},q}N^{2q}\varepsilon_{\rm F}^q
\frac{(1+|m|)^q}{|m|^2}
\langle m/N\rangle^{-[q-2]_+}(1+\varepsilon_N)
\le C_{{r_{\rm x}},q}N^{2q+[q-2]_+}\varepsilon_{\rm F}^q.
\]
Here
\[
\sup_{m\ne0}\frac{(1+|m|)^q}{|m|^2}
\langle m/N\rangle^{-[q-2]_+}\le C_qN^{[q-2]_+}.
\]
Taking the supremum over $s$ for each $k$ before summing the columns gives
\[
\begin{aligned}
\|\mathbb T_q\rho\|_{\mathscr X_N^{r_{\rm x}}}
&\le C_{r_{\rm x},q}N^{2q+[q-2]_+}\varepsilon_{\rm F}^q
\|\rho\|_{\mathscr X_N^{r_{\rm x}+[q-2]_+}},\\
\|\mathbb V^{\rm C,\nu}\mathbb T_1\rho\|_{\mathscr X_N^{r_{\rm x}}}
&\le C_{r_{\rm x}}N\,\|\mathbb T_1\rho\|_{\mathscr X_N^{r_{\rm x}}}
\le C_{r_{\rm x}}N^3\varepsilon_{\rm F}\|\rho\|_{\mathscr X_N^{r_{\rm x}}}.
\end{aligned}
\]
The absolutely summable columns extend these bounds from finite
sequences and preserve modewise time continuity.
\par
\end{proof}
\par\endgroup
\subsubsection{Spatial norms and reconstruction of the iterates}\label{app:physicalFieldDerivatives}
The density recursion determines $\rho^{(n)}$ in the electric source $\sum_{m\ne0}\mathcal F_m^{\rm b}\rho_m^{(n)}$ of the collision equation.
We first bound the reconstructed distributions in spatial Sobolev norms;
the energy construction then supplies their integrated dissipation.

\begin{lem}\label{lem:fieldSpatialInputs}
Under the hypotheses of Proposition~\ref{thm:finiteFieldResponse},
fix an integer $r\ge0$, a moment $\ell\ge0$, and a finite integer $q\ge1$.
Choose the moments and radii required by Proposition~\ref{prop:finiteFieldTrace} for the estimate at velocity weight $\max\{\ell,4\}$.
Then, for $m\ne0$ and $t_0\le\tau\le s\le T_1$,
\begin{equation}\label{eq:quantitativeElectricDistribution}
\|\mathbb U_{\rm b}^{(q)}(s,\tau)\mathcal F_m^{\rm b}(\tau)\|_{H_x^rL_\ell^2}
\le C_{r,\ell,q}N^{2q+2}\varepsilon_{\rm F}^q\langle m\rangle^{r+q+1}.
\end{equation}
For $q=0$ the bound is $C_{r,\ell}N^2L_N\langle m\rangle^{r+1}$,
where $\log L_N=o(N)$.
\end{lem}
\begin{proof}
The zeroth velocity-derivative term of the centered norm gives
\[
\|f\|_{H_x^rL_\ell^2}
\le C_{r,\varrho}\langle j\rangle^r
\|f\|_{\mathfrak G_{c,j}^{\varrho,\lambda,\ell}(s)},
\qquad
\sup_k\langle k-j\rangle^r e^{-\varrho|k-j|^\chi}<\infty.
\]
Use the electric profiles $V_{m,l}$ of
\eqref{eq:electricOrderedInputBound}, with $j=m+l$ and $c=m(\tau-t_0)$.
Equations \eqref{eq:finiteFieldDistribution} and
\eqref{eq:electricOrderedInputBound} yield
\[
\begin{aligned}
\|\mathbb U_{\rm b}^{(q)}\mathcal F_m^{\rm b}\|_{H_x^rL_\ell^2}
&\le C_{r,\ell,q}\sum_l\langle m+l\rangle^r
\Xi_q(m+l,-l(\tau-t_0);s,\tau)
\|V_{m,l}(\tau)\|_{\mathscr A_{\ell+q}^{\lambda_q,\chi}}\\
&\le \frac{C_{r,\ell,q}N^{2q}\varepsilon_{\rm F}^q
\langle m\rangle^{r+q+1}}{|m|}
e^{-c_{{\rm dec},q}(\tau-t_0)^\chi}
\sum_l\langle l\rangle^{r+q+2}
\big(\mathbf1_{l=0}+\varepsilon_Ne^{-\varrho_{\rm b}|l|^\chi}\big).
\end{aligned}
\]
The sum is uniformly bounded. For $q=0$, use
\eqref{eq:physicalCollisionEnergy} and
\eqref{eq:completeElectricInputBound}; the same spatial sum gives
\[
\|\mathbb U_{\rm b}^{(0)}(s,\tau)\mathcal F_m^{\rm b}(\tau)\|_{H_x^rL_\ell^2}
\le C_{r,\ell}\langle m\rangle^{r+1}
[1+\varepsilon_N(1+T_1)]
\le C_{r,\ell}N^2L_N\langle m\rangle^{r+1}.
\]
For $\ell<4$, the estimate at moment $4$ implies the stated bound.
\end{proof}

We now reconstruct the distributions and estimate the last field
source by \eqref{eq:electricL2SourceBound}.
\par

In the following estimates $A_{r,\ell,J}(N)$ denotes a function of the form
\[
A_{r,\ell,J}(N)=C_{r,\ell,J}N^{C_{r,\ell,J}}L_N.
\]
Here $L_N=C N^P e^{CN^\chi}$ uses only the finitely many coefficient orders needed for the estimate. These constants may depend on $r,\ell,J$; for every fixed choice $\log A_{r,\ell,J}(N)=o(N)$.

\begingroup
\begin{prop}\label{prop:fieldQuantitativeInputs}
Under the hypotheses of Proposition~\ref{thm:finiteFieldExpansion},
fix a polynomial density weight $r_{\rm x}\ge0$.
Choose the radius constant in \eqref{eq:echoDensityNorm} sufficiently small
depending on $J$. The coefficients constructed there satisfy
\begin{equation}\label{eq:fieldMainCoefficients}
\|\rho[H_n]\|_{\mathscr X_N^{r_{\rm x}}}
+\frac{\sqrt N}{B_*}|\rho_{m_1}[H_n]({T_1})|
\le\nu^{n{\sigma_{\rm int}}-o(1)},\qquad 1\le n\le J.
\end{equation}
For each fixed integer $r\ge0$ and $\ell\ge0$, the same coefficients satisfy, for $0\le n\le J$,
\begin{equation}\label{eq:quantitativeFieldDistribution}
\sup_{s\in I}\|H_n(s)/\sqrt\mu\|_{H_x^rL_\ell^2}
\le A_{r,\ell,J}(N)\widetilde B(1+T_1)(\zeta\varepsilon_{\rm F})^n.
\end{equation}
For each fixed integer $r\ge2$, the final source $S_J=-\mathcal T_{\rm b}H_J$ satisfies
\begin{equation}\label{eq:quantitativeFieldResidual}
\nu^{-1/2}\sum_{j=0}^{r}\left\|\partial_{x_1}^j\!\left(\mu^{-1/2}\bigl(S_J\bigr)\right)\right\|_{L^2_{s,x}(I)\Sigma_{2r+6-2j}^*}
\le A_{r,J}(N)\widetilde B\frac{1+T_1}{\sqrt\nu\,\zeta}
(\zeta\varepsilon_{\rm F})^{J+1}.
\end{equation}

\par
Consequently,
\begin{equation}\label{eq:logQuantitativeFieldResidual}
\nu^{-1/2}\sum_{j=0}^{r}\left\|\partial_{x_1}^j\!\left(\mu^{-1/2}\bigl(S_J\bigr)\right)\right\|_{L^2_{s,x}(I)\Sigma_{2r+6-2j}^*}
\le\nu^{(J+1){\sigma_{\rm int}}-1/2-KC_{\rm B}-o(1)}.
\end{equation}
The principal loss $1/2+KC_{\rm B}$ is independent of $J$.
Each additional correction is generated by one further field operation
and contributes another factor $\zeta\varepsilon_{\rm F}
\le\nu^{\sigma_{\rm int}-o(1)}$ in
\eqref{eq:quantitativeFieldDistribution}--\eqref{eq:quantitativeFieldResidual}.
All these estimates concern fixed finite $J$; no monotonicity of the
actual errors or convergence as $J\to\infty$ is asserted.
\end{prop}
\par
\begin{proof}[Proofs of Propositions~\ref{thm:finiteFieldExpansion} and~\ref{prop:fieldQuantitativeInputs}]
\emph{1. Density equations and spatial regularity requirements.} Start with $H_0=H^{\rm C,\nu}$ and $\rho^{(0)}=\rho^{\rm C,\nu}$. For $1\le n\le J$, define the density successively by
\begin{equation}\label{eq:fieldDensityRecursion}
\rho^{(n)}=\mathbb V^{\rm C,\nu}
\left(\mathfrak g_n+\sum_{q=1}^n\mathbb T_q\rho^{(n-q)}\right).
\end{equation}
Fix a density regularity order ${r_{\rm x}}$ and choose the spatial regularity orders
\begin{equation}\label{eq:fieldRegularityOrders}
\beta_n={r_{\rm x}}+J-n,\qquad 0\le n\le J.
\end{equation}
At $n=0$, \eqref{eq:logExactCollisionInverse} gives
$\rho^{(0)}\in\mathscr X_N^{\beta_0}$. For $1\le q\le n$,
\[
[q-2]_+=\max\{q-2,0\}\le q,
\qquad \beta_n+[q-2]_+\le\beta_n+q=\beta_{n-q}.
\]
Proposition~\ref{thm:finiteFieldResponse} and the monotonicity of the
polynomial weights therefore give the maps
\[
\mathbb T_q:\mathscr X_N^{\beta_{n-q}}
\longrightarrow\mathscr X_N^{\beta_n},
\qquad
\mathbb V^{\rm C,\nu}:\mathscr X_N^{\beta_n}
\longrightarrow\mathscr X_N^{\beta_n}.
\]
Induction in \eqref{eq:fieldDensityRecursion} constructs
$\rho^{(n)}\in\mathscr X_N^{\beta_n}$ through $\beta_J={r_{\rm x}}$.
Increasing ${r_{\rm x}}$ supplies any further fixed spatial orders required below;
uniqueness of the same density recursion identifies these constructions.

Choose all moments and radii for the finitely many estimates for the initial datum and electric coefficients
in Lemma~\ref{lem:fieldSpatialInputs} before $\nu\to0$.
The Gaussian datum has every required fixed moment and derivative.

\emph{2. Energy construction and Duhamel reconstruction.}
For any desired energy order $r\ge2$, carry out the construction at
orders $r+J-n$, $0\le n\le J$, increasing the density orders above if necessary.
At step $n\ge1$, $\rho^{(n)}$ and $H_{n-1}$ are already known.
Solve the equation with prescribed $\rho^{(n)}$ and $H_{n-1}$
\begin{equation}\label{eq:fieldPrescribedReconstruction}
\begin{aligned}
\partial_sH_n+v_1\partial_{x_1}H_n
-\nu\big(Q(F^{\rm b,\nu},H_n)+Q(H_n,F^{\rm b,\nu})\big)
&=-E[\rho^{(n)}]\partial_{v_1}F^{\rm b,\nu}
-E[F^{\rm b,\nu}]\partial_{v_1}H_{n-1},\\
H_n(t_0)&=0.
\end{aligned}
\end{equation}
The background coefficient bounds and \eqref{eq:electricL2SourceBound} give
\[
\begin{aligned}
&\int_{I}\mathcal E_{r+J-n}[\mu^{-1/2}\bigl(E[\rho^{(n)}]\partial_{v_1}F^{\rm b,\nu}\bigr)]^{1/2}\,ds \le C_{\nu,N,r,J}\|\rho^{(n)}\|_{\ell^1_{r+J-n}L^\infty}<\infty,\\
&\nu^{-1/2}\sum_{j=0}^{r+J-n}\left\|\partial_{x_1}^j\!\left(\mu^{-1/2}\bigl(E[F^{\rm b,\nu}]\partial_{v_1}H_{n-1}\bigr)\right)\right\|_{L^2_{s,x}(I)\Sigma_{2(r+J-n)+6-2j}^*}\\
&\quad\le C_{r,J}\nu^{-1/2}
\|E[F^{\rm b,\nu}]\|_{L^2_sC_x^{r+J-n}}
\|H_{n-1}/\sqrt\mu\|_{\mathcal N_{r+J-n+1}}<\infty.
\end{aligned}
\]
The second bound uses exactly the energy order of the preceding
iterate. Starting from $H_0$ at order $r+J$,
Lemma~\ref{lem:collisionSourceEstimates} constructs successively
\[
H_n/\sqrt\mu\in C(I;\mathscr E_{r+J-n}),\qquad
\|H_n/\sqrt\mu\|_{\mathcal N_{r+J-n}}<\infty,
\qquad 0\le n\le J.
\]
All constants in this construction may depend on the fixed $\nu$;
Step 5 supplies the quantitative bounds.

To justify the ordered Duhamel formula, first increase the density
orders so that $\rho^{(n)}\in\mathscr X_N^{r+J+2}$ for $0\le n\le J$.
By \eqref{eq:echoDensityEmbeddings} and modewise time continuity,
spatial truncation followed by time smoothing after continuous
extension at the endpoints gives smooth functions $\rho_\epsilon^{(n)}$ whose Fourier transforms in $x$ have finite support such that
\par
\[
\sum_{n=0}^J
\|\rho_\epsilon^{(n)}-\rho^{(n)}\|_{\ell^1_{r+J+2}L^\infty(I)}
\longrightarrow0\qquad(\epsilon\downarrow0).
\]
Let $H_{n,\epsilon}$ solve \eqref{eq:fieldPrescribedReconstruction}
with the prescribed functions $\rho_\epsilon^{(n)},H_{n-1,\epsilon}$.
For $n=0$ omit the second source and impose
$H_{0,\epsilon}(t_0)=H_{\rm in}$; for $n\ge1$ the initial datum is zero.
{The continuous dependence in
Lemma~\ref{lem:collisionSourceEstimates}, with the two source bounds
above, gives}
\[
\|(H_{0,\epsilon}-H_0)/\sqrt\mu\|_{\mathcal N_{r+J}}
\le C_{\nu,N,r,J}
\|\rho_\epsilon^{(0)}-\rho^{(0)}\|_{\ell^1_{r+J}L^\infty}
\longrightarrow0.
\]
For $1\le n\le J$ it gives
\begin{equation}\label{eq:fieldReconstructionEnergyLimit}
\begin{aligned}
&\|(H_{n,\epsilon}-H_n)/\sqrt\mu\|_{\mathcal N_{r+J-n}} \le C_{\nu,N,r,J}
\|\rho_\epsilon^{(n)}-\rho^{(n)}\|_{\ell^1_{r+J-n}L^\infty}\\
&\qquad+C_{r,J}\nu^{-1/2}
\|E[F^{\rm b,\nu}]\|_{L^2_sC_x^{r+J-n}}
\|(H_{n-1,\epsilon}-H_{n-1})/\sqrt\mu\|_{\mathcal N_{r+J-n+1}}
\longrightarrow0.
\end{aligned}
\end{equation}
The last limit follows by induction. In particular,
\[
\begin{aligned}
&\nu^{-1/2}\sum_{j=0}^{r+J-n}\left\|\partial_{x_1}^j\!\left(\mu^{-1/2}\bigl(E[F^{\rm b,\nu}]\partial_{v_1}(H_{n-1,\epsilon}-H_{n-1})\bigr)\right)\right\|_{L^2_{s,x}(I)\Sigma_{2(r+J-n)+6-2j}^*}\\
&\quad\le C_{r,J}\nu^{-1/2}
\|E[F^{\rm b,\nu}]\|_{L^2_sC_x^{r+J-n}}
\|(H_{n-1,\epsilon}-H_{n-1})/\sqrt\mu\|_{\mathcal N_{r+J-n+1}}
\longrightarrow0.
\end{aligned}
\]

For the smooth densities $\rho_\epsilon^{(n)}$ whose Fourier transforms in $x$ have finite support, Duhamel's formula reads
\[
\begin{aligned}
\frac{H_{n,\epsilon}(s)}{\sqrt\mu}
={}&\sum_{m\ne0}\int_{t_0}^s
\mathbb U_{\rm b}^\nu(s,\tau)\mathcal F_m^{\rm b}(\tau)
\rho_{\epsilon,m}^{(n)}(\tau)\,d\tau\\
&+\int_{t_0}^s\mathbb U_{\rm b}^\nu(s,\tau)
\widetilde{\mathcal T}_{\rm b}(\tau)\frac{H_{n-1,\epsilon}(\tau)}{\sqrt\mu}\,d\tau,
\qquad n\ge1.
\end{aligned}
\]
At $n=0$ replace the last integral by
$\mathbb U_{\rm b}^\nu(s,t_0)(H_{\rm in}/\sqrt\mu)$.
Substitution of the preceding iterate uses the defining identity
\[
\int_\tau^s\mathbb U_{\rm b}^\nu(s,u)\widetilde{\mathcal T}_{\rm b}(u)
\mathbb U_{\rm b}^{(q)}(u,\tau)\,du
=\mathbb U_{\rm b}^{(q+1)}(s,\tau),\qquad q\ge0.
\]
For $0\le q<n\le J$, the ordered bounds
\eqref{eq:fieldLevelBounds}--\eqref{eq:fieldOrderedInduction}, applied
to the electric profiles in \eqref{eq:electricOrderedInputBound}, give
\[
\begin{aligned}
&\sum_{m\ne0}\int_{t_0}^s\int_\tau^s
\|\mathbb U_{\rm b}^\nu(s,u)\widetilde{\mathcal T}_{\rm b}(u)
\mathbb U_{\rm b}^{(q)}(u,\tau)\mathcal F_m^{\rm b}(\tau)\|_{H_x^rL_\ell^2}
\,du\,|\rho_{\epsilon,m}^{(n-1-q)}(\tau)|\,d\tau\\
&\quad\le C_{\nu,N,r,\ell,J}\sum_{m\ne0}
\langle m\rangle^{r+q+2}\sup_I|\rho_{\epsilon,m}^{(n-1-q)}|<\infty.
\end{aligned}
\]
The datum term has the finite ordered bound
\eqref{eq:fieldOrderedInduction}. Fubini and induction therefore yield
\par
\[
\frac{H_{n,\epsilon}(s)}{\sqrt\mu}
=\mathbb U_{\rm b}^{(n)}(s,t_0)(H_{\rm in}/\sqrt\mu)
+\sum_{q=0}^n\sum_{m\ne0}\int_{t_0}^s
\mathbb U_{\rm b}^{(q)}(s,\tau)\mathcal F_m^{\rm b}(\tau)
\rho_{\epsilon,m}^{(n-q)}(\tau)\,d\tau.
\]
For every fixed $\ell\ge0$, Lemma~\ref{lem:fieldSpatialInputs} bounds
the difference of the electric sums by
\[
\begin{aligned}
&\sup_{s\in I}\left\|
\sum_{q=0}^n\sum_{m\ne0}\int_{t_0}^s
\mathbb U_{\rm b}^{(q)}(s,\tau)\mathcal F_m^{\rm b}(\tau)
\big(\rho_{\epsilon,m}^{(n-q)}-\rho_m^{(n-q)}\big)(\tau)\,d\tau
\right\|_{H_x^rL_\ell^2}\\
&\quad\le C_{\nu,N,r,\ell,J}
\sum_{q=0}^n\|\rho_\epsilon^{(n-q)}-\rho^{(n-q)}
\|_{\ell^1_{r+q+1}L^\infty(I)}\longrightarrow0.
\end{aligned}
\]
{The $C(I;H_x^rL_\ell^2)$ limit therefore equals
$H_n/\sqrt\mu$ by \eqref{eq:fieldReconstructionEnergyLimit}, proving}
\begin{equation}\label{eq:fieldDistributionRecursion}
\frac{H_n(s)}{\sqrt\mu}
=\mathbb U_{\rm b}^{(n)}(s,t_0)(H_{\rm in}/\sqrt\mu)
+\sum_{q=0}^n\sum_{m\ne0}\int_{t_0}^s
\mathbb U_{\rm b}^{(q)}(s,\tau)\mathcal F_m^{\rm b}(\tau)
\rho_m^{(n-q)}(\tau)\,d\tau.
\end{equation}
{The field-source convergence displayed after
\eqref{eq:fieldReconstructionEnergyLimit} uses
\eqref{eq:electricL2SourceBound} and no additional estimate
of an ordinary velocity derivative of $H_n$.}
\par

\emph{3. Density consistency and the correction equations.}
Taking the bounded velocity density pairing in the convergent reconstruction gives
\[
\rho[H_n]=\mathfrak g_n+\sum_{q=1}^n\mathbb T_q\rho^{(n-q)}
+(\mathbb K+\mathbb P_\downarrow+\mathbb P_\uparrow+\Delta\mathbb T)\rho^{(n)}
=\rho^{(n)},
\]
where the last equality is \eqref{eq:fieldDensityRecursion}.
For the zero density mode,
\[
\int E[\rho^{(n)}]\partial_{v_1}F^{\rm b,\nu}\,dv
=\int E[F^{\rm b,\nu}]\partial_{v_1}H_{n-1}\,dv=0,
\]
so local collision mass conservation gives
\[
\partial_s\rho_0[H_n]=0,\qquad \rho_0[H_n](t_0)=0,
\qquad \rho_0[H_n](s)=0.
\]
Homogeneous distribution modes remain included in the reconstruction.
\par
Substituting the density identity into
\eqref{eq:fieldPrescribedReconstruction} gives
\[
\mathcal L_{\rm C}H_n=-\mathcal T_{\rm b}H_{n-1},
\qquad H_n(t_0)=0.
\]
For two solutions with the same preceding iterate and initial datum,
subtraction gives
\[
\mathcal L_{\rm C}(H_n-\widetilde H_n)=0,\qquad
(H_n-\widetilde H_n)(t_0)=0.
\]
The finite-energy uniqueness proved in Step 3 of the proof of
Proposition~\ref{thm:exactCollisionEcho} gives $H_n=\widetilde H_n$.
In particular, the higher-order constructions agree.
\par
\par

\emph{4. Quantitative density bounds and cancellation.}
\par
 Proposition~\ref{thm:finiteFieldResponse} and the common base inverse \eqref{eq:logExactCollisionInverse} give
\[
\|\rho^{(n)}\|_{\mathscr X_N^{r_{\rm x}}}
\le C_{r_{\rm x}}N\left[
C_{{r_{\rm x}},n}(\zeta\varepsilon_{\rm F})^n
+\sum_{q=1}^n C_{{r_{\rm x}},q}N^{2q+[q-2]_+}\varepsilon_{\rm F}^{q}
\|\rho^{(n-q)}\|_{\mathscr X_N^{{r_{\rm x}}+[q-2]_+}}\right].
\]

Retain $\beta_n$ from \eqref{eq:fieldRegularityOrders} and set
\[
X_n=(\zeta\varepsilon_{\rm F})^{-n}
\|\rho^{(n)}\|_{\mathscr X_N^{\beta_n}}.
\]
The available source order and amplitude ratio satisfy
\[
\beta_n+[q-2]_+\le\beta_{n-q},\qquad
\frac{\varepsilon_{\rm F}^q(\zeta\varepsilon_{\rm F})^{n-q}}
{(\zeta\varepsilon_{\rm F})^n}=\zeta^{-q}\le1.
\]
For $0\le n\le J$, these bounds give
\begin{align*}
X_0&\le C_{r_{\rm x},J}N^{3J+1},\\
X_n&\le C_{r_{\rm x},J}N^{3J+1}
\left(1+\sum_{h=0}^{n-1}X_h\right)
\le\big(2C_{r_{\rm x},J}N^{3J+1}\big)^{n+1}.
\end{align*}
The last inequality follows by summing the geometric upper bounds for
$X_h$, with $C_{r_{\rm x},J}N^{3J+1}\ge1$. Since $n\le J$ is fixed,
\[
\|\rho^{(n)}\|_{\mathscr X_N^{r_{\rm x}}}
\le\|\rho^{(n)}\|_{\mathscr X_N^{\beta_n}}
\le C_{r_{\rm x},J}N^{C_{r_{\rm x},J}}
(\zeta\varepsilon_{\rm F})^n.
\]
This proves
\begin{equation}\label{eq:fieldCoefficientDensityBound}
\|\rho^{(n)}\|_{\mathscr X_N^{r_{\rm x}}}
\le C_{{r_{\rm x}},J}N^{C_{{r_{\rm x}},J}}
(\zeta\varepsilon_{\rm F})^n
\le\nu^{n{\sigma_{\rm int}}-o(1)},\qquad 1\le n\le J.
\end{equation}
The last inequality uses \eqref{eq:fieldExpansionSizes} and
$N=K\log(1/\nu)$, with $J$ fixed. Since $B_{m_1}=B_*$,
${w_{m_1,N}}({T_1})\ge1$ and $\langle m_1/N\rangle^{r_{\rm x}}\ge1$, the terminal bound follows directly:
\[
\frac{\sqrt N}{B_*}|\rho_{m_1}^{(n)}({T_1})|
\le\sqrt N\,\|\rho^{(n)}\|_{\mathscr X_N^{r_{\rm x}}}
\le\nu^{n\sigma_{\rm int}-o(1)}.
\]
Here $\sqrt N=\nu^{-o(1)}$, proving \eqref{eq:fieldMainCoefficients}.

For $1\le n\le J$, {\eqref{eq:fieldCorrectionRecursion}} gives
\[
\begin{aligned}
\mathcal L_{\rm C}H^{[n]}
&=\sum_{j=0}^n\mathcal L_{\rm C}H_j
=-\sum_{j=1}^n\mathcal T_{\rm b}H_{j-1}
=-\mathcal T_{\rm b}H^{[n-1]},\\
\mathcal A_\nu H^{[n]}
&=-\mathcal T_{\rm b}H^{[n-1]}+\mathcal T_{\rm b}H^{[n]}
=\mathcal T_{\rm b}H_n.
\end{aligned}
\]
For $n=0$ the second identity follows from $\mathcal L_{\rm C}H_0=0$.
The initial values give {$H^{[n]}(t_0)=H_{\rm in}$}, proving
\eqref{eq:fieldPartialSumIdentities}. {Subtracting the identity for $n=J$}
from {$\mathcal A_\nu H^\nu=0$}, with the same initial datum, proves
\eqref{eq:fieldFiniteUnknownError}. At the terminal mode,
\[
\frac{\sqrt N}{B_*}
|\rho_{m_1}[H^{[J]}-H_0]({T_1})|
\le\sum_{n=1}^J\frac{\sqrt N}{B_*}|\rho_{m_1}[H_n]({T_1})|
\le\nu^{\sigma_{\rm int}-o(1)}.
\]
Subtracting this bound from \eqref{eq:logExactCollisionLower} proves
\eqref{eq:fieldApproximationToBase}.

\emph{5. Spatial norms and the final residual.}
For the initial term in \eqref{eq:fieldDistributionRecursion},
$j=\pm m_0$ and $c=\pm(\zeta-m_0t_0)$.
Since $m_0T_1\le C\zeta$, every field time satisfies
$1+|c-j(t_i-t_0)|+t_i-t_0\le C\zeta$. Thus
\[
\Xi_n\le\frac{(C\zeta)^n}{n!}
\left(\int_{t_0}^s\mathfrak f(t)\,dt\right)^n
\le C_n(\zeta\varepsilon_{\rm F})^n.
\]
Equation \eqref{eq:finiteFieldDistribution}, using only the zeroth
centered derivative and the spatial exponential weight, bounds this
term by $C_{r,\ell,n}N^r(\zeta\varepsilon_{\rm F})^n$ in $H_x^rL_\ell^2$.
For $n=0$ use \eqref{eq:physicalCollisionEnergy}.

For each coefficient $\mathcal F_m^{\rm b}(\tau)$, Lemma~\ref{lem:fieldSpatialInputs} gives
\[
\|\mathbb U_{\rm b}^{(q)}(s,\tau)\mathcal F_m^{\rm b}(\tau)\|_{H_x^rL_\ell^2}
\le C_{r,\ell,q}N^{2q+2}L_N\varepsilon_{\rm F}^q
\langle m\rangle^{r+q+1},\qquad 0\le q\le J.
\]
By \eqref{eq:echoDensityEmbeddings},
\eqref{eq:fieldCoefficientDensityBound}, and the order-zero bound
\eqref{eq:logExactCollisionInverse}, each $0\le q\le n$ contributes
\[
\begin{aligned}
&\sup_{s\in I}\left\|\sum_{m\ne0}\int_{t_0}^s
\mathbb U_{\rm b}^{(q)}(s,\tau)\mathcal F_m^{\rm b}(\tau)
\rho_m^{(n-q)}(\tau)\,d\tau\right\|_{H_x^rL_\ell^2}\\
&\quad\le C_{r,\ell,q}N^{2q+2}L_N\varepsilon_{\rm F}^q T_1
\sum_{m\ne0}\langle m\rangle^{r+q+1}\sup_I|\rho_m^{(n-q)}|\\
&\quad\le C_{r,\ell,q}N^{r+3q+3}L_N\varepsilon_{\rm F}^q T_1\widetilde B
\|\rho^{(n-q)}\|_{\mathscr X_N^{r+q+1}}\\
&\quad\le A_{r,\ell,J}(N)T_1\widetilde B
\varepsilon_{\rm F}^q(\zeta\varepsilon_{\rm F})^{n-q}
\le A_{r,\ell,J}(N)T_1\widetilde B(\zeta\varepsilon_{\rm F})^n.
\end{aligned}
\]
Adding the datum term and the $n+1$ electric terms proves
\eqref{eq:quantitativeFieldDistribution}.
\par
The weighted Fourier estimate for the background field \eqref{eq:fieldFourierEnvelope} implies
\[
\|E[F^{\rm b,\nu}]\|_{L^2(I;C_x^r)}
\le C_r\varepsilon_{\rm F}
\left(\int_0^\infty e^{-2c_{\rm dec}t^\chi}\,dt\right)^{1/2}
\le C_r\varepsilon_{\rm F}.
\]
The required energy weights are controlled by one fixed moment:
\[
\mathcal E_{r+1}[f]
=\sum_{j=0}^{r+1}\|\langle v\rangle^{2r+8-2j}\partial_{x_1}^jf\|_2^2
\le\|f\|_{H_x^{r+1}L^2_{2r+8}}^2.
\]
Thus \eqref{eq:quantitativeFieldDistribution}, with spatial order
$r+1$ and moment $2r+8$, implies
\[
\sup_I\mathcal E_{r+1}[H_J/\sqrt\mu]^{1/2}
\le A_{r+1,2r+8,J}(N)\widetilde B(1+T_1)
(\zeta\varepsilon_{\rm F})^J.
\]
Substitution in \eqref{eq:electricL2SourceBound} gives
\[
\begin{aligned}
\nu^{-1/2}\sum_{j=0}^{r}\left\|\partial_{x_1}^j\!\left(\mu^{-1/2}\bigl(S_J\bigr)\right)\right\|_{L^2_{s,x}(I)\Sigma_{2r+6-2j}^*}
&\le C_r\nu^{-1/2}\varepsilon_{\rm F}
\sup_I\mathcal E_{r+1}[H_J/\sqrt\mu]^{1/2}\\
&\le C_r A_{r+1,2r+8,J}(N)\widetilde B(1+T_1)
\nu^{-1/2}\varepsilon_{\rm F}(\zeta\varepsilon_{\rm F})^J\\
&\le A_{r,J}(N)\widetilde B\frac{1+T_1}{\sqrt\nu\,\zeta}
(\zeta\varepsilon_{\rm F})^{J+1}.
\end{aligned}
\]
Here $\varepsilon_{\rm F}(\zeta\varepsilon_{\rm F})^J
=\zeta^{-1}(\zeta\varepsilon_{\rm F})^{J+1}$, and the fixed-order
factor is included in $A_{r,J}$. This proves
\eqref{eq:quantitativeFieldResidual}.
\par
Since $T_1=\zeta/m_1$, $m_1\sim N$ and $\zeta=\nu^{-1/3+o(1)}$,
\[
\frac{1+T_1}{\sqrt\nu\,\zeta}=\nu^{-1/2+o(1)},
\qquad \widetilde B\le\nu^{-KC_{\rm B}}.
\]
Together with \eqref{eq:fieldExpansionSizes}, this gives
\eqref{eq:logQuantitativeFieldResidual}.
The source has zero velocity integral by Lemma~\ref{lem:electricFormDual}.
Proposition~\ref{thm:fullForcedCoarse}, with zero initial datum,
$S_0=0$ and $S_1=S_J$, now yields
\begin{equation}\label{eq:logFullLinearRemainder}
\begin{aligned}
\|\rho[H^\nu-H^{[J]}]\|_{\ell^1_{r-1}L^\infty}
&\le C_rN^4e^{C_*N}\nu^{-1/2}\sum_{j=0}^{r}\left\|\partial_{x_1}^j\!\left(\mu^{-1/2}\bigl(S_J\bigr)\right)\right\|_{L^2_{s,x}(I)\Sigma_{2r+6-2j}^*}\\
&\le\nu^{(J+1)\sigma_{\rm int}-1/2-K(C_*+C_{\rm B})-o(1)}
=\nu^{\iota_{\rm F}-o(1)}.
\end{aligned}
\end{equation}
The polynomial $N^4$ is absorbed in $A_{r,J}$; the exponential
constants are independent of all fixed orders.
This proves \eqref{eq:commonFieldCompletionError}.
At $J=0$ the same calculation gives
\eqref{eq:overviewDirectFieldError}; the obstruction to using that
bound without corrections was established in
Lemma~\ref{lem:coarseLossVersusEcho}.
\par
\end{proof}
\par\endgroup
The common completion estimate now transfers the finite expansion's
terminal lower bound to the complete linear solution. {The positive}
power of $\nu$ is smaller than the logarithmic collision-comparison error.

\begingroup
\begin{proof}[Proof of Theorem~\ref{thm:logFullLinearEcho}]
First apply Proposition~\ref{thm:fullForcedCoarse} with {${{\bf H}}={{\bf W}}=0$} and $S_0=S_1=0$ and datum {$H_{\rm in}$}, whose weighted energy is finite at every fixed order. It constructs a solution of the complete equation. The solutions obtained at different fixed spatial orders agree by finite-energy uniqueness; denote this solution by {$H^\nu$}.

Use the single $J$ of \eqref{eq:logFieldOrderChoice}.
Proposition~\ref{thm:finiteFieldExpansion}, at order $r=4$, gives
\[
\begin{aligned}
\frac{\sqrt N}{B_*}
|\rho_{m_1}[H^{[J]}-H^{\rm C,\nu}]({T_1})|
&\le\nu^{\sigma_{\rm int}-o(1)},\\
\frac{\sqrt N}{B_*}
|\rho_{m_1}[H^\nu-H^{[J]}]({T_1})|
&\le\frac{\sqrt N}{B_*}
\|\rho[H^\nu-H^{[J]}]\|_{\ell^1_{3}L^\infty}
\le\nu^{\iota_{\rm F}-o(1)}.
\end{aligned}
\]
The second inequality uses $B_*\ge1$ and
\eqref{eq:commonFieldCompletionError}; since $\iota_{\rm F}>1>
\sigma_{\rm int}$, it is bounded by the first error scale.

The exact collision comparison and the two field estimates give
\[
\begin{aligned}
\frac{\sqrt N}{B_*}|\rho_{m_1}[H^\nu](T_1)-\rho_{m_1}^\downarrow(T_1)|
&\le C N^{9/2}\nu T_1^3
+C N^{7/2}aL_N
+\nu^{\sigma_{\rm int}-o(1)}+\nu^{\iota_{\rm F}-o(1)}\\
&\le C[\log(1/\nu)]^{3{\rm s}_0+21/2-3p}
+\nu^{1/3-o(1)}\longrightarrow0.
\end{aligned}
\]
Here $\iota_{\rm F}>1>\sigma_{\rm int}>1/2$, and
$N^{9/2}e^{-cN^2}=o(\nu^{1/3})$.
Proposition~\ref{thm:chain} and Lemma~\ref{lem:logarithmicScales} imply
\[
\begin{aligned}
|\rho_{m_1}[H^\nu](T_1)|
&\ge|\rho_{m_1}^\downarrow(T_1)|
-|\rho_{m_1}[H^\nu](T_1)-\rho_{m_1}^\downarrow(T_1)|
\ge(c-o(1))B_*/\sqrt N,\\
\nu^{1/3}(T_1-t_0)
&\le C[\log(1/\nu)]^{{\rm s}_0+2-p}\longrightarrow0.
\end{aligned}
\]
This proves \eqref{eq:logFullLinearComparison}--\eqref{eq:logFullLinearLower}.
\par
The initial norm in \eqref{eq:logLinearObservable} follows from \eqref{eq:maininitialnorm}. Since $m_1({T_1}-t_0)\ge\zeta/2$ and $B_*=\exp(c_{\rm e}N+O(1))$, the terminal observable is bounded below by
\[
c\epsilon_\nu\zeta^{-3}B_*/\sqrt N
=\nu^{4/3-Kc_{\rm e}+o(1)}.
\]
The assumed inequality on $K$ makes this tend to infinity. All quantities in this calculation belong to the linear problem.
\end{proof}
\par\endgroup
We finish with the spatial energy and density upper bounds needed for
the nonlinear construction. The electric source estimate in
Lemma~\ref{lem:electricFormDual} uses these bounds at one higher fixed spatial order.
\par

\begin{proof}[Proof of Proposition~\ref{prop:linearNonlinearInput}]
Write {$H=H^\nu$} and use the same $J$ fixed in
\eqref{eq:logFieldOrderChoice}. The common completion estimate
\eqref{eq:commonFieldCompletionError}, applied at order ${\rm s}$, gives
\[
\|\rho[H-H^{[J]}]\|_{\ell^1_{{\rm s}-1}L^\infty}
\le\nu^{\iota_{\rm F}-o(1)},\qquad \iota_{\rm F}>1.
\]
The finite field coefficients satisfy \eqref{eq:fieldCoefficientDensityBound} {and the corresponding order-zero estimate.} Removing the spatial gain weight yields
\[
\|\rho[H^{[J]}]\|_{\ell^1_{{\rm s}-1}L^\infty}
\le C_{\rm s}N^{{\rm s}-1}\widetilde{B}
\|\rho[H^{[J]}]\|_{\mathscr X_N^{{\rm s}-1}}
\le A_{\rm s}(N)\widetilde{B}.
\]
The preceding embedding is \eqref{eq:echoDensityEmbeddings}.
The Gaussian datum \eqref{eq:highdatum} satisfies
$\mathcal E_{\rm s}[H_{\rm in}/\sqrt\mu]^{1/2}\le C_{\rm s}N^{\rm s}$.
Hence \eqref{eq:densityNormComparison} and
Proposition~\ref{prop:fullLinearEnergy} give
\[
\begin{aligned}
\|\rho[H]\|_{L^\infty(I;H_x^{{\rm s}-1})}
&\le\sqrt{2\pi}\|\rho[H]\|_{\ell^1_{{\rm s}-1}L^\infty(I)}
\le A_{\rm s}(N)\widetilde B,\\
\|H/\sqrt\mu\|_{\mathcal N_{\rm s}}
&\le C_{\rm s}\big(N^{\rm s}
+T_1[1+\varepsilon_N(1+T_1)]A_{\rm s}(N)\widetilde B\big) \le A_{\rm s}(N)(1+T_1)\widetilde B.
\end{aligned}
\]
Enlarging the subexponential factor $A_{\rm s}$ therefore yields
\par

\begin{equation}\label{eq:linearInputBase}
\|\rho[H]\|_{\ell^1_{{\rm s}-1}L^\infty}
\le A_{\rm s}(N)\widetilde{B},\qquad
\|H/\sqrt \mu\|_{\mathcal N_{\rm s}}\le A_{\rm s}(N)(1+{T_1})\widetilde{B}.
\end{equation}
The construction at this higher fixed order agrees with the order-four solution by uniqueness in the energy class.

This proves \eqref{eq:linearHierarchyInput}. All additional factors are
subexponential in $N$; the coefficient $C_{\rm B}$ is unchanged.
\end{proof}
\par

\begingroup
\section{Nonlinear estimates and proof of the main theorem}\label{sec:nonlinear}\label{part:nonlinear}
\begingroup
We prove Proposition~\ref{thm:exactNonlinearSolution} and then complete
Theorem~\ref{thm:logNonlinearInflation}. The remainder equation and the
initial temperature choice are \eqref{eq:nonlinearExpansionRoadmap} and
\eqref{eq:initialTemperatureCorrection}. We first estimate the quadratic
source and construct the exact nonlinear remainder by iteration.
We then identify the energy solution with a local classical solution;
the mixed-derivative estimates continue this identification to $T_1$.
Conservation and the density-to-Sobolev estimate complete the proof.
All quantitative comparisons use the exact linear solution of
Theorem~\ref{thm:logFullLinearEcho}; the classical continuation is carried
out for each fixed $\nu$.
\par\endgroup
\par\endgroup
\subsection{The quadratic residual and coefficient smallness}\label{sec:quadraticResponse}
Proposition~\ref{prop:linearNonlinearInput} supplies the density and spatial
energy of $H^\nu$. 
Lemma~\ref{lem:electricFormDual} and the collision form estimate
give the quadratic source bounds needed in
Proposition~\ref{thm:fullForcedCoarse}.
\par We use the energy and comparison
norms in \eqref{eq:fullLinearEnergyNorms}--\eqref{eq:errorComparisonNorm}
{and the source integrals in \eqref{eq:quadraticSourceSpaces}.}
\par

Use $A_{\rm s}$ and $\eta_{\nu,{\rm s}}$ from
\eqref{eq:additionalSpatialFactors}.

\begin{prop}\label{prop:quadraticApproximation}\label{lem:approximationCoefficients}\label{lem:nonlinearPowerEstimates}
Use \eqref{eq:comparisonLogParameters} and \eqref{eq:nonlinearTwoScales}
with fixed integers ${\rm s}\ge4$, ${\rm s}_0\ge11$, and fixed $p>{\rm s}_0+7/2$, $K>0$.
The factors in \eqref{eq:additionalSpatialFactors} can be chosen so that
the following assertions hold, with $\Theta$ given by
\eqref{eq:nonlinearRoadmapMargin}. The quantity $A_{\rm s}(N)\widetilde{B}$ bounds the
linear solution, $2\delta_\nu A_{\rm s}(N)\widetilde{B}$ is the radius of the remainder ball, and
$\eta_{\nu,{\rm s}}$ controls both the forced residual and the iteration.
\emph{(i) Quadratic defect and relative source bound.}
The defect of the first-order approximation is exactly
\begin{equation}\label{eq:quadraticResidualIdentity}
\mathcal R_{\rm a}:=\mathcal A_\nu H_{\rm a}-\mathbb Q_\nu(H_{\rm a},H_{\rm a})
=E[H_{\rm a}]\partial_{v_1}H_{\rm a}-\nu Q(H_{\rm a},H_{\rm a}).
\end{equation}
\par

{The residual satisfies}
\begin{equation}\label{eq:quadraticResidualSource}
\begin{aligned}
&\nu^{-1/2}\sum_{j=0}^{\rm s}\left\|\partial_{x_1}^j\!\left(\mu^{-1/2}\bigl(\mathcal R_{\rm a}\bigr)\right)\right\|_{L^2_{s,x}(I)\Sigma_{2{\rm s}+6-2j}^*}\\
&\quad\le C_{\rm s}\big(\nu^{-1/2}(1+T_1)^{3/2}+(1+T_1)^2\big)
(\delta_\nu A_{\rm s}(N)\widetilde B)^2.
\end{aligned}
\end{equation}

\par

The prescribed residual satisfies
\begin{equation}\label{eq:residualRelativeEstimate}
\frac{A_{\rm s}e^{C_*N}\nu^{-1/2}\sum_{j=0}^{\rm s}\left\|\partial_{x_1}^j\!\left(\mu^{-1/2}\bigl(\mathcal R_{\rm a}\bigr)\right)\right\|_{L^2_{s,x}(I)\Sigma_{2{\rm s}+6-2j}^*}}
{\delta_\nu B_*/\sqrt N}\le\nu^{\Theta-o(1)}.
\end{equation}

\par

\emph{(ii) Coefficients of the exact error equation.}
The coefficient bounds are
\begin{equation}\label{eq:approximationCoefficientNorms}
\|\rho[H_{\rm a}]\|_{\ell^1_{\rm s}L^\infty}
\le\delta_\nu A_{\rm s}(N)\widetilde B,\qquad
\|H_{\rm a}/\sqrt\mu\|_{\mathcal N_{{\rm s}+1}}
\le(1+T_1)\delta_\nu A_{\rm s}(N)\widetilde B.
\end{equation}
\par

Furthermore,
\begin{equation}\label{eq:approximationContinuity}
\int_{\mathbb R^3}\mathcal R_{\rm a}\,dv=0,
\qquad \partial_s\rho[H_{\rm a}]+\partial_{x_1}j[H_{\rm a}]=0.
\end{equation}
\emph{(iii) Smallness conditions.}
The comparison factors satisfy
\begin{equation}\label{eq:quadraticComparisonPowers}
\eta_{\nu,{\rm s}}\le\nu^{\Theta+K(C_{\rm B}-c_{\rm e})-o(1)},\qquad
\frac{\sqrt N A_{\rm s}(N)\widetilde{B}}{B_*}\le\nu^{-K(C_{\rm B}-c_{\rm e})-o(1)}{.}
\end{equation}
If $\Theta>0$, these bounds imply the prescribed-coefficient
smallness required in Proposition~\ref{thm:fullForcedCoarse}
with radius $2\delta_\nu A_{\rm s}(N)\widetilde B$.
\par
All $o(1)$ terms are taken with the displayed orders and parameters fixed.
\end{prop}

\begin{proof}
\emph{1. The residual and the coefficients of the approximation.}
Since {$\mathcal A_\nu H^\nu=0$} and $\delta_\nu$ is constant,
{$\mathcal A_\nu H_{\rm a}=0$}. Expanding the definition of
$\mathbb Q_\nu$ gives \eqref{eq:quadraticResidualIdentity}.

Proposition~\ref{prop:linearNonlinearInput} at order ${\rm s}+1$ gives
\eqref{eq:approximationCoefficientNorms} after multiplication by $\delta_\nu$.
Lemma~\ref{lem:electricFormDual} and the collision form estimate
\eqref{eq:mixedSpatialLandauForm} yield
\[
\begin{aligned}
&\nu^{-1/2}\sum_{j=0}^{\rm s}\left\|\partial_{x_1}^j\!\left(\mu^{-1/2}\bigl(E[H_{\rm a}]\partial_{v_1}H_{\rm a}\bigr)\right)\right\|_{L^2_{s,x}(I)\Sigma_{2{\rm s}+6-2j}^*}\\
&\quad\le C_{\rm s}\sqrt{(T_1-t_0)/\nu}\,
\|\rho[H_{\rm a}]\|_{\ell^1_{{\rm s}-1}L^\infty}
\sup_I\mathcal E_{{\rm s}+1}[H_{\rm a}/\sqrt\mu]^{1/2}\\
&\quad\le C_{\rm s}\nu^{-1/2}(1+T_1)^{3/2}(\delta_\nu A_{\rm s}\widetilde B)^2,\\
&\nu^{-1/2}\sum_{j=0}^{\rm s}\left\|\partial_{x_1}^j\!\left(\mu^{-1/2}\bigl(\nu Q(H_{\rm a},H_{\rm a})\bigr)\right)\right\|_{L^2_{s,x}(I)\Sigma_{2{\rm s}+6-2j}^*}\\
&\quad\le C_{\rm s}\|H_{\rm a}/\sqrt\mu\|_{\mathcal N_{\rm s}}^2
\le C_{\rm s}(1+T_1)^2(\delta_\nu A_{\rm s}\widetilde B)^2.
\end{aligned}
\]
Adding these inequalities proves \eqref{eq:quadraticResidualSource}.
\par

The electric residual and the Landau residual are velocity divergences. Their velocity integrals vanish by Gaussian integrability and the form estimate, using velocity cutoffs to justify the integrations. The same calculation in the linear equation for {$H_{\rm a}$} {proves the local mass conservation of {$H_{\rm a}$}} in \eqref{eq:approximationContinuity}.

\emph{2. Powers of $\nu$ and the terminal comparison.}
Write $A_{\rm s}=e^{o(N)}=\nu^{-o(1)}$. Lemma~\ref{lem:logarithmicScales} gives
\begin{equation}\label{eq:logNonlinearScales}
\delta_\nu=\nu^{({\rm s}_0+1)/3+o(1)},\qquad
\zeta=\nu^{-1/3+o(1)},\qquad 1+{T_1}=\nu^{-1/3+o(1)}.
\end{equation}
Hence
\[
\begin{aligned}
(1+T_1)^2&=\nu^{-2/3+o(1)},&
\nu^{-1/2}(1+T_1)^{3/2}&=\nu^{-1+o(1)},\\
A_{\rm s}(N)\widetilde B&\le\nu^{-KC_{\rm B}-o(1)},&
e^{C_*N}&=\nu^{-KC_*}.
\end{aligned}
\]
Substitution into \eqref{eq:additionalSpatialFactors} gives
\[
\eta_{\nu,{\rm s}}
\le\nu^{({\rm s}_0+1)/3-K(C_*+C_{\rm B})-1-o(1)}
=\nu^{\Theta+K(C_{\rm B}-c_{\rm e})-o(1)}.
\]
\par

Furthermore,
\[
B_*=e^{c_{\rm e}N+O(1)}=\nu^{-Kc_{\rm e}+o(1)},\qquad
\frac{\sqrt N A_{\rm s}(N)\widetilde{B}}{B_*}\le\nu^{-K(C_{\rm B}-c_{\rm e})-o(1)},
\]
This proves \eqref{eq:quadraticComparisonPowers}. The same substitution gives
\begin{equation}\label{eq:nonlinearCoefficientPower}
2(1+T_1)\delta_\nu A_{\rm s}(N)\widetilde B
\le\nu^{{\rm s}_0/3-KC_{\rm B}-o(1)},\qquad
\frac{{\rm s}_0}{3}-KC_{\rm B}
=\Theta+K(C_{\rm B}-c_{\rm e})+\frac23+KC_*.
\end{equation}
By \eqref{eq:quadraticResidualSource},
\par

\[
A_{\rm s}e^{C_*N}{\nu^{-1/2}\sum_{j=0}^{\rm s}\left\|\partial_{x_1}^j\!\left(\mu^{-1/2}\bigl(\mathcal R_{\rm a}\bigr)\right)\right\|_{L^2_{s,x}(I)\Sigma_{2{\rm s}+6-2j}^*}}
\le C_{\rm s}\delta_\nu A_{\rm s}(N)\widetilde{B}\eta_{\nu,{\rm s}}.
\]
Division by $\delta_\nu B_*/\sqrt N$ proves \eqref{eq:residualRelativeEstimate}.
The definition of $\eta_{\nu,{\rm s}}$ also bounds all additional
integrated energy coefficients at once:
\[
\begin{aligned}
&2(1+T_1)^2\delta_\nu A_{\rm s}(N)\widetilde B
+4\big(\delta_\nu A_{\rm s}(N)\widetilde B\big)^2
\big[\nu^{-1}(1+T_1)+(1+T_1)^2\big]\\
&\qquad\le\eta_{\nu,{\rm s}}+2\eta_{\nu,{\rm s}}^2=o(1)
\qquad(\Theta>0).
\end{aligned}
\]
\par

The terminal comparison cancels the extra $K(C_{\rm B}-c_{\rm e})$
in the contraction exponent against the maximal-to-terminal gain ratio,
leaving precisely $\Theta$. The fixed field order $J$ enters only the
subexponential factors.
\end{proof}

\subsection{Construction and comparison of the exact solution}\label{sec:nonlinearConstruction}\label{sec:exactNonlinearRemainder}

We prove Proposition~\ref{thm:exactNonlinearSolution}, stated in
Section~\ref{sec:roadmapNonlinear}:
the first-order approximation generates an exact solution, {and the remainder}
is smaller than the terminal linear density. The construction has three steps:
bound the iterates at order ${\rm s}$, contract their differences at order ${\rm s}-1$,
and compare the limiting remainder with $\delta_\nu B_*/\sqrt N$.

The common estimate used to solve each iterate is
Proposition~\ref{thm:fullForcedCoarse}, with {${{\bf H}}=H_{\rm a}$} and
$\mathfrak r=2\delta_\nu A_{\rm s}(N)\widetilde{B}$.

\begingroup
\begin{proof}[Proof of Proposition~\ref{thm:exactNonlinearSolution}]
\emph{1. Initial remainder and the exact equation.}
Proposition~\ref{prop:quadraticApproximation} and $\Theta>0$ give
\[
\eta_{\nu,{\rm s}}\le\nu^{\Theta+K(C_{\rm B}-c_{\rm e})-o(1)}\to0,
\qquad 2(1+{T_1})\delta_\nu A_{\rm s}(N)\widetilde{B}=o(1).
\]
The initial remainder $W_{\rm in}$ is defined in
\eqref{eq:initialTemperatureCorrection}. It has zero total mass and
satisfies \eqref{eq:initialTemperatureSpatialBound}.
\par
Since $\zeta-m_0t_0\ge\zeta/2$ and $\zeta=\nu^{-1/3+o(1)}$,
\begin{equation}\label{eq:smallInitialRemainder}
{
A_{\rm s}e^{C_*N}\mathcal E_{\rm s}[W_{\rm in}/\sqrt \mu]^{1/2} =o(\delta_\nu A_{\rm s}(N)\widetilde{B}\eta_{\nu,{\rm s}}) =o(2\delta_\nu A_{\rm s}(N)\widetilde{B}).
}
\end{equation}
Substituting {$F^\nu=F^{\rm b,\nu}+H_{\rm a}+W$} and using
\eqref{eq:quadraticResidualIdentity} gives
\begin{equation}\label{eq:exactAmplitudeError}
{
\mathcal A_\nu W-\mathbb Q_\nu(H_{\rm a}+W,W)
-\mathbb Q_\nu(W,H_{\rm a})=-\mathcal R_{\rm a},\qquad W(t_0)=W_{\rm in}.
}
\end{equation}
\emph{2. Iteration at the higher order.}

Set {$W^{(0)}=0$} and define
\begin{equation}\label{eq:nonlinearPicard}
{
\begin{aligned}
&\mathcal A_\nu W^{(n+1)}-\mathbb Q_\nu(H_{\rm a}+W^{(n)},W^{(n+1)})
-\mathbb Q_\nu(W^{(n+1)},H_{\rm a})=-\mathcal R_{\rm a},\\
&W^{(n+1)}(t_0)=W_{\rm in}.
\end{aligned}
}
\end{equation}
Equations \eqref{eq:approximationCoefficientNorms} and
\eqref{eq:approximationContinuity} verify the prescribed-coefficient hypotheses
of Proposition~\ref{thm:fullForcedCoarse} with {${{\bf H}}=H_{\rm a}$}.
{The required smallness condition is exactly} $\eta_{\nu,{\rm s}}=o(1)$.
If {$\|W^{(n)}\|_{\mathscr V_{\rm s}}\le2\delta_\nu A_{\rm s}(N)\widetilde{B}$} and {$W^{(n)}$}
satisfies local mass conservation, that proposition with {${{\bf W}}=W^{(n)}$}
{and sources $S_0=0$, $S_1=-\mathcal R_{\rm a}$ gives a solution
on all of $I$ and}
\begin{align*}
\|W^{(n+1)}\|_{\mathscr V_{\rm s}}
&\le A_{\rm s}e^{C_*N}\left(\mathcal E_{\rm s}[W_{\rm in}/\sqrt \mu]^{1/2}
+{\nu^{-1/2}\sum_{j=0}^{\rm s}\left\|\partial_{x_1}^j\!\left(\mu^{-1/2}\bigl(\mathcal R_{\rm a}\bigr)\right)\right\|_{L^2_{s,x}(I)\Sigma_{2{\rm s}+6-2j}^*}}\right)\\
&\le C_{\rm s}\delta_\nu A_{\rm s}(N)\widetilde{B}\eta_{\nu,{\rm s}}
\le\delta_\nu A_{\rm s}(N)\widetilde{B}.
\end{align*}
Here \eqref{eq:smallInitialRemainder}, \eqref{eq:quadraticResidualSource}
and $\eta_{\nu,{\rm s}}=o(1)$ give the last two inequalities. Also
\[
\int\mathcal R_{\rm a}\,dv=0
\quad\Longrightarrow\quad
\partial_s\rho[W^{(n+1)}]+\partial_{x_1}j[W^{(n+1)}]=0.
\]
Thus induction closes. Reality follows from real coefficients and prescribed-equation uniqueness.

\emph{3. The derivative loss and contraction at the adjacent order.}
For $n\ge1$, subtracting consecutive equations gives
\begin{equation}\label{eq:nonlinearDifferenceEquation}
{
\begin{aligned}
&\mathcal A_\nu (W^{(n+1)}-W^{(n)})-\mathbb Q_\nu(H_{\rm a}+W^{(n)},(W^{(n+1)}-W^{(n)}))\\
&\quad-\mathbb Q_\nu((W^{(n+1)}-W^{(n)}),H_{\rm a})=\mathbb Q_\nu(W^{(n)}-W^{(n-1)},W^{(n)}),\\
&(W^{(n+1)}-W^{(n)})(t_0)=0.
\end{aligned}
}
\end{equation}
All positive-index iterates have the same initial value. We estimate
the difference source at order ${\rm s}-1$ to prove contraction.
In Step~4 we pass to the equation at order ${\rm s}-2$ and then
recover strong continuity at order ${\rm s}$.
For real zero-total-mass $H_1,W$, Lemma~\ref{lem:electricFormDual} with
$r={\rm s}-1$ and \eqref{eq:mixedSpatialLandauForm} give
\begin{equation}\label{eq:lowerOrderDifferenceSource}
\begin{aligned}
&\nu^{-1/2}\sum_{j=0}^{{\rm s}-1}\left\|\partial_{x_1}^j\!\left(\mu^{-1/2}\bigl(\mathbb Q_\nu(H_1,W)\bigr)\right)\right\|_{L^2_{s,x}(I)\Sigma_{2{\rm s}+4-2j}^*}\\
&\quad\le C_{\rm s}\big((1+T_1)^2+\nu^{-1/2}(1+T_1)^{3/2}\big)
\|H_1\|_{\mathscr V_{{\rm s}-1}}\|W\|_{\mathscr V_{\rm s}}.
\end{aligned}
\end{equation}
Indeed, the electric part uses
${\sup_I\mathcal E_{\rm s}[W/\sqrt\mu]^{1/2}
\le(1+T_1)\|W\|_{\mathscr V_{\rm s}}}$, while the collision part satisfies
\[
\begin{aligned}
&\nu^{-1/2}\sum_{j=0}^{{\rm s}-1}\left\|\partial_{x_1}^j\!\left(\mu^{-1/2}\bigl(\nu Q(H_1,W)\bigr)\right)\right\|_{L^2_{s,x}(I)\Sigma_{2{\rm s}+4-2j}^*}\\
&\quad\le C_{\rm s}\|H_1/\sqrt\mu\|_{\mathcal N_{{\rm s}-1}}
\|W/\sqrt\mu\|_{\mathcal N_{{\rm s}-1}} \le C_{\rm s}(1+T_1)^2\|H_1\|_{\mathscr V_{{\rm s}-1}}\|W\|_{\mathscr V_{\rm s}}.
\end{aligned}
\]

\par

For $\|W\|_{\mathscr V_{\rm s}}\le 2\delta_\nu A_{\rm s}(N)\widetilde{B}$ and a coefficient {${{\bf W}}$} satisfying
\eqref{eq:prescribedRemainderClass} with {${{\bf H}}=H_{\rm a}$},
{the solution {$H$} of \eqref{eq:fullForcedEquation} with zero initial
datum and {${S_0=0,\ S_1=\mathbb Q_\nu(H_1,W)}$}} obeys, by
{Proposition~\ref{thm:fullForcedCoarse} and \eqref{eq:additionalSpatialFactors},}
\begin{equation}\label{eq:lowerOrderContraction}
\|H\|_{\mathscr V_{{\rm s}-1}}
\le C_{\rm s}\eta_{\nu,{\rm s}}\|H_1\|_{\mathscr V_{{\rm s}-1}}.
\end{equation}
Taking {$H_1=W^{(n)}-W^{(n-1)}$}, {$W=W^{(n)}$} and {${{\bf W}}=W^{(n)}$} gives
\[
\begin{aligned}
\|W^{(n+1)}-W^{(n)}\|_{\mathscr V_{{\rm s}-1}}
&\le C_{\rm s}\eta_{\nu,{\rm s}}
\|W^{(n)}-W^{(n-1)}\|_{\mathscr V_{{\rm s}-1}} \le\tfrac12\|W^{(n)}-W^{(n-1)}\|_{\mathscr V_{{\rm s}-1}},\\
\sum_{n\ge1}\|W^{(n)}-W^{(n-1)}\|_{\mathscr V_{{\rm s}-1}}&<\infty.
\end{aligned}
\]
\emph{4. Passage to the exact solution and uniqueness.}
Consequently, there is a physical remainder {$W$} such that
\[
\begin{aligned}
\frac{W^{(n)}}{\sqrt \mu }&\longrightarrow\frac{W}{\sqrt \mu }
 &&\text{in }C(I;\mathscr E_{{\rm s}-1})
   \text{ and the integrated form norm},\\
\frac{W^{(n)}}{\sqrt \mu }&\rightharpoonup\frac{W}{\sqrt \mu }
 &&\text{in the integrated order-${\rm s}$ form space},\\
\rho[W^{(n)}]&\longrightarrow\rho[W]
 &&\text{in }\ell^1_{{\rm s}-2}L^\infty.
\end{aligned}
\]
The uniform order-${\rm s}$ bound passes to the limit:
\[
\sup_I\mathcal E_{\rm s}[W/\sqrt\mu]
\le\liminf_{n\to\infty}\sup_I\mathcal E_{\rm s}[W^{(n)}/\sqrt\mu]<\infty.
\]
For a smooth periodic test $\varphi$ compactly supported in velocity,
integrate the highest spatial derivative once by parts. Then
\[
\begin{aligned}
&\sum_{j=0}^{\rm s}
\left\langle\partial_{x_1}^j\frac{W(s)-W(s')}{\sqrt\mu},
\partial_{x_1}^j\varphi\right\rangle_{2{\rm s}+6-2j}\\
&\quad=\sum_{j=0}^{{\rm s}-1}
\left\langle\partial_{x_1}^j\frac{W(s)-W(s')}{\sqrt\mu},
\partial_{x_1}^j\varphi\right\rangle_{2{\rm s}+6-2j}
-\left\langle\partial_{x_1}^{{\rm s}-1}\frac{W(s)-W(s')}{\sqrt\mu},
\partial_{x_1}^{{\rm s}+1}\varphi\right\rangle_6,\\
&\left|\sum_{j=0}^{\rm s}
\left\langle\partial_{x_1}^j\frac{W(s)-W(s')}{\sqrt\mu},
\partial_{x_1}^j\varphi\right\rangle_{2{\rm s}+6-2j}\right|
\le C_\varphi\mathcal E_{{\rm s}-1}
\left[\frac{W(s)-W(s')}{\sqrt\mu}\right]^{1/2}\longrightarrow0.
\end{aligned}
\]
Here $s\to s'$, and compact velocity support controls the difference
between the order-${\rm s}$ and order-${\rm s}-1$ weights. For
$\varphi\in\mathscr E_{\rm s}$, choose smooth tests $\varphi_n$ with
$\mathcal E_{\rm s}[\varphi-\varphi_n]\to0$. The remaining pairing is bounded by
\[
\left|\sum_{j=0}^{\rm s}
\left\langle\partial_{x_1}^j\frac{W(s)-W(s')}{\sqrt\mu},
\partial_{x_1}^j(\varphi-\varphi_n)\right\rangle_{2{\rm s}+6-2j}\right|
\le2\sup_I\mathcal E_{\rm s}[W/\sqrt\mu]^{1/2}
\mathcal E_{\rm s}[\varphi-\varphi_n]^{1/2}.
\]
First let $s\to s'$ and then $n\to\infty$. This proves
$W/\sqrt\mu\in C_{\rm w}(I;\mathscr E_{\rm s})$.
\par
Lower semicontinuity and modewise uniform convergence imply
\[
\sum_{k\ne0}\langle k\rangle^{{\rm s}-1}\sup_I|\rho_k[W]|
\le\liminf_n\sum_{k\ne0}\langle k\rangle^{{\rm s}-1}\sup_I|\rho_k[W^{(n)}]|,
\]
and hence
\begin{equation}\label{eq:nonlinearWeakLimitBound}
\|W\|_{\mathscr V_{\rm s}}
\le C_{\rm s}\delta_\nu A_{\rm s}(N)\widetilde{B}\eta_{\nu,{\rm s}}.
\end{equation}

To pass to the equation, use the single bilinear identity
\[
{
\mathbb Q_\nu(W^{(n)},W^{(n+1)})-\mathbb Q_\nu(W,W) =\mathbb Q_\nu(W^{(n)}-W,W^{(n+1)}) +\mathbb Q_\nu(W,W^{(n+1)}-W).
}
\]
Estimate \eqref{eq:lowerOrderDifferenceSource}, with ${\rm s}$ replaced
by ${\rm s}-1$, gives
\[
{\nu^{-1/2}\sum_{j=0}^{{\rm s}-2}\left\|\partial_{x_1}^j\!\left(\mu^{-1/2}\bigl(\mathbb Q_\nu(U,V)\bigr)\right)\right\|_{L^2_{s,x}(I)\Sigma_{2{\rm s}+2-2j}^*}}
\le C_{\nu,{\rm s}}\|U\|_{\mathscr V_{{\rm s}-2}}
\|V\|_{\mathscr V_{{\rm s}-1}}.
\]
Apply it to the two pairs in the preceding identity. Since
\[
{
\|W^{(n)}-W\|_{\mathscr V_{{\rm s}-2}} \|W^{(n+1)}\|_{\mathscr V_{{\rm s}-1}} +\|W\|_{\mathscr V_{{\rm s}-2}} \|W^{(n+1)}-W\|_{\mathscr V_{{\rm s}-1}}\longrightarrow0,
}
\]
{the preceding product bounds give
\[
\sum_{j=0}^{{\rm s}-2}
\left\|\partial_{x_1}^j\left\{\mu^{-1/2}
\bigl(\mathbb Q_\nu(W^{(n)},W^{(n+1)})-\mathbb Q_\nu(W,W)\bigr)\right\}\right\|_
{L^2_{s,x}(I)\Sigma_{2({\rm s}-2)+6-2j}^*}\longrightarrow0.
\]
}
This order leaves the velocity-weight margin required by
Lemma~\ref{lem:electricFormDual} for the strongly convergent second factor in $\mathbb Q_\nu(W^{(n)},W^{(n+1)})$.
The convergence suffices for the distributional equation; here $\nu$ is fixed.
\par
 The linear
terms pass against smooth compactly supported tests by energy and form
convergence. Thus {$W$} solves \eqref{eq:exactAmplitudeError}, and integration
in velocity gives
\begin{equation}\label{eq:errorContinuity}
\partial_s\rho[W]+\partial_{x_1}j[W]=0.
\end{equation}

Now treat {$W$} as a prescribed coefficient. {Weak continuity of {$W$},} \eqref{eq:nonlinearWeakLimitBound} and \eqref{eq:errorContinuity} verify Proposition~\ref{thm:fullForcedCoarse}. Let
\[
{
\begin{aligned}
&\mathcal A_\nu \widetilde W-\mathbb Q_\nu(H_{\rm a}+W,\widetilde W)
-\mathbb Q_\nu(\widetilde W,H_{\rm a})=-\mathcal R_{\rm a},\\
&\widetilde W(t_0)=W_{\rm in},\qquad \widetilde W/\sqrt \mu \in C(I;\mathscr E_{\rm s}).
\end{aligned}
}
\]
At order ${\rm s}-1$, both {$\widetilde W$} and {$W$} solve this same equation. Uniqueness gives {$\widetilde W=W$}, proving strong order-${\rm s}$ continuity and \eqref{eq:nonlinearComparisonClosed}.

If {$W,\widehat W$} are two remainders in the order-${\rm s}$ ball, then
\[
{
\begin{aligned}
&\mathcal A_\nu (W-\widehat W)-\mathbb Q_\nu(H_{\rm a}+W,(W-\widehat W))
-\mathbb Q_\nu((W-\widehat W),H_{\rm a})=\mathbb Q_\nu(W-\widehat W,\widehat W),\\
&\|W-\widehat W\|_{\mathscr V_{{\rm s}-1}}\le C_{\rm s}\eta_{\nu,{\rm s}}\|W-\widehat W\|_{\mathscr V_{{\rm s}-1}}.
\end{aligned}
}
\]
Their difference has zero initial datum, so {$W=\widehat W$} for $C_{\rm s}\eta_{\nu,{\rm s}}<1$.

\emph{5. Compatibility of orders and terminal comparison.}
Constructions at different fixed orders agree because
\[
\mathscr V_{\rm s}\hookrightarrow\mathscr V_4,
\qquad
\frac{A_{\rm s}(N)\eta_{\nu,{\rm s}}}{2A_4(N)}
\le\nu^{\Theta+K(C_{\rm B}-c_{\rm e})-o(1)}\to0,
\qquad A_{\rm s}/A_4=e^{o(N)}.
\]
They therefore lie in the common order-four uniqueness ball. Adding the background, linear-solution and remainder estimates gives
\[
\sup_I\mathcal E_{\rm s}[f^\nu]^{1/2}
\le C_{\rm s}\varepsilon_N+2C_{\rm s}(1+{T_1})\delta_\nu A_{\rm s}(N)\widetilde{B}=o(1),\qquad \int_I\mathcal D_{\rm s}[f^\nu]\,ds<\infty,
\]
by \eqref{eq:nonlinearCoefficientPower}. Finally,
\[
\frac{|\rho_{m_1}[W]({T_1})|}{\delta_\nu B_*/\sqrt N}
\le C\frac{\sqrt N A_4(N)\widetilde{B}}{B_*}\eta_{\nu,4}
\le\nu^{\Theta-o(1)}\le C_\omega\nu^\omega,\qquad 0<\omega<\Theta.
\]
Since $\Theta>0$, subtraction from the linear lower bound proves \eqref{eq:exactNonlinearEcho}. {The second lower bound in \eqref{eq:exactNonlinearEcho} is}
\[
\zeta^{-3}B_*/\sqrt N=\nu^{1-Kc_{\rm e}+o(1)}\to\infty
\qquad(Kc_{\rm e}>1).
\]

\end{proof}
\par\endgroup
\begingroup
\subsection{Local classical solutions and continuation}\label{sec:classicalConclusion}\label{sec:nonlinearMixedRegularity}

The preceding subsection constructs the energy solution and proves the
quantitative comparison. We now identify that solution with a local
classical VPL solution supplied by the existing local theory. The
mixed-derivative equations give bounds in the original Maxwellian
weight, which verify the continuation criterion. Constants in this subsection may depend on the fixed
$\nu$; they do not enter the terminal error estimate.

We use the Gaussian normalizations of Appendix~\ref{sec:prescribedSolvability},
\eqref{eq:gaussianEnergyPullback}--\eqref{eq:broaderGaussianNorms}.
The original weight corresponds to $\varpi=1$; the single
weaker weight $\varpi=5/4$ is used only for uniqueness.

\begin{lem}\label{lem:localClassicalIdentification}
Fix $\nu>0$ and an integer $k\ge6$. Let $G_\tau\ge0$ have mean
density one and satisfy
\begin{equation}\label{eq:localGaussianDatum}
\|e^{15|v|^2/64}G_\tau\|_{H^k_{x_1,v}}<\infty.
\end{equation}
The periodic VPL equation has a unique local nonnegative classical
solution $G$ with $G(\tau)=G_\tau$ and
\begin{equation}\label{eq:localGaussianClass}
e^{7|v|^2/32}G\in C([\tau,\tau+h];H^k_{x_1,v})
\end{equation}
for some $h>0$. At fixed $\nu$ and $k$, the local existence time has a
positive lower bound when the norm in \eqref{eq:localGaussianDatum}
ranges over a bounded set.

Suppose, in addition, that a real VPL energy solution $F$ on $I$
satisfies $F(\tau)=G_\tau$ and
\begin{equation}\label{eq:localIdentificationCoefficients}
\begin{gathered}
(F-\mu)/\sqrt\mu\in C(I;\mathscr E_k),\\
\sup_I\mathcal E_k[(F-\mu)/\sqrt\mu]^{1/2}\le\epsilon,
\qquad \int_I\mathcal D_k[(F-\mu)/\sqrt\mu]\,ds<\infty.
\end{gathered}
\end{equation}
Assume also local mass conservation. If $\epsilon$ is sufficiently
small, then $G=F$ on $[\tau,\tau+h]\cap I$.
\end{lem}

\begin{proof}
\emph{1. The local theory and the Gaussian margin.}
We use the VPL local existence, uniqueness and continuation result
stated in \cite[Remark~3.2, equation~(3.7)]{ref4}, based on the
Landau construction in \cite[Theorem~1.1 and Section~3]{HendersonSnelsonTarfulea2019}.
The latter construction works at each finite order $k\ge4$, with a
decreasing Gaussian exponent. In the present normalization, choose
that exponent initially as $15/64$ and restrict the local interval so that
\[
15/64-(s-\tau)\ge7/32,
\qquad 0\le s-\tau\le1/64,
\]
with decrease rate one in that local construction.
The local norm and the time of existence depend on the fixed $\nu$;
no uniform lower bound as $\nu\downarrow0$ is used.
The VPL assertion is the one stated in \cite{ref4}, rather than a
claim that the Landau theorem itself contains the Poisson coupling.
For the finite-order extension of that local argument, the field obeys
\[
\|E[G]\|_{H_x^{k+1}}
\le C\|\rho[G]-1\|_{H_x^k}
\le C\bigl(1+\|e^{7|v|^2/32}G\|_{H^k_{x_1,v}}\bigr).
\]
In particular, the commutator with the electric transport has no loss
of the total differentiation order:
\[
[\partial_{x_1}^j\partial_v^\alpha,E[G]\partial_{v_1}]G
=\sum_{p=1}^j\binom jp(\partial_{x_1}^pE[G])
\partial_{x_1}^{j-p}\partial_v^{\alpha+e_1}G,
\qquad j-p+|\alpha|+1\le k.
\]
For the decreasing Gaussian weight, its additional transport term is
bounded by the negative quadratic-velocity term in the local energy
identity: for $0<\lambda\le15/64$ and any test $u$,
\[
2\lambda\left|\iint E[G]v_1|u|^2\right|
\le\frac12\|vu\|_2^2
+2\lambda^2\|E[G]\|_\infty^2\|u\|_2^2.
\]
These are the additional terms in the local VPL argument. The
finite-order construction and uniqueness give
\eqref{eq:localGaussianClass} and the stated dependence of $h$.
Applying the three-dimensional periodic result to data independent
of $x_2,x_3$, uniqueness under translations in these variables gives
a solution depending only on $x_1$.

\emph{2. Identification with the energy solution.}
The margins
\[
\frac15<\frac7{32}<\frac{15}{64}<\frac14,
\qquad \mu_{5/4}^{-1/2}=C e^{|v|^2/5}
\]
imply, for every fixed polynomial weight and $j+|\alpha|\le k$,
\[
\sup_{[\tau,\tau+h]}
\|\langle v\rangle^\ell\mu_{5/4}^{-1/2}
\partial_{x_1}^j\partial_v^\alpha G\|_2
\le C_{k,\ell}\sup_{[\tau,\tau+h]}
\|e^{7|v|^2/32}G\|_{H^k}<\infty.
\]
At $s=\tau$, multiplication by $\sqrt{\mu/\mu_{5/4}}$ and
\eqref{eq:localIdentificationCoefficients} give
\[
\mathcal E_1[(G_\tau-\mu)/\sqrt{\mu_{5/4}}]^{1/2}
\le C\epsilon.
\]
Continuity therefore makes both coefficients $F,G$ satisfy
\eqref{eq:prescribedGaussianCoefficients} at $r=0$ on a short interval
$I'=[\tau,T]$. For $G$ the integrated coefficient condition follows from
\[
\int_{I'}\mathcal D_1[(G-\mu)/\sqrt{\mu_{5/4}}]\,ds
\le C\nu|I'|\left(1+\sup_{I'}\|e^{7|v|^2/32}G\|_{H^2}^2\right)<\infty;
\]
for $F$ it follows from \eqref{eq:localIdentificationCoefficients}.
The Poisson equation and local mass conservation also give
\[
\|E[F]\|_{C_x^1}+\|\phi[F]\|_{C_x^2}
+\|\partial_s\phi[F]\|_\infty
\le C\mathcal E_k[(F-\mu)/\sqrt\mu]^{1/2}.
\]
Subtracting the two VPL equations yields the exact linear equation
\begin{equation}\label{eq:localEnergyIdentification}
\begin{aligned}
&(\partial_s+v_1\partial_{x_1}+E[F]\partial_{v_1})(G-F)\\
&\quad=\nu\bigl(Q(F,G-F)+Q(G-F,G)\bigr)
-E[G-F]\partial_{v_1}G,
\qquad (G-F)(\tau)=0.
\end{aligned}
\end{equation}
For any $U$, Gaussian integration and the one-dimensional Poisson estimate give
\[
\begin{aligned}
\|E[U]\|_\infty&\le C\|\rho[U]\|_{L_x^2}
\le C\mathcal E_0[U/\sqrt{\mu_{5/4}}]^{1/2},\\
\mathcal E_0[\mu_{5/4}^{-1/2}E[U]\partial_{v_1}G]^{1/2}
&\le C\mathcal E_0[\mu_{5/4}^{-1/2}\partial_{v_1}G]^{1/2}
\mathcal E_0[U/\sqrt{\mu_{5/4}}]^{1/2}.
\end{aligned}
\]
Thus Lemma~\ref{lem:physicalDerivativeEnergy}, with
$\varpi=5/4$, $r=0$, $F^{\rm c}_1=F$, $F^{\rm c}_2=G$
and $\mathcal J_sU=-E[U]\partial_{v_1}G$, gives
\[
\|(G-F)/\sqrt{\mu_{5/4}}\|_{\mathcal N_0^{I'}}
\le C_{\nu,I'}\mathcal E_0[(G-F)(\tau)/\sqrt{\mu_{5/4}}]^{1/2}=0.
\]
This equality cannot stop while $G$ exists: at an endpoint of an
interval of equality,
\[
\mathcal E_1[(G-\mu)/\sqrt{\mu_{5/4}}]^{1/2}
=\mathcal E_1[(F-\mu)/\sqrt{\mu_{5/4}}]^{1/2}
\le C\epsilon<\epsilon_*/4,
\]
so the same uniqueness argument applies on a further interval.
This proves the identification on the whole common interval.
\end{proof}

\begin{prop}\label{thm:nonlinearMixedRegularity}
Fix ${\rm s}\ge4$ and let $F^\nu$ be the energy solution of
Proposition~\ref{thm:exactNonlinearSolution}. Take $\nu$ sufficiently
small for that proposition at order $\max\{{\rm s},6\}$.
For $|\alpha|=q\le{\rm s}$, the derivative
$U_\alpha=\partial_v^\alpha F^\nu$ satisfies
\begin{equation}\label{eq:finiteMixedRegularity}
\mu^{-1/2}U_\alpha\in C(I;\mathscr E_{{\rm s}-q}),
\qquad \|U_\alpha\|_{\mathcal N_{{\rm s}-q,1}}<\infty.
\end{equation}
Consequently, for each fixed $\ell\ge0$,
\begin{equation}\label{eq:commonGaussianMixedRegularity}
\sum_{j+|\alpha|\le{\rm s}}
\|e^{|v|^2/8}\langle v\rangle^\ell
\partial_{x_1}^j\partial_v^\alpha F^\nu(s)\|_2
\end{equation}
is finite and continuous on $I$. If $0\le n\le{\rm s}$ and
$\ell+2n\le2{\rm s}+6$, then
\begin{equation}\label{eq:normalizedMixedRegularity}
f^\nu\in C(I;\mathscr H^n_\ell),\qquad
\sup_I\|f^\nu\|_{\mathscr H^n_\ell}
\le C_{{\rm s},n,\ell}\left(1+
\sum_{q=0}^n\sum_{|\alpha|=q}\|U_\alpha\|_{\mathcal N_{{\rm s}-q,1}}\right)<\infty.
\end{equation}
These norms use the physical variables $(x_1,v)$. The solution is
classical and nonnegative on $I$. The bounds are for each fixed
$\nu$ and need not be uniform as $\nu\downarrow0$.
\end{prop}

\begin{proof}
It suffices to work at ${\rm s}\ge6$: orders $4$ and $5$ follow by
restriction from order $6$, using the compatibility of energy solutions
in Proposition~\ref{thm:exactNonlinearSolution}. Write $F=F^\nu$ and
$f=(F-\mu)/\sqrt\mu$.

\emph{1. A finite hierarchy on the full interval.}
Set $U_0=F$. For $|\alpha|=q\ge1$, define $U_\alpha$ successively by
\begin{equation}\label{eq:physicalMixedDerivative}
(\partial_s+v_1\partial_{x_1}+E[F]\partial_{v_1})U_\alpha
=\nu\bigl(Q(F,U_\alpha)+Q(U_\alpha,F)\bigr)+\mathcal S_\alpha,
\qquad U_\alpha(t_0)=\partial_v^\alpha F_{\rm in},
\end{equation}
where
\begin{equation}\label{eq:physicalMixedSource}
\mathcal S_\alpha=-\alpha_1\partial_{x_1}U_{\alpha-e_1}
+\nu\sum_{\substack{\beta\le\alpha\\0<|\beta|<q}}
\binom\alpha\beta Q(U_\beta,U_{\alpha-\beta}).
\end{equation}
The transport source is omitted if $\alpha_1=0$.
For this finite hierarchy put
\[
\mathcal M_q=\sum_{|\alpha|=q}\|U_\alpha\|_{\mathcal N_{{\rm s}-q,1}},
\qquad
\mathcal M_{{\rm in},q}=\sum_{|\alpha|=q}
\mathcal E_{{\rm s}-q}[\mu^{-1/2}\partial_v^\alpha F_{\rm in}]^{1/2}.
\]
The explicit Gaussian datum gives $\mathcal M_{{\rm in},q}<\infty$.
The energy construction gives $\mathcal M_0<\infty$ and
\begin{equation}\label{eq:mixedOriginalCoefficients}
\sup_I\mathcal E_{\rm s}[f]^{1/2}=o(1),\qquad
\int_I\mathcal D_{\rm s}[f]\,ds<\infty.
\end{equation}
We will prove, before identifying the derivatives,
\begin{equation}\label{eq:finiteMixedRecurrence}
\mathcal M_q\le C_{\nu,{\rm s}}\left(
\mathcal M_{{\rm in},q}+(1+T_1)\mathcal M_{q-1}
+\sum_{q_1=1}^{q-1}\mathcal M_{q_1}\mathcal M_{q-q_1}\right),
\qquad 1\le q\le{\rm s}.
\end{equation}

Put $r={\rm s}-q$. The derivative and weight allocations are
\[
2(r+1)+6-2(j+1)=2r+6-2j,
\qquad 1\le|\beta|\le q-1\ \Longrightarrow\
{\rm s}-|\beta|,\ {\rm s}-q+|\beta|\ge r+1.
\]
Thus the two source integrals satisfy
\begin{equation}\label{eq:mixedSourceBound}
\begin{aligned}
&\int_I\mathcal E_r[\mu^{-1/2}\alpha_1\partial_{x_1}U_{\alpha-e_1}]^{1/2}\,ds
\le C_{\rm s}(1+T_1)\mathcal M_{q-1},\\
&\nu^{-1/2}\sum_{j=0}^r
\|\partial_{x_1}^j\{\mu^{-1/2}\nu Q(U_\beta,U_{\alpha-\beta})\}\|_{L^2_{s,x}(I)\Sigma_{2r+6-2j}^*}\\
&\quad=\sqrt\nu\sum_{j=0}^r
\|\partial_{x_1}^j\Gamma(\mu^{-1/2}U_\beta,\mu^{-1/2}U_{\alpha-\beta})\|_{L^2_{s,x}(I)\Sigma_{2r+6-2j}^*}\\
&\quad\le C_{\rm s}\sup_I\mathcal E_{r+1}[\mu^{-1/2}U_\beta]^{1/2}
\left(\int_I\mathcal D_{r+1}[\mu^{-1/2}U_{\alpha-\beta}]\,ds\right)^{1/2}\\
&\qquad+C_{\rm s}\left(\int_I\mathcal D_{r+1}[\mu^{-1/2}U_\beta]\,ds\right)^{1/2}
\sup_I\mathcal E_{r+1}[\mu^{-1/2}U_{\alpha-\beta}]^{1/2}\\
&\quad\le C_{\rm s}\|U_\beta\|_{\mathcal N_{{\rm s}-|\beta|,1}}
\|U_{\alpha-\beta}\|_{\mathcal N_{{\rm s}-q+|\beta|,1}}.
\end{aligned}
\end{equation}
For $r\ge1$ this is \eqref{eq:mixedSpatialLandauForm}. At $r=0$,
both $U_\beta$ and $U_{\alpha-\beta}$ still have one spatial derivative, and the endpoint is
\begin{equation}\label{eq:mixedLastLevelEstimate}
\begin{aligned}
&\|\Gamma(\mu^{-1/2}U_\beta,\mu^{-1/2}U_{\alpha-\beta})\|_{L_x^2\Sigma_6^*}\\
&\quad\le C\|\mu^{-1/2}U_\beta\|_{L_x^\infty L_v^2}
\|\mu^{-1/2}U_{\alpha-\beta}\|_{L_x^2\Sigma_6}
+C\|\mu^{-1/2}U_\beta\|_{L_x^\infty\Sigma_0}
\|\langle v\rangle^6\mu^{-1/2}U_{\alpha-\beta}\|_2\\
&\quad\le C\nu^{-1/2}\left(
\mathcal E_1[\mu^{-1/2}U_\beta]^{1/2}\mathcal D_1[\mu^{-1/2}U_{\alpha-\beta}]^{1/2}
+\mathcal D_1[\mu^{-1/2}U_\beta]^{1/2}\mathcal E_1[\mu^{-1/2}U_{\alpha-\beta}]^{1/2}\right).
\end{aligned}
\end{equation}
Here $H_x^1\subset L_x^\infty$ and the weights at spatial orders
$0,1$ in $\mathcal E_1$ are $8,6$.

For the coefficient hypotheses, Poisson's equation and local mass conservation give
\[
\phi_k[F]=\rho_k[F-\mu]/k^2,\qquad
E_k[F]=-i\rho_k[F-\mu]/k,\qquad
\partial_s\phi_k[F]=-ij_k[F]/k\quad(k\ne0),
\]
and hence
\[
\sup_I\bigl(\|E[F]\|_{C_x^{{\rm s}-1}}+\|\phi[F]\|_{C_x^{\rm s}}
+\|\partial_s\phi[F]\|_{C_x^{{\rm s}-1}}\bigr)
\le C_{\rm s}\sup_I\mathcal E_{\rm s}[f]^{1/2}<\infty.
\]
Since $r_{\rm c}=\max(r,1)\le{\rm s}-1$,
\[
2\sup_I\mathcal E_{r_{\rm c}}[f]^{1/2}\le\epsilon_*,
\qquad 2\int_I\mathcal D_{r_{\rm c}}[f]\,ds<\infty.
\]
These inequalities and \eqref{eq:mixedSourceBound} verify
Lemma~\ref{lem:physicalDerivativeEnergy} with
$F^{\rm c}_1=F^{\rm c}_2=F$, $\varpi=1$ and $\mathcal J=0$.
It constructs the entire hierarchy on $I$ and proves
\eqref{eq:finiteMixedRecurrence}. Every coefficient is already known
from the energy solution, and the source at level $q$ uses only
levels strictly below $q$.

\emph{2. Identification on the classical interval.}
The datum in \eqref{eq:positivefamily} is nonnegative and satisfies
\eqref{eq:localGaussianDatum} at every fixed order. Apply
Lemma~\ref{lem:localClassicalIdentification} at order ${\rm s}$.
On the resulting interval the local solution $G$ equals $F$.
Differentiation of the local equation gives
\[
\partial_v^\alpha Q(G,G)
=\sum_{\beta\le\alpha}\binom\alpha\beta
Q(\partial_v^\beta G,\partial_v^{\alpha-\beta}G),\qquad
[\partial_v^\alpha,v_1\partial_{x_1}]G
=\alpha_1\partial_{x_1}\partial_v^{\alpha-e_1}G.
\]
The first identity follows by differentiating simultaneous velocity
translations of both collision arguments. For finite Sobolev
regularity these identities hold in distributions, by testing against
smooth functions with compact velocity support.
If $\partial_v^\beta G=U_\beta$ for $|\beta|<q$, then
$\partial_v^\alpha G$ satisfies \eqref{eq:physicalMixedDerivative}
with exactly the source \eqref{eq:physicalMixedSource} and the same
initial datum as $U_\alpha$.

We check the energy class required for uniqueness, including the
top derivative. On a compact classical interval $I'=[t_0,T]$,
\eqref{eq:localGaussianClass} implies
\[
\sup_{I'}\mathcal E_r[\mu_{5/4}^{-1/2}\partial_v^\alpha G]<\infty,
\qquad r={\rm s}-q.
\]
Write $U=\partial_v^\alpha G$ in this step. We first estimate
$Q(U,F)$ without a velocity derivative of $U$, and then justify the
localized energy estimate for $Q(F,U)$.

The divergence form is
\[
Q(U,F)=\partial_{v_i}\bigl((\Phi_{ik}*U)\partial_{v_k}F
-(\partial_{v_k}\Phi_{ik}*U)F\bigr).
\]
For $0\le j\le r$, $\ell\ge0$ and a smooth compactly supported
velocity test $\varphi$, spatial differentiation and velocity
integration by parts give, at each fixed $x_1$,
\begin{align*}
&\int_{\mathbb R^3}\langle v\rangle^{2\ell}
\mu_{5/4}^{-1/2}\partial_{x_1}^jQ(U,F)\,\overline\varphi\,dv\\
&=-\sum_{p=0}^j\binom jp\int_{\mathbb R^3}
\langle v\rangle^{2\ell}\mu_{5/4}^{-1/2}
\Big[(\Phi_{ik}*\partial_{x_1}^pU)\partial_{v_k}\partial_{x_1}^{j-p}F
-(\partial_{v_k}\Phi_{ik}*\partial_{x_1}^pU)\partial_{x_1}^{j-p}F\Big]\\
&\hspace{35mm}\times
\left(\partial_{v_i}+\frac{2v_i}{5}
+\frac{2\ell v_i}{\langle v\rangle^2}\right)\overline\varphi\,dv.
\end{align*}
The three terms in the last operator come respectively from the test,
the Gaussian normalization and the polynomial weight. The collision
matrix bounds imply
\[
\left\|\langle v\rangle^{\ell-3/2}
\left(\nabla_v+\frac{2v}{5}
+\frac{2\ell v}{\langle v\rangle^2}\right)\varphi\right\|_2
\le C_\ell\|\varphi\|_{\Sigma_\ell}.
\]
For the two convolution factors, Cauchy--Schwarz near and away from
$v=v_*$, and then the Hardy--Littlewood--Sobolev inequality, give
\[
\begin{aligned}
|\Phi*V(v)|
&\le C\left(\int_{|v-v_*|\le1}|v-v_*|^{-2}\,dv_*\right)^{1/2}\|V\|_2\\
&\quad+C\left(\int_{|v-v_*|>1}|v-v_*|^{-2}
\langle v_*\rangle^{-6}\,dv_*\right)^{1/2}
\|\langle v\rangle^3V\|_2
\le C\|\langle v\rangle^3V\|_2,\\
\|\nabla\Phi*V\|_{L_v^6}&\le C\|V\|_2.
\end{aligned}
\]
Consequently, H\"older's inequality with exponents $(2,2)$ for the
first flux and $(6,3,2)$ for the second flux yields
\begin{align*}
&\left|\int\langle v\rangle^{2\ell}\mu_{5/4}^{-1/2}
\partial_{x_1}^jQ(U,F)\,\overline\varphi\,dv\right|\\
&\le C_{r,\ell}\sum_{p=0}^j
\|\langle v\rangle^3\partial_{x_1}^pU\|_{L_v^2}
\Big(\|\langle v\rangle^{\ell+3/2}\mu_{5/4}^{-1/2}
\nabla_v\partial_{x_1}^{j-p}F\|_{L_v^2}\\
&\hspace{48mm}
+\|\langle v\rangle^{\ell+3/2}\mu_{5/4}^{-1/2}
\partial_{x_1}^{j-p}F\|_{L_v^3}\Big)\|\varphi\|_{\Sigma_\ell}.
\end{align*}
Taking the supremum over $\|\varphi\|_{\Sigma_\ell}\le1$ gives the
pointwise $\Sigma_\ell^*$ estimate. To take the $L_x^2$ norm, use
\[
\begin{aligned}
p=j:&\quad
\|\langle v\rangle^3\partial_{x_1}^jU\|_{L_x^2L_v^2}
\|\langle v\rangle^{\ell+5/2}\mu_{5/4}^{-1/2}F\|_{L_x^\infty H_v^1},\\
p<j:&\quad
\|\langle v\rangle^3\partial_{x_1}^pU\|_{L_x^\infty L_v^2}
\|\langle v\rangle^{\ell+5/2}\mu_{5/4}^{-1/2}
\partial_{x_1}^{j-p}F\|_{L_x^2H_v^1},\\
&\quad
\|\langle v\rangle^3\partial_{x_1}^pU\|_{L_x^\infty L_v^2}
\le C\sum_{i=0}^1\|\langle v\rangle^3
\partial_{x_1}^{p+i}U\|_2,\qquad p+1\le j\le r.
\end{aligned}
\]
Here $H_v^1\subset L_v^3$. The extra factor $\langle v\rangle$
also bounds the term $2v/5$ when the derivative is moved from $F$
to $\mu_{5/4}^{-1/2}F$. More explicitly,
$F=\mu+\sqrt\mu f$ and
$\sqrt{\mu/\mu_{5/4}}=C e^{-|v|^2/20}$ give
\[
\begin{aligned}
\sum_{h=0}^{r+1}
\|\langle v\rangle^{2r+17/2}\mu_{5/4}^{-1/2}
\partial_{x_1}^hF\|_{L_x^2H_v^1}
&\le C_r\bigl(1+\mathcal E_{r+1}[f]^{1/2}
+\nu^{-1/2}\mathcal D_{r+1}[f]^{1/2}\bigr),\\
\sum_{p=0}^r\|\langle v\rangle^3\partial_{x_1}^pU\|_2
&\le C_r\mathcal E_r[\mu_{5/4}^{-1/2}U]^{1/2}.
\end{aligned}
\]
For the first line, the weighted gradient bound in $\Sigma_{2r+8-2h}$
controls $\nabla_v\partial_{x_1}^hf$; multiplication by
$e^{-|v|^2/20}$ absorbs the difference of polynomial weights.
Applying $H_x^1\subset L_x^\infty$ to the coefficient when $p=j$,
and setting $\ell=2r+6-2j$, therefore proves
\begin{equation}\label{eq:mixedIdentificationNonlocalSource}
\begin{aligned}
&\sum_{j=0}^r\|\partial_{x_1}^j
\{\mu_{5/4}^{-1/2}Q(U,F)\}\|_{L_x^2\Sigma_{2r+6-2j}^*}\\
&\quad\le C_r\bigl(1+\mathcal E_{r+1}[f]^{1/2}
+\nu^{-1/2}\mathcal D_{r+1}[f]^{1/2}\bigr)
\mathcal E_r[\mu_{5/4}^{-1/2}U]^{1/2}.
\end{aligned}
\end{equation}
Since $r+1\le{\rm s}$, the coefficient on the right is square
integrable in time. No velocity derivative of $U$ occurs on the right.

We next make the velocity boundary terms explicit. Choose
$0\le\psi_L(v)=\psi(v/L)\le1$, equal to one on $|v|\le L$ and zero
on $|v|\ge2L$, with $|\nabla^i\psi_L|\le C_iL^{-i}$ for $i=1,2$.
Use the potential-conjugated unknown from
\eqref{eq:prescribedPotentialConjugation}, now
$y=e^{2\phi[F]/5}\mu_{5/4}^{-1/2}U$.
The potential multiplier bounds identify its energy and form norms
with those of $\mu_{5/4}^{-1/2}U$. The local collision term is
\[
\mu_{5/4}^{-1/2}Q(F,\sqrt{\mu_{5/4}}y)
=\left(\partial_{v_i}-\frac{2v_i}{5}\right)
\left[\sigma_{ik}[F]\left(\partial_{v_k}-\frac{2v_j}{5}\right)y
-(\partial_{v_k}\Phi_{ik}*F)y\right].
\]
The coefficient bounds used in the collision form, with
$F=\mu+\sqrt\mu f$, give, for $0\le p\le r$,
\[
\begin{aligned}
\left|X^{\mathsf T}(\partial_{x_1}^p\sigma[F])Y\right|
&\le C_r\bigl(1+\mathcal E_{r+1}[f]^{1/2}\bigr)
|\sigma^{1/2}X|\,|\sigma^{1/2}Y|,\\
\left|\partial_{x_1}^p\left(\tfrac25\sigma[F]v
+\nabla_v\cdot\sigma[F]\right)\right|
&\le C_r\bigl(1+\mathcal E_{r+1}[f]^{1/2}
+\nu^{-1/2}\mathcal D_{r+1}[f]^{1/2}\bigr).
\end{aligned}
\]
The bounds are uniform in $x_1,v$; one spatial derivative supplies
$L_x^\infty$, and the Gaussian in $\sqrt\mu f$ supplies the velocity
moments. For the drift convolution, one uses
\[
\partial_{v_k}\Phi_{ik}*V=\Phi_{ik}*\partial_{v_k}V,
\qquad \|\nabla_v\Phi*V\|_\infty
\le C\|\langle v\rangle^3\nabla_vV\|_2;
\]
this accounts for the dissipation term in the second bound. Also
\[
|\sigma^{1/2}\nabla\psi_L|\le C L^{-1}\langle v\rangle^{-1/2},
\qquad
\psi_L\nabla_v\varphi=\nabla_v(\psi_L\varphi)-\varphi\nabla_v\psi_L.
\]
Thus the boundary term from the differentiated principal part obeys,
for any $\eta>0$,
\begin{align*}
&2\nu\sum_{j=0}^r\sum_{p=0}^j\binom jp
\left|\iint\langle v\rangle^{2(2r+6-2j)}\psi_L
(\partial_{x_1}^p\sigma[F])\nabla_v\partial_{x_1}^{j-p}y
\cdot\nabla_v\psi_L\,\overline{\partial_{x_1}^jy}\,dx_1dv\right|\\
&\quad\le \eta\mathcal D_r[\psi_Ly]
+\frac{C_{r,\eta}\nu}{L^2}
\bigl(1+\mathcal E_{r+1}[f]\bigr)\mathcal E_r[y].
\end{align*}
Indeed Cauchy--Schwarz pairs
$\sigma^{1/2}\nabla_v(\psi_L\partial_{x_1}^{j-p}y)$ with
$\sigma^{1/2}\nabla\psi_L\,\partial_{x_1}^jy$; the identity above
leaves a second term containing two cutoff derivatives. The lower
spatial component has the larger available weight, since
$2r+6-2(j-p)\ge2r+6-2j$.
The collision drift and the electric transport boundary terms satisfy
\begin{align*}
&\nu\sum_{j=0}^r\sum_{p=0}^j\binom jp
\left|\iint\langle v\rangle^{2(2r+6-2j)}\psi_L\nabla_v\psi_L
\cdot\partial_{x_1}^p\left(\tfrac25\sigma[F]v
+\nabla_v\cdot\sigma[F]\right)
\partial_{x_1}^{j-p}y\,\overline{\partial_{x_1}^jy}\,dx_1dv\right|\\
&\quad\le\frac{C_r\nu}{L}
\bigl(1+\mathcal E_{r+1}[f]^{1/2}
+\nu^{-1/2}\mathcal D_{r+1}[f]^{1/2}\bigr)\mathcal E_r[y],\\
&\sum_{j=0}^r\sum_{p=0}^j\binom jp
\left|\iint\langle v\rangle^{2(2r+6-2j)}
(\partial_{x_1}^pE[F])\psi_L\partial_{v_1}\psi_L
\partial_{x_1}^{j-p}y\,\overline{\partial_{x_1}^jy}\,dx_1dv\right|\\
&\quad\le\frac{C_r}{L}\|E[F]\|_{C_x^r}\mathcal E_r[y].
\end{align*}
These are the terms where a velocity derivative hits the cutoff in
the local collision and transport pairings. Terms where it hits the
polynomial weight are the ones already estimated in
\eqref{eq:prescribedGaussianEnergy}. Taking $\eta$ smaller than the
coercivity constant and integrating gives
\[
\int_{I'}|\mathrm{Err}_L|\,ds
\le\frac c4\int_{I'}\mathcal D_r[\psi_Ly]\,ds
+C_{\nu,I'}L^{-1}\sup_{I'}\mathcal E_r[y].
\]
All coefficient norms in $C_{\nu,I'}$ are finite by
\eqref{eq:mixedOriginalCoefficients}.

For the nonlocal ordering and the prescribed source, the cutoff test
is controlled by
\[
\|\psi_L^2\varphi\|_{\Sigma_\ell}
\le C_\ell\bigl(\|\psi_L\varphi\|_{\Sigma_\ell}
+L^{-1}\|\langle v\rangle^\ell\varphi\|_2\bigr).
\]
In particular, \eqref{eq:mixedIdentificationNonlocalSource},
Cauchy--Schwarz in time and Young's inequality give
\begin{align*}
&\nu\int_{I'}\left|\sum_{j=0}^r
\left\langle\partial_{x_1}^j
\{e^{2\phi[F]/5}\mu_{5/4}^{-1/2}Q(U,F)\},
\psi_L^2\partial_{x_1}^jy\right\rangle_{2r+6-2j}\right|\,ds\\
&\quad\le\frac c4\int_{I'}\mathcal D_r[\psi_Ly]\,ds
+C_{\nu,I'}\sup_{I'}\mathcal E_r[y],
\end{align*}
because
\[
\nu\int_{I'}\bigl(1+\mathcal E_{r+1}[f]^{1/2}
+\nu^{-1/2}\mathcal D_{r+1}[f]^{1/2}\bigr)^2ds
\le C\nu|I'|\bigl(1+\sup_{I'}\mathcal E_{r+1}[f]\bigr)
+C\int_{I'}\mathcal D_{r+1}[f]\,ds<\infty.
\]
The two parts of $\mathcal S_\alpha$ have the source bounds
\eqref{eq:mixedSourceBound}; the potential multiplier preserves those
bounds. Pairing them with $\psi_L^2\partial_{x_1}^jy$ gives,
uniformly in $L$,
\begin{align*}
\int_{I'}\mathcal D_r[\psi_Ly]\,ds
\le C_{\nu,I'}\Bigg(&\sup_{I'}\mathcal E_r[y]\\
&+\left[(1+T_1)\mathcal M_{q-1}
+\sum_{q_1=1}^{q-1}\mathcal M_{q_1}\mathcal M_{q-q_1}\right]^2\Bigg)<\infty.
\end{align*}
Use the regularized tests in
\eqref{eq:prescribedRegularizedEnergyIdentity} on each fixed velocity
ball, first removing the mollifier and then letting $L\to\infty$.
The estimates above bound the localized dissipation before this last
limit. Fatou's lemma therefore yields
\[
\int_{I'}\mathcal D_r[\mu_{5/4}^{-1/2}U]\,ds
\le C_{\nu,I'}\liminf_{L\to\infty}
\int_{I'}\mathcal D_r[\psi_Ly]\,ds<\infty.
\]
Thus both $\partial_v^\alpha G$ and $U_\alpha$ belong to the same
weaker Gaussian energy class. Their difference satisfies
\[
\begin{aligned}
&(\partial_s+v_1\partial_{x_1}+E[F]\partial_{v_1})
(\partial_v^\alpha G-U_\alpha)\\
&\quad=\nu\bigl(Q(F,\partial_v^\alpha G-U_\alpha)
+Q(\partial_v^\alpha G-U_\alpha,F)\bigr),
\qquad (\partial_v^\alpha G-U_\alpha)(t_0)=0.
\end{aligned}
\]
Lemma~\ref{lem:physicalDerivativeEnergy} now gives
\[
\|\mu_{5/4}^{-1/2}(\partial_v^\alpha G-U_\alpha)\|_{\mathcal N_r^{I'}}=0.
\]
Induction identifies every $U_\alpha$ with the actual derivative on
the classical interval. The original-weight bounds come from Step~1,
not from an inverse estimate between the two Gaussian weights.

\emph{3. Continuation to the observation time.}
Let $T_{\rm cl}$ be the supremum of the endpoints up to which this
classical identification holds. For $t_0\le\tau<T_{\rm cl}$,
Steps~1--2 give
\begin{equation}\label{eq:classicalContinuationBound}
\begin{aligned}
\|e^{15|v|^2/64}F(\tau)\|_{H^{\rm s}_{x_1,v}}
&\le C_{\rm s}\sum_{q=0}^{\rm s}\mathcal M_q<\infty,\\
\sum_{j+|\alpha|\le4}
\|e^{|v|^2/8}\partial_{x_1}^j\partial_v^\alpha F(\tau)\|_2
&\le C_{\rm s}\sum_{q=0}^{4}\mathcal M_q<\infty.
\end{aligned}
\end{equation}
Indeed every polynomial multiple of
$e^{15|v|^2/64}\sqrt\mu=C e^{-|v|^2/64}$ is bounded; Leibniz's
formula also controls the derivatives of the exponential weight.
The second line is the continuation norm of
\cite[Remark~3.2]{ref4}; the first also supplies a uniform bound for
restarting the local construction at the chosen order ${\rm s}$.
Both bounds are independent of $\tau$ for the fixed $\nu$.
The nonnegative local solutions therefore have a common positive
existence time $h_\nu$ when restarted at these times.
If $T_{\rm cl}<T_1$, choose
\[
\max\{t_0,T_{\rm cl}-h_\nu/2\}<\tau<T_{\rm cl}.
\]
Lemma~\ref{lem:localClassicalIdentification} identifies the restarted
solution with $F$ up to $\min\{\tau+h_\nu,T_1\}>T_{\rm cl}$.
Its differentiated equations identify the same hierarchy $U_\alpha$
by the zero-data uniqueness estimate of Step~2, now with initial
time $\tau$. This contradicts the definition of $T_{\rm cl}$.
Thus $T_{\rm cl}=T_1$; restarting shortly before $T_1$ also supplies
the classical endpoint and its one-sided time derivative.
Nonnegativity follows from the local theory on every overlapping
interval. This proves \eqref{eq:finiteMixedRegularity} on $I$ without
any additional smallness condition involving $\mathcal M_q$.

\emph{4. The normalized Sobolev norms.}
The bounded multiplier $e^{|v|^2/8}\langle v\rangle^\ell\sqrt\mu$
gives \eqref{eq:commonGaussianMixedRegularity}, including continuity.
For the normalized perturbation, Leibniz's formula is
\begin{equation}\label{eq:normalizedMixedLeibniz}
\partial_{x_1}^j\partial_v^\alpha f
=\sum_{\beta\le\alpha}\binom\alpha\beta
(\partial_v^\beta\mu^{-1/2})\partial_{x_1}^jU_{\alpha-\beta}
-\mathbf1_{j=0}\partial_v^\alpha\sqrt\mu.
\end{equation}
For $j+|\alpha|\le n$ and $\ell+2n\le2{\rm s}+6$,
\[
|\partial_v^\beta\mu^{-1/2}|
\le C_\beta\langle v\rangle^{|\beta|}\mu^{-1/2},\qquad
\ell+|\beta|\le2({\rm s}-|\alpha-\beta|)+6-2j.
\]
Each term in \eqref{eq:normalizedMixedLeibniz} is therefore bounded by
$C\mathcal E_{{\rm s}-|\alpha-\beta|}[\mu^{-1/2}U_{\alpha-\beta}]^{1/2}$.
Applying the same bounded linear maps to time differences proves
\eqref{eq:normalizedMixedRegularity}. For ${\rm s}\ge6$, Sobolev
embedding in the four variables $(x_1,v)$ gives local $C^2$ regularity,
and the equation gives
\[
\partial_sF=-v_1\partial_{x_1}F-E[F]\partial_{v_1}F
+\nu\bigl((\Phi_{ij}*F)\partial_{v_iv_j}F+8\pi F^2\bigr)\in C.
\]
The Gaussian bounds justify the convolutions and the endpoint traces.
\end{proof}
\par\endgroup

\begingroup
\subsection{Nonnegativity, conservation, and norm inflation}\label{sec:nonlinearNormInflation}

The local identification has already proved nonnegativity and classical
regularity. We now verify the conserved quantities, then use
Lemma~\ref{lem:nonlinearSobolevTrace} to convert the terminal density
bound into norm inflation. The final proof combines the linear and
nonlinear comparisons with background damping and returns to physical time.

\begin{prop}\label{prop:nonlinearPositivity}
{Assume the hypotheses of Proposition~\ref{thm:exactNonlinearSolution},
and let $F^\nu$ be the solution constructed there. Take $\nu$
sufficiently small for the construction at order $\max\{{\rm s},6\}$.
Then $F^\nu$ is
nonnegative and has the mass, total momentum, and total kinetic plus
electric energy of $\mu$ throughout $I$.}
\end{prop}

\begin{proof}

Proposition~\ref{thm:nonlinearMixedRegularity}, applied at order
$\max\{{\rm s},6\}$, identifies $F^\nu$ on overlapping intervals
with nonnegative classical solutions $G$:
\[
F^\nu(s)=G(s)\ge0,\qquad
s\in[\tau,\tau+h_\nu]\cap I.
\]
These intervals cover $I$, so no separate negative-part construction
is required. For conservation, the collision invariants give

\begin{equation}\label{eq:conservationCollisionWeakForm}
\int_{\mathbb R^3}Q(F^\nu,F^\nu)
\begin{pmatrix}1\\v\\|v|^2/2\end{pmatrix}\,dv=0.
\end{equation}
This follows by symmetrizing the collision integral: the derivative
differences for the three tests are $0,0,v-v_*$, and
$\Phi(v-v_*)(v-v_*)=0$. Gaussian integrability and the collision
form bound justify removal of the velocity cutoffs. Write
$\rho=\int F^\nu\,dv$, $j=\int v_1F^\nu\,dv$ and $E=E[F^\nu]$.
Then
\begin{equation}\label{eq:conservationContinuityField}
\partial_s\rho+\partial_{x_1}j=0,\qquad
\partial_{x_1}E=\rho-1,\qquad
\partial_sE=-j+\frac1{2\pi}\int_0^{2\pi}j\,dx_1.
\end{equation}
Periodicity and the zero mean of $E$ imply
\begin{equation}\label{eq:conservationMassMomentum}
\frac d{ds}\iint F^\nu\,dx_1dv=0,\qquad
\frac d{ds}\iint v_iF^\nu\,dx_1dv
=\delta_{i1}\int E(1+\partial_{x_1}E)\,dx_1=0.
\end{equation}
The kinetic and electric energy derivatives cancel:
\begin{equation}\label{eq:conservationEnergyCancellation}
\frac d{ds}\iint\frac{|v|^2}{2}F^\nu\,dx_1dv=\int Ej\,dx_1,
\qquad
\frac d{ds}\frac12\int|E|^2\,dx_1=-\int Ej\,dx_1.
\end{equation}
Continuity of these moments and of the field energy extends the
integrated identities to both endpoints.

\end{proof}
\par\endgroup

We now convert the terminal density lower bound into a lower bound for a finite Sobolev norm of the exact solution. We use the physical free-transport coordinates based at the initial time. We then strengthen the gain condition to make the terminal norm itself tend to infinity, {rather than merely the ratio} to the initial size.

\begin{lem}\label{lem:nonlinearSobolevTrace}
Let $n\ge0$ be an integer and $\ell>3/2$. With the Fourier convention
\eqref{eq:fourier}, every $H\in\mathscr H^n_\ell$ satisfies
\begin{equation}\label{eq:nonlinearSobolevTrace}
\langle k,\eta\rangle^n|\widehat H_k(\eta)|
\le C_{n,\ell}\|H\|_{\mathscr H^n_\ell},
\qquad k\in\mathbb Z,\quad\eta\in\mathbb R^3.
\end{equation}
Multiplication by {$\sqrt \mu $} satisfies
\begin{equation}\label{eq:profileGaussianMultiplier}
\|\sqrt \mu \,g\|_{\mathscr H^n_\ell}
\le C_{n,\ell}\|g\|_{\mathscr H^n_\ell}.
\end{equation}
For the exact solution and the normalized profiles in
\eqref{eq:shiftedProfileCorrespondence}, for $0\le t\le t_\nu$,
\begin{equation}\label{eq:nonlinearPhysicalTrace}
{
\rho_k[\widetilde F^\nu-\mu ](t) =\rho_k[F^\nu-\mu ](t_0+t) =\widehat{\sqrt \mu \,g^\nu}_k(t_0+t,kte_1) =\widehat{\sqrt \mu \,\widetilde g^\nu}_k(t,kte_1).
}
\end{equation}
If Proposition~\ref{thm:nonlinearMixedRegularity} holds at order
${\rm s}\ge n$, then
{$\sqrt \mu \,g^\nu\in C(I;\mathscr H^n_\ell)$} and
{$\sqrt \mu \,\widetilde g^\nu\in C([0,t_\nu];\mathscr H^n_\ell)$}.
If also $\ell+2n\le2{\rm s}+6$, then
\begin{equation}\label{eq:normalizedProfileContinuity}
g^\nu\in C(I;\mathscr H^n_\ell),\qquad
\widetilde g^\nu\in C([0,t_\nu];\mathscr H^n_\ell).
\end{equation}
The constants in \eqref{eq:nonlinearSobolevTrace} and
\eqref{eq:profileGaussianMultiplier} are independent of time and $\nu$.
The continuity assertions are for each fixed $\nu$.
\end{lem}

\begin{proof}
For $j+|\alpha|\le n$, Fourier transformation and
Cauchy--Schwarz yield
\[
|k|^j|\eta^\alpha|\,|\widehat H_k(\eta)|
\le\|\langle v\rangle^{-\ell}\|_2
\|\langle v\rangle^\ell(\partial_z^j\partial_v^\alpha H)_k\|_2
\le C_\ell\|H\|_{\mathscr H^n_\ell}.
\]
Summing the monomials proves \eqref{eq:nonlinearSobolevTrace}.
Leibniz's rule and boundedness of every fixed derivative of {$\sqrt \mu $}
give \eqref{eq:profileGaussianMultiplier}.

The transport change of variables and
{$F^\nu-\mu =\sqrt \mu \,f^\nu$} yield
\[
\widehat{\sqrt \mu \,g^\nu}_k(s,\eta)
=\widehat{(F^\nu-\mu )}_k\big(s,\eta-k(s-t_0)e_1\big).
\]
Taking $s=t_0+t$ and $\eta=kte_1$ proves
\eqref{eq:nonlinearPhysicalTrace}.
For the normalized profile the chain rule is
\[
\begin{aligned}
\partial_z^j\partial_v^\alpha g^\nu(s,z,v)
={}&\sum_{m=0}^{\alpha_1}\binom{\alpha_1}{m}(s-t_0)^m\,
\big(\partial_{x_1}^{j+m}\partial_v^{\alpha-me_1}f^\nu\big)
(s,z+(s-t_0)v_1,v),\\
\|g^\nu(s)\|_{\mathscr H^n_\ell}
\le{}& C_n(1+{T_1}-t_0)^n\|f^\nu(s)\|_{\mathscr H^n_\ell}.
\end{aligned}
\]
Under $\ell+2n\le2{\rm s}+6$, equation
\eqref{eq:normalizedMixedRegularity} and strong continuity of the
transport change of variables prove \eqref{eq:normalizedProfileContinuity}.
The latter continuity follows first for smooth compactly supported
functions, and then by density and the displayed uniform bound on $I$.
Applying the same chain rule to {$F^\nu-\mu $}, using
\eqref{eq:commonGaussianMixedRegularity}, proves the unnormalized
continuity for every fixed $\ell$ when ${\rm s}\ge n$.
Time translation gives both physical-time assertions.
\end{proof}

We now prove the logarithmic-amplitude result stated in the introduction. The parameter choice leaves a strictly positive power margin for the nonlinear comparison. We then pass from the density of the perturbation about the exact background to that of the full solution about {$\mu $}.

\begingroup
\begin{proof}[Proof of Theorem~\ref{thm:logNonlinearInflation}]
Fix the mode parameters and choose $K$ once, independently of the
Sobolev order and of $\nu$, so that
\begin{equation}\label{eq:logAbsoluteGainChoice}
Kc_{\rm e}>2.
\end{equation}
Since $C_*$, $C_{\rm B}$ and $c_{\rm e}$ do not depend on the Sobolev order or the fixed logarithmic power $p$, there is a finite integer ${\rm s}_{\log}\ge11$ such that
\[
{\rm s}_{\log}>2+3K(C_*+2C_{\rm B}-c_{\rm e}).
\]
For every fixed ${\rm s}_0\ge {\rm s}_{\log}$ this guarantees \eqref{eq:nonlinearRoadmapMargin}. Fix $p>{\rm s}_0+7/2$, use \eqref{eq:comparisonLogParameters}, and put $\delta_\nu=c_1\nu^{1/3}[\log(1/\nu)]^p\zeta^{-{\rm s}_0}$, with fixed small $c_0,c_1>0$. Choose $\gamma_{\rm a}>1/(3K)$ as specified after \eqref{eq:comparisonLogParameters}, and fix $A_0\ge1$, as allowed by Proposition~\ref{thm:bgcomparison}.

Choose an additional fixed integer ${\rm s}$ such that
\[
{\rm s}\ge {\rm s}_0+2,\qquad \ell+2{\rm s}_0\le2{\rm s}+6.
\]
Apply Proposition~\ref{thm:exactNonlinearSolution} at this order. It
constructs the exact solution for the datum $F_{\rm in}$
and identifies it with the order-four comparison solution. This higher
order is used only for qualitative regularity in the prescribed
weight: it changes neither ${\rm s}_0$, $p$, nor the terminal comparison;
it only decreases the admissible upper bound on $\nu$. Proposition~\ref{thm:nonlinearMixedRegularity} gives classical regularity and finite continuous mixed norms, and Proposition~\ref{prop:nonlinearPositivity} gives nonnegativity and conservation. The initial mass, momentum and energy are exactly those of {$\mu $} by Lemma~\ref{lem:data}. These steps concern the exact nonlinear equation, including both collision orderings {and the self-consistent electric field.}

Set $t_\nu={T_1}-t_0$ and return to physical time by
\[
\widetilde F^\nu(t,x,v)=F^\nu(t_0+t,x_1,v),
\qquad 0\le t\le t_\nu.
\]
The time translation in Section~\ref{sec:mainIdeas} shows that this is
a nonnegative classical solution of \eqref{eq:vplPhysicalTime}, with
the same conserved quantities. {The initial datum, transport profile}
and density satisfy the exact identities
\[
\begin{aligned}
\widetilde F^\nu(0)&=F_{\rm in},\qquad \widetilde f^\nu(t)=f^\nu(t_0+t),\qquad
\widetilde g^\nu(t)=g^\nu(t_0+t),\\
\rho_k[\widetilde F^\nu-\mu ](t)&=\rho_k[F^\nu-\mu ](t_0+t).
\end{aligned}
\]

We verify the size of this datum in the requested Gaussian-weighted norm. Bounds \eqref{eq:maininitialnorm} and \eqref{eq:initialTemperatureMixedBound} give
\begin{equation}\label{eq:logNonlinearInitialUpper}
\begin{split}
\left\|(F_{\rm in}-\mu )/\sqrt \mu \right\|_{H^{{\rm s}_0}_{q,\ell}}
\le{}&C\left[a(1+t_0)^{{\rm s}_0}+\delta_\nu\zeta^{{\rm s}_0}\right]\\
&+C(1+t_0)^{{\rm s}_0}\delta_\nu^2m_0^{-2}e^{-(\zeta-m_0t_0)^2/2}
\le C\nu^{1/3}[\log(1/\nu)]^p.
\end{split}
\end{equation}
For the matching lower bound use spatial mode $1$. The high-frequency perturbation has modes $\pm m_0$, with $m_0>1$, so it cannot cancel this mode. {The normalized coefficient in mode $1$ is exactly}
{$a\mu_{\theta_{\delta_\nu}}e^{-it_0v_1}/(2\sqrt \mu )$}.
The proof of Lemma~\ref{lem:initialGaussianBounds} applies uniformly to $\theta_{\delta_\nu}\to1$ and yields a lower bound $ca t_0^{{\rm s}_0}\ge c'\nu^{1/3}[\log(1/\nu)]^p$. This proves the norm assertion in \eqref{eq:logNonlinearInitialTime} with the same ${\rm s}_0,q,\ell$ on both sides. Also
\[
{T_1}\mathrel{{\sim}}\nu^{-1/3}[\log(1/\nu)]^{{\rm s}_0+2-p},\qquad t_0=O(\log(1/\nu))=o({T_1}).
\]
Thus $t_\nu={T_1}-t_0$ has the stated size and is $o(\nu^{-1/3})$.

Proposition~\ref{thm:exactNonlinearSolution} gives
$|\rho_{m_1}[F^\nu-F^{\rm b,\nu}]({T_1})|\ge c\delta_\nu B_*/\sqrt N$.
The background contribution is negligible by Proposition~\ref{prop:bguniform} and the field decay in \eqref{eq:bguniformfield}:
\[
|\rho_{m_1}[F^{\rm b,\nu}-\mu ]({T_1})|
\le Cm_1\varepsilon_Ne^{-c({T_1}-t_0)^\chi}
=o(\delta_\nu B_*/\sqrt N).
\]
Indeed ${T_1}-t_0=\nu^{-1/3+o(1)}$, whereas $\delta_\nu B_*/\sqrt N$ has a fixed power of $\nu$ up to subpower factors. The triangle inequality therefore gives
\begin{equation}\label{eq:logFullNonlinearDensity}
|\rho_{m_1}[F^\nu-\mu ]({T_1})|\ge c\delta_\nu B_*/\sqrt N.
\end{equation}

Set $k_\nu=m_1$. Since $m_1({T_1}-t_0)\ge\zeta/2$ for small $\nu$,
\begin{equation}\label{eq:logAbsoluteNonlinearGrowth}
\begin{split}
\langle m_1,m_1({T_1}-t_0)e_1\rangle^{{\rm s}_0-3}
|\rho_{m_1}[F^\nu-\mu ]({T_1})|
&\ge c\nu^{1/3}[\log(1/\nu)]^p\zeta^{-3}B_*/\sqrt N\\
&\ge c\nu^{4/3-Kc_{\rm e}}
[\log(1/\nu)]^{4p-3{\rm s}_0-19/2}
\ge c\nu^{-2/3}.
\end{split}
\end{equation}
Here the construction scales and the descending-chain gain give
\[
\zeta\mathrel{{\sim}}\nu^{-1/3}[\log(1/\nu)]^{{\rm s}_0+3-p},
\qquad B_*\mathrel{{\sim}}\nu^{-Kc_{\rm e}},\qquad
N\mathrel{{\sim}}\log(1/\nu).
\]
Since $p>{\rm s}_0+7/2$, the logarithmic exponent satisfies
$4p-3{\rm s}_0-19/2>{\rm s}_0+9/2>0$.
The last inequality follows from \eqref{eq:logAbsoluteGainChoice}.
Since ${T_1}=t_0+t_\nu$ and
{$\rho_{k_\nu}[\widetilde F^\nu-\mu ](t_\nu)=\rho_{k_\nu}[F^\nu-\mu ]({T_1})$},
this proves \eqref{eq:logNonlinearDensity}.

For the stronger profile conclusion, apply Lemma~\ref{lem:nonlinearSobolevTrace} at order ${\rm s}_0-4$ to the same density lower bound. It gives
\begin{equation}\label{eq:logTransportNormGrowth}
\begin{aligned}
\|\sqrt \mu \,g^\nu({T_1})\|_{\mathscr H^{{\rm s}_0-4}_\ell}
&\ge c\langle m_1,m_1({T_1}-t_0)e_1\rangle^{{\rm s}_0-4}
|\rho_{m_1}[F^\nu-\mu ]({T_1})|\\
&\ge c\nu^{1/3}[\log(1/\nu)]^p\zeta^{-4}B_*/\sqrt N\\
&\ge c\nu^{5/3-Kc_{\rm e}}
[\log(1/\nu)]^{5p-4{\rm s}_0-25/2}
\ge c\nu^{-1/3}.
\end{aligned}
\end{equation}
Indeed $5p-4{\rm s}_0-25/2>{\rm s}_0+5>0$ and
$Kc_{\rm e}>2$. The exact profile identity
\[
\sqrt \mu (v)\,g^\nu({T_1},z,v)
=\widetilde F^\nu(t_\nu,z+t_\nu v_1,v)-\mu (v)
\]
therefore proves \eqref{eq:logNonlinearInflation}. The multiplier estimate
\eqref{eq:profileGaussianMultiplier} then gives
\[
\|\widetilde g^\nu(t_\nu)\|_{\mathscr H^{{\rm s}_0-4}_\ell}
\ge C^{-1}\|\sqrt \mu \,\widetilde g^\nu(t_\nu)\|_{\mathscr H^{{\rm s}_0-4}_\ell}
\ge c\nu^{-1/3}.
\]
For each fixed $\nu$, the choice $\ell+2{\rm s}_0\le2{\rm s}+6$
allows us to apply \eqref{eq:normalizedMixedRegularity} and
\eqref{eq:normalizedProfileContinuity}, giving
\[
\widetilde f^\nu\in C([0,t_\nu];\mathscr H^{{\rm s}_0}_\ell),\qquad
\widetilde g^\nu\in C([0,t_\nu];\mathscr H^{{\rm s}_0}_\ell).
\]
Multiplication by {$\sqrt \mu $} gives the asserted continuity of the
unnormalized profile. No uniform-in-$\nu$ continuity bound is claimed.
At $t=0$,
\[
\widetilde g^\nu(0)=\widetilde f^\nu(0)
=\frac{F_{\rm in}-\mu }{\sqrt \mu },\qquad
\|\widetilde g^\nu(0)\|_{\mathscr H^{{\rm s}_0}_\ell}
\le\|\widetilde f^\nu(0)\|_{H^{{\rm s}_0}_{q,\ell}}
\longrightarrow0.
\]
The multiplier estimate gives the same initial smallness for
{$\sqrt \mu \,\widetilde g^\nu$}. This completes the proof for the original
initial-value problem \eqref{eq:vplPhysicalTime}.
\end{proof}
\par\endgroup
\appendix

\section{Auxiliary normalization and coefficient estimates}\label{app:coefficients}

This appendix proves the auxiliary estimates invoked in the background and energy
arguments. We begin with the initial norms and fixed normalization, then establish
Gaussian coefficient bounds and differentiated collision forms. All constants may
depend on the stated fixed orders; the dependence on $N,\nu$ is recorded explicitly.

\subsection{Initial norms and a fixed normalization}\label{app:initialNormalization}

\begingroup
We prove the initial-data properties of Lemma~\ref{lem:data} and the scale calculation of Lemma~\ref{lem:logarithmicScales}. The following normalization lemma then checks the convention and spatial period used in the cited background theorem.

\begin{proof}[Proof of Lemma~\ref{lem:data}]\begingroup
\emph{1. Temperature selection, positivity and the initial remainder.}
Spatial averaging gives mass one and momentum zero for
$F^{\rm b}_{\rm in}$. Equation \eqref{eq:initialFullEnergy} at $\delta=0$
{shows that the total energy of $F^{\rm b}_{\rm in}$ equals} {$\mathcal E(\mu )$} exactly when
\eqref{eq:bgtemp} holds.
The left-hand side of \eqref{eq:bgtemp} changes sign between
$1-a^2/6$ and $1$, and
\[
3-\frac{a^2t_0^2}{2}e^{-\theta t_0^2}
\ge3-\frac{a^2}{2e\theta}>2
\qquad(1-a^2/6\le\theta\le1).
\]
This proves existence and uniqueness. The equation then gives
\[
{
(1-\theta)t_0^2 =\frac{a^2}{6}t_0^2e^{-\theta t_0^2}\le Ca^2,\qquad 1-\theta =\frac{a^2}{6}e^{-t_0^2}e^{(1-\theta)t_0^2} \le Ca^2e^{-t_0^2}.
}
\]
The resulting datum has the required conserved quantities, and
{$F^{\rm b}_{\rm in}-F^{\rm b}_{\rm app}(t_0)
=(\mu_\theta-\mu )(1+a\cos(x_1-t_0v_1))$}.

For the full datum, the two field modes are orthogonal since $m_0\ge2$:
\begin{align*}
E[F^\delta_{\rm in}]
&=a e^{-\theta_\delta t_0^2/2}\sin x_1
+\frac{\delta}{\sqrt\pi m_0}e^{-(\zeta-m_0t_0)^2/4}\sin(m_0x_1),\\
2\bigl(\mathcal E(F^\delta_{\rm in})-\mathcal E(\mu )\bigr)
&=3(\theta_\delta-1)+\frac{a^2}{2}e^{-\theta_\delta t_0^2}
+\frac{\delta^2}{2\pi m_0^2}e^{-(\zeta-m_0t_0)^2/2}.
\end{align*}
Consequently, energy matching is precisely \eqref{eq:thetadelta}.
The preceding derivative bound and the implicit function theorem
{give a unique root of \eqref{eq:thetadelta} near} $\theta$ for small $|\delta|$, and the
mean value theorem gives
\[
|\theta_\delta-\theta|
\le C\delta^2m_0^{-2}e^{-(\zeta-m_0t_0)^2/2}.
\]
For small $|\delta|$ this root lies in $[3/4,1]$.
The oscillatory terms have zero spatial mean, so the mass and total
momentum remain one and zero, respectively.

For $3/4\le\theta_\delta\le1$,
\[
\frac{\mu_{1/2}(v)}{\mu_{\theta_\delta}(v)}
=(2\theta_\delta)^{3/2}
\exp\left(-\left(1-\frac1{2\theta_\delta}\right)|v|^2\right)
\le2\sqrt2.
\]
Since $a\le1/2$, sufficiently small $|\delta|$ gives
\[
F^\delta_{\rm in}\ge\mu_{\theta_\delta}
\left(1-a-\frac{2\sqrt2}{\sqrt\pi}|\delta|\right)
\ge\frac14\mu_{\theta_\delta}.
\]
The temperature is smooth and even in $\delta$; hence
\[
\left.\partial_\delta\theta_\delta\right|_{\delta=0}=0,
\qquad
\left.\partial_\delta F^\delta_{\rm in}\right|_{\delta=0}
=\frac1{\sqrt\pi}\mu_{1/2}
\cos\big(m_0x_1+(\zeta-m_0t_0)v_1\big)=H_{\rm in}.
\]
\par

To prove the correction bounds, use
\[
\mu_{\theta_\delta}-\mu_\theta
=(\theta_\delta-\theta)\int_0^1
\partial_{\theta'} \mu_{\theta'}\big|_{\theta'=\theta+u(\theta_\delta-\theta)}\,du.
\]
For each fixed multi-index $\beta$, Gaussian differentiation gives
\[
\sup_{3/4\le\theta'\le1}
\left\|e^{q|v|^2}\langle v\rangle^\ell
\partial_v^\beta\left(\frac{\partial_{\theta'}\mu_{\theta'}}{\sqrt\mu}\right)\right\|_2
\le C_{\beta,q,\ell},\qquad 0<q<1/4.
\]
The mean-value formula and the product rule therefore yield
\[
\begin{aligned}
\partial_{x_1}^j\partial_v^\alpha
\left(\frac{W^\delta_{\rm in}}{\sqrt\mu}\right)
&=\sum_{\beta\le\alpha}\binom\alpha\beta
\partial_v^\beta\left(\frac{\mu_{\theta_\delta}-\mu_\theta}{\sqrt\mu}\right)
\partial_{x_1}^j\partial_v^{\alpha-\beta}
\big[1+a\cos(x_1-t_0v_1)\big],\\
\left\|e^{q|v|^2}\langle v\rangle^\ell
\partial_{x_1}^j\partial_v^\alpha
\left(\frac{W^\delta_{\rm in}}{\sqrt\mu}\right)\right\|_2
&\le C_{j,\alpha,q,\ell}|\theta_\delta-\theta|
\big[1+a(1+t_0)^{|\alpha|}\big].
\end{aligned}
\]
Summing over $j+|\alpha|\le n$ and inserting the temperature bound gives
\[
\left\|\frac{W^\delta_{\rm in}}{\sqrt\mu}\right\|_{H^n_{q,\ell}}
\le C_{n,q,\ell}\delta^2m_0^{-2}e^{-(\zeta-m_0t_0)^2/2}
\big[1+a(1+t_0)^n\big].
\]
Taking $\alpha=0$ in the same product estimate gives
\[
\sum_{j=0}^n\left\|\langle v\rangle^\ell\partial_{x_1}^j
\left(\frac{W^\delta_{\rm in}}{\sqrt\mu}\right)\right\|_2
\le C_{n,\ell}|\theta_\delta-\theta|
\le C_{n,\ell}\delta^2m_0^{-2}e^{-(\zeta-m_0t_0)^2/2}.
\]
\par
\par\endgroup
\emph{2. Sobolev norms.}
Put {${\mathfrak m_\theta}=\mu_\theta/\sqrt \mu $}. Choose a compact interval $J_q$ about $1$ such that
\[
\inf_{\theta\in J_q}\left(\frac1{2\theta}-\frac14-q\right)>c_q>0.
\]
Part (i) places $\theta$ in this interval for large $N$.
Gaussian differentiation and the mean-value theorem give, uniformly
for $|\alpha|\le{\rm s}_0$ and $\theta\in J_q$,
\[
\begin{aligned}
&\|e^{q|v|^2}\langle v\rangle^\ell\partial_v^\alpha {\mathfrak m_\theta}\|_2
+\|e^{q|v|^2}\langle v\rangle^\ell
\partial_\theta\partial_v^\alpha {\mathfrak m_\theta}\|_2\le C_{{\rm s}_0,q,\ell},\\
&\|e^{q|v|^2}\langle v\rangle^\ell
\partial_v^\alpha({\mathfrak m_\theta}-{\mathfrak m_1})\|_2
\le C_{{\rm s}_0,q,\ell}|\theta-1|.
\end{aligned}
\]
For $j+|\alpha|\le {\rm s}_0$,
\[
\partial_{x_1}^j\partial_v^\alpha({\mathfrak m_\theta} e^{ix_1-it_0v_1})
=i^je^{ix_1-it_0v_1}
\sum_{j_1=0}^{\alpha_1}\binom{\alpha_1}{j_1}(-it_0)^{\alpha_1-j_1}
\partial_1^{j_1}\partial_{v_\perp}^{\alpha_\perp}{\mathfrak m_\theta}.
\]
Applying this identity to the two conjugate modes, and to the high-frequency datum with
Gaussian {$\mu_{1/2}/\sqrt \mu $} and frequencies $m_0,\zeta-m_0t_0$, yields
\begin{align*}
\|(F^{\rm b}_{\rm in}-\mu )/\sqrt \mu \|_{H^{{\rm s}_0}_{q,\ell}}
&\le C\big(|\theta-1|+a(1+t_0)^{{\rm s}_0}\big)
\le Ca(1+t_0)^{{\rm s}_0},\\
\|H_{\rm in}/\sqrt \mu \|_{H^{{\rm s}_0}_{q,\ell}}
&\le C(1+m_0+|\zeta-m_0t_0|)^{{\rm s}_0}\le C\zeta^{{\rm s}_0}.
\end{align*}
Here $|\theta-1|\le Ca^2e^{-t_0^2}$, and
$a(1+t_0)^{{\rm s}_0}\le CN^{-{\rm s}_0-2}(1+A_0N)^{{\rm s}_0}\le CN^{-2}$.

For the lower bound, uniform Gaussian bounds give fixed $c,C>0$ with
\begin{align*}
0<c&\le\inf_{\theta\in J_q}\|e^{q|v|^2}\langle v\rangle^\ell {\mathfrak m_\theta}\|_2,\\
\sum_{j_1=1}^{{\rm s}_0}\binom{{\rm s}_0}{j_1}
\sup_{\theta\in J_q}\|e^{q|v|^2}\langle v\rangle^\ell\partial_1^{j_1}{\mathfrak m_\theta}\|_2&\le C.
\end{align*}
The $k=1$ mode of the normalized background is
$(a/2){\mathfrak m_\theta} e^{-it_0v_1}$. For $t_0\ge\max\{1,2C/c\}$,
\begin{align*}
&\left\|e^{q|v|^2}\langle v\rangle^\ell\partial_1^{{\rm s}_0}
\left(\frac a2{\mathfrak m_\theta} e^{-it_0v_1}\right)\right\|_2\\
&\quad\ge\frac a2t_0^{{\rm s}_0}\|e^{q|v|^2}\langle v\rangle^\ell {\mathfrak m_\theta}\|_2\\
&\qquad-\frac a2\sum_{j_1=1}^{{\rm s}_0}\binom{{\rm s}_0}{j_1}t_0^{{\rm s}_0-j_1}
\|e^{q|v|^2}\langle v\rangle^\ell\partial_1^{j_1}{\mathfrak m_\theta}\|_2\\
&\quad\ge\frac a2(ct_0^{{\rm s}_0}-Ct_0^{{\rm s}_0-1})
\ge\frac{ac}4t_0^{{\rm s}_0}.
\end{align*}
Parseval gives the required lower bound.

\emph{3. Gevrey norms and comparison with the cited convention.}
Put {${\mathfrak m_\theta}=\mu_\theta/\sqrt \mu $} and choose an even integer $L_1\ge L$.
Uniformly for $\theta$ near $1$, {the Fourier transform and any fixed}
{number of $\theta$ and $\eta$ derivatives} have Gaussian decay.
The only spatial modes are
\begin{equation}\label{eq:initialBackgroundFourierModes}
\widehat f^{\rm b}_{{\rm in},0}=\widehat {\mathfrak m_\theta}-\widehat {\mathfrak m_1},
\qquad
\widehat f^{\rm b}_{{\rm in},\pm1}(\eta)
=\frac a2\widehat {\mathfrak m_\theta}(\eta\pm t_0e_1).
\end{equation}
For $|\beta|\le L_1$, $j+|\alpha|\le{\rm s}$ and $k=\pm1$,
Leibniz's rule and subadditivity of the Gevrey exponent give
\[
\left\|\partial_\eta^\beta\left[(ik)^j(i\eta)^\alpha
e^{\lambda\langle k,\eta\rangle^\chi}
\widehat {\mathfrak m_\theta}(\eta\pm t_0e_1)\right]\right\|_2
\le C(1+t_0)^De^{C(1+t_0)^\chi}
\]
for some fixed $D$. The zero mode is bounded by $C|\theta-1|$ by
the mean-value theorem. Plancherel and \eqref{eq:bgtemp} therefore yield
\[
\|f^{\rm b}_{\rm in}\|_{\mathscr G^{\lambda,\chi}_{{\rm s},L}}
\le C|\theta-1|+Ca(1+t_0)^De^{C(1+t_0)^\chi}
\le aC(1+N)^{r_L}e^{CN^\chi}.
\]

For use here and in Proposition~\ref{prop:bguniform}, the standard
Stirling and maximum estimates are, for $0<\lambda'<\lambda$ and
fixed $D_0$,
\begin{equation}\label{eq:gevreyRadiusComparison}
\begin{gathered}
c e^{c(\lambda|\xi|)^\chi}
\le\left(\sum_{n\ge0}\frac{\lambda^{2n}|\xi|^{2n}}{(n!)^{2/\chi}}\right)^{1/2}
\le C e^{C(\lambda|\xi|)^\chi},\\
\sup_\xi|\xi|^n\langle\xi\rangle^{D_0}
e^{-(\lambda-\lambda')\langle\xi\rangle^\chi}
\le C_{\lambda,\lambda',D_0}^{n+1}(n!)^{1/\chi}.
\end{gathered}
\end{equation}
The first follows from Stirling's inequalities, and the second by
maximizing $x^{n+D_0}e^{-(\lambda-\lambda')x^\chi}$.
At time zero the six reference vector fields are ordinary derivatives.
In \cite[(1.12), (2.6)]{ref6}, the multinomial coefficients cost at most
$6^n$, there are at most $(n+1)^5$ multi-indices, and the remaining
coefficient is $(n+1)^{12}/(n!)^{1/\chi}$. Thus, after the fixed
coordinate changes of Lemma~\ref{lem:fixedPeriodNormalization},
\[
{
\big\|(2\pi)^{3/4}f^{\rm b}_{\rm in}(c_xx'_1,\sqrt2v') \big\|_{\mathcal G_{c_x;\chi,L}^{\lambda_{\rm der},{\rm s}}}^2 \le C\|f^{\rm b}_{\rm in}\|_{\mathscr G^{\lambda,\chi}_{{\rm s}',L_1}}^2 \sum_{n\ge0}(n+1)^{29}(C\lambda_{\rm der})^{2n} \le Ca^2L_N^2,
}
\]
where $c_x=\sqrt2\pi^{3/4}$, the larger indices ${\rm s}',L_1$ are fixed,
and $\lambda_{\rm der}>0$ is sufficiently small but independent of $\nu$.
Polynomial velocity moments follow by Fourier differentiation,
{absorbing the fixed polynomial factors} in the radius gap in
\eqref{eq:gevreyRadiusComparison}.

Finally, $\chi<1$ and the parameter choice give
\[
\log L_N=O(\log N+N^\chi)=o(N),\qquad
\log(aL_N)=-\tfrac13\log(1/\nu)+o(\log(1/\nu)),
\]
which proves \eqref{eq:epsNsmall}.
\end{proof}

\begin{proof}[Proof of Lemma~\ref{lem:logarithmicScales}]
Since $N=K\log(1/\nu)$, $t_0=A_0N$ and $m_1\mathrel{{\sim}} N$,
\[
\begin{aligned}
a&\mathrel{{\sim}}\nu^{1/3}[\log(1/\nu)]^{p-{\rm s}_0},\qquad \zeta=N^3/a\mathrel{{\sim}}\nu^{-1/3}[\log(1/\nu)]^{{\rm s}_0+3-p},\\
{T_1}=\zeta/m_1&\mathrel{{\sim}}\nu^{-1/3}[\log(1/\nu)]^{{\rm s}_0+2-p}.
\end{aligned}
\]
The initial restrictions follow from
\[
{
ae^{\gamma_{\rm a}N} \mathrel{{\sim}}\nu^{1/3-K\gamma_{\rm a}}[\log(1/\nu)]^{p-{\rm s}_0}\to\infty,\qquad aN^{{\rm s}_0+2}\mathrel{{\sim}}\nu^{1/3}[\log(1/\nu)]^{p+2}\to0.
}
\]
Also
\[
\begin{aligned}
\log {T_1}&=\tfrac13\log(1/\nu)+O(\log\log(1/\nu)),\\
{T_1}/t_0&\mathrel{{\sim}}\nu^{-1/3}[\log(1/\nu)]^{{\rm s}_0+1-p}\to\infty,
\qquad \zeta\mathrel{{\sim}} N{T_1}\to\infty.
\end{aligned}
\]
Consequently
\[
{
\nu {T_1}^3\mathrel{{\sim}}[\log(1/\nu)]^{3{\rm s}_0+6-3p},\qquad \nu {T_1}(1+\zeta^2)\mathrel{{\sim}}\nu N^2{T_1}^3 \mathrel{{\sim}}[\log(1/\nu)]^{3{\rm s}_0+8-3p}.
}
\]
Multiplication by $\sqrt N$ and $N^{9/2}$ proves \eqref{eq:logObjectiveCollision}. Since $\log(1/\nu)\to\infty$, the three limits are equivalent, respectively, to
\[
{\rm s}_0+2-p<0,\qquad 3{\rm s}_0+8-3p<0,\qquad
\max\{3{\rm s}_0+17/2-3p,3{\rm s}_0+21/2-3p\}<0.
\]
This proves \eqref{eq:preparationTimeRange}--\eqref{eq:retainedComparisonRange}.
{Substitution of $p={\rm s}_0+7/2+\eta$ into the second line of
\eqref{eq:logObjectiveCollision} gives \eqref{eq:introTerminalCollisionEstimates}.} Finally,
\[
\sup_{0\le t\le {T_1}}|\xi-kt|\le C(\zeta+N{T_1})\le C\zeta
\quad\Longrightarrow\quad
\nu\int_0^{T_1}(1+|\xi-kt|^2)\,dt\le C\nu {T_1}(1+\zeta^2).
\]
\end{proof}
\par\endgroup
\begin{lem}\label{lem:fixedPeriodNormalization}
Set $c_x=\sqrt2\pi^{3/4}$ and $c_t=\pi^{3/4}$. For the physical-time
unknowns in \eqref{eq:vplPhysicalTime}, the change of variables
\begin{equation}\label{eq:fullMaxwellianscaling}
\begin{aligned}
t&=c_tt',\qquad x=c_xx',\qquad v=\sqrt2v',\qquad F'(t',x',v')=(2\pi)^{3/2}\widetilde F(t,x,v),\\
E'(t',x')&=\frac{c_t}{\sqrt2}\widetilde E(t,x),\qquad
\nu'=\frac{\nu}{2^{3/2}\pi^{3/4}}
\end{aligned}
\end{equation}
maps \eqref{eq:vplPhysicalTime}, with equilibrium {$\mu $} on the torus of side
$2\pi$, to the convention with equilibrium {$\mu'(v')=e^{-|v'|^2}$} and
background density $\pi^{3/2}$ on the torus of side $2\pi/c_x$.
For the shifted solution of \eqref{eq:vpl}, the same transformation reads
\[
{
F'(t',x',v')=(2\pi)^{3/2}F(t_0+c_tt',c_xx'_1,\sqrt2v'),\qquad E'(t',x')=\frac{c_t}{\sqrt2}E[F](t_0+c_tt',c_xx'_1)e_1.
}
\]
Thus $t'=0$ corresponds to $t=0$ and $s=t_0$.
Both the Gevrey estimates in \cite[Theorem 2.1]{ref6} and the finite-Sobolev
estimates in \cite[Theorem 3.1]{ref4} admit this fixed period, with constants
and smallness bounds depending on $c_x$ but not on $\nu'$. Their powers of
$\nu'$ and their derivative counts are unchanged.
\end{lem}

\begingroup
\begin{proof}
Use $v=\sqrt2v'$ and $\widetilde F=(2\pi)^{-3/2}F'$,
suppressing the corresponding time and spatial arguments.
The velocity measure, Coulomb kernel and two derivatives contribute
$(\sqrt2)^3$, $(\sqrt2)^{-1}$ and $(\sqrt2)^{-2}$, whose product is one.
The two distribution factors therefore give
\[
Q(\widetilde F,\widetilde F)(\sqrt2v')=(2\pi)^{-3}Q(F',F')(v').
\]
The transformed kinetic equation has transport coefficient
$\sqrt2c_t/c_x=1$, field coefficient $c_t/\sqrt2$, and collision
coefficient $\nu c_t(2\pi)^{-3/2}=\nu'$. Writing
$\widetilde\rho=\int\widetilde F\,dv$ and $\rho'=\int F'\,dv'$, we obtain
\[
\widetilde\rho=\pi^{-3/2}\rho',\qquad
\operatorname{div}_{x'}E'
=\frac{c_x c_t}{\sqrt2}(\widetilde\rho-1)=\rho'-\pi^{3/2}.
\]
The field remains a mean-zero gradient. This proves the transformed Poisson equation and the asserted normalization. Spatial averages use the volume of the corresponding torus; hence equality of the conserved quantities with those of the equilibrium is preserved. The hypotheses transform as follows:
\[
\begin{array}{c|c}
\text{current convention}&\text{reference convention after \eqref{eq:fullMaxwellianscaling}}\\ \hline
\mu ,\quad \int \mu \,dv=1&\mu'=e^{-|v'|^2},\quad \int\mu'\,dv'=\pi^{3/2}\\
k\in\mathbb Z^3&\xi=c_xk\in c_x\mathbb Z^3\\
\nu&\nu'=\nu/(2^{3/2}\pi^{3/4})\\
t\le {T_1}-t_0=o(\nu^{-1/2})&t'\le({T_1}-t_0)/c_t=o((\nu')^{-1/2}).
\end{array}
\]
Polynomial weights and fixed derivative orders are equivalent under this fixed scaling. A fixed Gevrey radius becomes another fixed positive radius; the Gaussian datum in Lemma~\ref{lem:bginitialsize} permits that choice before $\nu$ is made small. The remaining issue is therefore the fixed spatial lattice, which we check next.

We specify why the period change is permitted in the cited estimates. For a fixed side $2\pi/c_x$ the normalized Fourier coefficients are indexed by $\xi\in c_x\mathbb Z^3$. Products still give sums with the exact relation $\xi=\eta+(\xi-\eta)$, Parseval's identity has only the fixed volume factor, and
\begin{equation}\label{eq:fixedLatticeBounds}
\begin{gathered}
|\xi|\ge c_x>1\quad(\xi\ne0),\qquad
\sum_{\xi\in c_x\mathbb Z^3}\langle\xi\rangle^{-p}<\infty\quad(p>3),\\
\langle\xi\rangle^\chi\le
\langle\eta\rangle^\chi+\langle\xi-\eta\rangle^\chi,
\qquad 0<\chi<1.
\end{gathered}
\end{equation}

Use normalized spatial measure on the fixed lattice
$\Lambda_{c_x}=c_x\mathbb Z^3$. Parseval's identity and the product
convolution retain their usual form. By \eqref{eq:fixedLatticeBounds}
and $\ell^2_{\rm s}(\Lambda_{c_x})\subset\ell^1(\Lambda_{c_x})$,
for ${\rm s}>3/2$,
\begin{equation}\label{eq:fixedLatticeOperators}
\begin{aligned}
\|f\|_{L^\infty_x}&\le C_{c_x,{\rm s}}\|f\|_{H^{\rm s}_x},\\
\|e^{\lambda\langle\xi\rangle^\chi}\langle\xi\rangle^{\rm s}(f*g)_\xi\|_{\ell^2}
&\le C_{c_x,{\rm s}}
\|e^{\lambda\langle\xi\rangle^\chi}\langle\xi\rangle^{\rm s}f_\xi\|_{\ell^2}
\|e^{\lambda\langle\xi\rangle^\chi}\langle\xi\rangle^{\rm s}g_\xi\|_{\ell^2},\\
\|\nabla\Delta^{-1}{\rho}\|_{H^{{\rm s}+1}_x}
&\le C_{c_x}\|{\rho}\|_{H^{\rm s}_x},\qquad {\rho}_0=0.
\end{aligned}
\end{equation}
The commuting field $\partial_{v_j}+it\xi_j$ and the identities
\[
\int_0^\infty\langle\xi t\rangle^{-m}\,dt=C_m/|\xi|\le C_m/c_x
\quad(m>1),\qquad
\nu'\int_0^t|\xi\tau|^2\,d\tau=\nu'|\xi|^2t^3/3
\]
show that these changes affect fixed constants, not the collision
powers or derivative counts.
\par

We also check the Maxwellian resolvent, since the cited theorems are stated on the integer lattice. Write $K_\xi^{(\nu')}$ for the scalar density kernel in the convention of \cite{ref4}, and use $\operatorname{Lap}[K](z)=\int_0^\infty e^{-zt}K(t)\,dt$. The proof of \cite[equation (7.13)]{ref4} uses the mode equation and the physical magnitude $|\xi|$, and gives
\[
\sup_{\operatorname{Re}z\ge0}
|\operatorname{Lap}[K_\xi^{(\nu')}](z)|\le C|\xi|^{-2}.
\]
Thus the Penrose denominator is bounded below by $1/2$ for $|\xi|\ge R$, with a fixed $R>c_x$. On the finite set
$\{\xi\in c_x\mathbb Z^3:0<|\xi|<R\}$, the collisionless Maxwellian Penrose condition gives
\[
c_{\rm Pen}(c_x):=\min_{\substack{\xi\in c_x\mathbb Z^3\\0<|\xi|<R}}
\inf_{\operatorname{Re}z\ge0}
|1+\operatorname{Lap}[K_\xi^{(0)}](z)|>0.
\]

The collisionless Gaussian Penrose argument in
\cite[(7.17)]{ref4} applies to each nonzero physical frequency.
Let $\Lambda_{c_x}^R=\{\xi\in\Lambda_{c_x}:0<|\xi|<R\}$.
The uniform tail bound in \cite[Lemma 7.5]{ref4} and the fixed-time
comparison in \cite[Proposition 6.5]{ref4}, applied to this finite set,
allow first ${T}$ and then $\nu_{c_x}$ to be chosen so that
\begin{align*}
\sup_{\xi\in\Lambda_{c_x}^R,\ 0\le\nu'\le\nu_{c_x}}
\int_{T}^\infty|K_\xi^{(\nu')}(t)|\,dt&\le c_{\rm Pen}(c_x)/8,\\
\sup_{\xi\in\Lambda_{c_x}^R}
\int_0^{T}|K_\xi^{(\nu')}(t)-K_\xi^{(0)}(t)|\,dt&\le c_{\rm Pen}(c_x)/4.
\end{align*}
Consequently,
\[
\sup_{\xi\in\Lambda_{c_x}^R}\|K_\xi^{(\nu')}-K_\xi^{(0)}\|_{L^1(0,\infty)}
\le c_{\rm Pen}(c_x)/2,
\]
and for $\operatorname{Re}z\ge0$,
\[
|1+\operatorname{Lap}[K_\xi^{(\nu')}](z)|
\ge |1+\operatorname{Lap}[K_\xi^{(0)}](z)|
-\|K_\xi^{(\nu')}-K_\xi^{(0)}\|_1\ge c_{\rm Pen}(c_x)/2.
\]
Combining this with the large-frequency bound yields
\[
\inf_{\xi\in c_x\mathbb Z^3\setminus\{0\}}
\inf_{\operatorname{Re}z\ge0}
|1+\operatorname{Lap}[K_\xi^{(\nu')}](z)|
\ge \min\{1/2,c_{\rm Pen}(c_x)/2\}>0.
\]
The derivative estimates for the Laplace transform and the inversion argument in \cite[Section 7]{ref4} now apply on this lattice. The same lower bound verifies the Penrose condition in \cite[Section 6.4, equation (6.30)]{ref6}; {the proof of Proposition 6.4 in Section 6.5 of that reference} then gives the required Gevrey resolvent estimate. This is a fixed-lattice adaptation of those proofs, with period-dependent constants.

We next check the frequency-dependent estimates used for the nonlinear density
estimates in \cite[Section 7.2]{ref6}, {retaining the time and mode}
decomposition.
Use physical frequencies $\xi=c_xk$, $\eta=c_xl$ in $\Lambda_{c_x}$ and a
continuous velocity frequency $\eta_{\rm v}\in\mathbb R^3$. Their phase and
Poisson factors are
\begin{equation}\label{eq:fixedLatticeEchoFactors}
\begin{aligned}
\xi t-\eta\tau
&=\big(\eta_{\rm v}+(\xi-\eta)\tau\big)
+\big(\xi(t-\tau)-\eta_{\rm v}\big),\\
\langle\xi t-\eta\tau\rangle^2
&\le2\langle\eta_{\rm v}+(\xi-\eta)\tau\rangle^2
\langle\xi(t-\tau)-\eta_{\rm v}\rangle^2,\\
\frac{\eta}{|\eta|^2}&=c_x^{-1}\frac{l}{|l|^2},
\qquad |\xi-\eta|\ge c_x\quad(k\ne l).
\end{aligned}
\end{equation}
For nonzero $\xi,\eta$, write $e=\eta/|\eta|$ and
$\xi t=(\xi t\cdot e)e+(\xi t)_\perp$. Extension of the integral to
the real line gives the two time-Schur bounds
\begin{equation}\label{eq:fixedLatticeEchoSchur}
{
\sup_{\xi,\eta,t}\int_0^\infty \frac{|\eta|\,d\tau}{1+|\xi t-\eta\tau|^2} \le\int_{\mathbb R}\frac{du}{1+u^2}=\pi,\qquad \sup_{\xi,\eta,\tau}\int_0^\infty \frac{|\xi|\,dt}{1+|\xi t-\eta\tau|^2}\le\pi.
}
\end{equation}
Thus neither the Poisson gain nor the resonant time integral incurs
a loss depending on $\nu'$ or on the length of the time interval.

Let $q_{\rm ref}$ be the fixed time-weight exponent in
\cite[Proposition 2.2]{ref6}. {Choose the radius parameters in that proposition} so that
\[
\lambda(t)=\lambda_\infty+\widetilde\delta(1+t)^{-\alpha},
\qquad
0<\alpha<\min\{\chi/10,\chi-1/2,\chi/(q_{\rm ref}+4)\}.
\]
These choices satisfy \cite[(2.36), (7.12)]{ref6}. For $t\ge1$,
split the time integral at $2t/3$, and split the late part into
$\xi=\eta$ and $\xi\ne\eta$.

On $0\le\tau\le2t/3$, the mean value formula gives
\begin{align*}
\frac{\lambda(\tau)-\lambda(t)}{\lambda(\tau)}
&=\frac{\alpha\widetilde\delta}{\lambda(\tau)}
\int_\tau^t(1+u)^{-1-\alpha}\,du
\ge c\,t(1+t)^{-1-\alpha}\ge c'\langle t\rangle^{-\alpha},\\
\left(\frac{\lambda(t)}{\lambda(\tau)}\right)^n
&\le e^{-c'n\langle t\rangle^{-\alpha}}
\le C_\gamma n^{-\gamma}\langle t\rangle^{\alpha\gamma},
\qquad n\ge1,\ \gamma>0.
\end{align*}
Thus the derivative allocation in \cite[(7.18)--(7.20)]{ref6}
applies with $k,l$ replaced by $\xi,\eta$; {the inverse frequency}
powers are bounded because $|\xi|,|\eta|\ge c_x$. This estimate uses
time separation and does not require a lower bound on
$|\xi t-\eta\tau|$. The terms of derivative order zero are retained
as in that reference estimate.

For the late off-diagonal part, put
\[
\Omega_{c_x}=\{t\ge1,\ 2t/3\le\tau\le t,\ \xi\ne\eta,
\ |\xi t-\eta\tau|\le c_xt/2\}.
\]
On $\Omega_{c_x}$,
\begin{align*}
|\eta|(t-\tau)
&=|(\eta-\xi)t+(\xi t-\eta\tau)|\ge c_xt/2,\\
\frac{\lambda(\tau)-\lambda(t)}{\lambda(\tau)}
&\ge c(t-\tau)\langle t\rangle^{-1-\alpha}
\ge c_{c_x}\langle\tau\rangle^{-\alpha}|\eta|^{-1},\\
\frac{\lambda(t)}{\lambda(\tau)}
&\le\exp\{-c_{c_x}\langle\tau\rangle^{-\alpha}|\eta|^{-1}\}.
\end{align*}
Our choice of $\alpha$ also satisfies
$\alpha\le\min\{(3\chi-1)/2,\chi/2\}$. Hence
\cite[(7.23)--(7.24)]{ref6} permits a fixed $\Theta_{\rm ref}$ with
\[
\frac1{\chi-\alpha}\le\Theta_{\rm ref}
\le\min\{\alpha^{-1},3/(1-\alpha)\}.
\]
For an integer derivative order $n\ge1$, use
$e^{-x}\le C_\Theta x^{-\Theta}$ to obtain
\begin{equation}\label{eq:fixedLatticeRadiusGain}
{
\left(\frac{\lambda(t)}{\lambda(\tau)}\right)^n \frac{|\eta\tau|}{|\eta|^3}\mathbf1_{\Omega_{c_x}} \le\frac{C_{c_x}\langle\tau\rangle^{\alpha\Theta_{\rm ref}} |\eta|^{\Theta_{\rm ref}}}{n^{\Theta_{\rm ref}}} \frac{|\eta\tau|}{|\eta|^3} \le\frac{C_{c_x}\langle\eta\tau\rangle^{1+\alpha\Theta_{\rm ref}}} {n^{\Theta_{\rm ref}}}.
}
\end{equation}
The last step uses $\tau\ge2/3$, $|\eta|\ge c_x$ and
$\Theta_{\rm ref}(1-\alpha)\le3$. At $n=0$ the factor
$|\eta\tau|/|\eta|^3$ is retained, exactly as in the reference
weight $w_1$. On the remaining late off-diagonal region,
$|\xi t-\eta\tau|>c_xt/2$ and
\[
\frac{t}{\langle\xi t-\eta\tau\rangle^3}
\le\frac{2/c_x}{\langle\xi t-\eta\tau\rangle^2}.
\]
Together with \eqref{eq:fixedLatticeEchoSchur}, these give the kernel bounds needed for the two off-diagonal estimates \cite[(7.29)--(7.30)]{ref6}.

For the late diagonal term $\xi=\eta$, the distribution factor is
{the spatial mean $f_0$.} In this mode,
\[
(\nabla_v+\tau\nabla_x)f_0=\nabla_v f_0,\qquad
\xi t-\eta\tau=\xi(t-\tau),\qquad
\int_0^t\frac{|\xi|\,d\tau}
 {1+|\xi|^2(t-\tau)^2}\le\frac{\pi}{2}.
\]
Thus the extra velocity derivative produces no factor $\tau$.
The same allocation of derivatives to $f_0$ as in
\cite[(7.31)]{ref6} applies, using the displayed time integral.
This covers the diagonal contribution without a radius-gap argument.
For $t<1$, use the fixed product bounds \eqref{eq:fixedLatticeOperators}.

These estimates cover all time and mode regions in
\cite[(7.20), (7.29)--(7.31)]{ref6}. The remaining derivative
allocations and density/current estimates {in Sections 7--8 of that reference}
use the same phase identities, time-Schur bounds and convolutions:
the index relations remain $(k,l,k-l)$, all factorial weights
remain unchanged, and the physical frequencies are $(c_xk,c_xl,c_x(k-l))$.
Equations \eqref{eq:fixedLatticeOperators} and
\eqref{eq:fixedLatticeEchoFactors}--\eqref{eq:fixedLatticeRadiusGain}
therefore transfer those estimates with fixed period-dependent
constants. The finite-Sobolev estimates in \cite[Sections 8 and 11]{ref4}
use the same phase bounds and time-Schur estimates with polynomial weights.

For precision, the estimates obtained from \cite[Theorem 2.1]{ref6}
are the following. Let $\mathcal G_{c_x;\chi,L}^{\lambda,N_0}$ denote
{the norm (2.6) in that reference} on $\mathbb T^3_{2\pi/c_x}$, formed with
$Z=(\nabla_{x'},\nabla_{v'}+t'\nabla_{x'})$.
Put {$f'=(F'-\mu')/\sqrt{\mu'}$}, {let $f'_0$ be the normalized spatial}
mean and set $f'_{\ne}=f'-f'_0$. Fix admissible derivative and moment
orders $N_0,L_0$ and an initial radius $\lambda_{\rm in}>0$.
For initial data obeying the reference constraints (1.11) and
$\|f'(0)\|_{\mathcal G_{c_x;\chi,L_0}^{\lambda_{\rm in},N_0}}
\le\epsilon$ sufficiently small, the fixed-lattice estimates read
\begin{equation}\label{eq:fixedLatticeReferenceOutputs}
\begin{aligned}
\|f'_0(t')\|_{\mathcal G_{c_x;\chi,L_0}^{\lambda_\infty,N_0}}
+\langle(\nu')^{1/3}t'\rangle^{L_0/6}
\|f'_{\ne}(t')\|_{\mathcal G_{c_x;\chi,L_0/4}^{\lambda_\infty,N_0}}
&\le C_{c_x}\epsilon,\\
\|e^{c_{c_x}\lambda_\infty|t'\nabla_{x'}|^\chi}E'(t')\|_2
&\le C_{c_x}\epsilon,
\qquad 0\le t'\le(\nu')^{-1/2}.
\end{aligned}
\end{equation}
Equivalently, the nonzero-mode norm acts by ordinary derivatives on
$f'_{\ne}(t',z'+t'v',v')$. The velocity collision estimates are
unchanged by the spatial period, and the energy proofs use only
spatial integration by parts, the commuting fields $Z$, and the
product bounds already checked. The continuity arguments in
\cite[Section 2.3]{ref6} and \cite[Sections 10 and 13]{ref4} can
therefore choose their smallness constants after the fixed constants
$C_{c_x}$ and their admissible radii. This retains the powers of $\nu'$
and the stated derivative and moment counts. The global continuation
in \cite[Theorem 1.1]{ref6} uses the same finite-Sobolev estimates.
Together with \eqref{eq:fullMaxwellianscaling}, these are exactly the
background bounds of Proposition~\ref{prop:bguniform}.

\end{proof}
\par\endgroup
\begingroup
\subsection{The loss in the coarse density estimate}\label{app:coarseLossVerification}
We compare the norm-conversion cost of the coarse inverse with the descending-chain growth.
\begin{proof}[Proof of Lemma~\ref{lem:coarseLossVersusEcho}]
We first check the cost of the contraction used in that proof.
Take $m=\lfloor\sqrt N\rfloor$ and the descending column indexed by $m$,
$k=m-1$. The explicit kernel \eqref{eq:densityexact} gives
\[
\begin{aligned}
\int_{T_{\rm a}}^{T_1}|P_{m-1,m}(T_1,\tau)|\,d\tau
&=\frac{a(m-1)}{2m^3}
\int_{-T_1}^{(m-1)T_1-mT_{\rm a}}(T_1+z)e^{-z^2/2}\,dz\\
&=(\Lambda+o(1))\frac{aT_1}{m^2}.
\end{aligned}
\]
Here $T_{\rm a}/T_1\to0$, $m\to\infty$ and
$aT_1/m^2\to\infty$. The kernel has one sign, so the constant function ${\rho}_m=1$ attains this absolute time integral. Also,
\[
\|K_{m-1}^\nu\|_{L^1(0,T_1)}=O(m^{-2}),\qquad
\|({\mathrm{Id}}-K_{m-1}^\nu*)^{-1}-{\mathrm{Id}}\|_{L^\infty\to L^\infty}=O(m^{-2}).
\]
The second identity follows by the Neumann series at this large mode;
it does not replace the uniform scalar inverse used earlier.
In particular,
\[
\begin{aligned}
&\|[(\mathrm{Id}-K_{m-1}^\nu*)^{-1}-\mathrm{Id}]
\int_{T_{\rm a}}^{\,\cdot}P_{m-1,m}(\cdot,\tau)\,d\tau\|_{L^\infty(I_+)}\\
&\qquad\le Cm^{-2}\sup_{s\in I_+}
\left|\int_{T_{\rm a}}^sP_{m-1,m}(s,\tau)\,d\tau\right|
\le C\frac{aT_1}{m^4}=o\left(\frac{aT_1}{m^2}\right).
\end{aligned}
\]
\par
Taking $q=0$ in the weighted norm, the corresponding
column norm after scalar inversion is therefore bounded below by
\[
\frac{(\Lambda+o(1))aT_1/m^2}
{1+C_{\rm maj}a(1+T_1)/m^2}
=\frac{\Lambda}{C_{\rm maj}}+o(1).
\]
The contraction bound $1/4$ imposed in the proof of Proposition~\ref{prop:coarseDensityInverse} thus
requires $C_{\rm maj}\ge4\Lambda$.

Next, the product defining $G_1$ satisfies
\[
\begin{aligned}
\log G_1
&=\sum_{m=2}^{\infty}\log\left(1+
\frac{C_{\rm maj}a(1+T_1)}{m^2}\right)\\
&=\int_0^\infty\log\left(1+
\frac{C_{\rm maj}a(1+T_1)}{y^2}\right)\,dy+O(\log N)\\
&=\pi\sqrt{C_{\rm maj}a(1+T_1)}+O(\log N),\qquad
\frac{\log G_1}{N}\longrightarrow\pi\sqrt{C_{\rm maj}/{\bf b}_1}.
\end{aligned}
\]
The integral follows by scaling and integration by parts; monotonicity
of the summand bounds the sum--integral error by $O(\log N)$.
We used $aT_1=N^3/m_1=N^2/{\bf b}_1+O(N)$.
For $\rho_k(s)=\mathbf1_{k=1}$,
\[
\|\rho\|_{\ell^1_0L^\infty(I_+)}=1,\qquad
\|\rho\|_{G,0}=G_1^{-1},\qquad
\frac{\|\rho\|_{\ell^1_0L^\infty(I_+)}}{\|\rho\|_{G,0}}=G_1.
\]
Thus the exponent delivered by weight removal is at least
$\lim_{N\to\infty}N^{-1}\log G_1$; polynomial weight conversion cannot
reduce this bound.
\par
Finally, the chain gain and the elementary maximum of
$\sqrt x(1-x)$ on $[0,1]$ give
\[
\begin{aligned}
c_{\rm e}
&=3\left({\bf b}_0-{\bf b}_1
+{\bf b}_1\log\frac{{\bf b}_1}{{\bf b}_0}\right)
<3({\bf b}_0-{\bf b}_1),\\
\sqrt{\frac{{\bf b}_1}{{\bf b}_0}}
\left(1-\frac{{\bf b}_1}{{\bf b}_0}\right)
&\le\max_{0\le x\le1}\sqrt x(1-x)=\frac{2}{3\sqrt3},\qquad \frac{C_*}{c_{\rm e}}
>\frac{2\pi\sqrt{\Lambda/{\bf b}_1}}
{3({\bf b}_0-{\bf b}_1)}
\ge\pi\sqrt3,
\end{aligned}
\]
where $\Lambda={\bf b}_0^3$. If $Kc_{\rm e}\ge1/3$, then
\[
\sigma_{\rm int}-\tfrac12-K(C_*+C_{\rm B}-c_{\rm e})
<\tfrac12-KC_*<\tfrac12-\frac{\pi}{\sqrt3}<0.
\]
This proves both assertions. The obstruction is a loss of the coarse
argument, and is compatible with sharper estimates for particular
oscillatory functions.
\end{proof}
\par\endgroup

\subsection{Gaussian coefficient bounds and velocity moments}\label{app:gaussianCoefficients}

We establish the cancellation and convolution estimates used in the amplitude-sensitive
comparison. The subsequent interpolation argument supplies polynomial moments at
a reduced Gevrey radius, and then converts the background profile bounds into
coefficient estimates for the physical-time derivative $\widehat Z_{k,c}$. The collision norm is
\eqref{eq:backgroundCollisionNorm}, and the profiles are defined in
\eqref{eq:physicalCollisionProfiles}.

All coefficient estimates use the convention that $\sigma[\cdot]$ and $b[\cdot]$ act on physical distributions in
\eqref{eq:normalizedBilinearNotation}, with {$\sigma=\sigma[\mu ]$}.

\begin{lem}
\label{lem:LandauCancellation}\label{lem:gaussianCollisionCoefficients}
For Schwartz functions {${F_1},{F_2}$} on $\mathbb R^3$,
\begin{equation}\label{eq:LandauCancellation}
Q({F_1},{F_2})=(\Phi_{ij}*{F_1})\partial_{v_iv_j}{F_2}+8\pi {F_1} {F_2}.
\end{equation}
Fix $c>0$. Suppose that the possibly complex function {${F_1}$} satisfies
{$|{F_1}(v)|\le\varepsilon_{\rm c} e^{-c|v|^2}$}, and put
\[
\sigma_{ij}[{F_1}]=\Phi_{ij}*{F_1},\qquad
b_{i}[{F_1}]=(\partial_j\Phi_{ij})*{F_1}.
\]
For $X,Y\in\mathbb C^3$ and $v\in\mathbb R^3$,
\begin{equation}\label{eq:gaussianCollisionCoefficients}
\begin{aligned}
|\overline X_i\sigma_{ij}[{F_1}](v)Y_j|
&\le C_c\varepsilon_{\rm c}(\overline X_i\sigma_{ij}(v)X_j)^{1/2}
(\overline Y_i\sigma_{ij}(v)Y_j)^{1/2},\\
|\sigma(v)^{-1/2}b[{F_1}](v)|
&\le C_c\varepsilon_{\rm c}\langle v\rangle^{-1/2}.
\end{aligned}
\end{equation}
In addition, let $\ell_{\rm ker}>3$, $C_{\rm ker}\ge0$, and let $K$ be a measurable kernel such that
\begin{equation}\label{eq:weightedCoulombEnvelope}
|K(v,v_*)|\le C_{\rm ker}\langle v\rangle^{-\ell_{\rm ker}}|v-v_*|^{-1}e^{-c|v_*|^2}.
\end{equation}
Then, for every $f\in L^2(\mathbb R^3)$,
\begin{equation}\label{eq:weightedCoulombSchur}
\begin{aligned}
\sup_v\int|K(v,v_*)|\,dv_*
+\sup_{v_*}\int|K(v,v_*)|\,dv
&\le C_{c,\ell_{\rm ker}}C_{\rm ker},\\
\left\|\int K(\cdot,v_*)f(v_*)\,dv_*\right\|_2
&\le C_{c,\ell_{\rm ker}}C_{\rm ker}\|f\|_2.
\end{aligned}
\end{equation}
\end{lem}

\begin{proof}
First,
\[
\partial_i\sigma_{ij}[{F_1}]=b_j[{F_1}],\qquad
\partial_i b_i[{F_1}]=-8\pi {F_1},
\]
so expansion of
{$\partial_i(\sigma_{ij}[{F_1}]\partial_j{F_2}-b_i[{F_1}]{F_2})$} proves
\eqref{eq:LandauCancellation}. Kernel regularization justifies the
identities distributionally.

Put ${f_c}(v)=e^{-c|v|^2}$. For $r=|v|\ge2$, $e=v/r$ and
$P_e={\mathrm{Id}}-e\otimes e$, radial symmetry diagonalizes both
$\sigma[{f_c}]$ and $\sigma$ in the radial and tangential directions.
The estimates
\[
\begin{gathered}
e\cdot\Phi(v-v_*)e=\frac{|P_ev_*|^2}{|v-v_*|^3},\qquad
\|\Phi(v-v_*)\|\le |v-v_*|^{-1},\\
\int_{|v_*|>r/2}|v-v_*|^{-j}{f_c}(v_*)\,dv_*
\le C_c e^{-cr^2/8}\quad(j=1,2)
\end{gathered}
\]
give radial and tangential upper bounds $C_cr^{-3}$ and $C_cr^{-1}$.
For the matching lower bounds on $\sigma$, integrate over
$|v_*|\le1/2$: the transverse Maxwellian second moment is positive,
and $|v-v_*|\mathrel{{\sim}} r$. Consequently,
\[
\sigma[{f_c}](v)\le C_c(r^{-3}e\otimes e+r^{-1}P_e),\qquad
\sigma(v)\ge c(r^{-3}e\otimes e+r^{-1}P_e).
\]
On $|v|\le2$, positivity of $\sigma$ and boundedness of $\sigma[{f_c}]$
give the same comparison. Hence $\sigma[{f_c}]\le C_c\sigma$ globally.
Since $\Phi$ is nonnegative, matrix Cauchy--Schwarz gives the mixed
direction estimate even when {${F_1}$} is complex:
\[
\begin{aligned}
|\overline X \sigma[{F_1}] Y|
&\le\varepsilon_{\rm c}\int
(\overline X\Phi(v-v_*)X)^{1/2}
(\overline Y\Phi(v-v_*)Y)^{1/2}{f_c}(v_*)\,dv_*\\
&\le\varepsilon_{\rm c}(\overline X \sigma[{f_c}]X)^{1/2}(\overline Y \sigma[{f_c}]Y)^{1/2}
\le C_c\varepsilon_{\rm c}(\overline X\sigma X)^{1/2}(\overline Y\sigma Y)^{1/2}.
\end{aligned}
\]
Finally $\partial_j\Phi_{ij}(z)=-2z_i/|z|^3$. The same split,
now with $j=2$, gives
\[
|b[{F_1}](v)|\le C_c\varepsilon_{\rm c}\langle v\rangle^{-2},\qquad
\|\sigma^{-1/2}(v)\|\le C\langle v\rangle^{3/2}.
\]
Multiplication proves the second bound in
\eqref{eq:gaussianCollisionCoefficients}.

For the kernel assertion, split $|v-v_*|\le1$ {from the complementary region:}
\begin{align*}
\sup_v\int|K(v,v_*)|\,dv_*
&\le C_{\rm ker}\left(\int_{|z|\le1}|z|^{-1}\,dz
+\int e^{-c|v_*|^2}\,dv_*\right)\le C_c C_{\rm ker},\\
\sup_{v_*}\int|K(v,v_*)|\,dv
&\le C_{\rm ker}\left(\int_{|z|\le1}|z|^{-1}\,dz
+\int\langle v\rangle^{-\ell_{\rm ker}}\,dv\right)
\le C_{\ell_{\rm ker}}C_{\rm ker}.
\end{align*}
Schur's inequality now gives the operator bound in
\eqref{eq:weightedCoulombSchur}.

\end{proof}

\begin{lem}\label{lem:weightedGevreyInterpolation}
Fix integers $L\ge0$ and $\ell_{\rm mom}\ge2L+2$, a radius
$\lambda>0$, and $0<\chi<1$. Suppose that a function $f(z,v)$ satisfies
\[
\|e^{\lambda\langle\nabla_z,\nabla_v\rangle^\chi}f\|_2\le C_{\rm amp}\varepsilon_{\rm int},
\qquad
\sum_{|\alpha|\le \ell_{\rm mom}}\|v^\alpha f\|_2\le C_{\rm amp},
\qquad 0<\varepsilon_{\rm int}\le1,
\]
for some $C_{\rm amp}\ge0$. There are constants $C,c,\lambda'>0$,
depending only on $L,\ell_{\rm mom},\lambda,\chi$ and independent of
$f,C_{\rm amp},\varepsilon_{\rm int}$, such that, for every
$l\in\mathbb Z$ and $|\alpha|\le L$,
\begin{equation}\label{eq:weightedGevreyInterpolation}
\|e^{\lambda'\langle\nabla_v\rangle^\chi}(v^\alpha f_l)\|_2
\le C C_{\rm amp}\varepsilon_{\rm int}^{1/2}e^{-c\langle l\rangle^\chi}.
\end{equation}
After a further fixed radius reduction, this also controls \eqref{eq:velocityInputGevrey} {for $f_l$, $\partial_{v_1}f_l$, and $v_1f_l$}, when $L\ge\ell+2$.
\end{lem}

\begin{proof}
Cover velocity frequency space by unit cubes $Q_n$, $n\in\mathbb Z^3$, with fixed larger cubes $Q_n'$ of bounded overlap. The first assumption gives
\[
\|\widehat f_l\|_{L^2(Q_n')}\le C C_{\rm amp}\varepsilon_{\rm int} e^{-c\langle l,n\rangle^\chi}.
\]
The moment assumption bounds $\|\widehat f_l\|_{H^{\ell_{\rm mom}}(Q_n')}$ by $C C_{\rm amp}$.
Local Sobolev interpolation, with $q=|\alpha|\le L$, yields
\[
\|\partial_\eta^\alpha \widehat f_l\|_{L^2(Q_n)}
\le C\|\widehat f_l\|_{L^2(Q_n')}^{1-q/\ell_{\rm mom}}
\|\widehat f_l\|_{H^{\ell_{\rm mom}}(Q_n')}^{q/\ell_{\rm mom}}
\le C C_{\rm amp}\varepsilon_{\rm int}^{1/2}e^{-c\langle l,n\rangle^\chi/2}.
\]
For $q=0$ the first inequality is simply the restriction estimate. The local interpolation follows by a bounded extension on the enlarged cube and the $L^2$ interpolation inequality for integer derivatives; {the constants are uniform in} $n$.

Use $\langle l\rangle^\chi+\langle n\rangle^\chi\le2\langle l,n\rangle^\chi$.
For sufficiently small $\lambda'>0$, summing over the cubes gives
\[
\|e^{\lambda'\langle\eta\rangle^\chi}\partial_\eta^\alpha \widehat f_l\|_2^2 \le C C_{\rm amp}^2\varepsilon_{\rm int} e^{-(c/2)\langle l\rangle^\chi} \sum_{n\in\mathbb Z^3}e^{-(c/2)\langle n\rangle^\chi/2} \le C C_{\rm amp}^2\varepsilon_{\rm int} e^{-(c/2)\langle l\rangle^\chi}.
\]
\par
Plancherel proves \eqref{eq:weightedGevreyInterpolation}.

For the series norm,
$\lceil\ell\rceil+1\le\ell+2\le L$ and
$\langle v\rangle^\ell\le C_\ell\sum_{|\alpha|\le \lceil\ell\rceil}|v^\alpha|$.

Use the derivative estimate \eqref{eq:gevreyRadiusComparison}
and the exact commutation identity below.
\par
\begin{align*}
v^\alpha\partial_1^nf_l
&=\sum_{j=0}^{\min(n,\alpha_1)}(-1)^j\binom nj
\frac{\alpha_1!}{(\alpha_1-j)!}
\partial_1^{n-j}(v^{\alpha-je_1}f_l).
\end{align*}
Apply these identities to $|\alpha|\le \lceil\ell\rceil$. They imply
\[
\begin{aligned}
\|\langle v\rangle^\ell\partial_1^nf_l\|_2
&\le C_{\rm amp}\varepsilon_{\rm int}^{1/2}e^{-c\langle l\rangle^\chi} C^{n+1}(n+1)^{\lceil\ell\rceil}(n!)^{1/\chi},\\
\|f_l\|_{\mathscr A_\ell^{\lambda_{\rm v},\chi}}^2
&\le C C_{\rm amp}^2\varepsilon_{\rm int} e^{-2c\langle l\rangle^\chi}\sum_{n\ge0}(C\lambda_{\rm v})^{2n}(n+1)^{2\lceil\ell\rceil+24} \le C C_{\rm amp}^2\varepsilon_{\rm int} e^{-2c\langle l\rangle^\chi}
\end{aligned}
\]
if $C\lambda_{\rm v}<1/2$. One extra derivative changes the polynomial
factor by $(n+1)^{1/\chi}$; for $v_1f_l$ use moments through degree
$\lceil\ell\rceil+1\le L$. The same convergent series proves
the two remaining assertions.
\end{proof}
\begin{lem}\label{lem:physicalCoefficientDerivatives}
For any fixed integer $L$, there are fixed $C_0,\varrho_{\rm b}>0$ such that, for every $n\ge0$, $l\in\mathbb Z$, and $|\alpha|\le2$,
\begin{equation}\label{eq:physicalCoefficientDerivatives}
\sup_v\langle v\rangle^L
|\partial_{v_1}^n\partial_v^\alpha g_l(t,v)|
\le C_0^{n+1}g_\chi(n)\varepsilon_{\rm c} e^{-\varrho_{\rm b}|l|^\chi}.
\end{equation}
One may take $\varepsilon_{\rm c}=\varepsilon_N$ for {$g\in\{g^{\rm b,\nu},g^{\rm b}_{\rm app},{g^{\rm b}_{\rm err}}\}$}. These bounds are uniform for $0\le t\le {T_1}-t_0$.
\end{lem}

\begin{proof}
Choose an even integer $L_*\ge\max\{L,0\}$ so that $\langle v\rangle^{L_*}$ is a polynomial of fixed degree $L_*$.
Proposition~\ref{prop:bguniform} and
Lemma~\ref{lem:weightedGevreyInterpolation}, with $\varepsilon_{\rm int}=1$,
$C_{\rm amp}=C\varepsilon_N$ and sufficiently many fixed moments, give
a radius $\lambda'>0$ and spatial decay $\varrho_{\rm b}>0$, both
independent of $\nu$, such that
\begin{equation}\label{eq:coefficientWeightedFourier}
\sum_{|\beta|\le L_*}
\|e^{\lambda'\langle\eta\rangle^\chi}\partial_\eta^\beta
\widehat g_l^{\rm b,\nu}(t,\eta)\|_{L^2_\eta}
\le C\varepsilon_N e^{-\varrho_{\rm b}|l|^\chi}.
\end{equation}
Applying \eqref{eq:gevreyRadiusComparison} with four additional
velocity derivatives gives
\par
\begin{equation}\label{eq:coefficientFactorialAbsorption}
\begin{aligned}
\sup_\eta\langle\eta\rangle^{n+4}e^{-\lambda'\langle\eta\rangle^\chi/2}
&\le C^{n+1}((n+4)!)^{1/\chi},\\
(n+1)^{L_*}((n+4)!)^{1/\chi}
&\le C_{L_*}^{n+1}(n!)^{1/\chi}(n+1)^{L_*+4/\chi} \le \widetilde C_{L_*}^{n+1}g_\chi(n).
\end{aligned}
\end{equation}
The constants may increase at each fixed radius reduction.
For $|\alpha|\le2$ and $|\gamma|\le2$, each monomial in
$\partial_v^\gamma(\langle v\rangle^{L_*}\partial_{v_1}^n\partial_v^\alpha g_l^{\rm b,\nu})$
has Fourier transform bounded by a finite sum of terms
\[
C_{L_*}(n+1)^{L_*}\langle\eta\rangle^{n+4}
|\partial_\eta^\beta\widehat g_l^{\rm b,\nu}|,
\qquad |\beta|\le L_*.
\]
Indeed Leibniz's rule differentiates $\eta_1^n\eta^\alpha$ at most
$L_*$ times; {these derivatives lower the polynomial degree and contribute} at most
$C_{L_*}(n+1)^{L_*}$. Thus \eqref{eq:coefficientWeightedFourier},
\eqref{eq:coefficientFactorialAbsorption}, and
$H^2(\mathbb R^3)\subset L^\infty(\mathbb R^3)$ yield
\[
{
\|\langle v\rangle^L\partial_{v_1}^n\partial_v^\alpha g_l^{\rm b,\nu}\|_\infty \le C\|\langle v\rangle^{L_*}\partial_{v_1}^n\partial_v^\alpha g_l^{\rm b,\nu}\|_{H^2_v} \le C_0^{n+1}g_\chi(n)\varepsilon_N e^{-\varrho_{\rm b}|l|^\chi}.
}
\]
For $g^{\rm b}_{\rm app}$ the only modes are $l=\pm1$, and Gaussian Fourier
translation gives, at a fixed reduced radius,
\[
\sum_{|\beta|\le L_*}
\|e^{\lambda'\langle\eta\rangle^\chi}\partial_\eta^\beta
\widehat g_{\rm app,l}^{\rm b}(\eta)\|_2
\le C_{L_*}a e^{C t_0^\chi}\mathbf1_{|l|=1}
\le C\varepsilon_N e^{-\varrho_{\rm b}|l|^\chi}.
\]
The preceding derivation therefore applies to $g^{\rm b}_{\rm app}$ as well.
The triangle inequality for {${g^{\rm b}_{\rm err}}=g^{\rm b,\nu}-g^{\rm b}_{\rm app}$} completes
\eqref{eq:physicalCoefficientDerivatives}. All moments and radii are
fixed before $\nu\to0$.
\end{proof}

\subsection{Collision forms, weights and the background potential}\label{app:energyCoefficients}

We first estimate both collision orderings in the centered norms of
\eqref{eq:physicalCenteredNorms}. We then prove the coefficient and multiplier
bounds needed for the finite-order energy norms in
\eqref{eq:fullLinearEnergyNorms}--\eqref{eq:errorComparisonNorm} and {the source
bounds associated with \eqref{eq:quadraticDualNorm} and
\eqref{eq:quadraticSourceSpaces},}
together with their Gaussian-weighted variants in
\eqref{eq:gaussianEnergyPullback}--\eqref{eq:broaderGaussianNorms}, including the
{potential conjugation and the bilinear collision estimate.} These results are used
in the complete forced problem and in the nonlinear comparison.

\begin{lem}\label{lem:physicalCollisionForm}
Fix $\ell\ge0$. Suppose that $g$ satisfies \eqref{eq:physicalCoefficientDerivatives}, with $L$ sufficiently large depending on $\ell$. For small fixed $\varrho_0,\lambda_{\rm v}$, the associated {$\mathcal C_{F^{\rm c}}$} satisfies
\begin{equation}\label{eq:physicalCollisionForm}
\begin{aligned}
&\left|\sum_k\kappa_{k,j}^2\sum_{n\ge0}a_n^2
\left\langle\widehat Z_{k,c}^n(\mathcal C_{F^{\rm c}}f)_k,
\widehat Z_{k,c}^nw_k\right\rangle_{\ell}\right|\\
&\quad\le C\varepsilon_{\rm c}\left[
\mathcal D^{\rm ph}_{c,j}(s;f)^{1/2}
\mathcal D^{\rm ph}_{c,j}(s;w)^{1/2}
+(1+(s-t_0)^2)
\mathcal E^{\rm ph}_{c,j}(s;f)^{1/2}
\mathcal E^{\rm ph}_{c,j}(s;w)^{1/2}\right].
\end{aligned}
\end{equation}
The estimate holds with the same constant for finite derivative sums and finite spatial projections. It contains no derivative of $f$ outside the displayed dissipation.
The weights satisfy \begin{equation}\label{eq:physicalSpatialRatio}
\kappa_{k,j}/\kappa_{m,j}\le e^{\varrho_0|k-m|^\chi}.
\end{equation}
Also $[\partial_s+ikv_1,\widehat Z_{k,c}(s)]=0$.
\end{lem}

\begin{proof}
The weight ratio and the transport commutator follow from
\[
|k-j|^\chi-|m-j|^\chi\le|k-m|^\chi,\qquad
[\partial_s,\widehat Z_{k,c}]=ik,
\qquad[ikv_1,\widehat Z_{k,c}]=-ik.
\]
Fix $t=s-t_0$, a coefficient mode $l$, and $k=m+l$.
For either velocity collision ordering, denote the corresponding operator by
$\mathcal A_{l,0}$ and define
\[
\mathcal A_{l,r+1}=[\partial_{v_1},\mathcal A_{l,r}]+ilt\mathcal A_{l,r}.
\]
Since $\widehat Z_{k,c}\mathcal A_{l,0}-\mathcal A_{l,0}\widehat Z_{m,c}=\mathcal A_{l,1}$,
induction and Pascal's identity give
\begin{equation}\label{eq:physicalCollisionCommutator}
\widehat Z_{k,c}^n\mathcal A_{l,0}y
=\sum_{r=0}^n\binom nr\mathcal A_{l,r}\widehat Z_{m,c}^{n-r}y.
\end{equation}
The operators $\mathcal A_{l,r}$ do not contain $c$.

Use the Gaussian-conjugated derivatives defined in \eqref{eq:gaussianConjugatedDerivatives}.
For the first collision ordering put {$\sigma[{F^{\rm c}}_l]=\Phi*{F^{\rm c}}_l$}, {$b_{i}[{F^{\rm c}}_l]=(\partial_q\Phi_{iq})*{F^{\rm c}}_l$}, and
use the superscript $[r]$ for application of $(\partial_{v_1}+ilt)^r$
to a coefficient.
Because {${F^{\rm c}}_l=e^{-iltv_1}\sqrt \mu g_l$},
\[
(\sigma[{F^{\rm c}}_l])^{[r]}=\Phi*\big(e^{-iltv_1}\partial_{v_1}^r(\sqrt \mu g_l)\big),
\qquad
(b_{i}[{F^{\rm c}}_l])^{[r]}=(\partial_q\Phi_{iq})*
\big(e^{-iltv_1}\partial_{v_1}^r(\sqrt \mu g_l)\big).
\]
Gaussian differentiation and \eqref{eq:physicalCoefficientDerivatives} imply
\[
|\partial_{v_1}^r(\sqrt \mu g_l)(v)|
\le C^{r+1}g_\chi(r)\varepsilon_{\rm c} e^{-\varrho_{\rm b}|l|^\chi}e^{-c|v|^2}.
\]
Lemma~\ref{lem:gaussianCollisionCoefficients} therefore yields
\begin{equation}\label{eq:physicalCoefficientMatrix}
\begin{aligned}
|\overline X(\sigma[{F^{\rm c}}_l])^{[r]}Y|
&\le C^{r+1}g_\chi(r)\varepsilon_{\rm c} e^{-\varrho_{\rm b}|l|^\chi}
(\overline X\sigma X)^{1/2}(\overline Y\sigma Y)^{1/2},\\
|\sigma^{-1/2}(b[{F^{\rm c}}_l])^{[r]}|
&\le C^{r+1}g_\chi(r)\varepsilon_{\rm c} e^{-\varrho_{\rm b}|l|^\chi}
\langle v\rangle^{-1/2}.
\end{aligned}
\end{equation}
Here and below {$\sigma=\Phi*\mu $}.
The first ordering is {$\mathcal A_{l,0}=D_{1,i}^-(\sigma_{iq}[{F^{\rm c}}_l]D_{1,q}^--b_{i}[{F^{\rm c}}_l])$}.
Using $[\partial_{v_1},D_{1,i}^-]=-\delta_{i1}/2$ {gives the complete expansion}
\begin{equation}\label{eq:physicalFirstCommutatorExpansion}
\begin{aligned}
\mathcal A_{l,r}={}&D_{1,i}^-\big((\sigma_{iq}[{F^{\rm c}}_l])^{[r]}D_{1,q}^--(b_{i}[{F^{\rm c}}_l])^{[r]}\big)\\
&-\frac r2\delta_{i1}
\big((\sigma_{iq}[{F^{\rm c}}_l])^{[r-1]}D_{1,q}^--(b_{i}[{F^{\rm c}}_l])^{[r-1]}\big)
-\frac r2D_{1,i}^-\big((\sigma_{i1}[{F^{\rm c}}_l])^{[r-1]}\,\cdot\big)\\
&+\frac{r(r-1)}4(\sigma_{11}[{F^{\rm c}}_l])^{[r-2]}\,\cdot.
\end{aligned}
\end{equation}
Terms with a negative coefficient-derivative index are omitted.
In particular, the expression with $D_{1,i}^-$ on the second line is a
composition with multiplication, not just {differentiation of the coefficient alone.}

{The undifferentiated pairing is}
\begin{equation}\label{eq:physicalFirstPair}
\langle\mu^{-1/2}Q({F^{\rm c}}_l,\sqrt \mu y),z\rangle_{\ell}
=-\int(\sigma_{iq}[{F^{\rm c}}_l]D_{1,q}^-y-b_{i}[{F^{\rm c}}_l]y)\overline{D_{1,i}^+(\langle v\rangle^{2\ell}z)}\,dv.
\end{equation}
The definition of the form norm and the lower-order weighted bound give
\begin{equation}\label{eq:physicalFormFactors}
\begin{aligned}
\|\sigma^{1/2}D_1^-y\|_{L_\ell^2}
+\|\langle v\rangle^{-\ell}\sigma^{1/2}D_1^+(\langle v\rangle^{2\ell}y)\|_2
+\|\langle v\rangle^{\ell-1/2}y\|_2
\le C_\ell\|y\|_{\Sigma_\ell}.
\end{aligned}
\end{equation}
{Indeed,
\[
D_1^+(\langle v\rangle^{2\ell}y)
=\langle v\rangle^{2\ell}\left(D_1^+y+\frac{2\ell v}{\langle v\rangle^2}y\right),
\qquad
\left|\frac{\ell\sigma^{1/2}v}{\langle v\rangle^2}\right|
\le C_\ell\langle v\rangle^{-3/2},
\]
and}
$|\sigma^{1/2}e_1|\le C\langle v\rangle^{-1/2}$.
Apply \eqref{eq:physicalCoefficientMatrix} to the four terms in
\eqref{eq:physicalFirstCommutatorExpansion}, integrating an outer $D_1^-$
once by parts. Their respective bounds are
\[
C\varepsilon_{\rm c} e^{-\varrho_{\rm b}|l|^\chi}
\big(C^rg_\chi(r)+rC^rg_\chi(r-1)
+r(r-1)C^rg_\chi(r-2)\big)
\|y\|_{\Sigma_\ell}\|z\|_{\Sigma_\ell}.
\]
The terms with {$b_{i}[{F^{\rm c}}_l]$} use the third norm in
\eqref{eq:physicalFormFactors}; those with {$\sigma_{i1}[{F^{\rm c}}_l]$} {use the displayed}
bound for $\sigma^{1/2}e_1$. For $d=0,1,2$ and $r\ge d$,
\[
\frac{r!}{(r-d)!}g_\chi(r-d)
\le C_d^{r+1}g_\chi(r).
\]
Consequently,
\begin{equation}\label{eq:physicalFirstOrdering}
|\langle\mathcal A_{l,r}y,z\rangle_{\ell}|
\le C^{r+1}g_\chi(r)\varepsilon_{\rm c} e^{-\varrho_{\rm b}|l|^\chi}
\|y\|_{\Sigma_\ell}\|z\|_{\Sigma_\ell}.
\end{equation}

For the second ordering, the cancellation formula gives
\[
\mu^{-1/2}Q(\sqrt \mu y,{F^{\rm c}}_l)
=(\Phi_{iq}*(\sqrt \mu y))e^{-iltv_1}H_{l,iq}+8\pi {F^{\rm c}}_ly,
\]
where
\[
H_{l,iq}=(D_{1,i}^--ilt\delta_{i1})(D_{1,q}^--ilt\delta_{q1})g_l.
\]
{The integral kernel is}
{$K_l(v,v_*)=\Phi_{iq}(v-v_*)\sqrt \mu (v_*)e^{-iltv_1}H_{l,iq}(v)$}.
Integration by parts in $v_*$ shows that the $r$th operator commutator
has kernel $(\partial_{v_1}+\partial_{v_{*1}}+ilt)^rK_l$. Since
$(\partial_{v_1}+\partial_{v_{*1}})\Phi(v-v_*)=0$,
\begin{equation}\label{eq:physicalSecondKernelExpansion}
(\partial_{v_1}+\partial_{v_{*1}}+ilt)^rK_l
=\Phi_{iq}(v-v_*)e^{-iltv_1}
\sum_{a=0}^r\binom ra
(\partial_1^aH_{l,iq})(v)(\partial_1^{r-a}\sqrt \mu )(v_*).
\end{equation}
Leibniz's rule in $H_{l,iq}$ and the coefficient hypothesis give
\begin{align*}
|\partial_1^aH_{l,iq}(v)|
&\le C^{a+1}g_\chi(a)\varepsilon_{\rm c} e^{-\varrho_{\rm b}|l|^\chi}
(1+l^2t^2)\langle v\rangle^{-L+2},\\
|\partial_1^{r-a}\sqrt \mu (v_*)|
&\le C^{r-a+1}g_\chi(r-a)e^{-c|v_*|^2},\\
\sum_{a=0}^r\binom rag_\chi(a)g_\chi(r-a)
&\le(r+1)g_\chi(r)\le C^{r+1}g_\chi(r).
\end{align*}
{After conjugation by $\langle v\rangle^\ell$,} the absolute kernel is bounded by
\[
C^{r+1}g_\chi(r)\varepsilon_{\rm c} e^{-\varrho_{\rm b}|l|^\chi}(1+l^2t^2)
\frac{\langle v\rangle^{\ell-L+2}e^{-c|v_*|^2}
\langle v_*\rangle^{-\ell}}{|v-v_*|}.
\]
Choose $L>\ell+8$, so $L-\ell-2>3$.
Equation \eqref{eq:weightedCoulombSchur} bounds the operator by the
displayed amplitude $C^{r+1}g_\chi(r)\varepsilon_{\rm c} e^{-\varrho_{\rm b}|l|^\chi}(1+l^2t^2)$.
The multiplication term has commutator
{$8\pi e^{-iltv_1}\partial_1^r(\sqrt \mu g_l)$} and obeys the same bound.
Schur's inequality therefore gives
\begin{equation}\label{eq:physicalSecondOrdering}
\|\mathcal A_{l,r}y\|_{L_\ell^2}
\le C^{r+1}g_\chi(r)\varepsilon_{\rm c} e^{-\varrho_{\rm b}|l|^\chi}
(1+l^2t^2)\|y\|_{L_\ell^2}.
\end{equation}

It remains to sum the two indices. For the first ordering set
\[
U_{k,n}=\kappa_{k,j}a_n\|\widehat Z_{k,c}^nf_k\|_{\Sigma_\ell},
\qquad V_{k,n}=\kappa_{k,j}a_n\|\widehat Z_{k,c}^nw_k\|_{\Sigma_\ell}.
\]
Extend these sequences by zero for $n<0$. Equations
\eqref{eq:centeredfactorialratio} and \eqref{eq:physicalSpatialRatio}
give the summable convolution coefficients
\[
\sum_l\sum_{r\ge0}
e^{-(\varrho_{\rm b}-\varrho_0)|l|^\chi}(C\lambda_{\rm v})^r\le C
\]
when $\varrho_0<\varrho_{\rm b}$ and $C\lambda_{\rm v}<1/2$.
Discrete Young's inequality therefore bounds the absolute first-ordering sum by
\begin{align*}
&C\varepsilon_{\rm c}\sum_{k,n}\sum_l\sum_{r=0}^n
e^{-(\varrho_{\rm b}-\varrho_0)|l|^\chi}(C\lambda_{\rm v})^r
U_{k-l,n-r}V_{k,n}\\
&\quad\le C\varepsilon_{\rm c}\left(\sum_l\sum_{r\ge0}
e^{-(\varrho_{\rm b}-\varrho_0)|l|^\chi}(C\lambda_{\rm v})^r\right)
\|U\|_2\|V\|_2\\
&\quad\le C\varepsilon_{\rm c}\mathcal D^{\rm ph}_{c,j}(f)^{1/2}
\mathcal D^{\rm ph}_{c,j}(w)^{1/2}.
\end{align*}
For the second ordering replace $\Sigma_\ell$ by ${L_\ell^2}$ and use
\[
1+l^2t^2\le\langle l\rangle^2(1+t^2),\qquad
\sum_l\sum_{r\ge0}\langle l\rangle^2
e^{-(\varrho_{\rm b}-\varrho_0)|l|^\chi}(C\lambda_{\rm v})^r\le C.
\]
The same convolution bound gives the energy term in
\eqref{eq:physicalCollisionForm}. Restricting either sequence to finite
sets only decreases these absolute sums, proving the uniform finite-sum
and finite-projection statements.
\end{proof}

\begin{lem}\label{lem:fullBackgroundCoefficients}
For every fixed velocity weight $L$, spatial order $r$, and velocity derivative order $r_{\rm v}$, the constants in
$L_N$ can be enlarged so that, with $\varepsilon_N=aL_N$,
\begin{equation}\label{eq:fullBackgroundFourierEnvelope}
\sum_{l\in\mathbb Z}\langle l\rangle^{r+2}
\sup_v\langle v\rangle^L
|\partial_v^\alpha(f^{\rm b,\nu})_l(s,v)|
\le C_{r,L,r_{\rm v}}\varepsilon_N\langle s\rangle^{|\alpha|},
\qquad |\alpha|\le r_{\rm v}.
\end{equation}
All powers and radii used to obtain these bounds are fixed independently of $\nu$.
\end{lem}

\begin{proof}
Put $t=s-t_0$. {The shear transformation and the corresponding Fourier formula are}
\[
\begin{aligned}
f^{\rm b,\nu}(s,x,v)&=g^{\rm b,\nu}(t,x-tv_1,v),\\
\partial_v^\alpha f^{\rm b,\nu}(s,x,v)
&=\sum_{a=0}^{\alpha_1}\binom{\alpha_1}{a}(-t)^a
(\partial_z^a\partial_v^{\alpha-ae_1}g^{\rm b,\nu})(t,x-tv_1,v),\\
(\partial_v^\alpha f^{\rm b,\nu})_l(s,v)
&=e^{-iltv_1}\sum_{a=0}^{\alpha_1}\binom{\alpha_1}{a}
(-ilt)^a\partial_v^{\alpha-ae_1}g_l^{\rm b,\nu}(t,v).
\end{aligned}
\]
The last phase has modulus one and is removed before applying velocity
Sobolev embedding. For every fixed $d\ge0$ and multi-index $\beta$,
$H^2_v\subset L^\infty_v$, sequence Cauchy--Schwarz and Parseval give
\[
\begin{aligned}
\sum_l\langle l\rangle^{d+2}
\sup_v\langle v\rangle^L|\partial_v^\beta g_l^{\rm b,\nu}|
&\le C\sum_l\langle l\rangle^{d+2}
\|\langle v\rangle^L\partial_v^\beta g_l^{\rm b,\nu}\|_{H^2_v}\\
&\le C\left(\sum_l\langle l\rangle^{-4}\right)^{1/2}
\left(\sum_l\langle l\rangle^{2(d+4)}
\|\langle v\rangle^L\partial_v^\beta g_l^{\rm b,\nu}\|_{H^2_v}^2\right)^{1/2}\\
&\le C_{d,L,\beta}\varepsilon_N.
\end{aligned}
\]
The last inequality is \eqref{eq:bguniformsobolev} at the finitely many
orders $d+4$ in $z$ and $|\beta|+2$ in velocity, with moment $L$.
Apply it with $d=r+a$ and $\beta=\alpha-ae_1$ in the Fourier formula:
\[
\begin{aligned}
\sum_l\langle l\rangle^{r+2}
\sup_v\langle v\rangle^L|(\partial_v^\alpha f^{\rm b,\nu})_l|
&\le C_\alpha\sum_{a=0}^{\alpha_1}t^a
\sum_l\langle l\rangle^{r+a+2}
\sup_v\langle v\rangle^L|\partial_v^{\alpha-ae_1}g_l^{\rm b,\nu}|\\
&\le C_{r,L,\alpha}\varepsilon_N\langle s\rangle^{|\alpha|}.
\end{aligned}
\]
This proves \eqref{eq:fullBackgroundFourierEnvelope} without an additional
power of $t$ from velocity embedding.

\end{proof}

\begin{lem}\label{lem:fullCollisionForm}
{Let $\mathcal C_{\rm b}(s)$ be the normalized
background collision operator defined in
\eqref{eq:backgroundNormalizedLinearization}.}
For sufficiently small $\varepsilon_N$,
\begin{equation}\label{eq:fullCollisionForm}
\operatorname{Re}\sum_{j=0}^r
\left\langle\partial_{x_1}^j\mathcal C_{\rm b}y,
\partial_{x_1}^jy\right\rangle_{2r+6-2j}
\le-c_r\nu^{-1}\mathcal D_r[y]
+C_r(1+\varepsilon_N\langle s\rangle^2)\mathcal E_r[y].
\end{equation}
\end{lem}

\begin{proof}
Write {${F^{\rm c}}=F^{\rm b,\nu}-\mu =\sqrt \mu  f^{\rm b,\nu}$} and set
\[
y_j=\partial_{x_1}^jy.
\]
The Maxwellian contribution, by \eqref{eq:allWeightMaxwellianCoercivity}, is
\[
-\operatorname{Re}\sum_{j=0}^r\langle\mathcal Ly_j,y_j\rangle_{2r+6-2j}
\le-2c_r\nu^{-1}\mathcal D_r[y]+C_r\mathcal E_r[y].
\]
For $0\le q\le r$, \eqref{eq:fullBackgroundFourierEnvelope} and
Lemma~\ref{lem:gaussianCollisionCoefficients} give
\begin{align*}
|\overline X\partial_{x_1}^q \sigma[{F^{\rm c}}]Y|
&\le C_r\varepsilon_N(\overline X\sigma X)^{1/2}(\overline Y\sigma Y)^{1/2},\\
|\sigma^{-1/2}\partial_{x_1}^qb[{F^{\rm c}}]|
&\le C_r\varepsilon_N\langle v\rangle^{-1/2},
\end{align*}
where {$\sigma[{F^{\rm c}}]=\Phi*{F^{\rm c}}$} and {$b_{i}[{F^{\rm c}}]=(\partial_a\Phi_{ia})*{F^{\rm c}}$}.
The first remainder is exactly
\[
S_1=\sum_{j=0}^r\sum_{q=0}^j\binom jq
\langle\mu^{-1/2}Q(\partial_{x_1}^q {F^{\rm c}},\sqrt \mu  y_{j-q}),y_j\rangle_{2r+6-2j}.
\]
Apply \eqref{eq:physicalFirstPair} and \eqref{eq:physicalFormFactors}
with resulting velocity weight $2r+6-2j$. Since $(2r+6-2j+2q)-(2r+6-2j)=2q\ge0$,
\begin{align*}
\|y_{j-q}\|_{L^2_x\Sigma_{2r+6-2j}}&\le\|y_{j-q}\|_{L^2_x\Sigma_{2r+6-2j+2q}},\\
|S_1|&\le C_r\varepsilon_N\sum_{j=0}^r\sum_{q=0}^j\|y_{j-q}\|_{L^2_x\Sigma_{2r+6-2j+2q}}\|y_j\|_{L^2_x\Sigma_{2r+6-2j}}
\le C_r\varepsilon_N\nu^{-1}\mathcal D_r[y].
\end{align*}

For the second ordering the exact formula is
\[
\mu^{-1/2}Q(\sqrt \mu  f,H)=(\Phi_{ia}*(\sqrt \mu  f))\mu^{-1/2}\partial_{v_iv_a}H+8\pi Hf.
\]
Expanding the spatial derivatives gives
\[
S_2=\sum_{j=0}^r\sum_{q=0}^j\binom jq
\langle\mu^{-1/2}Q(\sqrt \mu  y_{j-q},\partial_{x_1}^q {F^{\rm c}}),y_j\rangle_{2r+6-2j}.
\]
The derivative bounds for $F^{\rm c}$ obtained from \eqref{eq:physicalCoefficientDerivatives} yield
\[
|\mu^{-1/2}\partial_{v_iv_a}\partial_{x_1}^q {F^{\rm c}}|
\le C_r\varepsilon_N\langle s\rangle^2\langle v\rangle^{-L+2}.
\]
Choose $L>2r+14$. After conjugation by the weights $\langle v_*\rangle^{-(2r+6-2j+2q)}$ and $\langle v\rangle^{2r+6-2j}$, the kernel has the pointwise bound
\[
C_r\varepsilon_N\langle s\rangle^2
\frac{\langle v\rangle^{2r+8-2j-L}\sqrt \mu (v_*)
\langle v_*\rangle^{-2r-6+2j-2q}}{|v-v_*|}.
\]
Since $L-2r-8+2j>3$ and $2r+6-2j+2q\ge0$,
\eqref{eq:weightedCoulombSchur} applies with decay exponent $L-2r-8+2j$.
The resulting operator norm is at most $C_{r,L}\varepsilon_N\langle s\rangle^2$.
The multiplier {$8\pi\partial_{x_1}^q {F^{\rm c}}$} is bounded by $C_r\varepsilon_N$.
Schur's inequality and $\langle v\rangle^{2r+6-2j}\le \langle v\rangle^{2r+6-2j+2q}$ therefore imply
\begin{align*}
&\|\langle v\rangle^{2r+6-2j}\mu^{-1/2}Q(\sqrt \mu  y_{j-q},\partial_{x_1}^q {F^{\rm c}})\|_2\\
&\quad\le C_r\varepsilon_N\langle s\rangle^2
\|\langle v\rangle^{2r+6-2j+2q}y_{j-q}\|_2,\\
|S_2|&\le C_r\varepsilon_N\langle s\rangle^2
\sum_{j=0}^r\sum_{q=0}^j
\|\langle v\rangle^{2r+6-2j+2q}y_{j-q}\|_2
\|\langle v\rangle^{2r+6-2j}y_j\|_2\\
&\le C_r\varepsilon_N\langle s\rangle^2\mathcal E_r[y].
\end{align*}
Choose $\varepsilon_N$ so that $C_r\varepsilon_N\le c_r$.
Adding the three contributions proves \eqref{eq:fullCollisionForm}.
\end{proof}
\begin{lem}\label{lem:weightedDualMultipliers}
For $0\le\ell_-\le\ell_+$, the velocity form spaces satisfy
\begin{equation}\label{eq:weightedDualEmbeddingGeneral}
\begin{aligned}
\|\psi\|_{\Sigma_{\ell_-}}&\le\|\psi\|_{\Sigma_{\ell_+}},\qquad \|\langle v\rangle^{2(\ell_--\ell_+)}\psi\|_{\Sigma_{\ell_+}}
\le C_{\ell_-,\ell_+}\|\psi\|_{\Sigma_{\ell_-}},\\
\|g\|_{\Sigma_{\ell_-}^*}
&\le C_{\ell_-,\ell_+}\|g\|_{\Sigma_{\ell_+}^*}.
\end{aligned}
\end{equation}
Fix an integer $r\ge0$, and let $\mathfrak m\in C^r(\mathbb T)$ depend only on
$x_1$. Then
\begin{equation}\label{eq:spatialMultiplierNorms}
\begin{aligned}
\mathcal E_r[\mathfrak mf]^{1/2}&\le C_r\|\mathfrak m\|_{C^r_x}\mathcal E_r[f]^{1/2},\qquad \mathcal D_r[\mathfrak mf]^{1/2}\le C_r\|\mathfrak m\|_{C^r_x}\mathcal D_r[f]^{1/2},\\
\sum_{j=0}^{r}\left\|\partial_{x_1}^j\!\left(\mathfrak mg\right)\right\|_{L_x^2\Sigma_{2r+6-2j}^*}&\le C_r\|\mathfrak m\|_{C^r_x}\sum_{j=0}^{r}\left\|\partial_{x_1}^j\!\left(g\right)\right\|_{L_x^2\Sigma_{2r+6-2j}^*}.
\end{aligned}
\end{equation}
In particular, for real $\phi\in C^r(\mathbb T)$ and $1\le{{\varpi}}\le3/2$,
multiplication by $e^{\phi/(2{{\varpi}})}$ and its inverse is
bounded in these three spaces, uniformly when $\|\phi\|_{C^r}$ is bounded.
\end{lem}

\begin{proof}
The first inequality in \eqref{eq:weightedDualEmbeddingGeneral} follows
from the nonnegative integrand defining $\Sigma_\ell$. For the second,
put $h=\langle v\rangle^{2(\ell_--\ell_+)}$. The
Maxwellian matrix bounds and the weighted coercive norm give
\[
|\sigma^{1/2}\nabla\log\langle v\rangle|\le C\langle v\rangle^{-3/2},
\qquad
\|\langle v\rangle^{\ell_--1/2}\psi\|_2
\le C\|\psi\|_{\Sigma_{\ell_-}},
\]
as in Lemma~\ref{lem:centeredcoercivity}. Since
$\nabla(h\psi)=h\nabla\psi-2(\ell_+-\ell_-) h\psi\nabla\log\langle v\rangle$,
\begin{align*}
\|\langle v\rangle^{\ell_+}\sigma^{1/2}\nabla(h\psi)\|_2
&\le\|\langle v\rangle^{\ell_-}\sigma^{1/2}\nabla\psi\|_2
+C(\ell_+-\ell_-)\|\langle v\rangle^{\ell_--3/2}\psi\|_2
\le C_{\ell_-,\ell_+}\|\psi\|_{\Sigma_{\ell_-}},\\
\|\langle v\rangle^{\ell_+}|\sigma^{1/2}v|h\psi\|_2
&\le\|\langle v\rangle^{\ell_-}|\sigma^{1/2}v|\psi\|_2.
\end{align*}
These are the two components of $\|h\psi\|_{\Sigma_{\ell_+}}$.
For a smooth source $g$, the weighted pairing then yields
\[
|\langle g,\psi\rangle_{\ell_-}|=|\langle g,h\psi\rangle_{\ell_+}| \le\|g\|_{\Sigma_{\ell_+}^*}\|h\psi\|_{\Sigma_{\ell_+}} \le C_{\ell_-,\ell_+}\|g\|_{\Sigma_{\ell_+}^*} \|\psi\|_{\Sigma_{\ell_-}}.
\]

\par
Taking the test supremum proves the dual embedding; completion extends
it to distributional sources with the stated finite norm.

For $0\le j_1\le j\le r$, use the single product rule
\[
\partial_{x_1}^j(\mathfrak mf)=\sum_{j_1=0}^j\binom j{j_1}
(\partial_{x_1}^{j_1}\mathfrak m)\partial_{x_1}^{j-j_1}f,
\qquad 2r+6-2j+2j_1\ge2r+6-2j.
\]
Monotonicity of the weighted $L^2$ and form norms gives the first two
bounds in \eqref{eq:spatialMultiplierNorms}. For the dual norm the
embedding just proved gives explicitly
\[
\|\partial_{x_1}^j(\mathfrak mg)\|_{L^2_x\Sigma_{2r+6-2j}^*}
\le C_r\sum_{j_1=0}^j\binom j{j_1}\|\partial_{x_1}^{j_1}\mathfrak m\|_\infty
\|\partial_{x_1}^{j-j_1}g\|_{L^2_x\Sigma_{2r+6-2j+2j_1}^*}.
\]
Finite-sum Cauchy--Schwarz and summation in $j$ prove the third bound;
{the distributional product rule extends it to sources whose spatial derivatives have the displayed weighted dual norms.}
Finally, the finite chain rule gives
\[
\|e^{\pm\phi/(2{{\varpi}})}\|_{C^r_x}
\le C_r e^{\|\phi\|_\infty/(2{{\varpi}})}(1+\|\phi\|_{C^r_x})^r,
\qquad 1\le{{\varpi}}\le3/2.
\]
Apply \eqref{eq:spatialMultiplierNorms} to each sign. This proves both
boundedness assertions and the resulting equivalence of the norms.
\end{proof}

\begin{lem}\label{lem:potentialTransport}
Fix an integer $r\ge0$, $\nu>0$, $1\le{{\varpi}}\le3/2$, and a real potential
$\phi\in W^{1,1}(I;C^r_x)\cap L^\infty(I;C^{r+1}_x)$, with
$E=-\partial_{x_1}\phi\in L^\infty(I;C^{\max(r,1)}_x)$.
Write {$y=e^{\phi/(2{\varpi})}\mu_{{\varpi}}^{-1/2}U$}.
For smooth $U$, and then in distributions whenever the terms exist,
\begin{equation}\label{eq:generalPotentialConjugation}
{
e^{\phi/(2{\varpi})}\mu_{{\varpi}}^{-1/2}(\partial_s+v_1\partial_{x_1}+E\partial_{v_1})U ={}(\partial_s+v_1\partial_{x_1}+E\partial_{v_1})y -\frac{\partial_s\phi}{2{{\varpi}}}y.
}
\end{equation}
For smooth $y$, put $y_j=\partial_{x_1}^jy$ and
\[
\mathcal T_r=\operatorname{Re}\sum_{j=0}^r
\langle\partial_{x_1}^j
(v_1\partial_{x_1}+E\partial_{v_1})y,y_j\rangle_{2r+6-2j}.
\]
For every $\eta>0$,
\begin{equation}\label{eq:generalFieldEnergy}
|\mathcal T_r|
\le\eta\mathcal D_r[y]
+C_{r,\eta}\big(\|E\|_\infty+\nu^{-1}\|E\|_{C^r_x}^2\big)
\mathcal E_r[y].
\end{equation}
The norms of $y$ and {$\mu_{{\varpi}}^{-1/2}U$} are equivalent in energy
and collision form, with constants depending only on $r$ and the
displayed potential bound.
\end{lem}

\begin{proof}
Let {$f=\mu_{{\varpi}}^{-1/2}U$}. Direct differentiation gives
\[
\begin{aligned}
\mu_{{\varpi}}^{-1/2}(\partial_s+v_1\partial_{x_1}+E\partial_{v_1})U
&=(\partial_s+v_1\partial_{x_1}+E\partial_{v_1})f-\frac{v_1E}{2{\varpi}}f,\\
{e^{-\phi/(2{\varpi})}}(\partial_s+v_1\partial_{x_1}+E\partial_{v_1})e^{\phi/(2{\varpi})}
&=\frac{\partial_s\phi-v_1E}{2{\varpi}}.
\end{aligned}
\]
Multiplying the first identity by $e^{\phi/(2{\varpi})}$ and applying the product rule
to $y=e^{\phi/(2{\varpi})}f$ cancels the two $v_1Ey/(2{\varpi})$ terms and proves
\eqref{eq:generalPotentialConjugation}.
Lemma~\ref{lem:weightedDualMultipliers} gives the two norm equivalences.

The product rule is
\[
\partial_{x_1}^j(E\partial_{v_1}y)
=E\partial_{v_1}y_j+
\sum_{q=1}^j\binom jq(\partial_{x_1}^qE)\partial_{v_1}y_{j-q}.
\]
Periodicity and velocity integration by parts give
\begin{align*}
\operatorname{Re}\langle v_1\partial_{x_1}y_j,y_j\rangle_{2r+6-2j}&=0,\\
\operatorname{Re}\langle E\partial_{v_1}y_j,y_j\rangle_{2r+6-2j}
&=-\tfrac12\int E\partial_{v_1}(\langle v\rangle^{4r+12-4j})|y_j|^2\,dx_1dv,\\
\sum_{j=0}^r\left|\operatorname{Re}
\langle E\partial_{v_1}y_j,y_j\rangle_{2r+6-2j}\right|
&\le C_r\|E\|_\infty\mathcal E_r[y].
\end{align*}
For $1\le q\le j\le r$, the precise weight margin is
\begin{equation}\label{eq:fieldCommutatorWeight}
\begin{aligned}
(2r+6-2j+2q)-\tfrac32-(2r+6-2j)&=2q-\tfrac32\ge\tfrac12,\\
\|\langle v\rangle^{2r+6-2j}\partial_{v_1}y_{j-q}\|_2
&\le C\|y_{j-q}\|_{L^2_x\Sigma_{2r+6-2j+2q}}.
\end{aligned}
\end{equation}
Consequently
\begin{equation}\label{eq:fieldCommutatorPairing}
\begin{aligned}
&\sum_{j=1}^r\sum_{q=1}^j\binom jq
\left|\langle(\partial_{x_1}^qE)
\partial_{v_1}y_{j-q},y_j\rangle_{2r+6-2j}\right|\\
&\quad\le C_r\|E\|_{C^r_x}\sum_{j=1}^r\sum_{q=1}^j\|y_{j-q}\|_{L^2_x\Sigma_{2r+6-2j+2q}}\|\langle v\rangle^{2r+6-2j}y_j\|_2\\
&\quad\le C_r\|E\|_{C^r_x}\nu^{-1/2}\mathcal D_r[y]^{1/2}\mathcal E_r[y]^{1/2}
\le\eta\mathcal D_r[y]
+C_{r,\eta}\nu^{-1}\|E\|_{C^r_x}^2\mathcal E_r[y].
\end{aligned}
\end{equation}
For $r=0$ the double sum is empty. Adding the principal transport
pairing proves \eqref{eq:generalFieldEnergy}. Cutoff and smooth
approximation extend the pairings to the finite energy and form class.
\end{proof}

We extend the bilinear operator $\Gamma$ from \eqref{eq:normalizedGammaExplicit}
to the Gaussian widths used in Lemma~\ref{lem:physicalDerivativeEnergy}.
For $1\le{\varpi}\le3/2$, set
\[
\Gamma_{{\varpi}}(f,g)
=\mu_{{\varpi}}^{-1/2}
 Q(\sqrt{\mu_{{\varpi}}}f,\sqrt{\mu_{{\varpi}}}g).
\]
At ${\varpi}=1$, $\Gamma_{{\varpi}}=\Gamma$.
The collision coefficients are {$\sigma[\sqrt{\mu_{{\varpi}}}f]$} and
{$b[\sqrt{\mu_{{\varpi}}}f]$}, with the same convention
\eqref{eq:normalizedBilinearNotation}.
All form norms below still use the physical matrix {$\sigma=\Phi*\mu $}.
The estimates allow signed and complex functions.

\begin{lem}\label{lem:quadraticLandauForm}
For every fixed $\ell\ge0$ and $1\le{\varpi}\le3/2$,
\begin{equation}\label{eq:quadraticVelocityForm}
\|\Gamma_{{\varpi}}(f,g)\|_{\Sigma_\ell^*}
\le C_\ell\big(
\|f\|_2\|g\|_{\Sigma_\ell}
+\|f\|_{\Sigma_0}\|\langle v\rangle^\ell g\|_2\big).
\end{equation}
Consequently, for the energies and dissipations of \eqref{eq:fullLinearEnergyNorms} and every fixed $r\ge1$,
\begin{equation}\label{eq:quadraticSpatialForm}
\sum_{j=0}^{r}\left\|\partial_{x_1}^j\!\left(\Gamma_{{\varpi}}(g,g)\right)\right\|_{L_x^2\Sigma_{2r+6-2j}^*}
\le C_r\mathcal E_r[g]^{1/2}\nu^{-1/2}\mathcal D_r[g]^{1/2}.
\end{equation}
For two functions $f$ and $g$ one has
\begin{equation}\label{eq:mixedSpatialLandauForm}
\sum_{j=0}^{r}\left\|\partial_{x_1}^j\!\left(\Gamma_{{\varpi}}(f,g)\right)\right\|_{L_x^2\Sigma_{2r+6-2j}^*}
\le C_r\left(
\mathcal E_r[f]^{1/2}\nu^{-1/2}\mathcal D_r[g]^{1/2}
+\nu^{-1/2}\mathcal D_r[f]^{1/2}\mathcal E_r[g]^{1/2}\right).
\end{equation}
For $r=0$, one spatial derivative is required only on the prescribed function $f$:
\begin{equation}\label{eq:broaderGaussianZeroOrder}
\begin{aligned}
\left\|\Gamma_{{\varpi}}(f,g)\right\|_{L_x^2\Sigma_{6}^*}
&\le C\left(\mathcal E_1[f]^{1/2}\nu^{-1/2}\mathcal D_0[g]^{1/2}
+\nu^{-1/2}\mathcal D_1[f]^{1/2}\mathcal E_0[g]^{1/2}\right),\\
\left\|\Gamma_{{\varpi}}(g,f)\right\|_{L_x^2\Sigma_{6}^*}
&\le C\left(\mathcal E_0[g]^{1/2}\nu^{-1/2}\mathcal D_1[f]^{1/2}
+\nu^{-1/2}\mathcal D_0[g]^{1/2}\mathcal E_1[f]^{1/2}\right).
\end{aligned}
\end{equation}
All constants are uniform for $1\le{\varpi}\le3/2$.
No second velocity derivative occurs on these right-hand sides.
\end{lem}

\begin{proof}
Start with smooth rapidly decaying $f,g$. The compact range of widths gives,
for a fixed $c>0$ and $|\alpha|\le1$,
\[
\sup_{1\le{\varpi}\le3/2}
 |\partial_v^\alpha\sqrt{\mu_{{\varpi}}}(v)|
\le C_\alpha e^{-c|v|^2}.
\]
The two uniform coefficient estimates needed for the divergence pairing are
\begin{equation}\label{eq:quadraticCoefficientBounds}
\begin{aligned}
|\overline X \sigma[\sqrt{\mu_{{\varpi}}}f](v)Y|
&\le C\|f\|_2(\overline X \sigma(v)X)^{1/2}
(\overline Y \sigma(v)Y)^{1/2},\\
|\sigma(v)^{-1/2}b[\sqrt{\mu_{{\varpi}}}f](v)|
&\le C\|f\|_{\Sigma_0}\langle v\rangle^{-1/2}.
\end{aligned}
\end{equation}
Here is a verification including the Coulomb singularity. For $|v|\le2$, Cauchy--Schwarz bounds the convolution defining {$\sigma[\sqrt{\mu_{{\varpi}}}f]$}, since $|v-v_*|^{-1}$ is locally square integrable and the other factor is Gaussian. On this set $\sigma$ is uniformly positive definite. For $|v|>2$, split $|v_*|\le |v|/2$ {from the complementary region.} With $e=v/|v|$, the first part has tangential size at most $C|v|^{-1}\|f\|_2$ and radial size at most
\[
C|v|^{-3}\int |v_*|^2\sqrt{\mu_{{\varpi}}}(v_*)|f(v_*)|\,dv_*
\le C|v|^{-3}\|f\|_2.
\]
In the complementary region the Gaussian and local square integrability give an exponentially small bound. Apply Cauchy--Schwarz to the positive matrix integral with $|f|$ to control mixed directions. The radial and tangential eigenvalue bounds for $\sigma$ prove the first assertion of \eqref{eq:quadraticCoefficientBounds}.

For the drift, the dissipation controls both $\langle v\rangle^{-3/2}\nabla f$ and $\langle v\rangle^{-1/2}f$. Hence, for a fixed sufficiently small $c>0$,
\[
\|e^{c|v|^2}\sqrt{\mu_{{\varpi}}} f\|_{H^1_v}
+\|\sqrt{\mu_{{\varpi}}} f\|_{L^1_v}
\le C\|f\|_{\Sigma_0}.
\]
On $|v-v_*|\le1$, use $H^1(\mathbb R^3)\subset L^6$ and
$|v-v_*|^{-2}\in L^{6/5}(|v-v_*|\le1)$.
On the complementary set use the $L^1$ bound. For large $|v|$, the part with $|v_*|\le|v|/2$ is bounded by $C|v|^{-2}\|f\|_{\Sigma_0}$, and the other part has Gaussian decay by the weighted $H^1$ bound. Thus {$|b[\sqrt{\mu_{{\varpi}}}f](v)|\le C\langle v\rangle^{-2}\|f\|_{\Sigma_0}$}. Since $\|\sigma^{-1/2}(v)\|\le C\langle v\rangle^{3/2}$, the second assertion follows.

The exact normalized divergence identity is
\begin{equation}\label{eq:quadraticDivergence}
\Gamma_{{\varpi}}(f,g)
=(\partial_{v_i}-\frac{v_i}{2{\varpi}})
\left[\sigma_{ij}[\sqrt{\mu_{{\varpi}}}f](\partial_{v_j}-\frac{v_j}{2{\varpi}})g-b_{i}[\sqrt{\mu_{{\varpi}}}f]g\right].
\end{equation}
{Using the derivatives in \eqref{eq:gaussianConjugatedDerivatives},
integration by parts gives}
\begin{equation}\label{eq:quadraticWeightedPairing}
\langle\Gamma_{{\varpi}}(f,g),\varphi\rangle_\ell
=-\int\big(\sigma_{ij}[\sqrt{\mu_{{\varpi}}}f]D_{{\varpi},j}^-g-b_{i}[\sqrt{\mu_{{\varpi}}}f]g\big)
\overline{D_{{\varpi},i}^+(\langle v\rangle^{2\ell}\varphi)}\,dv.
\end{equation}
Since $1/3\le \frac1{2{\varpi}}\le1/2$, the weight derivative and
the definition of $\Sigma_\ell$ imply, uniformly in ${\varpi}$,
\[
{
\|\sigma^{1/2}D_{{\varpi}}^-g\|_{L_\ell^2}\le C\|g\|_{\Sigma_\ell},\qquad \|\langle v\rangle^{-\ell}\sigma^{1/2}D_{{\varpi}}^+(\langle v\rangle^{2\ell}\varphi)\|_2 \le C_\ell\|\varphi\|_{\Sigma_\ell}.
}
\]
Consequently \eqref{eq:quadraticCoefficientBounds} yields
\[
{
\left|\langle\Gamma_{{\varpi}}(f,g),\varphi\rangle_\ell\right| \le C_\ell\|f\|_2\|g\|_{\Sigma_\ell}\|\varphi\|_{\Sigma_\ell} +C_\ell\|f\|_{\Sigma_0} \|\langle v\rangle^{\ell-1/2}g\|_2 \|\varphi\|_{\Sigma_\ell}.
}
\]
Taking the supremum over $\|\varphi\|_{\Sigma_\ell}\le1$ proves
\eqref{eq:quadraticVelocityForm}.

For the spatial estimate write $f_q=\partial_{x_1}^qf$ and
$g_h=\partial_{x_1}^hg$. Then
\[
\partial_{x_1}^j\Gamma_{{\varpi}}(f,g)
=\sum_{q=0}^j\binom jq\Gamma_{{\varpi}}(f_q,g_{j-q}),\qquad 0\le j\le r.
\]
For $q<r$, the velocity bound and H\"older's inequality in $x_1$ give
\begin{equation}\label{eq:quadraticSpatialLowerSplit}
\begin{aligned}
\|\Gamma_{{\varpi}}(f_q,g_{j-q})\|_{L^2_x\Sigma_{2r+6-2j}^*}
\le C_r\big(&\|f_q\|_{L^\infty_xL^2_v}
\|g_{j-q}\|_{L^2_x\Sigma_{2r+6-2j}}\\
&+\|f_q\|_{L^\infty_x\Sigma_0}
\|g_{j-q}\|_{L^2_xL^2_{2r+6-2j}}\big).
\end{aligned}
\end{equation}
Since $q+1\le r$ and $2r+6-2j+2q\ge 2r+6-2j$, the Hilbert-valued embedding
$H^1_x\subset L^\infty_x$ gives
\[
\begin{aligned}
\|f_q\|_{L^\infty_xL^2_v}&\le C_r\mathcal E_r[f]^{1/2},
&\|f_q\|_{L^\infty_x\Sigma_0}&\le C_r\nu^{-1/2}\mathcal D_r[f]^{1/2},\\
\|g_{j-q}\|_{L^2_x\Sigma_{2r+6-2j}}&\le\nu^{-1/2}\mathcal D_r[g]^{1/2},
&\|g_{j-q}\|_{L^2_xL^2_{2r+6-2j}}&\le\mathcal E_r[g]^{1/2}.
\end{aligned}
\]
For $q=r$ one necessarily has $j=r$. Reverse the spatial placement:
\begin{equation}\label{eq:quadraticSpatialTopSplit}
{
\|\Gamma_{{\varpi}}(f_r,g)\|_{L^2_x\Sigma_{6}^*} \le C_r\big(\|f_r\|_{L^2_xL^2_v}\|g\|_{L^\infty_x\Sigma_{6}} +\|f_r\|_{L^2_x\Sigma_0}\|g\|_{L^\infty_xL^2_{6}}\big).
}
\end{equation}
Here $r\ge1$ and $2r+6\ge2r+4\ge6$, so
\[
\begin{aligned}
\|g\|_{L^\infty_x\Sigma_{6}}
&\le C\sum_{h=0}^1\|g_h\|_{L^2_x\Sigma_{6}}
\le C_r\nu^{-1/2}\mathcal D_r[g]^{1/2},\\
\|g\|_{L^\infty_xL^2_{6}}
&\le C\sum_{h=0}^1\|g_h\|_{L^2_xL^2_{6}}
\le C_r\mathcal E_r[g]^{1/2}.
\end{aligned}
\]
Summing \eqref{eq:quadraticSpatialLowerSplit} and
\eqref{eq:quadraticSpatialTopSplit} over the finite set $q\le j\le r$
proves \eqref{eq:mixedSpatialLandauForm}; setting $f=g$ proves
\eqref{eq:quadraticSpatialForm}. Only spatial derivatives through $r$
occur.

For \eqref{eq:broaderGaussianZeroOrder}, the velocity weight is $6$ in
$\mathcal E_0$; in $\mathcal E_1$, the undifferentiated and differentiated
terms have weights $8$ and $6$, respectively.
The same Hilbert-valued spatial embedding gives
\begin{align*}
\|f\|_{L^\infty_xL^2_v}
 +\|\langle v\rangle^6f\|_{L^\infty_xL^2_v}
&\le C\mathcal E_1[f]^{1/2},\\
\|f\|_{L^\infty_x\Sigma_0}
 +\|f\|_{L^\infty_x\Sigma_6}
&\le C\nu^{-1/2}\mathcal D_1[f]^{1/2}.
\end{align*}
Apply \eqref{eq:quadraticVelocityForm} with $\ell=6$ in each ordering of its two factors:
\[
\begin{aligned}
\left\|\Gamma_{{\varpi}}(f,g)\right\|_{L_x^2\Sigma_{6}^*}
&\le C\big(\|f\|_{L^\infty_xL^2_v}\|g\|_{L^2_x\Sigma_6}
 +\|f\|_{L^\infty_x\Sigma_0}\|g\|_{L^2_xL^2_{6}}\big),\\
\left\|\Gamma_{{\varpi}}(g,f)\right\|_{L_x^2\Sigma_{6}^*}
&\le C\big(\|g\|_{L^2_xL^2_v}\|f\|_{L^\infty_x\Sigma_6}
 +\|g\|_{L^2_x\Sigma_0}\|f\|_{L^\infty_xL^2_{6}}\big).
\end{aligned}
\]
These inequalities prove \eqref{eq:broaderGaussianZeroOrder}, with no
spatial derivative on $g$. Approximation in the energy and form norms extends these
identities and bounds to the stated spaces.
\end{proof}

\section{Construction of the collision and prescribed-coefficient evolutions}\label{app:evolution}

\begingroup
{We construct the solutions used in the main text and prove their time
continuity and source estimates. We first treat the localized collision
equation without electric term in the density equation, then the prescribed-coefficient
equation in finite energy spaces. Finally, we bound the density produced
by sources estimated directly in the energy norm or against the collision
dissipation. The arguments use the coefficient and collision estimates
already established and require no nonlinear existence or instability result.}
\par\endgroup
\subsection{A common energy-limit argument}\label{app:commonEnergyLimits}

The constructions below have different transport operators and
regularizations. We first record the common passage from regularized
pairings to an energy identity. {Each construction will separately verify the}
transport limit; the collision and source limits then use this lemma.

\begin{lem}\label{lem:commonEnergyLimits}
Let $V\subset H\subset V^*$ be a dense Hilbert triple on a compact
time interval $I=[s_-,s_+]$, with $s_-<s_+$. If $z\in L^2(I;V)$ and
$\partial_sz\in L^2(I;V^*)$, then $z$ has a representative in
$C(I;H)$ and
\begin{equation}\label{eq:localizedTimeChainRule}
\frac12\|z(s_2)\|_H^2-\frac12\|z(s_1)\|_H^2
=\operatorname{Re}\int_{s_1}^{s_2}\langle\partial_sz,z\rangle\,ds.
\end{equation}
The bilinear identity follows by polarization.

Suppose instead that $z\in C_{\rm w}(I;H)\cap L^2(I;V)$, and let
$\mathfrak a_s$ be measurable sesquilinear forms satisfying
\begin{equation}\label{eq:commonPairingHypothesis}
|\mathfrak a_s(u,w)|\le C_{\rm form}\|u\|_V\|w\|_V
+c_{\rm mix}(s)(\|u\|_H\|w\|_V+\|u\|_V\|w\|_H)
+c_{\rm zer}(s)\|u\|_H\|w\|_H,
\end{equation}
where $C_{\rm form}<\infty$, $c_{\rm mix}\in L^2(I)$ and $c_{\rm zer}\in L^1(I)$ are nonnegative.
Let $w_n\to w$ in $L^2(I;V)$ and pointwise in $H$, with
$\sup_n\|w_n\|_{L^\infty H}<\infty$.
For $f_0\in L^1(I;H)$ and $f_1\in L^2(I;V^*)$,
\begin{equation}\label{eq:commonPairingLimit}
\int_I|\mathfrak a_s(z,w_n-w)|\,ds
+\int_I|\langle f_0+f_1,w_n-w\rangle|\,ds\longrightarrow0.
\end{equation}
The form convergence also holds {with the two form arguments interchanged.}
Consequently, if regularized energy identities have strongly
convergent endpoint values and these form and source pairings, their
limits give the energy identity for $z$. {If the right-hand side of that identity is}
integrable, then $z\in C(I;H)$.
\end{lem}

\begin{proof}

For the first assertion, time mollification after reflection across the
endpoints gives $z_\epsilon\to z$ in $L^2(I;V)$ and
$z_\epsilon'\to z'$ in $L^2(I;V^*)$. The inequality
\[
\|w\|_{C(I;H)}^2\le |I|^{-1}\|w\|_{L^2H}^2
+2\|w'\|_{L^2V^*}\|w\|_{L^2V}
\]
applied to $z_\epsilon-z_\delta$ gives convergence in $C(I;H)$;
Cauchy--Schwarz passes the smooth energy identity to
\eqref{eq:localizedTimeChainRule}. Polarization gives the bilinear identity.

For the second assertion, $z,w_n,w$ are uniformly bounded in
$L^\infty(I;H)$. The terms containing $\|w_n-w\|_V$ tend to zero by
\[
\int_I\big(C_{\rm form}\|z\|_V+c_{\rm mix}\|z\|_H+\|f_1\|_{V^*}\big)
\|w_n-w\|_V\,ds\longrightarrow0.
\]
The remaining integrands tend pointwise to zero and are bounded by
$C(c_{\rm mix}\|z\|_V+c_{\rm zer}+\|f_0\|_H)\in L^1(I)$.
Dominated convergence proves \eqref{eq:commonPairingLimit}, also with
the form arguments interchanged. The same limits on subintervals
pass the regularized energy identities to their endpoints.
The resulting integrable identity makes $\|z(s)\|_H^2$ continuous;
together with weak continuity this gives $z\in C(I;H)$.
\par

\end{proof}

\subsection{The localized collision evolution}\label{app:localizedEvolution}

We regularize the collision equation, prove estimates uniform in the approximation,
and pass to the energy solution. The same construction accommodates each spatial
weight and the centered derivative sums used in Section~\ref{sec:exactCollisionEvolution}.

\begin{lem}\label{lem:localizedConstruction}\label{app:localizedConstruction}
{Let $I=[\tau,T]$ be a finite interval, $\nu>0$, and $\ell\ge0$.}
Let {${F^{\rm c}}(s,x_1,v)$} be real, with
\[
{F^{\rm c}}_l(s,v)=\sqrt \mu (v)e^{-il(s-s_*)v_1}g_l(s-s_*,v),
\qquad \mathcal C(s)=-\mathcal L+\mathcal C_{F^{\rm c}}(s).
\]
Assume that $g_l$ is continuous in time and smooth in velocity and obeys
\eqref{eq:physicalCoefficientDerivatives}, with $t_0$ replaced by $s_*$,
for one fixed integer $J_v>\ell+8$. The constants may depend on this
fixed moment order. Let
\[
\mathcal Z_k(s)=\partial_{v_1}+i\big(k(s-s_*)-c\big),
\qquad a_n=\lambda_{\rm v}^n/g_\chi(n),
\]
where $s_*,c$ are fixed. For a derivative cutoff $J_{\rm der}\ge0$
and positive spatial weights $\kappa_k$, set
\begin{equation}\label{eq:constructionFiniteNorms}
\begin{aligned}
E_{J_{\rm der}}(s;f)&=\sum_k\kappa_k^2\sum_{n=0}^{J_{\rm der}} a_n^2
\|\langle v\rangle^\ell\mathcal Z_k(s)^nf_k\|_2^2,\\
D_{J_{\rm der}}(s;f)&=\sum_k\kappa_k^2\sum_{n=0}^{J_{\rm der}} a_n^2
\|\mathcal Z_k(s)^nf_k\|_{\Sigma_\ell}^2.
\end{aligned}
\end{equation}
Assume that, for every finite derivative sum and spatial projection, the coefficient satisfies
\begin{equation}\label{eq:constructionFormHypothesis}
\begin{split}
\left|\sum_k\kappa_k^2\sum_{n=0}^{J_{\rm der}}a_n^2
\langle\mathcal Z_k^n(\mathcal C_{F^{\rm c}} f)_k,
\mathcal Z_k^nw_k\rangle_{\ell}\right|
\le C\varepsilon_{\rm c}\big(&D_{J_{\rm der}}(f)^{1/2}D_{J_{\rm der}}(w)^{1/2}\\
&+c_{\rm col}(s)E_{J_{\rm der}}(f)^{1/2}E_{J_{\rm der}}(w)^{1/2}\big),
\end{split}
\end{equation}
where $c_{\rm col}\ge1$ is bounded on $I$, and $C$ is independent of the projections, ${J_{\rm der}},c,s_*$. Assume also the corresponding zeroth-order estimate for each prescribed polynomial spatial weight, with constants allowed to depend on that weight. The spatial weights used here are $e^{\varrho_0|k-j_0|^\chi}$ with
$\varrho_0<\varrho_{\rm b}$, or prescribed polynomial weights
$\langle k\rangle^r$ or $\langle k/N\rangle^r$. Their ratio under a mode
shift is absorbed by the coefficient spatial decay. The fixed moment
$J_v$ is not increased with the number of velocity derivatives.

{Fix a spatial order $r\ge0$, and let
$f_\tau\in H^r_{x_1}L^2_\ell$ and
$S\in L^1(I;H^r_{x_1}L^2_\ell)$.
For sufficiently small $\varepsilon_{\rm c}$ and $\lambda_{\rm v}$,
the initial-value problem}
\begin{equation}\label{eq:constructionEquation}
\partial_sf+v_1\partial_{x_1}f=\nu\mathcal C(s)f+S,
\qquad f(\tau)=f_\tau,
\end{equation}
{admits a unique energy solution $f\in C(I;H^r_{x_1}L^2_\ell)$
with finite integrated collision dissipation.
Let $E_\infty,D_\infty$ denote the nonnegative sums in
\eqref{eq:constructionFiniteNorms}. If, in addition,
\[
E_\infty(\tau;f_\tau)^{1/2}
+\int_\tau^T E_\infty(s;S)^{1/2}\,ds<\infty,
\]
then}
\begin{equation}\label{eq:constructionEnergy}
\begin{aligned}
&E_\infty(s;f)^{1/2} \le e^{C\nu\int_\tau^s(1+\varepsilon_{\rm c} c_{\rm col}(\tau'))\,d\tau'}
\left(E_\infty(\tau;f_\tau)^{1/2}
+\int_\tau^sE_\infty(\tau';S)^{1/2}\,d\tau'\right),\\
&\sup_{\tau\le s'\le s}E_\infty(s';f)
+\nu\int_\tau^sD_\infty(\tau';f)\,d\tau'\\
&\quad\le C e^{C\nu\int_\tau^s(1+\varepsilon_{\rm c} c_{\rm col}(\tau'))\,d\tau'}
\left(E_\infty(\tau;f_\tau)^{1/2}
+\int_\tau^sE_\infty(\tau';S)^{1/2}\,d\tau'\right)^2.
\end{aligned}
\end{equation}
The constants are independent of $c,{J_{\rm der}}$ and of $0<\lambda_{\rm v}\le\lambda_*$, for a sufficiently small fixed $\lambda_*$. For the spatial weights used below they are also independent of the spatial mode $j$ of the datum. The homogeneous solutions define an evolution family, and the inhomogeneous solution is {given by the associated Duhamel integral.}
\end{lem}

\begin{proof}

\emph{1. Regularization and velocity Galerkin construction.}
Let $P_{J_x}$ project onto $|k|\le J_x$ and set
$\mathcal H_v=-\Delta_v+|v|^2$. For finite-mode Schwartz data and
smooth Schwartz sources, first solve
\begin{equation}\label{eq:constructionRegularization}
\partial_sf^{J_x,\epsilon}+v_1\partial_{x_1}f^{J_x,\epsilon}
=\nu P_{J_x}\mathcal C(s)P_{J_x}f^{J_x,\epsilon}
-\epsilon\mathcal H_vf^{J_x,\epsilon}+P_{J_x}S,
\qquad f^{J_x,\epsilon}(\tau)=P_{J_x}f_\tau .
\end{equation}
At fixed $J_x$ and $\epsilon>0$, use the finite product of
$H_\ell=L^2_\ell$ and the reference form domain
\[
\|y\|_{V_\ell}^2
=\|\langle v\rangle^\ell y\|_2^2
+\|\langle v\rangle^\ell\nabla_vy\|_2^2
+\|\langle v\rangle^\ell v y\|_2^2.
\]
The embedding $V_\ell\Subset H_\ell$ is compact. Spectral Galerkin
projections for this form give
\[
\sup_I\|f_n\|_{H_\ell}^2+\epsilon\int_I\|f_n\|_{V_\ell}^2\,ds
+\|\partial_sf_n\|_{L^2(I;V_\ell^*)}^2\le C_{J_x,\epsilon,I}.
\]
Here transport is skew and bounded $V_\ell\to V_\ell^*$ at fixed
$J_x$, and the regularizer supplies coercivity. Weak compactness
and Lemma~\ref{lem:commonEnergyLimits} give the energy solution of
\eqref{eq:constructionRegularization}; the difference energy estimate
gives uniqueness. The next steps justify differentiation and obtain
bounds independent of $J_x,\epsilon$.
\par

\emph{2. Coefficients for the difference-quotient equation.} We next justify finite differentiations of this weak solution.
This step uses the specified Landau coefficient class, in addition to the
{form inequality. Put ${F^{\rm c,tot}}=\mu +{F^{\rm c}}$.}
The first collision ordering has the exact local expression
\begin{equation}\label{eq:constructionLocalOperator}
\mu^{-1/2}Q({F^{\rm c,tot}},\sqrt \mu  y)
=(\partial_{v_i}-v_i/2)
\big(\sigma_{iq}[{F^{\rm c,tot}}](\partial_{v_q}-v_q/2)y-b_{i}[{F^{\rm c,tot}}]y\big).
\end{equation}
For small $\varepsilon_{\rm c}$, Gaussian matrix comparison gives {$\sigma[{F^{\rm c,tot}}]\ge c\sigma$}.
This remains true for the quadratic form after spatial projection:
the solution is in the range of $P_{J_x}$, so $P_{J_x}$ can be removed from
{the pairing with that solution.}

Expanding \eqref{eq:constructionLocalOperator} gives
\begin{equation}\label{eq:constructionLocalCoefficients}
\begin{aligned}
\mathcal B_{F^{\rm c,tot}}y&=\operatorname{div}_v(\sigma[{F^{\rm c,tot}}]\nabla_vy)
-(\sigma[{F^{\rm c,tot}}]v+b[{F^{\rm c,tot}}])\cdot\nabla_vy+c^{F^{\rm c,tot}}y,\\
c^{F^{\rm c,tot}}&=-\tfrac12\operatorname{div}_v(\sigma[{F^{\rm c,tot}}]v)-\operatorname{div}_vb[{F^{\rm c,tot}}]
+\tfrac14v\cdot \sigma[{F^{\rm c,tot}}]v+\tfrac12v\cdot b[{F^{\rm c,tot}}].
\end{aligned}
\end{equation}
All fixed $v_1$ derivatives of these coefficients are bounded on $I$.
Here is the velocity bound needed for this assertion. The coefficient
hypothesis makes {${F^{\rm c,tot}}$}, every fixed $v_1$ derivative, {and the required}
additional velocity derivatives of order at most two
Gaussian-decaying, uniformly on $I$. Splitting the Coulomb convolution
into $|v_*|\le |v|/2$, a neighborhood of $v_*=v$, and the remaining
Gaussian tail gives
\[
|\partial_{v_1}^h \sigma[{F^{\rm c,tot}}]|\le C_{h,I}\langle v\rangle^{-1-h},\qquad
|\partial_{v_1}^h b[{F^{\rm c,tot}}]|\le C_{h,I}\langle v\rangle^{-2-h}.
\]
In the neighborhood of the singularity, derivatives are transferred
to {${F^{\rm c,tot}}$}; outside it the kernel is smooth and homogeneous. Cutoff
derivatives are supported in the Gaussian tail. For the undifferentiated
potential, matrix comparison also gives
{$|v_i\sigma_{iq}[{F^{\rm c,tot}}]v_q|\le C\langle v\rangle^{-1}$}.
For $h\ge1$ the displayed derivative bounds control the factors of
degree at most two in $v$. They prove the claimed boundedness of all
differentiated drift and multiplication coefficients. At fixed $J_x,I$
their spatial norms are finite as well. These auxiliary constants need
not be uniform in derivative order, $J_x$, or $I$.

The second ordering has a nonlocal part with kernel
\begin{equation}\label{eq:constructionNonlocalKernel}
K_{F^{\rm c,tot}}(s,x_1;v,v_*)
=\mu (v)^{-1/2}\partial_{v_iv_q}{F^{\rm c,tot}}(s,x_1,v)
\Phi_{iq}(v-v_*)\sqrt \mu (v_*),
\end{equation}
and the multiplication term {$8\pi {F^{\rm c,tot}}y$}. {For the corresponding integral operator}
{$\mathcal K_{F^{\rm c,tot}}$}, integration by parts in $v_{*1}$ gives
\begin{equation}\label{eq:constructionKernelDifferentiation}
\partial_{v_1}^h(\mathcal K_{F^{\rm c,tot}} y)
=\sum_{j=0}^h\binom hj\mathcal K_{{F^{\rm c,tot}},j}\partial_{v_1}^{h-j}y,
\qquad K_{{F^{\rm c,tot}},j}=(\partial_{v_1}+\partial_{v_{*1}})^jK_{F^{\rm c,tot}}.
\end{equation}
The simultaneous derivative annihilates $\Phi(v-v_*)$. The coefficient
{hypothesis consequently bounds the kernel after conjugation by $\langle v\rangle^\ell$}
by a finite constant $C_{j,I}$ times
\[
\langle v\rangle^{\ell-J_v+2}|v-v_*|^{-1}
e^{-c|v_*|^2}\langle v_*\rangle^{-\ell}.
\]
Both Schur integrals are finite because $J_v>\ell+8$. Thus
\begin{equation}\label{eq:constructionFiniteKernelBound}
\|\langle v\rangle^\ell\partial_{v_1}^h\mathcal K_{F^{\rm c,tot}}y\|_2
\le C_{h,I}\sum_{j=0}^h\|\langle v\rangle^\ell\partial_{v_1}^jy\|_2.
\end{equation}
The multiplication term has the same estimate. In particular this
step does not require additional moments as $h$ increases.

\emph{3. Finite velocity regularity.}
At fixed $J_x,\epsilon$, write $y_h=\partial_{v_1}^hf^{J_x,\epsilon}$,
\[
E_h=\sum_{|k|\le J_x}\|y_{h,k}\|_{L_\ell^2}^2,\qquad
\mathfrak R_h=\sum_{|k|\le J_x}
\big(\|\nabla_vy_{h,k}\|_{L_\ell^2}^2+\|v y_{h,k}\|_{L_\ell^2}^2\big).
\]
Step 1 starts the induction in $h$. Set
$\mathcal T_\eta y(v)=y(v+\eta e_1)$,
$\Delta_\eta=(\mathcal T_\eta-{\mathrm{Id}})/\eta$ and
$z_\eta=\Delta_\eta y_{h-1}$, with $0<|\eta|\le1$.
The local principal part satisfies
\begin{equation}\label{eq:constructionDifferencePrincipal}
\Delta_\eta\operatorname{div}_v(\sigma[{F^{\rm c,tot}}]\nabla_vy)
=\operatorname{div}_v\big((\mathcal T_\eta\sigma[{F^{\rm c,tot}}])\nabla_v\Delta_\eta y
+(\Delta_\eta\sigma[{F^{\rm c,tot}}])\nabla_vy\big).
\end{equation}
For $m=h-1$, {the differentiated lower-order terms are}
\[
\begin{aligned}
&\partial_{v_1}^m(\mathcal B_{F^{\rm c,tot}}f)-\mathcal B_{F^{\rm c,tot}}y_m =\sum_{j=1}^m\binom mj\big[\operatorname{div}_v((\partial_{v_1}^j\sigma[F^{\rm c,tot}])\nabla_vy_{m-j})\\
&\qquad\qquad-\partial_{v_1}^j(\sigma[F^{\rm c,tot}]v+b[F^{\rm c,tot}])\cdot\nabla_vy_{m-j}
+(\partial_{v_1}^jc^{F^{\rm c,tot}})y_{m-j}\big].
\end{aligned}
\]
Let $E_\eta,\mathfrak R_\eta$ be the preceding energy and regularization
forms evaluated at $z_\eta$. Step 2 and the fundamental theorem of
calculus give, for $0\le j<h$,
\[
\begin{aligned}
\Delta_\eta\partial_{v_1}^j\sigma[F^{\rm c,tot}]
&=\int_0^1\mathcal T_{\xi\eta}
\partial_{v_1}^{j+1}\sigma[F^{\rm c,tot}]\,d\xi,\\
\sup_{0<|\eta|\le1}
\|\Delta_\eta\partial_{v_1}^j\sigma[F^{\rm c,tot}]\|_{L^\infty_{s,x,v}}
&\le\|\partial_{v_1}^{j+1}\sigma[F^{\rm c,tot}]\|_{L^\infty_{s,x,v}}
\le C_{h,J_x,I}.
\end{aligned}
\]
The highest-order coefficient difference in
\eqref{eq:constructionDifferencePrincipal} therefore satisfies
\[
\begin{aligned}
&\nu\left|\left\langle
\nabla_v\!\cdot\!\big((\Delta_\eta\sigma[F^{\rm c,tot}])
\nabla_v y_{h-1}\big),z_\eta\right\rangle_\ell\right|\\
&\quad=\nu\left|\int_{\mathbb T\times\mathbb R^3}
\langle v\rangle^{2\ell}
(\Delta_\eta\sigma[F^{\rm c,tot}])\nabla_v y_{h-1}\cdot
\left(\nabla_v\overline{z_\eta}
+\frac{2\ell v}{\langle v\rangle^2}\overline{z_\eta}\right)\,dx_1\,dv\right|\\
&\quad\le C_{h,J_x,I}\nu\mathfrak R_{h-1}^{1/2}
\big(\mathfrak R_\eta^{1/2}+E_\eta^{1/2}\big) \le\frac{\epsilon}{16}\mathfrak R_\eta
+C_{h,J_x,\epsilon,I}(E_\eta+\mathfrak R_{h-1}).
\end{aligned}
\]
The spatial projections can be omitted from this pairing because
they are orthogonal and commute with the velocity weight.
For $0\le p<h-1$, bounded translations in the polynomial weight give
\[
\Delta_\eta y_p=\int_0^1\mathcal T_{\xi\eta}y_{p+1}\,d\xi,\qquad \|\Delta_\eta y_p\|_{L_x^2L_\ell^2}
+\|\nabla_v\Delta_\eta y_p\|_{L_x^2L_\ell^2}
\le C_\ell\big(E_{p+1}^{1/2}+\mathfrak R_{p+1}^{1/2}\big).
\]
Apply the difference product rule to the remaining divergence terms
in the displayed commutator and use these bounds. Since
$p=h-1-j$ implies $p+1=h-j<h$ for $j\ge1$, their sum satisfies
\[
\begin{aligned}
&\nu\sum_{j=1}^{h-1}\binom{h-1}{j}
\left|\left\langle\Delta_\eta\nabla_v\!\cdot\!
\big((\partial_{v_1}^j\sigma[F^{\rm c,tot}])
\nabla_vy_{h-1-j}\big),z_\eta\right\rangle_\ell\right|\\
&\quad\le\frac{\epsilon}{16}\mathfrak R_\eta
+C_{h,J_x,\epsilon,I}
\left(E_\eta+\sum_{j<h}(E_j+\mathfrak R_j)\right).
\end{aligned}
\]
The differentiated drift and multiplication coefficients are bounded
by Step 2; the same product rule bounds their pairings by
$C_{h,J_x,I}(E_\eta+\sum_{j<h}(E_j+\mathfrak R_j))$.
Every constant here is independent of the difference step $\eta$.
\par

For the nonlocal term, simultaneous translation of the two kernel
variables gives the exact identity
\begin{equation}\label{eq:constructionDifferenceKernel}
\begin{aligned}
\Delta_\eta(\mathcal K_{F^{\rm c,tot}}y)(v)
={}&\int K_{F^{\rm c,tot}}(v+\eta e_1,v_*+\eta e_1)\Delta_\eta y(v_*)\,dv_*\\
&+\int_0^1\int (\partial_{v_1}+\partial_{v_{*1}})
K_{F^{\rm c,tot}}(v+\xi\eta e_1,v_*+\xi\eta e_1)y(v_*)\,dv_*d\xi.
\end{aligned}
\end{equation}
The Coulomb kernel is unchanged. The Schur bound
\eqref{eq:constructionFiniteKernelBound}, {also applied to the differentiated}
kernels, therefore gives
\[
\|\langle v\rangle^\ell\Delta_\eta\partial_{v_1}^{h-1}(\mathcal K_{F^{\rm c,tot}}f)\|_2
\le C_{h,I}\left(\|\langle v\rangle^\ell z_\eta\|_2+\sum_{j<h}E_j^{1/2}\right).
\]
The remaining transport and regularizer commutators are
$[\Delta_\eta,ikv_1]=ik\mathcal T_\eta$,
$[\Delta_\eta,\mathcal H_v]=(2v_1+\eta)\mathcal T_\eta$ and
$[\partial_{v_1}^m,\mathcal H_v]f
=2mv_1y_{m-1}+m(m-1)y_{m-2}$, with negative indices zero.
They give, together with the positive regularizer,
\begin{equation}\label{eq:constructionDifferenceEnergy}
E_\eta'+c\epsilon\mathfrak R_\eta
\le C_{h,J_x,\epsilon,I}\left(E_\eta+\sum_{j<h}(E_j+\mathfrak R_j)\right)
+C_\ell E_\eta^{1/2}\|\langle v\rangle^\ell\partial_{v_1}^hP_{J_x}S\|_2.
\end{equation}
Lemma~\ref{lem:commonEnergyLimits} justifies this pairing at fixed
$\eta$. Gronwall and weak compactness as $\eta\to0$ construct $y_h$
with $\sup_I E_h+\epsilon\int_I\mathfrak R_h<\infty$ {and the corresponding}
differentiated equation. {The same velocity weight $\langle v\rangle^\ell$ suffices at every}
finite order. These auxiliary constants may depend on
$h,J_x,\epsilon,I$; the next step proves the uniform centered estimate.
\par

\emph{4. The uniform centered estimate.} Put $y_{k,n}=\mathcal Z_k^nf_k^{J_x,\epsilon}$. The spatial projection commutes with $\mathcal Z_k$, and
$[\partial_s+ikv_1,\mathcal Z_k]=0$. The added velocity operator has the exact commutator
\begin{equation}\label{eq:constructionRegularizerCommutator}
\mathcal Z_k^n\mathcal H_v f_k
=\mathcal H_v y_{k,n}
+2n v_1y_{k,n-1}+n(n-1)y_{k,n-2},
\end{equation}
where negative derivative indices mean zero. In particular this commutator is independent of $c$. Weighted integration by parts gives
\[
\operatorname{Re}\langle\mathcal H_v y,y\rangle_{\ell}
\ge \tfrac12\|\langle v\rangle^\ell\nabla_vy\|_2^2
+\|\langle v\rangle^\ell v y\|_2^2
-C_\ell\|\langle v\rangle^\ell y\|_2^2.
\]
The factorial coefficients obey, uniformly in $n$,
\[
n\,a_n/a_{n-1}\le C\lambda_{\rm v},\qquad
n(n-1)a_n/a_{n-2}\le C\lambda_{\rm v}^2.
\]
{Take the weighted inner product of
\eqref{eq:constructionRegularizerCommutator} with
${\kappa_k^2a_n^2y_{k,n}}$ and sum in ${k,n}$. With negative indices set to zero,}
\begin{align*}
2\epsilon\sum_{k,n\le {J_{\rm der}}}\kappa_k^2a_n^2n
|\langle v_1y_{k,n-1},y_{k,n}\rangle_{\ell}|
&\le C\epsilon\lambda_{\rm v}
\left(\sum_{k,n\le {J_{\rm der}}}\kappa_k^2a_n^2\|v y_{k,n}\|_{L_\ell^2}^2\right)^{1/2}E_{J_{\rm der}}^{1/2}\\
&\le\tfrac14\epsilon \mathfrak R_{J_{\rm der}}+C\epsilon\lambda_{\rm v}^2E_{J_{\rm der}},\\
\epsilon\sum_{k,n\le {J_{\rm der}}}\kappa_k^2a_n^2n(n-1)
|\langle y_{k,n-2},y_{k,n}\rangle_{\ell}|
&\le C\epsilon\lambda_{\rm v}^2E_{J_{\rm der}}.
\end{align*}
Here $\mathfrak R_{J_{\rm der}}$ is the positive regularization expression below. Both estimates are independent of $c,J_x,{J_{\rm der}}$; only $\lambda_{\rm v}\le\lambda_*$ is used.

The Maxwellian estimate and \eqref{eq:constructionFormHypothesis} consequently yield
\begin{equation}\label{eq:constructionUniformEstimate}
\frac12 E_{J_{\rm der}}'+c\nu D_{J_{\rm der}}+c\epsilon \mathfrak R_{J_{\rm der}}
\le C\big(\nu[1+\varepsilon_{\rm c} c_{\rm col}(s)]+\epsilon\big)E_{J_{\rm der}}
+E_{J_{\rm der}}^{1/2}E_{J_{\rm der}}(S)^{1/2},
\end{equation}
where
\[
\mathfrak R_{J_{\rm der}}=\sum_k\kappa_k^2\sum_{n=0}^{J_{\rm der}}a_n^2
\left(\|\langle v\rangle^\ell\nabla_v y_{k,n}\|_2^2
+\|\langle v\rangle^\ell v y_{k,n}\|_2^2\right).
\]
Choose $\varepsilon_{\rm c}$ small enough to {absorb the relative dissipation contribution from {$\mathcal C_{F^{\rm c}}$}.} Gronwall gives an estimate uniform in $J_x,\epsilon,{J_{\rm der}}$ except for the harmless factor $e^{C\epsilon(s-\tau)}$. The same calculation at derivative order zero with polynomial spatial weights gives the energy bounds needed for the construction independently of centered regularity.

\emph{5. Removal of the regularization and strong time traces.}
Use the polynomial spatial weights $\langle k\rangle^r$ and write
\[
\|f\|_H^2=\sum_k\langle k\rangle^{2r}\|\langle v\rangle^\ell f_k\|_2^2,
\qquad \|f\|_V^2=\|f\|_H^2+\sum_k\langle k\rangle^{2r}\|f_k\|_{\Sigma_\ell}^2.
\]
Constants used only to justify the following limits may depend on fixed
$\nu,I,r$; the a priori estimate \eqref{eq:constructionUniformEstimate}
{retains the stated uniformity.}
The collision form $\mathfrak c_s(f,w)=\nu\langle\mathcal C(s)f,w\rangle_H$
satisfies, by Maxwellian boundedness and the zeroth-order coefficient hypothesis,
\begin{equation}\label{eq:localizedZerothFormContinuity}
|\mathfrak c_s(f,w)|
\le C\nu\left(\|f\|_V\|w\|_V+
\varepsilon_{\rm c} c_{\rm col}(s)\|f\|_H\|w\|_H\right),\qquad c_{\rm col}\in L^\infty(I).
\end{equation}
This extends the form from smooth functions to $V\times V$.

First let $\epsilon\downarrow0$ at fixed $J_x$, and then let $J_x\to\infty$.
The zeroth-order energy bound gives, along subsequences,
\[
f^{J_x,\epsilon}\overset{*}{\rightharpoonup}f^{J_x}
\text{ in }L^\infty(I;H),\quad
f^{J_x,\epsilon}\rightharpoonup f^{J_x}\text{ in }L^2(I;V),
\]
and the same convergences from $f^{J_x}$ to $f$ in the second limit.
For a smooth compactly supported test $\varphi$,
\begin{align*}
\left|\epsilon\int_I\langle f^{J_x,\epsilon},\mathcal H_v\varphi\rangle\,ds\right|
&\le C_\varphi\epsilon\|f^{J_x,\epsilon}\|_{L^\infty H}\longrightarrow0,\\
\left|\int_I\mathfrak c_s(f^{J_x},(P_{J_x}-{\mathrm{Id}})\varphi)\,ds\right|
&\le C_{I,\nu}\|f^{J_x}\|_{L^2V}
\|(P_{J_x}-{\mathrm{Id}})\varphi\|_{L^2V}\longrightarrow0.
\end{align*}
The remaining collision pairing passes by weak convergence in $L^2V$.
Transport is tested as $-\langle f,v_1\partial_{x_1}\varphi\rangle_H$;
{the test coefficient is} smooth with compact velocity support. Thus the limit
solves \eqref{eq:constructionEquation}, {including the prescribed initial term.}
For each fixed such test, and every measurable time set $\mathfrak I\subset I$, the weak equation gives
\[
\int_{\mathfrak I}|\partial_s(f,\varphi)_H|\,ds
\le C_\varphi\left(|\mathfrak I|^{1/2}\|f\|_{L^2V}
+|\mathfrak I|\|f\|_{L^\infty H}+\|\mathbf1_{\mathfrak I} S\|_{L^1H}\right).
\]
Applied first to the approximations, this bound yields a weakly continuous
representative $f\in C_{\rm w}(I;H)$ with $f(\tau)=f_\tau$.

{We derive the energy identity for $f$} without treating transport as a bounded
map $V\to V^*$. Fix ${K_{\rm cut}}<\infty$, put $f_{K_{\rm cut}}=P_{K_{\rm cut}}f$, and set
\[
(\mathcal W_{K_{\rm cut}}(s)f_{K_{\rm cut}})_k=e^{ik(s-\tau)v_1}(f_{K_{\rm cut}})_k,
\qquad z_{K_{\rm cut}}=\mathcal W_{K_{\rm cut}}f_{K_{\rm cut}}.
\]
On these finitely many modes,
\begin{align*}
\|\mathcal W_{K_{\rm cut}}(s)y\|_H&=\|y\|_H,\\
\nabla_v(\mathcal W_{K_{\rm cut}}(s)y)_k
&=e^{ik(s-\tau)v_1}\big(\nabla_vy_k+ik(s-\tau)e_1y_k\big),\\
\|\mathcal W_{K_{\rm cut}}(s)^{\pm1}y\|_V&\le C_{{K_{\rm cut}},I}\|y\|_V.
\end{align*}
Let $V_{K_{\rm cut}}=P_{K_{\rm cut}}V$, $H_{K_{\rm cut}}=P_{K_{\rm cut}}H$. The equation becomes
\[
\partial_sz_{K_{\rm cut}}=\nu\mathcal W_{K_{\rm cut}}P_{K_{\rm cut}}\mathcal C(s)f+\mathcal W_{K_{\rm cut}}P_{K_{\rm cut}}S
\in L^2(I;V_{K_{\rm cut}}^*),\qquad z_{K_{\rm cut}}\in L^2(I;V_{K_{\rm cut}}).
\]
Here the source is still the smooth source chosen in Step 1. Applying
\eqref{eq:localizedTimeChainRule} in $V_{K_{\rm cut}}\subset H_{K_{\rm cut}}\subset V_{K_{\rm cut}}^*$ gives,
{for every $\tau\le s_1\le s_2$ in $I$,}
\begin{equation}\label{eq:localizedProjectedEnergyIdentity}
\frac12\|P_{K_{\rm cut}}f(s_2)\|_H^2-\frac12\|P_{K_{\rm cut}}f(s_1)\|_H^2
=\operatorname{Re}\int_{s_1}^{s_2}
\left(\mathfrak c_s(f,P_{K_{\rm cut}}f)+(S,P_{K_{\rm cut}}f)_H\right)\,ds.
\end{equation}
The endpoint values coincide with the previously constructed weak traces.
The projections satisfy
\[
P_{K_{\rm cut}}f\to f\text{ in }L^2(I;V),\qquad
P_{K_{\rm cut}}f(s)\to f(s)\text{ in }H,\qquad
\sup_{K_{\rm cut}}\|P_{K_{\rm cut}}f\|_{L^\infty H}\le\|f\|_{L^\infty H}.
\]
By \eqref{eq:localizedZerothFormContinuity}, Lemma~\ref{lem:commonEnergyLimits}
applies with $c_{\rm mix}=0$ and
$c_{\rm zer}=C\nu\varepsilon_{\rm c}c_{\rm col}\in L^\infty(I)$.
Passing to the limit in \eqref{eq:localizedProjectedEnergyIdentity} gives,
for every pair of endpoints,
\begin{equation}\label{eq:localizedExactEnergyIdentity}
\frac12\|f(s_2)\|_H^2-\frac12\|f(s_1)\|_H^2
=\operatorname{Re}\int_{s_1}^{s_2}
\left(\mathfrak c_s(f,f)+(S,f)_H\right)\,ds.
\end{equation}
Lemma~\ref{lem:commonEnergyLimits} gives $f\in C(I;H)$.
Applied to a difference, this identity and Gronwall prove uniqueness,
so every approximation subsequence has the same limit.

General data and $S\in L^1(I;H)$ follow by smooth approximation and
the difference estimate
\begin{equation}\label{eq:localizedSourceCompletion}
\begin{aligned}
&\|f_n-f_m\|_{C(I;H)}
+\left(\nu\int_I\sum_k\langle k\rangle^{2r}
\|(f_n-f_m)_k\|_{\Sigma_\ell}^2\,ds\right)^{1/2}\\
&\qquad\le C\big(\|f_{\tau,n}-f_{\tau,m}\|_H+\|S_n-S_m\|_{L^1(I;H)}\big).
\end{aligned}
\end{equation}
Completion and Lemma~\ref{lem:commonEnergyLimits} preserve
\eqref{eq:localizedExactEnergyIdentity} and the strong time traces.
\par

For each fixed ${J_{\rm der}}$, the centered norm $E_{J_{\rm der}}$ is a finite-order weighted
Sobolev norm after a finite spatial projection. Approximate the initial
datum in $E_{J_{\rm der}}(\tau;\cdot)$ and the source in $L^1_sE_{J_{\rm der}}(s;\cdot)^{1/2}$.
The zeroth-order convergence is \eqref{eq:localizedSourceCompletion};
distributional differentiation identifies all weak centered-derivative limits.
At each fixed time, weak energy convergence and the uniform $E_{J_{\rm der}}$
bound identify the same centered derivatives by testing in velocity.
Lower semicontinuity therefore gives the $E_{J_{\rm der}}$ estimate at every time,
as well as the integrated $D_{J_{\rm der}}$ estimate in
\eqref{eq:constructionEnergy}. Uniqueness {identifies the resulting solution with} $f$.
Finally,
\[
E_{J_{\rm der}}(s;f)\uparrow E_\infty(s;f),\qquad
\int_\tau^sD_{J_{\rm der}}(t;f)\,dt\uparrow\int_\tau^sD_\infty(t;f)\,dt.
\]
Monotone convergence proves \eqref{eq:constructionEnergy}. No density
assertion in an infinite derivative norm is needed.

\emph{6. Evolution and Duhamel formula.}
Uniqueness gives $\mathbb U(s,r)\mathbb U(r,\tau)=\mathbb U(s,\tau)$
and $\mathbb U(\tau,\tau)={\mathrm{Id}}$. The energy bound gives local uniform
boundedness, and for finite-mode Schwartz $y$,
\[
\|\mathbb U(s,\tau)y-y\|_H
\le C_I\int_\tau^s\|(\nu\mathcal C(t)-v_1\partial_{x_1})y\|_H\,dt
\longrightarrow0\quad(s-\tau\to0).
\]
Density and the composition law give joint strong continuity.
For $S\in L^1(I;H)$, the energy bound and Fubini's theorem verify
the Duhamel integral in the weak equation; uniqueness yields
\[
f(s)=\mathbb U(s,\tau)f_\tau+\int_\tau^s\mathbb U(s,r)S(r)\,dr.
\]
For $\nu=0$ the formula follows from free transport; the dissipation
is asserted {only after multiplication by $\nu$.}
\par

\end{proof}

\begingroup
\subsection{Energy solutions with prescribed coefficients}\label{sec:prescribedSolvability}

We construct the prescribed-coefficient energy solutions used in the
complete linear equation, the quantitative nonlinear comparison and
the identification of local classical solutions. The proof derives
the energy estimate and establishes existence, uniqueness and strong
time continuity. The collision matrix remains $\sigma=\Phi*\mu$.

For $1\le\varpi\le3/2$ and a physical distribution $U$, use
\begin{equation}\label{eq:gaussianEnergyPullback}
\|U\|_{\mathcal N_{r,\varpi}}
=\|\mu_{\varpi}^{-1/2}U\|_{\mathcal N_r},\qquad
\|U\|_{\mathcal N_{r,1}}=\|U/\sqrt\mu\|_{\mathcal N_r}.
\end{equation}
The corresponding source conditions are
\begin{equation}\label{eq:broaderGaussianNorms}
\int_I\mathcal E_r[\mu_{\varpi}^{-1/2}S_0]^{1/2}\,ds
+\nu^{-1/2}\sum_{j=0}^r
\|\partial_{x_1}^j(\mu_{\varpi}^{-1/2}S_1)\|_{L^2_{s,x}(I)\Sigma_{2r+6-2j}^*}<\infty.
\end{equation}
The scalar $\varpi$ changes only the normalization, not the
physical Maxwellian $\mu$ or the collision matrix $\sigma=\Phi*\mu$.
At $\varpi=1$ these are the source conditions in
\eqref{eq:quadraticSourceSpaces}. The local-solution identification in
Section~\ref{sec:classicalConclusion} uses $\varpi=5/4$;
the quantitative estimates use $\varpi=1$.
\par\endgroup

\begin{lem}\label{lem:physicalDerivativeEnergy}\label{app:physicalDerivativeEnergy}
{Fix a finite interval $I=[t_0,T]$, $\nu>0$,
$1\le{\varpi}\le3/2$ and an integer $r\ge0$.} Set $r_{\rm c}=\max\{r,1\}$. Let {${F^{\rm c}_1},{F^{\rm c}_2}$} be prescribed real distributions and put
\[
f_i=\mu_{{\varpi}}^{-1/2}({F^{\rm c}_i}-\mu ),\qquad i=1,2.
\]
{Assume}
\begin{equation}\label{eq:prescribedGaussianCoefficients}
\sup_I\sum_{i=1}^2\mathcal E_{r_{\rm c}}[f_i]^{1/2}\le\epsilon_*,
\qquad
\int_I\sum_{i=1}^2\mathcal D_{r_{\rm c}}[f_i]\,ds<\infty.
\end{equation}
Let $E=-\partial_{x_1}\phi$ be real, with mean-zero potential, and assume
\[
E\in L^\infty(I;C^{r_{\rm c}}_x),\qquad
\phi\in W^{1,1}(I;C^r_x)\cap L^\infty(I;C^{r+1}_x).
\]
Let $\mathcal J_s=\mathcal J_s^0+\mathcal J_s^1$ be a prescribed linear
operator satisfying
\begin{equation}\label{eq:prescribedBoundedResponse}
\begin{aligned}
\mathcal E_r[\mu_{{\varpi}}^{-1/2}\mathcal J_s^0U]^{1/2}
&\le c_{\mathcal J,0}(s)\mathcal E_r[\mu_{{\varpi}}^{-1/2}U]^{1/2},
&c_{\mathcal J,0}&\in L^1(I),\\
\sum_{j=0}^{r}\left\|\partial_{x_1}^j\!\left(\mu_{{\varpi}}^{-1/2}\mathcal J_s^1U\right)\right\|_{L_x^2\Sigma_{2r+6-2j}^*}
&\le c_{\mathcal J,1}(s)\mathcal E_r[\mu_{{\varpi}}^{-1/2}U]^{1/2},
&c_{\mathcal J,1}&\in L^2(I),
\end{aligned}
\end{equation}

\par
{Assume that $s\mapsto\mathcal J_s f$ is measurable for each fixed $f$; $\mathcal J_s=0$ is
allowed. Let $U_{\rm in}$ and $S_0,S_1$ satisfy
\[
\mu_{{\varpi}}^{-1/2}U_{\rm in}\in\mathscr E_r,
\qquad
\int_{I}\mathcal E_r[\mu_{{\varpi}}^{-1/2}S_0]^{1/2}\,ds
+\nu^{-1/2}\sum_{j=0}^{r}\left\|\partial_{x_1}^j\!\left(\mu_{{\varpi}}^{-1/2}S_1\right)\right\|_{L^2_{s,x}(I)\Sigma_{2r+6-2j}^*}<\infty,
\]
{with the source conditions in \eqref{eq:broaderGaussianNorms}.}
For $\epsilon_*$ sufficiently small depending on the fixed order,
weight and potential bound, the initial-value problem}
\begin{equation}\label{eq:translatedDerivativeEquation}
(\partial_s+v_1\partial_{x_1}+E\partial_{v_1})U
=\nu\big(Q({F^{\rm c}_1},U)+Q(U,{F^{\rm c}_2})\big)
+\mathcal J_sU+S_0+S_1
\end{equation}
{with $U(t_0)=U_{\rm in}$ admits a unique energy solution $U$.
This solution satisfies}
\begin{equation}\label{eq:physicalDerivativeEnergyBound}
\begin{aligned}
\|U\|_{\mathcal N_{r,{\varpi}}}
\le C\Bigg(&\mathcal E_r[\mu_{{\varpi}}^{-1/2}U(t_0)]^{1/2}
+\int_I\mathcal E_r[\mu_{{\varpi}}^{-1/2}S_0]^{1/2}\,ds\\
&+\nu^{-1/2}\sum_{j=0}^r
\|\partial_{x_1}^j(\mu_{{\varpi}}^{-1/2}S_1)\|_{L^2_{s,x}(I)\Sigma_{2r+6-2j}^*}\Bigg).
\end{aligned}
\end{equation}
The normalized solution is strongly continuous in $\mathscr E_r$. The constant is finite and uniform on bounded sets of the stated coefficient norms, with a fixed smallness margin in \eqref{eq:prescribedGaussianCoefficients}. There is no zero-mass condition. 
\end{lem}

\begingroup
\begin{proof}
We first prove the energy inequality for smooth solutions, then construct
a weak solution and justify the inequality {at the regularity of that solution.}
Use $\Gamma_{\varpi}$
from Lemma~\ref{lem:quadraticLandauForm}. {The estimates}
\eqref{eq:mixedSpatialLandauForm} for $r\ge1$ and
\eqref{eq:broaderGaussianZeroOrder} for $r=0$ supply all bilinear
collision bounds below, uniformly for $1\le\varpi\le3/2$.

The coercive part remains the physical Maxwellian {$\mu $}. Since {$b[\mu ]=-\sigma v$},
\[
\mu_{\varpi}^{-1/2}Q(\mu ,\sqrt{\mu_{\varpi}}f)
=(\partial_{v_i}-\frac{v_i}{2\varpi})
\left[\sigma_{ij}(\partial_{v_j}+(1-\frac1{2\varpi})v_j)f\right].
\]
{Integration of the outer derivative and then of the
two mixed terms gives the exact identity}
\[
\begin{aligned}
\operatorname{Re}\langle \mu_{\varpi}^{-1/2}Q(\mu ,\sqrt{\mu_{\varpi}}f),f\rangle_{\ell}
={}&-\|\sigma^{1/2}\nabla f\|_{L_\ell^2}^2
-\frac{2\varpi-1}{4\varpi^2}\|\sigma^{1/2}vf\|_{L_\ell^2}^2\\
&-\operatorname{Re}\int(\nabla\langle v\rangle^{2\ell})^{\mathsf T}\sigma\nabla f\,\bar f
+\frac12\langle\operatorname{div}(\sigma v)f,f\rangle_\ell\\
&+\left(\frac1{2\varpi}-\frac12\right)\int\nabla\langle v\rangle^{2\ell}\cdot \sigma v\,|f|^2.
\end{aligned}
\]
The Maxwellian matrix estimates and the fixed polynomial weight imply
\[
\frac{4\ell^2|\sigma^{1/2}v|^2}{\langle v\rangle^4}
+|\operatorname{div}(\sigma v)|
+\frac{2\ell|v\cdot\sigma v|}{\langle v\rangle^2}\le C_\ell.
\]
Consequently,
\begin{align*}
\left|\int(\nabla\langle v\rangle^{2\ell})^{\mathsf T}\sigma\nabla f\,\bar f\right|
&\le \|\sigma^{1/2}\nabla f\|_{L_\ell^2}\left(\int\frac{|\sigma^{1/2}\nabla\langle v\rangle^{2\ell}|^2}{\langle v\rangle^{2\ell}}|f|^2\right)^{1/2}
\\
&\le\tfrac14\|\sigma^{1/2}\nabla f\|_{L_\ell^2}^2+C_\ell\|f\|_{L_\ell^2}^2,\\
\frac1{2\varpi}\in[1/3,1/2]&\quad\Longrightarrow\quad \frac{2\varpi-1}{4\varpi^2}\ge2/9.
\end{align*}
It follows that
\[
\operatorname{Re}\langle \mu_{\varpi}^{-1/2}Q(\mu ,\sqrt{\mu_{\varpi}}f),f\rangle_{\ell}
\le-c\|f\|_{\Sigma_\ell}^2+C_\ell\|f\|_{L_\ell^2}^2.
\]

For the other Maxwellian ordering, cancellation gives
\[
\mu_{\varpi}^{-1/2}Q(\sqrt{\mu_{\varpi}}f,\mu )
=\int K_{\varpi}(v,v_*)f(v_*)\,dv_*+8\pi \mu f,
\]
where
\[
K_{\varpi}(v,v_*)
=\frac{\partial_{v_iv_j}\mu (v)}{\sqrt{\mu_{\varpi}(v)}}
\Phi_{ij}(v-v_*)\sqrt{\mu_{\varpi}(v_*)}.
\]
{For $1\le\varpi\le3/2$, conjugation by $\langle v\rangle^\ell$ gives, for a}
fixed $c>0$ and $C_\ell$ uniform in this interval,
\[
|\langle v\rangle^\ell K_{\varpi}(v,v_*)\langle v_*\rangle^{-\ell}|
\le C_\ell e^{-c|v|^2}|v-v_*|^{-1}e^{-c|v_*|^2}.
\]
Indeed {$\mu /\sqrt{\mu_{\varpi}}=C_{\varpi}
e^{-(1/2-1/(4\varpi))|v|^2}$}, and the exponent is at least $1/4$.
For any fixed $d>3$, $e^{-c|v|^2}\le C_{c,d}\langle v\rangle^{-d}$.
Equation \eqref{eq:weightedCoulombSchur} therefore gives
\[
\left\|\int_{\mathbb R^3}
K_{\varpi}(v,v_*)f(v_*)\,dv_*\right\|_{L_\ell^2}
\le C_\ell\|f\|_{L_\ell^2}.
\]
The multiplication term obeys the same weighted $L^2$ bound. Applying
the first-ordering coercivity and this bound to each spatial derivative,
with $\ell=2r+6-2j$, yields
\[
\begin{aligned}
&\operatorname{Re}\sum_{j=0}^r
\left\langle\partial_{x_1}^j \mu_{\varpi}^{-1/2}
\big(Q(\mu ,\sqrt{\mu_{\varpi}}f)+Q(\sqrt{\mu_{\varpi}}f,\mu )\big),
\partial_{x_1}^jf\right\rangle_{2r+6-2j}\\
&\hspace{20mm}\le-c_r\nu^{-1}\mathcal D_r[f]+C_r\mathcal E_r[f].
\end{aligned}
\]
Equations \eqref{eq:mixedSpatialLandauForm} and \eqref{eq:broaderGaussianZeroOrder} bound the remaining collision terms by
\[
C_{r,{\varpi}}\sum_{i=1}^2
\left(\mathcal E_{r_{\rm c}}[f_i]^{1/2}\mathcal D_r[f]
+\mathcal D_{r_{\rm c}}[f_i]^{1/2}
\mathcal E_r[f]^{1/2}\mathcal D_r[f]^{1/2}\right).
\]
Choose $\epsilon_*$ to absorb the first terms. Young's inequality absorbs another fixed fraction of the dissipation and leaves
${C\sum_i\mathcal D_{r_{\rm c}}[f_i]\mathcal E_r[f]}$.

To treat the field, put
\[
f=\mu_{\varpi}^{-1/2}U,\qquad
y=e^{\phi/(2\varpi)}f.
\]
Define the normalized operators by
\[
\begin{aligned}
\mathcal C_{\varpi} f
&=\mu_{\varpi}^{-1/2}\big(Q(F^{\rm c}_1,\sqrt{\mu_{\varpi}}f)
+Q(\sqrt{\mu_{\varpi}}f,F^{\rm c}_2)\big),\\
\widehat{\mathcal J}^{\,i}f
&=\mu_{\varpi}^{-1/2}\mathcal J_s^i(\sqrt{\mu_{\varpi}}f),
\qquad i=0,1,\qquad
\widehat{\mathcal J}=\widehat{\mathcal J}^{\,0}+\widehat{\mathcal J}^{\,1}.
\end{aligned}
\]
\par

Lemma~\ref{lem:potentialTransport} gives
\begin{equation}\label{eq:prescribedPotentialConjugation}
\begin{aligned}
&e^{\phi/(2\varpi)}\mathcal C_{\varpi}
(e^{-\phi/(2\varpi)}y)=\mathcal C_{\varpi}y,\qquad (\partial_s+v_1\partial_{x_1}+E\partial_{v_1})y\\
&\quad=\nu\mathcal C_{\varpi}y+\frac{\partial_s\phi}{2\varpi}y
+e^{\phi/(2\varpi)}\widehat{\mathcal J}(e^{-\phi/(2\varpi)}y) +e^{\phi/(2\varpi)}\mu_{\varpi}^{-1/2}(S_0+S_1).
\end{aligned}
\end{equation}
The first identity uses locality in $x$ of each collision ordering.
The collision estimates above therefore apply to $y$ without a
commutator. The common field estimate \eqref{eq:generalFieldEnergy}
{and the spatial commutator calculation} \eqref{eq:fieldCommutatorPairing}
bound the transport contribution by
\[
\eta\mathcal D_r[y]
+C_{r,\eta}\big(\|E\|_\infty+
\nu^{-1}\|E\|_{C^r_x}^2\big)\mathcal E_r[y],\qquad \eta>0.
\]
This also covers $r=0$, where the spatial commutator sum is empty.
Lemma~\ref{lem:weightedDualMultipliers} and
\eqref{eq:prescribedBoundedResponse} yield
\[
\begin{aligned}
\mathcal E_r[\frac{\partial_s\phi}{2\varpi}y
+e^{\phi/(2\varpi)}\widehat{\mathcal J}^{\,0}(e^{-\phi/(2\varpi)}y)]^{1/2}
&\le C\big(\|\partial_s\phi\|_{C^r_x}+c_{\mathcal J,0}(s)\big)\mathcal E_r[y]^{1/2},\\
\left|\sum_{j=0}^r\left\langle
\partial_{x_1}^j\{e^{\phi/(2\varpi)}
\widehat{\mathcal J}^{\,1}(e^{-\phi/(2\varpi)}y)\},\partial_{x_1}^jy\right\rangle_{2r+6-2j}\right|
&\le Cc_{\mathcal J,1}(s)\mathcal E_r[y]^{1/2}\nu^{-1/2}\mathcal D_r[y]^{1/2}\\
&\le\eta\mathcal D_r[y]+C_\eta\nu^{-1}c_{\mathcal J,1}(s)^2\mathcal E_r[y],\\
\mathcal E_r[e^{\phi/(2\varpi)}\mu_{\varpi}^{-1/2}S_0]^{1/2}
&\le C\mathcal E_r[\mu_{\varpi}^{-1/2}S_0]^{1/2},\\
\sum_{j=0}^{r}\left\|\partial_{x_1}^j\!\left(e^{\phi/(2\varpi)}\mu_{\varpi}^{-1/2}S_1\right)\right\|_{L_x^2\Sigma_{2r+6-2j}^*}
&\le C\sum_{j=0}^{r}\left\|\partial_{x_1}^j\!\left(\mu_{\varpi}^{-1/2}S_1\right)\right\|_{L_x^2\Sigma_{2r+6-2j}^*}.
\end{aligned}
\]

\par

The same lemma gives equivalence of the energy and form norms of $f$ and
$y$. Choose $\eta$ smaller than a fixed fraction of the retained collision
coercivity constant in \eqref{eq:fieldCommutatorPairing}.

In particular an admissible integrable energy coefficient is
\begin{equation}\label{eq:mixedEnergyCoefficient}
c_{\rm en}(s)=C\bigg(1+\nu+\sum_{i=1}^2\mathcal D_{r_{\rm c}}[f_i]
+\nu^{-1}\|E\|_{C^{r_{\rm c}}_x}^2 +\|\partial_s\phi\|_{C^r_x}+c_{\mathcal J,0}(s)
+\nu^{-1}c_{\mathcal J,1}(s)^2\bigg).
\end{equation}

\par

We have proved
\begin{equation}\label{eq:prescribedGaussianEnergy}
\begin{aligned}
\frac{d}{ds}\mathcal E_r[y]+c\mathcal D_r[y]
\le{}&c_{\rm en}(s)\mathcal E_r[y]
+C\mathcal E_r[\mu_{\varpi}^{-1/2}S_0(s)]^{1/2}\mathcal E_r[y]^{1/2}\\
&+C\nu^{-1}\sum_{j=0}^r
\|\partial_{x_1}^j(\mu_{\varpi}^{-1/2}S_1(s))\|_{L_x^2\Sigma_{2r+6-2j}^*}^2.
\end{aligned}
\end{equation}

Integrating \eqref{eq:prescribedGaussianEnergy} with the factor
$\exp(-\int_{t_0}^sc_{\rm en})$ and applying Young's inequality to
the $S_0$ term gives
\[
\begin{aligned}
&\sup_{I}\mathcal E_r[y]^{1/2}
+\left(\int_I\mathcal D_r[y]\,ds\right)^{1/2} \le Ce^{\int_Ic_{\rm en}}
\Bigl(\mathcal E_r[y(t_0)]^{1/2}\\
&\qquad+\int_{I}\mathcal E_r[\mu_{\varpi}^{-1/2}S_0]^{1/2}\,ds +\nu^{-1/2}\sum_{j=0}^{r}\left\|\partial_{x_1}^j\!\left(\mu_{\varpi}^{-1/2}S_1\right)\right\|_{L^2_{s,x}(I)\Sigma_{2r+6-2j}^*}\Bigr).
\end{aligned}
\]
The potential multiplier bounds yield
\eqref{eq:physicalDerivativeEnergyBound} for smooth solutions.
\par

For the construction, use the form domain
\[
V_r=\{y\in\mathscr E_r:\mathcal D_r[y]<\infty\},\qquad
\|y\|_{V_r}^2=\mathcal E_r[y]+\nu^{-1}\mathcal D_r[y].
\]
Project the conjugated equation onto nested spans of smooth periodic
functions with compact velocity support, orthonormal in $\mathscr E_r$
and dense in $V_r$ and in the smooth test topology. Each finite
system has integrable coefficients. {The energy estimate}
\eqref{eq:prescribedGaussianEnergy} gives
\[
y_n\overset{*}{\rightharpoonup}y\text{ in }L^\infty(I;\mathscr E_r),
\qquad y_n\rightharpoonup y\text{ in }L^2(I;V_r).
\]
The limit is taken against fixed smooth tests. For such a test $\varphi$,
integration by parts gives
\begin{equation}\label{eq:fixedTestTransportPairing}
\begin{aligned}
\langle-(v_1\partial_{x_1}+E\partial_{v_1})z,\varphi\rangle_{\mathscr E_r}
={}&\sum_{j=0}^r\int\langle v\rangle^{4r+12-4j}\partial_{x_1}^jz\,
v_1\partial_{x_1}^{j+1}\overline\varphi\,dx_1dv\\
&+\sum_{j=0}^r\sum_{q=0}^j\binom jq
\int\partial_{x_1}^{j-q}z\,(\partial_{x_1}^qE)
\partial_{v_1}(\langle v\rangle^{4r+12-4j}\partial_{x_1}^j\overline\varphi)\,dx_1dv.
\end{aligned}
\end{equation}
Compact velocity support bounds this by
${C_\varphi(1+\|E\|_{C^r_x})\mathcal E_r[z]^{1/2}}$.
The collision form satisfies
\[
|\langle\nu\mathcal C_{\varpi}z,\varphi\rangle_{\mathscr E_r}|
\le C_\varphi\nu\left(1+\sum_i\mathcal E_{r_{\rm c}}[f_i]^{1/2}\right)\|z\|_{V_r}
+C_\varphi\sqrt\nu\sum_i\mathcal D_{r_{\rm c}}[f_i]^{1/2}\mathcal E_r[z]^{1/2}.
\]
Thus the collision pairing defines an $L^2(I;V_r^*)$ functional;
for each fixed smooth test, the transport, potential and bounded-electric-operator
pairings have $L^1(I;\mathscr E_r^*)$ coefficients. The additional electric term satisfies
\[
\sum_{j=0}^{r}\left\|\partial_{x_1}^j\!\left(\widehat{\mathcal J}^{\,1}(e^{-\phi/(2\varpi)}y_n)\right)\right\|_{L^2_{s,x}(I)\Sigma_{2r+6-2j}^*}
\le C\|c_{\mathcal J,1}\|_{L^2(I)}\sup_I\mathcal E_r[y_n]^{1/2}.
\]
The adjoint of this electric operator applied to a fixed smooth test lies in
$L^1(I;\mathscr E_r^*)$, so weak-star convergence passes this term as well.
The two weak convergences pass all pairings separately.
\par
 Their
absolute continuity on measurable time sets also yields weak time
traces, as in Step 5 of Lemma~\ref{lem:localizedConstruction}.
Consequently $y\in C_{\rm w}(I;\mathscr E_r)$ with the prescribed
initial datum. This uses neither a bounded transport map
$V_r\to V_r^*$ nor bounded Galerkin projections on $V_r$.
\par

To recover the physical weak equation from the differentiated pairing,
let $\psi$ be a smooth test compactly supported in velocity. The formula
\[
\widehat\varphi_k(v)=
\frac{\widehat\psi_k(v)}{\displaystyle\sum_{j=0}^r
\langle v\rangle^{4r+12-4j}k^{2j}}
\]
defines a test with the same properties, and ${\langle y,\varphi\rangle_{\mathscr E_r}=\langle y,\psi\rangle}$. The limiting equation is therefore the conjugated distributional equation. Multiplication by {$e^{-\phi/(2\varpi)}\sqrt{\mu_{\varpi}}$} recovers \eqref{eq:translatedDerivativeEquation}.

It remains to justify the energy identity for this limit. Set
$z_j=\langle v\rangle^{2r+6-2j}\partial_{x_1}^jy$ and let $J_\epsilon$
be an even mollifier in $(x_1,v)$. On each velocity ball,
\begin{equation}\label{eq:kineticFriedrichsCommutator}
\begin{aligned}
&[v_1\partial_{x_1}+E\partial_{v_1},J_\epsilon]z(x_1,v)\\
&=\int_{\mathbb R\times\mathbb R^3}
\Big[w_1\partial_{y_1}\eta_\epsilon(y_1,w)
+\big(E(x_1)-E(x_1-y_1)\big)\partial_{w_1}\eta_\epsilon(y_1,w)\Big]\\
&\qquad\times\big(z(x_1-y_1,v-w)-z(x_1,v)\big)\,dy_1dw,\\
&\|[v_1\partial_{x_1}+E\partial_{v_1},J_\epsilon]z\|_2
\le C(1+\|\partial_{x_1}E\|_\infty)\|z\|_2,\\
&[v_1\partial_{x_1}+E\partial_{v_1},J_\epsilon]z\longrightarrow0
\quad\hbox{in }L^2.
\end{aligned}
\end{equation}
The subtracted term vanishes because $\partial_{x_1}v_1+\partial_{v_1}E=0$;
translation continuity proves convergence.
Since $E\in L^\infty C^1_x$, the bound is time integrable.

Choose $0\le\psi_L\le1$ equal to one on $|v|\le L$ and zero on
$|v|\ge2L$, with $|\nabla\psi_L|\le C/L$. Set
\[
z_{j,\epsilon,L}=J_\epsilon(\psi_Lz_j),\qquad
\psi_{j,\epsilon,L}=\psi_LJ_\epsilon^2(\psi_Lz_j).
\]
At fixed $L,\epsilon$, the mollified equation has an $L^1(I;L^2)$
right-hand side. {The weak time traces agree with} those already constructed.
Self-adjointness of $J_\epsilon$ therefore gives the exact identity
\begin{equation}\label{eq:prescribedRegularizedEnergyIdentity}
\frac12\sum_{j=0}^r
\left(\|z_{j,\epsilon,L}(s_2)\|_2^2-
\|z_{j,\epsilon,L}(s_1)\|_2^2\right)
=\operatorname{Re}\sum_{j=0}^r\int_{s_1}^{s_2}
\langle\partial_sz_j,\psi_{j,\epsilon,L}\rangle\,ds.
\end{equation}
To pass to the limit, apply the differentiated equation in each component.
{The collision, field-commutator, potential and electric pairings have the common bound}
\begin{equation}\label{eq:energyPairingContinuity}
\begin{aligned}
&C\mathcal D_r[y]^{1/2}\mathcal D_r[z]^{1/2}
+c_{\rm zer}(s)\mathcal E_r[y]^{1/2}\mathcal E_r[z]^{1/2}\\
&\quad+c_{\rm mix}(s)\left(\mathcal E_r[y]^{1/2}\nu^{-1/2}\mathcal D_r[z]^{1/2}
+\nu^{-1/2}\mathcal D_r[y]^{1/2}\mathcal E_r[z]^{1/2}\right),\\
&\hspace{25mm} c_{\rm mix}\in L^2(I),\qquad c_{\rm zer}\in L^1(I).
\end{aligned}
\end{equation}

One can take, with the fixed potential multiplier bound included in $C$,
\begin{align*}
c_{\rm mix}(s)&=C\left(\sqrt\nu\sum_i\mathcal D_{r_{\rm c}}[f_i(s)]^{1/2}
+\|E(s)\|_{C^{r_{\rm c}}_x}+{c_{\mathcal J,1}(s)}\right),\\
c_{\rm zer}(s)&=C\left(\nu+\|E(s)\|_\infty+
\|\partial_s\phi(s)\|_{C^r_x}+{c_{\mathcal J,0}(s)}\right).
\end{align*}
Then $c_{\rm mix}\in L^2(I)$ and $c_{\rm zer}\in L^1(I)$ at fixed $\nu$.
This is \eqref{eq:commonPairingHypothesis} for the finite vector
of differentiated components, with energy norm $\mathcal E_r^{1/2}$
{and form norm ${(\mathcal E_r+\nu^{-1}\mathcal D_r)^{1/2}}$.}
The source coefficients have bounds {${\mathcal E_r[\mu_{\varpi}^{-1/2}S_0]^{1/2}\in L^1(I)}$} in the energy
dual and {${\sum_{j=0}^{r}\left\|\partial_{x_1}^j\!\left(\mu_{\varpi}^{-1/2}S_1\right)\right\|_{L_x^2\Sigma_{2r+6-2j}^*}\in L^2(I)}$} in the form dual.

Use the finite vector with components
$\langle v\rangle^{-2r-6+2j}\psi_{j,\epsilon,L}-\partial_{x_1}^jy$, $0\le j\le r$.
The proof of \eqref{eq:energyPairingContinuity} is a sum of component
pairings, so it applies to this vector without requiring the components
to be derivatives of one regularized function. Local comparability of
$\sigma(v)$ and $\sigma(v+z)$ for $|z|\le1$ bounds the mollifiers on these form spaces.
For the cutoff derivative,
\[
\|\langle v\rangle^{2r+6-2j} \sigma^{1/2}\nabla\psi_L\,\partial_{x_1}^jy\|_2^2
\le\frac{C}{L^2}\int_{|v|\ge L}
\langle v\rangle^{4r+11-4j}|\partial_{x_1}^jy|^2\,dx_1dv.
\]
The form and energy tails therefore give, successively as $\epsilon\to0$
and $L\to\infty$,
\[
\begin{aligned}
&\int_I\sum_{j=0}^r
\|\langle v\rangle^{-2r-6+2j}\psi_{j,\epsilon,L}-\partial_{x_1}^jy\|_{L^2_x\Sigma_{2r+6-2j}}^2\,ds\to0,\\
&\sum_{j=0}^r\|\psi_{j,\epsilon,L}(s)-\langle v\rangle^{2r+6-2j}\partial_{x_1}^jy(s)\|_2^2
\to0\quad\hbox{for every }s,\\
&\sum_{j=0}^r\|\psi_{j,\epsilon,L}(s)-\langle v\rangle^{2r+6-2j}\partial_{x_1}^jy(s)\|_2^2
\le C\mathcal E_r[y(s)].
\end{aligned}
\]
The constants are uniform for $L\ge1$ and $0<\epsilon\le1$.
By \eqref{eq:commonPairingLimit}, all nontransport and source pairings
converge in this product space, first as $\epsilon\to0$ and then as
$L\to\infty$. No mollifier is commuted through a collision coefficient.

For transport, \eqref{eq:kineticFriedrichsCommutator} gives a vanishing
commutator as $\epsilon\to0$ at fixed $L$. {The remaining transport boundary term is}
bounded by
\[
\frac{C\|E\|_{L^\infty(I\times\mathbb T)}}{L}
\int_I\sum_{j=0}^r\|\mathbf1_{|v|\ge L}z_j(s)\|_2^2\,ds
\longrightarrow0.
\]
At each endpoint, $J_\epsilon(\psi_Lz_j(s_i))\to z_j(s_i)$ in $L^2$
in the same order of limits; this only uses $y(s_i)\in \mathscr E_r$, not strong
time continuity. Hence \eqref{eq:prescribedRegularizedEnergyIdentity}
passes to the integrated energy identity between arbitrary endpoints.
Lemma~\ref{lem:commonEnergyLimits} now gives $y\in C(I;\mathscr E_r)$
from the integrable energy identity and the previously fixed weak traces.
This proves \eqref{eq:prescribedGaussianEnergy} for the actual solution. Applying it to differences gives uniqueness and continuous dependence in the energy supremum and integrated form norms, hence solvability for the completed source spaces.
\end{proof}
\par\endgroup
\subsection{\texorpdfstring{{Density estimates for inhomogeneous collision equations}}{Density estimates for inhomogeneous collision equations}}\label{app:sourceTraces}

The same spatially localized energy estimate controls initial data,
integrable $L^2_v$ sources and sources estimated by the collision form.
We prove their density bounds together, keeping the time supremum
inside the sum over the spatial indices $k$.

\begingroup
\begin{lem}\label{lem:collisionSourceEstimates}
Use \eqref{eq:comparisonLogParameters} with $p>{\rm s}_0+7/2$.
Fix an integer $r\ge2$ and $I'=[\tau_*,\tau^*]\subset I$.
{Suppose that ${\mathcal E_r[f_{\rm in}]<\infty}$ and that the
normalized sources $S_0,S_1$ satisfy
\[
\int_{I'}\mathcal E_r[S_0(s)]^{1/2}\,ds
+\sum_{j=0}^r\|\partial_{x_1}^jS_1\|_{L^2_{s,x}(I')\Sigma_{2r+6-2j}^*}<\infty.
\]
}
Let $f$ be an energy solution of
\begin{equation}\label{eq:collisionSourceForward}
(\partial_s+v_1\partial_{x_1})f
=\nu\mathcal C_{\rm b}(s)f+S_0+S_1,
\qquad f(\tau_*)=f_{\rm in}.
\end{equation}
Then
\begin{equation}\label{eq:collisionSourceEnergy}
\|f\|_{\mathcal N_r^{I'}}
\le C_r\big(\mathcal E_r[f_{\rm in}]^{1/2}
+\int_{I'}\mathcal E_r[S_0]^{1/2}\,ds
+\nu^{-1/2}\sum_{j=0}^{r}\left\|\partial_{x_1}^jS_1\right\|_{L^2_{s,x}(I')\Sigma_{2r+6-2j}^*}\big).
\end{equation}
For every such datum and sources a unique solution
$f\in C(I';\mathscr E_r)$ with finite integrated collision form norm
exists. Its density satisfies the sharper density estimate
\begin{equation}\label{eq:collisionCombinedTrace}
\begin{aligned}
\sum_{k\ne0}\langle k\rangle^{r-1}\sup_{s\in I'}
\left|\int\sqrt\mu f_k(s)\,dv\right|
\le C_r\big(&\|f_{\rm in}\|_{H_x^rL_v^2}
+\|S_0\|_{L^1(I';H_x^rL_v^2)}\\
&+\nu^{-1/2}\sum_{j=0}^{r}\left\|\partial_{x_1}^jS_1\right\|_{L^2_{s,x}(I')\Sigma_{2r+6-2j}^*}\big).
\end{aligned}
\end{equation}
If only $f_{\rm in}\in H_x^rL_v^2$ and
$S_0\in L^1(I';H_x^rL_v^2)$ are assumed, with
{$S_1$ satisfying the {same source condition} as above, then}
\eqref{eq:collisionSourceForward} has a unique solution in
$C(I';H_x^rL_v^2)$ with finite integrated collision form norm.
This solution satisfies \eqref{eq:collisionCombinedTrace} and
\begin{equation}\label{eq:collisionUnweightedEnergy}
\begin{aligned}
&\sup_{s\in I'}\|f(s)\|_{H_x^rL_v^2}
+\left(\nu\sum_{a=0}^r\int_{I'}
\|\partial_{x_1}^af(s)\|_{L_x^2\Sigma_0}^2\,ds\right)^{1/2}\\
&\qquad\le C_r\big(\|f_{\rm in}\|_{H_x^rL_v^2}
+\|S_0\|_{L^1(I';H_x^rL_v^2)}
+\nu^{-1/2}\sum_{j=0}^{r}\left\|\partial_{x_1}^jS_1\right\|_{L^2_{s,x}(I')\Sigma_{2r+6-2j}^*}\big).
\end{aligned}
\end{equation}
For the homogeneous equation, any real $q\ge0$, $\tau\in I$ and
$S\in H_x^{q+1}L_v^2$ satisfy
\par
\begin{equation}\label{eq:coarseCollisionTrace}
\sum_{k\ne0}\langle k\rangle^q\sup_{\tau\le s\le T_1}
\left|\int\sqrt\mu[\mathbb U_{\rm b}^\nu(s,\tau)S]_k\,dv\right|
\le C_q\sum_j\langle j\rangle^q\|S_j\|_2
\le C_q\|S\|_{H_x^{q+1}L_v^2}.
\end{equation}
{For a normalized source $S$ satisfying ${\sum_{j=0}^r\|\partial_{x_1}^jS\|_{L^2_{s,x}(I')\Sigma_{2r+6-2j}^*}<\infty}$, the solution with zero initial data satisfies}
\begin{equation}\label{eq:dualCollisionTrace}
\sum_{k\ne0}\langle k\rangle^{r-1}\sup_{s\in I'}
\left|\int\sqrt\mu\left[\int_{\tau_*}^s\mathbb U_{\rm b}^\nu(s,\tau)S(\tau)\,d\tau\right]_k\,dv\right|
\le C_r\nu^{-1/2}\sum_{j=0}^{r}\left\|\partial_{x_1}^jS\right\|_{L^2_{s,x}(I')\Sigma_{2r+6-2j}^*}.
\end{equation}
Duhamel expressions for form sources mean limits
of solutions with smooth sources. All distribution modes are retained;
the constants depend only on fixed orders.
\end{lem}
\par\endgroup

\begingroup
\begin{proof}
\emph{1. Energy solutions.}
Apply Lemma~\ref{lem:physicalDerivativeEnergy} with
$\varpi=1$, $F_1^{\rm c}=F_2^{\rm c}=F^{\rm b,\nu}$,
$E=0$, $\mathcal J=0$ and physical sources $\sqrt\mu S_0,\sqrt\mu S_1$.
Proposition~\ref{lem:fixedCoefficientOrders} verifies the coefficients.
The collision form and Young's inequality give
\[
\begin{aligned}
\partial_s\mathcal E_r[f]+c_r\mathcal D_r[f]
\le{}&C_r\nu[1+\varepsilon_N(1+(s-t_0)^2)]\mathcal E_r[f]\\
&+C_r\mathcal E_r[S_0(s)]^{1/2}\mathcal E_r[f]^{1/2}
+C_r\nu^{-1}\sum_{j=0}^{r}\left\|\partial_{x_1}^j\!\left(S_1(s)\right)\right\|_{L_x^2\Sigma_{2r+6-2j}^*}^2,\\
\nu\int_{I'}[1+\varepsilon_N(1+(s-t_0)^2)]\,ds
&\le C(\nu T_1+\nu\varepsilon_N(1+T_1)^3)=o(1).
\end{aligned}
\]
An integrating factor proves \eqref{eq:collisionSourceEnergy};
source completion and uniqueness in that lemma give the stated solution.

\emph{2. One spatial mode $j$ for $f_{\rm in}$, $S_0$ and $S_1$.}
Initially take smooth data and sources with finitely many modes.
For each $j$, let $f^{[j]}$ solve \eqref{eq:collisionSourceForward}
with datum $\mathbf e_j f_{{\rm in},j}$ and sources
$\mathbf e_j S_{0,j}$ and $\mathbf e_j S_{1,j}$.
Fix a sufficiently small $\varrho_0>0$ below the background spatial
radius. The weight equals one at mode $j$.
The zero-velocity-derivative collision form at weight zero gives
\[
\begin{aligned}
&\frac{d}{ds}\sum_k e^{2\varrho_0|k-j|^\chi}\|f_k^{[j]}\|_2^2
+c\nu\sum_k e^{2\varrho_0|k-j|^\chi}\|f_k^{[j]}\|_{\Sigma_0}^2\\
&\quad\le C\nu[1+\varepsilon_N(1+(s-t_0)^2)]
\sum_k e^{2\varrho_0|k-j|^\chi}\|f_k^{[j]}\|_2^2
+2\|S_{0,j}\|_2\|f_j^{[j]}\|_2
+C\nu^{-1}\|S_{1,j}\|_{\Sigma_0^*}^2.
\end{aligned}
\]
Here $2|\langle S_{1,j},f_j^{[j]}\rangle|
\le(c\nu/2)\|f_j^{[j]}\|_{\Sigma_0}^2
+C\nu^{-1}\|S_{1,j}\|_{\Sigma_0^*}^2$.
The small integrated coefficient from Step 1 implies
\begin{equation}\label{eq:forwardLocalizedSourceTrace}
\begin{aligned}
\sup_{I'}\|f_k^{[j]}\|_2
+\left(\nu\int_{I'}\|f_k^{[j]}\|_{\Sigma_0}^2\,ds\right)^{1/2}
\le Ce^{-\varrho_0|k-j|^\chi}\big(&\|f_{{\rm in},j}\|_2
+\|S_{0,j}\|_{L^1(I';L_v^2)}\\
&+\nu^{-1/2}\|S_{1,j}\|_{L^2(I';\Sigma_0^*)}\big).
\end{aligned}
\end{equation}
{The energy construction in Lemma~\ref{lem:localizedConstruction}
justifies these estimates for smooth initial data $f_{\rm in}$ and sources $S_0,S_1$. The passage to the
completed source spaces is given in Step 4 below.}

\emph{3. Density summation.}
By uniqueness $f=\sum_j f^{[j]}$.
Pair \eqref{eq:forwardLocalizedSourceTrace} with $\sqrt\mu$ and use
\[
\langle k\rangle^{r-1}\le C_r\langle j\rangle^{r-1}\langle k-j\rangle^{r-1},
\qquad \sum_k\langle k-j\rangle^{r-1}e^{-\varrho_0|k-j|^\chi}\le C_r.
\]
Thus
\[
\begin{aligned}
&\sum_{k\ne0}\langle k\rangle^{r-1}\sup_{I'}
\left|\int\sqrt\mu f_k(s)\,dv\right|\\
&\quad\le C_r\sum_j\langle j\rangle^{r-1}
\left(\|f_{{\rm in},j}\|_2+\int_{I'}\|S_{0,j}(s)\|_2\,ds
+\nu^{-1/2}\|S_{1,j}\|_{L^2(I';\Sigma_0^*)}\right).
\end{aligned}
\]
For the first two terms, Cauchy--Schwarz gives
\[
\begin{aligned}
&\sum_j\langle j\rangle^{r-1}
\left(\|f_{{\rm in},j}\|_2+\int_{I'}\|S_{0,j}(s)\|_2\,ds\right)\\
&\quad\le
\left(\sum_j\langle j\rangle^{-2}\right)^{1/2}
\left[
\left(\sum_j\langle j\rangle^{2r}\|f_{{\rm in},j}\|_2^2\right)^{1/2}
+\int_{I'}\left(\sum_j\langle j\rangle^{2r}\|S_{0,j}(s)\|_2^2\right)^{1/2}ds
\right]\\
&\quad\le C_r\big(\|f_{\rm in}\|_{H_x^rL_v^2}
+\|S_0\|_{L^1(I';H_x^rL_v^2)}\big).
\end{aligned}
\]
For the form source, Lemma~\ref{lem:weightedDualMultipliers} gives
\[
\|\langle v\rangle^{-2\ell}\varphi\|_{\Sigma_\ell}
\le C_\ell\|\varphi\|_{\Sigma_0},\qquad
\|S_{1,j}\|_{\Sigma_0^*}\le C_\ell\|S_{1,j}\|_{\Sigma_\ell^*}.
\]
Hence
\[
\begin{aligned}
\sum_j\langle j\rangle^{r-1}\|S_{1,j}\|_{L^2\Sigma_0^*}
&\le C\left(\sum_j\langle j\rangle^{2r}\|S_{1,j}\|_{L^2\Sigma_0^*}^2\right)^{1/2}\\
&\le C_r\left(\sum_{a=0}^r
\|\partial_{x_1}^aS_1\|_{L^2(I';L_x^2\Sigma_{2r+6-2a}^*)}^2\right)^{1/2}
\le C_r\sum_{a=0}^{r}\|\partial_{x_1}^aS_1\|_{L^2_{s,x}(I')\Sigma_{2r+6-2a}^*}.
\end{aligned}
\]
This proves \eqref{eq:collisionCombinedTrace} with the time supremum
inside the mode sum. For every real $q\ge0$,
\[
\sum_k\langle k\rangle^qe^{-\varrho_0|k-j|^\chi}
\le C_q\langle j\rangle^q.
\]
With zero sources, summing \eqref{eq:forwardLocalizedSourceTrace}
using this bound gives the first inequality in
\eqref{eq:coarseCollisionTrace}. Cauchy--Schwarz gives the second:
\[
\sum_j\langle j\rangle^q\|S_j\|_2
\le\left(\sum_j\langle j\rangle^{-2}\right)^{1/2}
\left(\sum_j\langle j\rangle^{2q+2}\|S_j\|_2^2\right)^{1/2}
\le C_q\|S\|_{H_x^{q+1}L_v^2}.
\]
Taking $f_{\rm in}=S_0=0$
in \eqref{eq:collisionCombinedTrace} gives \eqref{eq:dualCollisionTrace}.

\emph{4. Unweighted energy and source completion.}
The inequalities
\[
\langle k\rangle^r\le C_r\langle j\rangle^r\langle k-j\rangle^r,
\qquad
\sum_{h\in\mathbb Z}\langle h\rangle^re^{-\varrho_0|h|^\chi}<\infty
\]
allow discrete Young's inequality to be applied to
\eqref{eq:forwardLocalizedSourceTrace}. Minkowski's inequality in time
and the form norm gives
\[
\begin{aligned}
&\left(\sum_k\langle k\rangle^{2r}
\sup_{I'}\|f_k(s)\|_2^2\right)^{1/2}
+\left(\nu\sum_k\langle k\rangle^{2r}
\int_{I'}\|f_k(s)\|_{\Sigma_0}^2\,ds\right)^{1/2}\\
&\quad\le C_r\left[
\left(\sum_j\langle j\rangle^{2r}\|f_{{\rm in},j}\|_2^2\right)^{1/2}
+\left(\sum_j\langle j\rangle^{2r}
\left(\int_{I'}\|S_{0,j}(s)\|_2\,ds\right)^2\right)^{1/2}
\right.\\
&\hspace{35mm}\left.
+\nu^{-1/2}
\left(\sum_j\langle j\rangle^{2r}
\|S_{1,j}\|_{L^2(I';\Sigma_0^*)}^2\right)^{1/2}
\right].
\end{aligned}
\]
The strong-source term satisfies
\[
\left(\sum_j\langle j\rangle^{2r}
\left(\int_{I'}\|S_{0,j}(s)\|_2\,ds\right)^2\right)^{1/2}
\le\int_{I'}\left(\sum_j\langle j\rangle^{2r}
\|S_{0,j}(s)\|_2^2\right)^{1/2}ds.
\]
The form-source bound in Step 3 controls the last term. Fourier
norm equivalence therefore proves \eqref{eq:collisionUnweightedEnergy}.

Choose smooth data and sources
$f_{{\rm in},\epsilon},S_{0,\epsilon},S_{1,\epsilon}$
{such that
\[
\|f_{{\rm in},\epsilon}-f_{\rm in}\|_{H_x^rL_v^2}
+\|S_{0,\epsilon}-S_0\|_{L^1(I';H_x^rL_v^2)}
+\sum_{j=0}^r\|\partial_{x_1}^j(S_{1,\epsilon}-S_1)\|_
{L^2_{s,x}(I')\Sigma_{2r+6-2j}^*}\longrightarrow0.
\]
For their solutions, apply}
\eqref{eq:collisionUnweightedEnergy} and
\eqref{eq:collisionCombinedTrace} to differences:
\[
\begin{aligned}
&\sup_{I'}\|f_\epsilon-f_{\epsilon'}\|_{H_x^rL_v^2}
+\left(\nu\sum_{a=0}^r\int_{I'}
\|\partial_{x_1}^a(f_\epsilon-f_{\epsilon'})\|_{L_x^2\Sigma_0}^2\,ds\right)^{1/2}\\
&\quad+\sum_{k\ne0}\langle k\rangle^{r-1}\sup_{I'}
\left|\int\sqrt\mu(f_{\epsilon,k}-f_{\epsilon',k})\,dv\right|\\
&\le C_r\big(\|f_{{\rm in},\epsilon}-f_{{\rm in},\epsilon'}\|_{H_x^rL_v^2}
+\|S_{0,\epsilon}-S_{0,\epsilon'}\|_{L^1(I';H_x^rL_v^2)}\\
&\hspace{35mm}
+\nu^{-1/2}\sum_{j=0}^{r}\left\|\partial_{x_1}^j\!\left(S_{1,\epsilon}-S_{1,\epsilon'}\right)\right\|_{L^2_{s,x}(I')\Sigma_{2r+6-2j}^*}\big)
\longrightarrow0.
\end{aligned}
\]
The limit belongs to $C(I';H_x^rL_v^2)$ and has the asserted
integrated form norm. The collision terms converge in the weak
equation by the collision-form bounds; the sources converge against
smooth test functions by the displayed norms of $f_{\rm in}$, $S_0$ and $S_1$. Moreover,
\[
\sup_{I'}\left|\int\sqrt\mu(f_{\epsilon,k}-f_k)\,dv\right|
\le\|\sqrt\mu\|_2\sup_{I'}\|f_{\epsilon,k}-f_k\|_2\longrightarrow0.
\]
Thus the limit in $\ell^1_{r-1}C(I')$ is precisely the density of $f$.
The difference of two energy solutions with the same initial datum and sources has
zero initial datum and zero source. Uniqueness follows from
Lemma~\ref{lem:localizedConstruction} with $\ell=0$ and
the polynomial spatial weight $\kappa_k=\langle k\rangle^r$.
For $f_{\rm in}$, $S_0$ and $S_1$ in the stronger spaces of the lemma, choose the
approximations to converge also in $\mathscr E_r$ and
$L^1(I';\mathscr E_r)$. Step 1 then gives
\[
\begin{aligned}
\|f_\epsilon-f_{\epsilon'}\|_{\mathcal N_r^{I'}}
\le C_r\big(&\mathcal E_r[f_{{\rm in},\epsilon}-f_{{\rm in},\epsilon'}]^{1/2}
+\int_{I'}\mathcal E_r[S_{0,\epsilon}-S_{0,\epsilon'}]^{1/2}\,ds\\
&+\nu^{-1/2}\sum_{j=0}^{r}\left\|\partial_{x_1}^j\!\left(S_{1,\epsilon}-S_{1,\epsilon'}\right)\right\|_{L^2_{s,x}(I')\Sigma_{2r+6-2j}^*}\big)
\longrightarrow0.
\end{aligned}
\]
The two limits coincide in $C(I';H_x^rL_v^2)$, so the completed
Duhamel integral is the energy solution constructed in Step 1.
\par
\end{proof}
\par
\par\endgroup
% Native TikZ road map: no external image files, colors, or version identifiers.
\clearpage
\begingroup
\newcommand{\RMbg}{F^{\rm b,\nu}}
\newcommand{\RML}{\mathcal L_{\rm C}}
\newcommand{\RMA}{\mathcal A_\nu}
\newcommand{\RMT}{\mathcal T_{\rm b}}
\newcommand{\RMQ}{\mathbb Q_\nu}
\newcommand{\RMeq}[1]{\overset{\text{\scriptsize\eqref{#1}}}{=}}
\newcommand{\RMhead}[1]{{\sffamily\bfseries\fontsize{9}{10.5}\selectfont #1}\par\smallskip}
\newcommand{\RMnote}[1]{{\sffamily\fontsize{7}{8.5}\selectfont #1}}
\tikzset{
 RMcard/.style={draw=black,fill=white,rounded corners=1.5mm,line width=.45pt,
 inner sep=1.5mm,align=left,anchor=north,text width=7cm,
 font=\fontsize{8}{9.5}\selectfont,text=black},
 RMwide/.style={RMcard,text width=22.7cm},
 RMarr/.style={-{Stealth[length=1.7mm]},line width=.55pt,draw=black}}
\begin{landscape}
\section{Road map of the proof}\label{app:roadMap}
Figures~\ref{fig:roadMapDistributions}--\ref{fig:roadMapDensityComparison}
collect the density decompositions and the corresponding distribution equations;
$s=t+t_0\in[t_0,T_1]$ and $\rho[H]=\int_{\mathbb R^3}H\,dv$.
\begin{figure}[!ht]
\centering
\resizebox{.90\linewidth}{!}{%
\begin{tikzpicture}[x=1cm,y=1cm]
\node[RMwide] (root) at (11.75,0) {
\RMhead{Density decomposition of the exact nonlinear solution}
\[
\begin{gathered}
\rho[F^\nu-\mu]\RMeq{eq:nonlinearFirstOrderExpansion}
\rho[\RMbg-\mu]+\delta_\nu\rho[H^\nu]+\rho[W],\\
\partial_sF^\nu+v_1\partial_{x_1}F^\nu+E[F^\nu]\partial_{v_1}F^\nu
\RMeq{eq:vplPhysicalTime}\nu Q(F^\nu,F^\nu).
\end{gathered}
\]
};
\node[RMcard,minimum height=2.6cm] (bg) at ($(root.south)+(-7.85,-.35)$) {
\RMhead{Background density $\rho[\RMbg-\mu]$}
\[
\begin{aligned}
(\partial_s+v_1\partial_{x_1})\RMbg
&+E[\RMbg]\partial_{v_1}\RMbg\\
&\RMeq{eq:exactVPLbackground}\nu Q(\RMbg,\RMbg),\\
\RMbg(t_0)&\RMeq{eq:maininitial}\mu_\theta(1+a\cos(x_1-t_0v_1)).
\end{aligned}
\]
\RMnote{Section~\ref{sec:referenceEstimates}: construct the background and estimate its electric field and coefficients.}
};
\node[RMcard,minimum height=2.6cm] (lin) at ($(root.south)+(0,-.35)$) {
\RMhead{Complete linear density $\rho[H^\nu]$}
\[
\begin{gathered}
\RMA H^\nu\RMeq{eq:nonlinearFirstOrderExpansion}0,
\qquad H^\nu(t_0)\RMeq{eq:fullRangePhysicalEquation}H_{\rm in},\\
H_{\rm in}\RMeq{eq:maininitial}\frac{\mu_{1/2}}{\sqrt\pi}
\cos\!\big(m_0x_1+(\zeta-m_0t_0)v_1\big).
\end{gathered}
\]
\RMnote{Section~\ref{sec:completeLinear}: compare the complete linear density with the growing density $\rho^\downarrow$.}
};
\node[RMcard,minimum height=2.6cm] (nonlin) at ($(root.south)+(7.85,-.35)$) {
\RMhead{Nonlinear error density $\rho[W]$}
\[
\begin{aligned}
\RMA W&\RMeq{eq:nonlinearExpansionRoadmap}\delta_\nu^2\RMQ(H^\nu,H^\nu)\\
&\quad+\RMQ(\delta_\nu H^\nu,W)+\RMQ(W,\delta_\nu H^\nu)\\
&\quad+\RMQ(W,W),\qquad W(t_0)\RMeq{eq:roadmapInitialDecomposition}W_{\rm in}.
\end{aligned}
\]
\RMnote{Section~\ref{sec:nonlinear}: $W_{\rm in}$ is the initial temperature correction.}
};
\coordinate (topbranch) at ($(root.south)!.5!(lin.north)$);
\draw[RMarr] (root.south) -- (topbranch) -| (bg.north);
\draw[RMarr] (root.south) -- (lin.north);
\draw[RMarr] (root.south) -- (topbranch) -| (nonlin.north);
\node[fit=(bg)(lin)(nonlin),inner sep=0pt] (rowone) {};
\node[RMwide] (split) at ($(rowone.south)+(0,-.35)$) {
\RMhead{Decomposition of the complete linear density}
\[
\rho[H^\nu]\RMeq{eq:fieldFiniteUnknownError}\rho[H^{[J]}]+\rho[W_J]
\RMeq{eq:fieldPartialSumIdentities}\rho[H_0]+\sum_{n=1}^J\rho[H_n]+\rho[W_J].
\]
};
\draw[RMarr] (lin.south) -- (split.north);
\node[RMcard,minimum height=1.9cm] (h0) at ($(split.south)+(-7.85,-.35)$) {
\RMhead{Main density $\rho[H_0]$}
\[
\begin{gathered}
\rho[H_0]\RMeq{eq:exactCollisionModel}\rho^{\rm C,\nu},\\
\RML H_0\RMeq{eq:exactCollisionModel}0,\qquad
H_0(t_0)\RMeq{eq:exactCollisionModel}H_{\rm in}.
\end{gathered}
\]
};
\node[RMcard,minimum height=1.9cm] (hn) at ($(split.south)+(0,-.35)$) {
\RMhead{Correction densities $\rho[H_n]$, $1\le n\le J$}
\[
\begin{gathered}
\rho[H_n]\RMeq{eq:fieldDistributionRecursion}\rho^{(n)},\\
\RML H_n\RMeq{eq:fieldCorrectionRecursion}-\RMT H_{n-1},
\qquad H_n(t_0)\RMeq{eq:fieldCorrectionRecursion}0.
\end{gathered}
\]
};
\node[RMcard,minimum height=1.9cm] (wj) at ($(split.south)+(7.85,-.35)$) {
\RMhead{Linear error density $\rho[W_J]$}
\[
\begin{gathered}
\RMA W_J\RMeq{eq:fieldFiniteUnknownError}-\RMT H_J,\\
W_J(t_0)\RMeq{eq:fieldFiniteUnknownError}0.
\end{gathered}
\]
};
\coordinate (branch) at ($(split.south)!.5!(hn.north)$);
\draw[RMarr] (split.south) -- (branch) -| (h0.north);
\draw[RMarr] (split.south) -- (hn.north);
\draw[RMarr] (split.south) -- (branch) -| (wj.north);
\node[fit=(h0)(hn)(wj),inner sep=0pt] (rowtwo) {};
\node[RMwide] at ($(rowtwo.south)+(0,-.35)$) {
\RMhead{Distribution operators}
\[
\begin{aligned}
\RML H&\RMeq{eq:fieldExpansionPhysicalOperator}
\partial_sH+v_1\partial_{x_1}H+E[H]\partial_{v_1}\RMbg
-\nu\big(Q(\RMbg,H)+Q(H,\RMbg)\big),\\
\RMT H&\RMeq{eq:fullPhysicalOperatorDefinition}E[\RMbg]\partial_{v_1}H,
\quad \RMA\RMeq{eq:fullPhysicalOperatorDefinition}\RML+\RMT,
\quad \RMQ(h,w)\RMeq{eq:nonlinearBilinearOperator}-E[h]\partial_{v_1}w+\nu Q(h,w).
\end{aligned}
\]
};
\end{tikzpicture}%
}
\caption{Density decompositions and the corresponding distribution equations.}
\label{fig:roadMapDistributions}
\end{figure}
\clearpage
\begin{figure}[!ht]
\centering
\resizebox{.90\linewidth}{!}{%
\begin{tikzpicture}[x=1cm,y=1cm]
\node[RMwide] (root) at (11.75,0) {
\RMhead{The complete linear density: growing term and five errors}
\[
\rho[H^\nu]\RMeq{eq:linearConstructionThreeErrors}\rho^\downarrow
+\mathbb V^{\rm C,\nu}(\mathfrak g^{\rm C,\nu}-\mathfrak g^0)
+\mathbb V^{\rm C,\nu}\Delta\mathbb T\rho^\downarrow
+\mathbb V^{\rm C,\nu}\mathbb P_\uparrow\rho^\downarrow
+\sum_{n=1}^J\rho[H_n]+\rho[W_J].
\]
};
\node[RMcard,minimum height=6cm] (main) at ($(root.south)+(-7.85,-.35)$) {
\RMhead{Main density and the descending chain}
\[
\begin{aligned}
\rho^\downarrow&\RMeq{eq:overviewDescendingDensityEquation}\mathfrak g^0+
(\mathbb K+\mathbb P_\downarrow)\rho^\downarrow,\\
\rho^{\rm C,\nu}&\RMeq{eq:directChainComparisonIdentity}\rho^\downarrow\\
&\quad+\mathbb V^{\rm C,\nu}(\mathfrak g^{\rm C,\nu}-\mathfrak g^0)\\
&\quad+\mathbb V^{\rm C,\nu}\Delta\mathbb T\rho^\downarrow
+\mathbb V^{\rm C,\nu}\mathbb P_\uparrow\rho^\downarrow.
\end{aligned}
\]
\[
\begin{gathered}
\mathfrak g^0_{m_0}(s)\RMeq{eq:freeGaussianDensity}
\frac{e^{-(m_0s-\zeta)^2/4}}{\sqrt{4\pi}},\\
|\rho^\downarrow_{m_1}(T_1)|\gtrsim\frac{B_*}{\sqrt N}.
\end{gathered}
\]
\RMnote{Propositions~\ref{thm:chain} and~\ref{thm:exactCollisionEcho}: compare the initial-data evolutions, the collision and background coefficients, and the upward interaction. Each error is $o(B_*/\sqrt N)$ at $(m_1,T_1)$.}
};
\node[RMcard,minimum height=6cm] (corr) at ($(root.south)+(0,-.35)$) {
\RMhead{Correction densities and ordered propagation}
\[
\begin{aligned}
\rho[H_n]&\RMeq{eq:fieldDistributionRecursion}\rho^{(n)}\\
&\RMeq{eq:fieldDensityRecursion}\mathbb V^{\rm C,\nu}
\left(\mathfrak g_n+\sum_{q=1}^n\mathbb T_q\rho^{(n-q)}\right).
\end{aligned}
\]
\[
\begin{aligned}
\frac{H_n(s)}{\sqrt\mu}
&\RMeq{eq:fieldDistributionRecursion}\mathbb U_{\rm b}^{(n)}(s,t_0)
\frac{H_{\rm in}}{\sqrt\mu}\\
&\quad+\sum_{q=0}^n\sum_{m\ne0}\int_{t_0}^s
\mathbb U_{\rm b}^{(q)}(s,\tau)
\mathcal F_m^{\rm b}(\tau)\rho_m^{(n-q)}(\tau)\,d\tau.
\end{aligned}
\]
\[
\frac{\sqrt N}{B_*}|\rho_{m_1}[H_n](T_1)|
\le\nu^{n\sigma_{\rm int}-o(1)}.
\]
\RMnote{Lemma~\ref{lem:oneFieldOperation} and Proposition~\ref{prop:finiteFieldTrace}: retain the insertion times and frequency shifts. Use the finite hierarchy of spatial orders $r_{\rm x}+J-n$.}
};
\node[RMcard,minimum height=6cm] (err) at ($(root.south)+(7.85,-.35)$) {
\RMhead{Linear and nonlinear error densities}
\[
\RMA W_J\RMeq{eq:fieldFiniteUnknownError}-\RMT H_J,
\qquad W_J(t_0)\RMeq{eq:fieldFiniteUnknownError}0.
\]
\RMnote{First estimate the source $-\RMT H_J$. Proposition~\ref{thm:fullForcedCoarse} then gives the density and energy bounds for $W_J$. Choose fixed finite $J$ to compensate for the loss $e^{C_*N}$.}
\[
|\rho_{m_1}[W_J](T_1)|=o(B_*/\sqrt N).
\]
\smallskip
\RMhead{Nonlinear comparison}
\RMnote{Use the equation for $W$ in Figure~\ref{fig:roadMapDistributions}. Its known quadratic source is $\delta_\nu^2\mathbb Q_\nu(H^\nu,H^\nu)$, and $W(t_0)=W_{\rm in}$. Section~\ref{sec:nonlinear} controls the remaining nonlinear terms.}
\[
|\rho_{m_1}[W](T_1)|=o(\delta_\nu B_*/\sqrt N).
\]
};
\coordinate (branch) at ($(root.south)!.5!(corr.north)$);
\draw[RMarr] (root.south) -- (branch) -| (main.north);
\draw[RMarr] (root.south) -- (corr.north);
\draw[RMarr] (root.south) -- (branch) -| (err.north);
\node[fit=(main)(corr)(err),inner sep=0pt] (rowthree) {};
\node[RMwide] at ($(rowthree.south)+(0,-.35)$) {
\RMhead{Normalized transport and ordered solution operators}
\[
\begin{aligned}
\widetilde{\mathcal T}_{\rm b}(s)f
&\RMeq{eq:fieldNormalizedOperator}-E[\RMbg](s)(\partial_{v_1}-v_1/2)f,\\
\mathbb U_{\rm b}^{(0)}(s,\tau)&\RMeq{eq:fieldEvolutionCoefficients}\mathbb U_{\rm b}^\nu(s,\tau),\qquad
\mathbb U_{\rm b}^{(q)}(s,\tau)f\RMeq{eq:fieldEvolutionCoefficients}
\int_\tau^s\mathbb U_{\rm b}^\nu(s,u)\widetilde{\mathcal T}_{\rm b}(u)
\mathbb U_{\rm b}^{(q-1)}(u,\tau)f\,du.
\end{aligned}
\]
\RMnote{The operators $\mathbb U_{\rm b}^\nu$ and $\mathbb U_{\rm b}^{(q)}$ act on normalized functions; $\mathbb T_q$ acts on densities.}
};
\end{tikzpicture}%
}
\caption{Density decompositions and the corresponding distribution equations: density comparison and finite corrections.}
\label{fig:roadMapDensityComparison}
\end{figure}
\end{landscape}
\endgroup
\clearpage

\bigbreak

\noindent{\bf Acknowledgments}
 Ruizhao Zi was   supported by NSF of China under  Grant 12571235.

\noindent{\bf Conflict of interest}  The authors declare that they have no conflict of interest. 

\noindent{\bf Data availability} Data sharing is not applicable to this article as no datasets
were generated or analyzed during the current study.

\noindent{\bf Disclosure of AI assistance} Some of the key ideas underlying the proofs in this paper emerged from discussions between the authors and an AI system(ChatGPT 6.0). The author also used the AI tools to help revise the manuscript. The authors take full responsibility for the accuracy of the entire paper, including the correctness of all mathematical statements and proofs.

\par
\end{document}